\documentclass[11pt,leqno]{amsart}

\usepackage{amsmath}
\usepackage{amssymb}
\usepackage{a4wide}
\usepackage{bbm}
\usepackage{graphics}
\usepackage{epsfig}
\usepackage{color}
\usepackage{mathrsfs}
\usepackage[sans]{dsfont}

\usepackage[breaklinks]{hyperref}

\newcommand{\Subsection}[1]{\subsection{ #1} ${}^{}$}

\newtheorem{theorem}{Theorem}[section]
\newtheorem{lemma}[theorem]{Lemma}
\newtheorem{proposition}[theorem]{Proposition}
\newtheorem{definition}[theorem]{Definition}
\newtheorem{remark}[theorem]{Remark}
\newtheorem{corollary}[theorem]{Corollary}
\newtheorem{example}[theorem]{Example}

\newcounter{hypo}
\newcounter{hypoa}

\newenvironment{hyp}{  \begin{enumerate} \setcounter{enumi}{\value{hypo}} \item}{\stepcounter{hypo} \end{enumerate}}

\def\C{{\mathbb C}}

\def\N{{\mathbb N}} 
\def\R{{\mathbb R}} 

\def\Z{{\mathbb Z}}

\def\S{{\mathbb S}}

\def\CC{\mathcal {C}}
\def\CA{\mathcal {A}}
\def\CB{\mathcal {B}}
\def\CD{\mathcal {D}}
\def\CE{\mathcal {E}}

\def\CH{\mathcal {H}}
\def\CI{{\mathcal I}}
\def\CJ{{\mathcal J}}
\def\CK{\mathcal {K}}
\def\CL{\mathcal {L}}
\def\CM{\mathcal {M}}
\def\CN{\mathcal {N}}
\def\CO{\mathcal {O}}
\def\CP{\mathcal {P}}
\def\CQ{\mathcal {Q}}
\def\CR{\mathcal {R}}
\def\CS{\mathcal {S}}
\def\CT{\mathcal {T}}
\def\CU{\mathcal {U}}
\def\CV{\mathcal {V}}

\def\CZ{\mathcal {Z}}

\def\SA{\mathscr {A}}
\def\SB{\mathscr {B}}
\def\SC{\mathscr {C}}
\def\SD{\mathscr {D}}
\def\SE{\mathscr {E}}
\def\SF{\mathscr {F}}
\def\SG{\mathscr {G}}

\def\SJ{\mathscr {J}}

\def\SM{\mathscr {M}}

\def\SP{\mathscr {P}}
\def\SQ{\mathscr {Q}}
\def\SR{\mathscr {R}}

\def\SV{\mathscr {V}}

\def\ker{\mathop{\rm Ker}\nolimits}

\def\one{\mathds{1}}
\def\re{\mathop{\rm Re}\nolimits}
 \def\im{\mathop{\rm Im}\nolimits}

\def\Op{\mathop{\rm Op}\nolimits}

\newcommand{\tr}{\operatorname{tr}}
\newcommand{\ang}{\operatorname{angle}}
\newcommand{\card}{\operatorname{card}}

\newcommand{\ag}{\operatorname{Ag}}

\newcommand{\adj}{\operatorname{adj}}

\newcommand{\mes}{\operatorname{mes}}
\newcommand{\supp}{\operatorname{supp}}
\newcommand{\spe}{\operatorname{sp}}
\newcommand{\spr}{\operatorname{spr}}
 
\def\arg{\mathop{\rm arg}\nolimits}
\def\dist{\mathop{\rm dist}\nolimits}
\def\diag{\mathop{\rm diag}\nolimits}
\def\sgn{\mathop{\rm sgn}\nolimits}
\def\Hess{\mathop{\rm Hess}\nolimits}
\def\ad{\mathop{\rm ad}\nolimits}
\def\<{\langle}
\def\>{\rangle}

\def\rank{\mathop{\rm Rank}\nolimits}

\def\res{\mathop{\rm Res}\nolimits}

\def\Res{{\rm Res}}

\def\ds{\displaystyle}

\newcommand{\fract}[2]{\genfrac{}{}{0pt}{}{\scriptstyle #1}{\scriptstyle #2}}

\makeatletter
 \@addtoreset{equation}{section}
 \makeatother
 
\makeatletter
\def\@tocline#1#2#3#4#5#6#7{\relax
  \ifnum #1>\c@tocdepth 
  \else
    \par \addpenalty\@secpenalty\addvspace{#2}%
    \begingroup \hyphenpenalty\@M
    \@ifempty{#4}{%
      \@tempdima\csname r@tocindent\number#1\endcsname\relax
    }{%
      \@tempdima#4\relax
    }%
    \parindent\z@ \leftskip#3\relax \advance\leftskip\@tempdima\relax
    \rightskip\@pnumwidth plus4em \parfillskip-\@pnumwidth
    {\ifnum #1<2\medskip\fi}#5\leavevmode\hskip-\@tempdima
      \ifcase #1
       \or\or \hskip 2em \or \hskip 2em \else \hskip 3em \fi%
       {\ifnum #1<2 \bf #6 \else #6\fi}\nobreak\relax
    {\ifnum #1<2 \hfill\else\dotfill\fi}\hbox to\@pnumwidth{{\ifnum
#1<2\null\qquad\@tocpagenum{\textbf{#7}}\else\@tocpagenum{#7}\fi}}\par
    \nobreak
    \endgroup
  \fi}
\makeatother

\title{Resonances for homoclinic trapped sets}

\author[J.-F. Bony]{Jean-Fran\c{c}ois Bony}
\address{Jean-Fran\c{c}ois Bony, IMB, CNRS (UMR 5251), Universit\'e de Bordeaux, 33405 Talence, France}
\email{bony@math.u-bordeaux.fr}
\author[S. Fujii\'e]{Setsuro Fujii\'e}
\address{Setsuro Fujii\'e, Department of Mathematical Sciences, Ritsumeikan University, 1-1-1 Noji-Higashi, Kusatsu, 525-8577 Japan}
\email{fujiie@fc.ritsumei.ac.jp}
\author[T. Ramond]{Thierry Ramond}
\address{Thierry Ramond, Laboratoire de Math\'ematiques d'Orsay, Univ. Paris-Sud, CNRS, Universit\'e Paris-Saclay, 91405 Orsay, France}
\email{thierry.ramond@math.u-psud.fr}
\author[M. Zerzeri]{Maher Zerzeri}
\address{Maher Zerzeri, Universit\'e Paris 13, Sorbonne Paris Cit\'e, LAGA, CNRS (UMR 7539), 93430 Villetaneuse, France}
\email{zerzeri@math.univ-paris13.fr}


\thanks{\textbf{Acknowledgments:} This work was partially supported by the JSPS KAKENHI Grant 15K04971 and the ANR project NOSEVOL 2011 BS 010119-01. The third and fourth authors would like to thank the mathematical department of Ritsumeikan University for its kind hospitality.}

\begin{document}

\begin{abstract}
We study semiclassical resonances generated by homoclinic trapped sets. First, under some general assumptions, we prove that there is no resonance in a region below the real axis. Then, we obtain a quantization rule and the asymptotic expansion of the resonances when there is a finite number of homoclinic trajectories. The same kind of results is proved for homoclinic sets of maximal dimension. Next, we generalize to the case of homoclinic/heteroclinic trajectories and we study the three bump case. In all these settings, the resonances may either accumulate on curves or form clouds. We also describe the corresponding resonant states.
\end{abstract}

\maketitle

\section{Introduction}

In this work, we are interested in quantum resonances, that is poles of the resolvent meromorphically extended through the real axis of self-adjoint operators $P$. Mostly, we consider the semiclassical Schr\"odinger operator $P = - h^{2} \Delta + V$ on $L^{2} ( \R^{n} )$ for potentials $V$ decaying at infinity, but the discussion often applies to more general pseudodifferential operators provided they are close enough to the Laplace operator $- h^{2} \Delta$ out of a compact set.

In physics, the notion of quantum resonances goes back to the very beginnings of quantum mechanics (see e.g. Gamow \cite{Ga28_01}). It was proposed to explain results in scattering experiments, and it is associated to the notion of pseudo-particles with finite lifetime. To some extent, it is reasonable to think of resonances as complex numbers which are generalized eigenvalues of quantum observable. The imaginary part of a resonance is interpreted as the inverse of the lifetime of the pseudo-particle it is associated to. Thus, the closer they are to the real axis, the more meaningful they should be.

We consider here resonances that are close to the real axis, namely with imaginary part of size $\CO ( h )$, in the semiclassical regime $h \to 0$. It is well known that they are related to the existence of bounded classical trajectories for the Hamiltonian vector field $H_{p}$ on $T^{*} \R^{n}$. Here, $p$ denotes the classical observable corresponding to, or the symbol of the operator $P$. The set $K ( E_{0} )$ of bounded trajectories in $p^{- 1} ( E_{0} )$ is usually called the {\it trapped set} at energy $E_{0} > 0$. For example, it has been proved by Martinez \cite{Ma02_01} in the $C^{\infty}$ setting that for a given real $E_{0}$ such that the trapped set at energy $E_{0}$ is empty, there are no resonance is any complex disk of size $h \vert \ln h \vert$ around $E_{0}$ (see also Va{\u\i}nberg \cite{Va89_01} in the classical case).

As far as the asymptotic of resonances is concerned, trapped sets of different kinds, yet relatively simple, have been considered. First of all, Helffer and Sj\"ostrand \cite{HeSj86_01} have studied the so-called ``well in the island'' situation, where the potential $V$ have the shape of a well at energy $E_{0}$. They have proved that the resonances near $E_{0}$ are exponentially close to the eigenvalues of a corresponding confining configuration, for example a suitable Dirichlet realization $P_{D}$ of the differential operator $P$. Therefore, their imaginary part is exponentially small. Moreover, in the case of a punctual well, that is when $V$ has a non-degenerate local minimum with critical value $E_{0}$, they have given the asymptotic with respect to $h$ of the imaginary part of these resonances.

When the potential $V$ has a general non-degenerate critical point with critical value $E_{0}$, Sj\"ostrand \cite{Sj87_01} has obtained the asymptotic of the resonances close to $E_{0}$. In particular, in the case of a maximum, they satisfy
\begin{equation} \label{n1}
z = E_{0} - i h \sum_{j = 1}^{n} \Big( \alpha_{j} + \frac{1}{2} \Big) \lambda_{j} + o ( h ) ,
\end{equation}
for some $\alpha \in \N^{n}$. Here, the $\lambda_{j}$'s are the positive eigenvalues of the linearized Hamiltonian vector field $H_{p}$ at the critical point. This result has also been proved by Briet, Combes and Duclos \cite{BrCoDu87_01,BrCoDu87_02} (see \cite{Ra96_01} in dimension $1$). The interested reader may find in \cite{BoFuRaZe11_01} a precise description of the spectral projection associated to these resonances.

Then, G\'{e}rard and Sj\"ostrand \cite{GeSj87_01} have studied the case where the trapped set at energy $E_0$ consists of a single closed orbit of hyperbolic type. They exhibit a Bohr--Sommerfeld like quantization rule along the orbit, and show that there is a bijection between the set of resonances and the set of roots of the quantization rule. Up to lower order terms, the resonances in that case form a lattice at distance $h$ of the real axis. Let us mention that resonances for hyperbolic situations have also been studied in slightly different contexts. They have been first investigated by Ikawa \cite{Ik83_01} and G\'{e}rard \cite{Ge88_01} in the case of scattering by obstacles. There is also a wide related literature in general relativity (see e.g. S{\'a} Barreto and Zworski \cite{SaZw97_01}) as well as for hyperbolic surfaces (see e.g. Dyatlov, Faure and Guillarmou \cite{DyFaGu14_01}). For weaker trappings, we would like to mention the results by Burq \cite{Bu97_01} for resonances generated by a corner and by Alexandrova and Tamura \cite{AlTa14_01} in the context of Aharonov--Bohm effect.

\Subsection{Brief description of the results}

Here, we consider homoclinic situations. More precisely, we assume that the trapped set at energy $E_{0}$ consists of a hyperbolic fixed point, say $0 \in T^{*} \R^{n}$, and homoclinic trajectories. We recall that a non-constant integral curve of a vector field is called homoclinic when it converges to a fixed point of the vector field as the time tends to both $\pm \infty$. Some configurations of this kind have already been studied in dimension $1$ by the second and third authors \cite{FuRa98_01,FuRa03_01}, and by Servat \cite{Se04_01}. Let $\Lambda_{\pm}$ denote the incoming/outgoing Lagrangian manifolds associated to the fixed point $0$. The homoclinic set is then $\CH = \Lambda_{-} \cap \Lambda_{+} \setminus \{ 0 \}$. We also denote $0 < \lambda_{1} \leq \cdots \leq \lambda_{n}$ the positive eigenvalues of the linearization of the Hamiltonian vector field at $0$.

We first obtain resonance free domains under fairly general assumptions. Besides their intrinsic interest, these results should be seen as a preamble to the asymptotic of the resonances. The domains depend on the $\lambda_{j}$'s and on the nature of the contact between $\Lambda_{-}$ and $\Lambda_{+}$ along $\CH$. This last notion is more relevant than the dimension of the homoclinic set, as can be seen comparing previous results for other kinds of trappings. We distinguish between two cases, the precise corresponding results being stated in Section \ref{s22}. Note that, if we do not add any hypothesis, there may exist resonances exponentially close to the real axis.

We assume first that the hyperbolic fixed point is anisotropic (i.e. $\lambda_{1} < \lambda_{n}$) or that $\Lambda_{-}$ and $\Lambda_{+}$ intersect transversally in at least one direction. In the last case, the dimension of $\CH$ is at most $n - 1$. In this weak trapping situation, we show that $P$ has no resonance and that its truncated resolvent has a polynomial estimate in a neighborhood of size $h$ of the real axis. Such results have already been obtained for hyperbolic trapping (see e.g. Ikawa \cite{Ik88_01}, Nonnenmacher and Zworski \cite{NoZw09_01} or Petkov and Stoyanov \cite{PeSt10_01}).

Second, we suppose that the $n$ dimensional measure of $\CH$ is not too large. In this situation which allows strong trappings, we show that $P$ has no resonance and that its truncated resolvent has a polynomial estimate in a neighborhood of size $h \vert \ln h \vert^{- 1}$ of the real axis. A single homoclinic trajectory in dimension $1$ enters into this setting.

We then give in Section \ref{s6} the asymptotic of the resonances in the most natural weak trapping situation: we suppose that $\CH$ consists of $K$ trajectories along which $\Lambda_{-}$ and $\Lambda_{+}$ intersect transversally. We form a $K \times K$ matrix $\CQ ( z , h )$, depending only on dynamical quantities along $\CH$ and near $0$, which mixes the contributions of the different homoclinic trajectories. In this case, the {\it quantization rule} reads
\begin{equation} \label{m93}
1 \in \spe \Big( h^{\sum_{j \geq 2} \frac{\lambda_{j}}{2 \lambda_{1}} - i \frac{z - E_{0}}{\lambda_{1} h }} \CQ ( z , h ) \Big) ,
\end{equation}
where $\spe$ stands for spectrum of. Writing this in terms of the determinant leads to a scalar equation, whose solutions are called {\it pseudo-resonances}. The first main theorem of the paper states that the resonances are close to the pseudo-resonances. Notice that the quantization rule is not necessarily exact. More precisely, we have
\begin{equation} \label{m97}
\text{distance} \big( \text{resonances of } P , \text{pseudo-resonances of } P \big) = o \Big( \frac{h}{\vert \ln h \vert} \Big) ,
\end{equation}
in the region
\begin{equation} \label{m96}
0 \geq \im z \geq - \sum_{j = 2}^{n} \frac{\lambda_{j}}{2} h - C \frac{h}{\vert \ln h \vert} ,
\end{equation}
for any $C > 0$. Combining \eqref{m93} with \eqref{m97}, the asymptotic of the resonances $z$ satisfying $\re z \approx E_{0} + \tau h$ is given by
\begin{equation} \label{m94}
z = E_{0} + 2 q \pi \lambda_{1} \frac{h}{\vert \ln h \vert} - i h \sum_{j = 2}^{n} \frac{\lambda_{j}}{2} + i \ln ( \mu_{k} ( \tau , h ) ) \lambda_{1} \frac{h}{\vert \ln h \vert} + o \Big( \frac{h}{\vert \ln h \vert} \Big) ,
\end{equation}
where $q \in \Z$ and the $\mu_{k}$'s are the eigenvalues of $\CQ ( E_{0} + h \tau - i h \sum_{j \geq 2} \lambda_{j} / 2 , h )$. In particular, the $\mu_{k} ( \tau , h )$'s are continuous functions of $\tau , h$. We also obtain a polynomial estimate of the distorted resolvent at distance $h \vert \ln h \vert^{- 1}$ from the resonances. Comparing \eqref{n1} with \eqref{m94}, we notice that the resonances are closer to the real axis here than for the barrier-top alone and that there are way more resonances (typically at least $\vert \ln h \vert$ in domains of size $h$).

\begin{figure}
\begin{center}
\begin{picture}(0,0)%
\includegraphics{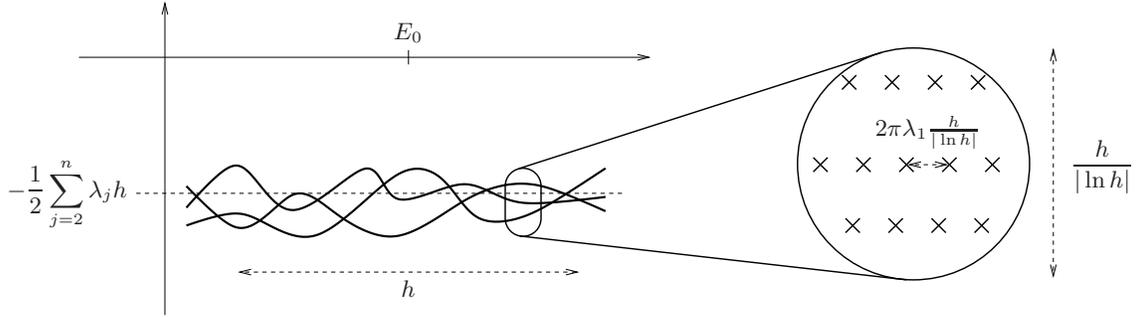}%
\end{picture}%
\setlength{\unitlength}{1184sp}%
\begingroup\makeatletter\ifx\SetFigFont\undefined%
\gdef\SetFigFont#1#2#3#4#5{%
  \reset@font\fontsize{#1}{#2pt}%
  \fontfamily{#3}\fontseries{#4}\fontshape{#5}%
  \selectfont}%
\fi\endgroup%
\begin{picture}(22305,6644)(-2114,-6383)
\put(-2099,-3886){\makebox(0,0)[lb]{\smash{{\SetFigFont{9}{10.8}{\rmdefault}{\mddefault}{\updefault}$\ds - \frac{1}{2} \sum_{j = 2}^{n} \lambda_{j} h$}}}}
\put(6301,-586){\makebox(0,0)[b]{\smash{{\SetFigFont{9}{10.8}{\rmdefault}{\mddefault}{\updefault}$E_{0}$}}}}
\put(17176,-2611){\makebox(0,0)[b]{\smash{{\SetFigFont{9}{10.8}{\rmdefault}{\mddefault}{\updefault}$2 \pi \lambda_{1} \frac{h}{\vert \ln h \vert}$}}}}
\put(20176,-3361){\makebox(0,0)[lb]{\smash{{\SetFigFont{9}{10.8}{\rmdefault}{\mddefault}{\updefault}$\ds \frac{h}{\vert \ln h \vert}$}}}}
\put(6301,-5986){\makebox(0,0)[b]{\smash{{\SetFigFont{9}{10.8}{\rmdefault}{\mddefault}{\updefault}$h$}}}}
\end{picture}%
\end{center}
\caption{The two scale asymptotic of the resonances.} \label{f64}
\end{figure}

The main feature in \eqref{m94} is that the resonances satisfy a {\it two scale asymptotic}, illustrated in Figure \ref{f64}. At the macroscopic scale $h$, they concentrate on continuous {\it accumulation curves} given by
\begin{equation} \label{m95}
\im \sigma = - i \sum_{j = 2}^{n} \frac{\lambda_{j}}{2} + \ln \big( \vert \mu_{k} ( \re \sigma , h ) \vert \big) \frac{\lambda_{1}}{\vert \ln h \vert} ,
\end{equation}
for $k = 1 , \ldots , K$ where $\sigma = ( z - E_{0} ) / h$ is the natural macroscopic parameter. Since some of the $\mu_{k}$'s can vanish, there may be less than $K$ of these curves in the region \eqref{m96}. On the other hand, at the microscopic scale $h \vert \ln h \vert^{- 1}$, the resonances are regularly distributed on horizontal lines. This particular distribution of the resonances has been observed in other settings by physicists and quantum chemists (see e.g. Korsch \cite{Ko89_01} or Burghardt and Gaspard \cite{BuGa97_01}) and in numerical computations (see e.g. Barkhofen, Faure and Weich \cite{BaFaWe14_01}). Moreover, these accumulation curves enjoy remarkable properties that we study in Section \ref{s18}. First, they are quasi-periodical functions of $h^{- 1}$. We call this phenomenon {\it vibration} in the sequel. We also describe their asymptotic behavior as $\re \sigma \to \pm \infty$ and show that they are related to the nature of the trapping for the energies above and below $E_{0}$. In other words, \eqref{m95} may be seen as the {\it transition} between different types of trappings. Eventually, we prove that \eqref{m94} is {\it stable} under small perturbations of the homoclinic set $\CH$. As a curiosity, we provide situations where the resonances move away from the real axis when the trapping increases.

We also study the distribution of resonances deeper in the complex. Typically, there is a band of size $h$ without resonances below the region \eqref{m96}, and then $( n - 1 ) K$ accumulation curves near the line
\begin{equation*}
\im z \approx - i h \sum_{j = 2}^{n} \frac{\lambda_{j}}{2} - i h \lambda_{1} ,
\end{equation*}
corresponding to the $n - 1$ transverse directions to the $K$ homoclinic trajectories. Thus, the situation seems similar to the one of hyperbolic trappings (see e.g. G{\'e}rard and Sj{\"o}strand \cite{GeSj87_01} or Faure and Tsujii \cite{FaTs15_01}). Nevertheless, the mechanisms behind these phenomena differ. Moreover, under additional assumptions, the asymptotic of resonances can be sharpened. More precisely, there exists in this case a $K \times K$ matrix $\CQ_{\rm tot} ( z , h ) = \CQ ( z , h ) + o ( 1 )$ such that the remainder term $o ( h \vert \ln h \vert^{- 1} )$ is replaced by $\CO ( h^{\infty} )$ in \eqref{m97}, if the pseudo-resonances are defined by $\CQ_{\rm tot}$ instead of $\CQ$ in \eqref{m93}. Eventually, we obtain the repartition of the resonances when $\Lambda_{-}$ and $\Lambda_{+}$ have a tangential intersection of finite order along a finite number of trajectories. The corresponding results are similar to \eqref{m93}--\eqref{m94} but the resonances are closer to the real axis.

\begin{figure}
\begin{center}
\begin{picture}(0,0)%
\includegraphics{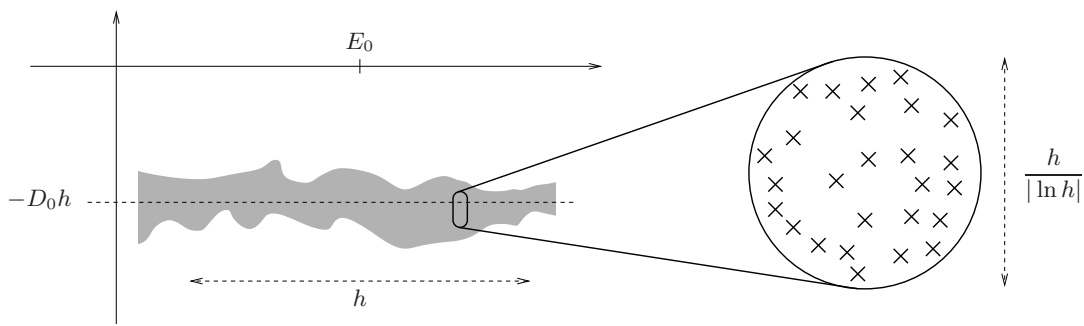}%
\end{picture}%
\setlength{\unitlength}{1184sp}%
\begingroup\makeatletter\ifx\SetFigFont\undefined%
\gdef\SetFigFont#1#2#3#4#5{%
  \reset@font\fontsize{#1}{#2pt}%
  \fontfamily{#3}\fontseries{#4}\fontshape{#5}%
  \selectfont}%
\fi\endgroup%
\begin{picture}(21285,6644)(-1094,-6383)
\put(6301,-586){\makebox(0,0)[b]{\smash{{\SetFigFont{9}{10.8}{\rmdefault}{\mddefault}{\updefault}$E_{0}$}}}}
\put(20176,-3361){\makebox(0,0)[lb]{\smash{{\SetFigFont{9}{10.8}{\rmdefault}{\mddefault}{\updefault}$\ds \frac{h}{\vert \ln h \vert}$}}}}
\put(6301,-5986){\makebox(0,0)[b]{\smash{{\SetFigFont{9}{10.8}{\rmdefault}{\mddefault}{\updefault}$h$}}}}
\put(-1079,-3886){\makebox(0,0)[lb]{\smash{{\SetFigFont{9}{10.8}{\rmdefault}{\mddefault}{\updefault}$- D_{0} h$}}}}
\end{picture}%
\end{center}
\caption{The cloud of resonances.} \label{f65}
\end{figure}

In Section \ref{s42}, we generalize the previous results to several maximums. We allow for a finite number of hyperbolic fixed points and of transversal homoclinic/heteroclinic trajectories between them. In this situation, we see the trapped set as a graph where the vertices are the fixed points and the edges are the homoclinic/heteroclinic trajectories. To each cycle of this graph, we associate its {\it damping}, a non-negative quantity measuring the dissipation of the energy along this cycle. The resonances closest to the real axis are generated by the cycles with the minimal damping, noted $D_{0}$. More precisely, the pseudo-resonances are defined here as the solutions of quantization rule
\begin{equation} \label{m98}
1 \in \spe ( \SQ ( z , h ) ) ,
\end{equation}
where the matrix $\SQ ( z , h )$ is given in terms of dynamical quantities. By comparison with \eqref{m93}, the exponent of $h$ depends on the coefficient and can no longer be factored out of $\SQ$. We then prove that \eqref{m97} holds true up to the zone $\im z \approx - D_{0} h$. The situation seems similar to \eqref{m93}--\eqref{m94}, but the asymptotic of resonances can be different. Namely, the resonances $z$ with $\re z \approx E_{0} + \tau h$ satisfy
\begin{equation}
z = E_{0} + \tau h - i D_{0} h + Z \frac{h}{\vert \ln h \vert} + o \Big( \frac{h}{\vert \ln h \vert} \Big) ,
\end{equation}
where $Z$ is a zero of some exponential sum $f_{\tau} ( \cdot , h )$. Depending on the arithmetic properties of the eigenvalues of the Hessian of $V$ at the vertices, these resonances either concentrate on accumulations curves as in Figure \ref{f64}, or form a {\it cloud} as illustrated in Figure \ref{f65}. Note that the proof is the same in these two cases. Some of the phenomena of the two previous paragraphs occur also in the present setting. As an application, we give the explicit asymptotic of the resonances in the three bump case. In dimension $1$, the asymptotic of the eigenvalues has been obtained by Colin de Verdi\`ere and Parisse \cite{CoPa99_01} under similar assumptions.

We provide in Section \ref{s26} the asymptotic of the resonances in the most natural strong trapping situation. We assume that the hyperbolic fixed point is isotropic (i.e. $\lambda_{1} = \cdots = \lambda_{n} = : \lambda$). Let $\CH^{\infty}$ denote the compact set of the normalized asymptotic directions of $\CH$ at the barrier-top. Under some geometric assumptions, the pseudo-resonances are defined by the quantization rule
\begin{equation}
1 \in \spe \big( h^{- i \frac{z - E_{0}}{\lambda h }} \CT ( z , h ) \big) .
\end{equation}
Here, $\CT$ is no longer a matrix but an operator on the space $L^{2} ( \CH^{\infty} )$. Near the real axis, the resonances are again approximated by the pseudo-resonances in the sense of \eqref{m97}. Moreover, the resonances $z$ such that $\re z \approx E_{0} + \tau h$ verify the asymptotic
\begin{equation} \label{m99}
z = E_{0} + 2 q \pi \lambda \frac{h}{\vert \ln h \vert} + i \ln ( \mu_{k} ( \tau , h ) ) \lambda \frac{h}{\vert \ln h \vert} + o \Big( \frac{h}{\vert \ln h \vert} \Big) ,
\end{equation}
where $q \in \Z$ and the $\mu_{k}$'s are the eigenvalues of $\CT ( E_{0} + h \tau , h )$. Hence, the resonances concentrate on accumulation curves enjoying properties similar to those of \eqref{m95}. There can be an infinite number of such accumulation curves. No regularity hypothesis (except compactness which always holds) is made on the set $\CH^{\infty}$, which can be a Cantor set for instance. Lastly, one can formally recover \eqref{m94} from \eqref{m99}, taking $\CH$ small.

Finally, we give the asymptotic of the {\it resonant states} in Section \ref{s32}. These results are stated for a finite number of transversal homoclinic/heteroclinic trajectories. We prove that the resonant states are Lagrangian distributions near each point of the homoclinic/heteroclinic trajectories. Moreover, they are naturally associated to the eigenvectors of the quantization operator $\SQ$ corresponding to the eigenvalue $1$ (see \eqref{m97} and \eqref{m98}). We also give some examples showing that these functions can be unexpectedly distributed. For instance, we sometimes observe the {\it delocalization} of the resonants states. By this, we mean that the resonances are generated by a part of the trapped set (i.e. the dynamical quantities appearing in the asymptotic of the resonances are localized near this part), whereas the corresponding resonant states are small near this part and large near another part of the trapped set.

\Subsection{General strategy for the proofs} \label{s36}

Let us now describe the approach we use throughout the paper. Let $H$ be a closed operator on some Banach space $B$ having only discrete spectrum near $z_{0} \in \C$. We want to show that $H$ has an eigenvalue close to $z_{0}$. We proceed in two steps: \parskip = 0 in

$i)$ First, we show that $H$ has no spectrum in the ring $\SC ( \varepsilon ) : = B ( z_{0} , \varepsilon_{0} ) \setminus B ( z_{0} , \varepsilon )$ and that the resolvent $( H - z )^{- 1}$ satisfies a nice estimate in this region, say
\begin{equation} \label{m85}
\forall z \in \SC ( \varepsilon ) , \qquad \big\Vert ( H - z )^{- 1} \big\Vert \leq M .
\end{equation}
Mostly, we obtain this inequality by a contradiction argument: if \eqref{m85} does not hold true, there exist $u \in B$ and $z \in \SC ( \varepsilon )$ such that
\begin{equation} \label{m86}
\Vert ( H - z ) u \Vert < M^{- 1} \qquad \text{and} \qquad \Vert u \Vert = 1 .
\end{equation}
All the work consists then in proving that this is impossible.

$ii)$ Second, we show that $H$ has at least one eigenvalue in $\SD ( \delta ) : = B ( z_{0} , \delta )$, with $\varepsilon < \delta < \varepsilon_{0}$. To this aim, we construct a test vector $v \in B$ and we consider
\begin{equation} \label{m87}
u ( \widetilde{z} ) : = ( H - \widetilde{z} )^{- 1} v ,
\end{equation}
for $\widetilde{z} \in \partial \SD ( \delta )$. Since $\partial \SD ( \delta ) \subset \SC ( \varepsilon )$, the resolvent estimate \eqref{m85} yields that $u ( \widetilde{z} )$ is well-defined, holomorphic and bounded by $M$ on $\partial \SD ( \delta )$. Then we prove that
\begin{equation} \label{m91}
U : = \int_{\partial \SD ( \delta )} u ( \widetilde{z} ) \, d \widetilde{z} = \int_{\partial \SD ( \delta )} ( H - \widetilde{z} )^{- 1} v \, d \widetilde{z} \neq 0 ,
\end{equation}
for a well chosen $v$. This is generally done computing first $u ( \widetilde{z} )$, the solution of $( H - \widetilde{z} ) u ( \widetilde{z} ) = v$. If $H$ had no eigenvalue in $\CD ( \delta )$, we would have $U = 0$. Therefore, $H$ has at least one eigenvalue in $\CD ( \delta )$.

\noindent
Summing up, the two previous points imply that \parskip = 0.05 in
\begin{equation*}
\dist \big( \spe ( H ) , \{ z_{0} \} \big) < \varepsilon ,
\end{equation*}
in $B ( z_{0} , \varepsilon_{0} )$ (in the sense Definition \ref{g80}), and that the resolvent of $H$ has an nice upper bound in $\SC ( \varepsilon )$. Note that to get an eigenvalue free domain, it is enough to prove the point $i)$. Moreover, to carry out $ii)$, the first point must already have been verified. Eventually, to refine the asymptotic of the spectrum, it is enough to refine $i)$ (i.e. to take $\varepsilon$ smaller). This approach is not specific to the study of the eigenvalues of an operator. It can be used to give the asymptotic of the characteristic values of holomorphic operator-valued functions (see e.g. Gohberg and Leiterer \cite{GoLe09_01}).

For the study of the semiclassical resonances, we take $H = P_{\theta}$, the distorted operator of angle $\theta$ defined in Section \ref{s2}. For the type of problem considered in this paper, it is natural to work with operators which are $C^{\infty}$ everywhere and analytic only at infinity (in order to define the resonances). Thus, we will use the $C^{\infty}$ semiclassical microlocal calculus. It means that the remainder terms will be of order $\CO ( h^{\infty} )$ (and not of order $\CO ( e^{- \delta / h} )$ as in the analytic category). The same way, it is natural to consider $M = h^{- N}$ in \eqref{m85}. This choice is also based on the belief that the resolvent of $P_{\theta}$ can not behave worse than $h^{- N}$ away from the resonances (there is no pseudo-spectrum). Then, \eqref{m86} leads to equations of the form
\begin{equation} \label{m89}
( P_{\theta} - z ) u = \CO ( h^{\infty} ) ,
\end{equation}
with $\Vert u \Vert = 1$, whereas \eqref{m87} writes
\begin{equation} \label{m90}
( P_{\theta} - z ) u = v .
\end{equation}
Note that \eqref{m86} guaranties that the unique solution $u$ of \eqref{m90} satisfies $\Vert u \Vert \leq h^{- N}$, which allows us to apply the $C^{\infty}$ microlocal analysis.

Using the ellipticity of $P_{\theta}$ at infinity, this type of equations can be reduced to some microlocal Cauchy problems near the trapped set:
\begin{equation} \label{m88}
\left\{ \begin{aligned}
&( P - z ) u = f &&\text{microlocally near the trapped set,}  \\
&u = g &&\text{microlocally near the incoming region.}
\end{aligned} \right.
\end{equation}
The problem is said to be homogeneous when $f = 0$, and $g$ is called the initial data. These microlocal Cauchy problems must be seen as usual evolution equations. But the propagation must be understood in the spirit of H{\"o}rmander and not using the time evolution $e^{i t P / h}$. The intuition is that $z$ is a resonance if and only if there is no uniqueness in the microlocal Cauchy problem \eqref{m88}. This idea exists in the folklore of the theory of resonances. More precisely, we show in Section \ref{s3} that the uniqueness in a quantitative sense in \eqref{m88} implies that \eqref{m89} has no solution and then that $z$ is not a resonance from $i)$. On the other hand, we give the asymptotic of the solution (existence) of the microlocal Cauchy problem \eqref{m88} in order to compute the solution of \eqref{m90} and realize $ii)$. Remark that it is generally enough to be able to solve the homogeneous problem. It means that one may take $v$ microlocalized outside the (problematic part of the) trapped set.

Most of the previous arguments hold in the general setting of semiclassical resonances. To the contrary, the way we solve \eqref{m88} is specific to the geometry treated here. The idea is to decompose \eqref{m88} into small microlocal Cauchy problems near the hyperbolic fixed points or along the homoclinic/heteroclinic trajectories. For the first ones, we use our previous work \cite{BoFuRaZe07_01}. For the second ones, we apply the usual propagation of singularities of H{\"o}rmander. To achieve this decomposition, the microlocal Cauchy problem near the fixed points must be well posed. This is why we require that $z$ avoids some discrete set $\Gamma ( h )$.

The main advantage of our approach is the flexibility. We do not need to understand precisely what happens at $z_{0}$, since it is enough to work at respectable distance. Thus, it is possible to consider complicated geometries such as trappings with different types of trajectories, to concentrate the study near the to pertinent part of the trapped set, to treat crossings of resonances, $\ldots$ In particular, the quantization rules of this paper approximate the resonances only in very specific regions of $\C$. By comparison with the Grushin method (see e.g. Helffer and Sj\"{o}strand \cite{HeSj86_01}), we do not have to construct quasimodes nor to inverse the so called Grushin problem. Instead, we only have to solve the microlocal Cauchy problem \eqref{m88} for one data $( f , g )$ of our choice. Moreover, we only have to show that $U \neq 0$ in \eqref{m91}. Furthermore, using the standard propagation of singularities instead of conjugating the operator $P$ near the trapped set allows to consider easily $C^{\infty}$ operators instead of analytic ones. Of course, it is possible to achieved the Grushin approach in the $C^{\infty}$ category (see e.g. the work of Lahmar-Benbernou, Martinez and the second author \cite{FuLaMa11_01}), as it is possible to work in the analytic category with our method. This strategy also provides polynomial estimates of the cut-off resolvent and the asymptotic of the resonant states on the base space. On the other hand, the main weakness of our approach is that we only give the asymptotic of the resonances as a set. In other words, we have no upper bound on the multiplicity. If we were able to solve the general inhomogeneous microlocal Cauchy problem \eqref{m88}, it might be possible to overcome this difficulty. Lastly, the Grushin method also provides a formula for the resolvent even near the spectrum that we do not have obtained here.

For the reader's convenience, we provide at the end of the paper an index of notations and the table of contents.

\section{General setting} \label{s2}

In this paper, we consider the Schr\"odinger operator on $\R^{n}$, $n \geq 1$,
\begin{equation}\label{a5}
P = - h^{2} \Delta + V (x),
\end{equation}
where $V$ is a smooth real-valued function and $h$ is a positive number. We work in the semiclassical regime $h \to 0$. We denote by $p ( x , \xi ) = \xi^{2} + V ( x )$ the associated classical Hamiltonian. The vector field
\begin{equation} \label{d27}
H_{p} = \partial_{\xi} p \cdot \partial_{x} - \partial_{x} p \cdot \partial_{\xi} = 2 \xi \cdot \partial_{x} - \nabla V ( x ) \cdot \partial_{\xi} ,
\end{equation}
is the Hamiltonian vector field associated to $p$. Integral curves $t \mapsto \exp (t H_p )( x , \xi )$ of $H_p$ are called classical trajectories or bicharacteristic curves, and $p$ is constant along such curves.
The trapped set at a real energy $E$ for $P$ is defined as 
\begin{equation*}
K(E) = \big\{ (x,\xi) \in p^{-1}(E); \ \exp(tH_p)(x,\xi)\not\to \infty \mbox{ as } t\to \pm \infty \big\} .
\end{equation*}
We shall suppose that $V$ satisfies the following assumptions
\begin{hyp} \label{h1}
$V \in C^{\infty} ( \R^{n} ; \R )$ extends holomorphically in the sector 
\begin{equation*}
\CS = \big\{ x \in \C^{n} ; \ \vert \im x \vert \leq \delta \< x \> \text{ and } \vert x \vert \geq C \big\} ,
\end{equation*}
for some $C , \delta > 0$. Moreover $V(x)\to 0$ as $x\to\infty$ in $\CS$.
\end{hyp}

Under the previous assumption, the operator $P$ is self-adjoint with domain $H^2(\R^n)$, and we define the resonances of $P$ as follows (see Aguilar and Combes \cite{AgCo71_01}, Hunziker \cite{Hu86_01} or Sj{\"o}strand and Zworski \cite{SjZw91_01} for an alternative approach). Let $R_0 > 0$ be a large constant, and let $F:\R^n\to\R^n$ be a smooth vector field, such that $F(x)=0$ for $\vert x\vert\leq R_0$ and $F(x)=x$ for $\vert x\vert\geq R_{0} +1$. For $\mu \in \R$ small enough, we denote $U_\mu:L^2(\R^n)\to L^2(\R^n)$ the unitary operator defined by
\begin{equation} \label{j70}
U_\mu \varphi (x) = \big\vert \det(1+\mu d F (x)) \big\vert^{1/2}\varphi(x+\mu F(x)) ,
\end{equation}
for $\varphi\in C^{\infty}_{0}(\R^n)$. Then, the operator $U_{\mu}P(U_{\mu})^{-1}$ is a differential operator with analytic coefficients with respect to $\mu$, and can be analytically continued to small enough complex values of $\mu$. For $\theta \in \R$ small enough, we denote
\begin{equation}\label{a6}
P_{\theta}=U_{i\theta}P(U_{i\theta})^{-1}.
\end{equation}
The spectrum of $P_{\theta}$ is discrete in $\CE_{\theta}=\{z\in \C ; \ -2\theta<\arg z\leq 0\}$, and the resonances of $P$ are by definition the eigenvalues of $P_{\theta}$ in $\CE_{\theta}$. The resonances do not really depend on $\theta$ and $F$, and we denote their set by $\Res(P)$.

As a matter of fact, the resonances are also the poles of the meromorphic extension from the upper complex half-plane of the resolvent $(P-z)^{-1}: L^{2}_{\rm comp} ( \R^{n} )\to L^{2}_{\rm loc} ( \R^{n} )$ (see e.g. Helffer and Martinez \cite{HeMa87_01}). Furthermore, if $\chi \in C^{\infty}_{0} ( \R^{n} )$ is supported outside of the complex dilation, we have (see Sj{\"o}strand and Zworski \cite{SjZw91_01})
\begin{equation} \label{c12}
\chi ( P - z)^{-1} \chi = \chi ( P_{\theta} - z)^{-1} \chi .
\end{equation}
Moreover, Proposition \ref{j66} shows that the truncated resolvent and the distorted resolvent have essentially the same norm. We send the reader to Sj\"{o}strand \cite{Sj07_01} or Dyatlov and Zworski \cite{DyZw16_01} for more details on the theory of resonances.

\begin{figure}
\begin{center}
\begin{picture}(0,0)%
\includegraphics{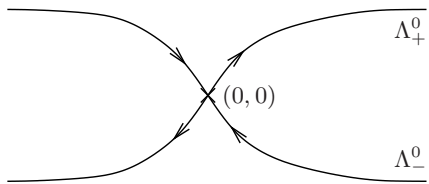}%
\end{picture}%
\setlength{\unitlength}{1184sp}%
\begingroup\makeatletter\ifx\SetFigFont\undefined%
\gdef\SetFigFont#1#2#3#4#5{%
  \reset@font\fontsize{#1}{#2pt}%
  \fontfamily{#3}\fontseries{#4}\fontshape{#5}%
  \selectfont}%
\fi\endgroup%
\begin{picture}(8841,3666)(1168,-5794)
\put(9301,-2761){\makebox(0,0)[lb]{\smash{{\SetFigFont{9}{10.8}{\rmdefault}{\mddefault}{\updefault}$\Lambda_{+}^{0}$}}}}
\put(5701,-4111){\makebox(0,0)[lb]{\smash{{\SetFigFont{9}{10.8}{\rmdefault}{\mddefault}{\updefault}$( 0 , 0 )$}}}}
\put(9301,-5461){\makebox(0,0)[lb]{\smash{{\SetFigFont{9}{10.8}{\rmdefault}{\mddefault}{\updefault}$\Lambda_{-}^{0}$}}}}
\end{picture}%
\end{center}
\caption{The incoming/outgoing manifolds $\Lambda_{\pm}^{0}$.} \label{f66}
\end{figure}

We now make some assumptions on the geometry of the Hamiltonian flow.
\begin{hyp} \label{h2}
$V$ has a non-degenerate maximum at $x=0$ and
\begin{equation*}
V (x) = E_{0} - \sum_{j=1}^{n} \frac{\lambda_{j}^{2}}{4} x_{j}^{2} + \CO ( x^{3} ) ,
\end{equation*}
with $E_{0}> 0$ and $0 < \lambda_{1} \leq \lambda_{2} \leq \dots \leq \lambda_{n}$.
\end{hyp}
The linearization $F_{p}$ at $(0,0)$ of the Hamilton vector field $H_{p}$ is given by
\begin{equation*}
F_{p} =
\left( \begin{array}{cc}
0 & 2 Id \\
\frac{1}{2} \diag ( \lambda_{1}^{2} , \ldots , \lambda_{n}^{2} ) & 0
\end{array} \right) ,
\end{equation*}
and has eigenvalues $- \lambda_{n} , \dots , - \lambda_{1} , \lambda_{1} , \dots , \lambda_{n}$. Thus, $(0,0)$ is a hyperbolic fixed point for $H_{p}$ and the local stable/unstable manifold theorem gives the existence of a local incoming Lagrangian manifold $\Lambda_{-}^{0}$ and a local outgoing Lagrangian manifold $\Lambda_{+}^{0}$ characterized by
\begin{align*}
\Lambda_{\pm}^{0} = \big\{ ( x , \xi ) \in T^{*} \R^{n} \text{ near } ( 0 , 0 ) ; \ & \exp(t H_{p})( x , \xi ) \to (0,0) \\
\mbox{ as } & t \to \mp \infty \text{ and staying near } ( 0 , 0 ) \big\} \subset p^{-1} (E_{0}) .
\end{align*}
Moreover, there exist two smooth functions $\varphi_{\pm}$, defined in a vicinity of $0$, satisfying
\begin{equation} \label{a57}
\varphi_{\pm} (x) = \pm \sum_{j=1}^{n} \frac{\lambda_{j}}{4} x_{j}^{2} + \CO (x^{3}) ,
\end{equation}
and such that $\Lambda_{\pm}^{0} = \Lambda_{\varphi_{\pm}} : = \{ (x , \xi ) ; \ \xi = \nabla \varphi_{\pm} (x) \}$ near 
$(0,0)$. Since $P$ is a Schr\"odinger operator, we have $\varphi_{-} = - \varphi_{+}$. We also define the global incoming/outgoing Lagrangian manifold $\Lambda_{\pm}$ by
\begin{equation*}
\Lambda_{\pm} = \bigcup_{t \in \R} \exp ( t H_{p} ) ( \Lambda_{\pm}^{0} ) ,
\end{equation*}
which satisfy
\begin{equation} \label{j49}
\Lambda_{\pm} = \big\{ ( x , \xi ) \in T^{*} \R^{n} ; \ \exp(t H_{p})( x , \xi ) \to (0,0) \mbox{ as } t \to \mp \infty \big\} .
\end{equation}
We assume
\begin{hyp} \label{h3}
The trapped set at energy $E_{0}$ satisfies
\begin{equation*}
K ( E_{0} ) \subset \Lambda_{-} \cap \Lambda_{+} .
\end{equation*}
\end{hyp}
We denote by $\CH : = \Lambda_{-} \cap \Lambda_{+} \setminus \{ ( 0 ,0 ) \}$ the set of homoclinic curves. Recall that a homoclinic curve is a non-constant bicharacteristic curve which converges to the fixed point $( 0 ,0 )$ as the time goes to $- \infty$ and to $+ \infty$.

For $\rho \in \Lambda_{\pm}$, let $( x (t) , \xi (t) ) = \exp ( t H_{p} ) ( \rho )$ be the corresponding Hamiltonian curve. From Helffer and Sj\"{o}strand \cite[Equation (2.7)]{HeSj85_01}, we know that
\begin{equation} \label{d2}
x (t) = g_{\pm} ( \rho ) e^{\pm \lambda_{1} t} + o \big( e^{\pm \lambda_{1} t} \big) ,
\end{equation}
as $t \to \mp \infty$ and for some $g_{\pm} ( \rho ) \in \R^{n}$. Therefore, $g_{\pm} ( \rho )$ is the asymptotic direction of the curve when this quantity does not vanish. Actually, from \cite{BoFuRaZe07_01,HeSj85_01}, the functions $g_{\pm}$ are $C^{\infty}$ on $\Lambda_{\pm}^{0}$, $g_{\pm} ( \rho ) \in \ker ( \Hess V ( 0 ) + \lambda_{1}^{2} / 2 )$ and $g_{\pm}$ vanishes on a submanifold of $\Lambda_{\pm}^{0}$ of dimension $\card \{ j ; \ \lambda_{j} \neq \lambda_{1} \}$. In the sequel, we suppose
\begin{hyp} \label{h4}
For all $\rho , \widetilde{\rho} \in \CH$, we have $g_{-} ( \rho ) \cdot g_{+} ( \widetilde{\rho} ) \neq 0$.
\end{hyp}
This hypothesis can be reformulate in a more symplectic way. Using $\Lambda_{\pm}^{0} = \Lambda_{\varphi_{\pm}}$, \eqref{a57} and \eqref{d2}, we deduce $\exp ( t H_{p} ) ( \rho ) = \rho_{\pm} e^{\pm \lambda_{1} t} + o ( e^{\pm \lambda_{1} t} )$ as $t \to \mp \infty$ for some $\rho_{\pm} \in T_{( 0 , 0 )} ( T^{*} \R^{n} )$ with $g_{\pm} ( \rho ) = d \pi_{x} ( \rho_{\pm} )$. Then, \ref{h4} is equivalent to
\begin{equation*}
\sigma ( \rho_{-} , \widetilde{\rho}_{+} ) \neq 0 ,
\end{equation*}
for all $\rho , \widetilde{\rho} \in \CH$. In the previous equation, $\sigma$ denotes the canonical 2--form on $T^{*} \R^{n}$.

The assumption \ref{h4} has some geometric consequences which are described in Section \ref{s11}. In particular, the manifolds $\Lambda_{\pm}$ can not be wound on themselves from Proposition \ref{a44} (this means that, in a neighborhood of $\exp ( t H_{p} ) ( \rho )$ with $\rho \in \Lambda_{\pm}^{0}$, the manifold $\Lambda_{\pm}$ coincides with $\exp ( t H_{p} ) ( \Lambda_{\pm}^{0} )$).

\begin{remark}\sl \label{c13}
$i)$ The results of this paper may remain valid for more general operators. For example, one could consider metric perturbations, obstacles (see Example \ref{b86} below), pseudodifferential operators (as in Example \ref{b81}), $\ldots$ What is important is that the resonances can be defined by complex distortion, that the semiclassical microlocal analysis can be used and that the geometric setting induced by the assumptions \ref{h2}--\ref{h4} holds true. In this direction, the framework of ``black box'' operators, due to Sj\"{o}strand and Zworski \cite{SjZw91_01}, can be used to define the resonances in many situations.

$ii)$ It is also possible to add a damping potential of the form $- i \alpha (h) W (x)$ (or more general) with $h \leq \alpha ( h ) \leq 1$, $W \in C^{\infty}_{0} ( \R^{n} )$ and $W \geq 0$ to ``remove'' some part of the trapped set. Indeed, if $\alpha$ is not too small, the contribution of the Hamiltonian curves passing trough $W$ vanishes. This allows to construct examples more easily. Note that the usual definition of the resonances by complex distortion (see Hunziker \cite{Hu86_01}) can be adapted to this setting. Such an idea has been developed by Royer \cite{Ro10_01} for the limiting absorption principle and by Datchev and Vasy \cite[Section 5.3]{DaVa12_01} in a situation close to ours.

$iii)$ Eventually, note that, thanks to Datchev and Vasy \cite{DaVa10_01}, some of the present results can directly be extended to manifolds with more complicated ends.
\end{remark}

\section{Resonance free domains} \label{s22}

In this section, we show the existence of resonance free zones below the real axis. We distinguish between two regimes that we call weak and strong trapping. They give rise to rather different zones. The corresponding proofs are gathered in Section \ref{s4} and Section \ref{s5} respectively.

\Subsection{Weak trapping}

In addition to the general hypotheses of the previous section, we make in this part the following geometric assumption which, roughly speaking, means that the trapping induced by the homoclinic curves is not very strong.
\begin{hyp} \label{h5}
We assume that one of the two conditions below is fulfilled. \smallskip
\begin{enumerate}  \renewcommand{\theenumii}{{{\rm (\alph{enumii})}}} \renewcommand{\labelenumii}{{(\alph{enumii})}}
\item {\rm Anisotropic case:} \quad $\lambda_{1} < \lambda_{n}$, \label{h5a}  \smallskip
\item {\rm Transversal case:} \quad $\forall \rho \in \CH , \qquad T_{\rho} \Lambda_{-} \neq T_{\rho} \Lambda_{+}$. \label{h5b}
\end{enumerate}
\end{hyp}
Note that, in dimension $n=1$, this assumption implies that $\CH = \emptyset$ (i.e. the trapped set reduces to the hyperbolic fixed point $( 0 , 0)$). Indeed, in the other case, $\Lambda_{-}$ and $\Lambda_{+}$ (and then their tangent space) will coincide along $\CH$. The following remark explains that the dimension of the homoclinic set (in a particular sense) is not maximal if \ref{h5} is verified.

\begin{remark}\sl \label{c16}
From \ref{h4}, we have $g_{\pm} ( \rho ) \neq 0$ for $\rho \in \CH$. Then, we can define
\begin{equation*}
\CH^{\pm \infty} = \Big\{ \frac{g_{\pm} ( \rho )}{\vert g_{\pm} ( \rho ) \vert} ; \ \rho \in \CH \Big\} \subset \S^{n-1} ,
\end{equation*}
the asymptotic directions of $\CH$ at the critical point $( 0 , 0 )$. Note that $\CH^{+ \infty} = \CH^{- \infty}$ since $P$ is a Schr\"{o}dinger operator. If \ref{h5a} is satisfied, $\CH^{\pm \infty}$ is a subset of $\{ x \in \S^{n - 1} ; \ x_{j} = 0 \text{ if } \lambda_{j} \neq \lambda_{1} \}$ which is a submanifold of $\S^{n-1}$ of dimension $\card \{ j ; \ \lambda_{j} = \lambda_{1} \} \leq n - 1$. The same way, if \ref{h5b} holds true, $\CH^{\pm \infty}$ is locally contained in a hypersurface of $\S^{n-1}$. Thus, \ref{h5} implies that the trapped set is not too big in the sense that the set of asymptotic directions $\CH^{\pm \infty}$ is at most of dimension $n-1$ (even if $\CH$ can be of dimension $n$ near some points).
\end{remark}

Under the previous assumption, our main result is the following.

\begin{figure}
\begin{center}
\begin{picture}(0,0)%
\includegraphics{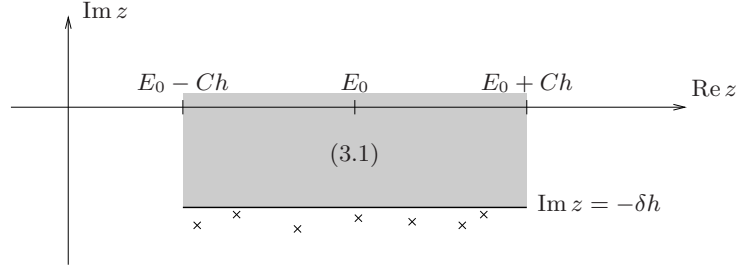}%
\end{picture}%
\setlength{\unitlength}{1184sp}%
\begingroup\makeatletter\ifx\SetFigFont\undefined%
\gdef\SetFigFont#1#2#3#4#5{%
  \reset@font\fontsize{#1}{#2pt}%
  \fontfamily{#3}\fontseries{#4}\fontshape{#5}%
  \selectfont}%
\fi\endgroup%
\begin{picture}(14287,5320)(-2421,-8483)
\put(4801,-4786){\makebox(0,0)[b]{\smash{{\SetFigFont{9}{10.8}{\rmdefault}{\mddefault}{\updefault}$E_{0}$}}}}
\put(1201,-4786){\makebox(0,0)[b]{\smash{{\SetFigFont{9}{10.8}{\rmdefault}{\mddefault}{\updefault}$E_{0} - C h$}}}}
\put(8401,-4786){\makebox(0,0)[b]{\smash{{\SetFigFont{9}{10.8}{\rmdefault}{\mddefault}{\updefault}$E_{0} + C h$}}}}
\put(11851,-4936){\makebox(0,0)[lb]{\smash{{\SetFigFont{9}{10.8}{\rmdefault}{\mddefault}{\updefault}$\re z$}}}}
\put(-899,-3286){\makebox(0,0)[lb]{\smash{{\SetFigFont{9}{10.8}{\rmdefault}{\mddefault}{\updefault}$\im z$}}}}
\put(8626,-7336){\makebox(0,0)[lb]{\smash{{\SetFigFont{9}{10.8}{\rmdefault}{\mddefault}{\updefault}$\im z = - \delta h$}}}}
\put(4801,-6286){\makebox(0,0)[b]{\smash{{\SetFigFont{9}{10.8}{\rmdefault}{\mddefault}{\updefault}\eqref{m74}}}}}
\end{picture}%
\end{center}
\caption{The resonance free zone given by Theorem \ref{a1}.} \label{f12}
\end{figure}

\begin{theorem}\sl \label{a1}
Assume \ref{h1}--\ref{h4} and \ref{h5}. Then, there exists $\delta > 0$ such that, for all $C >0$, $P$ has no resonance in
\begin{equation} \label{m74}
[ E_{0} - C h , E_{0} + C h] + i [ - \delta h , h ] ,
\end{equation}
for $h$ small enough. Moreover, for all $\chi \in C^{\infty}_{0} ( \R^{n} )$, there exists $M > 0$ such that
\begin{equation*}
\big\Vert \chi ( P -z )^{-1} \chi \big\Vert \lesssim h^{- M} ,
\end{equation*}
uniformly for $h$ small enough and $z \in \eqref{m74}$.
\end{theorem}

From Section \ref{s19} {\rm (A)}, it is not possible to hope a larger resonance free region under the assumptions of Theorem \ref{a1}. Using the semiclassical maximum principle of Tang and Zworski \cite[Lemma 2]{TaZw98_01} (see also Burq \cite[Lemma 4.7]{Bu04_01}), one can give a more accurate bound of the truncated resolvent on the real axis. More precisely, it is proved in Appendix \ref{b73} that Theorem \ref{a1} implies

\begin{corollary}\sl \label{a7}
Under the conditions of Theorem \ref{a1}, there exists $N > 0$ such that
\begin{equation*}
\big\Vert \chi ( P -z )^{-1} \chi \big\Vert \lesssim\frac{\vert \ln h \vert}{h} e^{N \vert \im z \vert \vert \ln h \vert / h} ,
\end{equation*}
uniformly for $h$ small enough and $z \in \eqref{m74}$.
\end{corollary}

For $z$ on the real axis, one can also replace the cut-off $\chi \in C^{\infty}_{0} ( \R^{n} )$ by a weight $\< x \>^{- s}$, with $s > 1/2$, as in the work of Robert and Tamura \cite{RoTa87_01}. In this regard, if we are only interested in estimating the weighted resolvent on the real axis, the analytic extension of the potential stated in \ref{h1} is not necessary. Assuming only that the potential is a symbol which decays at infinity, one can use the original approach of Burq \cite{Bu02_01} (see also Castella and Jecko \cite{CaJe06_01}) and modify the proof of Theorem \ref{a1} to obtain an upper bound for the weighted resolvent of order $h^{- N}$ on the real axis.

\begin{remark}\sl \label{j63}
$i)$ The assumption \ref{h5}, which guaranties that the trapping is not too strong, seems technical but the conclusions of Theorem \ref{a1} do not hold in general if this hypothesis is removed. Indeed, Section \ref{s19} {\rm (B)} provides an operator in dimension $n = 1$ which satisfies \ref{h1}--\ref{h4} but not \ref{h5} and whose resonances have an imaginary part of size $h \vert \ln  h \vert^{- 1}$.

$ii)$ Nevertheless, \ref{h5b} may be weakened. More precisely, the conclusions of Theorem \ref{a1} could probably be obtained under the hypothesis that $\Lambda_{-}$ and $\Lambda_{+}$ have a contact of finite order at each point of their intersection. This is done in Section \ref{s81} in dimension $n = 2$.
\end{remark}

One can be more explicit about the constant $\delta$ in Theorem \ref{a1}. Following the proof of this result (more precisely, the estimates \eqref{a39} and \eqref{a40} which quantify the gain after a ``turn''), it seems possible to take $\delta$ arbitrary close to
\begin{equation} \label{c14}
\delta_{1} = \sum_{j=1}^{n} \frac{\lambda_{j}}{2} - \frac{\lambda_{1}}{2} \sup_{\rho \in \CH} \dim \big( T_{\rho} \Lambda_{-} \cap T_{\rho} \Lambda_{+} \big) .
\end{equation}
Unfortunately, since (the uniqueness part of) Theorem \ref{a32} is only valid for $\im z > - \delta_{0} h$ for some $\delta_{0} > 0$, we have to restrict ourself to $\delta < \min ( \delta_{0} , \delta_{1} )$. Remark also that $\delta < \delta_{1}$ is optimal in the case $\CH = \emptyset$ (see Sj\"{o}strand \cite{Sj87_01}).

Finally, we give a typical example of operator $P$ which satisfies the assumptions of Theorem \ref{a1}.

\begin{figure}
\begin{center}
\begin{picture}(0,0)%
\includegraphics{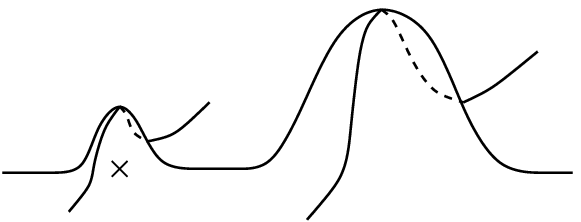}%
\end{picture}%
\setlength{\unitlength}{987sp}%
\begingroup\makeatletter\ifx\SetFigFont\undefined%
\gdef\SetFigFont#1#2#3#4#5{%
  \reset@font\fontsize{#1}{#2pt}%
  \fontfamily{#3}\fontseries{#4}\fontshape{#5}%
  \selectfont}%
\fi\endgroup%
\begin{picture}(11038,4120)(-11443,-4905)
\put(-6599,-1486){\makebox(0,0)[lb]{\smash{{\SetFigFont{9}{10.8}{\rmdefault}{\mddefault}{\updefault}$V ( x )$}}}}
\put(-9299,-4561){\makebox(0,0)[lb]{\smash{{\SetFigFont{9}{10.8}{\rmdefault}{\mddefault}{\updefault}$0$}}}}
\end{picture} $\qquad \qquad \qquad$ \begin{picture}(0,0)%
\includegraphics{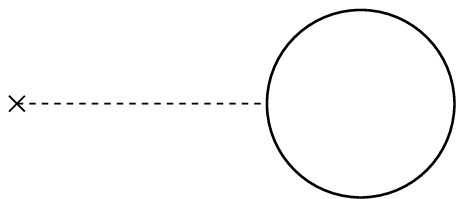}%
\end{picture}%
\setlength{\unitlength}{987sp}%
\begingroup\makeatletter\ifx\SetFigFont\undefined%
\gdef\SetFigFont#1#2#3#4#5{%
  \reset@font\fontsize{#1}{#2pt}%
  \fontfamily{#3}\fontseries{#4}\fontshape{#5}%
  \selectfont}%
\fi\endgroup%
\begin{picture}(9195,3658)(2836,-4740)
\put(5401,-1336){\makebox(0,0)[lb]{\smash{{\SetFigFont{9}{10.8}{\rmdefault}{\mddefault}{\updefault}$\{ V ( x ) = E_{0} \}$}}}}
\put(2851,-3061){\makebox(0,0)[lb]{\smash{{\SetFigFont{9}{10.8}{\rmdefault}{\mddefault}{\updefault}$0$}}}}
\put(5251,-3511){\makebox(0,0)[lb]{\smash{{\SetFigFont{9}{10.8}{\rmdefault}{\mddefault}{\updefault}$\pi_{x} ( \CH )$}}}}
\end{picture}%
\end{center}
\bigskip \bigskip
\begin{center}
\begin{picture}(0,0)%
\includegraphics{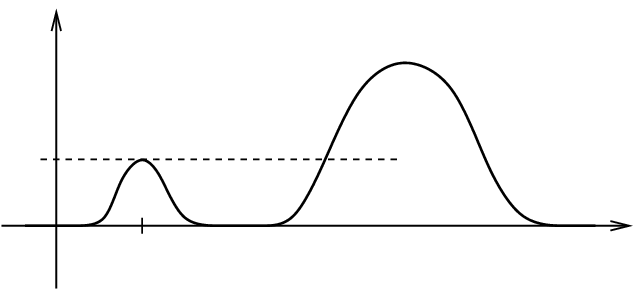}%
\end{picture}%
\setlength{\unitlength}{987sp}%
\begingroup\makeatletter\ifx\SetFigFont\undefined%
\gdef\SetFigFont#1#2#3#4#5{%
  \reset@font\fontsize{#1}{#2pt}%
  \fontfamily{#3}\fontseries{#4}\fontshape{#5}%
  \selectfont}%
\fi\endgroup%
\begin{picture}(12216,5466)(-11882,-5194)
\put(-9149,-4711){\makebox(0,0)[b]{\smash{{\SetFigFont{9}{10.8}{\rmdefault}{\mddefault}{\updefault}$0$}}}}
\put(-2549,-1711){\makebox(0,0)[lb]{\smash{{\SetFigFont{9}{10.8}{\rmdefault}{\mddefault}{\updefault}$V_{1} ( x_{1} )$}}}}
\put(-11249,-2836){\makebox(0,0)[rb]{\smash{{\SetFigFont{9}{10.8}{\rmdefault}{\mddefault}{\updefault}$E_{0}$}}}}
\end{picture} $\qquad \qquad$ \begin{picture}(0,0)%
\includegraphics{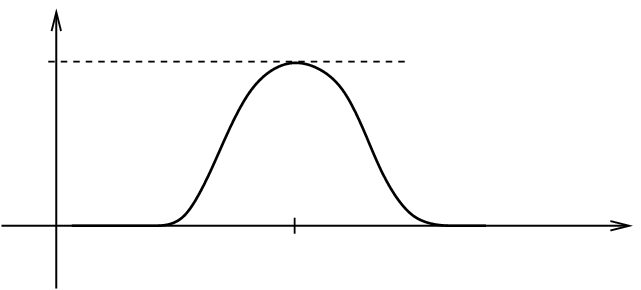}%
\end{picture}%
\setlength{\unitlength}{987sp}%
\begingroup\makeatletter\ifx\SetFigFont\undefined%
\gdef\SetFigFont#1#2#3#4#5{%
  \reset@font\fontsize{#1}{#2pt}%
  \fontfamily{#3}\fontseries{#4}\fontshape{#5}%
  \selectfont}%
\fi\endgroup%
\begin{picture}(12216,5466)(-11882,-5194)
\put(-6224,-4711){\makebox(0,0)[b]{\smash{{\SetFigFont{9}{10.8}{\rmdefault}{\mddefault}{\updefault}$0$}}}}
\put(-11249,-886){\makebox(0,0)[rb]{\smash{{\SetFigFont{9}{10.8}{\rmdefault}{\mddefault}{\updefault}$1$}}}}
\put(-4499,-2161){\makebox(0,0)[lb]{\smash{{\SetFigFont{9}{10.8}{\rmdefault}{\mddefault}{\updefault}$V_{2} ( x_{2} )$}}}}
\end{picture}%
\end{center}
\caption{The potential $V ( x )$ of Example \ref{c15}, the base space projection of $\CH$ and the choice of $V_{1}$ and $V_{2}$.} \label{f20}
\end{figure}

\begin{example}\rm \label{c15}
In dimension $n = 2$, we consider
\begin{equation*}
V (x) = V_{1} ( x_{1} ) V_{2} ( x_{2} ) ,
\end{equation*}
where the functions $V_{\bullet} \in C^{\infty}_{0} ( \R )$ are such as in Figure \ref{f20}. In particular, they satisfy
\begin{equation*}
V_{1} ( x_{1} ) = E_{0} - \frac{\lambda_{1}^{2}}{2} x_{1}^{2} + \CO ( x_{1}^{3} )   \qquad \text{and} \qquad   V_{2} ( x_{2} ) = 1 - \frac{\lambda_{2}^{2}}{2 E_{0}} x_{2}^{2} + \CO ( x_{2}^{3} ) ,
\end{equation*}
near $0$. Then, \ref{h1} and \ref{h2} are verified. Moreover, $K ( E_{0} ) = \{ ( 0 , 0 ) \} \cup \CH$ where $\CH \subset \{ x_{2} = 0 \}$ consists of a unique homoclinic curve. Thus, the assumption \ref{h3} holds true. For $\rho \in \CH$, we have $g_{\pm} ( \rho ) = ( 1 , 0 )$ and \ref{h4} follows (see Figure \ref{f20}). Finally, one can verify that \ref{h5a} is satisfied.

Changing the form of the ``croissant'' in the forthcoming Example \ref{b80}, it is also possible to construct a potential $V$ satisfying \ref{h1}--\ref{h4} and \ref{h5b} but not \ref{h5a}.
\end{example}

\Subsection{Strong trapping} \label{s77}

In this part, we prove a smaller resonance free zone when the trapping is stronger. Since the setting \ref{h5a} has been treated in Theorem \ref{a1}, we now assume that
\begin{hyp} \label{h6}
We are in the isotropic case: $\lambda_{1} = \dots = \lambda_{n} = : \lambda$.
\end{hyp}

We next define an operator $\CT_{0}$, acting on the classical quantities, which measures the ``decay'' of the solutions of $( P - z ) u = 0$ after a turn around the trapped set (through the fixed point and along the homoclinic curves). The proof of the following assertions are given in Section \ref{a71}. We set
\begin{equation*}
\CH_{\rm tang} = \{ \rho \in \CH ; \ T_{\rho} \Lambda_{-} = T_{\rho} \Lambda_{+} \} ,
\end{equation*}
and
\begin{equation*}
\CH_{\rm tang}^{\pm \infty} = \Big\{ \frac{g_{\pm} ( \rho )}{\vert g_{\pm} ( \rho ) \vert} ; \ \rho \in \CH_{\rm tang} \Big\} \subset \S^{n-1} ,
\end{equation*}
the asymptotic directions of $\CH_{\rm tang}$ at the critical point $( 0 , 0 )$. In some sense, $\CH_{\rm tang}^{\pm \infty}$ is a quotient of $\CH_{\rm tang}$ by the Hamiltonian flow. The set $\CH_{\rm tang}^{\pm \infty}$ is compact and, since we consider a Schr\"{o}dinger operator, we have $\CH_{\rm tang}^{+ \infty} = \CH_{\rm tang}^{- \infty}$.

There exists a smooth function $x_{+} ( t , \alpha ) : \R \times \S^{n-1} \longrightarrow \R^{n}$ such that
\begin{equation*}
x_{+} ( t , \alpha ) = e^{\lambda t} \alpha + o \big( e^{\lambda t} \big) ,
\end{equation*}
as $t \to - \infty$, and such that, for all $\alpha \in \S^{n-1}$,
\begin{equation} \label{c17}
t \longmapsto ( x_{+} ( t , \alpha ) , \nabla \varphi_{+} ( x_{+} ( t , \alpha ) ) ) ,
\end{equation}
is a Hamiltonian curve. In particular, note that $\Lambda_{+} \setminus \{ 0 \}$ is the disjoint union over $\alpha \in \S^{n-1}$ of these curves.

\begin{figure}
\begin{center}
\begin{picture}(0,0)%
\includegraphics{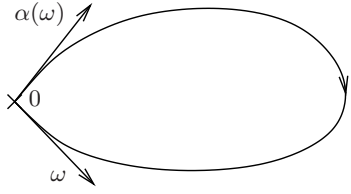}%
\end{picture}%
\setlength{\unitlength}{1184sp}%
\begingroup\makeatletter\ifx\SetFigFont\undefined%
\gdef\SetFigFont#1#2#3#4#5{%
  \reset@font\fontsize{#1}{#2pt}%
  \fontfamily{#3}\fontseries{#4}\fontshape{#5}%
  \selectfont}%
\fi\endgroup%
\begin{picture}(7168,3966)(3418,-5794)
\put(4351,-5611){\makebox(0,0)[lb]{\smash{{\SetFigFont{9}{10.8}{\rmdefault}{\mddefault}{\updefault}$\omega$}}}}
\put(3901,-4036){\makebox(0,0)[lb]{\smash{{\SetFigFont{9}{10.8}{\rmdefault}{\mddefault}{\updefault}$0$}}}}
\put(3601,-2236){\makebox(0,0)[lb]{\smash{{\SetFigFont{9}{10.8}{\rmdefault}{\mddefault}{\updefault}$\alpha ( \omega )$}}}}
\end{picture}%
\end{center}
\caption{The definition of $\alpha ( \omega )$.} \label{f9}
\end{figure}

For $\alpha \in \CH_{\rm tang}^{+ \infty}$, we define $\omega ( \alpha ) \in \CH_{\rm tang}^{- \infty}$ as the normalized asymptotic direction of $x_{+} ( t ,\alpha )$ as $t \to + \infty$. Then, $\omega ( \cdot )$ is a homeomorphism from $\CH_{\rm tang}^{+ \infty}$ onto $\CH_{\rm tang}^{- \infty}$ and $\alpha ( \cdot ) : \CH_{\rm tang}^{- \infty} \longrightarrow \CH_{\rm tang}^{+ \infty}$ denotes its inverse (see Figure \ref{f9}).

For all $\alpha \in \CH_{\rm tang}^{+ \infty}$ and $0 < \varepsilon \leq \varepsilon_{0}$ (with $\varepsilon_{0}$ independent of $\alpha$), the characteristic curve $x_{+} ( t , \alpha )$ meets $\Lambda_{\pm}^{0} \cap \{ \vert x \vert = \varepsilon \}$ only once, at time $t_{\pm}^{\varepsilon} ( \alpha )$. Moreover, the Maslov determinant
\begin{equation} \label{a81}
\CM_{\varepsilon} ( \alpha ) = \frac{\sqrt{\Big\vert \det \frac{\partial x_{+} ( t , \beta )}{\partial ( t , \beta )} \vert_{t = t_{+}^{\varepsilon} ( \alpha ) , \ \beta = \alpha} \Big\vert}}{\sqrt{\Big\vert \det \frac{\partial x_{+} ( t , \beta)}{\partial ( t , \beta )} \vert_{t = t_{-}^{\varepsilon} ( \alpha ) , \ \beta = \alpha} \Big\vert}} ,
\end{equation}
is well-defined for all $( \varepsilon , \alpha ) \in ] 0 , \varepsilon_{0} ] \times \CH_{\rm tang}^{+ \infty}$. The limit
\begin{equation} \label{i49}
\CM_{0} ( \alpha ) = \lim_{\varepsilon \to 0} \CM_{\varepsilon} ( \alpha ) ,
\end{equation}
exists uniformly in $\alpha \in \CH_{\rm tang}^{+ \infty}$. Furthermore, the map $( \varepsilon , \alpha ) \longmapsto \CM_{\varepsilon} ( \alpha )$ is continuous on $[ 0 , \varepsilon_{0} ] \times \CH_{\rm tang}^{+ \infty}$ and there exists $C > 0$ such that $1/ C \leq \CM_{\varepsilon} ( \alpha ) \leq C$ on this set. Roughly speaking, $\CM_{0}$ measures the amplification along the homoclinic curves.

\begin{remark}\sl \label{a82}
In some sense, we have
\begin{equation} \label{a83}
\CM_{0} ( \alpha ) = \sqrt{\Big\vert \det \frac{\partial \omega ( \alpha )}{\partial \alpha} \Big\vert^{- 1}} .
\end{equation}
But, since $\CH_{\rm tang}^{+ \infty}$ is only a compact set, the right hand side in the previous expression has, in general, no meaning. However, this formula holds true in the interior of $\CH_{\rm tang}^{+ \infty}$. See the discussion at the end of Appendix \ref{a71} for more details.
\end{remark}

We are now in position to define the limit operator which governs the geometry. For $\tau \in \R$, $\CT_{0} ( \tau )$ denotes the operator on $L^{\infty} ( \CH_{\rm tang}^{- \infty} )$, endowed with the Lesbegue measure on $\S^{n-1}$, with kernel
\begin{equation} \label{b70}
\CT_{0} ( \tau ) ( \omega , \widetilde{\omega} ) : = ( 2 \pi )^{- \frac{n}{2}} \CM_{0} ( \alpha ( \omega ) ) \big\vert \alpha ( \omega ) \cdot \widetilde{\omega} \big\vert^{- \frac{n}{2}} e^{- \frac{\pi \tau}{2 \lambda} \sgn ( \alpha ( \omega ) \cdot \widetilde{\omega} )} \Big\vert \Gamma \Big( \frac{n}{2} - i \frac{\tau}{\lambda} \Big) \Big\vert .
\end{equation}
From \ref{h4}, the properties of $\CM_{0}$ and the compactness of $\CH_{\rm tang}^{\pm \infty}$, this kernel is continuous on $\CH_{\rm tang}^{- \infty} \times \CH_{\rm tang}^{- \infty}$ and $\CT_{0} ( \tau )$ is bounded. Measuring the decay after a turn around the trapped set of the solutions of $( P - z ) u = 0$ with $\re z \approx E_{0} + \tau h$, this operator is essentially the modulus of a quantization relation. Eventually, we define
\begin{equation} \label{b71}
\CA_{0} ( \tau ) = \spr ( \CT_{0} ( \tau ) ),
\end{equation}
where
\begin{equation*}
\spr ( T ) = \lim_{k \to + \infty} \Vert T^{k} \Vert^{\frac{1}{k}} = \max_{\mu \in \spe ( T )} \vert \mu \vert ,
\end{equation*}
is the spectral radius of an operator $T$ on a Banach space. Recall that the spectral radius is upper semicontinuous (i.e. if $T_{n} \to T$ in norm then $\limsup_{n \to + \infty} \spr ( T_{n} ) \leq \spr (T)$) but, in general, not continuous (see the example of Kakutani in Halmos \cite[Problem 87]{Ha67_01}). However in our case, $\tau \longmapsto \CA_{0} ( \tau )$ is continuous on $\R$ from Proposition \ref{c7}. Note also that $\CT_{0}$ is a compact operator on all the Banach spaces $L^{q} ( \CH_{\rm tang}^{- \infty} )$ for $1 \leq q \leq + \infty$. Moreover, its spectrum and its spectral radius are independent of $q$.

\begin{remark}\sl \label{b72}
It may not be easy to compute $\CA_{0}$ in the applications. Nevertheless, one can easily estimate this quantity. Let
\begin{gather*}
\CM_{0} : = \max_{\alpha \in \CH_{\rm tang}^{+ \infty}} \CM_{0} ( \alpha ) \qquad \text{and} \qquad \CJ_{0} ( \tau ) : = \max_{\alpha \in \CH_{\rm tang}^{+ \infty}} \int_{\CH_{\rm tang}^{- \infty}} e^{- \frac{\pi \tau}{2 \lambda} \sgn ( \alpha \cdot \omega )} \vert \alpha \cdot \omega \vert^{- \frac{n}{2}} d \omega .
\end{gather*}
which are bounded from the previous discussion. Then,
\begin{equation*}
\CA_{0} ( \tau ) \leq \big\Vert \CT_{0} ( \tau ) \big\Vert \leq ( 2 \pi )^{- \frac{n}{2}} \Big\vert \Gamma \Big( \frac{n}{2} - i \frac{\tau}{\lambda} \Big) \Big\vert \CM_{0} \CJ_{0} ( \tau ).
\end{equation*}
\end{remark}

Under the previous assumptions, our main result is the following.

\begin{figure}
\begin{center}
\begin{picture}(0,0)%
\includegraphics{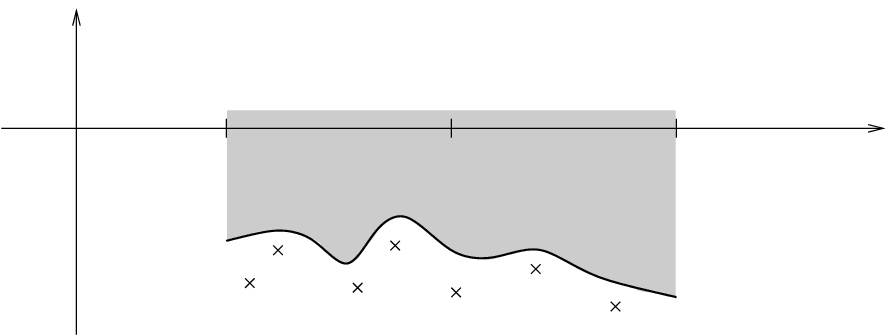}%
\end{picture}%
\setlength{\unitlength}{1184sp}%
\begingroup\makeatletter\ifx\SetFigFont\undefined%
\gdef\SetFigFont#1#2#3#4#5{%
  \reset@font\fontsize{#1}{#2pt}%
  \fontfamily{#3}\fontseries{#4}\fontshape{#5}%
  \selectfont}%
\fi\endgroup%
\begin{picture}(14287,5320)(-2421,-8483)
\put(4801,-6136){\makebox(0,0)[b]{\smash{{\SetFigFont{9}{10.8}{\rmdefault}{\mddefault}{\updefault}\eqref{b68}}}}}
\put(4801,-4786){\makebox(0,0)[b]{\smash{{\SetFigFont{9}{10.8}{\rmdefault}{\mddefault}{\updefault}$E_{0}$}}}}
\put(1201,-4786){\makebox(0,0)[b]{\smash{{\SetFigFont{9}{10.8}{\rmdefault}{\mddefault}{\updefault}$E_{0} - C h$}}}}
\put(8401,-4786){\makebox(0,0)[b]{\smash{{\SetFigFont{9}{10.8}{\rmdefault}{\mddefault}{\updefault}$E_{0} + C h$}}}}
\put(-899,-3286){\makebox(0,0)[lb]{\smash{{\SetFigFont{9}{10.8}{\rmdefault}{\mddefault}{\updefault}$\im z$}}}}
\put(11851,-4936){\makebox(0,0)[lb]{\smash{{\SetFigFont{9}{10.8}{\rmdefault}{\mddefault}{\updefault}$\re z$}}}}
\put(8551,-7906){\makebox(0,0)[lb]{\smash{{\SetFigFont{9}{10.8}{\rmdefault}{\mddefault}{\updefault}$\im z = \lambda \ln \Big( \CA_{0} \Big( \frac{\re z - E_{0}}{h} \Big) \Big) \frac{h}{\vert \ln h \vert}$}}}}
\end{picture}%
\end{center}
\caption{The resonance free zone given by Theorem \ref{a2}.} \label{f11}
\end{figure}

\begin{theorem}\sl \label{a2}
Assume \ref{h1}--\ref{h4} and \ref{h6}. Then, for all $C , \delta > 0$, $P$ has no resonance in the domain
\begin{equation} \label{b68}
\left\{ \begin{aligned}
&E_{0} - C h \leq \re z \leq E_{0} + C h ,    \\
&\Big( \lambda \ln \Big( \CA_{0} \Big( \frac{\re z - E_{0}}{h} \Big) \Big) + \delta \Big) \frac{h}{\vert \ln h \vert} \leq \im z \leq \frac{h}{\vert \ln h \vert} ,
\end{aligned} \right.
\end{equation}
for $h$ small enough. Moreover, for all $\chi \in C^{\infty}_{0} ( \R^{n} )$, there exists $M > 0$ such that
\begin{equation*}
\big\Vert \chi ( P -z )^{-1} \chi \big\Vert \lesssim h^{- M} ,
\end{equation*}
uniformly for $h$ small enough and $z \in \eqref{b68}$.
\end{theorem}

The domain \eqref{b68} is drawn in Figure \ref{f11}. We now give a typical example of an operator illustrating Theorem \ref{a2}.

\begin{figure}
\begin{center}
\begin{picture}(0,0)%
\includegraphics{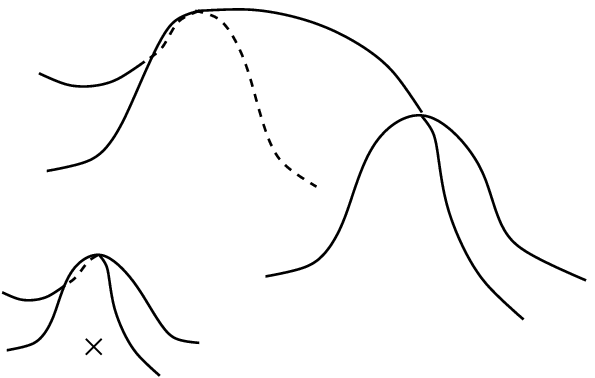}%
\end{picture}%
\setlength{\unitlength}{987sp}%
\begingroup\makeatletter\ifx\SetFigFont\undefined%
\gdef\SetFigFont#1#2#3#4#5{%
  \reset@font\fontsize{#1}{#2pt}%
  \fontfamily{#3}\fontseries{#4}\fontshape{#5}%
  \selectfont}%
\fi\endgroup%
\begin{picture}(11294,7122)(-10199,-6661)
\put(-2174,-361){\makebox(0,0)[lb]{\smash{{\SetFigFont{9}{10.8}{\rmdefault}{\mddefault}{\updefault}$V (x)$}}}}
\put(-8699,-6211){\makebox(0,0)[rb]{\smash{{\SetFigFont{9}{10.8}{\rmdefault}{\mddefault}{\updefault}$0$}}}}
\end{picture} $\qquad \qquad$ \begin{picture}(0,0)%
\includegraphics{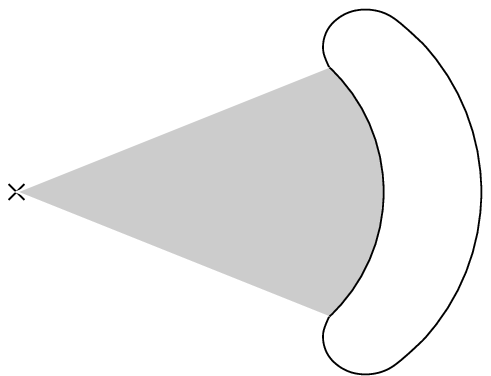}%
\end{picture}%
\setlength{\unitlength}{987sp}%
\begingroup\makeatletter\ifx\SetFigFont\undefined%
\gdef\SetFigFont#1#2#3#4#5{%
  \reset@font\fontsize{#1}{#2pt}%
  \fontfamily{#3}\fontseries{#4}\fontshape{#5}%
  \selectfont}%
\fi\endgroup%
\begin{picture}(9339,7070)(-8789,-8696)
\put(-2699,-2161){\makebox(0,0)[rb]{\smash{{\SetFigFont{9}{10.8}{\rmdefault}{\mddefault}{\updefault}$\{ V(x) = E_{0} \}$}}}}
\put(-8774,-5236){\makebox(0,0)[rb]{\smash{{\SetFigFont{9}{10.8}{\rmdefault}{\mddefault}{\updefault}$0$}}}}
\put(-6899,-5236){\makebox(0,0)[lb]{\smash{{\SetFigFont{9}{10.8}{\rmdefault}{\mddefault}{\updefault}$2 \theta_{0}$}}}}
\put(-3299,-5236){\makebox(0,0)[rb]{\smash{{\SetFigFont{9}{10.8}{\rmdefault}{\mddefault}{\updefault}$\pi_{x} ( {\mathcal H} )$}}}}
\end{picture}%
\end{center}
\caption{The potential of Example \ref{b80} and the base space projection of $\CH$.} \label{f13}
\end{figure}

\begin{example}\rm \label{b80}
In dimension $n = 2$, let $(r , \theta )$ be the polar coordinates. We consider
\begin{equation*}
V (x) = V_{0} ( r ) + V_{1} ( r - a ) \psi ( \theta ) ,
\end{equation*}
where the $V_{\bullet} ( r )$'s are even functions in $C_{0}^{\infty} ( \R )$ satisfying $r V_{\bullet}^{\prime} (r)<0$ for $r$ in the interior of $\supp V_{\bullet} \setminus \{ 0 \}$ and $E_{0} = V_{0} ( 0 ) < V_{1} ( 0 )$. Moreover, $V_{0}$ verifies \ref{h2}. The constant $a > 0$ is chosen sufficiently large such that $\supp V_{0} (r) \cap \supp V_{1} (r-a) = \emptyset$. Eventually, $\psi ( \theta ) \in C_{0}^{\infty} ( [ - \theta_{1} , \theta_{1} ] )$ is equal to $1$ for $\vert \theta \vert \leq \theta_{0}$ and $\theta \psi^{\prime} ( \theta ) < 0$ for $\theta_{0} < \vert \theta \vert < \theta_{1}$ for some $0 \leq \theta_{0} < \theta_{1} \leq \pi$. The setting is illustrated in Figure \ref{f13}.

If $\theta_{1} < \pi / 4$, it can be checked that the conditions \ref{h1}--\ref{h4} and \ref{h6} are all satisfied and that $\CH^{\infty} = \CH^{\pm \infty}_{\rm tang} = [ - \theta_{0} , \theta_{0} ]$. On the contrary, \ref{h4} is clearly not verified for $\theta_{0} \geq \pi / 4$. Moreover, using the notations of Remark \ref{b72}, we have
\begin{equation*}
\CM_{0} = 1 \qquad \text{and} \qquad \CJ_{0} ( \tau ) = e^{- \frac{\pi \tau}{2 \lambda}} \ln \tan \Big( \theta_{0} + \frac{\pi}{4} \Big) .
\end{equation*}
Thus, Theorem \ref{a2} provides an explicit resonance free zone in this example. Note also that $g_{-} \cdot g_{+} > 0$ for all $g_{\pm} \in \CH^{\pm \infty}_{\rm tang}$. The non-trapping potential $V_{1} ( r - a ) \psi ( \theta )$ has been considered for the Helmholtz equation by Castella and Klak in \cite[Section 1.3]{CaKl14_01}.
\end{example}

\begin{remark}\sl \label{b69}
$i)$ Since the kernel of $\CT_{0} ( \tau )$ is bounded from below by a positive constant, $\CA_{0}$ vanishes at some point if and only if $\mes_{\S^{n - 1}} ( \CH_{\rm tang}^{- \infty} ) = 0$. In this case, $\CA_{0}$ vanishes  on the whole real axis and the previous theorem must be read as follows: For all $C > 0$, $P$ has no resonance in the set
\begin{equation*}
[ E_{0} - C h , E_{0} + C h ] + i \Big[ - C \frac{h}{\vert \ln h \vert} , \frac{h}{\vert \ln h \vert} \Big] ,
\end{equation*}
for $h$ small enough. Moreover, the truncated resolvent of $P$ satisfies a polynomial estimate in this set. Note that the Lebesgue measure of $\CH_{\rm tang}^{- \infty}$ is zero, in particular, when $\CH$ is the union of an at most countable number of Hamiltonian curves in dimension $n \geq 2$. Example \ref{b86} provides an instance of such situation.

$ii)$ On the contrary, the statement of the above theorem is empty in the region where $\CA_{0} ( ( \re z - E_{0} ) / h ) \geq 1$. Section \ref{s19} {\rm (C)} provides an example of such a situation. Remark moreover that, in this example, the resonances have exponentially small imaginary part. Thus, the assumptions of Theorem \ref{a2} are not enough to guarantee a resonance free zone of size $h \vert \ln h \vert^{- 1}$ in general.
\end{remark}

\begin{example}\rm \label{b86}
We construct here an example of operator $P$ satisfying the assumptions \ref{h2}--\ref{h4} and such that $\CH = \CH_{\rm tang}$ is a countable set of homoclinic curves. It will not be of the form \eqref{a5}, but, using Remark \ref{c13} and the propagation of singularities for boundary value problems, one can hope to extend Theorem \ref{a2} to this case.

\begin{figure}
\begin{center}
\begin{picture}(0,0)%
\includegraphics{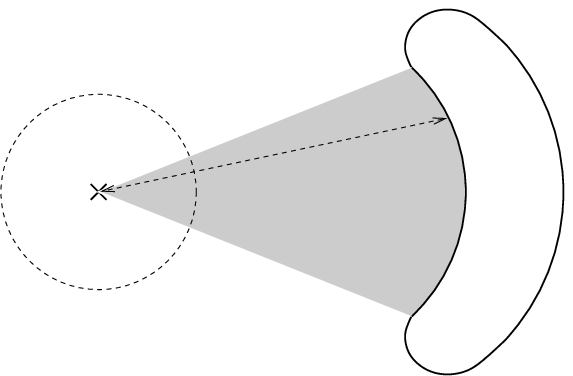}%
\end{picture}%
\setlength{\unitlength}{987sp}%
\begingroup\makeatletter\ifx\SetFigFont\undefined%
\gdef\SetFigFont#1#2#3#4#5{%
  \reset@font\fontsize{#1}{#2pt}%
  \fontfamily{#3}\fontseries{#4}\fontshape{#5}%
  \selectfont}%
\fi\endgroup%
\begin{picture}(10839,7070)(-10289,-8696)
\put(-2549,-3736){\makebox(0,0)[rb]{\smash{{\SetFigFont{9}{10.8}{\rmdefault}{\mddefault}{\updefault}$R_{0}$}}}}
\put(-8399,-2986){\makebox(0,0)[b]{\smash{{\SetFigFont{9}{10.8}{\rmdefault}{\mddefault}{\updefault}$\supp V$}}}}
\put(-8774,-5236){\makebox(0,0)[rb]{\smash{{\SetFigFont{9}{10.8}{\rmdefault}{\mddefault}{\updefault}$0$}}}}
\put(-479,-5236){\makebox(0,0)[b]{\smash{{\SetFigFont{9}{10.8}{\rmdefault}{\mddefault}{\updefault}$\CO_{0}$}}}}
\put(-2999,-5236){\makebox(0,0)[rb]{\smash{{\SetFigFont{9}{10.8}{\rmdefault}{\mddefault}{\updefault}$\pi_{x} ( {\mathcal H} )$}}}}
\end{picture} $\qquad \qquad$ \begin{picture}(0,0)%
\includegraphics{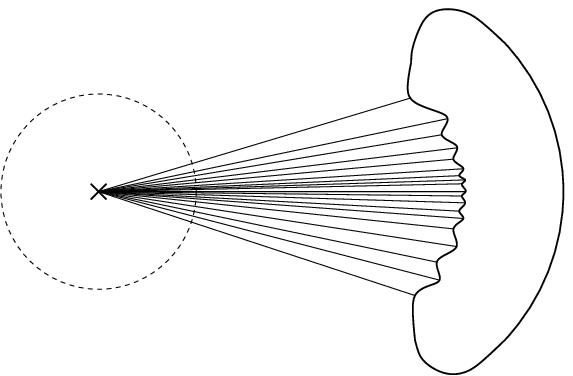}%
\end{picture}%
\setlength{\unitlength}{987sp}%
\begingroup\makeatletter\ifx\SetFigFont\undefined%
\gdef\SetFigFont#1#2#3#4#5{%
  \reset@font\fontsize{#1}{#2pt}%
  \fontfamily{#3}\fontseries{#4}\fontshape{#5}%
  \selectfont}%
\fi\endgroup%
\begin{picture}(10839,7073)(-10289,-8699)
\put(-8399,-2986){\makebox(0,0)[b]{\smash{{\SetFigFont{9}{10.8}{\rmdefault}{\mddefault}{\updefault}$\supp V$}}}}
\put(-479,-5236){\makebox(0,0)[b]{\smash{{\SetFigFont{9}{10.8}{\rmdefault}{\mddefault}{\updefault}$\CO$}}}}
\put(-8774,-5236){\makebox(0,0)[rb]{\smash{{\SetFigFont{9}{10.8}{\rmdefault}{\mddefault}{\updefault}$0$}}}}
\put(-4124,-3511){\makebox(0,0)[rb]{\smash{{\SetFigFont{9}{10.8}{\rmdefault}{\mddefault}{\updefault}$\pi_{x} ( {\mathcal H} )$}}}}
\end{picture}%
\end{center}
\caption{The obstacles of Example \ref{b86} and the corresponding homoclinic sets.} \label{f16}
\end{figure}

In dimension $n = 2$, let $(r , \theta )$ be the polar coordinates. First let $V (r)$ be a smooth potential as in Section \ref{s21}. We then consider
\begin{equation} \label{i79}
P_{0} = - h^{2} \Delta_{\R^{2} \setminus \CO_{0}} + V (r) ,
\end{equation}
with Dirichlet condition at the boundary of the obstacle $\CO_{0}$ which is a piece of a ring as illustrated in Figure \ref{f16}. As in Example \ref{b80}, the homoclinic set at energy $E_{0}$ consists of the rays whose angle is in a certain sector. We now perturb $\CO_{0}$ to obtain the required geometric setting. Let $\CO$ be the smooth non-trapping obstacle which coincides with $\CO_{0}$ except that the part of the boundary $r = R_{0}$ is replaced by
\begin{equation*}
r = R_{0} + F ( \theta ) ,
\end{equation*}
where $F \in C^{\infty}$ is a well chosen function such that $F ( 0 ) = 0$ and $F$ is a primitive of
\begin{equation} \label{b87}
F^{\prime} ( \theta ) = \Big( \sin \Big( \frac{1}{\theta} \Big) e^{- \frac{1}{\theta^{2}}} \Big)^{3} ,
\end{equation}
near $0$ (see Figure \ref{f16}). We then define
\begin{equation*}
P = - h^{2} \Delta_{\R^{2} \setminus \CO} + V (r) .
\end{equation*}
Using Proposition \ref{a14}, one can prove that the trapped set of $P$ at energy $E_{0}$ verifies \ref{h3} and that the homoclinic set consists of the radial rays whose angle $\theta$ satisfies $F^{\prime} ( \theta ) = 0$. Thus, \eqref{b87} implies that $\CH$ is a countable set of homoclinic curves which accumulate only at $\theta = 0$. In particular, $\mes_{\S^{n - 1}} ( \CH_{\rm tang}^{\pm \infty} ) = 0$. Moreover, since $F^{\prime} ( \theta ) = 0$ implies $F^{\prime \prime} ( \theta ) = 0$, $\Lambda_{-}$ and $\Lambda_{+}$ intersect transversally along each homoclinic curve. In other words, $\CH = \CH_{\rm tang}$.
\end{example}

As stated before, $\CT_{0} ( \tau ) ( \omega , \widetilde{\omega} )$ is the modulus of the kernel of some ``quantization operator'' $\CT ( \tau )$ defined in \eqref{c19}. Under some additional assumptions, it is possible to use this operator, to obtain sharp resonance free zones and even the asymptotic of the resonances in any vicinity of size $h \vert \ln h \vert^{- 1}$ of the real axis. This is done in Section \ref{s26}.

\begin{remark}\sl \label{b77}
Under the assumptions of Theorem \ref{a2}, there is no hope to have in general a resonance free zone of order greater than $h \vert \ln h \vert^{- 1}$. Indeed, Section \ref{s19} {\rm (B)} provides an operator, satisfying the required hypotheses, which has resonances with imaginary part of order $h \vert \ln h \vert^{- 1}$. Moreover, using the computation of $\CA_{0}$ made in \eqref{b75}, Theorem \ref{a2} shows that this operator has no resonance $z$ with
\begin{equation*}
\Big( - \lambda \frac{\ln 2}{2} + \delta \Big) \frac{h}{\vert \ln h \vert} \leq \im z \leq 0 ,
\end{equation*}
and $\re z = E_{0} + o ( h )$ for any $\delta > 0$. Otherwise, from \eqref{e39}, the resonances $z$ in the set $B ( E_{0} , o ( h ) )$ satisfy
\begin{equation*}
\im z =  - \lambda \frac{\ln 2}{2} \frac{h}{\vert \ln h \vert} + o \Big( \frac {h}{\vert \ln h \vert} \Big) .
\end{equation*}
Thus, the resonance free zone given by our result is sharp in this case.
\end{remark}

As in Corollary \ref{a7}, the behavior of the norm of the truncated resolvent can be specified near the real axis.

\begin{corollary}\sl \label{a8}
Assume that $P$ has no resonance in
\begin{equation*}
[ a (h) , b (h) ] + i \Big[ - c \frac{h}{\vert \ln h \vert} , \frac{h}{\vert \ln h \vert} \Big] + B \Big( 0 , d \frac{h}{\vert \ln h \vert} \Big) ,
\end{equation*}
with $a(h) < b (h)$, $c, d > 0$ and that the norm of $\chi ( P -z )^{-1} \chi$ satisfies a polynomial estimate in this set. Then, there exists $N > 0$ such that
\begin{equation*}
\forall z \in [ a (h) , b (h) ] + i \Big[ - c \frac{h}{\vert \ln h \vert} , \frac{h}{\vert \ln h \vert} \Big] , \qquad \big\Vert \chi ( P -z )^{-1} \chi \big\Vert \lesssim\frac{\vert \ln h \vert^{2}}{h} e^{N \vert \im z \vert \vert \ln h \vert^{2} / h} ,
\end{equation*}
for $h$ small enough.
\end{corollary}

\section{Asymptotic of the resonances generated by a finite number of homoclinic trajectories} \label{s6}

In this part, we give the precise localization of the resonances in the most natural case where the homoclinic set consists in a finite number of trajectories. First, we state the main results in Section \ref{s61}. Their proofs can be found in Section \ref{s12}. The rest of Section  \ref{s6} is devoted to the description of particular phenomena, additional results and examples. The corresponding proofs are given in Section \ref{s23}.

\Subsection{Main results} \label{s61}

In addition to the assumptions of Section \ref{s2}, we suppose that
\begin{hyp} \label{h8}
The homoclinic set $\CH$ consists of a finite number of trajectories on which $\Lambda_{-}$ and $\Lambda_{+}$ intersect transversally.
\end{hyp}
Here, transversally means that $T_{\rho} ( \Lambda_{-} \cap \Lambda_{+} ) = T_{\rho} \Lambda_{-} \cap T_{\rho} \Lambda_{+}$ (which coincides with $\R H_{p} ( \rho )$) for all $\rho \in \CH$. Remark that the transversality of the intersection implies that the number of homoclinic trajectories is finite. We denote by $\gamma_{k}$, $k = 1 , \ldots , K$, these Hamiltonian curves. To each of them, we associate its action
\begin{equation} \label{r5}
A_{k} = \int_{\gamma_{k}} \xi \cdot d x ,
\end{equation}
and the Maslov index $\nu_{k}$ of $\Lambda_{+}$ along $\gamma_{k}$. Note that these two quantities are well-defined (and finite) thanks to \cite[Proposition C.1]{ALBoRa08_01}. We also set the notations for the asymptotic directions of $\gamma_{k}$. Let $\gamma_{k} (t) = ( x_{k} ( t ) , \xi_{k} ( t ) )$ be a parametrization of the Hamiltonian trajectory $\gamma_{k}$. Following \eqref{d2}, the vectors $g_{\pm}^{k} \in \R^{n}$ are defined by
\begin{equation} \label{d6}
x_{k} (t) = g_{\pm}^{k} e^{\pm \lambda_{1} t} + o \big( e^{\pm \lambda_{1} t} \big) ,
\end{equation}
as $t \to \mp \infty$. Eventually, we define the asymptotic Maslov determinants associated to these curves. Let $\gamma_{k} ( t , y ) = ( x_{k} ( t , y) , \xi_{k} ( t , y ) ) : \R \times \R^{n-1} \longrightarrow T^{*} \R^{n}$, defined from a neighborhood of $\R \times \{ 0 \}$ to a neighborhood of $\gamma_{k}$, be a smooth parametrization of $\Lambda_{+}$ by Hamiltonian curves such that $\gamma_{k} ( t , 0 ) = \gamma_{k} ( t )$. Then, the limits
\begin{equation} \label{d7}
\begin{aligned}
\CM_{k}^{+} & = \lim_{s \to - \infty} \sqrt{\Big\vert \det \frac{\partial x_{k} ( t , y )}{\partial ( t , y )} \vert_{t = s , \ y = 0} \Big\vert} e^{- s \sum_{j} \lambda_{j} / 2} ,   \\
\CM_{k}^{-} & = \lim_{s \to + \infty} \sqrt{\Big\vert \det \frac{\partial x_{k} ( t , y )}{\partial ( t , y )} \vert_{t = s , \ y = 0} \Big\vert} e^{- s ( \sum_{j} \lambda_{j} - 2 \lambda_{1} ) / 2} ,
\end{aligned}
\end{equation}
exist and belong to $] 0 , + \infty [$ (see \cite{ALBoRa08_01,BoFuRaZe07_01}). Remark that the present setting enters the framework of Theorem \ref{a1} when the dimension is greater than $1$.

As we will see below, the asymptotic of the resonances is governed by the $K \times K$ matrix $\CQ$ whose entries are given by
\begin{equation} \label{d4}
\CQ_{k , \ell} ( z , h ) = e^{i A_{k} / h} \Gamma \big( S ( z , h ) / \lambda_{1} \big) \sqrt{\frac{\lambda_{1}}{2 \pi}} \frac{\CM_{k}^{+}}{\CM_{k}^{-}} e^{- \frac{\pi}{2} ( \nu_{k} + \frac{1}{2} ) i} \big\vert g_{-}^{\ell} \big\vert \big( i \lambda_{1} g_{+}^{k} \cdot g_{-}^{\ell} \big)^{- S ( z , h ) / \lambda_{1}} ,
\end{equation}
with the convention $( i a )^{b} = \vert a \vert^{b} e^{i \frac{\pi}{2} \sgn ( a ) b}$ for $a \in \R$, $b \in \C$ and
\begin{equation} \label{d5}
S ( z , h ) = \sum_{j = 1}^{n} \frac{\lambda_{j}}{2} - i \frac{z - E_{0}}{h} .
\end{equation}
It is important to note that $\CQ ( z , h )$ can be written as
\begin{equation} \label{d10}
\CQ ( z , h ) = \widetilde{\CQ} ( \rho , \sigma ) = \sum_{k = 1}^{K} \rho_{k} \widetilde{\CQ}_{k} ( \sigma ) ,
\end{equation}
where $\rho = ( \rho_{1} , \ldots , \rho_{K} ) \in \C^{K}$ and $\rho_{k} = e^{i A_{k} / h} \in \S^{1}$ is periodic with respect to $h^{- 1}$. The matrix $\widetilde{\CQ}_{k}$ is independent of $h$ and depends meromorphically on $z$ only through the rescaled spectral parameter
\begin{equation} \label{d92}
\sigma = \frac{z - E_{0}}{h} ,
\end{equation}
which will belong to some compact subset of $\C$. Note that all the coefficients of the matrix $\widetilde{\CQ}_{k}$ are zero except those on the $k$-th line.

\begin{remark}\sl \label{d3}
$i)$ The quantities $\CM_{\bullet}^{\pm}$ and $g_{\pm}^{\bullet}$ may depend on the choice of the parametrization of the curves $\gamma_{\bullet}$, but we prove in Section \ref{s23} that the eigenvalues (as well as their multiplicities) of the matrix $\CQ$ are independent of this choice.

$ii)$ The matrix $\CQ$ is not defined when the scalar product between incoming and outgoing directions vanishes. This justifies \ref{h4} and shows that this hypothesis is not technical. For some related results (concerning the scattering amplitude) where \ref{h4} is not supposed to hold, we send the reader to \cite{ALBoRa08_01}.

$iii)$ We could also have written the matrix $\CQ$ as the operators $\CT_{0} , \CT$ in Section \ref{s77} and Section \ref{s26}, using the time delay $T$ and geometric quantities which measure the amplification (of a WKB solution for example) along the curves $\gamma_{\bullet}$. We did not in order to avoid the introduction of some definitions which are useless for a finite number of curves.
\end{remark}

We have defined all the quantities we need to give the approximate quantization rule for resonances.

\begin{definition}[Quantization rule]\sl \label{d1}
We say that $z$ is a pseudo-resonance when
\begin{equation*}
1 \in \spe \big( h^{S ( z , h ) / \lambda_{1} - 1 / 2} \CQ ( z , h ) \big) .
\end{equation*}
The set of pseudo-resonances is denoted by $\res_{0} ( P )$.
\end{definition}

Let us first describe the asymptotic of the pseudo-resonances. From the previous definition, the spectrum of $\CQ ( z , h )$ plays a central role in this question. Thus, let $\mu_{1} ( \tau , h ) , \ldots , \mu_{K} ( \tau , h )$ denote, for $\tau \in \R$, the eigenvalues of
\begin{equation} \label{e18}
\widehat{\CQ} ( \tau , h ) : = \CQ \Big( E_{0} + h \tau - i h \sum_{j = 2}^{n} \frac{\lambda_{j}}{2} , h \Big) = \widetilde{\CQ} \Big( \rho , \tau - i \sum_{j = 2}^{n} \frac{\lambda_{j}}{2} \Big) .
\end{equation}
Since $[- C , C ] \times ] 0 , 1 ] \ni ( \tau , h ) \mapsto \widehat{\CQ} ( \tau , h )$ is analytic, the perturbation theory in finite dimensional spaces (see Chapter II of Kato \cite{Ka76_01}) ensures that, with an appropriate labeling of the eigenvalues, $( \tau , h ) \mapsto \mu_{k} ( \tau , h )$ are continuous and that, for $h$ fixed, $\tau \mapsto \mu_{k} ( \tau , h )$ are analytic functions with only algebraic singularities at some exceptional points. The pseudo-resonances satisfy the following property.

\begin{proposition}[Asymptotic of the pseudo-resonances]\sl \label{d9}
Assume \ref{h1}--\ref{h4}, \ref{h8}. Let $C > 0$ and let $\delta (h)$ be a function which goes to $0$ as $h \to 0$. Then, uniformly for $\tau \in [ - C , C ]$, the pseudo-resonances $z$ in
\begin{equation} \label{d12}
E_{0} + [ - C h , C  h ] + i \Big[ - \sum_{j = 2}^{n} \frac{\lambda_{j}}{2} h - C \frac{h}{\vert \ln h \vert} , h \Big] ,
\end{equation}
with $\re z \in E_{0} + \tau h + h \delta (h) [ - 1 , 1 ]$ satisfy
\begin{equation} \label{d94}
z = z_{q , k} ( \tau ) + o \Big( \frac{h}{\vert \ln h \vert} \Big) ,
\end{equation}
with
\begin{equation} \label{d95}
z_{q , k} ( \tau ) = E_{0} + 2 q \pi \lambda_{1} \frac{h}{\vert \ln h \vert} - i h \sum_{j = 2}^{n} \frac{\lambda_{j}}{2} + i \ln ( \mu_{k} ( \tau , h ) ) \lambda_{1} \frac{h}{\vert \ln h \vert} ,
\end{equation}
for some $q \in \Z$ and $k \in \{ 1 , \ldots , K\}$. On the other hand, for each $\tau \in [ - C , C ]$, $q \in \Z$ and $k \in \{ 1 , \ldots , K\}$ such that $z_{q , k} ( \tau )$ belongs to \eqref{d12} with a real part lying in $E_{0} + \tau h + h \delta (h) [ - 1 , 1 ]$, there exists a pseudo-resonance $z$ satisfying \eqref{d94} uniformly with respect to $q , k , \tau$.
\end{proposition}

We do not need to specify the determination of the logarithm of $\mu_{k} ( \tau , h )$ in \eqref{d95}, since a change of determination is balanced by a change of $q \in \Z$. Note also that $z_{q , k} ( \tau )$ is in \eqref{d12} only for eigenvalues $\mu_{k} ( \tau , h )$ outside a vicinity of $0$.

To compare the set of resonances and the set of pseudo-resonances, we will use the following definition. Its flexibility avoids the problems which may occur at the boundary of the domain of study.

\begin{definition}\sl \label{g80}
Let $A , B , C$ be subsets of $\C$ and $\varepsilon \geq 0$. We say that
\begin{equation*}
\dist ( A , B ) \leq \varepsilon \text{ in } C ,
\end{equation*}
if and only if
\begin{align*}
&\forall a \in A \cap C , \quad \exists b \in B , \qquad \vert a - b \vert \leq \varepsilon ,   \\
\text{and} \quad &\forall b \in B \cap C , \quad \exists a \in A , \qquad \vert a - b \vert \leq \varepsilon .
\end{align*}
\end{definition}

\begin{figure}
\begin{center}
\begin{picture}(0,0)%
\includegraphics{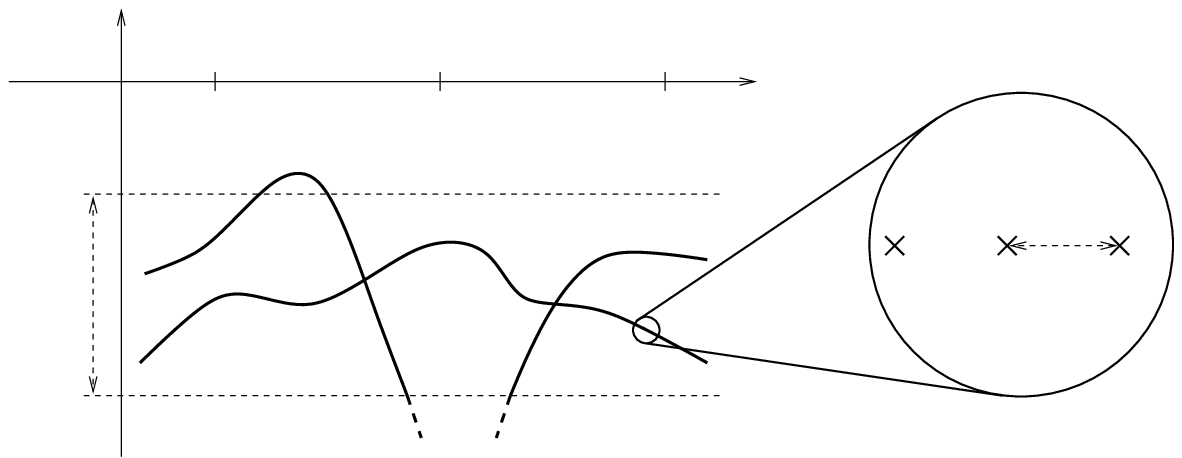}%
\end{picture}%
\setlength{\unitlength}{1184sp}%
\begingroup\makeatletter\ifx\SetFigFont\undefined%
\gdef\SetFigFont#1#2#3#4#5{%
  \reset@font\fontsize{#1}{#2pt}%
  \fontfamily{#3}\fontseries{#4}\fontshape{#5}%
  \selectfont}%
\fi\endgroup%
\begin{picture}(20166,7244)(-2114,-6983)
\put(-2099,-2836){\makebox(0,0)[lb]{\smash{{\SetFigFont{9}{10.8}{\rmdefault}{\mddefault}{\updefault}$\ds - \frac{1}{2} \sum_{j = 2}^{n} \lambda_{j} h$}}}}
\put(6301,-586){\makebox(0,0)[b]{\smash{{\SetFigFont{9}{10.8}{\rmdefault}{\mddefault}{\updefault}$E_{0}$}}}}
\put(2701,-586){\makebox(0,0)[b]{\smash{{\SetFigFont{9}{10.8}{\rmdefault}{\mddefault}{\updefault}$E_{0} - C h$}}}}
\put(9901,-586){\makebox(0,0)[b]{\smash{{\SetFigFont{9}{10.8}{\rmdefault}{\mddefault}{\updefault}$E_{0} + C h$}}}}
\put(16351,-2986){\makebox(0,0)[b]{\smash{{\SetFigFont{9}{10.8}{\rmdefault}{\mddefault}{\updefault}$2 \pi \lambda_{1} \frac{h}{\vert \ln h \vert}$}}}}
\put(-1349,-4486){\makebox(0,0)[lb]{\smash{{\SetFigFont{9}{10.8}{\rmdefault}{\mddefault}{\updefault}$C \frac{h}{\vert \ln h \vert}$}}}}
\end{picture}%
\end{center}
\caption{Two scale asymptotic of resonances of Theorem \ref{d8}.} \label{f23}
\end{figure}

For a finite number of homoclinic trajectories, our main result is the following.

\begin{theorem}[Asymptotic of resonances]\sl \label{d8}
Assume \ref{h1}--\ref{h4}, \ref{h8} and let $C , \delta > 0$. In the domain
\begin{equation} \label{d90}
E_{0} + [ - C h , C  h ] + i \Big[ - \sum_{j = 2}^{n} \frac{\lambda_{j}}{2} h - C \frac{h}{\vert \ln h \vert} , h \Big] \setminus \big( \Gamma (h) + B ( 0 , \delta h ) \big),
\end{equation}
we have
\begin{equation*}
\dist \big( \res (P) , \res_{0} (P) \big) = o \Big( \frac {h}{\vert \ln h \vert} \Big) ,
\end{equation*}
as $h$ goes to $0$. Moreover, for all $\chi \in C^{\infty}_{0} ( \R^{n} )$, there exists $M > 0$ such that
\begin{equation*}
\big\Vert \chi ( P -z )^{-1} \chi \big\Vert \lesssim h^{- M} ,
\end{equation*}
uniformly for $h$ small enough and $z \in \eqref{d90}$ with $\dist ( z , \res_{0} ( P ) ) \geq \delta h \vert \ln h \vert^{- 1}$.
\end{theorem}

The result above describes the semiclassical behavior of the resonances as a set. It shows that the resonances are close to the pseudo-resonances and that there exists at least one resonance near any pseudo-resonance. But it does not give the number of resonances that are close to each pseudo-resonance, neither it specifies the multiplicity of these resonances. This is not necessarily due to our general strategy (see Section \ref{s36}) but is more likely a consequence of the fact that only homogeneous Cauchy problem can be treated in Theorem \ref{a32}. Nevertheless, it is possible to give a lower bound on the number of resonances close to a pseudo-resonance by its multiplicity. More precisely,

\begin{proposition}[Lower bound on the multiplicity]\sl \label{d91}
In the setting of Theorem \ref{d8}, let $s > 0$ be small enough. Then, for any pseudo-resonance $z$ in \eqref{d90}, we have
\begin{equation*}
\card \Big( \res (P) \cap B \Big( z , 2 s \frac{h}{\vert \ln h \vert} \Big) \Big) \geq \card \Big\{ ( q , k ) \in \Z \times \{ 1 , \ldots , K \} ; \ z_{q , k} ( \re \sigma ) \in B \Big( z , s \frac{h}{\vert \ln h \vert} \Big) \Big\} ,
\end{equation*}
where $\sigma$ and $z_{q , k} ( \cdot )$ are given by \eqref{d92} and \eqref{d95} respectively. In the previous expression, the resonances are counted with their multiplicity.
\end{proposition}

\Subsection{Remarkable phenomena} \label{s18}

From the previous results, we deduce four important properties of the resonances generated by homoclinic trajectories: they accumulate on curves, these curves vibrate with $h$, their large behavior illustrates the transition of the trapping at the critical energy and they are stable by small perturbations of the set of homoclinic trajectories.

\subsubsection{Accumulation on curves} \label{s17}

First, we show that the resonances generated by homoclinic trajectories concentrate on curves. More precisely, Proposition \ref{d9} and Theorem \ref{d8} imply

\begin{remark}[Two scale asymptotic]\sl \label{e2}
The resonances satisfy a two scale asymptotic in this homoclinic setting. This phenomenon is illustrated in Figure \ref{f23}.

At the macroscopic scale $h$, the resonances accumulate on the curves
\begin{equation} \label{d93}
\im \sigma = - \sum_{j = 2}^{n} \frac{\lambda_{j}}{2} + \ln \big( \vert \mu_{k} ( \re \sigma , h ) \vert \big)  \frac{\lambda_{1}}{\vert \ln h \vert} ,
\end{equation}
where $k \in \{ 1 , \ldots , K \}$ and
\begin{equation} \label{i70}
\sigma = \frac{z - E_{0}}{h} ,
\end{equation}
with $\re \sigma \in [ - C , C ]$, is the associated rescaled spectral parameter. For $h$ fixed, these curves are continuous and analytic outside of their crossings. Note that this scale is the typical gap between two resonances generated by a periodic trajectory of hyperbolic type (see G{\'e}rard and Sj{\"o}strand \cite{GeSj87_01}) and corresponds to an interval of size $1$ in the non-semiclassical setting (see Ikawa \cite{Ik83_01}).

However, at the microscopic scale $h \vert \ln h \vert^{- 1}$ (at a matter of fact, any $o ( h )$ scale would give the same picture), these accumulation curves become horizontal lines on which the packets of resonances (see the discussion above Proposition \ref{d91}) are regularly distributed and spaced out by
\begin{equation*}
2 \pi \lambda_{1} \frac{h}{\vert \ln h \vert} + o \Big( \frac{h}{\vert \ln h \vert} \Big) .
\end{equation*}
Nevertheless, some of this lines sometimes coincide and the distance between two consecutive packets of resonances can be different (see Example \ref{e6}).
\end{remark}

The results of the previous section imply that there are at most $K$ accumulation curves in the set \eqref{d12}. In general, the number of these curves can be smaller since some of the $\mu_{k}$ may  vanish (see Section \ref{s16}). It can even happen that the set of pseudo-resonances in \eqref{d12} is empty (see Example \ref{e48} {\rm (B)}).

Note that two accumulation curves cross if and only if the modulus (and not the value) of the corresponding eigenvalues of $\CQ$ coincide. Thus, these curves tend to be more regular (i.e. analytic with respect to $\tau , h$) at their crossings than one may think.

We now apply our results to the simplest possible case: when $\CH$ consists of a single trajectory. More precisely, we replace \ref{h8} in this discussion by
\begin{hyp} \label{h9}
The homoclinic set $\CH$ consists of a unique trajectory $\gamma$ on which $\Lambda_{-}$ and $\Lambda_{+}$ intersect transversally.
\end{hyp}
Example \ref{c15} provides operators satisfying \ref{h1}--\ref{h4} and \ref{h9}. Under the assumption \ref{h9}, we remove the subscript $k = 1$ which was used to indicate the number of the trajectory. Note that the matrix $\CQ$ is just a scalar in this case. Moreover, for Schr\"{o}dinger operators, the usual symmetry of the Hamiltonian trajectories $( x (t) , \xi (t) ) \longmapsto ( x ( - t ) ,  - \xi ( - t ) )$ and the uniqueness of the homoclinic trajectory imply that $g_{-}$ and $g_{+}$ are the same up to some positive constant. From Proposition \ref{d9}, we then deduce that
\begin{equation} \label{e8}
z_{q} ( \tau ) =  E_{0} - \frac{A \lambda_{1}}{\vert \ln h \vert} + 2 q \pi \lambda_{1} \frac{h}{\vert \ln h \vert} - i h \sum_{j = 2}^{n} \frac{\lambda_{j}}{2} + i \ln ( \mu ( \tau ) ) \lambda_{1} \frac{h}{\vert \ln h \vert} ,
\end{equation}
with $q \in \Z$ and
\begin{equation} \label{g4}
\mu ( \tau ) = \Gamma \Big( \frac{1}{2} - i \frac{\tau}{\lambda_{1}} \Big) \sqrt{\frac{\lambda_{1}}{2 \pi}} \frac{\CM^{+}}{\CM^{-}} e^{- \frac{\pi}{2} ( \nu + 1 ) i} \vert g_{-} \vert \big( \lambda_{1} \vert g_{+} \vert \vert g_{-} \vert \big)^{- \frac{1}{2} + i \frac{\tau}{\lambda_{1}} } e^{- \frac{\pi \tau}{2 \lambda_{1}}} .
\end{equation}
Combining with Theorem \ref{d8}, we deduce the localization of resonances under \ref{h9}. In particular, there is a unique accumulation curve in the sense of Remark \ref{e2}. Using the formula
\begin{equation} \label{e36}
\forall y \in \R , \qquad \Big\vert \Gamma \Big( \frac{1}{2} + i y \Big) \Big\vert^{2} = \frac{\pi}{\cosh ( \pi y )} ,
\end{equation}
(see \cite[(6.1.30)]{AbSt64_01}), the asymptotic behavior of this curve for large $\re \sigma$ is as follows.

\begin{figure}
\begin{center}
\begin{picture}(0,0)%
\includegraphics{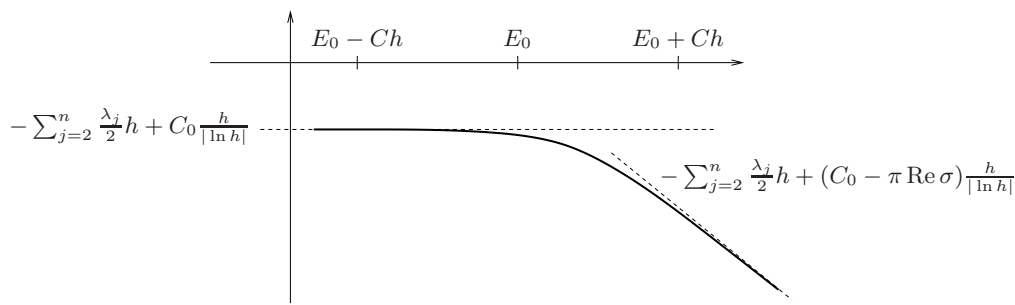}%
\end{picture}%
\setlength{\unitlength}{1105sp}%
\begingroup\makeatletter\ifx\SetFigFont\undefined%
\gdef\SetFigFont#1#2#3#4#5{%
  \reset@font\fontsize{#1}{#2pt}%
  \fontfamily{#3}\fontseries{#4}\fontshape{#5}%
  \selectfont}%
\fi\endgroup%
\begin{picture}(15930,6644)(-2414,-6383)
\put(13501,-3661){\makebox(0,0)[b]{\smash{{\SetFigFont{9}{10.8}{\rmdefault}{\mddefault}{\updefault}$- \sum_{j = 2}^{n} \frac{\lambda_{j}}{2} h + ( C_{0} - \pi \re \sigma ) \frac{h}{\vert \ln h \vert}$}}}}
\put(6301,-586){\makebox(0,0)[b]{\smash{{\SetFigFont{9}{10.8}{\rmdefault}{\mddefault}{\updefault}$E_{0}$}}}}
\put(2701,-586){\makebox(0,0)[b]{\smash{{\SetFigFont{9}{10.8}{\rmdefault}{\mddefault}{\updefault}$E_{0} - C h$}}}}
\put(9901,-586){\makebox(0,0)[b]{\smash{{\SetFigFont{9}{10.8}{\rmdefault}{\mddefault}{\updefault}$E_{0} + C h$}}}}
\put(-2399,-2536){\makebox(0,0)[b]{\smash{{\SetFigFont{9}{10.8}{\rmdefault}{\mddefault}{\updefault}$- \sum_{j = 2}^{n} \frac{\lambda_{j}}{2} h + C_{0} \frac{h}{\vert \ln h \vert}$}}}}
\end{picture}%
\end{center}
\caption{The accumulation curve generated by one homoclinic trajectory.} \label{f24}
\end{figure}

\begin{corollary}\sl \label{e4}
Under \ref{h1}--\ref{h4} and \ref{h9}, the accumulation curve given by \eqref{d93} verifies
\begin{align}
\im \sigma &= - \sum_{j = 2}^{n} \frac{\lambda_{j}}{2} + \frac{C_{0}}{\vert \ln h \vert} - \frac{\lambda_{1}}{2 \vert \ln h \vert} \ln \big( e^{\frac{2 \pi}{\lambda_{1}} \re \sigma} + 1 \big)    \nonumber \\
&= - \sum_{j = 2}^{n} \frac{\lambda_{j}}{2} + \frac{1}{\vert \ln h \vert}
\left\{ \begin{aligned}
&C_{0} + \CO \big( e^{- \frac{2 \pi}{\lambda_{1}} \vert \re \sigma \vert} \big) &&\text{ as } \re \sigma \to - \infty , \\
&- \pi \re \sigma + C_{0} + \CO \big( e^{- \frac{2 \pi}{\lambda_{1}} \vert \re \sigma \vert} \big) &&\text{ as } \re \sigma \to + \infty , \\
\end{aligned} \right.  \label{e3}
\end{align}
with the constant
\begin{equation} \label{e7}
C_{0} = \lambda_{1} \ln \bigg( \frac{\CM^{+}}{\CM^{-}} \sqrt{\frac{\vert g_{-} \vert}{\vert g_{+} \vert}} \bigg) .
\end{equation}
\end{corollary}

This result is summarized in Figure \ref{f24}. From a formal point of view, the asymptotic of the resonances obtained here can be understood as the combination of \eqref{n1} and \eqref{e37} at least for potentials of product type as in Example \ref{c15}. Indeed, the operator has a homoclinic orbit on the $x_{1}$-axis as in Section \ref{s19} {\rm (B)}, whereas it has a barrier-top in the other directions. The difference is that $C_{0}$ vanishes in dimension one and not necessarily in higher dimensions.

From \eqref{e3}, the unique accumulation curve is independent of $h$. Thus, the asymptotic of the resonances is frozen at the macroscopic scale. The situation is different at the microscopic scale. For $\tau$ fixed, \eqref{e8} shows that the $z_{q} ( \tau )$'s with real part close to $E_{0} + \tau h$ verify
\begin{equation*}
\frac{\vert \ln h \vert}{h} \Big( z_{q} ( \tau ) - \Big( E_{0} + \tau h  - i h \sum_{j = 2}^{n} \frac{\lambda_{j}}{2} \Big) \Big) = - A \lambda_{1} h^{- 1} + 2 q \pi \lambda_{1} - \tau \vert \ln h \vert + i \ln ( \mu ( \tau ) ) \lambda_{1} .
\end{equation*}
Hence, the resonances seem to ``march past'' from the right to the left as $h$ goes to $0$. This phenomenon is an example of microscopic dynamic on the accumulations curves.

If we relax the assumption that $P$ is a Schr\"{o}dinger operator, it is possible to construct pseudodifferential operators satisfying assumptions similar (as explained in Remark \ref{c13}) to \ref{h1}--\ref{h4} and \ref{h9} and for which $g_{-} \cdot g_{+} < 0$. In this case, $\re \sigma$ is replaced by $- \re \sigma$ in Corollary \ref{e4}. This geometric setting can not occur for a Schr\"{o}dinger operator $- h^{2} \Delta + V (x)$ since $g_{-} / \vert g_{-} \vert = g_{+} / \vert g_{+} \vert$ in that case. This last property follows from the usual symmetry of the Hamiltonian curves $( x (t) , \xi (t) ) \longmapsto ( x ( - t ) ,  - \xi ( - t ) )$ which holds for the Schr\"{o}dinger operators (and more generally if the symbol $p ( x , \xi )$ of $P$ is even with respect to $\xi$). We now give an example of such an operator.

\begin{figure}
\begin{center}
\begin{picture}(0,0)%
\includegraphics{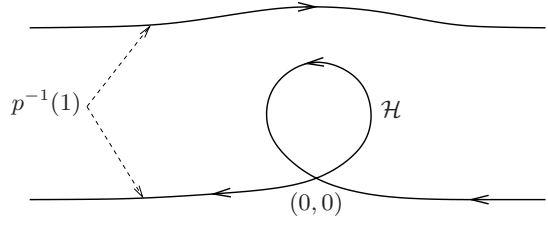}%
\end{picture}%
\setlength{\unitlength}{1184sp}%
\begingroup\makeatletter\ifx\SetFigFont\undefined%
\gdef\SetFigFont#1#2#3#4#5{%
  \reset@font\fontsize{#1}{#2pt}%
  \fontfamily{#3}\fontseries{#4}\fontshape{#5}%
  \selectfont}%
\fi\endgroup%
\begin{picture}(11223,4425)(811,-5746)
\put(7201,-5686){\makebox(0,0)[b]{\smash{{\SetFigFont{9}{10.8}{\rmdefault}{\mddefault}{\updefault}$( 0 , 0 )$}}}}
\put(8551,-3736){\makebox(0,0)[lb]{\smash{{\SetFigFont{9}{10.8}{\rmdefault}{\mddefault}{\updefault}$\CH$}}}}
\put(826,-3586){\makebox(0,0)[lb]{\smash{{\SetFigFont{9}{10.8}{\rmdefault}{\mddefault}{\updefault}$p^{-1} ( 1 )$}}}}
\end{picture}%
\end{center}
\caption{The energy level $p = 1$ in Example \ref{b81}.} \label{f14}
\end{figure}

\begin{example}\rm \label{b81}
In dimension $n = 1$, we set
\begin{equation*}
P = - h^{2} \Delta + \chi (x) \Op ( V ( x , \xi ) ) \chi (x) ,
\end{equation*}
where
\begin{equation*}
V ( x , \xi ) = e^{- x^{2}}  \big( 1 - \xi^{2} ( \xi + 1 / 2 ) \big) e^{- \xi^{4} / 5} \in S (1) ,
\end{equation*}
and $\chi \in C^{\infty}_{0} ( \R ; [ 0 , 1 ] )$ satisfies $\chi = 1$ on a sufficiently large neighborhood of $0$. The pseudodifferential operator $P$ does not satisfy the assumptions of the paper (it is not even of the form \eqref{a5}). Nevertheless, the results can be adapted to the present operator. In particular, since the cut-off function $\chi$ guaranties that the perturbation is compactly supported, the resonances of $P$ can be defined using the abstract theory of resonances for operators of ``black box'' type developed by Sj\"{o}strand and Zworski \cite{SjZw91_01}.

From the semiclassical pseudodifferential calculus, we have
\begin{equation*}
P = \Op ( p ( x, \xi ) ) + \Psi ( h^{2} ) ,
\end{equation*}
where $p ( x , \xi ) = \xi^{2} (x) + \chi^{2} V ( x , \xi )$. The energy level $p = 1$ is described in Figure \ref{f14}. Moreover, near $( 0 , 0 )$, we have
\begin{equation*}
p ( x , \xi ) = 1 + \xi^{2} / 2 - x^{2} + \CO \big( (x , \xi )^{3} \big) .
\end{equation*}
Thus, $( 0 , 0 )$ is a hyperbolic fixed point. One can see in Figure \ref{f14} that the trapped set at energy $E_{0} = 1$ satisfies \ref{h3}. Furthermore, \ref{h4} and \ref{h9} are automatically satisfied in dimension $n = 1$, as explained in Section \ref{s19} {\rm (B)}. Nevertheless, $g_{\pm} \in \R \setminus \{ 0 \}$ have opposite sign and so $g_{-} \cdot g_{+} < 0$. Then, the accumulation curve of resonances is described by Corollary \ref{e4} with $\re \sigma$ replaced by $- \re \sigma$. In this case, one can show that $C_{0} = 0$ which illustrates the conservation of energy along $\CH$ (as always for self-adjoint situations in dimension $n = 1$).

Eventually, we discuss the resonance free regions of Section \ref{s22} for this example. Using the notation $\S^{0} = \{ - 1 , 1 \}$, we have $\CH_{\rm tang}^{- \infty} = \{ - 1 \}$, $\CH_{\rm tang}^{+ \infty} = \{ + 1 \}$ and
\begin{equation} \label{b82}
\CA_{0} ( \tau ) = \frac{e^{\frac{\pi}{2} \frac{\tau}{\lambda_{1}}}}{\sqrt{2 \cosh ( \pi \tau / \lambda_{1} )}} ,
\end{equation}
which can be compared with \eqref{b75}. In particular, \ref{h6} holds true but not \ref{h5}, as always for $n = 1$. Then, as in Section \ref{s19} {\rm (B)}, Theorem \ref{a2} (but not Theorem \ref{a1}) can be apply to give the resonance free zone
\begin{equation}
\im \sigma \geq - \frac{\lambda_{1}}{2} \ln \big( e^{- \frac{2 \pi}{\lambda_{1}} \re \sigma} + 1 \big) \frac{1}{\vert \ln h \vert} + \delta \frac{1}{\vert \ln h \vert} ,
\end{equation}
below any interval $\re \sigma \in [ - C , C]$ with $C > 0$. Comparing with Corollary \ref{e4}, we notice that this example provides another situation where the resonance free zone of Theorem \ref{a2} is optimal.
\end{example}

\subsubsection{Vibration phenomena}

Coming back to the general assumptions of Theorem \ref{d8}, Proposition \ref{d9} also shows that the accumulation curves of Remark \ref{e2} oscillate as $h$ goes to $0$. This follows from the quasi-periodicity of $h^{- 1} \longmapsto \widehat{\CQ} ( \tau , h )$. More precisely,

\begin{remark}[Vibration of the accumulation curves]\sl \label{e5}
The second main property verified by the resonances is that the curves on which they accumulate (see Remark \ref{e2}) depend on $h$. Here, we do not refer to the factor $\vert \ln h \vert^{- 1}$ in \eqref{d93}, but to the fact that even the shape of these curves changes with $h$. Indeed the $\mu_{k} ( \tau , h )$'s, the eigenvalues of $\widehat{\CQ} ( \tau , h )$, may depend on $h$. Formula \eqref{d10} shows that this dependence only comes from the actions $\rho_{k} = e^{i A_{k} / h}$. From usual perturbation theory (see Chapter II.5.7 of \cite{Ka76_01} for instance), the accumulation curves are continuous with respect to $h$ and analytic outside of the crossings but can have rather complicated singularities.

If all the actions are equal, one can factor $\widehat{\CQ} ( \tau , h )$ by $e^{i A_{1} / h} = \cdots = e^{i A_{K} / h}$, so that $\mu_{k} ( \tau , h ) = e^{i A_{1} / h} e^{- i A_{1}} \mu_{k} ( \tau , 1 )$. Since only the modulus of the $\mu_{k} ( \tau , h )$'s appears in \eqref{d93}, the accumulation curves do not depend on $h$ in this case. Of course, this assumption is satisfied when $\CH$ consists of an a single trajectory. In this particular setting, one can verify from Corollary \ref{e4} that the unique accumulation curve does not move with $h$.

Assume now that the actions take exactly two different values, denoted $A_{1}$ and $A_{2}$. As before, we can factor out $\widehat{\CQ} ( \tau , h )$ by $e^{i A_{1} / h}$. Then the relevant matrix, as well as the $\vert \mu_{k} ( \tau , h ) \vert$, depend on $h$ through the quantity $e^{i ( A_{2} - A_{1} ) / h}$. Thus, the accumulation curves are periodic as functions of $h^{- 1}$ with period $2 \pi \vert A_{2} - A_{1} \vert^{- 1}$. Example \ref{e6} below provides an operator which satisfies this condition.

Eventually, if there exist at least three actions $A_{k}$ which are $\Z$-independent, the accumulation curves are no longer periodic with respect to $h^{- 1}$ but are smooth (at least continuous) functions of $e^{i ( A_{2} - A_{1} ) / h} , \ldots , e^{i ( A_{K} - A_{1} ) / h}$. This is the vibration phenomenon.
\end{remark}

We now construct an example of operator for which the vibration phenomenon described above occurs. It is a Schr\"{o}dinger operator with two homoclinic trajectories with opposite asymptotic directions.

\begin{example}\rm \label{e6}
In dimension $n = 2$, we consider
\begin{equation*}
V (x) = V_{0} ( x ) + V_{1} ( x - x_{1} ) + V_{2} ( x - x_{2} ) ,
\end{equation*}
where $V_{0}$ is as in Section \ref{s21}. For $k = 1 , 2$, the potential $V_{k} \in C_{0}^{\infty} ( \R^{2} )$ is a radial function satisfying $x \cdot \nabla V_{k} (x) < 0$ for $x$ in the interior of $\supp V_{k} \setminus \{ 0 \}$ and $E_{0} < V_{k} ( 0 )$. One can take $V_{1} = V_{2}$. Eventually, the vectors $x_{k} \in \R^{2}$ are chosen in such way that their norm is sufficiently large and such that the angle $( x_{1} , x_{2} )$ is close enough to $\pi$. The geometric setting as well as the asymptotic of resonances are illustrated in Figure \ref{f25}.

\begin{figure}
\begin{center}
\begin{picture}(0,0)%
\includegraphics{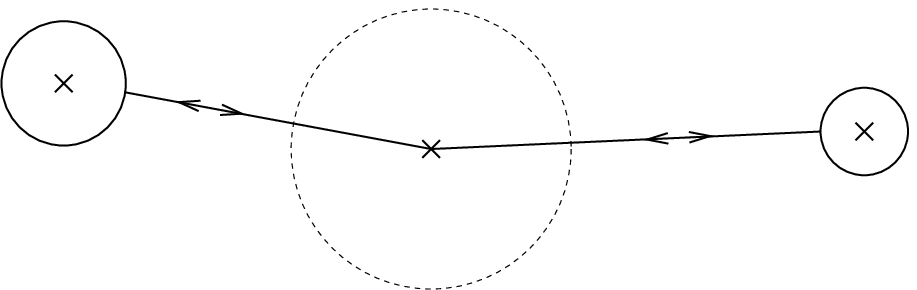}%
\end{picture}%
\setlength{\unitlength}{1105sp}%
\begingroup\makeatletter\ifx\SetFigFont\undefined%
\gdef\SetFigFont#1#2#3#4#5{%
  \reset@font\fontsize{#1}{#2pt}%
  \fontfamily{#3}\fontseries{#4}\fontshape{#5}%
  \selectfont}%
\fi\endgroup%
\begin{picture}(15587,5336)(-15788,-7575)
\put(-12224,-5161){\makebox(0,0)[b]{\smash{{\SetFigFont{9}{10.8}{\rmdefault}{\mddefault}{\updefault}$\pi_{x} ( \gamma_{2} )$}}}}
\put(-8399,-4786){\makebox(0,0)[b]{\smash{{\SetFigFont{9}{10.8}{\rmdefault}{\mddefault}{\updefault}$0$}}}}
\put(-14699,-4636){\makebox(0,0)[b]{\smash{{\SetFigFont{9}{10.8}{\rmdefault}{\mddefault}{\updefault}$x_{2}$}}}}
\put(-899,-5386){\makebox(0,0)[b]{\smash{{\SetFigFont{9}{10.8}{\rmdefault}{\mddefault}{\updefault}$x_{1}$}}}}
\put(-899,-3697){\makebox(0,0)[b]{\smash{{\SetFigFont{9}{10.8}{\rmdefault}{\mddefault}{\updefault}$\{ V_{1} ( x - x_{1} ) = E_{0} \}$}}}}
\put(-14699,-2611){\makebox(0,0)[b]{\smash{{\SetFigFont{9}{10.8}{\rmdefault}{\mddefault}{\updefault}$\{ V_{2} ( x - x_{2} ) = E_{0} \}$}}}}
\put(-8399,-2386){\makebox(0,0)[b]{\smash{{\SetFigFont{9}{10.8}{\rmdefault}{\mddefault}{\updefault}$\supp V_{0}$}}}}
\put(-4124,-5686){\makebox(0,0)[b]{\smash{{\SetFigFont{9}{10.8}{\rmdefault}{\mddefault}{\updefault}$\pi_{x} ( \gamma_{1} )$}}}}
\end{picture}%
\end{center}
\bigskip \bigskip
\begin{center}
\begin{picture}(0,0)%
\includegraphics{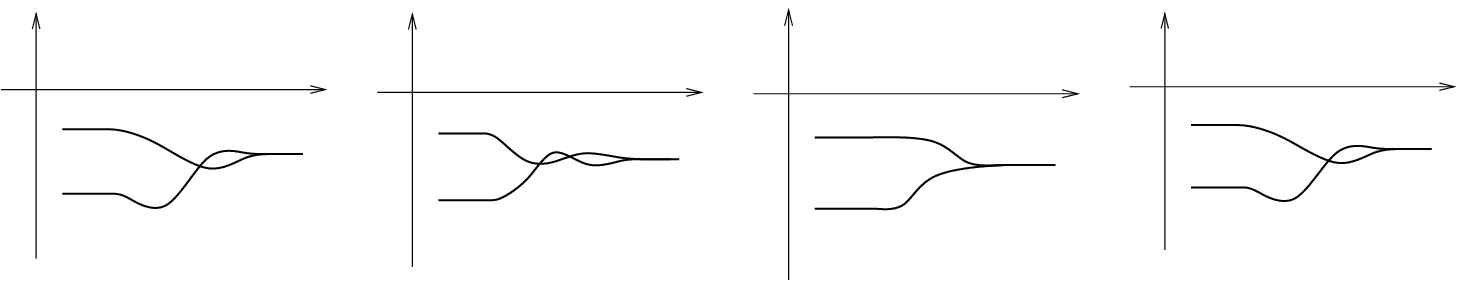}%
\end{picture}%
\setlength{\unitlength}{1105sp}%
\begingroup\makeatletter\ifx\SetFigFont\undefined%
\gdef\SetFigFont#1#2#3#4#5{%
  \reset@font\fontsize{#1}{#2pt}%
  \fontfamily{#3}\fontseries{#4}\fontshape{#5}%
  \selectfont}%
\fi\endgroup%
\begin{picture}(25019,4735)(-621,-5091)
\put(22051,-936){\makebox(0,0)[b]{\smash{{\SetFigFont{9}{10.8}{\rmdefault}{\mddefault}{\updefault}$\frac{A_{2} - A_{1}}{h} \equiv 2 \pi \text{ mod } 2 \pi$}}}}
\put(9151,-977){\makebox(0,0)[b]{\smash{{\SetFigFont{9}{10.8}{\rmdefault}{\mddefault}{\updefault}$\frac{A_{2} - A_{1}}{h} \equiv \frac{2}{3} \pi \text{ mod } 2 \pi$}}}}
\put(2701,-957){\makebox(0,0)[b]{\smash{{\SetFigFont{9}{10.8}{\rmdefault}{\mddefault}{\updefault}$\frac{A_{2} - A_{1}}{h} \equiv 0 \text{ mod } 2 \pi$}}}}
\put(15601,-941){\makebox(0,0)[b]{\smash{{\SetFigFont{9}{10.8}{\rmdefault}{\mddefault}{\updefault}$\frac{A_{2} - A_{1}}{h} \equiv \frac{4}{3} \pi \text{ mod } 2 \pi$}}}}
\end{picture}%
\end{center}
\caption{The geometry of the potential $V$ and the vibration of the accumulation curves in Example \ref{e6}.} \label{f25}
\end{figure}

For this operator, the assumptions \ref{h1}--\ref{h3} are satisfied, $\lambda_{1} = \lambda_{2} = : \lambda$ and the homoclinic set $\CH$ is the union of two radial curves $\gamma_{k}$ whose base space projection is included in $\R^{+} x_{k}$ for $k = 1 , 2$. To determine the trapped set, we remark that no Hamiltonian curve coming from one of the ``barriers'' (say the support of $V_{1} ( x - x_{1} )$) can touch the other ``barrier'' (the support of $V_{2} ( x - x_{2} )$). Since $g_{\pm}^{k}$ is collinear to $x_{k}$ and $x_{1} \cdot x_{2} \neq 0$, \ref{h4} holds true. Moreover, one can verify \ref{h8} as in Example \ref{c15}.

Then, we can apply Theorem \ref{d8} to obtain the asymptotic of the resonances. More precisely, the operator $\CQ$ is a $2 \times 2$ matrix in this case. We now use that the two eigenvalues of a general $2 \times 2$ matrix
\begin{equation} \label{i44}
\left( \begin{array}{cc}
a & b \\
c & d
\end{array} \right) ,
\end{equation}
are given by
\begin{equation} \label{e21}
\frac{a + d \pm \sqrt{a^{2} - 2 a d + d^{2} + 4 b c}}{2} .
\end{equation}
Combining with the explicit form of $\CQ$ defined in \eqref{d4}, we deduce the explicit form of the leading term \eqref{d95} of the pseudo-resonances as well as the accumulation curves \eqref{d93}. In particular, we observe that these accumulation curves are periodic with respect to $h^{- 1}$. Furthermore, the $z_{q , k} ( \tau )$'s satisfy the asymptotic
\begin{equation} \label{e9}
z_{q , k} ( \tau ) =  E_{0} - \frac{A_{k} \lambda}{\vert \ln h \vert} - i \frac{\lambda}{2} h + \big( 2 q \pi \lambda + i \ln ( \mu_{k} ( \tau ) ) \lambda + \CO ( \tau^{-1} ) \big) \frac{h}{\vert \ln h \vert} ,
\end{equation}
as $\tau \to - \infty$ with
\begin{equation} \label{l11}
\mu_{k} ( \tau ) = \Gamma \Big( \frac{1}{2} - i \frac{\tau}{\lambda} \Big) \sqrt{\frac{\lambda}{2 \pi}} \frac{\CM^{+}_{k}}{\CM^{-}_{k}} e^{- \frac{\pi}{2} ( \nu_{k} + 1 ) i} \vert g_{-}^{k} \vert \big( \lambda \vert g_{+}^{k} \vert \vert g_{-}^{k} \vert \big)^{- \frac{1}{2} + i \frac{\tau}{\lambda} } e^{- \frac{\pi \tau}{2 \lambda}} .
\end{equation}
On the contrary, we have
\begin{equation} \label{e16}
\left\{ \begin{aligned}
z_{q , 1} ( \tau ) &= E_{0} - \frac{( A_{1} + A_{2} ) \lambda}{2 \vert \ln h \vert} - i \frac{\lambda}{2} h + \big( 2 q \pi \lambda + i \ln ( \mu_{\widehat{12}} ( \tau ) ) \lambda + \CO ( \tau^{-1} ) \big) \frac{h}{\vert \ln h \vert} ,   \\
z_{q , 2} ( \tau ) &= z_{q , 1} ( \tau ) + \big( \pi \lambda + \CO ( \tau^{-1} ) \big) \frac{h}{\vert \ln h \vert} ,
\end{aligned} \right.
\end{equation}
as $\tau \to + \infty$ with
\begin{equation*}
\mu_{\widehat{12}} ( \tau ) = \Gamma \Big( \frac{1}{2} - i \frac{\tau}{\lambda} \Big) \sqrt{\frac{\lambda}{2 \pi}} \sqrt{\frac{\CM^{+}_{1} \CM^{+}_{2}}{\CM^{-}_{1} \CM^{-}_{2}}} e^{- \frac{\pi}{4} ( \nu_{1} + \nu_{2} ) i} \sqrt{\vert g_{-}^{1} \vert \vert g_{-}^{2} \vert} \Big( \lambda \sqrt{\vert g_{-}^{1} \cdot g_{+}^{2} \vert \vert g_{-}^{2} \cdot g_{+}^{1} \vert} \Big)^{- \frac{1}{2} + i \frac{\tau}{\lambda} } e^{\frac{\pi \tau}{2 \lambda}} .
\end{equation*}
Since \eqref{e8} and \eqref{e9} have the same behavior, we deduce that, in the limit $\tau$ goes to $- \infty$, the distribution of resonances is similar to the superposition of the resonances generated by the two trajectories $\gamma_{k}$ separately. On the contrary, as $\tau$ goes to $+ \infty$, the resonances seem to be generated by a unique trajectory (as in Example \ref{b81}) which could be some recombination of the $\gamma_{k}$'s, except that the packets of resonances are now spaced by $\pi \lambda h \vert \ln h \vert^{- 1}$.

Moreover, mimicking Corollary \ref{e4}, we deduce from \eqref{e36}, \eqref{e9} and \eqref{e16} that the two accumulation curves (in the sense of Remark \ref{e2}) verify
\begin{equation} \label{e33}
\im \sigma = - \frac{\lambda}{2} + \frac{1}{\vert \ln h \vert}
\left\{ \begin{aligned}
&C_{k} + \CO \big( \vert \re \sigma \vert^{- 1} \big) &&\text{ as } \re \sigma \to - \infty , \\
&C_{\widehat{12}} + \CO \big( \vert \re \sigma \vert^{- 1} \big) \ &&\text{ as } \re \sigma \to + \infty , \\
\end{aligned} \right. 
\end{equation}
with the constants
\begin{equation} \label{e17}
C_{k} = \lambda \ln \bigg( \frac{\CM^{+}_{k}}{\CM^{-}_{k}} \sqrt{\frac{\vert g_{-}^{k} \vert}{\vert g_{+}^{k} \vert}} \bigg) \qquad \text{and} \qquad C_{\widehat{12}} = \frac{C_{1} + C_{2}}{2} + \frac{\lambda}{2} \big\vert \ln \vert \cos ( x_{1} , x_{2} ) \vert \big\vert .
\end{equation}
Note that we always have $C_{\widehat{12}} \geq ( C_{1} + C_{2} ) / 2$. It means that the (unique) right asymptote of the accumulation curves is above the middle of the two left asymptotes (see Figure \ref{f25}).

\begin{figure}
\begin{center}
\begin{picture}(0,0)%
\includegraphics{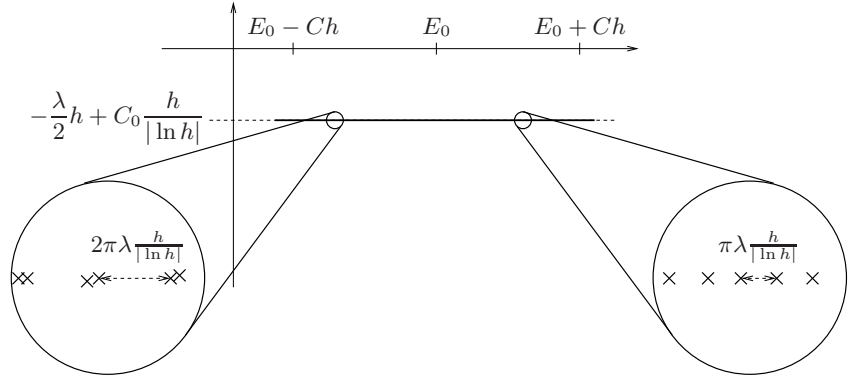}%
\end{picture}%
\setlength{\unitlength}{987sp}%
\begingroup\makeatletter\ifx\SetFigFont\undefined%
\gdef\SetFigFont#1#2#3#4#5{%
  \reset@font\fontsize{#1}{#2pt}%
  \fontfamily{#3}\fontseries{#4}\fontshape{#5}%
  \selectfont}%
\fi\endgroup%
\begin{picture}(21029,9433)(-4402,-9172)
\put(-1199,-6061){\makebox(0,0)[b]{\smash{{\SetFigFont{9}{10.8}{\rmdefault}{\mddefault}{\updefault}$2 \pi \lambda \frac{h}{\vert \ln h \vert}$}}}}
\put(6301,-586){\makebox(0,0)[b]{\smash{{\SetFigFont{9}{10.8}{\rmdefault}{\mddefault}{\updefault}$E_{0}$}}}}
\put(2701,-586){\makebox(0,0)[b]{\smash{{\SetFigFont{9}{10.8}{\rmdefault}{\mddefault}{\updefault}$E_{0} - C h$}}}}
\put(9901,-586){\makebox(0,0)[b]{\smash{{\SetFigFont{9}{10.8}{\rmdefault}{\mddefault}{\updefault}$E_{0} + C h$}}}}
\put(-1724,-2836){\makebox(0,0)[b]{\smash{{\SetFigFont{9}{10.8}{\rmdefault}{\mddefault}{\updefault}$\ds - \frac{\lambda}{2} h+ C_{0} \frac{h}{\vert \ln h \vert}$}}}}
\put(14401,-6061){\makebox(0,0)[b]{\smash{{\SetFigFont{9}{10.8}{\rmdefault}{\mddefault}{\updefault}$\pi \lambda \frac{h}{\vert \ln h \vert}$}}}}
\end{picture}%
\end{center}
\caption{The microscopic dynamic of the resonances in the symmetric case.} \label{f26}
\end{figure}

Moreover, if the two barriers are symmetric (i.e. $V_{1} = V_{2}$ and $x_{1} = - x_{2}$), a direct calculus shows that the two accumulation curves coincide and satisfy
\begin{equation}
\im \sigma = - \frac{\lambda}{2} + \frac{C_{0}}{\vert \ln h \vert} ,
\end{equation}
where $C_{0} : = C_{1} = C_{2} = C_{\widehat{12}}$ is given by \eqref{e17}. Thus, the asymptotic of the resonances is particularly basic at the macroscopic scale. On the other hand, the asymptotic of the $z_{q , k} ( \tau )$ stated above give in this symmetric case
\begin{equation*}
z_{q , 1} ( \tau ) = z_{q , 2} ( \tau ) + \frac{h}{\vert \ln h \vert}
\left\{ \begin{aligned}
&\CO ( \tau^{-1} )  &&\text{ as } \tau \to - \infty , \\
&\pi \lambda + \CO ( \tau^{-1} ) &&\text{ as } \tau \to + \infty . \\
\end{aligned} \right.
\end{equation*}
Then, for $\tau \ll - 1$, the packets of resonances (see Remark \ref{e2}) have multiplicity at least $2$ from Proposition \ref{d91} and are spaced out by $2 \pi \lambda h \vert \ln h \vert^{- 1}$. On the contrary, for $\tau \gg 1$, the packets of resonances have multiplicity at least $1$ and are spaced out by $\pi \lambda h \vert \ln h \vert^{- 1}$. Thus, the microscopic distribution of the resonances changes along the accumulation curves. We call this phenomenon the dynamic on the accumulations curves. Such properties have already been observed by Colin de Verdi\`ere and Parisse \cite{CoPa94_02} for the eigenvalues in the presence of a degenerated double-well potential in dimension $n = 1$. Finally, note that the trapped set consists of two (resp. one) hyperbolic periodic trajectories for the energies below (resp. above) $E_{0}$.
\end{example}

\subsubsection{Transition phenomena} \label{s14}

The conclusions of Theorem \ref{d8} hold true for $\re z \in [ E_{0} - C h , E_{0} + C h ]$ where $C$ is any positive constant. Thus, it is natural to examine the asymptotic of the accumulation curves given by \eqref{d93} when $\re \sigma$ goes to $\pm \infty$. In this manner, we can hope to see the transition from the (non-)trapping situation below $E_{0}$ to the (non-)trapping situation above $E_{0}$. Indeed, these two situations (as well as the homoclinic regime at energy $E_{0}$) are completely different in general. This will be an illustration of the instability (which occurs generally) of the homoclinic regime.

We first investigate the behavior of the accumulation curves above $E_{0}$. From \ref{h2}, the potential $V (x)$ is a barrier of height $E_{0}$ near $x = 0$. Thus, it can not stop the Hamiltonian trajectories of energy higher than $E_{0}$ and these trajectories can pass through the origin. In other words, the transmission is favored in this case.

The simplest case for $E > E_{0}$ is then when the transmission can not occur in $\CH$. Stated another way, it means that $g_{+}^{k} \cdot g_{-}^{\ell} > 0$ for all $k , \ell \in \{ 1 ,\ldots , K \}$.  Since $g_{\pm}^{k}$ is the asymptotic directions of the base space projection of $\gamma_{k}$, we shall say that all the homoclinic trajectories are ``on the same side of $0$'' in this situation. Note that this is the case in Example \ref{c15}, Corollary \ref{e4}, Example \ref{e29} and Section \ref{s19} {\rm (A)}. To explain the link between the hypothesis $g_{+} \cdot g_{-} > 0$ and the fact the energies are non-trapped, we state a result on the classical dynamic under an additional assumption. Its proof can be found in Section \ref{s86}.

\begin{proposition}\sl \label{e24}
Assume \ref{h1}--\ref{h4}, \ref{h8}, $g_{+}^{k} \cdot g_{-}^{\ell} > 0$ for all $k , \ell \in \{ 1 ,\ldots , K \}$ and $V (x) = E_{0} - \sum \lambda_{j}^{2} x_{j}^{2} / 4$ near $0$. Then, there exists $\delta > 0$ such that the energies in $] E_{0} , E_{0} + \delta ]$ are non-trapped.
\end{proposition}

On the other hand, if all the homoclinic trajectories are on the same side of $0$, Theorem \ref{d8} together with \eqref{e36} imply the following result on the distribution of resonances. This resonance free region is the widest possible in view of \eqref{e3}.

\begin{remark}[Transition to the non-trapping regime]\sl \label{e23}
Assume $g_{+}^{k} \cdot g_{-}^{\ell} > 0$ for all $k , \ell \in \{ 1 ,\ldots , K \}$. Then, the accumulation curves given by \eqref{d93} verify
\begin{equation} \label{e26}
\im \sigma \leq - \sum_{j = 2}^{n} \frac{\lambda_{j}}{2} - \frac{\pi \re \sigma}{\vert \ln h \vert} + \CO \Big( \frac{1}{\vert \ln h \vert} \Big) ,
\end{equation}
as $\re \sigma \to + \infty$.
\end{remark}

Thus, the accumulation curves of resonances move away from the real axis. More precisely, the resonance free zone increases linearly, with a universal rate, as $\re \sigma$ goes to $+ \infty$. The growth of this region is consistent with the localization of the resonances for the energies higher than $E_{0}$. Indeed, near non-trapped energies, the imaginary part of the resonances is larger that $h \vert \ln h \vert$ from Martinez \cite{Ma02_01}, whereas the imaginary part of the resonances provided by Theorem \ref{d8} is of size $h$. Finally, note that the boundary of the region \eqref{e26} coincides with the curve $B$ of \cite[Th\'eor\`eme 0.1]{FuRa98_01}, outside of which there is no resonance in Example \ref{k42} in dimension $n = 1$.

If we now assume that all the trajectories are not on the same side of $0$, the transmission between incoming paths $\gamma_{\ell}$ and outgoing paths $\gamma_{k}$ such that $g_{+}^{k} \cdot g_{-}^{\ell} < 0$ plays the central role. Technically, it means that the main part of the coefficients of $\widehat{\CQ} ( \tau , h )$ behaves like
\begin{equation} \label{e27}
\Big\vert \Gamma \Big( \frac{1}{2} - i \frac{\tau}{\lambda_{1}} \Big) \big( i \lambda_{1} g_{+}^{k} \cdot g_{-}^{\ell} \big)^{- \frac{1}{2} + i \frac{\tau}{\lambda_{1}}} \Big\vert \sim \sqrt{\frac{2 \pi}{\lambda_{1} \vert g_{+}^{k} \cdot g_{-}^{\ell} \vert}}
\left\{ \begin{aligned}
&e^{- \frac{\pi \tau}{\lambda_{1}}} &&\text{ if } g_{+}^{k} \cdot g_{-}^{\ell} > 0 ,  \\
&1 &&\text{ if } g_{+}^{k} \cdot g_{-}^{\ell} < 0 ,
\end{aligned} \right.
\end{equation}
as $\tau \to + \infty$. Here, we have used again the formula \eqref{e36}. This implies that, for the asymptotic behavior of the accumulation curves in the limit $\re \sigma \to + \infty$, one can neglect the coefficients $( k , \ell )$ of $\widehat{\CQ}$ such that $g_{+}^{k} \cdot g_{-}^{\ell} > 0$.

Thus, in \eqref{e16} of Example \ref{e6}, we have seen that the accumulation curves become horizontal and coincide as $\re \sigma \to + \infty$. This is consistent with the localization of the resonances for the energies higher than $E_{0}$. Indeed, at least when $x_{1}$ and $x_{2}$ are collinear, one can verify that the trapped set consists here of a hyperbolic periodic trajectory for the energies above $E_{0}$. For such a trapping, G\'erard and Sj\"{o}strand \cite{GeSj87_01} have proved that the resonances form a lattice at distance $h$ from the real axis. This explains that the first accumulation curves have a horizontal asymptote as $\re \sigma \to + \infty$. That the two accumulation curves coincide could be explained by the non-vanishing of the symbol of the resonant states on the curves $\gamma_{\bullet}$ (see \eqref{j47}) which is typically a property that is true only for the first line of resonances generated by hyperbolic trappings (see Theorem 4.1 $iv)$ of \cite{BoFuRaZe11_01} for barrier-top resonances).

We now investigate the behavior of the resonances below $E_{0}$. Of course, the geometric setting is opposite to the previous one. Indeed, for Hamiltonian trajectories of energy lower than $E_{0}$, the potential $V (x)$ is an impassable barrier near $x = 0$. Thus, the reflexion is favored in this case.

In particular, if the incoming directions and outgoing directions are always ``on the opposite side of $0$'' (i.e. $g_{+}^{k} \cdot g_{-}^{\ell} < 0$ for all $k , \ell \in \{ 1 ,\ldots , K \}$), the accumulation curves of resonances move away from the real axis as $\re \sigma \to - \infty$. More precisely,  the conclusions of Remark \ref{e23} hold true with $\re \sigma$ replaced by $- \re \sigma$. As explained before, the parity of the symbol $p ( x , \xi )$ with respect to $\xi$ for Schr\"{o}dinger operators prevents the asymptotic directions from being on the opposite side of $0$, but in the setting of \eqref{n1} where $\CH = \emptyset$. Nevertheless, Example \ref{b81} provides a pseudodifferential operator satisfying this hypothesis. Note eventually that the energies below $E_{0}$ are non-trapped in these examples. Again, the behavior of the accumulation curves can be interpreted as a transition to a non-trapping situation.

If now some incoming directions and outgoing directions are on the same side of $0$, we expect a situation similar to the one of \eqref{e27}, mutatis mutandis. It means that the reflexion between incoming paths $\gamma_{\ell}$ and outgoing paths $\gamma_{k}$ such that $g_{+}^{k} \cdot g_{-}^{\ell} > 0$ are dominant since the main part of the coefficients of $\widehat{\CQ} ( \tau , h )$ behaves like
\begin{equation} \label{e28}
\Big\vert \Gamma \Big( \frac{1}{2} - i \frac{\tau}{\lambda_{1}} \Big) \big( i \lambda_{1} g_{+}^{k} \cdot g_{-}^{\ell} \big)^{- \frac{1}{2} + i \frac{\tau}{\lambda_{1}}} \Big\vert \sim \sqrt{\frac{2 \pi}{\lambda_{1} \vert g_{+}^{k} \cdot g_{-}^{\ell} \vert}}
\left\{ \begin{aligned}
&1 &&\text{ if } g_{+}^{k} \cdot g_{-}^{\ell} > 0 ,  \\
&e^{\frac{\pi \tau}{\lambda_{1}}} &&\text{ if } g_{+}^{k} \cdot g_{-}^{\ell} < 0 ,
\end{aligned} \right.
\end{equation}
as $\tau \to - \infty$.

The simplest situation with an incoming and an outgoing direction on the same side of $0$ is given by Example \ref{c15} where $\CH$ consists of a single trajectory. In this case, we have proved in Corollary \ref{e4} that the accumulation curve of resonance has a horizontal asymptote in the limit $\re \sigma \to - \infty$. Once more, this is in agreement with the nature of the trapped set for energies below $E_{0}$ which consists of a hyperbolic periodic trajectory. Indeed, as stated before, G\'erard and Sj\"{o}strand \cite{GeSj87_01} have proved that the resonances have an imaginary part of size $h$ for such trappings. Formula \eqref{e3} also shows that the accumulation curve reaches its asymptotes exponentially fast.

In Example \ref{e6}, the case of a unique trajectory on each side of $0$ is considered. The asymptotic behavior as $\re \tau \to - \infty$ of the resonances is given in \eqref{e9} in this situation. Comparing with \eqref{e8}, the set of resonances seems to be the union of the sets of resonances generated by each homoclinic trajectory separately. This illustrates the fact that the potential is an impassable barrier for energies below $E_{0}$: the communication between the two trajectories is impossible since they are on the opposite side of $0$. Moreover, this phenomenon can also be explained by the structure of the trapped set for energies below $E_{0}$. Indeed, it consists here of two separated hyperbolic periodic trajectories.

Eventually, we discuss the transition phenomenon in an example of operator with several homoclinic curves on the same side of $0$.

\begin{example}\rm \label{e29}
We construct here an example of operator $P$ satisfying the hypotheses of Theorem \ref{d8} with two or three trajectories on the same side of $0$. As in Example \ref{b86}, it will not be of the form \eqref{a5} but will enter in the setting of Remark \ref{c13}.

\begin{figure}
\begin{center}
\begin{picture}(0,0)%
\includegraphics{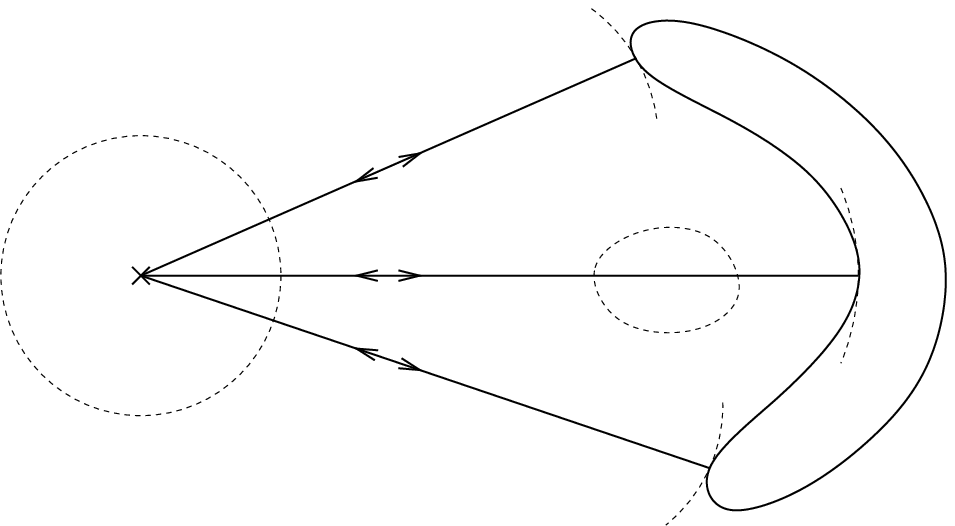}%
\end{picture}%
\setlength{\unitlength}{1105sp}%
\begingroup\makeatletter\ifx\SetFigFont\undefined%
\gdef\SetFigFont#1#2#3#4#5{%
  \reset@font\fontsize{#1}{#2pt}%
  \fontfamily{#3}\fontseries{#4}\fontshape{#5}%
  \selectfont}%
\fi\endgroup%
\begin{picture}(16246,8880)(-10814,-9451)
\put(-4049,-7336){\makebox(0,0)[b]{\smash{{\SetFigFont{9}{10.8}{\rmdefault}{\mddefault}{\updefault}$\pi_{x} ( \gamma_{2} )$}}}}
\put(-8399,-4786){\makebox(0,0)[b]{\smash{{\SetFigFont{9}{10.8}{\rmdefault}{\mddefault}{\updefault}$0$}}}}
\put(-8399,-2386){\makebox(0,0)[b]{\smash{{\SetFigFont{9}{10.8}{\rmdefault}{\mddefault}{\updefault}$\supp V$}}}}
\put(676,-3961){\makebox(0,0)[b]{\smash{{\SetFigFont{9}{10.8}{\rmdefault}{\mddefault}{\updefault}$\supp W$}}}}
\put(4426,-5236){\makebox(0,0)[b]{\smash{{\SetFigFont{9}{10.8}{\rmdefault}{\mddefault}{\updefault}$x_{3}$}}}}
\put(601,-1486){\makebox(0,0)[b]{\smash{{\SetFigFont{9}{10.8}{\rmdefault}{\mddefault}{\updefault}$x_{1}$}}}}
\put(1876,-8686){\makebox(0,0)[b]{\smash{{\SetFigFont{9}{10.8}{\rmdefault}{\mddefault}{\updefault}$x_{2}$}}}}
\put(2326,-2236){\makebox(0,0)[b]{\smash{{\SetFigFont{9}{10.8}{\rmdefault}{\mddefault}{\updefault}$\CO$}}}}
\put(-4049,-5761){\makebox(0,0)[b]{\smash{{\SetFigFont{9}{10.8}{\rmdefault}{\mddefault}{\updefault}$\pi_{x} ( \gamma_{3} )$}}}}
\put(-4049,-3961){\makebox(0,0)[b]{\smash{{\SetFigFont{9}{10.8}{\rmdefault}{\mddefault}{\updefault}$\pi_{x} ( \gamma_{1} )$}}}}
\end{picture}%
\end{center}
\bigskip \bigskip
\begin{center}
\begin{picture}(0,0)%
\includegraphics{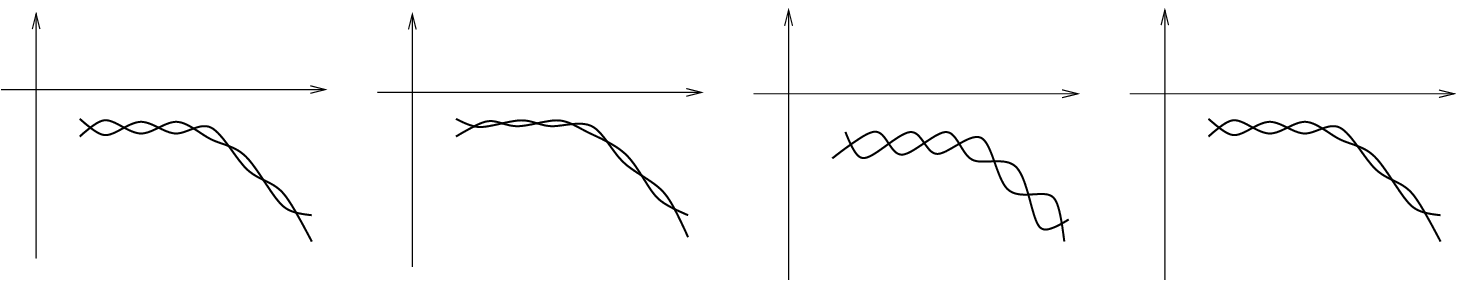}%
\end{picture}%
\setlength{\unitlength}{1105sp}%
\begingroup\makeatletter\ifx\SetFigFont\undefined%
\gdef\SetFigFont#1#2#3#4#5{%
  \reset@font\fontsize{#1}{#2pt}%
  \fontfamily{#3}\fontseries{#4}\fontshape{#5}%
  \selectfont}%
\fi\endgroup%
\begin{picture}(25019,4735)(-621,-5091)
\put(22051,-936){\makebox(0,0)[b]{\smash{{\SetFigFont{9}{10.8}{\rmdefault}{\mddefault}{\updefault}$\frac{A_{2} - A_{1}}{h} \equiv 2 \pi \text{ mod } 2 \pi$}}}}
\put(9151,-977){\makebox(0,0)[b]{\smash{{\SetFigFont{9}{10.8}{\rmdefault}{\mddefault}{\updefault}$\frac{A_{2} - A_{1}}{h} \equiv \frac{2}{3} \pi \text{ mod } 2 \pi$}}}}
\put(2701,-957){\makebox(0,0)[b]{\smash{{\SetFigFont{9}{10.8}{\rmdefault}{\mddefault}{\updefault}$\frac{A_{2} - A_{1}}{h} \equiv 0 \text{ mod } 2 \pi$}}}}
\put(15601,-941){\makebox(0,0)[b]{\smash{{\SetFigFont{9}{10.8}{\rmdefault}{\mddefault}{\updefault}$\frac{A_{2} - A_{1}}{h} \equiv \frac{4}{3} \pi \text{ mod } 2 \pi$}}}}
\end{picture}%
\end{center}
\caption{The geometry of Example \ref{e29} and the accumulation curves in the case of two symmetric trajectories.} \label{f28}
\end{figure}

In dimension $n = 2$, we consider the semiclassical operator
\begin{equation*}
P = - h^{2} \Delta_{\R^{2} \setminus \CO} + V (x) ,
\end{equation*}
with Dirichlet condition at the boundary of the obstacle $\CO$. Here, $V (x) \in C^{\infty}_{0} ( \R^{2} )$ is a potential as in Section \ref{s21} and $\CO$ is a small non-trapping ``croissant'' far away from the origin as illustrated in Figure \ref{f28}. In particular, the hypotheses \ref{h1} and \ref{h2} are satisfied and $\lambda_{1} = \lambda_{2} = : \lambda$.

On the other hand, Proposition \ref{a14} implies that the trapped set of $P$ at energy $E_{0}$ satisfies \ref{h3} and that the base space projection of the homoclinic trajectories consists of the radial rays which are normal to the boundary of $\CO$. Since $\CO$ is small, \ref{h4} holds true. Moreover, \ref{h8} is automatically satisfied if, for each point $x \in \partial \CO \cap \pi_{x} ( \CH )$, the curves $\partial B ( 0 , x )$ and $\partial \CO$ have a contact of order $1$. Thus, one can construct a situation with three homoclinic trajectories on the same side of $0$ as in Figure \ref{f28}. It may be more tricky to provide an example with only two homoclinic curves; but, as explained in Remark \ref{c13} $ii)$, on can add a potential $- i h \vert \ln h \vert W (x)$ to $P$ in order to ``remove'' a trajectory. Eventually, it should be possible to replace $\CO$ by a potential and to provide a Schr\"{o}dinger version of this example.

Using \eqref{d4}, Proposition \ref{d4}, Theorem \ref{d8} and \eqref{e21}, which gives the eigenvalues of a generic $2 \times 2$ matrix, we immediately obtain the explicit form of the resonances when there are two homoclinic trajectories. In particular, as stated in Remark \ref{e5}, the accumulation curves are periodic as function of $h^{- 1}$. Following the approach of Example \ref{e6}, we will now investigate the asymptotic of the accumulation curves in the limits $\re \sigma \to \pm \infty$.

From \eqref{e18}, we can factor $\widehat{\CQ} ( \tau , h )$ as
\begin{equation} \label{e30}
\widehat{\CQ} ( \tau , h ) = - i \frac{\lambda^{i \frac{\tau}{\lambda}}}{\sqrt{2 \pi}} \Gamma \Big( \frac{1}{2} - i \frac{\tau}{\lambda} \Big) e^{- \frac{\pi \tau}{2 \lambda}} \breve{\CQ} ( \tau , h ) ,
\end{equation}
where the coefficients of the $2 \times 2$ matrix $\breve{\CQ} ( \tau , h )$ are given by
\begin{equation*}
\breve{\CQ}_{k , \ell} ( \tau , h ) = e^{i A_{k} / h} \frac{\CM_{k}^{+}}{\CM_{k}^{-}} e^{- i \nu_{k} \frac{\pi}{2}} \big\vert g_{-}^{\ell} \big\vert \big( g_{+}^{k} \cdot g_{-}^{\ell} \big)^{- \frac{1}{2} + i \frac{\tau}{\lambda}} .
\end{equation*}
Recall that the accumulation curves only depend on the modulus of the eigenvalues of $\widehat{\CQ}$. As in \eqref{e27} and \eqref{e28}, the common factor satisfies
\begin{equation} \label{e31}
\Big\vert \frac{\lambda^{i \frac{\tau}{\lambda}}}{\sqrt{2 \pi}}  \Gamma \Big( \frac{1}{2} - i \frac{\tau}{\lambda} \Big) e^{- \frac{\pi \tau}{2 \lambda}} \Big\vert \sim
\left\{ \begin{aligned}
&e^{- \frac{\pi \tau}{\lambda}} &&\text{ as } \tau \to - \infty ,  \\
&1 &&\text{ as } \tau \to + \infty .
\end{aligned} \right.
\end{equation}
For $k = 1 , 2$, let $\breve{\mu}_{k} ( \tau , h )$ denote the eigenvalues of $\breve{\CQ} ( \tau , h )$. A direct computation shows that
\begin{equation*}
\vert \breve{\mu}_{1} ( \tau , h ) \vert \vert \breve{\mu}_{2} ( \tau , h ) \vert = \big\vert \det \breve{\CQ} ( \tau , h ) \big\vert = \frac{\CM_{1}^{+} \CM_{2}^{+}}{\CM_{1}^{-} \CM_{2}^{-}} \sqrt{\frac{\vert g_{-}^{1} \vert \vert g_{-}^{2} \vert}{\vert g_{+}^{1} \vert \vert g_{+}^{2} \vert}} \big\vert 1 - \cos ( x_{1} , x_{2} )^{- 1 + i \frac{2 \tau}{\lambda}} \big\vert \gtrsim 1 ,
\end{equation*}
since $\vert \cos ( x_{1} , x_{2} ) \vert < 1$ (the two trajectories $\gamma_{1}$ and $\gamma_{2}$ are different). On the other hand, the $\breve{\mu}_{\bullet} ( \tau , h )$'s are uniformly bounded since the matrix $\breve{\CQ}$ is uniformly bounded. Combining these estimates, there exists $R > 1$ such that
\begin{equation} \label{e32}
R^{-1} \leq \vert \breve{\mu}_{k} ( \tau , h ) \vert \leq R ,
\end{equation}
for $k = 1 , 2$, $\tau \in \R$ and $h \in ] 0 , 1 ]$.

From \eqref{d93} and the relations \eqref{e30} and \eqref{e31}, the two accumulation curves satisfy
\begin{equation*}
\im \sigma = -\frac{\lambda}{2} + \frac{1}{\vert \ln h \vert}
\left\{ \begin{aligned}
&\lambda \ln \big( \vert \breve{\mu}_{k} ( \re \sigma , h ) \vert \big) + o (1) &&\text{ as } \re \sigma \to - \infty , \\
&- \pi \re \sigma + \lambda \ln \big( \vert \breve{\mu}_{k} ( \re \sigma , h ) \vert \big) + o (1) &&\text{ as } \re \sigma \to + \infty , \\
\end{aligned} \right.
\end{equation*}
where the logarithm of the modulus of the $\breve{\mu}_{k}$'s is of order $1$ from \eqref{e32}. Even if \eqref{e21} provides explicit formulas for the eigenvalues $\breve{\mu}_{k}$, the expressions can be rather complicated. Thus, we make an additional symmetry assumption in order to obtain simple formulas.

We assume further that the geometric quantities ($A_{\bullet}$, $g_{\pm}^{\bullet}$, $\nu_{\bullet}$, \ldots) are the same for the two homoclinic trajectories. This can be achieved by taking $\CO$ symmetric with respect to $x_{3}$ in Figure \ref{f28}. In this case, we remove the subscript $k$ (except for $x_{1}$ and $x_{2}$). Here, the eigenvalues of $\breve{\CQ}$ can easily be calculated and the asymptotic of the accumulation curves becomes
\begin{equation*}
\im \sigma = -\frac{\lambda}{2} + \frac{1}{\vert \ln h \vert}
\left\{ \begin{aligned}
&C + \lambda \ln \big( \big\vert 1 \pm \cos ( x_{1} , x_{2} )^{- \frac{1}{2} + i \frac{\re \sigma}{\lambda}} \big\vert \big) + o (1) &&\text{as } \re \sigma \to - \infty , \\
&- \pi \re \sigma + C + \lambda \ln \big( \big\vert 1 \pm \cos ( x_{1} , x_{2} )^{- \frac{1}{2} + i \frac{\re \sigma}{\lambda}} \big\vert \big) + o (1) &&\text{as } \re \sigma \to + \infty , \\
\end{aligned} \right.
\end{equation*}
where, as in \eqref{e7} or \eqref{e17}, the constant $C$ is given by
\begin{equation*}
C = \lambda \ln \bigg( \frac{\CM^{+}}{\CM^{-}} \sqrt{\frac{\vert g_{-} \vert}{\vert g_{+} \vert}} \bigg) .
\end{equation*}

Since all the trajectories are on the same side of $0$, we observe again a transition to non-trapping energies as $\re \sigma \to + \infty$. But we see here that the accumulation curves generally do not have asymptote in the setting of Remark \ref{e23}. On the contrary, as $\re \sigma \to - \infty$, the accumulation curves do not escape to infinity neither tend to $0$. Moreover, in the symmetric case, these curves are (asymptotically) periodic with respect to $\re \sigma$. In some sense, everything happens as if the two homoclinic trajectories form a stable system of interactions. Once more, this example proves that the accumulation curves generally do not have asymptote as $\re \sigma \to - \infty$.

One can also consider the case of three homoclinic trajectories. Using Cardano's formula, one can again give explicit expressions for the pseudo-resonances which are of course more complicated. Nevertheless, one can verify that the accumulation curves satisfy similar properties.
\end{example}

Until now, we only have considered the transition phenomena at the macroscopic scale. But such phenomena occur also at the microscopic scale. For instance, we have already pointed out the dynamic of the resonances on the accumulation curves at the end of Example \ref{e6}. Nevertheless, some transition properties do not occur in the intervals of size $h$. For example, the space between packets of resonances under the assumption \ref{h9} is always less than $2 \pi \lambda_{1} h \vert \ln h \vert^{- 1}$ in the limit $\tau \to + \infty$ (see \eqref{e8}). On the other hand, the space between resonances associated to hyperbolic trappings is of order $h$ (see \cite{GeSj87_01}). Thus, we do not observe the expected increase in the size of the gaps which means that this transition does not hold in intervals like $[ E_{0} - C h , E_{0} + C h ]$. This is not surprising since the eigenvalues for homoclinic orbits in dimension $n = 1$ satisfy similar properties (see Labl{\'e}e \cite{La10_01}).

\subsubsection{Stability phenomena} \label{s24}

As we have seen in the previous part, the resonances are very sensitive to the energy. Their imaginary part can change radically between $E_{0}$ and $E_{0} + \varepsilon$. In this sense, the resonances are unstable. Nevertheless, we will show that they are stable in the class of homoclinic trapped set.

More precisely, let $P_{0}$ be an operator satisfying \ref{h1}--\ref{h4}, \ref{h8}. It will be our reference operator. We denote by $\CH_{0} = \gamma_{1} \cup \cdots \cup \gamma_{K_{0}}$ the homoclinic set of $P_{0}$. We now consider another operator $P$ satisfying also \ref{h1}--\ref{h4}, \ref{h8} and such that the symbol of $P$ coincides with the symbol of $P_{0}$ in a neighborhood of $\CH_{0}$. In particular, $\CH_{0} \subset \CH = \gamma_{1} \cup \cdots \cup \gamma_{K}$ and the exceptional set $\Gamma (h)$ is the same for $P_{0}$ and $P$. This last operator can be seen as a perturbation of $P_{0}$. The following result shows that the resonances of $P$ are close to those of $P_{0}$ if the perturbation (corresponding to the additional homoclinic trajectories $\gamma_{K_{0} + 1} , \ldots , \gamma_{K}$) is ``small enough''. Section \ref{s25} collects the proofs of the results of this part.

\begin{proposition}[Stability with respect to the trapped set]\sl \label{i38}
Let $P_{0}$ be fixed as before and let $C , \delta > 0$ be any given positive numbers. There exists $\varepsilon > 0$ such that, if $P$ is of the previous form with
\begin{equation} \label{i40}
\max_{k \in \{ K_{0} + 1 , \ldots , K \}} \frac{\CM^{+}_{k}}{\CM^{-}_{k}} \sqrt{\frac{\vert g_{-}^{k} \vert}{\vert g_{+}^{k} \vert}} \leq \varepsilon \qquad \text{and} \qquad \min_{\fract{k \in \{ K_{0} + 1 , \ldots , K \}}{\ell \in \{ 1 , \ldots , K \}}} \big\vert \cos ( g_{+}^{k} , g_{-}^{\ell} ) \big\vert \geq \delta ,
\end{equation}
then, in the domain
\begin{equation*}
E_{0} + [ - C h , C  h ] + i \Big[ - \sum_{j = 2}^{n} \frac{\lambda_{j}}{2} h - C \frac{h}{\vert \ln h \vert} , h \Big] \setminus \big( \Gamma (h) + B ( 0 , \delta h ) \big) ,
\end{equation*}
we have
\begin{equation} \label{i48}
\dist \big( \res ( P ) , \res ( P_{0} ) \big) \leq \delta \frac {h}{\vert \ln h \vert} ,
\end{equation}
for $h$ small enough.
\end{proposition}

\begin{figure}
\begin{center}
\begin{picture}(0,0)%
\includegraphics{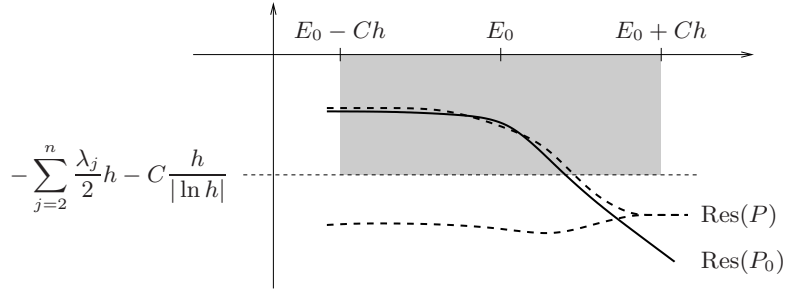}%
\end{picture}%
\setlength{\unitlength}{1105sp}%
\begingroup\makeatletter\ifx\SetFigFont\undefined%
\gdef\SetFigFont#1#2#3#4#5{%
  \reset@font\fontsize{#1}{#2pt}%
  \fontfamily{#3}\fontseries{#4}\fontshape{#5}%
  \selectfont}%
\fi\endgroup%
\begin{picture}(14287,6494)(-2264,-6233)
\put(6301,-586){\makebox(0,0)[b]{\smash{{\SetFigFont{9}{10.8}{\rmdefault}{\mddefault}{\updefault}$E_{0}$}}}}
\put(2701,-586){\makebox(0,0)[b]{\smash{{\SetFigFont{9}{10.8}{\rmdefault}{\mddefault}{\updefault}$E_{0} - C h$}}}}
\put(9901,-586){\makebox(0,0)[b]{\smash{{\SetFigFont{9}{10.8}{\rmdefault}{\mddefault}{\updefault}$E_{0} + C h$}}}}
\put(10801,-4711){\makebox(0,0)[lb]{\smash{{\SetFigFont{9}{10.8}{\rmdefault}{\mddefault}{\updefault}$\Res ( P )$}}}}
\put(10801,-5761){\makebox(0,0)[lb]{\smash{{\SetFigFont{9}{10.8}{\rmdefault}{\mddefault}{\updefault}$\Res ( P_{0} )$}}}}
\put(-2249,-3811){\makebox(0,0)[b]{\smash{{\SetFigFont{9}{10.8}{\rmdefault}{\mddefault}{\updefault}$\ds - \sum_{j = 2}^{n} \frac{\lambda_{j}}{2} h - C \frac{h}{\vert \ln h \vert}$}}}}
\end{picture}%
\end{center}
\caption{The stability of the resonances in Example \ref{i41}.} \label{f35}
\end{figure}

Thus, the resonances of $P$ near the real axis are close to those of $P_{0}$. This phenomenon is illustrated in Figure \ref{f35} for $P_{0}$ and $P$ as in Example \ref{i41}. Note that the stability is not verified below the line $\im z = - h \sum_{j = 2}^{n} \lambda_{j} / 2 - C h \vert \ln h \vert^{- 1}$. This is natural since the resonances (just) below this line are given by the logarithm of the small eigenvalues of $\CQ$ which is unstable. Note that the condition \eqref{i40} is independent of the parametrization of the curves $\gamma_{\bullet}$ (see the proof of Remark \ref{d3} $i)$).

\begin{remark}\sl
Proposition \ref{i38} still holds true if we do not assume that $P$ coincides with $P_{0}$ microlocally near $\CH_{0}$ but suppose that the first $K_{0}$ homoclinic trajectories of $P$ are ``sufficiently close'' to the ones of $P_{0}$. In particular, the classical action \eqref{r5} for $P_{0}$ and $P$ must be the same along these trajectories.
\end{remark}

\begin{example}\rm \label{i41}
We go back to Example \ref{e6}. Consider
\begin{equation*}
P_{0} = - h^{2} \Delta + V_{0} ( x ) + V_{1} ( x - x_{1} ) \qquad \text{and} \qquad P = P_{0} + V_{2} ( x - x_{2} ) ,
\end{equation*}
where the potentials $V_{\bullet}$ are as in Example \ref{e6}. These operators satisfy \ref{h1}--\ref{h4} and \ref{h8}. Moreover, $\CH_{0} = \gamma_{1}$ and $\CH = \gamma_{1} \cup \gamma_{2}$. Since $\CH_{0}$ consists of a unique trajectory, the asymptotic of the resonances of $P_{0}$ is given below Remark \ref{e2}. On the other hand, Example \ref{e6} provides the resonances of $P$.

We now explain how to construct $V_{2}$ such that \eqref{i40} holds true. The first idea is to consider a well-chosen family of potentials $V_{2}$ such that $\{ V_{2} ( x ) = E_{0} \}$ ``goes to $\{ 0 \}$''. This can be done by taking $V_{2} ( x ) = ( 1 + \varepsilon ) \widetilde{V}_{2} ( x )$ where $\widetilde{V}_{2}$ is a non-degenerate bump of height $\widetilde{V}_{2} ( 0 ) = E_{0}$. In this case, one can prove that \eqref{i40} is satisfied. This can be more easily done replacing $V_{2} ( x - x_{2} )$ by a small obstacle contracting to $x_{2}$ (one can take $B ( x_{2} , \varepsilon )$ for instance). Another approach is to add to $P$ an absorbing potential $i h W (x)$ supported on $\pi_{x} ( \gamma_{2} )$ in order to ``weaken the contribution'' of $\gamma_{2}$ (see Remark \ref{c13} $ii)$). Figure \ref{f35} summarizes the distribution of resonances when \eqref{i40} is verified.
\end{example}

In Proposition \ref{i38}, the trapped set at energy $E_{0}$ of $P$ contains that of $P_{0}$. In other words, $P$ is ``more geometrically trapping'' than $P_{0}$. By Bohr's correspondence principle, it is natural to think that $P$ will be ``more quantum trapping'' than $P_{0}$, namely that the resonances will be closer to the real axis. This definition of quantum trapping is justified for instance by the resonance expansion of the propagator (see e.g. Lax and Phillips \cite[Theorem III.5.4]{LaPh67_01}). Nevertheless, this intuition lead by the correspondence principle is inexact.

\begin{figure}
\begin{center}
\begin{picture}(0,0)%
\includegraphics{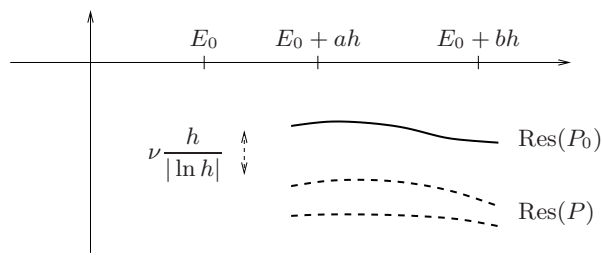}%
\end{picture}%
\setlength{\unitlength}{1105sp}%
\begingroup\makeatletter\ifx\SetFigFont\undefined%
\gdef\SetFigFont#1#2#3#4#5{%
  \reset@font\fontsize{#1}{#2pt}%
  \fontfamily{#3}\fontseries{#4}\fontshape{#5}%
  \selectfont}%
\fi\endgroup%
\begin{picture}(12644,5519)(-621,-5258)
\put(10801,-4486){\makebox(0,0)[lb]{\smash{{\SetFigFont{9}{10.8}{\rmdefault}{\mddefault}{\updefault}$\Res ( P )$}}}}
\put(3751,-586){\makebox(0,0)[b]{\smash{{\SetFigFont{9}{10.8}{\rmdefault}{\mddefault}{\updefault}$E_{0}$}}}}
\put(6301,-586){\makebox(0,0)[b]{\smash{{\SetFigFont{9}{10.8}{\rmdefault}{\mddefault}{\updefault}$E_{0} + a h$}}}}
\put(9901,-586){\makebox(0,0)[b]{\smash{{\SetFigFont{9}{10.8}{\rmdefault}{\mddefault}{\updefault}$E_{0} + b h$}}}}
\put(4201,-3136){\makebox(0,0)[rb]{\smash{{\SetFigFont{9}{10.8}{\rmdefault}{\mddefault}{\updefault}$\ds \nu \frac{h}{\vert \ln h \vert}$}}}}
\put(10801,-2836){\makebox(0,0)[lb]{\smash{{\SetFigFont{9}{10.8}{\rmdefault}{\mddefault}{\updefault}$\Res ( P_{0} )$}}}}
\end{picture}%
\end{center}
\caption{The accumulations curves of Remark \ref{i42}.} \label{f36}
\end{figure}

\begin{remark}[More geometric trapping does not always imply more quantum trapping]\sl \label{i42}
Let $P_{0}$ be as in Example \ref{i41}. For all $\delta > 0$, there exist $a < b$, $\nu > 0$, a perturbation $P$ of $P_{0}$ with $K > K_{0}$ and a sequence of positive numbers $h$ which converges to $0$ such that $P_{0}$ and $P$ satisfy the assumptions of Proposition \ref{i38} and
\begin{align*}
\max \big\{ \im z ; \ z & \in \res ( P ) \text{ with } \re z \in E_{0} + [ a h , b h ] \big\} \\
&< \max \big\{ \im z ; \ z \in \res ( P_{0} ) \text{ with } \re z \in E_{0} + [ a h , b h ] \big\} - \nu \frac{h}{\vert \ln h \vert} ,
\end{align*}
for $h$ in this sequence (see Figure \ref{f36}).
\end{remark}

This property is the consequence of the vibration phenomena. Due to the presence of the action factor $e^{i A_{k} / h}$ in the expression of $\CQ$, the additional homoclinic trajectories $\CH \setminus \CH_{0}$ increase or decrease (depending on $h$) the imaginary part of the resonances. Thus, in other situations, more geometric trapping gives more quantum trapping (see e.g. \eqref{e35}). More generally, these phenomena show that the nature of the trapped set (i.e. the geometric setting) is not enough to characterize the distribution of resonances and that some quantities (i.e. the dynamical setting) are needed.

Remark \ref{i42} is stated in the setting of small perturbations. But this result also holds true for perturbation of size of the original trapping (see e.g. Example \ref{e48} {\rm (B)}). In this case, the gap between the resonance of $P_{0}$ and $P$ can be even larger (since the stability given by Proposition \ref{i38} is not verified). Eventually, note that $a$ and $b$ can be chosen arbitrarily in Example \ref{e48} {\rm (B)}.

\Subsection{The one dimensional situation} \label{s19}

\begin{figure}
\begin{center}
\begin{picture}(0,0)%
\includegraphics{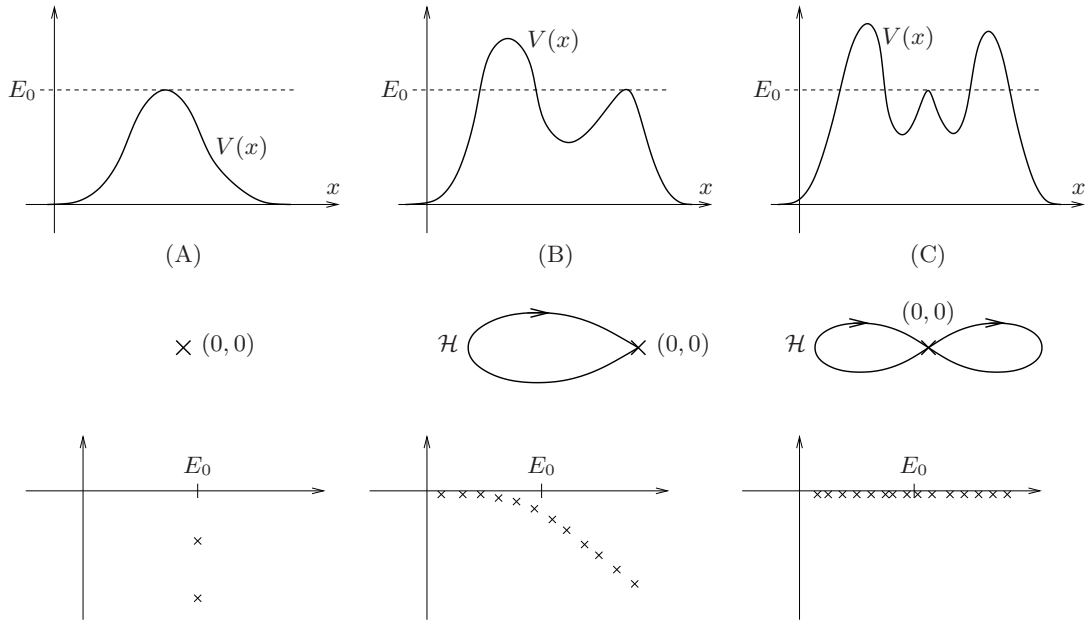}%
\end{picture}%
\setlength{\unitlength}{1184sp}%
\begingroup\makeatletter\ifx\SetFigFont\undefined%
\gdef\SetFigFont#1#2#3#4#5{%
  \reset@font\fontsize{#1}{#2pt}%
  \fontfamily{#3}\fontseries{#4}\fontshape{#5}%
  \selectfont}%
\fi\endgroup%
\begin{picture}(22612,12944)(211,-13883)
\put(4201,-10786){\makebox(0,0)[b]{\smash{{\SetFigFont{9}{10.8}{\rmdefault}{\mddefault}{\updefault}$E_{0}$}}}}
\put(8026,-2911){\makebox(0,0)[lb]{\smash{{\SetFigFont{9}{10.8}{\rmdefault}{\mddefault}{\updefault}$E_{0}$}}}}
\put(226,-2911){\makebox(0,0)[lb]{\smash{{\SetFigFont{9}{10.8}{\rmdefault}{\mddefault}{\updefault}$E_{0}$}}}}
\put(15826,-2911){\makebox(0,0)[lb]{\smash{{\SetFigFont{9}{10.8}{\rmdefault}{\mddefault}{\updefault}$E_{0}$}}}}
\put(19501,-6361){\makebox(0,0)[b]{\smash{{\SetFigFont{9}{10.8}{\rmdefault}{\mddefault}{\updefault}(C)}}}}
\put(3901,-6361){\makebox(0,0)[b]{\smash{{\SetFigFont{9}{10.8}{\rmdefault}{\mddefault}{\updefault}(A)}}}}
\put(11701,-6361){\makebox(0,0)[b]{\smash{{\SetFigFont{9}{10.8}{\rmdefault}{\mddefault}{\updefault}(B)}}}}
\put(11101,-1861){\makebox(0,0)[lb]{\smash{{\SetFigFont{9}{10.8}{\rmdefault}{\mddefault}{\updefault}$V (x)$}}}}
\put(4576,-4111){\makebox(0,0)[lb]{\smash{{\SetFigFont{9}{10.8}{\rmdefault}{\mddefault}{\updefault}$V (x)$}}}}
\put(18526,-1786){\makebox(0,0)[lb]{\smash{{\SetFigFont{9}{10.8}{\rmdefault}{\mddefault}{\updefault}$V (x)$}}}}
\put(14701,-4936){\makebox(0,0)[lb]{\smash{{\SetFigFont{9}{10.8}{\rmdefault}{\mddefault}{\updefault}$x$}}}}
\put(22501,-4936){\makebox(0,0)[lb]{\smash{{\SetFigFont{9}{10.8}{\rmdefault}{\mddefault}{\updefault}$x$}}}}
\put(6901,-4936){\makebox(0,0)[lb]{\smash{{\SetFigFont{9}{10.8}{\rmdefault}{\mddefault}{\updefault}$x$}}}}
\put(4276,-8236){\makebox(0,0)[lb]{\smash{{\SetFigFont{9}{10.8}{\rmdefault}{\mddefault}{\updefault}$( 0 , 0 )$}}}}
\put(13801,-8236){\makebox(0,0)[lb]{\smash{{\SetFigFont{9}{10.8}{\rmdefault}{\mddefault}{\updefault}$( 0 , 0 )$}}}}
\put(9226,-8236){\makebox(0,0)[lb]{\smash{{\SetFigFont{9}{10.8}{\rmdefault}{\mddefault}{\updefault}$\CH$}}}}
\put(19501,-7561){\makebox(0,0)[b]{\smash{{\SetFigFont{9}{10.8}{\rmdefault}{\mddefault}{\updefault}$( 0 , 0 )$}}}}
\put(16501,-8236){\makebox(0,0)[lb]{\smash{{\SetFigFont{9}{10.8}{\rmdefault}{\mddefault}{\updefault}$\CH$}}}}
\put(11401,-10786){\makebox(0,0)[b]{\smash{{\SetFigFont{9}{10.8}{\rmdefault}{\mddefault}{\updefault}$E_{0}$}}}}
\put(19201,-10786){\makebox(0,0)[b]{\smash{{\SetFigFont{9}{10.8}{\rmdefault}{\mddefault}{\updefault}$E_{0}$}}}}
\end{picture}%
\end{center}
\caption{The three possible situations in dimension $1$ with the corresponding trapped sets and distributions of resonances.} \label{f6}
\end{figure}

In this part, we apply the results obtained previously in Section \ref{s22} and Section \ref{s6} in dimension $n = 1$. In this case, there exist only three possible geometric settings (modulo symmetries) whose trapped set consists of homoclinic curves (i.e. that satisfy the assumptions \ref{h1}--\ref{h3}). They are illustrated in Figure \ref{f6}.

In case {\rm (A)}, the trapped set at energy $E_{0}$ is reduced to the fixed point and then $\CH = \emptyset$. In particular, the assumptions \ref{h4} and \ref{h8} are automatically satisfied and Theorem \ref{d8} can be applied. As there is no homoclinic trajectory, $\CQ$ is the trivial operator on $\C^{0} = \{ 0 \}$ and Theorem \ref{d8} implies that $P$ has no resonance and a polynomial estimate of its resolvent in the set \eqref{d90}. This is in agreement with the computation of the resonances at the barrier-top made by Briet, Combes and Duclos \cite{BrCoDu87_02}, Sj\"{o}strand \cite{Sj87_01} and the third author \cite{Ra96_01}. In any domain of the form $B ( E_{0} , C h )$ with $C > 0$, these resonances are simple and verify
\begin{equation*}
z_{k} = E_{0} - i h \lambda_{1} \Big( \frac{1}{2} + k \Big) + \CO ( h^{2} ) ,
\end{equation*}
for $k \in \N$ and $h$ small enough. See \eqref{n1} in the general $n$ dimensional case. Thus, the resonances coincide with $\Gamma_{0} (h)$ up to lower order terms. This justifies that neighborhoods of these points (or of $\Gamma (h)$ for technical reasons) are removed from the set \eqref{d90}. Note also that the polynomial resolvent estimate in the complex plane has already been obtained by Michel and the first author \cite[Theorem 6.4]{BoMi04_01} and in a previous paper \cite[Theorem 3.1]{BoFuRaZe11_01}. Of course, the discussions of Section \ref{s18} are pointless here since $P$ has no accumulation curve. Finally, Theorem \ref{a1} and Theorem \ref{a2} (see also Remark \ref{b69} $i)$) can also be applied here since \ref{h5b} holds true.

In case {\rm (B)}, the homoclinic set consists of a single trajectory. In dimension $1$, $g_{\pm}$ do not vanish and are collinear to the unique spatial axis. Moreover, the intersection between $\Lambda_{-}$ and $\Lambda_{+}$ is automatically transversal since these manifolds are here curves. Then, \ref{h4} and \ref{h8} always hold true in dimension $1$. Hence, we can apply again Theorem \ref{d8} and also Corollary \ref{e4} since \ref{h9} is verified here. As in that case, we remove the subscript $k = 1$ labeling the trajectory. We first compute the geometric quantities $\CM^{\pm}$ defined in \eqref{d7}. Since \eqref{d6} is valid in the expandible sense (see \eqref{m92}), we obtain
\begin{equation} \label{e38}
\CM^{+} = \lim_{s \to - \infty} \sqrt{\big\vert \partial_{s} x_{1} ( s ) \big\vert} e^{- s \lambda_{1} / 2} = \lim_{s \to - \infty} \sqrt{\lambda_{1} \vert g_{+} \vert e^{s \lambda_{1}}} e^{- s \lambda_{1} / 2} = \sqrt{\lambda_{1} \vert g_{+} \vert} ,
\end{equation}
and the same way $\CM^{-} = \sqrt{\lambda_{1} \vert g_{-} \vert}$. Thus, the constant defined in \eqref{e7} satisfies $C_{0} = 0$. On the other hand, the Maslov's index of $\gamma$ is $\nu = - 1$. Combining the previous relations with \eqref{e8}, the points $z_{q} ( \tau )$ of Proposition \ref{d9} which approximate the (pseudo-)resonances modulo $o ( h \vert \ln h \vert^{- 1} )$ are given by
\begin{equation} \label{e37}
z_{q} ( \tau ) =  E_{0} - \frac{A \lambda_{1}}{\vert \ln h \vert} + 2 q \pi \lambda_{1} \frac{h}{\vert \ln h \vert} + i \ln ( \mu ( \tau ) ) \lambda_{1} \frac{h}{\vert \ln h \vert} ,
\end{equation}
with
\begin{equation*}
\mu ( \tau ) = \frac{1}{\sqrt{2 \pi}} \Gamma \Big( \frac{1}{2} - i \frac{\tau}{\lambda_{1}} \Big) e^{- \frac{\pi \tau}{2 \lambda_{1}}} \big( \lambda_{1} \vert g_{+} \vert \vert g_{-} \vert \big)^{i \frac{\tau}{\lambda_{1}} } .
\end{equation*}
In particular, the imaginary part of the resonances is of size $h \vert \ln h \vert^{- 1}$ whereas it is of size $h$ in dimension $n \geq 2$.

We now consider the resonances in the window $\re z = E_{0} + o (h)$. For that, we take $\tau = 0$ in \eqref{e37} and obtain
\begin{equation} \label{e39}
z_{q} ( 0 ) =  E_{0} - \frac{A \lambda_{1}}{\vert \ln h \vert} + 2 q \pi \lambda_{1} \frac{h}{\vert \ln h \vert} - i \lambda_{1} \frac{\ln 2}{2} \frac{h}{\vert \ln h \vert} .
\end{equation}
This is consistent with the computation of the resonances made by the second and third authors in dimension $n = 1$. Indeed, they have proved in \cite[Th\'eor\`eme 0.7]{FuRa98_01} that the resonances in $B ( E_{0} , C h \vert \ln h \vert^{- 1} )$ with $C > 0$ satisfy \eqref{e39} modulo $\CO ( h \vert \ln h \vert^{- 2} )$. As a matter of fact, this result concerns the heteroclinic curves but, transposed to the case of homoclinic curves, it gives the previous asymptotic. Note also that their theorem shows that $P$ has no resonance near $\Gamma_{0} (h)$.

On the other hand, Remark \ref{e2} provides a unique accumulation curve in this case. From \eqref{e36}, \eqref{e37} and Corollary \ref{e4}, it satisfies
\begin{align}
\im \sigma &= - \frac{\lambda_{1}}{2 \vert \ln h \vert} \ln \big( e^{\frac{2 \pi}{\lambda_{1}} \re \sigma} + 1 \big)    \nonumber \\
&= - \frac{1}{\vert \ln h \vert}
\left\{ \begin{aligned}
&\frac{\lambda_{1}}{2} e^{\frac{2 \pi}{\lambda_{1}} \re \sigma} + \CO \big( e^{\frac{4 \pi}{\lambda_{1}} \re \sigma} \big) &&\text{ as } \re \sigma \to - \infty , \\
&\pi \re \sigma + \CO \big( e^{- \frac{2 \pi}{\lambda_{1}} \re \sigma} \big) &&\text{ as } \re \sigma \to + \infty . \\
\end{aligned} \right.   \label{e40}
\end{align}
As already stated in Section \ref{s14}, the behavior of this curve as $\re \sigma \to + \infty$ can be explained by the transition to non-trapping. The asymptotic behavior of the accumulation curve as $\re \sigma \to - \infty$ is also consistent with the localization of the resonances for the energies lower than $E_{0}$. Indeed, we are in a well in the island situation below $E_{0}$. For such situations, the imaginary part of the resonances satisfies
\begin{equation} \label{e41}
\im z \approx h e^{- 2 \ag ( \re z ) / h} .
\end{equation}
where $\ag ( E ) = \int \sqrt{( V (x) - E )_{+}} d x$ denotes the Agmon distance between the well and the sea. We refer to Servat \cite[Th\'eor\`eme 2.6]{Se04_01} in this one dimensional case and to Helffer and Sj\"{o}strand \cite[Th\'eor\`eme 10.12]{HeSj86_01} for punctual wells in any dimension. Thus, that the resonances are exponentially small below $E_{0}$ justifies that the accumulation curve \eqref{e40} is asymptote to the real axis as $\re \sigma \to - \infty$. Moreover, a direct computation gives $\ag ( \re z ) \sim - \pi ( \re z - E_{0} ) / \lambda_{1}$ as $\re z \nearrow E_{0}$. So, if we formally apply \eqref{e41} for $\re z - E_{0}$ of size $h$, we obtain
\begin{equation} \label{e42}
\im \sigma \approx e^{\frac{2 \pi}{\lambda_{1}} \re \sigma} .
\end{equation}
Thus, the exponential factors (which are of order $1$ in the studied domain) of \eqref{e40} and \eqref{e42} coincide. The additional factor $\vert \ln h \vert^{- 1}$ in \eqref{e40} may be perhaps related to the concentration at the barrier-top of the eigenvalues generated by homoclinic orbits in dimension $n = 1$ (see Colin de Verdi\`ere and Parisse \cite[Th\'eor\`eme 20]{CoPa94_01}). Eventually, Th\'eor\`eme 0.1 of \cite{FuRa98_01} shows that the operator $P$ has no resonance at distance $1$ of the curve $B$ defined by
\begin{equation}
\im z =
\left\{ \begin{aligned}
&0 &&\text{ for } \re z < E_{0} ,  \\
&\pi \frac{\re z - E_{0}}{\ln ( \re z - E_{0} )} &&\text{ for } \re z > E_{0} ,
\end{aligned} \right.
\end{equation}
up to lower order terms. If we consider this curve for $z - E_{0}$ of order $h$, we obtain a good approximation of the accumulation curve \eqref{e40}.

Finally, we discuss the assumptions of Section \ref{s22}. In this particular geometric situation, neither \ref{h5a} (since there exists only one $\lambda_{\bullet}$) nor \ref{h5b} (since $\Lambda_{-}$ and $\Lambda_{+}$ coincide along the homoclinic set) is true. So, Theorem \ref{a1} can not be applied here. Moreover, its conclusions do not hold since the imaginary part of the resonances given in \eqref{e37} is of order $h \vert \ln h \vert^{- 1}$. With the notations of Section \ref{s77}, we have $\S^{0} = \{ - 1 , 1 \}$ and $\CH_{\rm tang}^{\pm \infty} = \{ - 1 \}$. Thus, $\CT_{0} ( \tau )$ is a linear operator on $\C$ and
\begin{equation*}
\CA_{0} ( \tau ) = ( 2 \pi )^{- \frac{1}{2}} \Big\vert \Gamma \Big( \frac{1}{2} - i \frac{\tau}{\lambda_{1}} \Big) \Big\vert \CM_{0} \CJ_{0} ( \tau ) ,
\end{equation*}
from \eqref{b70}. Moreover, $\CM_{0} = 1$ and $\CJ_{0} ( \tau ) = e^{- \frac{\pi}{2 \lambda_{1}} \tau}$ by definition. Combining with \eqref{e36}, it yields
\begin{equation} \label{b75}
\CA_{0} ( \tau ) = \frac{e^{- \frac{\pi}{2 \lambda_{1}} \tau}}{\sqrt{2 \cosh ( \pi \tau / \lambda_{1} )}} .
\end{equation}
Comparing with \eqref{e40}, we see that Theorem \ref{a2} provides in this setting the sharp resonance free domain below any interval $[ E_{0}- C h , E_{0} + C h ]$ with $C > 0$.

In case {\rm (C)}, the homoclinic set consists of two homoclinic trajectories which are not on the same side of $0$. As explained for the case {\rm (B)}, the assumptions \ref{h4} and \ref{h8} are always satisfied in dimension $n = 1$. Then, Theorem \ref{d8} provides again the asymptotic of the resonances. Furthermore, note that the present situation is the one dimensional version of Example \ref{e6}. Using the computation of the geometric quantities made in \eqref{e38}, the $2 \times 2$ matrix $\widehat{\CQ}$ defined in \eqref{e18} satisfies
\begin{align}
G_{-}^{\frac{1}{2}} \widehat{\CQ} ( \tau & , h ) G_{-}^{- \frac{1}{2}}   \nonumber \\
&= \frac{\lambda_{1}^{i \frac{\tau}{\lambda_{1}}}}{\sqrt{2 \pi}} \Gamma \Big( \frac{1}{2} - i \frac{\tau}{\lambda_{1}} \Big)
\left( \begin{array}{cc}
e^{- \frac{\pi \tau}{2 \lambda_{1}}} e^{i A_{1} / h} \big( \vert g_{+}^{1} \vert \vert g_{-}^{1} \vert \big)^{i \frac{\tau}{\lambda_{1}}} & i e^{\frac{\pi \tau}{2 \lambda_{1}}} e^{i A_{1} / h} \big( \vert g_{+}^{1} \vert \vert g_{-}^{2} \vert \big)^{i \frac{\tau}{\lambda_{1}}} \\
i e^{\frac{\pi \tau}{2 \lambda_{1}}} e^{i A_{2} / h} \big( \vert g_{+}^{2} \vert \vert g_{-}^{1} \vert \big)^{i \frac{\tau}{\lambda_{1}}} &  e^{- \frac{\pi \tau}{2 \lambda_{1}}} e^{i A_{2} / h} \big( \vert g_{+}^{2} \vert \vert g_{-}^{2} \vert \big)^{i \frac{\tau}{\lambda_{1}}}
\end{array} \right) , \label{e43}
\end{align}
with the $2 \times 2$ matrix $G_{-} = \diag ( \vert g_{-}^{k} \vert )$. Combining with Proposition \ref{d9}, Theorem \ref{d8} and \eqref{e21} which gives the eigenvalues of a generic $2 \times 2$ matrix, we can obtain the explicit asymptotic of the resonances.

We now study the imaginary part of the resonances. We first claim that
\begin{equation} \label{e44}
\forall \tau \in \R , \quad \forall h \in ] 0 , 1] , \quad \forall k \in \{ 1 , 2 \} , \qquad \vert \mu_{k} ( \tau , h ) \vert = 1 .
\end{equation}
For that, we first consider a generic matrix of the form
\begin{equation} \label{e45}
B = \left( \begin{array}{cc}
x a c & i \frac{a d}{x}  \\
i \frac{b c}{x} & x b d
\end{array} \right) ,
\end{equation}
with $x \in ] 0 , + \infty [$ and $a , b , c , d \in \S^{1} \subset \C$. In particular, a direct calculus shows that
\begin{equation*}
K : = \frac{a c + b d}{\sqrt{a b c d}} \in [ - 2 , 2 ]
\end{equation*}
Then, using \eqref{e21}, the two eigenvalues of $B$, denoted $\mu_{\pm} (B)$ satisfies
\begin{align}
4 \vert \mu_{\pm} ( B ) \vert^{2} &= \big\vert ( a c + b d ) x \pm \sqrt{( a c + b d )^{2} x^{2} - 4 a b c d x^{2} - 4 a b c d x^{- 2}} \big\vert^{2}    \nonumber \\
&= \big\vert K x \pm \sqrt{K^{2} x^{2} - 4 x^{2} - 4 x^{- 2}} \big\vert^{2}   \nonumber \\
&= K^{2} x^{2} + 4 x^{2} + 4 x^{- 2} - K^{2} x^{2}   \nonumber \\
&= 4 x^{2} + 4 x^{- 2} . \label{e46}
\end{align}
Since the right hand side of \eqref{e43} is of the form \eqref{e45} modulo a scalar term, \eqref{e46} together with \eqref{e36} give
\begin{equation*}
\vert \mu_{k} ( \tau , h ) \vert = \frac{1}{\sqrt{2 \pi}} \Big\vert \Gamma \Big( \frac{1}{2} - i \frac{\tau}{\lambda_{1}} \Big) \big\vert \sqrt{e^{\frac{\pi \tau}{\lambda_{1}}} + e^{- \frac{\pi \tau}{\lambda_{1}}}} = 1 ,
\end{equation*}
and \eqref{e44} is proved. From \eqref{d95} and \eqref{e44}, we deduce that the imaginary part of the resonances $z$ given by Theorem \ref{d8} satisfies
\begin{equation} \label{e47}
\im z = o \Big( \frac {h}{\vert \ln h \vert} \Big) .
\end{equation}
In particular, the two accumulation curves defined in Remark \ref{e2} coincide with the real axis. This is consistent with what is known about the resonances in this geometric setting. Indeed, the case {\rm (C)} can also be seen as a well in the island situation. Then, their imaginary part is exponentially small from Helffer and Sj\"{o}strand \cite[Proposition 9.6]{HeSj86_01}. We obtain here the weaker estimate \eqref{e47} because we only consider the symbol of the operator on the homoclinic trajectories and we do not consider the tunneling effect inside the island.

Thus, the macroscopic behavior of the resonances is very simple: they accumulate on the real axis. The situation is similar to the symmetric case in Example \ref{e6} (even if the two trajectories $\gamma_{1}$ and $\gamma_{2}$ can be different here). At the microscopic level, the distribution of the resonances follows also the properties obtained in Example \ref{e6}. As in \eqref{e16}, the packets of resonances are regularly spaced by $\pi \lambda_{1} h \vert \ln h \vert^{- 1}$ in the limit $\re \sigma \to + \infty$. On the contrary, as $\re \sigma \to - \infty$, the distribution of resonances is given by \eqref{e9} mutatis mutandis and seems to be the reunion of two copies of the resonances generated in case {\rm (B)}. One can also verify that the real part of the resonances coincide, up to lower order terms, to the eigenvalues computed by Colin de Verdi\`ere and Parisse \cite{CoPa94_02} in the case of the (symmetric) double-well potential.

As in {\rm (B)}, the assumption \ref{h5} is not satisfied and Theorem \ref{a1} can not be applied. In the present setting, we have $\S^{0} = \CH_{\rm tang}^{\pm \infty} = \{ - 1 , 1 \}$. Thus, computing as in \eqref{b75}, $\CT_{0} ( \tau )$ can be seen as the $2 \times 2$ matrix
\begin{equation*}
\frac{1}{\sqrt{2 \cosh ( \pi \tau / \lambda_{1} )}} \left( \begin{array}{cc}
e^{- \frac{\pi \tau}{2 \lambda_{1}}} & e^{\frac{\pi \tau}{2 \lambda_{1}}} \\
e^{\frac{\pi \tau}{2 \lambda_{1}}} & e^{- \frac{\pi \tau}{2 \lambda_{1}}}
\end{array} \right) ,
\end{equation*}
which satisfies
\begin{equation*}
\CA_{0} ( \tau ) = \spr ( \CT_{0} ( \tau ) ) = \frac{e^{- \frac{\pi \tau}{2 \lambda_{1}}} + e^{\frac{\pi \tau}{2 \lambda_{1}}}}{\sqrt{2 \cosh ( \pi \tau / \lambda_{1} )}} = \sqrt{1 + \frac{1}{\cosh ( \pi \tau / \lambda_{1} )}} > 1,
\end{equation*}
for all $\tau \in \R$. Note that $\CA_{0} ( \tau )$ exceeds $1$ because the kernel of $\CT_{0}$ is the modulus of the kernel of the quantization operator. Then, Theorem \ref{a2} gives no resonance free region below any interval. Of course, since the imaginary part of the resonances is exponentially small, there is no hope to have a resonance free zone of size $h \vert \ln h \vert^{- 1}$.

\Subsection{Resonances in deeper zones} \label{s43}

Theorem \ref{d8} provides the distribution of resonances up to the line $\im z = - h \sum_{j = 2}^{n} \lambda_{j} / 2$ and we investigate here what happens below this line. To that aim, we first consider formally the quantization rule of Definition \ref{d1}. In the region $\im z < - h \sum_{j = 2}^{n} \lambda_{j} / 2$, we have $\vert h^{S ( z , h ) / \lambda_{1} - 1 / 2} \vert \gg 1$. Thus, it is possible to have a pseudo-resonance only at the complex numbers $z$ such that $0$ is close to $\spe ( \CQ (z , h ) )$. Then, we will first prove that $P$ has no resonance if $0 \notin \spe ( \CQ (z , h ) )$. After that we will consider a general situation where $0 \in \spe ( \CQ (z , h ) )$ for all $z \in \C$ and investigate the distribution of resonance in this case. The proof of all the results stated in the present Section \ref{s43} can be found in Section \ref{s23}.

We note that $0 \in \spe ( \CQ (z , h ) )$ if and only if $\det \CQ ( z , h ) = 0$. Using \eqref{d4}, \eqref{d10} and the multilinearity of the determinant, we can write
\begin{equation} \label{e10}
\big\vert \det \widehat{\CQ} ( \tau , h ) \big\vert = \widehat{q} ( \tau ) \big\vert \det \widehat{\CZ} ( \tau ) \big\vert ,
\end{equation}
where the function
\begin{equation*}
\widehat{q} ( \tau ) : = \Big\vert \Gamma \Big( \frac{1}{2} - i \frac{\tau}{\lambda_{1}} \Big) \Big\vert^{K} ( 2 \pi )^{- \frac{K}{2}} \prod_{k=1}^{K} \frac{\CM^{+}_{k}}{\CM^{-}_{k}} \sqrt{\frac{\vert g_{-}^{k} \vert}{\vert g_{+}^{k} \vert}} ,
\end{equation*}
satisfies $1 \lesssim \widehat{q} ( \tau ) \lesssim 1$ uniformly for $\tau \in [ - C , C ]$ and the coefficients of the $K \times K$ matrix $\widehat{\CZ}$ are given by
\begin{equation}
\widehat{\CZ}_{k , \ell} ( \tau ) : = \big( i \widehat{g}_{+}^{k} \cdot \widehat{g}_{-}^{\ell} \big)^{- \frac{1}{2} + i \frac{\tau}{\lambda_{1}}} .
\end{equation}
Here, $\widehat{g}_{\pm}^{\bullet} : = g_{\pm}^{\bullet} / \vert g_{\pm}^{\bullet} \vert$ denotes the normalized asymptotic directions. From \ref{h4}, the function $\widehat{q} ( \tau )$ and the matrix $\widehat{\CZ} ( \tau )$ are analytic in $[ - C , C ]$. Thus, \eqref{e10} implies that $0 \in \spe ( \widehat{\CQ} ( \tau , h ) )$ if and only if $\det \widehat{\CZ} ( \tau ) = 0$. In particular, the zeros of $\det \widehat{\CQ} ( \tau , h )$ are independent of $h$, but we will see that the multiplicity of $0$ as an eigenvalue of $\widehat{\CQ}$ may depend on $h$. Note also that $\widehat{\CZ}$ depends only on geometric quantities (the normalized asymptotic directions of $\CH$) and not on dynamical quantities (like the constants $\CM^{\pm}_{\bullet}$, the Maslov indices $\nu_{\bullet}$, \ldots).

\subsubsection{Resonance free domains when $0 \notin \spe ( \CQ ( z , h ) )$} \label{s15}

We now prove that, other than the resonances given in Section \ref{s61}, there is no resonance near the points $z$ within a reasonable region such that $0 \notin \spe ( \CQ ( z , h ) )$. We first consider the situation near the line $\im z = - h \sum_{j = 2}^{n} \lambda_{j} / 2$ on which the accumulation curves of Theorem \ref{d8} concentrate. Let $a < b$ and assume that
\begin{hyp} \label{h10}
For all $\tau \in [ a , b ]$, we have $\det \widehat{\CZ} ( \tau ) \neq 0$.
\end{hyp}
In this case, the $K$ eigenvalues of $\widehat{\CQ} ( \tau , h )$ avoid a neighborhood of $0$ and all the pseudo-resonances (given by Definition \ref{d1}) below $[ E_{0} + a h , E_{0} + b h ]$ belong to the set \eqref{d90} for $C$ large enough. Under this hypothesis, we have the following result.

\begin{figure}
\begin{center}
\begin{picture}(0,0)%
\includegraphics{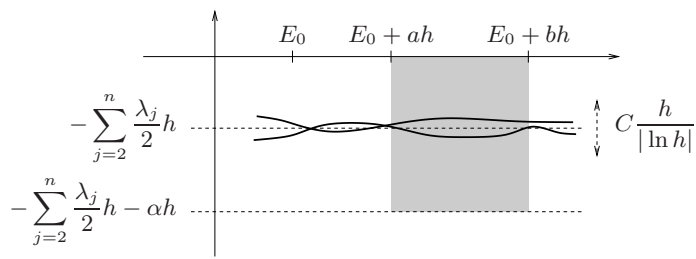}%
\end{picture}%
\setlength{\unitlength}{987sp}%
\begingroup\makeatletter\ifx\SetFigFont\undefined%
\gdef\SetFigFont#1#2#3#4#5{%
  \reset@font\fontsize{#1}{#2pt}%
  \fontfamily{#3}\fontseries{#4}\fontshape{#5}%
  \selectfont}%
\fi\endgroup%
\begin{picture}(12044,6269)(-621,-6008)
\put(3151,-586){\makebox(0,0)[b]{\smash{{\SetFigFont{9}{10.8}{\rmdefault}{\mddefault}{\updefault}$E_{0}$}}}}
\put(9076,-586){\makebox(0,0)[b]{\smash{{\SetFigFont{9}{10.8}{\rmdefault}{\mddefault}{\updefault}$E_{0} + b h$}}}}
\put(5626,-586){\makebox(0,0)[b]{\smash{{\SetFigFont{9}{10.8}{\rmdefault}{\mddefault}{\updefault}$E_{0} + a h$}}}}
\put(241,-2836){\makebox(0,0)[rb]{\smash{{\SetFigFont{9}{10.8}{\rmdefault}{\mddefault}{\updefault}$\ds - \sum_{j=2}^{n} \frac{\lambda_{j}}{2} h$}}}}
\put(241,-4936){\makebox(0,0)[rb]{\smash{{\SetFigFont{9}{10.8}{\rmdefault}{\mddefault}{\updefault}$\ds - \sum_{j=2}^{n} \frac{\lambda_{j}}{2} h - \alpha h$}}}}
\put(11251,-2836){\makebox(0,0)[lb]{\smash{{\SetFigFont{9}{10.8}{\rmdefault}{\mddefault}{\updefault}$\ds C \frac{h}{\vert \ln h \vert}$}}}}
\end{picture}%
\end{center}
\caption{The region studied in Theorem \ref{e11}.} \label{f27}
\end{figure}

\begin{theorem}[Resonance free domain below the first accumulation curves]\sl \label{e11}
Assume \ref{h1}--\ref{h4}, \ref{h8} and \ref{h10}. There exists $\alpha > 0$ such that, for all $\delta > 0$, the conclusions of Theorem \ref{d8} hold true with \eqref{d90} replaced by
\begin{equation} \label{e12}
E_{0} + [ a h , b  h ] + i \Big[ - \sum_{j = 2}^{n} \frac{\lambda_{j}}{2} h - \alpha h , h \Big] \setminus \big( \Gamma (h) + B ( 0 , \delta h ) \big) .
\end{equation}
\end{theorem}

In other words, comparing with Theorem \ref{d8}, this result just adds that $P$ has no resonance and that its truncated resolvent has a polynomial estimate in the domain
\begin{equation} \label{e84}
E_{0} + [ a h , b  h ] - i \sum_{j = 2}^{n} \frac{\lambda_{j}}{2} h + i \Big[ - \alpha h , - C \frac{h}{\vert \ln h \vert} \Big] \setminus \big( \Gamma (h) + B ( 0 , \delta h ) \big) ,
\end{equation}
for $C$ large enough and $h$ small enough. The setting is illustrated in Figure \ref{f27}. For example, the assumption \ref{h10} is satisfied for a single trajectory (see Corollary \ref{e4}), in Example \ref{e6}, in Example \ref{e29} and in the three possible situations in dimension $1$ (see Section \ref{s19}).

The natural question is then what is the optimal depth of this resonance free domain (i.e. $\alpha$). It is actually limited for two reasons: that $0 \in \spe ( \CQ (z , h ) )$ and because of the second set of accumulation curves which concentrates near the line $\im z = - h \sum_{j = 2}^{n} \lambda_{j} / 2 - h \lambda_{1}$ (see Section \ref{s83}). Thus we have the following result in domains that avoid these limitations. For $\sigma \in \C$, let $\widetilde{\CZ} ( \sigma )$ be the $K \times K$ matrix given by
\begin{equation} \label{g82}
\widetilde{\CZ}_{k , \ell} ( \sigma ) : = \big( i \widehat{g}_{+}^{k} \cdot \widehat{g}_{-}^{\ell} \big)^{i \frac{\sigma}{\lambda_{1}} - \sum_{j = 1}^{n} \frac{\lambda_{j}}{2 \lambda_{1}}} .
\end{equation}
In particular, as in \eqref{e18}, we have
\begin{equation*}
\widehat{\CZ} ( \tau ) = \widetilde{\CZ} \Big( \tau - i \sum_{j = 2}^{n} \frac{\lambda_{j}}{2} \Big) .
\end{equation*}
Let $\Omega$ be a compact subset of $\R - i \sum_{j = 2}^{n} \lambda_{j} / 2 + i ] - \lambda_{1} , 0 [$. We assume that
\begin{hyp} \label{h11}
For all $\sigma \in \Omega$, we have $\det \widetilde{\CZ} ( \sigma ) \neq 0$.
\end{hyp}
This assumption guaranties that $\CQ ( z , h )$ is invertible for $z \in E_{0} + h \Omega$. Note that $\Omega$ is indeed between the two first sets of accumulation curves.

\begin{proposition}\sl \label{e15}
Assume \ref{h1}--\ref{h4}, \ref{h8} and \ref{h11}. For all $\delta > 0$, $P$ has no resonance in $E_{0} + h \Omega \setminus ( \Gamma (h) + B ( 0 , \delta h ) )$ for $h$ small enough. Moreover, for all $\chi \in C^{\infty}_{0} ( \R^{n} )$, there exists $M > 0$ such that
\begin{equation*}
\big\Vert \chi ( P -z )^{-1} \chi \big\Vert \lesssim h^{- M} ,
\end{equation*}
uniformly for $h$ small enough and $z \in E_{0} + h \Omega \setminus ( \Gamma (h) + B ( 0 , \delta h ) )$.
\end{proposition}

In the case of a single trajectory (i.e. under the assumption \ref{h9}), this result implies that the conclusions of Theorem \ref{d8} hold true in any region above the optimal line $\im z = - h \sum_{j = 2}^{n} \lambda_{j} / 2 - h \lambda_{1}$ since the scalar function $\widetilde{\CZ} ( \sigma )$ never vanishes on $\C$. On the contrary, if $\CH$ has more than one trajectory, $\det \widetilde{\CZ} ( \sigma ) = 0$ at the point $\sigma_{0} = - i \sum_{j = 1}^{n} \lambda_{j} / 2$. Indeed, all the coefficients of $\widetilde{\CZ} ( \sigma_{0} )$ are equal to $1$. Note that $\sigma_{0} \in \Gamma_{0} (h)$ and that $E_{0} + \sigma_{0} h$ corresponds to (the leading term in the asymptotic of) the first resonance generated by the barrier-top without homoclinic trajectories (see \eqref{n1}).

\subsubsection{A general setting where $0 \in \spe ( \CQ ( z , h ) )$} \label{s16}

We now consider situations where $\det \widehat{\CZ}$ takes the value $0$. At a first sight, this case can seem to be artificial. But, we will see that $\det \widehat{\CZ}$ vanishes identically for a large class of operators.

Assume $n \geq 2$ and $\lambda_{1} < \lambda_{2}$. We denote by $e_{1} : = ( 1 , 0 , \ldots , 0 ) \in \R^{n}$ the first vector of the canonical basis. From \eqref{d2} and \ref{h4}, the non-zero asymptotic directions $g_{\pm}^{\bullet}$ are collinear to $e_{1}$ in this situation. In the sequel, we will say that two such vectors are ``at the same side of $0$'' if they are positively related. This definition is consistent with the one of Section \ref{s14}. Thus, if two incoming (or outgoing) asymptotic directions (say $g_{-}^{1}$ and $g_{-}^{2}$) are at the same side of $0$, the two first columns of $\widehat{\CZ} ( \tau )$ are the same for all $\tau \in \R$ and its determinant necessarily vanishes identically. More generally, $\det \widetilde{\CZ} ( \sigma ) = 0$ for all $\sigma  \in \C$.

Moreover, if all the asymptotic directions are at the same side of $0$, one can verify that $0$ is an eigenvalue of multiplicity $K - 1$ of $\widehat{\CQ} ( \tau )$ and that the last eigenvalue (which can also be $0$) is given by
\begin{align}
\mu ( \tau , h ) : ={}& \tr \big( \widehat{\CQ} ( \tau , h ) \big)  \nonumber \\
={}& \sum_{k = 1}^{K} e^{i A_{k} / h} \Gamma \Big( \frac{1}{2} - i \frac{\tau}{\lambda_{1}} \Big) \sqrt{\frac{\lambda_{1}}{2 \pi}} \frac{\CM_{k}^{+}}{\CM_{k}^{-}} e^{- \frac{\pi}{2} ( \nu_{k} + 1 ) i} \big\vert g_{-}^{k} \big\vert \big( \lambda_{1} \vert g_{+}^{k} \vert \vert g_{-}^{k} \vert \big)^{- \frac{1}{2} + i \frac{\tau}{\lambda_{1}}} e^{- \frac{\pi \tau}{2 \lambda_{1}}} . \label{e13}
\end{align}
Under these hypotheses, Theorem \ref{d8} provides at most one accumulation curve in the set \eqref{d90}. More precisely, it corresponds to the points
\begin{equation} \label{e34}
z_{q} ( \tau ) = E_{0} + 2 q \pi \lambda_{1} \frac{h}{\vert \ln h \vert} - i h \sum_{j = 2}^{n} \frac{\lambda_{j}}{2} + i \ln ( \mu ( \tau , h ) ) \lambda_{1} \frac{h}{\vert \ln h \vert} ,
\end{equation}
given by Proposition \ref{d9}. Note also that the vibration phenomenon described in Remark \ref{e5} takes place explicitly here as soon as two of the actions $A_{\bullet}$ are different.

\begin{proposition}\sl \label{e19}
Assume \ref{h1}--\ref{h4}, \ref{h8}, $\lambda_{1} < \lambda_{2}$, that all the asymptotic directions are at the same side of $0$ and that $\mu ( \tau , h )$ avoids a neighborhood of $0$ for all $\tau \in [ a , b ]$ and $h$ small enough. Then, there exists $\alpha > 0$ such that, for all $\delta > 0$, the conclusions of Theorem \ref{d8} hold true with \eqref{d90} replaced by
\begin{equation} \label{e20}
E_{0} + [ a h , b  h ] + i \Big[ - \sum_{j = 2}^{n} \frac{\lambda_{j}}{2} h - \alpha h , h \Big] \setminus \big( \Gamma (h) + B ( 0 , \delta h ) \big) .
\end{equation}
\end{proposition}

In other words, the previous result provides a situation where there actually is only one accumulation curve in \eqref{e20} (given by Theorem \ref{d8}) and a resonance free region below of size $h$. On the other hand, $\CH$ consists of $K$ different homoclinic trajectories. Thus, there is in general no equality between the number of accumulation curves of resonances and the number $K$ of homoclinic trajectories. To complete the exposition, we give an example of operator satisfying the assumptions of Proposition \ref{e19} with more than one homoclinic trajectories.

\begin{example}\rm \label{e22}
We adapt Example \ref{e29} in order that the new operator satisfies the hypotheses of Proposition \ref{e19} and $K \geq 2$. Note that it will not be of the form \eqref{a5} but will enter in the setting of Remark \ref{c13}.

\begin{figure}
\begin{center}
\begin{picture}(0,0)%
\includegraphics{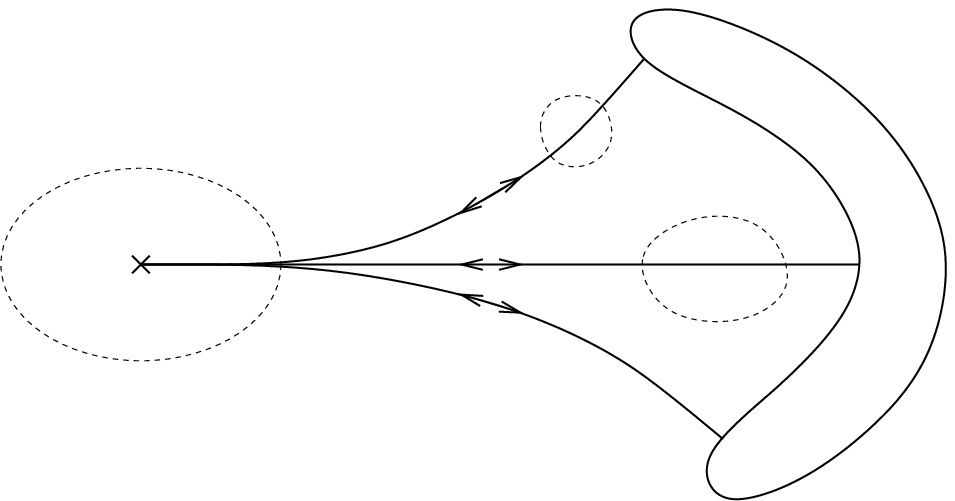}%
\end{picture}%
\setlength{\unitlength}{1105sp}%
\begingroup\makeatletter\ifx\SetFigFont\undefined%
\gdef\SetFigFont#1#2#3#4#5{%
  \reset@font\fontsize{#1}{#2pt}%
  \fontfamily{#3}\fontseries{#4}\fontshape{#5}%
  \selectfont}%
\fi\endgroup%
\begin{picture}(16246,8461)(-10814,-9217)
\put(-8399,-4786){\makebox(0,0)[b]{\smash{{\SetFigFont{9}{10.8}{\rmdefault}{\mddefault}{\updefault}$0$}}}}
\put(4426,-5236){\makebox(0,0)[b]{\smash{{\SetFigFont{9}{10.8}{\rmdefault}{\mddefault}{\updefault}$x_{3}$}}}}
\put(601,-1486){\makebox(0,0)[b]{\smash{{\SetFigFont{9}{10.8}{\rmdefault}{\mddefault}{\updefault}$x_{1}$}}}}
\put(2326,-2236){\makebox(0,0)[b]{\smash{{\SetFigFont{9}{10.8}{\rmdefault}{\mddefault}{\updefault}$\CO$}}}}
\put(-8399,-3061){\makebox(0,0)[b]{\smash{{\SetFigFont{9}{10.8}{\rmdefault}{\mddefault}{\updefault}$\supp V$}}}}
\put(1951,-8536){\makebox(0,0)[b]{\smash{{\SetFigFont{9}{10.8}{\rmdefault}{\mddefault}{\updefault}$x_{2}$}}}}
\put(-1049,-1936){\makebox(0,0)[b]{\smash{{\SetFigFont{9}{10.8}{\rmdefault}{\mddefault}{\updefault}$\supp \widetilde{W}$}}}}
\put(-2474,-3436){\makebox(0,0)[b]{\smash{{\SetFigFont{9}{10.8}{\rmdefault}{\mddefault}{\updefault}$\pi_{x} ( \gamma_{1} )$}}}}
\put(-2474,-4786){\makebox(0,0)[b]{\smash{{\SetFigFont{9}{10.8}{\rmdefault}{\mddefault}{\updefault}$\pi_{x} ( \gamma_{3} )$}}}}
\put(-2474,-6511){\makebox(0,0)[b]{\smash{{\SetFigFont{9}{10.8}{\rmdefault}{\mddefault}{\updefault}$\pi_{x} ( \gamma_{2} )$}}}}
\put(1501,-3961){\makebox(0,0)[b]{\smash{{\SetFigFont{9}{10.8}{\rmdefault}{\mddefault}{\updefault}$\supp W$}}}}
\end{picture}%
\end{center}
\bigskip \bigskip
\begin{center}
\begin{picture}(0,0)%
\includegraphics{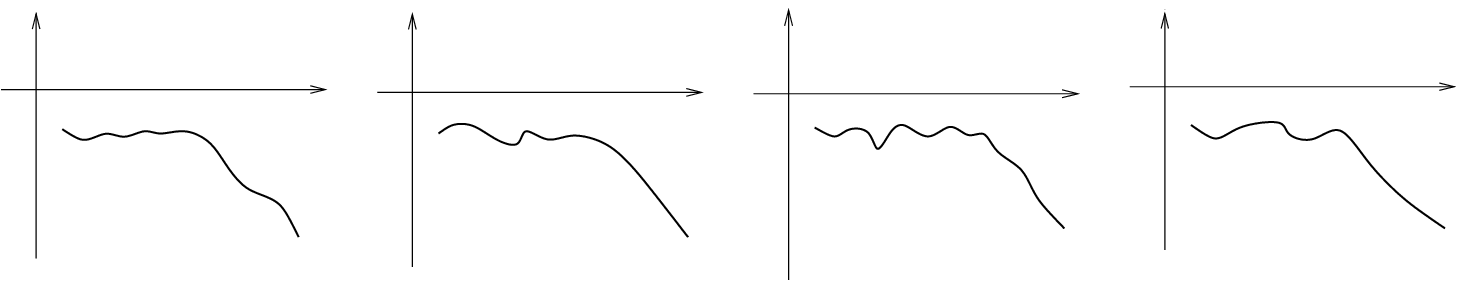}%
\end{picture}%
\setlength{\unitlength}{1105sp}%
\begingroup\makeatletter\ifx\SetFigFont\undefined%
\gdef\SetFigFont#1#2#3#4#5{%
  \reset@font\fontsize{#1}{#2pt}%
  \fontfamily{#3}\fontseries{#4}\fontshape{#5}%
  \selectfont}%
\fi\endgroup%
\begin{picture}(25019,4735)(-621,-5091)
\put(2701,-957){\makebox(0,0)[b]{\smash{{\SetFigFont{9}{10.8}{\rmdefault}{\mddefault}{\updefault}$h = h_{1}$}}}}
\put(15601,-941){\makebox(0,0)[b]{\smash{{\SetFigFont{9}{10.8}{\rmdefault}{\mddefault}{\updefault}$h = h_{3}$}}}}
\put(22051,-936){\makebox(0,0)[b]{\smash{{\SetFigFont{9}{10.8}{\rmdefault}{\mddefault}{\updefault}$h = h_{4}$}}}}
\put(9151,-977){\makebox(0,0)[b]{\smash{{\SetFigFont{9}{10.8}{\rmdefault}{\mddefault}{\updefault}$h = h_{2}$}}}}
\end{picture}%
\end{center}
\caption{The geometry of Example \ref{e22} and the corresponding accumulation curve in the case $K = 2 , 3$ without symmetry.} \label{f29}
\end{figure}

In dimension $n = 2$, we consider
\begin{equation*}
P = - h^{2} \Delta_{\R^{2} \setminus \CO} + V (x) ,
\end{equation*}
with Dirichlet condition at the boundary of the obstacle $\CO$. The potential $V$ is of the form
\begin{equation*}
V (x) = V_{1} ( x_{1} ) V_{2} ( x_{2} ) ,
\end{equation*}
where the functions $V_{\bullet} \in C^{\infty}_{0} ( \R )$ are simple barriers (like the case {\rm (A)} in Figure \ref{f6}) with
\begin{equation*}
V_{1} ( x_{1} ) = E_{0} - \frac{\lambda_{1}^{2}}{2} x_{1}^{2} + \CO ( x_{1}^{3} )  \qquad \text{and} \qquad   V_{2} ( x_{2} ) = 1 - \frac{\lambda_{2}^{2}}{2 E_{0}} x_{2}^{2} + \CO ( x_{2}^{3} ) ,
\end{equation*}
near $0$ with $\lambda_{1} < \lambda_{2}$. In particular, \ref{h2} holds true.

On the other hand, $\CO$ is a small enough obstacle near $( a , 0 ) \in \R^{2}$ with $a > 0$ fixed outside of the support of $V_{1}$. Using Lemma \ref{a43}, one can verify that no Hamiltonian trajectory of energy $E_{0}$ can start from the boundary of $\CO$, touch the support of $V$ and then come back to $\CO$. Thus, if $\CO$ is as illustrated in Figure \ref{f29}, the trapped set of $P$ at energy $E_{0}$ verifies \ref{h3} and $\CH$ consists of bicharacteristic curves touching $\partial \CO$ only one time. Furthermore, their asymptotic directions are on the same side of $0$. Moreover, as in Example \ref{c15}, \ref{h4} is automatically satisfied and \ref{h8} holds true under an assumption on the curvature of $\partial \CO$ at the points of $\partial \CO \cap \pi_{x} ( \CH )$. Figure \ref{f29} corresponds to $K = 3$. But, adding to $P$ a potential $i h \vert \ln h \vert W (x)$ with $0 \leq W \in C^{\infty}_{0} ( \R^{2} )$ allows to ``remove'' a trajectory and to consider the case $K = 2$.

If now we add a potential of the form $\nu h \widetilde{W} (x)$ with $\nu \in \C$ and $0 \leq \widetilde{W} \in C^{\infty}_{0} ( \R^{2} )$ supported on the trajectory $\gamma_{1}$, then the usual transport equations show that the contribution of this trajectory in \eqref{e13} is multiplied by the constant $e^{- i \nu \int_{\R} \widetilde{W} ( x_{1} ( t ) ) \, d t}$. In particular, we can adapt $\nu$ so that this constant is any number of $\C \setminus \{ 0 \}$. Thus, it is easy to guaranty that $\mu ( \tau , h )$ avoids a neighborhood of $0$ for $\tau$ in a fixed interval.

Consider now the case $K = 2$. Combining \eqref{e36}, \eqref{e13} and \eqref{e34}, the unique accumulation curve given by \eqref{d93} satisfies the following asymptotic
\begin{align}
\im \sigma = -\frac{\lambda_{2}}{2} + \frac{\lambda_{1}}{\vert \ln h \vert} \ln \bigg\vert \sum_{k = 1}^{2} e^{i A_{k} / h} \frac{\CM_{k}^{+}}{\CM_{k}^{-}} & \sqrt{\frac{\vert g_{-}^{k} \vert}{\vert g_{+}^{k} \vert}} e^{- i \nu_{k} \frac{\pi}{2}} \big( \vert g_{+}^{k} \vert \vert g_{-}^{k} \vert \big)^{i \frac{\re \sigma}{\lambda_{1}}} \bigg\vert  \nonumber \\
&+ \frac{1}{\vert \ln h \vert}
\left\{ \begin{aligned}
&o (1) &&\text{ as } \re \sigma \to - \infty , \\
&- \pi \re \sigma + o (1) &&\text{ as } \re \sigma \to + \infty , \\
\end{aligned} \right.
\end{align}
The behavior as $\re \sigma \to + \infty$ is in agreement with Remark \ref{e23}. Moreover, in the limit $\re \sigma \to  - \infty$, the accumulation curve is asymptotically periodic with respect to $\re \sigma$ since one can factorized by one term in the modulus of the sum. Such a phenomenon has already been observed in Example \ref{e29}.
\end{example}

Eventually, coming back to the general setting of Proposition \ref{e19}, assume further that all the trajectories have the same geometric quantities ($A_{\bullet}$, $g_{\pm}^{\bullet}$, $\nu_{\bullet}$, \ldots). In this setting, we remove the subscript $k$. The case $K = 2$ can be realized by taking $\CO$ symmetric in Example \ref{e22}. In this situation, the accumulation curve is given by
\begin{equation} \label{e35}
\im \sigma = - \frac{\lambda_{2}}{2} + \frac{C ( K )}{\vert \ln h \vert} - \frac{\lambda_{1}}{2 \vert \ln h \vert} \ln \big( e^{\frac{2 \pi}{\lambda_{1}} \re \sigma} + 1 \big) ,
\end{equation}
where
\begin{equation}
C ( K ) = \lambda_{1} \ln \bigg( K \frac{\CM^{+}}{\CM^{-}} \sqrt{\frac{\vert g_{-} \vert}{\vert g_{+} \vert}} \bigg) .
\end{equation}
Comparing \eqref{e3} with \eqref{e35}, we note that all happen as if there was only one homoclinic trajectory except that $C ( K ) = C_{0} + \lambda_{1} \ln K$. In particular, since $C ( K ) > C ( \widetilde{K} )$ when $K > \widetilde{K}$, the accumulation curve is closer to the real axis as $K$ increases. That several identical trajectories are ``more trapping'' than one trajectory can explain this property.

\subsubsection{An example of change of multiplicity} \label{s75}

In the previous parts, we have only considered situations where $0$ has a fixed multiplicity as an eigenvalue of $\widehat{\CQ} ( \tau , h )$: this multiplicity is $0$ in Section \ref{s15} and $K - 1$ in Section \ref{s16}. So far, we have not considered resonances (or accumulation curves) leaving the set \eqref{d90} where the conclusions of Theorem \ref{d8} are valid. In this case, the relevant question is of course how far these resonances go down in the complex plane. In particular, we can wonder if they reach to the second set of accumulation curves described in the Section \ref{s83} below.

We study here a situation where the multiplicity of $0$ as an eigenvalue of $\widehat{\CQ} ( \tau , h )$ changes and investigate the behavior of the resonances in deeper domains. For that, we come back to the setting of Section \ref{s16}. In this case, $0$ is an eigenvalue of multiplicity (at least) $K - 1$ and the last one is given by \eqref{e13}. Thus, the points $( \tau , h)$ at which $\mu ( \tau , h )$ vanishes correspond to the points where the multiplicity of $0$ as an eigenvalue of $\widehat{\CQ}$ changes. By definition, $\mu ( \tau , h)$ can be written
\begin{equation} \label{j61}
\mu ( \tau , h ) = \Gamma \Big( \frac{1}{2} - i \frac{\tau}{\lambda_{1}} \Big) e^{- \frac{\pi \tau}{2 \lambda_{1}}} \sum_{k = 1}^{K} e^{i A_{k} / h} B_{k} e^{i T_{k} \tau} ,
\end{equation}
for some $B_{k} \in \C \setminus \{ 0 \}$ and $T_{k} = \ln ( \lambda_{1} \vert g_{+}^{k} \vert \vert g_{-}^{k} \vert ) / \lambda_{1} \in \R$. The following example shows that, even in a simple situation, the changes of multiplicity can be rather complicated.

\begin{figure}
\begin{center}
\begin{picture}(0,0)%
\includegraphics{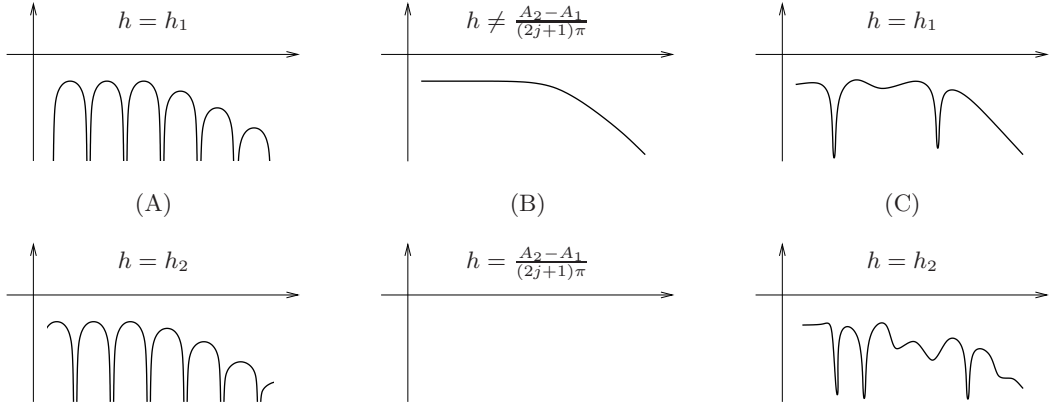}%
\end{picture}%
\setlength{\unitlength}{1105sp}%
\begingroup\makeatletter\ifx\SetFigFont\undefined%
\gdef\SetFigFont#1#2#3#4#5{%
  \reset@font\fontsize{#1}{#2pt}%
  \fontfamily{#3}\fontseries{#4}\fontshape{#5}%
  \selectfont}%
\fi\endgroup%
\begin{picture}(23444,9055)(-1821,-9394)
\put(1501,-957){\makebox(0,0)[b]{\smash{{\SetFigFont{9}{10.8}{\rmdefault}{\mddefault}{\updefault}$h = h_{1}$}}}}
\put(1501,-6357){\makebox(0,0)[b]{\smash{{\SetFigFont{9}{10.8}{\rmdefault}{\mddefault}{\updefault}$h = h_{2}$}}}}
\put(18301,-957){\makebox(0,0)[b]{\smash{{\SetFigFont{9}{10.8}{\rmdefault}{\mddefault}{\updefault}$h = h_{1}$}}}}
\put(18301,-6357){\makebox(0,0)[b]{\smash{{\SetFigFont{9}{10.8}{\rmdefault}{\mddefault}{\updefault}$h = h_{2}$}}}}
\put(9901,-957){\makebox(0,0)[b]{\smash{{\SetFigFont{9}{10.8}{\rmdefault}{\mddefault}{\updefault}$h \neq \frac{A_{2} - A_{1}}{( 2 j + 1 ) \pi}$}}}}
\put(9901,-6357){\makebox(0,0)[b]{\smash{{\SetFigFont{9}{10.8}{\rmdefault}{\mddefault}{\updefault}$h = \frac{A_{2} - A_{1}}{( 2 j + 1 ) \pi}$}}}}
\put(18301,-5086){\makebox(0,0)[b]{\smash{{\SetFigFont{9}{10.8}{\rmdefault}{\mddefault}{\updefault}${\rm (C)}$}}}}
\put(1501,-5086){\makebox(0,0)[b]{\smash{{\SetFigFont{9}{10.8}{\rmdefault}{\mddefault}{\updefault}${\rm (A)}$}}}}
\put(9901,-5086){\makebox(0,0)[b]{\smash{{\SetFigFont{9}{10.8}{\rmdefault}{\mddefault}{\updefault}${\rm (B)}$}}}}
\end{picture}%
\end{center}
\caption{The accumulation curve in the three cases of Example \ref{e48} for different values of $h$}\label{f30}
\end{figure}

\begin{example}\rm \label{e48}
Coming back to the setting of Example \ref{e22} and adding possibly some potentials supported on the characteristic trajectories $\pi_{x} ( \gamma_{\bullet} )$, we can realize the three different situations described below. The corresponding distribution of resonances is illustrated in Figure \ref{f30}. There may also exist more complicated situations.

In case {\rm (A)}, we assume $K = 2$, $B_{1} = B_{2}$ and $T_{1} \neq T_{2}$. In this situation, $\mu$ vanishes at the points $( \tau , h )$ such that
\begin{equation} \label{e50}
\tau = \frac{( 2 j + 1 ) \pi}{T_{1} - T_{2}} + \frac{A_{2} - A_{1}}{( T_{1} - T_{2} ) h} ,
\end{equation}
with $j \in \Z$. In particular, these points form a discrete set for $h$ fixed.

In case {\rm (B)}, we assume $K = 2$, $B_{1} = B_{2}$, $T_{1} = T_{2}$ and $A_{1} \neq A_{2}$. In this situation, the vanishing points of $\mu$ are given by
\begin{equation} \label{e49}
h = \frac{A_{2} - A_{1}}{( 2 j + 1 ) \pi} ,
\end{equation}
with $j \in \Z$. Then, for the values of $h$ satisfying \eqref{e49}, $\tau \mapsto \mu ( \tau , h )$ vanishes identically and there is no resonance (and a fortiori no accumulation curve) in the set \eqref{d90}. This shows that Theorem \ref{d8} may sometimes provide not a single resonance. In the present setting, one can even show that $P$ has no resonance (outside $\Gamma (h) + B ( 0 , \delta h)$) above the line $\im z = - h \sum_{j = 2}^{n} \lambda_{j} / 2 - h \alpha$, for some $\alpha > 0$. This is done in Lemma \ref{j59}. On the contrary, if $h$ avoid the exceptional values \eqref{e49}, one could show that the conclusions of Proposition \ref{e19} hold true.

In case {\rm (C)}, we assume $K = 3$, $2 \max_{k = 1 , 2 , 3} \vert B_{k} \vert \leq \vert B_{1} \vert + \vert B_{2} \vert + \vert B_{3} \vert$ (in other words, there exists a unique (modulo symmetries) triangle whose sides are of length $\vert B_{\bullet} \vert$), $T_{1} / A_{1} = T_{2} / A_{2} = T_{3} / A_{3}$ and that the $A_{\bullet}$ are $\Z$-independent. Then, for any non-empty interval $I \subset \R$, the function $\mu ( \tau , h )$ does not vanish for $\tau \in I$ and $h$ small enough. Nevertheless, there exists a sequence $( \tau_{j} , h_{j} )_{j \in \N} \in I \times ] 0 , 1 ]$ such that $h_{j} \to 0$ and $\mu ( \tau_{j} , h_{j} ) \to 0$ as $j \to + \infty$. Thus, the change of multiplicity never happens effectively but holds true in the asymptotic regime $h \to 0$.
\end{example}

We will not treat all the situations described in the previous example. Instead, we will concentrate our study on a significant setting. Note that $\tau \mapsto \mu ( \tau , h )$ is analytic on $\R$. In addition to the previous hypotheses, we assume
\begin{hyp} \label{h12}
The eigenvalue $\mu ( \tau , h )$ does not depend on $h$ and there exists $\tau_{0} \in \R$, $\ell \in \N \setminus \{ 0 \}$ and $\alpha \in \C \setminus \{ 0 \}$ such that
\begin{equation*}
\mu ( \tau ) = \alpha ( \tau - \tau_{0} )^{\ell} + \CO \big( ( \tau - \tau_{0} )^{\ell + 1} \big) \text{ near } \tau_{0} .
\end{equation*}
\end{hyp}
This assumption is satisfied in Example \ref{e48} {\rm (A)} if $A_{1} = A_{2}$ and $h$ is restricted to an appropriate sequence which goes to $0$.

Under the assumption \ref{h12}, we still define the pseudo-resonances as in Definition \ref{d1}. Concretely, it means that $z$ is a pseudo-resonance if and only if
\begin{equation*}
h^{S ( z , h ) / \lambda_{1} - 1 / 2} \mu \Big( \frac{z - E_{0}}{h} + i \sum_{j = 2}^{n} \frac{\lambda_{j}}{2} \Big) = 1 .
\end{equation*}
The asymptotic of the pseudo-resonances away from $E_{0} + \tau_{0} h$ (more precisely, such that $\vert \re z - ( E_{0} + \tau_{0} h ) \vert \geq \varepsilon h$ with $\varepsilon > 0$) is again given by Proposition \ref{d9}. But the pseudo-resonances with real part close to $E_{0} + \tau_{0} h$ satisfy a different asymptotic given by Lemma \ref{e52} below. In the sequel, the $q$-th branch of the Lambert function $W$ (i.e. the multivalued inverse of the complex function $x \mapsto x e^{x}$) will be denoted by $W_{q}$. We refer to Corless, Gonnet, Hare, Jeffrey and Knuth \cite{CoGoHaJeKn96_01} (see also \cite{Wr59_01,BeCo63_01,BuSi73_01}) for the precise definition and properties on the Lambert function. In particular, the asymptotic behavior of $W_{q} ( A )$ for $A \to \infty$ can be found in \cite[(4.20)]{CoGoHaJeKn96_01} and some estimates are stated in Lemma \ref{e55} below.

\begin{lemma}\sl \label{e52}
Assume \ref{h1}--\ref{h4} and \ref{h12}, let $\alpha > 0$ be small enough and let $\varepsilon (h)$ be a function which goes to $0$ as $h \to 0$. Then, the pseudo-resonances $z$ in
\begin{equation} \label{e53}
E_{0} + h \tau_{0} + h \varepsilon (h) [ - 1 , 1 ] + i h \Big[ - \sum_{j = 2}^{n} \frac{\lambda_{j}}{2} - \alpha , 1 \Big] ,
\end{equation}
satisfy uniformly
\begin{equation} \label{e54}
z = z_{q , \beta} + o \Big( \frac{h}{\vert \ln h \vert} \Big) ,
\end{equation}
with
\begin{equation} \label{e65}
z_{q , \beta} = E_{0} + h \tau_{0} - i h \sum_{j = 2}^{n} \frac{\lambda_{j}}{2} - i \ell \lambda_{1} W_{q} \Big( i \frac{\vert \ln h \vert}{\ell \lambda_{1} \alpha^{1 / \ell}}e^{- i \frac{\tau_{0}}{\ell \lambda_{1}} \vert \ln h \vert} e^{i 2 \pi \beta / \ell} \Big) \frac{h}{\vert \ln h \vert} ,
\end{equation}
for some $q \in \Z$ and $\beta \in \{ 0 , \ldots , \ell - 1 \}$. On the other hand, for each $q \in \Z$ and $\beta \in \{ 0 , \ldots , \ell - 1 \}$ such that $z_{q , \beta} ( \tau )$ belongs to \eqref{e53}, there exists a pseudo-resonances $z$ satisfying \eqref{e54} uniformly with respect to $q , \beta$.
\end{lemma}

Eventually, as in Theorem \ref{d8}, the resonances verify the following asymptotic.

\begin{proposition}\sl \label{e56}
Assume \ref{h1}--\ref{h4} and \ref{h12}, let $\alpha , \delta > 0$ be small enough and let $\varepsilon (h)$ be a function which goes to $0$ as $h \to 0$. In the domain
\begin{equation} \label{e70}
E_{0} + h \tau_{0} + h \varepsilon (h) [ - 1 , 1 ] + i h \Big[ - \sum_{j = 2}^{n} \frac{\lambda_{j}}{2} - \alpha , 1 \Big]  \setminus \big( \Gamma (h) + B ( 0 , \delta h ) \big),
\end{equation}
we have
\begin{equation*}
\dist \big( \res (P) , \res_{0} (P) \big) = o \Big( \frac {h}{\vert \ln h \vert} \Big) ,
\end{equation*}
as $h$ goes to $0$. Moreover, the truncated resolvent of $P$ satisfies a polynomial estimate at distance $h \vert \ln h \vert^{- 1}$ of $\res_{0} ( P )$ in \eqref{e70}.
\end{proposition}

\begin{figure}
\begin{center}
\begin{picture}(0,0)%
\includegraphics{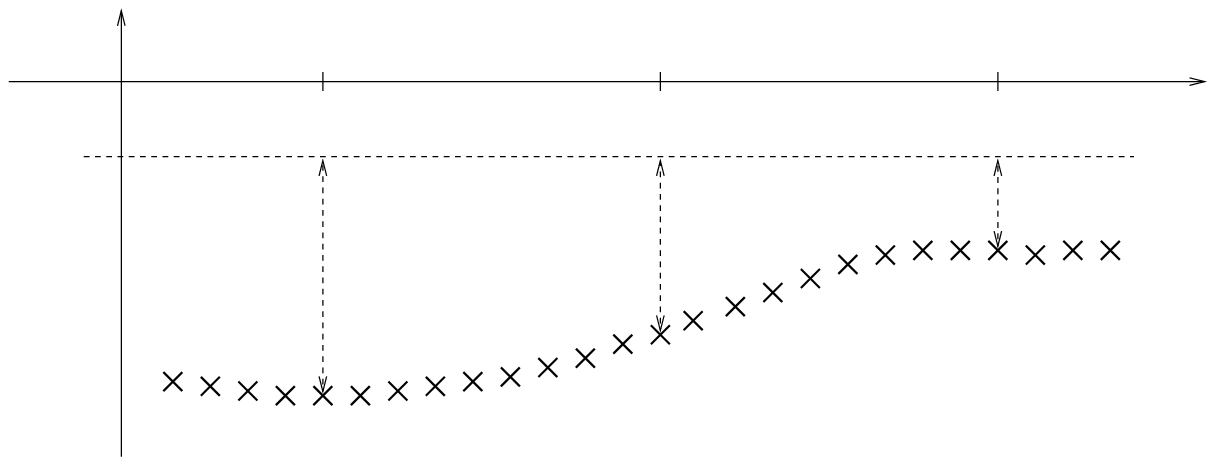}%
\end{picture}%
\setlength{\unitlength}{1184sp}%
\begingroup\makeatletter\ifx\SetFigFont\undefined%
\gdef\SetFigFont#1#2#3#4#5{%
  \reset@font\fontsize{#1}{#2pt}%
  \fontfamily{#3}\fontseries{#4}\fontshape{#5}%
  \selectfont}%
\fi\endgroup%
\begin{picture}(20737,7244)(-2114,-6983)
\put(-2099,-2236){\makebox(0,0)[lb]{\smash{{\SetFigFont{9}{10.8}{\rmdefault}{\mddefault}{\updefault}$\ds - \frac{1}{2} \sum_{j = 2}^{n} \lambda_{j} h$}}}}
\put(4426,-511){\makebox(0,0)[b]{\smash{{\SetFigFont{9}{10.8}{\rmdefault}{\mddefault}{\updefault}$E_{0}$}}}}
\put(10126,-3136){\makebox(0,0)[lb]{\smash{{\SetFigFont{9}{10.8}{\rmdefault}{\mddefault}{\updefault}$\frac{h}{\vert \ln h \vert} \ln \Big( \frac{\vert \ln h \vert}{q} \Big)$}}}}
\put(4726,-4186){\makebox(0,0)[lb]{\smash{{\SetFigFont{9}{10.8}{\rmdefault}{\mddefault}{\updefault}$\frac{h \ln \vert \ln h \vert}{\vert \ln h \vert}$}}}}
\put(9826,-511){\makebox(0,0)[b]{\smash{{\SetFigFont{9}{10.8}{\rmdefault}{\mddefault}{\updefault}$E_{0} + h \tau_{0} + 2 q \pi \lambda_{1} \frac{h}{\vert \ln h \vert}$}}}}
\put(15526,-2911){\makebox(0,0)[lb]{\smash{{\SetFigFont{9}{10.8}{\rmdefault}{\mddefault}{\updefault}$\frac{h}{\vert \ln h \vert}$}}}}
\put(15226,-511){\makebox(0,0)[b]{\smash{{\SetFigFont{9}{10.8}{\rmdefault}{\mddefault}{\updefault}$E_{0}+ h \tau_{0} + h$}}}}
\end{picture}%
\end{center}
\caption{The resonances of Proposition \ref{e56} in the case $\ell = 1$}\label{f31}
\end{figure}

The proof of these results can be found in Section \ref{s76}. Note that Proposition \ref{e56} is interesting only when $E_{0} + h \tau_{0} - i h \sum_{j = 2}^{n} \lambda_{j} / 2 \notin \Gamma (h) + B ( 0 , \delta h )$. This can be achieved as explained in Example \ref{e48}. In the rest of this part, we describe with more details the asymptotic behavior of the (pseudo-)resonances and assume for simplicity that $\ell = 1$. The setting is illustrated in Figure \ref{f31}. Using the rough estimate $\im W_{q} \in 2 q \pi + [ - 2 \pi , 2 \pi ]$ (see \cite[(4.5)]{CoGoHaJeKn96_01}), we deduce
\begin{equation} \label{e80}
\re z_{q , 0} \in E_{0} + h \tau_{0} + 2 q \pi \lambda_{1} \frac{h}{\vert \ln h \vert} + [ - 2 \pi , 2 \pi ] \lambda_{1} \frac{h}{\vert \ln h \vert} ,
\end{equation}
Thus, the real part of the $z_{q , 0}$'s is regularly distributed on average. Comparing with \eqref{d95}, we note that the resonances in the generic case satisfy also \eqref{e80}.

Let us now consider the asymptotic of the $z_{q , 0}$'s for $q$ fixed, which correspond to the resonances whose real part is the closest to $E_{0} + h \tau_{0}$ from the above discussion. Using the asymptotic of $W_{q} (A)$ for $A$ large stated in \cite[(4.20)]{CoGoHaJeKn96_01} (see also Bachelot and Motet-Bachelot \cite[(IV.23)]{BaMo93_01} for $q = 0$), \eqref{e65} yields
\begin{equation}
z_{q , 0} = E_{0} + h \tau_{0} - i h \sum_{j = 2}^{n} \frac{\lambda_{j}}{2} - i \lambda_{1} \big( \ln \vert \ln h \vert - \ln \ln \vert \ln h \vert + 2 q \pi i + c_{h} + o (1) \big) \frac{h}{\vert \ln h \vert} ,
\end{equation}
with
\begin{equation*}
c_{h} : = \ln \Big( \frac{i}{\lambda_{1} \alpha} e^{- i \frac{\tau_{0}}{\lambda_{1}} \vert \ln h \vert} \Big) .
\end{equation*}
In the previous equations, ``$\ln$'' denotes the principal branch of the logarithm (i.e. $\im \ln \in ] - \pi , \pi ]$). In particular, $c_{h}$ is bounded with respect to $h$. This implies that the imaginary part of the resonance below $E_{0} + h \tau_{0}$ satisfies
\begin{equation} \label{e81}
\im z_{q , 0} = - h \sum_{j = 2}^{n} \frac{\lambda_{j}}{2} - \lambda_{1} \frac{h \ln \vert \ln h \vert}{\vert \ln h \vert} + \lambda_{1} \frac{h \ln \ln \vert \ln h \vert}{\vert \ln h \vert} + \CO \Big( \frac {h}{\vert \ln h \vert} \Big) ,
\end{equation}
whereas the $z_{q , k} ( \tau )$'s given by \eqref{d95} satisfy
\begin{equation} \label{e83}
\im z_{q , k} ( \tau ) = - h \sum_{j = 2}^{n} \frac{\lambda_{j}}{2} + \CO \Big( \frac {h}{\vert \ln h \vert} \Big) .
\end{equation}
Thus, the behavior of the imaginary parts is different from the generic one. That the resonances are deeper in the complex plane here is in adequacy with Theorem \ref{d8}.

More generally, Lemma \ref{e55} shows that 
\begin{equation*}
\im z_{q , 0} \geq - h \sum_{j = 2}^{n} \frac{\lambda_{j}}{2} - \lambda_{1} \frac{h \ln \vert \ln h \vert}{\vert \ln h \vert} + \CO \Big( \frac {h}{\vert \ln h \vert} \Big) ,
\end{equation*}
and
\begin{equation*}
\Big\vert z_{q , 0} - \Big( E_{0} + h \tau_{0} - i h \sum_{j = 2}^{n} \frac{\lambda_{j}}{2} \Big) \Big\vert \gtrsim \frac{h \ln \vert \ln h \vert}{\vert \ln h \vert}  ,
\end{equation*}
for all $q \in \Z$ such that $z_{q , k} \in \eqref{e53}$. Eventually, for the large values of $q$, more precisely if $\ln \vert \ln h \vert \ll \vert q \vert  \ll \vert \ln h \vert$, Lemma \ref{e55} yields
\begin{align}
\im z_{q , 0} &= - h \sum_{j = 2}^{n} \frac{\lambda_{j}}{2} - \lambda_{1} \frac{h}{\vert \ln h \vert} \ln \Big( \frac{\vert \ln h \vert}{q} \Big) + \CO \Big( \frac {h}{\vert \ln h \vert} \Big)  \nonumber \\
&= - h \sum_{j = 2}^{n} \frac{\lambda_{j}}{2} + \lambda_{1} \frac{h}{\vert \ln h \vert} \ln \Big( \frac{\vert \re z_{q , 0} - ( E_{0} + h \tau_{0} ) \vert}{h} \Big) + \CO \Big( \frac {h}{\vert \ln h \vert} \Big) .  \label{e82}
\end{align}

Summing up, Proposition \ref{e56} (and more precisely \eqref{e82}) can be interpreted as a microscopic transitional regime. Corresponding to the small values of $q \in \Z$, the resonances just below $E_{0} + h \tau_{0}$ are furthest away from the real axis (see \eqref{e81}). Note that they belong in the region where $\mu$ is the smallest. But when $\vert q \vert$ increases, \eqref{e82} shows that the resonances get closer to the real axis. Eventually, if we formally apply \eqref{e82} with $q$ of order $\vert \ln h \vert$ (which correspond to the resonances at a distance $h$ from $E_{0} + h \tau_{0}$), we recognize \eqref{e83}.

\subsubsection{Other subsets of resonances} \label{s83}

In Proposition \ref{e15}, we have shown that $P$ has no resonance below the ones described by Theorem \ref{d8} up to the line
\begin{equation} \label{g84}
\im z  = - h \sum_{j = 2}^{n} \frac{\lambda_{j}}{2} - h \lambda_{1} ,
\end{equation}
if $0$ is away from the spectrum of $\CQ ( z, h )$. In this last part of Section \ref{s43}, we investigate the distribution of resonances near and below this line. The idea is that the quantization rule stated in Definition \ref{d1} only describes the ``first set of accumulation curves'' and that the resonances near \eqref{g84} are given by a new quantization rule. Roughly speaking, the resonances near $\im z  = - h \sum_{j = 2}^{n} \lambda_{j} / 2$ correspond to the $K$ ``values on the homoclinic trajectories'' whereas the resonances near \eqref{g84} correspond to the $(n-1) \times K$ ``first derivatives on the homoclinic trajectories''.

Let $a < b$ and assume that
\begin{hyp} \label{h14}
For all $\sigma \in [ a , b ] - i ( \sum_{j = 2}^{n} \lambda_{j} / 2 + \lambda_{1} )$, we have $\det \widetilde{\CZ} ( \sigma ) \neq 0$.
\end{hyp}
Recall that $\widetilde{\CZ}$ is defined in \eqref{g82}. This assumption guaranties that the possible resonances coming from Definition \ref{d1} will not perturb the setting near the line \eqref{g84}.

Under the previous assumption, there exists a $( n - 1 ) K \times ( n - 1 ) K$ matrix $\CQ^{2} ( z , h )$ governing the asymptotic of the resonances near the line \eqref{g84} and satisfying the following properties. As in \eqref{d10}, the matrix $\CQ^{2}$ can be written as
\begin{equation} \label{g91}
\CQ^{2} ( z , h ) = \widetilde{\CQ}^{2} ( \rho , \sigma ) = \sum_{k = 1}^{K} \rho_{k} \widetilde{\CQ}^{2}_{k} ( \sigma ) ,
\end{equation}
where $\rho_{k} = e^{i A_{k} / h}$ corresponds to the action $A_{k}$, $\sigma$ is the rescaled spectral parameter given by \eqref{d92} and $\widetilde{\CQ}^{2}_{k}$ is independent of $h$ and holomorphic in $\sigma$ near $\im \sigma = - \sum_{j = 2}^{n} \lambda_{j} / 2 - \lambda_{1}$. Contrary to \eqref{d4}, we do not give an explicit expression for the matrix $\CQ^{2}$ in terms of geometric quantities. Its construction is explained in Section \ref{s84}. In particular, the coefficients $\widetilde{\CQ}^{2}_{k}$ only depend on the symbol of $p ( x , \xi )$ in a neighborhood of $K ( E_{0} )$.  Mimicking Definition \ref{d1}, the quantization rule for resonances near \eqref{g84} is given by

\begin{definition}[Quantization rule for the second set of resonances]\sl \label{g85}
We say that $z$ is a pseudo-resonance of the second kind if and only if
\begin{equation}  \label{g83}
1 \in \spe \big( h^{S ( z , h ) / \lambda_{1} + 1 / 2} \CQ^{2} ( z , h ) \big) .
\end{equation}
The set of pseudo-resonances of the second kind is denoted by $\res_{0}^{2} ( P )$.
\end{definition}

As in Section \ref{s61}, we express the behavior of the pseudo-resonances of the second kind in terms of the spectrum of $\CQ^{2}$. In the present context, let $\mu_{1} ( \tau , h ) , \ldots , \mu_{( n - 1 ) K} ( \tau , h )$ denote the eigenvalues of
\begin{equation*}
\CQ^{2} \Big( E_{0} + h \tau - i h \sum_{j = 2}^{n} \frac{\lambda_{j}}{2} - i h \lambda_{1} , h \Big) .
\end{equation*}
The pseudo-resonances of the second kind satisfy the following two scale asymptotic. We omit the proof of this result since it is similar to the one of Proposition \ref{d9}.

\begin{proposition}\sl \label{g88}
Assume \ref{h1}--\ref{h4}, \ref{h8}, \ref{h14} and let $C > 0$. The pseudo-resonances of the second kind $z$ lying in
\begin{equation} \label{i27}
E_{0} + [ a h , b h ] - i h \sum_{j = 2}^{n} \frac{\lambda_{j}}{2} - i h \lambda_{1} + i \Big[ - C \frac{h}{\vert \ln h \vert} , C \frac{h}{\vert \ln h \vert} \Big] ,
\end{equation}
satisfy $z = z_{q , k}^{2} ( \tau ) + o ( h \vert \ln h \vert^{- 1} )$ in the sense of Proposition \ref{d9} with
\begin{equation*}
z_{q , k}^{2} ( \tau ) = E_{0} + 2 q \pi \lambda_{1} \frac{h}{\vert \ln h \vert} - i h \sum_{j = 2}^{n} \frac{\lambda_{j}}{2} - i h \lambda_{1} + i \ln ( \mu_{k} ( \tau , h ) ) \lambda_{1} \frac{h}{\vert \ln h \vert} ,
\end{equation*}
for some $q \in \Z$ and $k \in \{ 1 , \ldots , ( n - 1 ) K \}$.
\end{proposition}

\begin{figure}
\begin{center}
\begin{picture}(0,0)%
\includegraphics{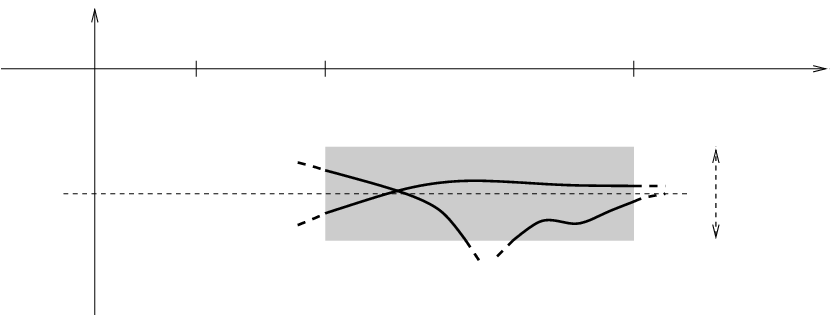}%
\end{picture}%
\setlength{\unitlength}{987sp}%
\begingroup\makeatletter\ifx\SetFigFont\undefined%
\gdef\SetFigFont#1#2#3#4#5{%
  \reset@font\fontsize{#1}{#2pt}%
  \fontfamily{#3}\fontseries{#4}\fontshape{#5}%
  \selectfont}%
\fi\endgroup%
\begin{picture}(15944,5969)(-621,-5708)
\put(3151,-586){\makebox(0,0)[b]{\smash{{\SetFigFont{9}{10.8}{\rmdefault}{\mddefault}{\updefault}$E_{0}$}}}}
\put(5626,-586){\makebox(0,0)[b]{\smash{{\SetFigFont{9}{10.8}{\rmdefault}{\mddefault}{\updefault}$E_{0} + a h$}}}}
\put(11551,-586){\makebox(0,0)[b]{\smash{{\SetFigFont{9}{10.8}{\rmdefault}{\mddefault}{\updefault}$E_{0} + b h$}}}}
\put(13576,-3436){\makebox(0,0)[lb]{\smash{{\SetFigFont{9}{10.8}{\rmdefault}{\mddefault}{\updefault}$\ds 2 C \frac{h}{\vert \ln h \vert}$}}}}
\put(241,-3436){\makebox(0,0)[rb]{\smash{{\SetFigFont{9}{10.8}{\rmdefault}{\mddefault}{\updefault}$\ds - \sum_{j=2}^{n} \frac{\lambda_{j}}{2} h - \lambda_{1} h$}}}}
\end{picture}%
\end{center}
\caption{The region studied in Theorem \ref{g86}.} \label{f34}
\end{figure}

With the notation of Definition \ref{g80}, the following result shows that the resonances near \eqref{g84} are close to the pseudo-resonances of the second kind. The setting is illustrated in Figure \ref{f34}.

\begin{theorem}[Asymptotic of the second set of resonances]\sl \label{g86}
Assume \ref{h1}--\ref{h4}, \ref{h8}, \ref{h14} and let $C , \delta > 0$. In the domain
\begin{equation} \label{g87}
E_{0} + [ a h , b h ]  - i h \sum_{j = 2}^{n} \frac{\lambda_{j}}{2} - i h \lambda_{1} + i \Big[ - C \frac{h}{\vert \ln h \vert} , C \frac{h}{\vert \ln h \vert} \Big] \setminus \big( \Gamma (h) + B ( 0 , \delta h ) \big),
\end{equation}
we have
\begin{equation*}
\dist \big( \res (P) , \res_{0}^{2} (P) \big) = o \Big( \frac {h}{\vert \ln h \vert} \Big) ,
\end{equation*}
as $h$ goes to $0$. Moreover, for all $\chi \in C^{\infty}_{0} ( \R^{n} )$, there exists $M > 0$ such that
\begin{equation*}
\big\Vert \chi ( P -z )^{-1} \chi \big\Vert \lesssim h^{- M} ,
\end{equation*}
uniformly for $h$ small enough and $z \in \eqref{g87}$ with $\dist ( z , \res_{0}^{2} ( P ) ) \geq \delta h \vert \ln h \vert^{- 1}$.
\end{theorem}

Since Theorem \ref{g86} is similar to Theorem \ref{d8}, the resonances lying in \eqref{g87} verify some of the phenomena described in Section \ref{s18}. Thus, they satisfy a two scale asymptotic (see Remark \ref{e2}) with at most $( n - 1 ) K$ accumulation curves which may vibrate as explained in Remark \ref{e5}. Roughly speaking, they correspond to the ``second set of accumulation curves''. Nevertheless, the transition and stability phenomena of Section \ref{s14} and Section \ref{s24} can not directly be observed here since we do not have an explicit expression of $\CQ^{2}$ in terms of dynamical quantities. One may also consider settings where $0 \in \spe ( \CQ^{2} ( z , h ) )$ as in Section \ref{s43}, prove asymptotic modulo $ \CO ( h^{\infty} )$ as in Section \ref{s79}, \ldots

As explained in the proof of Theorem \ref{g86}, the resonances of the second kind are generated by the $n - 1$ transversal directions to the $K$ homoclinic trajectories. Indeed, the quantization rule \eqref{i28} is about the derivatives of order $1$ in these directions. Moreover, it is explained in Proposition \ref{l28} that the leading term of the associated resonant states vanishes at order $1$ on the homoclinic trajectories. Thus, the situation seems to be similar to the case of a hyperbolic trajectory treated by G\'erard and Sj\"{o}strand \cite{GeSj87_01} or the case of the harmonic oscillator. Nevertheless, the mechanism is different here. Since the subprincipal term in the stationary phase expansion (arising from the transition through the hyperbolic fixed point) is given by a differential operator of order $2$, the derivatives of order $1$ and $2$ have to appear in the quantization rule (see \eqref{i1} and \eqref{i2}). But, thanks to some particular structure explained in Lemma \ref{g95}, the terms of order $2$ disappear.

In the one dimensional case, the matrix $\CQ^{2} ( z , h )$ is by definition the trivial matrix on $\C^{( n - 1 ) K} = \{ 0 \}$. Then, the set of pseudo-resonances of the second kind is empty and Theorem \ref{g86} means that $P$ has no resonance and a polynomial estimate of its truncated resolvent in \eqref{g87}. More generally, the remark below shows that there is no resonance when $n = 1$ away from the ones given by Theorem \ref{d8}, the zeros of $\det \widetilde{\CZ}$ and $\Gamma ( h )$. In some sense, it completes the resonance free zone of Proposition \ref{e15}. Note that this result is a consequence \cite[Th\'eor\`eme 0.7]{FuRa98_01} for two heteroclinic trajectories. Recall that \ref{h4} and \ref{h8} are always satisfied in dimension $n = 1$ (see Section \ref{s19}).

\begin{remark}\sl \label{i23}
Assume $n = 1$ and \ref{h1}--\ref{h3}. Let $\delta > 0$ and $\Omega$ be a compact subset of $\R + i ] - \infty , 0 [$ such that $\det \widetilde{\CZ} ( \sigma ) \neq 0$ for all $\sigma \in \Omega$. Then, for $h$ small enough, $P$ has no resonance and a polynomial estimate of its truncated resolvent in the domain $E_{0} + h \Omega \setminus ( \Gamma ( h ) + B ( 0 , \delta h ) )$.
\end{remark}

\begin{example}\rm \label{i29}
We now consider operators in dimension $n \geq 2$ for which Theorem \ref{g86} actually provides resonances in the region \eqref{g87}. Let
\begin{equation*}
V ( x ) = V_{1} ( x ) + V_{2} ( x - M e_{1} ) ,
\end{equation*}
where $e_{1}$ is the first vector of the canonical basis, $M > 1$ is large enough and $V_{\bullet} \in C^{\infty}_{0} ( \R^{n} )$ is a radial function satisfying $x \cdot \nabla V_{\bullet} (x) < 0$ for $x$ in the interior of $\supp V_{\bullet} \setminus \{ 0 \}$. Moreover, we assume that $V_{1} ( x ) = E_{0} - \lambda_{1}^{2} x^{2} /4$ near $0$ and $V_{2} ( 0 ) > E_{0}$. In this case, \ref{h1}--\ref{h4} are satisfied and $\CH \subset \{ x_{2} = \cdots = x_{n} = 0 \}$ consists of a single trajectory. Moreover, one can verify that \ref{h8} (and then \ref{h9}) holds true. Eventually, as noted below Proposition \ref{e15}, \ref{h14} is always verified under \ref{h9}. In dimension $n = 2$, the setting is close to the one of Example \ref{c15}. In this way, the first picture in Figure \ref{f20} represents such a potential.

\begin{figure}
\begin{center}
\begin{picture}(0,0)%
\includegraphics{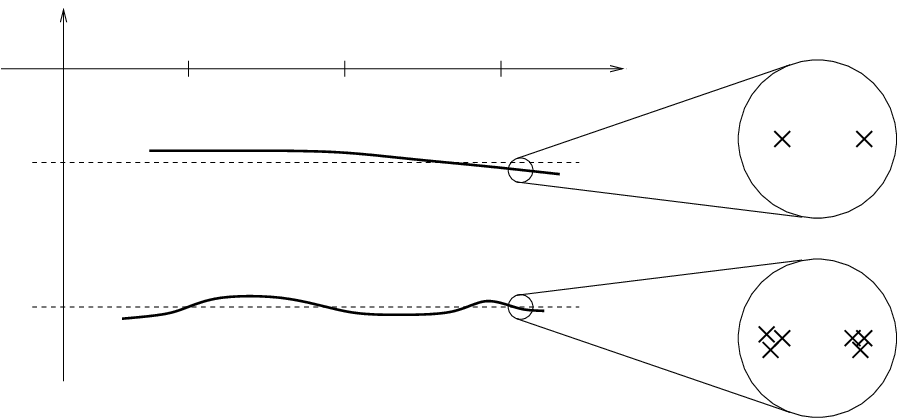}%
\end{picture}%
\setlength{\unitlength}{987sp}%
\begingroup\makeatletter\ifx\SetFigFont\undefined%
\gdef\SetFigFont#1#2#3#4#5{%
  \reset@font\fontsize{#1}{#2pt}%
  \fontfamily{#3}\fontseries{#4}\fontshape{#5}%
  \selectfont}%
\fi\endgroup%
\begin{picture}(17231,7931)(-621,-7670)
\put(-359,-5611){\makebox(0,0)[rb]{\smash{{\SetFigFont{9}{10.8}{\rmdefault}{\mddefault}{\updefault}$\ds - \sum_{j=2}^{n} \frac{\lambda_{j}}{2} h - \lambda_{1} h$}}}}
\put(3001,-586){\makebox(0,0)[b]{\smash{{\SetFigFont{9}{10.8}{\rmdefault}{\mddefault}{\updefault}$E_{0} - C h$}}}}
\put(6001,-586){\makebox(0,0)[b]{\smash{{\SetFigFont{9}{10.8}{\rmdefault}{\mddefault}{\updefault}$E_{0}$}}}}
\put(9001,-586){\makebox(0,0)[b]{\smash{{\SetFigFont{9}{10.8}{\rmdefault}{\mddefault}{\updefault}$E_{0} + C h$}}}}
\put(-359,-2836){\makebox(0,0)[rb]{\smash{{\SetFigFont{9}{10.8}{\rmdefault}{\mddefault}{\updefault}$\ds - \sum_{j=2}^{n} \frac{\lambda_{j}}{2} h$}}}}
\end{picture}%
\end{center}
\caption{The accumulation curves of Example \ref{i29}. The first one, of multiplicity one, is given by Theorem \ref{d8}. The second one, of multiplicity $n - 1$ is given by Theorem \ref{g86}.} \label{f33}
\end{figure}

In this symmetrical situation, the $( n - 1 ) \times ( n - 1 )$ matrix $\CQ^{2}$ takes the form
\begin{equation} \label{i30}
\CQ^{2} ( z , h ) = q^{2} ( z , h ) Id ,
\end{equation}
with $q^{2} ( z , h ) = e^{i A_{1} / h} \widetilde{q}^{2} ( \sigma )$ and $\widetilde{q}^{2} ( \sigma ) \neq 0$ for $\im \sigma = - \sum_{j = 2}^{n} \lambda_{j} / 2 - \lambda_{1}$. The proof of this formula can be found in Section \ref{s84}. Thus, the eigenvalues $\mu_{1} ( \tau , h ) , \ldots , \mu_{n - 1} ( \tau , h )$ coincide and are away from $0$. From Theorem \ref{g86}, the resonances $z$ of $P$ in \eqref{g87} with $\re z = E_{0} + \tau h + h \delta (h) [ - 1 , 1 ]$ satisfy
\begin{equation*}
z = z_{q}^{2} ( \tau ) + o \Big( \frac{h}{\vert \ln h \vert} \Big) ,
\end{equation*}
with $q \in \Z$ and
\begin{equation*}
z_{q}^{2} ( \tau ) = E_{0} - \frac{A_{1} \lambda_{1}}{\vert \ln h \vert} + 2 q \pi \lambda_{1} \frac{h}{\vert \ln h \vert} - i h \sum_{j = 2}^{n} \frac{\lambda_{j}}{2} - i h \lambda_{1} + i \ln \Big( \widetilde{q}^{2} \Big( \tau - i \sum_{j = 2}^{n} \frac{\lambda_{j}}{2} - i \lambda_{1} \Big) \Big) \lambda_{1} \frac{h}{\vert \ln h \vert} .
\end{equation*}
This formula can be compared with \eqref{e8} which gives the asymptotic of the pseudo-resonances of the first kind (i.e. defined by Definition \ref{d1}). Furthermore, Proposition \ref{d91}, adapted to the present setting, implies that there exist at least $n - 1$ resonances near each $z_{q}^{2} ( \tau )$. In particular, the resonances in \eqref{g87} concentrate on a unique accumulation curve (which formally has multiplicity $n - 1$). It is given by
\begin{equation*}
\im \sigma = - \sum_{j = 2}^{n} \frac{\lambda_{j}}{2} - \lambda_{1} + \ln \Big\vert \widetilde{q}^{2} \Big( \re \sigma - i \sum_{j = 2}^{n} \frac{\lambda_{j}}{2} - i \lambda_{1} \Big) \Big\vert \frac{\lambda_{1}}{\vert \ln h \vert} .
\end{equation*}
Note that this curve does not vibrate since there is only one action (see Remark \ref{e5}). The distribution of the resonances is illustrated in Figure \ref{f33}. To be precise, we have not yet proved that $P$ has no resonance in
\begin{equation*}
E_{0} + [ - C h , C h ] - i h \sum_{j = 2}^{n} \frac{\lambda_{j}}{2} - i h \lambda_{1} + i \Big[ C \frac{h}{\vert \ln h \vert} , \frac{h}{C} \Big] \setminus \big( \Gamma (h) + B ( 0 , \delta h ) \big) .
\end{equation*}
Indeed, this region is not covered by Proposition \ref{e15} and Theorem \ref{g86}. But it can be done following the proof of Theorem \ref{g86}.
\end{example}

\begin{remark}\sl
Assume $n \geq 2$. Theorem \ref{d8} and Theorem \ref{g86} give the asymptotic of the resonances near the lines
\begin{equation*}
\CL^{1} : = \Big\{ \im z  = - h \sum_{j = 2}^{n} \frac{\lambda_{j}}{2} \Big\}  \qquad \text{and} \qquad \CL^{2} : = \Big\{\im z  = - h \sum_{j = 2}^{n} \frac{\lambda_{j}}{2} - h \lambda_{1} \Big\} .
\end{equation*}
Moreover, under some hypotheses, there is no resonance between $\CL^{0} : = \R$ and $\CL^{1}$ (from Theorem \ref{d8}) and between $\CL^{1}$ and $\CL^{2}$ (see Section \ref{s15}). Using the methods developed in this part, one can probably give the asymptotic of the resonances in deeper zones. It is natural to imagine that most of them are close to some lines $\CL^{j}$ (see G\'erard and Sj\"{o}strand \cite{GeSj87_01}, Dyatlov \cite{Dy12_01} or Faure and Tsujii \cite{FaTs15_01} in the context of hyperbolic geometries).
\end{remark}

Formally, the approach to obtain the asymptotic of resonances in deeper zones is simple: one keeps more terms in the monodromy asymptotic (see Lemma \ref{g89}), expresses the first derivatives in the transversal directions in terms of higher derivatives (as in \eqref{g94}) and eventually gives the quantization rule (as in \eqref{i1} and \eqref{i2}). Unfortunately, some cancellations happen (see Lemma \ref{g95}) and the computation becomes complicated. Thus, the intuition that $\CL^{j}$ should be defined by
\begin{equation*}
\CL^{j} = \Big\{\im z  = - h \sum_{j = 2}^{n} \frac{\lambda_{j}}{2} - ( j - 1 ) h \lambda_{1} \Big\} ,
\end{equation*}
is probably inexact (the $\mu_{k}$'s, defined at the end of Section \ref{s35}, should play some role). As explained previously, this difficulty comes from the fact that the quantum monodromy does not behave like a harmonic oscillator in the transversal directions (as in the hyperbolic settings).

\Subsection{Asymptotic of higher order} \label{s79}

In Theorem \ref{d8} and more generally in the rest of Section \ref{s6}, the localization of the resonances is given modulo remainder terms of size $o ( h \vert \ln h \vert^{- 1} )$. The aim of this part is to obtain more accurate asymptotic and to give the localization of the resonances modulo $\CO ( h^{\infty} )$. This question can be tricky to deal with when the eigenvalues of $\CQ ( z , h )$ cross. Indeed, as it appears in Section \ref{s12} and Section \ref{s80}, the behavior of the inverse of $1 - h^{S ( z , h ) / \lambda_{1} - 1 / 2} \CQ ( z , h )$ near the resonances (and then the behavior of the resolvent of $\CQ$ near its eigenvalues) plays a central role here. To avoid this difficulty, we will work in a simple situation.

We assume here that $\CH$ consists of a unique trajectory (i.e. $K = 1$). More precisely, we suppose \ref{h9}. In this case, the matrix $\CQ$ is a scalar, the asymptotic of the pseudo-resonances is given by \eqref{e8} and Corollary \ref{e4} describes the accumulation curve of resonances.

There exists a function $\CQ_{\rm tot} ( z , h )$ governing the asymptotic of the resonances up to order $\CO ( h^{\infty} )$ and satisfying the following properties. It is defined for $z \in \eqref{d12}$ and $h$ small enough. In this domain, $\CQ_{\rm tot}$ verifies the asymptotic
\begin{equation} \label{e95}
\CQ_{\rm tot} ( z , h ) \simeq \sum_{a = 0}^{+ \infty} \sum_{b = 0}^{B_{a}} \sum_{c = 0}^{C_{a}} \CQ_{a , b , c} ( z , h ) \big( h^{S ( z , h ) / \lambda_{1} - 1 / 2} \big)^{b} ( \ln h )^{c} h^{\widehat{\mu}_{a} / \lambda_{1}} ,
\end{equation}
where $B_{a} , C_{a} \in \N$, $( \widehat{\mu}_{a} )_{a \geq 0}$ is the increasing sequence of the linear combinations over $\N$ of $\lambda_{1}$ and $\lambda_{j} - \lambda_{1}$, and $\CQ_{a , b , c} ( z , h )$ can be written (as in \eqref{d10})
\begin{equation} \label{g14}
\CQ_{a , b , c} ( z , h ) = e^{i A_{1} / h} \widetilde{\CQ}_{a , b , c} ( \sigma ) ,
\end{equation}
with $\sigma = ( z - E_{0} ) / h$. For $a = 0$, we have $B_{a} = C_{a} = 0$ and $\CQ_{0 , 0 , 0} ( z , h ) = \CQ ( z , h )$. Eventually, $\CQ_{\rm tot}$ as well as the $\CQ_{a , b , c}$'s are holomorphic in the domain \eqref{d12}. Note that $\vert h^{S ( z , h ) / \lambda_{1} - 1 / 2} \vert \lesssim 1$ for $z \in \eqref{d12}$.

As in \eqref{g91}, we do not give an explicit formula for $\CQ_{\rm tot}$ and $\CQ_{a , b , c}$. Their construction is explained in Section \ref{s80}. In particular, the coefficients $\CQ_{a , b , c}$ only depend on the symbol of $p ( x , \xi )$ in a neighborhood of $K ( E_{0} )$. Moreover, \eqref{e95} implies
\begin{equation} \label{e99}
\CQ_{\rm tot} ( z , h ) = \CQ ( z , h ) + \CO \big( ( \ln h )^{B_{1}} h^{\widehat{\mu}_{1} / \lambda_{1}} \big) .
\end{equation}
Thus, the $1 \times 1$ matrix $\CQ_{\rm tot}$ is a precise version of $\CQ$. Mimicking Definition \ref{d1}, the quantization rule modulo $\CO ( h^{\infty} )$ is given by

\begin{definition}\sl \label{e97}
We say that $z$ is a pseudo-resonance at infinite order if and only if
\begin{equation*}
h^{S ( z , h ) / \lambda_{1} - 1 / 2} \CQ_{\rm tot} ( z , h ) = 1 .
\end{equation*}
The set of pseudo-resonances at infinite order is denoted by $\res_{\infty} ( P )$.
\end{definition}

Since $\CQ_{\rm tot}$ is a perturbation of $\CQ$, the asymptotic of the pseudo-resonances at infinite order is still given by Proposition \ref{d9}. More precisely, we have the following result where $z_{q} ( \tau )$ is given by \eqref{e8}.

\begin{proposition}\sl \label{g1}
Assume \ref{h1}--\ref{h4}, \ref{h9}, let $C > 0$ and let $\delta (h)$ be a function which goes to $0$ as $h \to 0$. Then, uniformly for $\tau \in [ - C , C ]$, the pseudo-resonances at infinite order $z$ in \eqref{d12} with $\re z \in E_{0} + \tau h + h \delta (h) [ - 1 , 1 ]$ satisfy
\begin{equation} \label{g2}
z = z_{q} ( \tau ) + o \Big( \frac{h}{\vert \ln h \vert} \Big) ,
\end{equation}
for some (unique) $q \in \Z$.

On the other hand, for each $\tau \in [ - C , C ]$ and $q \in \Z$ such that $z_{q} ( \tau )$ belongs to \eqref{d12} with a real part lying in $E_{0} + \tau h + h \delta (h) [ - 1 , 1 ]$, there exists a unique pseudo-resonance at infinite order $z$ satisfying \eqref{g2} uniformly with respect to $q , \tau$.
\end{proposition}

Note that the remainder term in \eqref{g2} is $o ( h \vert \ln h \vert^{- 1} )$ and not $\CO ( h^{\infty} )$, as in \eqref{d94}. It could be tough to give an explicit formula for the pseudo-resonances at infinite order modulo $\CO ( h^{\infty} )$. Indeed, the asymptotic of the solutions of $h^{\sigma} f ( \sigma ) = 1$ is already typically in terms of $\vert \ln h \vert^{- N}$, $N \in \N$. We do not discuss this point. Nevertheless, Proposition \ref{g1} shows that $\res_{\infty} (P)$ is in bijection with the $z_{q} ( \tau )$'s. Finally, combining with Proposition \ref{d9}, Proposition \ref{g1} yields
\begin{equation} \label{g7}
\dist \big( \res_{0} (P) , \res_{\infty} (P) \big) = o \Big( \frac {h}{\vert \ln h \vert} \Big) ,
\end{equation}
in the domain \eqref{d12}.

\begin{theorem}[Asymptotic of resonances modulo $\CO ( h^{\infty} )$]\sl \label{g3}
Assume \ref{h1}--\ref{h4}, \ref{h9} and let $C , \delta > 0$. In the domain \eqref{d90}, we have
\begin{equation*}
\dist \big( \res (P) , \res_{\infty} (P) \big) = \CO ( h^{\infty} ) ,
\end{equation*}
as $h$ goes to $0$. Moreover, for all $\chi \in C^{\infty}_{0} ( \R^{n} )$, there exists $M > 0$ such that
\begin{equation} \label{g8}
\big\Vert \chi ( P -z )^{-1} \chi \big\Vert \lesssim h^{- M} ,
\end{equation}
uniformly for $h$ small enough and $z \in \eqref{d90}$ with $\dist ( z , \res_{\infty} ( P ) ) \geq h^{C}$.
\end{theorem}

Contrary to the $\CQ_{a , b , c}$'s, the total $1 \times 1$ matrix $\CQ_{\rm tot}$ is not unique. In fact, any $\CQ_{\rm tot}$ holomorphic in \eqref{d12} and satisfying \eqref{e95} will give the same set of pseudo-resonance at infinite order modulo $\CO ( h^{\infty} )$. As consequence, the conclusions of Theorem \ref{g3} are independent of this choice of $\CQ_{\rm tot}$.

\begin{remark}\sl
The asymptotic of the resonances modulo $\CO ( h^{\infty} )$ is stated and proved under the hypothesis that $\CH$ consists of a single trajectory. As explained in the beginning of the section, the general case seems difficult due to the possible presence of crossings of the eigenvalues of $\CQ$. Nevertheless, the results of this part may be extended outside of these crossings. 
\end{remark}

\Subsection{Tangential intersection of finite order} \label{s81}

Here, we relax \ref{h8} and assume that the manifolds $\Lambda_{-}$ and $\Lambda_{+}$ have an intersection of finite order along a finite number of trajectories. The situation is then intermediate between the transversal case treated in Theorem \ref{d8} and the fully tangential case treated in Theorem \ref{i55}. As explained in Remark \ref{j63} $ii)$, one could prove that the resonances are at least at distance of order $h$ from the real axis in dimension $n \geq 2$. The proof of the geometric assertions and results stated below can be found in Section \ref{s82}.

In the finite order intersection situation, the order of the contact between $\Lambda_{-}$ and $\Lambda_{+}$ may depend on the direction. It has no reason to be uniform. Thus, in large dimension, the asymptotic of the resonances can be rather complicated. To avoid this problem, we work in dimension $n = 2$ in the sequel. Then, near a point of the tangential intersection $\Lambda_{-} \cap \Lambda_{+}$, one direction of these $2$ dimensional manifolds corresponds to the Hamiltonian vector field and the order of contact is measured on the other ones. We assume that
\begin{hyp} \label{h13}
The homoclinic set $\CH$ consists of a finite number of trajectories on which $\Lambda_{-}$ and $\Lambda_{+}$ have an intersection of finite order.
\end{hyp}
As in Section \ref{s61}, we denote by $\gamma_{k}$, $k = 1 , \ldots , K$, these Hamiltonian trajectories. To each of theses curves, we associate its action $A_{k}$, its Maslov's index $\nu_{k}$ and its asymptotic direction $g_{\pm}^{k}$ defined in \eqref{d6}.

For $k \in \{ 1 , \ldots , K \}$, we make a linear change of coordinates such that $g_{-}^{k}$ is (positively) aligned with the base vector $( 1 , 0 )$. Let $\Lambda_{+}^{1}$ denote the evolution of $\Lambda_{+}^{0}$ by the Hamiltonian flow after a turn along $\CH$. Then, the manifold $\Lambda_{+}^{1}$ projects nicely on the $x$-space near $\pi_{x} ( \gamma_{k} \cap \Lambda_{-}^{0} )$. Let $\varphi_{+}^{k} \in C^{\infty} ( \R^{n} )$ be the unique generating function of $\Lambda_{+}^{1}$ (i.e. $\Lambda_{+}^{1} = \{ ( x , \nabla \varphi_{+}^{k} (x) ) \}$) defined near $\pi_{x} ( \gamma_{k} \cap \Lambda_{-}^{0} )$ with the normalization $\varphi_{+}^{k} = \varphi_{-}$ on the curve $\pi_{x} ( \gamma_{k} )$. The assumption \ref{h13} implies that, for all $x_{1} > 0$ small enough, one can write
\begin{equation} \label{g21}
\varphi_{+}^{k} ( x_{1} , x_{2} ) - \varphi_{-} ( x_{1} , x_{2} ) = \alpha_{k} ( x_{1} ) \big( x_{2} - x_{2}^{k} ( x_{1} ) \big)^{1 + m_{k}} + \CO \big( \big( x_{2} - x_{2}^{k} ( x_{1} ) \big)^{2 + m_{k}} \big) .
\end{equation}
for some $m_{k} \in \N \setminus \{ 0 \}$ and $\alpha_{k} ( x_{1} ) \neq 0$. Here, $x_{2}^{k} ( x_{1} )$ is the second spatial coordinate of the unique point of $\gamma_{k} \cap \Lambda_{-}^{0}$ with first spatial coordinate $x_{1}$. Since $\Lambda_{\pm}$ are stable by the Hamiltonian flow, $m_{k}$ is independent of $x_{1}$. The fact that $m_{k} \geq 1$ follows from the normalization $\varphi_{+}^{k} = \varphi_{-}$ on $\pi_{x} ( \gamma_{k} )$ and from $\gamma_{k} ( t ) = ( x_{k} ( t ) , \nabla \varphi_{-} ( x_{k} ( t ) ) ) = ( x_{k} ( t ) , \nabla \varphi_{+}^{k} ( x_{k} ( t ) ) ) \in \Lambda_{-} \cap \Lambda_{+}$. The number $m_{k}$ is called the order of the intersection of $\Lambda_{-}$ and $\Lambda_{+}$ along $\gamma_{k}$. In particular, $m_{k} = 1$ corresponds to a transversal intersection as in the assumption \ref{h8}. Contrary to $m_{k}$, the coefficient $\alpha_{k} ( x_{1} )$ depends on $x_{1}$. More precisely, if $m_{k} \geq 2$, we have
\begin{equation} \label{g22}
\alpha_{k} ( x_{1} ) = \alpha_{k}^{\infty} x_{1}^{- ( 1 + m_{k} ) \lambda_{2} / \lambda_{1}} + o \big( x_{1}^{- ( 1 + m_{k} ) \lambda_{2} / \lambda_{1}} \big) ,
\end{equation}
in the limit $x_{1}$ goes to $0$ for some $\alpha_{k}^{\infty} \neq 0$. That $\alpha_{k} ( x_{1} )$ diverges when $x_{1}$ goes to $0$ can be formally explained. Indeed, Lemma \ref{a43} shows that the Hamiltonian trajectories in $p^{- 1} ( E_{0} )$ close but not in $\Lambda_{-}$ must escape in the orthogonal of their incoming direction (see Section \ref{s11}). Thus, when $x_{1}$ is small, $\Lambda_{+}^{1}$ ``makes the transition'' between $\gamma_{k}$ which converges straight to $( 0 , 0 )$  and trajectories orthogonal to $\gamma_{k}$. Roughly speaking, \eqref{g22} quantifies this abrupt transition.

We now defined a matrix of interaction similar to the matrix $\CQ$ constructed in \eqref{d4}. We assume that the homoclinic trajectories $\gamma_{k}$ are labeled in such a way that $m_{k}$ is non-increasing. Let $K_{1} \in \{ 1 , \ldots , K \}$ be such that $m_{1} = \cdots = m_{K_{1}} > m_{K_{1} + 1}$. Thus, the $K_{1}$ first trajectories $\gamma_{k}$ are the curves along which $\Lambda_{-}$ and $\Lambda_{+}$ are the ``most tangential''. If $m_{1} = 1$, then $K_{1} = K$, all the intersections are transversal (i.e. \ref{h8} holds true) and the matrix $\CQ ( z , h)$ is given by \eqref{d4}. If $m_{1} > 1$, then the entries of the $K_{1} \times K_{1}$ matrix $\CQ ( z , h )$ are given by
\begin{align}
\CQ_{k , \ell} ( z , h ) = e^{i A_{k} / h} \vert \alpha_{k}^{\infty} & \vert^{-\frac{1}{1 + m_{1}}} \frac{\sqrt{\lambda_{1} \lambda_{2}}}{\pi ( 1 + m_{1} )} \Gamma \Big( \frac{1}{1 + m_{1}} \Big) \Gamma \big( S ( z , h ) / \lambda_{1} \big) \frac{\CM_{k}^{+}}{\CM_{k}^{-}} e^{- \frac{\pi}{2} ( \nu_{k} + 1 ) i} \nonumber \\
&\times \big\vert g_{-}^{\ell} \big\vert^{\frac{\lambda_{1} + \lambda_{2}}{\lambda_{1}}} \big( i \lambda_{1} g_{+}^{k} \cdot g_{-}^{\ell} \big)^{- S ( z , h ) / \lambda_{1}} \left\{ \begin{aligned}
&e^{i \frac{\pi}{2 + 2 m_{1}} \sgn ( \alpha_{k}^{\infty} )}  &&\text{ for odd } m_{1} ,  \\
&\cos \Big( \frac{\pi}{2 + 2 m_{1}} \Big) &&\text{ for even } m_{1} ,
\end{aligned} \right.   \label{g23}
\end{align}
where $\CM_{k}^{+}$ is defined in \eqref{d7} and $\CM_{k}^{-}$ is the limit
\begin{equation} \label{g44}
\CM_{k}^{-} = \lim_{s \to + \infty} \sqrt{\Big\vert \det \frac{\partial x_{k} ( t , y )}{\partial ( t , y )} \vert_{t = s , \ y = 0} \Big\vert} e^{s \sum_{j} \lambda_{j} / 2} ,
\end{equation}
which belongs to $] 0 , + \infty [$. The notation $x_{k} ( t , y )$ is defined above \eqref{d7}. Note that the ``original'' $\CM_{k}^{-}$ (i.e. the one given by \eqref{d7}) vanishes in the present situation.

Mimicking Section \ref{s61}, the quantization rule for the pseudo-resonances is given by

\begin{definition}\sl \label{g26}
We say that $z$ is a pseudo-resonance if and only if
\begin{equation*}
1 \in \spe \big( h^{S ( z , h ) / \lambda_{1} - \frac{m_{1}}{1 + m_{1}}} \CQ ( z , h ) \big) .
\end{equation*}
The set of pseudo-resonances is denoted by $\res_{0} ( P )$.
\end{definition}

The factor $h^{S ( z , h ) / \lambda_{1} - 1 / 2}$ appearing in Definition \ref{d1} is replaced by $h^{S ( z , h ) / \lambda_{1} - \frac{m_{1}}{1 + m_{1}}}$. This is the crucial factor in the quantization rule because it determines the imaginary part of the resonances at leading order. In the present context, let $\mu_{1} ( \tau , h ) , \ldots , \mu_{K_{1}} ( \tau , h )$ denote the eigenvalues of
\begin{equation*}
\CQ \Big( E_{0} + h \tau - i h \Big( \frac{1 - m_{1}}{2 + 2 m_{1}} \lambda_{1} + \frac{\lambda_{2}}{2} \Big) , h \Big) .
\end{equation*}
The pseudo-resonances satisfy the following two scale asymptotic. We omit the proof of this result since it is similar to the one of Proposition \ref{d9}.

\begin{proposition}\sl \label{g27}
Assume $n = 2$, \ref{h1}--\ref{h4}, \ref{h13} and let $C > 0$. The pseudo-resonances $z$ lying in
\begin{equation} \label{g28}
E_{0} + [ - C h , C  h ] + i \Big[ - \Big( \frac{1 - m_{1}}{2 + 2 m_{1}} \lambda_{1} + \frac{\lambda_{2}}{2} \Big) h - C \frac{h}{\vert \ln h \vert} , h \Big] ,
\end{equation}
satisfy $z = z_{q , k} ( \tau ) + o ( h \vert \ln h \vert^{- 1} )$ in the sense of Proposition \ref{d9} with
\begin{equation*}
z_{q , k} ( \tau ) = E_{0} + 2 q \pi \lambda_{1} \frac{h}{\vert \ln h \vert} - i h \Big( \frac{1 - m_{1}}{2 + 2 m_{1}} \lambda_{1} + \frac{\lambda_{2}}{2} \Big) + i \ln ( \mu_{k} ( \tau , h ) ) \lambda_{1} \frac{h}{\vert \ln h \vert} ,
\end{equation*}
for some $q \in \Z$ and $k \in \{ 1 , \ldots , K_{1} \}$.
\end{proposition}

With the notation of Definition \ref{g80}, the following result shows that the resonances are close to the pseudo-resonances.

\begin{theorem}\sl \label{g24}
Assume $n = 2$, \ref{h1}--\ref{h4}, \ref{h13} and let $C , \delta > 0$. In the domain
\begin{equation} \label{g25}
E_{0} + [ - C h , C  h ] + i \Big[ - \Big( \frac{1 - m_{1}}{2 + 2 m_{1}} \lambda_{1} + \frac{\lambda_{2}}{2} \Big) h - C \frac{h}{\vert \ln h \vert} , h \Big] \setminus \big( \Gamma (h) + B ( 0 , \delta h ) \big),
\end{equation}
we have
\begin{equation*}
\dist \big( \res (P) , \res_{0} (P) \big) = o \Big( \frac {h}{\vert \ln h \vert} \Big) ,
\end{equation*}
as $h$ goes to $0$. Moreover, for all $\chi \in C^{\infty}_{0} ( \R^{n} )$, there exists $M > 0$ such that
\begin{equation*}
\big\Vert \chi ( P -z )^{-1} \chi \big\Vert \lesssim h^{- M} ,
\end{equation*}
uniformly for $h$ small enough and $z \in \eqref{g25}$ with $\dist ( z , \res_{0} ( P ) ) \geq \delta h \vert \ln h \vert^{- 1}$.
\end{theorem}

This result shows that the asymptotic of the resonances generated by tangential intersections of finite order has the same structure than the one in the transversal setting (see Theorem \ref{d8}). Thus the resonances verify again the phenomena described in Section \ref{s18}: accumulation on curves, vibration, transition and stability. We do not develop more this point here.

\begin{remark}\sl \label{k62}
$i)$ At the leading order, the imaginary part of the resonances closest to the real axis behaves like
\begin{equation} \label{g71}
\im z \approx - \Big( \frac{\lambda_{1}}{1 + m_{1}} - \frac{\lambda_{1}}{2} + \frac{\lambda_{2}}{2} \Big) h ,
\end{equation}
except under some special circumstances (see Example \ref{e48}). Then, the resonances get closer to the real axis when $m_{1}$ increases. This is in agreement with the intuition that the more tangential the intersection is, the more trapping the situation is. Letting $m_{1}$ goes to $+ \infty$, \eqref{g71} becomes $\im z \approx ( \lambda_{1} - \lambda_{2} ) h /2$ which can be compared to \eqref{c14}.

$ii)$ The same phenomena can be observed for $m_{1}$ fixed. Assume for simplicity that $K = 1$ and $m_{1} \geq 2$. Theorem \ref{g24} shows that resonances move away from the real axis when $\vert \alpha_{1}^{\infty} \vert$ increases, the other parameters being fixed. This is natural since the intersections is less tangential and then the trapping is weaker when $\vert \alpha_{1}^{\infty} \vert$ increases.
\end{remark}

Only the $K_{1}$ most tangential trajectories contribute to $\CQ$ and then to the asymptotic of the resonances. This is an illustration of the stability principle developed in Section \ref{s24}. Indeed, these $K_{1}$ trajectories are the most trapped ones and the other $K - K_{1}$ curves must be seen as lower order perturbations. By comparison with Proposition \ref{i38}, the difference of strength of trapping is so high here that these $K - K_{1}$ trajectories do not even contribute in $\CQ$ (whereas they give small terms in \eqref{i48}). This phenomenon concerns the leading terms in the asymptotic of the resonances, but probably not more precise asymptotics. It means that a total quantization operator $\CQ_{\rm tot}$ similar to the one of Section \ref{s79} should take into account the contributions of all the $K$ homoclinic trajectories.

\begin{remark}\sl
One can probably also deal with contacts of infinite order between $\Lambda_{-}$ and $\Lambda_{+}$. The result will then depend on the nature of the contact. If $\lambda_{1} = \lambda_{2}$, the imaginary part of the resonances should verify $h \vert \ln h \vert^{- 1} \ll \im z \ll h$. Indeed, the situation is intermediate between the finite order case of Theorem \ref{g24} and the fully tangential case of Theorem \ref{i55}.
\end{remark}

One may ask why the coefficient $\alpha_{k}^{\infty}$ appears in the present quantization operator \eqref{g23} and not in the one of the transversal case \eqref{d4}. In fact, for tangential intersection of finite order, the nature of the contact between $\Lambda_{-}^{0}$ and $\Lambda_{+}^{1}$ is mainly described by \eqref{g21} where $\alpha_{k}$ is present. On the contrary, Proposition C.1 of \cite{ALBoRa08_01} shows that the asymptotic (as $x_{1}$ goes to $0$) of this contact is universal and uniformly transversal under the assumption \ref{h8}. In other words, in the transversal case, the corresponding coefficient $\alpha_{k}$ will only measure the difference between $\Lambda_{+}^{1}$ and its asymptote which is transversal to $\Lambda_{-}$. Thus, it is natural that it does not appear at the principal order (i.e. in the matrix $\CQ$ defined in \eqref{d4}) but would probably appear at lower orders (i.e. in the matrix $\CQ_{\rm tot}$ defined in \eqref{e95}).

When $\lambda_{1} = \lambda_{2}$, the assumptions of Section \ref{s77} or Section \ref{s53} are satisfied and $\CH_{\rm tang} = \bigcup_{m_{k} \geq 2} \gamma_{k}$ has Lebesgue measure zero in $\S^{1}$. Then, we already known from Theorem \ref{a2} (or Remark \ref{b69} $i)$ to be more precise) that $P$ has no resonance in the set $[ E_{0} - C h , E_{0} + C h ] + i [ - C h \vert \ln h \vert^{- 1} , 0 ]$ for any $C > 0$. One can also compare the notations in the quantization operators \eqref{b70}, \eqref{g23} and \eqref{c19} (see Remark \ref{d3} $iii)$). A direct computation shows
\begin{align}
\frac{\CM_{k}^{+}}{\CM_{k}^{-}} \big\vert g_{-}^{\ell} \big\vert^{\frac{\lambda_{1} + \lambda_{2}}{\lambda_{1}}} \big( i & \lambda_{1} g_{+}^{k} \cdot g_{-}^{\ell} \big)^{- S ( z , h ) / \lambda_{1}}    \nonumber \\
&= \CM_{0} ( \widehat{g}_{+}^{k} ) e^{i T ( \widehat{g}_{+}^{k} ) \frac{z - E_{0}}{h}} \big( i \lambda_{1} \widehat{g}_{+}^{k} \cdot \widehat{g}_{-}^{\ell} \big)^{- S ( z , h ) / \lambda_{1}} \frac{\vert g_{-}^{\ell} \vert^{\frac{\lambda_{1} + \lambda_{2}}{2 \lambda_{1}} + i \frac{z - E_{0}}{\lambda_{1} h}}}{\vert g_{-}^{k} \vert^{\frac{\lambda_{1} + \lambda_{2}}{2 \lambda_{1}} + i \frac{z - E_{0}}{\lambda_{1} h}}} , \label{g72}
\end{align}
with $\widehat{g}_{\pm}^{\bullet} : = g_{\pm}^{\bullet} / \vert g_{\pm}^{\bullet} \vert$. The last factor in the right hand side of \eqref{g72} plays no role in the asymptotic of the resonances and can formally be replaced by $1$. Indeed, only the eigenvalues of $\CQ$ are relevant (see Remark \ref{d3} $i)$). Thus, this part of the quantization operators looks quite similar. Nevertheless, the settings are rather different and this analogy can naturally not be complete: the exponents of $h$ in the quantization rules are not identical, some numerical factors change and the terms containing $m_{1}$ only appear in \eqref{g23}.

\begin{figure}
\begin{center}
\begin{picture}(0,0)%
\includegraphics{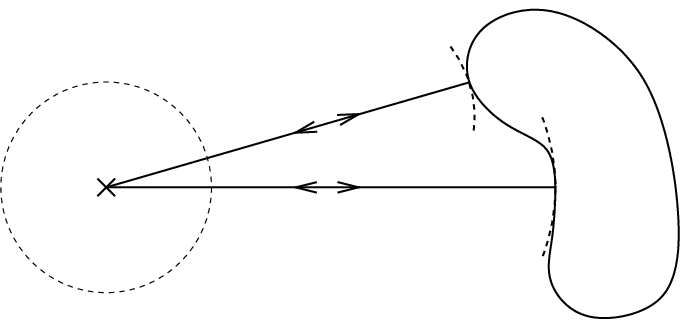}%
\end{picture}%
\setlength{\unitlength}{1105sp}%
\begingroup\makeatletter\ifx\SetFigFont\undefined%
\gdef\SetFigFont#1#2#3#4#5{%
  \reset@font\fontsize{#1}{#2pt}%
  \fontfamily{#3}\fontseries{#4}\fontshape{#5}%
  \selectfont}%
\fi\endgroup%
\begin{picture}(11670,5352)(-10220,-7435)
\put(-4574,-3586){\makebox(0,0)[b]{\smash{{\SetFigFont{9}{10.8}{\rmdefault}{\mddefault}{\updefault}$\pi_{x} ( \gamma_{2} )$}}}}
\put(-8399,-4786){\makebox(0,0)[b]{\smash{{\SetFigFont{9}{10.8}{\rmdefault}{\mddefault}{\updefault}$0$}}}}
\put(-674,-3361){\makebox(0,0)[b]{\smash{{\SetFigFont{9}{10.8}{\rmdefault}{\mddefault}{\updefault}$\CO$}}}}
\put(-8399,-2761){\makebox(0,0)[b]{\smash{{\SetFigFont{9}{10.8}{\rmdefault}{\mddefault}{\updefault}$\supp V$}}}}
\put(-4574,-5761){\makebox(0,0)[b]{\smash{{\SetFigFont{9}{10.8}{\rmdefault}{\mddefault}{\updefault}$\pi_{x} ( \gamma_{1} )$}}}}
\end{picture} $\qquad \qquad$ \begin{picture}(0,0)%
\includegraphics{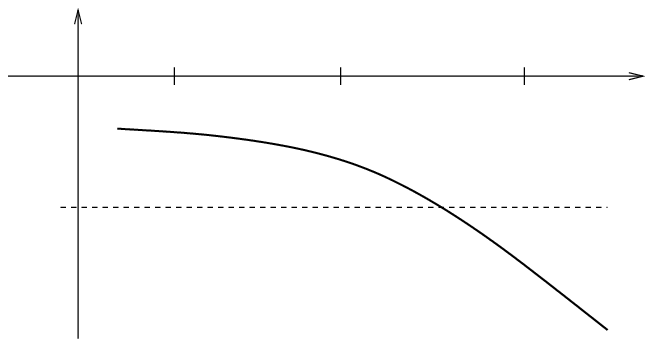}%
\end{picture}%
\setlength{\unitlength}{1105sp}%
\begingroup\makeatletter\ifx\SetFigFont\undefined%
\gdef\SetFigFont#1#2#3#4#5{%
  \reset@font\fontsize{#1}{#2pt}%
  \fontfamily{#3}\fontseries{#4}\fontshape{#5}%
  \selectfont}%
\fi\endgroup%
\begin{picture}(11212,5744)(961,-5708)
\put(976,-3511){\makebox(0,0)[b]{\smash{{\SetFigFont{9}{10.8}{\rmdefault}{\mddefault}{\updefault}$- \frac{\lambda_{1}}{1 + m_{1}} h$}}}}
\put(10051,-811){\makebox(0,0)[b]{\smash{{\SetFigFont{9}{10.8}{\rmdefault}{\mddefault}{\updefault}$E_{0} + C h$}}}}
\put(4051,-811){\makebox(0,0)[b]{\smash{{\SetFigFont{9}{10.8}{\rmdefault}{\mddefault}{\updefault}$E_{0} - C h$}}}}
\put(6901,-811){\makebox(0,0)[b]{\smash{{\SetFigFont{9}{10.8}{\rmdefault}{\mddefault}{\updefault}$E_{0}$}}}}
\end{picture}%
\end{center}
\caption{An example of operator satisfying \ref{h13} with $m_{1} > m_{2}$ and the corresponding accumulation curve.} \label{f32}
\end{figure}

\begin{example}\rm \label{g77}
We end this discussion by the construction of operators illustrating Theorem \ref{g24}. As in previous examples, they will not be of the form \eqref{a5} but will enter in the setting of Remark \ref{c13}.

In dimension $n = 2$, we consider the semiclassical operator
\begin{equation*}
P = - h^{2} \Delta_{\R^{2} \setminus \CO} + V (x) ,
\end{equation*}
with Dirichlet condition at the boundary of the obstacle $\CO$. Here, $V (x) \in C^{\infty}_{0} ( \R^{2} )$ is a potential which can be chosen as in Section \ref{s21}. In particular, we have $\lambda_{1} = \lambda_{2}$ (for examples with $\lambda_{1} < \lambda_{2}$, one can consider potentials $V$ as in Example \ref{e22}). On the other hand, $\CO$ is a small non-trapping obstacle far away from the origin as illustrated in Figure \ref{f32}.

As explained in Example \ref{e29}, the assumptions \ref{h1}--\ref{h4} are satisfied and the homoclinic trajectories are the radial rays which are normal to the boundary of $\partial \CO$. For $x \in \partial \CO \cap \pi_{x} ( \CH )$, one can verify that the order of the contact of the manifolds $\Lambda_{-}$ and $\Lambda_{+}$ along $[ 0 , x ]$ coincides with the order of contact of the curves $\partial B ( 0 , x )$ and $\partial \CO$ at $x$. Thus, one can easily construct operators satisfying \ref{h13}.

We precise the asymptotic of the resonances given by Theorem \ref{g24} in the simplest setting: when $\gamma_{1}$ is the unique homoclinic trajectory of order $m_{1}$ (i.e. $K_{1} = 1$ or $m_{1} > m_{2}$ if $K \geq 2$). Thus, the situation is similar to the one of the assumption \ref{h9} and the distribution of resonances will be quite close to \eqref{e8}. The symmetry of the Hamiltonian trajectories implies that $g_{-}^{1}$ and $g_{+}^{1}$ are the same up to some positive constant. Then, Proposition \ref{g27} gives
\begin{equation} \label{g79}
z_{q , 1} ( \tau ) = E_{0} - \frac{A_{1} \lambda_{1}}{\vert \ln h \vert} + 2 q \pi \lambda_{1} \frac{h}{\vert \ln h \vert} - i h \frac{\lambda_{1}}{1 + m_{1}} + i \ln ( \mu ( \tau ) ) \lambda_{1} \frac{h}{\vert \ln h \vert} ,
\end{equation}
with $q \in \Z$ and
\begin{align*}
\mu ( \tau ) : = \vert \alpha_{1}^{\infty} \vert^{-\frac{1}{1 + m_{1}}} & \frac{\lambda_{1}}{\pi ( 1 + m_{1} )} \Gamma \Big( \frac{1}{1 + m_{1}} \Big) \Gamma \Big( \frac{m_{1}}{1 + m_{1}} - i \frac{\tau}{\lambda_{1}} \Big) \frac{\CM_{1}^{+}}{\CM_{1}^{-}} e^{- \frac{\pi}{2} ( \nu_{1} + 1 ) i} \nonumber \\
&\times \big\vert g_{-}^{1} \big\vert^{2} \big( i \lambda_{1} \vert g_{+}^{1} \vert \vert g_{-}^{1} \vert \big)^{- \frac{m_{1}}{1 + m_{1}} + i \frac{\tau}{\lambda_{1}}} \left\{ \begin{aligned}
&e^{i \frac{\pi}{2 + 2 m_{1}} \sgn ( \alpha_{1}^{\infty} )}  &&\text{ for odd } m_{1} ,  \\
&\cos \Big( \frac{\pi}{2 + 2 m_{1}} \Big) &&\text{ for even } m_{1} .
\end{aligned} \right.
\end{align*}
The localization of the resonances is then given by Theorem \ref{g24} in the domain \eqref{g28}. We now study the unique corresponding accumulation curve in the sense of Remark \ref{e2}. For any fixed $x \in \R$, \cite[(5.11.9)]{OlLoBoCl10_01} shows that
\begin{equation} \label{k45}
\vert \Gamma ( x + i y ) \vert \sim \sqrt{2 \pi} \vert y \vert^{x - \frac{1}{2}} e^{- \frac{\pi \vert y \vert}{2}} .
\end{equation}
as $y$ goes to $\infty$. Combining with \eqref{g79}, the asymptotic behavior of the accumulation curve for large $\re \sigma$ satisfies
\begin{equation} \label{g78}
\im \sigma = - \frac{\lambda_{1}}{1 + m_{1}} + \frac{1}{\vert \ln h \vert}
\left\{ \begin{aligned}
&A_{1} \ln \vert \re \sigma \vert + A_{0} + o ( 1 ) &&\text{ as } \re \sigma \to - \infty , \\
&- \pi \re \sigma + A_{1} \ln \vert \re \sigma \vert + A_{0} + o ( 1 )  &&\text{ as } \re \sigma \to + \infty , \\
\end{aligned} \right.
\end{equation}
with the constants $A_{1} = \lambda_{1} \frac{m_{1} - 1}{2 + 2 m_{1}} > 0$, $A_{0} = \lambda_{1} \ln B_{0}$ and
\begin{align*}
B_{0} = \vert \alpha_{1}^{\infty} & \vert^{-\frac{1}{1 + m_{1}}} \sqrt{\frac{2}{\pi}} \frac{\lambda_{1}^{\frac{3 - m_{1}}{2 + 2 m_{1}}}}{1 + m_{1}} \Gamma \Big( \frac{1}{1 + m_{1}} \Big) \frac{\CM_{1}^{+}}{\CM_{1}^{-}} \frac{\vert g_{-}^{1} \vert^{\frac{2 + m_{1}}{1 + m_{1}}}}{\vert g_{+}^{1} \vert^{\frac{m_{1}}{1 + m_{1}}}} \left\{ \begin{aligned}
&1  &&\text{ for odd } m_{1} ,  \\
&\cos \Big( \frac{\pi}{2 + 2 m_{1}} \Big) &&\text{ for even } m_{1} .
\end{aligned} \right.
\end{align*}
By comparison with Corollary \ref{e4}, the novelty is the term $A_{1} \ln \vert \re \sigma \vert$. When $\re \sigma \to + \infty$, this term is secondary and the accumulation curve behaves essentially as in \eqref{e3}, mutatis mutandis. This is again an illustration of the transition to non-trapping principle described in Section \ref{s14}. Indeed, the energies above $E_{0}$ seems ``less trapping'' than $E_{0}$ (even if they can be trapping) since the asymptotic directions $g_{\pm}^{1}$ are on the same side of $0$. On the other hand, when $\re \sigma \to - \infty$, the new term $A_{1} \ln \vert \re \sigma \vert$ is dominant. Thus, the accumulation curve approaches slowly to the real axis. That the energies below $E_{0}$ seems more trapping than in the transversal case could explain this convergence. We send the reader to Section \ref{s14} for a general discussion on transition phenomena under the assumption \ref{h8}.
\end{example}

\section{Asymptotic of the resonances generated by nappes of homoclinic trajectories} \label{s26}

\Subsection{Main results} \label{s53}

We now give the asymptotic of the resonances when the intersection between $\Lambda_{-} $ and $\Lambda_{+}$ is maximal. Since we want to consider a strong trapping situation (as in Theorem \ref{a2}) and not a weak trapping situation (as in Theorem \ref{a1} and Section \ref{s6}), we come back to the assumptions and notations of Section \ref{s77}. Thus, we expect to see resonances at distance $h \vert \ln h \vert^{- 1}$ from the real axis.

We make the assumption \ref{h6}. This means that all the $\lambda_{j}$'s coincide and are denoted by $\lambda$ in the sequel. We send back the reader to Section \ref{s77} for the definition of the sets $\CH_{\rm tang} , \CH_{\rm tang}^{\pm \infty}$ and of the homeomorphisms $\alpha ( \cdot ) , \omega ( \cdot )$. We also define $\CH_{\rm trans} = \CH \setminus \CH_{\rm tang}$, the set of homoclinic trajectories along which $\Lambda_{-}$ and $\Lambda_{+}$ intersect transversally in at least one direction.

In order to give the asymptotic of the resonances, we need to define some dynamical quantities. The proofs of the following geometric assertions can be found in Section \ref{a71}. For $\alpha \in \CH^{+ \infty}_{\rm tang}$, let $\gamma_{\alpha}$ denote the Hamiltonian trajectory given by \eqref{c17}. Recall that the quantity $\CM_{0} ( \alpha )$, which measures the amplification along homoclinic trajectories, is defined in \eqref{i49}. The action
\begin{equation*}
A ( \alpha ) = \int_{\gamma_{\alpha}} \xi \cdot d x ,
\end{equation*}
and the Maslov's index of $\Lambda_{+}$ along the curve $\gamma_{\alpha}$, noted $\nu ( \alpha )$, are well-defined. Moreover, $A$ is continuous on $\CH^{+ \infty}_{\rm tang}$ and $\nu$ is locally constant. Eventually, let $T ( \alpha )$ denote the limit
\begin{equation} \label{i50}
T ( \alpha ) : = \lim_{\varepsilon \to 0} \big( t_{-}^{\varepsilon} ( \alpha ) - t_{+}^{\varepsilon} ( \alpha ) - 2 \vert \ln \varepsilon \vert / \lambda \big) ,
\end{equation}
which exists for all $\alpha \in \CH_{\rm tang}^{+ \infty}$ (the $t_{\pm}^{\varepsilon} ( \alpha )$'s are defined  below \eqref{a81}). In some sense, $T ( \alpha )$ can be seen as the time delay on the trajectory $\gamma_{\alpha}$. This function is continuous on $\CH_{\rm tang}^{+ \infty}$. In this part, we assume
\begin{hyp} \label{h7}
$\CH^{+ \infty}_{\rm tang}$ is the union of a finite number of isolated sets on which $A$ and $\nu$ are constant.
\end{hyp}
This assumption holds true in many situations. Indeed, the action and the Maslov's index are constant on any smooth path in $\CH^{+ \infty}_{\rm tang}$ since $\Lambda_{+}$ is a Lagrangian manifold. Thus, \ref{h7} is satisfied in Example \ref{b80}, in dimension $n = 1$ (see Section \ref{s19}) and in the examples of Section \ref{s54}. One should be able to relax this hypothesis by assuming that \ref{h7} holds outside a set of measure zero.

We can now define the quantization operator which governs the localization of the resonances in this strong trapping case. Let $\CT$ denote the operator on $L^{2} ( \CH^{- \infty}_{\rm tang} )$, endowed with the Lesbegue measure on $\S^{n-1}$, with kernel
\begin{align}
\CT ( z , h ) ( \omega , \widetilde{\omega} ) = e^{i A ( \alpha ( \omega ) ) / h} \Gamma & \Big( \frac{n}{2} - i \frac{z - E_{0}}{\lambda h} \Big) \Big( \frac{\lambda}{2 \pi} \Big)^{\frac{n}{2}} \CM_{0} ( \alpha ( \omega ) ) \nonumber  \\
&\times e^{- i ( \nu ( \alpha ( \omega ) ) \frac{\pi}{2} + n \frac{\pi}{4} - T ( \alpha ( \omega ) ) \frac{z - E_{0}}{h} )} \big( i \lambda \alpha ( \omega ) \cdot \widetilde{\omega} \big)^{- \frac{n}{2} + i \frac{z - E_{0}}{\lambda h}} .   \label{c19}
\end{align}
From \ref{h4}, the compactness of $\CH^{- \infty}_{\rm tang}$ and the previous properties on the geometric quantities, this kernel is continuous and then bounded on $\CH^{- \infty}_{\rm tang} \times \CH^{- \infty}_{\rm tang}$. In particular, $\CT ( z , h )$ is a compact operator. As noticed above Remark \ref{b72}, we could have defined equivalently $\CT$ on any $L^{q} ( \CH^{- \infty}_{\rm tang} )$ with $1 \leq q \leq + \infty$. We rather work on $L^{2} ( \CH^{- \infty}_{\rm tang} )$ in this part. From \ref{h7}, the action $A$ takes a finite number of values, say $A_{1} , \ldots , A_{K}$. Thus, as the matrix $\CQ$ in \eqref{d10}, the operator $\CT$ can be decomposed as
\begin{equation} \label{i81}
\CT ( z , h ) = \widetilde{\CT} ( \rho , \sigma ) = \sum_{k = 1}^{K} \rho_{k} \widetilde{\CT}_{k} ( \sigma ) , 
\end{equation}
where $\rho = ( \rho_{1} , \ldots , \rho_{K} )$ with $\rho_{k} = e^{i A_{k} / h}$ and $\widetilde{\CT}_{k}$ is independent of $h$ and meromorphic on the rescaled parameter $\sigma$ defined in \eqref{d92}. In this strong trapping case, the pseudo-resonances are defined by the

\begin{definition}[Quantization rule]\sl \label{i51}
We say that $z$ is a pseudo-resonance if and only if
\begin{equation*}
1 \in \spe \big( h^{- i \frac{z - E_{0}}{\lambda h}} \CT ( z , h ) \big) .
\end{equation*}
The set of pseudo-resonances is denoted by $\res_{0} ( P )$.
\end{definition}

This quantization rule looks similar to those of the transversal case (see Definition \ref{d1}) and of the tangential intersection of finite order (see Definition \ref{g26}). Nevertheless, there are some differences. First, the most important part, that is the power of $h$, changes. Second, some coefficients of $\CQ$ and $\CT$ do not have the same structure (the Maslov determinants $\CM$ of \eqref{i49} and \eqref{d7} for instance). Nevertheless, the two quantization rules are the same in dimension $n = 1$ (see Lemma \ref{i82} for the proof). This is natural since the assumptions of Section \ref{s6} and Section \ref{s26} are satisfied simultaneously when $n = 1$.

We now give the asymptotic of the pseudo-resonances which are implicitly defined. For all $\tau \in \R$, the operator
\begin{equation*}
\widehat{\CT} ( \tau , h ) : = \CT ( E_{0} + h \tau , h ) ,
\end{equation*}
is compact on $L^{2} ( \CH^{- \infty}_{\rm tang} )$. Its kernel is given by
\begin{align*}
\widehat{\CT} ( \tau , h ) ( \omega , \widetilde{\omega} ) = e^{i A ( \alpha ( \omega ) ) / h} \Gamma \Big( \frac{n}{2} - i \frac{\tau}{\lambda} & \Big) \Big( \frac{\lambda}{2 \pi} \Big)^{\frac{n}{2}} \CM_{0} ( \alpha ( \omega ) )   \\
&\times e^{- i ( \nu ( \alpha ( \omega ) ) \frac{\pi}{2} + n \frac{\pi}{4} - T ( \alpha ( \omega ) ) \tau )} \big( i \lambda \alpha ( \omega ) \cdot \widetilde{\omega} \big)^{- \frac{n}{2} + i \frac{\tau}{\lambda}} .
\end{align*}
Let $\mu_{1} ( \tau , h ) , \mu_{2} ( \tau , h ) , \ldots$ denote the non-zero eigenvalues of $\widehat{\CT} ( \tau , h )$ counted with their multiplicity. Depending on the situation, they can be in finite or infinite number and their cardinal may depend on $\tau , h$. Since $[- C , C ] \times ] 0 , 1 ] \ni ( \tau , h ) \mapsto \widehat{\CT} ( \tau , h )$ is analytic, the general perturbation theory (see Chapter VII of Kato \cite{Ka76_01}) shows that these non-zero eigenvalues (correctly labeled) are locally continuous functions of $\tau , h$ and locally analytic functions on $\tau$ for $h$ fixed with only algebraic singularities at some exceptional points. Eventually, \eqref{c19} yields that $( \tau , h ) \mapsto \widehat{\CT} ( \tau , h )$ is a locally uniformly bounded function of Hilbert--Schmidt operators. This means that
\begin{equation*}
\sup_{( \tau , h ) \in [- C , C ] \times ] 0 , 1 ]} \big\Vert \widehat{\CT} ( \tau , h ) \big\Vert_{\rm H S} < + \infty ,
\end{equation*}
for all $C , \delta > 0$. Combining with the Bienaym\'e--Tchebychev's inequality, it ensures
\begin{equation} \label{i52}
\sup_{( \tau , h ) \in [- C , C ] \times ] 0 , 1 ]} \card \{ k ; \ \vert \mu_{k} ( \tau , h ) \vert \geq \delta \} < + \infty ,
\end{equation}
for all $C , \delta > 0$. The pseudo-resonances of Definition \ref{i51} verify an asymptotic similar to the one of Proposition \ref{d9}.

\begin{proposition}[Asymptotic of the pseudo-resonances]\sl \label{i53}
Assume \ref{h1}--\ref{h4}, \ref{h6}, \ref{h7} and let $C > 0$. The pseudo-resonances $z$ lying in
\begin{equation} \label{i54}
E_{0} + [ - C h , C h ] + i \Big[ - C \frac{h}{\vert \ln h \vert} , \frac{h}{\vert \ln h \vert} \Big] ,
\end{equation}
satisfy $z = z_{q , k} ( \tau ) + o ( h \vert \ln h \vert^{- 1} )$ in the sense of Proposition \ref{d9} with
\begin{equation*}
z_{q , k} ( \tau ) = E_{0} + 2 q \pi \lambda \frac{h}{\vert \ln h \vert} + i \ln ( \mu_{k} ( \tau , h ) ) \lambda \frac{h}{\vert \ln h \vert} ,
\end{equation*}
for some $q \in \Z$ and $k$.
\end{proposition}

Note that only the eigenvalues $\mu_{k} ( \tau , h )$'s with $\vert \mu_{k} ( \tau , h ) \vert \geq e^{- 2 C / \lambda}$ provide pseudo-resonances in the domain \eqref{i54}. From \eqref{i52}, the numbers of those eigenvalues is then uniformly bounded with respect to $\tau , h$. With Definition \ref{g80} in mind, the asymptotic of the resonances in this strong trapping case is given by the following result.

\begin{figure}
\begin{center}
\begin{picture}(0,0)%
\includegraphics{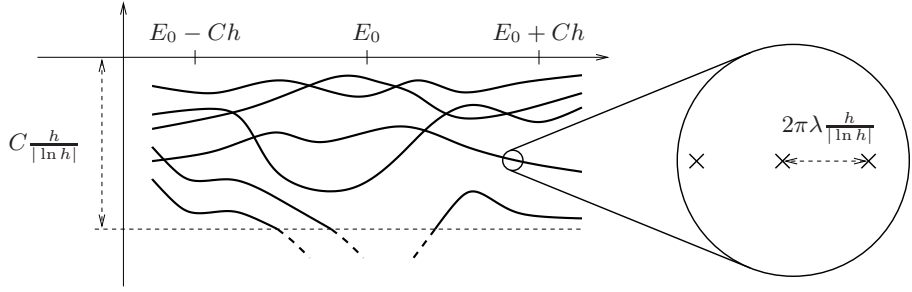}%
\end{picture}%
\setlength{\unitlength}{1184sp}%
\begingroup\makeatletter\ifx\SetFigFont\undefined%
\gdef\SetFigFont#1#2#3#4#5{%
  \reset@font\fontsize{#1}{#2pt}%
  \fontfamily{#3}\fontseries{#4}\fontshape{#5}%
  \selectfont}%
\fi\endgroup%
\begin{picture}(18891,6044)(-1214,-5783)
\put(-1199,-2911){\makebox(0,0)[lb]{\smash{{\SetFigFont{9}{10.8}{\rmdefault}{\mddefault}{\updefault}$C \frac{h}{\vert \ln h \vert}$}}}}
\put(15976,-2536){\makebox(0,0)[b]{\smash{{\SetFigFont{9}{10.8}{\rmdefault}{\mddefault}{\updefault}$2 \pi \lambda \frac{h}{\vert \ln h \vert}$}}}}
\put(6301,-586){\makebox(0,0)[b]{\smash{{\SetFigFont{9}{10.8}{\rmdefault}{\mddefault}{\updefault}$E_{0}$}}}}
\put(2701,-586){\makebox(0,0)[b]{\smash{{\SetFigFont{9}{10.8}{\rmdefault}{\mddefault}{\updefault}$E_{0} - C h$}}}}
\put(9901,-586){\makebox(0,0)[b]{\smash{{\SetFigFont{9}{10.8}{\rmdefault}{\mddefault}{\updefault}$E_{0} + C h$}}}}
\end{picture}%
\end{center}
\caption{Two scale asymptotic of resonances of Theorem \ref{i55}.} \label{f41}
\end{figure}

\begin{theorem}[Asymptotic of resonances]\sl \label{i55}
Assume \ref{h1}--\ref{h4}, \ref{h6}, \ref{h7} and let $C , \delta > 0$. In the domain \eqref{i54}, we have
\begin{equation*}
\dist \big( \res (P) , \res_{0} (P) \big) = o \Big( \frac {h}{\vert \ln h \vert} \Big) ,
\end{equation*}
as $h$ goes to $0$. Moreover, for all $\chi \in C^{\infty}_{0} ( \R^{n} )$, there exists $M > 0$ such that
\begin{equation*}
\big\Vert \chi ( P -z )^{-1} \chi \big\Vert \lesssim h^{- M} ,
\end{equation*}
uniformly for $h$ small enough and $z \in \eqref{i54}$ with $\dist ( z , \res_{0} ( P ) ) \geq \delta h \vert \ln h \vert^{- 1}$.
\end{theorem}

Thus, under the previous assumptions, the resonances are typically at distance $h \vert \ln h \vert^{- 1}$ from the real axis. However, this is not always the case. For instance, the imaginary part of the resonances can be exponential small in dimension $n = 1$ (see case {\rm (C)} in Section \ref{s19} for more details). On the contrary, there may be no resonance at distance $h \vert \ln h \vert^{- 1}$ from the real axis. We have already seen in Remark \ref{b69} $i)$ that this is the case when $\mes_{\S^{n - 1}} ( \CH_{\rm tang}^{\pm \infty} ) = 0$. In this situation, Theorem \ref{i55} gives no additional information.

This theorem shows that the resonances are close to the pseudo-resonances, but gives no information on the multiplicity. However, as in Proposition \ref{d91}, it is possible to give a lower bound for the number of resonances close to a pseudo-resonance by the multiplicity of this pseudo-resonance.

The asymptotic of the resonances given by Proposition \ref{i53} and Theorem \ref{i55} is partially implicit. Indeed, the eigenvalues of $\widehat{\CT}$ have to be computed. In the figures below, they are calculated numerically by discretizing the kernel \eqref{c19} on the compact $\CH^{- \infty}_{\rm tang}$. The situation of Section \ref{s6} was slightly different since the eigenvalues of a $K \times K$ matrix can be explicitly computed (at least for $K < 5$ or in particular cases like Section \ref{s16}). Roughly speaking, Theorem \ref{i55} (but also Theorem \ref{d8}) can be seen as the ``quotient by the Hamiltonian flow'' of the resonance problem and the ``reduction'' to the asymptotic directions at the hyperbolic fixed point of the homoclinic set.

In the results stated in Section \ref{s6}, the spectral parameter must avoid a neighborhood of size $h$ of the exceptional set $\Gamma (h)$. It is not the case in Theorem \ref{i55} where this exceptional set does not appear. This is natural since $\Gamma ( h )$, being at at distance $h$ of the real axis, does not meet the domain \eqref{i54}.

From Theorem \ref{i55}, the operator $\CT$ gives the resonances in any neighborhood of size $h \vert \ln h \vert^{- 1}$ of the real axis. There is no hope to generalize this result to larger zone. For instance, consider an operator satisfying the assumptions of Section \ref{s6} in dimension $n \geq 2$ as in Example \ref{c15}, Example \ref{e6} or Example \ref{g77}. In these cases, $L^{2} ( \CH^{- \infty}_{\rm tang} ) = \{ 0 \}$, $\CT$ has no eigenvalue and Definition \ref{i51} (even generalized to larger zone) provides no pseudo-resonance. Nevertheless, Theorem \ref{d8} and Theorem \ref{g24} show that in these situations there are  a lot of resonances in neighborhoods of size $h$ of the real axis. Thus, the small eigenvalues of $\CT$ do not provide, in general, the asymptotic of the resonances deeper in the complex plane. In some sense, $\CT$ contains only the ``pertinent'' dynamical informations for the asymptotic of the resonances at distance $h \vert \ln h \vert^{- 1}$ of the real axis. Understanding the distribution of resonances in larger zones would require to take into account the nature of the contact between $\Lambda_{-}$ and $\Lambda_{+}$ along $\CH$. Under the present assumptions, this seems rather difficult.

\begin{remark}\sl
$i)$ The assumption \ref{h7} is probably a technical hypothesis. In the proof, it is used to estimate the resolvent of $\CT$ and to invert the limits $h \to 0$ and $\varepsilon \to 0$ (see e.g. \eqref{j17}). Without \ref{h7}, the quantization operator $\CT$ can still be defined by \eqref{c19}. It can no longer be written like \eqref{i81}, but can be seen as a Fourier integral operator. In particular, it is not clear that the pseudo-resonances defined by Definition \ref{i51} satisfy the asymptotic stated in Proposition \ref{i53} and that they are close to the resonances.

$ii)$ More generally, one may think to give the asymptotic of the resonances without any assumption on the nature of the intersection of $\Lambda_{-}$ and $\Lambda_{+}$. In particular, such result will unify Section \ref{s6} and Section \ref{s26}. The idea will be to defined an abstract quantization operator whose kernel will be defined on all $\S^{n - 1}$. Its construction will be based on Theorem \ref{a32} and on the monodromy in a neighborhood of $\CH$. Unfortunately, this classical operator is too complicated to deal with, the crucial estimates on its resolvent can not be obtained and the asymptotic of the resonances is therefore out of reach. 
\end{remark}

\begin{figure}
\begin{center}
\begin{picture}(0,0)%
\includegraphics{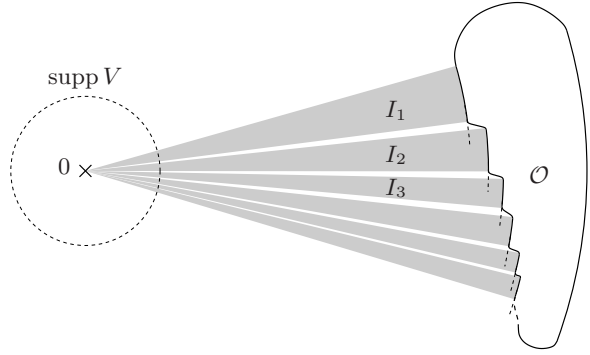}%
\end{picture}%
\setlength{\unitlength}{987sp}%
\begingroup\makeatletter\ifx\SetFigFont\undefined%
\gdef\SetFigFont#1#2#3#4#5{%
  \reset@font\fontsize{#1}{#2pt}%
  \fontfamily{#3}\fontseries{#4}\fontshape{#5}%
  \selectfont}%
\fi\endgroup%
\begin{picture}(14514,8766)(-10289,-9660)
\put(-8774,-5236){\makebox(0,0)[rb]{\smash{{\SetFigFont{9}{10.8}{\rmdefault}{\mddefault}{\updefault}$0$}}}}
\put(-8399,-2986){\makebox(0,0)[b]{\smash{{\SetFigFont{9}{10.8}{\rmdefault}{\mddefault}{\updefault}$\supp V$}}}}
\put(3001,-5461){\makebox(0,0)[b]{\smash{{\SetFigFont{9}{10.8}{\rmdefault}{\mddefault}{\updefault}$\CO$}}}}
\put(-599,-3811){\makebox(0,0)[b]{\smash{{\SetFigFont{9}{10.8}{\rmdefault}{\mddefault}{\updefault}$I_{1}$}}}}
\put(-599,-4936){\makebox(0,0)[b]{\smash{{\SetFigFont{9}{10.8}{\rmdefault}{\mddefault}{\updefault}$I_{2}$}}}}
\put(-599,-5761){\makebox(0,0)[b]{\smash{{\SetFigFont{9}{10.8}{\rmdefault}{\mddefault}{\updefault}$I_{3}$}}}}
\end{picture}%
\end{center}
\caption{A geometric setting where \ref{h7} does not hold.} \label{f44}
\end{figure}

Figure \ref{f44} provides a geometric situation where \ref{h7} is not satisfied. There are an infinite number of $I_{j} \subset \CH^{- \infty}_{\rm tang}$ with $\mes_{\S^{n - 1}} ( I_{j} ) > 0$ on which the action $A_{j}$ is constant and verifies $A_{j} \neq A_{k}$ for $j \neq k$. The corresponding operator can be rigorously constructed adapting the forthcoming Example \ref{i72} (see also Example \ref{i86}).

One can compare Theorem \ref{i55} with the resonance free domains obtained in Section \ref{s77}. Concerning the kernels, we have
\begin{equation} \label{i56}
\CT_{0} ( \tau ) ( \omega , \widetilde{\omega} ) = \big\vert \CT ( E_{0} + h \tau , h ) ( \omega , \widetilde{\omega} ) \big\vert ,
\end{equation}
from the explicit formulas \eqref{b70} and \eqref{c19}. This relation justifies the choice of $\CT_{0}$ in Section \ref{s77}. We will now study the resonance free domains given by Theorem \ref{a2} and Theorem \ref{i55}.  As in \eqref{b71}, we define
\begin{equation} \label{j50}
\CA ( \tau , h ) = \spr ( \CT ( E_{0} + h \tau , h ) ) ,
\end{equation}
the spectral radius of the operator $\CT ( E_{0} + h \tau , h )$. Using \eqref{i56}, we deduce
\begin{equation} \label{i57}
\CA ( \tau , h ) \leq \CA_{0} ( \tau ) ,
\end{equation}
for all $\tau , h$. The proof of this fact can be found at the end of Appendix \ref{s70}. On the other hand, Theorem \ref{i55} implies that $P$ has no resonance in the set
\begin{equation} \label{i58}
\left\{ \begin{aligned}
&E_{0} - C h \leq \re z \leq E_{0} + C h ,    \\
&\Big( \lambda \ln \Big( \CA \Big( \frac{\re z - E_{0}}{h}  , h \Big) \Big) + \delta \Big) \frac{h}{\vert \ln h \vert} \leq \im z \leq \frac{h}{\vert \ln h \vert} ,
\end{aligned} \right.
\end{equation}
for all $C , \delta > 0$ and then $h$ small enough. Moreover, this resonance free domain is optimal (at least when $\CA$ does not vanish). From \eqref{i57}, the zone \eqref{i58} is larger than \eqref{b68}, the resonance free domain provided by Theorem \ref{a2}. Thus, Theorem \ref{i55} is more precise under the assumption \ref{h7}, whereas Theorem \ref{a2} holds without this hypothesis.

Since $P$ is self-adjoint (except for some operators following Remark \ref{c13}), the resonances have a non-positive imaginary part. Using Proposition \ref{i53} and Theorem \ref{i55}, it gives
\begin{equation*}
\vert \mu_{k} ( \tau , h ) \vert \leq 1 + o_{h \to 0} ( 1 ) .
\end{equation*}
In other words, the spectral radius of $\CT$ satisfies $\CA ( \tau , h ) \leq 1 + o_{h \to 0} ( 1 )$ for all $\tau \in \R$ and $h \in ] 0 , 1 ]$. Moreover, $\CA ( \tau , h ) \leq 1$ if the action $A$ is constant since $\vert \mu_{k} ( \tau , h ) \vert$ does not depend on $h$ in that case. These estimates were not clear starting from \eqref{c19}. Such upper bounds are not necessarily verified for the matrices $\CQ$ of Section \ref{s6}.

\Subsection{Applications and examples} \label{s54}

Since the asymptotic stated in Theorem \ref{i55} is similar to that of Theorem \ref{d8}, the distribution of resonances in this strong trapping case verify the phenomena described in Section \ref{s18}. We follow the plan of this section and give examples corresponding to each phenomenon.

\subsubsection{Accumulation on curves}

From Theorem \ref{i55}, the resonances satisfy a two scale asymptotic as in Remark \ref{e2}. This is illustrated in Figure \ref{f41}. At the macroscopic scale $h$, they accumulate on the curves
\begin{equation} \label{i71}
\im \sigma =  \ln \big( \vert \mu_{k} ( \re \sigma , h ) \vert \big)  \frac{\lambda}{\vert \ln h \vert} ,
\end{equation}
where $k \in \N$ and $\sigma$ is the rescaled spectral parameter defined in \eqref{i70}. This expression is analogous to \eqref{d93} excepted that the leading term is now given by the eigenvalues of $\widehat{\CT}$. Note that the lower boundary of the region \eqref{i58} is given by the equation of the accumulation curve closest to the real axis. At the microscopic scale $h \vert \ln h \vert^{- 1}$, the packets of resonances (see Section \ref{s53} and the discussion above Proposition \ref{d91}) are horizontally spaced out by
\begin{equation*}
2 \pi \lambda \frac{h}{\vert \ln h \vert} + o \Big( \frac{h}{\vert \ln h \vert} \Big) .
\end{equation*}

\begin{example}\rm \label{i69}
We apply Theorem \ref{i55} in the simplest situation: the case of a single nappe. More precisely, we come back to Example \ref{b80} with $\theta_{1} < \pi / 4$. The geometric setting is described in Figure \ref{f13}. We already know that \ref{h1}--\ref{h4} and \ref{h6} hold true and that $\CH^{\pm \infty}_{\rm tang} = [ - \theta_{0} , \theta_{0} ]$, $\CH_{\rm trans} = \emptyset$. Moreover, $\alpha ( \omega ) = \omega$, $\CM_{0} ( \alpha ) = 1$ and $\nu ( \alpha ) = - 1$. Eventually, the action $A$ and the time delay $T$ are constant on $\CH^{+ \infty}_{\rm tang}$. Thus, all the assumptions of Theorem \ref{i55} are satisfied and $\widehat{\CT}$ takes the form
\begin{equation} \label{i73}
\widehat{\CT} ( \tau , h ) ( \omega , \widetilde{\omega} ) = e^{i A / h} \frac{\lambda}{2 \pi} \Gamma \Big( 1 - i \frac{\tau}{\lambda} \Big) e^{i T \tau} ( i \lambda \omega \cdot \widetilde{\omega} )^{- 1 + i \frac{\tau}{\lambda}} = e^{i A / h} \widetilde{\CT} ( \tau ) ( \omega , \widetilde{\omega} ) ,
\end{equation}
where $\widetilde{\CT}$ does not depend on $h$. Let $\widetilde{\mu}_{1} ( \tau ) , \widetilde{\mu}_{2} ( \tau ) , \ldots$ denote the non-zero eigenvalues of $\widetilde{\CT} ( \tau )$ counted with their multiplicity. Combining with Proposition \ref{i53}, the $z_{q , k}$'s which approximate the (pseudo-)resonances can be written
\begin{equation} \label{i61}
z_{q , k} ( \tau ) = E_{0} - \frac{A \lambda}{\vert \ln h \vert} + 2 q \pi \lambda \frac{h}{\vert \ln h \vert} + i \ln ( \widetilde{\mu}_{k} ( \tau ) ) \lambda \frac{h}{\vert \ln h \vert} .
\end{equation}
This expression is similar to \eqref{e8}, which concerns the resonances generated by a unique homoclinic trajectory. Nevertheless, the $\widetilde{\mu}_{k} ( \tau )$'s are now implicit. Otherwise, we deduce from \eqref{i61} that the equations for the accumulation curves are given by
\begin{equation} \label{i75}
\im \sigma = \ln \big( \vert \widetilde{\mu}_{k} ( \re \sigma ) \vert \big) \frac{\lambda}{\vert \ln h \vert} .
\end{equation}
Note that these curves do not depend on $h$ (except by the explicit factor $\vert \ln h \vert^{- 1}$). This is natural since the action is constant on $\CH^{+ \infty}_{\rm tang}$ (see the discussion below this example). Remark \ref{i60} shows that the number of these accumulation curves is infinite near $\sigma = 0$ when $\theta_{0} > 0$. They are illustrated in Figure \ref{f38}.

\begin{figure}
\begin{center}
\begin{picture}(370,170)
\includegraphics[width=350pt]{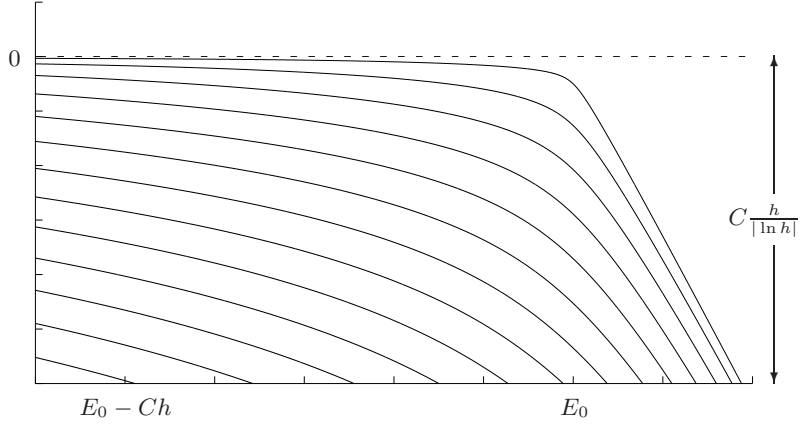}
\put(-25,92){\vector(0,1){52}}
\put(-25,72){\vector(0,-1){52}}
\put(-15,80){\makebox(0,0)[rb]{\smash{{\SetFigFont{9}{10.8}{\rmdefault}{\mddefault}{\updefault}$C \frac{h}{\vert \ln h \vert}$}}}}
\put(-310,140){\makebox(0,0)[rb]{\smash{{\SetFigFont{9}{10.8}{\rmdefault}{\mddefault}{\updefault}$0$}}}}
\put(-95,8){\makebox(0,0)[rb]{\smash{{\SetFigFont{9}{10.8}{\rmdefault}{\mddefault}{\updefault}$E_{0}$}}}}
\put(-253,8){\makebox(0,0)[rb]{\smash{{\SetFigFont{9}{10.8}{\rmdefault}{\mddefault}{\updefault}$E_{0} - C h$}}}}
\end{picture}
\end{center}
\caption{The accumulations curves in Example \ref{i69}. They are obtained through \eqref{i75} and a numerical computation of the eigenvalues of $\widetilde{\CT}$.} \label{f38}
\end{figure}

Finally, we study the distribution of the resonances when the ``opening angle'' (that is $2 \theta_{0}$) is small. For that, we investigate the behavior of the spectrum of $\widehat{\CT}$ in the limit $\theta_{0}$ goes to $0$. For $C > 0$ fixed, it is proved in Section \ref{s30} that
\begin{equation} \label{i65}
\sigma \big( \widehat{\CT} ( \tau , h ) \big) = \{ \mu_{\tau , h , \theta_{0}} \} \cup R_{\tau , h , \theta_{0}} ,
\end{equation}
for $\tau \in [ - C , C ]$. Here, $\mu_{\tau , h , \theta_{0}}$ is a simple eigenvalue satisfying
\begin{equation*}
\mu_{\tau , h , \theta_{0}} = e^{i A / h} q_{\tau} \theta_{0} + o_{\theta_{0} \to 0} ( \theta_{0} ) \qquad \text{with} \qquad q_{\tau} = - \frac{i}{\pi} \Gamma \Big( 1 - i \frac{\tau}{\lambda} \Big) e^{i T \tau} ( i \lambda )^{i \frac{\tau}{\lambda}} ,
\end{equation*}
and the rest of the spectrum verifies $R_{\tau , h , \theta_{0}} \subset B ( 0 , o_{\theta_{0} \to 0} ( \theta_{0} ) )$. These properties hold true uniformly with respect to $\tau \in [ - C , C ]$ and $h \in ] 0 , 1 ]$. Inserting these asymptotic in \eqref{i61}, we deduce
\begin{equation} \label{i67}
z_{q , 1} ( \tau ) = E_{0} - \frac{A \lambda}{\vert \ln h \vert} + 2 q \pi \lambda \frac{h}{\vert \ln h \vert} - i \lambda h \frac{\vert \ln \theta_{0} \vert}{\vert \ln h \vert} + i \ln ( q_{\tau} ) \lambda \frac{h}{\vert \ln h \vert} + o_{\theta_{0} \to 0} ( 1 ) \frac{h}{\vert \ln h \vert} ,
\end{equation}
whereas the other pseudo-resonances satisfy
\begin{equation} \label{i68}
\im z_{q , k} ( \tau ) \leq - \lambda h \frac{\vert \ln \theta_{0} \vert}{\vert \ln h \vert} - C \frac{h}{\vert \ln h \vert} ,
\end{equation}
for $k \geq 2$ and $\theta_{0}$ small enough. Note that, contrarily to \eqref{i61}, the expression for $z_{q , 1} ( \tau )$ is now explicit.

One can first deduce from \eqref{i67} and \eqref{i68} that the resonances go away from the real axis when $\theta_{0}$ goes to $0$. This is natural since the trapping ``decays'' in this regime (see Remark \ref{i76} for $\tau = 0$). On the other hand, the distribution of the resonances looks like to that generated by a unique transversal trajectory (see Corollary \ref{e4}): there is an isolated accumulation curve given by \eqref{i67} and the rest of the resonances is below in the complex plane from \eqref{i68}. This analogy can be specified. If we replace formally $\theta_{0}$ by an appropriate power of $h$ in \eqref{i67}, we fall back on \eqref{e8} (or Proposition \ref{g27}). Nevertheless, the constant $q_{\tau}$ differs from $\mu ( \tau )$. This link with Section \ref{s6} can be formally explained: due to the uncertainty principle, a function localized microlocally near some point of $T^{*} \R^{n}$ is sometimes seen in quantum mechanic as spread in a vicinity of size $h^{c}$ (for some appropriate $c > 0$) of this point. Thus, a transversal intersection as in \ref{h8} (or a tangential intersection of finite order as in \ref{h13}) could be seen as a tangential intersection whose opening angle is of size $h^{c}$.
\end{example}

Since the quantization rule (see Definition \ref{i51}) is implicit, it is not clear that $P$ has actually resonances in \eqref{i54}. We have already seen in \eqref{i67} that there exist operators with at least one accumulation curve. We now give conditions which guaranties an infinity of accumulation curves.

\begin{remark}[Infinity of accumulation curves]\sl \label{i60}
In addition to the assumptions of Theorem \ref{i55}, suppose that $n = 2$, $\mes_{\S^{n - 1}} ( \CH^{- \infty}_{\rm tang} ) > 0$, $\alpha \cdot \omega > 0$ for all $( \alpha , \omega ) \in \CH^{+ \infty}_{\rm tang} \times \CH^{- \infty}_{\rm tang}$, $\alpha ( \omega ) = \omega$ and that $A$, $T$, $\CM_{0}$, $\nu$ are constant on $\CH^{+ \infty}_{\rm tang}$. Then, the non-zero eigenvalues of $\widehat{\CT} ( 0 , h )$ are in infinite number and their modulus is independent of $h$. As consequence, $P$ has an infinity of accumulation curves below $\tau = 0$ and they do not depend on $h$. More precisely, the number of accumulation curves in 
\begin{equation*}
E_{0} + [ - h , h ] + i \Big[ - C \frac{h}{\vert \ln h \vert} , 0 \Big] ,
\end{equation*}
is finite for all $C > 0$ fixed (see the discussion below Proposition \ref{i53}) and goes to infinity as $C \to + \infty$.
\end{remark}

The setting is different from that of Section \ref{s61}. Indeed, the number of accumulation curves is bounded by $K$ in that case. Note that Example \ref{i69} satisfies the assumptions of Remark \ref{i60}. Under these hypotheses, one can also prove that the resonances approach the real axis when the trapping increases. More precisely,

\begin{remark}[Increase of the trapping]\sl \label{i76}
Let $P , Q$ be two operators satisfying the assumptions of Remark \ref{i60} with $\lambda ( P ) = \lambda ( Q )$, $\CM_{0} ( P ) = \CM_{0} ( Q )$ and $\CH^{\pm \infty}_{\rm tang} ( P ) \subset \CH^{\pm \infty}_{\rm tang} ( Q )$. Then, we have
\begin{equation}
\spr \big( \widehat{\CT} ( P ) ( 0 , h ) \big) \leq \spr \big( \widehat{\CT} ( Q ) ( 0 , h ) \big) ,
\end{equation}
for all $h \in ] 0 , 1 ]$. In particular, the resonances of the first accumulation curve of $Q$ are closest to the real axis than the ones of $P$, at least near $E_{0}$ and modulo $ o ( h \vert \ln h \vert^{- 1} )$.
\end{remark}

This phenomenon is natural since $Q$ is ``more geometrically trapping'' than $P$. Roughly speaking, the situation is similar to \eqref{e35} (where the contributions of the different trajectories add) and not to Remark \ref{i42} (where the contributions of the different trajectories conflict). On the other hand, it could be possible to adapt the result of Remark \ref{i42} to the setting of Theorem \ref{i55}.

\begin{example}\rm \label{i72}
We now construct operators such that $\CH^{\pm \infty}_{\rm tang}$ have a rough structure. In particular, the asymptotic tangential directions can form a Cantor set. For that, we adapt the constructions of Example \ref{b86}. Thus, $P$ will not be of the form \eqref{a5}, but will enter in the setting of Remark \ref{c13}.

\begin{figure}
\begin{center}
\begin{picture}(0,0)%
\includegraphics{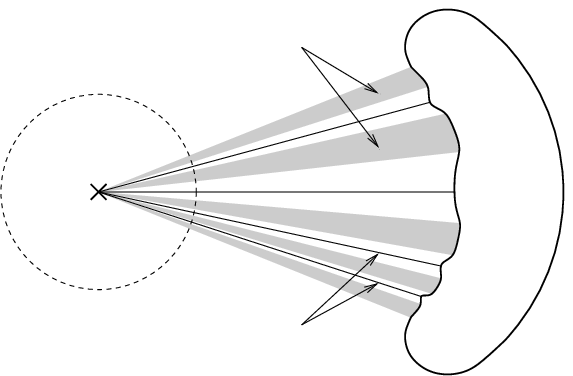}%
\end{picture}%
\setlength{\unitlength}{987sp}%
\begingroup\makeatletter\ifx\SetFigFont\undefined%
\gdef\SetFigFont#1#2#3#4#5{%
  \reset@font\fontsize{#1}{#2pt}%
  \fontfamily{#3}\fontseries{#4}\fontshape{#5}%
  \selectfont}%
\fi\endgroup%
\begin{picture}(10839,7070)(-10289,-8696)
\put(-8774,-5236){\makebox(0,0)[rb]{\smash{{\SetFigFont{9}{10.8}{\rmdefault}{\mddefault}{\updefault}$0$}}}}
\put(-8399,-2986){\makebox(0,0)[b]{\smash{{\SetFigFont{9}{10.8}{\rmdefault}{\mddefault}{\updefault}$\supp V$}}}}
\put(-479,-5236){\makebox(0,0)[b]{\smash{{\SetFigFont{9}{10.8}{\rmdefault}{\mddefault}{\updefault}$\CO$}}}}
\put(-5924,-2086){\makebox(0,0)[lb]{\smash{{\SetFigFont{9}{10.8}{\rmdefault}{\mddefault}{\updefault}$\pi_{x} ( {\mathcal H}_{\rm tang} )$}}}}
\put(-5924,-8236){\makebox(0,0)[lb]{\smash{{\SetFigFont{9}{10.8}{\rmdefault}{\mddefault}{\updefault}$\pi_{x} ( {\mathcal H}_{\rm trans} )$}}}}
\end{picture} $\qquad \qquad$
\begin{picture}(12000,9400)(1500,1400)
\includegraphics[width=220pt]{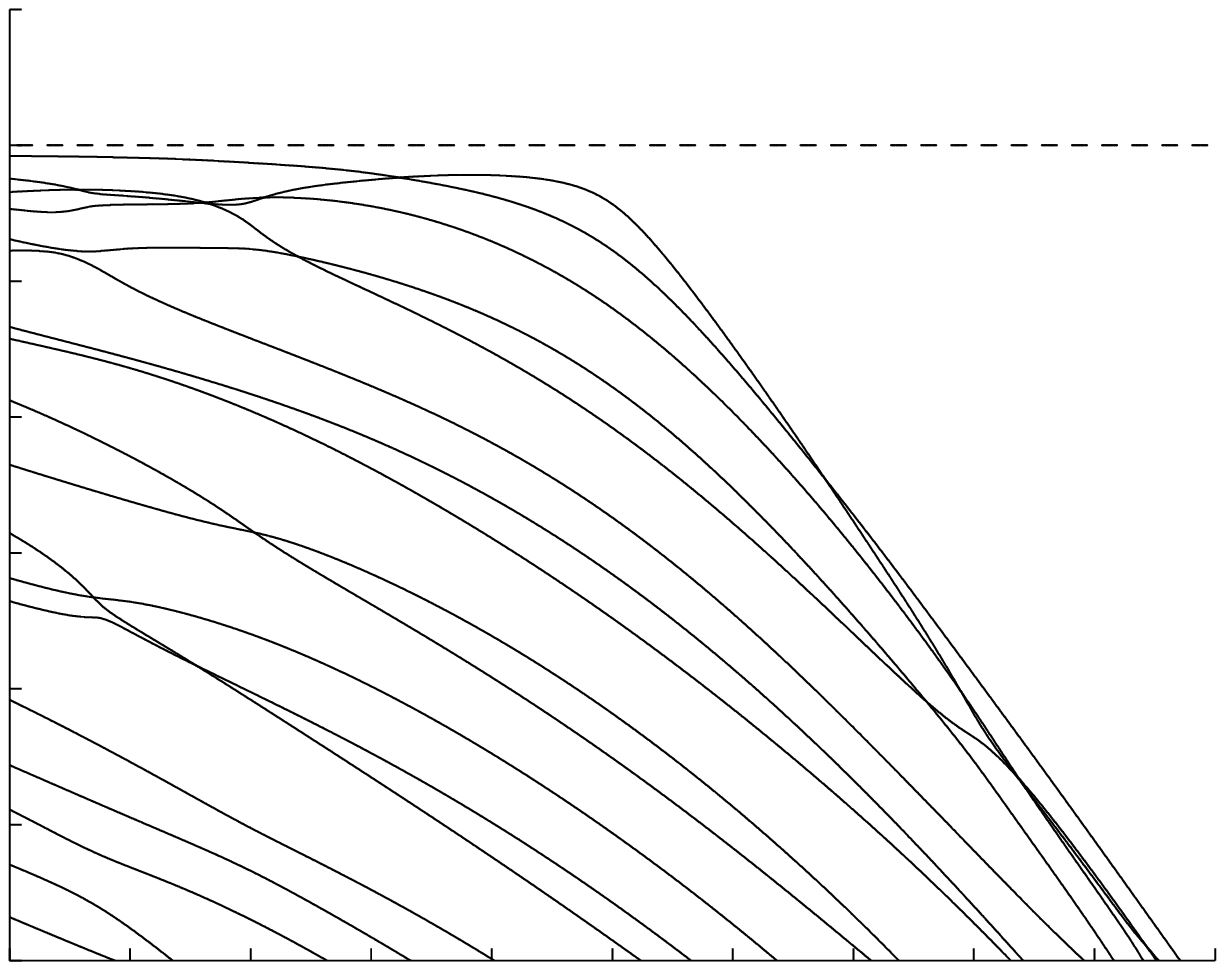}
\end{picture}
\end{center}
\caption{The geometry of Example \ref{i72} and the associated accumulation curves.} \label{f39}
\end{figure}

In dimension $n = 2$, let $(r , \theta )$ be the polar coordinates. As in \eqref{i79}, we first consider 
\begin{equation} \label{i87}
P_{0} = - h^{2} \Delta_{\R^{2} \setminus \CO_{0}} + V (r) ,
\end{equation}
with Dirichlet boundary conditions. Here, $V$ and $\CO_{0}$ are described in Example \ref{b86} (see also Figure \ref{f16}). We now perturb $\CO_{0}$ to change the homoclinic set. Let $\CO$ be the smooth non-trapping obstacle which coincides with $\CO_{0}$ except that the part of the boundary $r = R_{0}$ is replaced by
\begin{equation*}
r = R_{0} + F ( \theta ) ,
\end{equation*}
where $F \in C^{\infty}$ is a well chosen function such that $F ( \theta ) = F^{\prime} ( \theta ) = F^{\prime \prime} ( \theta ) = 0$ for $\theta \in K$ where $K$ is some compact set, and that $F^{\prime \prime} ( \theta ) \neq 0$ for all $\theta \notin K$ with $F^{\prime} ( \theta ) = 0$ (see Figure \ref{f39}). It is possible to construct such a function $F$ where $K$ is a Cantor set localized near $0$. Indeed, one can choose $F = F_{1} + F_{2} + \cdots$ where the $F_{\bullet}$'s have disjoint support and remove inductively some parts of the angular sector of the homoclinic trajectories of $P_{0}$.

We then define
\begin{equation*}
P = - h^{2} \Delta_{\R^{2} \setminus \CO} + V (r) .
\end{equation*}
Using Proposition \ref{a14}, one can prove that the trapped set of $P$ at energy $E_{0}$ verifies \ref{h3} and that the homoclinic set consists of the radial rays whose angle $\theta$ satisfies $F^{\prime} ( \theta ) = 0$. Thus, \ref{h4} holds true. From the construction of $F$, $\CH^{\pm \infty}_{\rm tang} = K$ and $\CH_{\rm trans}$ consists of the radial rays whose angle $\theta \notin K$ satisfies $F^{\prime} ( \theta ) = 0$. Moreover, the quantities $A$, $T$ are constant on $\CH^{+ \infty}_{\rm tang}$, $\CM_{0} ( \alpha ) = 1$, $\nu ( \alpha ) = - 1$, $\alpha ( \omega ) = \omega$ and $\alpha \cdot \omega > 0$ for all $( \alpha , \omega ) \in \CH^{+ \infty}_{\rm tang} \times \CH^{- \infty}_{\rm tang}$. In particular, \ref{h7} is verified.

Then, Proposition \ref{i53} and Theorem \ref{i55} provide the asymptotic of the resonances. The kernel of the operator $\widehat{\CT} ( \tau , h )$ on $L^{2} ( K )$ is given by \eqref{i73}. The $z_{q , k}$'s satisfy \eqref{i61} mutatis mutandis. Moreover, Remark \ref{b69} $i)$, Remark \ref{i60} and Remark \ref{i76} can be applied here. So, $P$ has an infinity of accumulation curves below $E_{0}$ (that is for $\tau = 0$) if $\mes_{\S^{n - 1}} ( K ) \neq 0$, and no resonance in \eqref{i54} if $\mes_{\S^{n - 1}} ( K ) = 0$. Since the action is constant on $K$, the accumulation curves are independent of $h$. They are drawn in Figure \ref{f39} using a computer calculation of the eigenvalues of $\widehat{\CT}$. That $K$ is irregular presents no difficulty for the numerical analysis. Indeed, one can replace $K$ by a finite union of intervals modulo lower order terms. This is a consequence of the stability result stated in Remark \ref{i89} $iii)$.
\end{example}

\begin{example}\rm \label{i80}
Using an absorbing potential, it is possible to construct an operator where $\CH_{\rm tang}^{\pm \infty}$ is ``formally'' any given compact subset $K$ of $\S^{n - 1}$ with sufficiently small diameter. Indeed, consider first $P_{0}$ the operator build in Example \ref{b86} (or the operator of Example \ref{b80}) but in dimension $n \geq 2$. The asymptotic homoclinic directions of $P_{0}$ form a sector $K_{0} \subset \S^{n - 1}$. Let also $F \in C^{\infty} ( \S^{n - 1} ; [ 0 , 1 ] )$ be such that $F ( \theta ) = 0$ for $\theta \in K$ and $F ( \theta ) > 0$ for $\theta \in K_{0} \setminus K$. The existence of such a function follows from the Whitney covering lemma. We set
\begin{equation*}
P = P_{0} - i h \vert \ln h \vert F ( \theta ) G ( r ) ,
\end{equation*}
where $G \in C^{\infty}_{0} ( \R ; [ 0 , 1 ] )$ is a cut-off function supported between the fixed point and the barrier (see Figure \ref{f40}). This operator enters in the setting of Remark \ref{c13}.

\begin{figure}
\begin{center}
\begin{picture}(0,0)%
\includegraphics{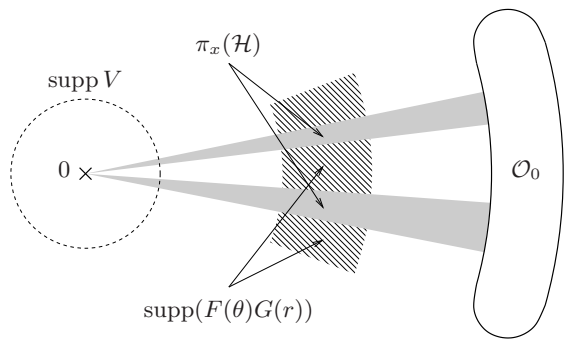}%
\end{picture}%
\setlength{\unitlength}{987sp}%
\begingroup\makeatletter\ifx\SetFigFont\undefined%
\gdef\SetFigFont#1#2#3#4#5{%
  \reset@font\fontsize{#1}{#2pt}%
  \fontfamily{#3}\fontseries{#4}\fontshape{#5}%
  \selectfont}%
\fi\endgroup%
\begin{picture}(13914,8328)(-10289,-9325)
\put(-8774,-5236){\makebox(0,0)[rb]{\smash{{\SetFigFont{9}{10.8}{\rmdefault}{\mddefault}{\updefault}$0$}}}}
\put(-8399,-2986){\makebox(0,0)[b]{\smash{{\SetFigFont{9}{10.8}{\rmdefault}{\mddefault}{\updefault}$\supp V$}}}}
\put(2677,-5236){\makebox(0,0)[b]{\smash{{\SetFigFont{9}{10.8}{\rmdefault}{\mddefault}{\updefault}$\CO_{0}$}}}}
\put(-4799,-8761){\makebox(0,0)[b]{\smash{{\SetFigFont{9}{10.8}{\rmdefault}{\mddefault}{\updefault}$\supp ( F ( \theta ) G ( r ) )$}}}}
\put(-4799,-2086){\makebox(0,0)[b]{\smash{{\SetFigFont{9}{10.8}{\rmdefault}{\mddefault}{\updefault}$\pi_{x} ( {\mathcal H} )$}}}}
\end{picture}%
\end{center}
\caption{The geometry of Example \ref{i80}.} \label{f40}
\end{figure}

In this situation, the principal symbol and the Hamiltonian trajectories of $P$ are those of $P_{0}$. Thus, the set of asymptotic tangential directions of $P$ is $K_{0}$, as for $P_{0}$. Nevertheless, the homoclinic trajectories which touch the support of $F ( \theta ) G ( r )$ play no role in the asymptotic of the resonances of $P$. Indeed, the microlocal solutions of $( P - z ) u = 0$ are absorbed in this support. Then, the kernel of the quantization operator $\CT$ is \eqref{c19} restricted to $K$ and Theorem \ref{i55} provides the asymptotic of the resonances of $P$. In other words, the ``effective''  set of asymptotic tangential directions is indeed $K$. Remark that the present operator satisfies also the assumptions of Remark \ref{i60}.
\end{example}

\subsubsection{Vibration phenomena}

Coming back to the general setting, the accumulation curves can vibrate in the sense of Remark \ref{e5}. It means that the shape of \eqref{i71} can change with $h$. This phenomenon occurs when the modulus of $\mu_{k} ( \re \sigma , h )$ depends on $h$. More precisely, thanks to \ref{h7}, the number of different actions $A_{1} , \ldots , A_{K}$ is finite. Thus, the accumulation curves in the domain \eqref{i54} are smooth (at least continuous) functions of $e^{i ( A_{2} - A_{1} ) / h} , \ldots , e^{i ( A_{K} - A_{1} ) / h}$.

Moreover, the particular cases of Remark \ref{e5} still hold true. First, if the action is constant on $\CH^{+ \infty}_{\rm tang}$ modulo a subset of measure zero, the accumulation curves do not depend on $h$. Second, if there exist exactly two different actions, denoted $A_{1}$ and $A_{2}$, modulo a subset of measure zero, the accumulation curves are periodic functions of $h^{- 1}$ with period $2 \pi \vert A_{2} - A_{1} \vert^{- 1}$.

We now construct an Schr\"{o}dinger operator with different actions. It will be similar to the one of Example \ref{i69} but with two ``croissants''.

\begin{example}\rm \label{b88}
In dimension $n = 2$, let $(r , \theta )$ be the polar coordinates. We consider
\begin{equation*}
V (x) = V_{0} (r) + V_{\infty} ( r - a_{1} ) \psi_{1} ( \theta ) + V_{\infty} ( r - a_{2} ) \psi_{2} ( \Theta + \theta ) ,
\end{equation*}
where $V_{0}$ is as in Section \ref{s21}. The potential $V_{\infty} \in C_{0}^{\infty} ( \R )$ is an even function satisfying $r V_{\infty}^{\prime} (r) < 0$ for $r$ in the interior of $\supp V_{\infty} \setminus \{ 0 \}$ and $E_{0} < V_{\infty} ( 0 )$. The constants $a_{\bullet} > 0$ are chosen sufficiently large (in particular, $\supp V_{0} (r) \cap \supp V_{\infty} ( r - a_{\bullet} ) = \emptyset$). Moreover, $\psi_{\bullet} ( \theta ) \in C_{0}^{\infty} ( [ - \theta^{1}_{\bullet} , \theta^{1}_{\bullet} ] )$ is equal to $1$ for $\vert \theta \vert \leq \theta^{0}_{\bullet}$ and $\theta \psi^{\prime}_{\bullet} ( \theta ) < 0$ for $\theta^{0}_{\bullet} < \vert \theta \vert < \theta^{1}_{\bullet}$ for some $0 \leq \theta^{0}_{\bullet} < \theta^{1}_{\bullet} \leq \pi$. Eventually, $\Theta \in \R$ is chosen near $\pi$. The geometric setting is illustrated in Figure \ref{f17}.

\begin{figure}
\begin{center}
\begin{picture}(0,0)%
\includegraphics{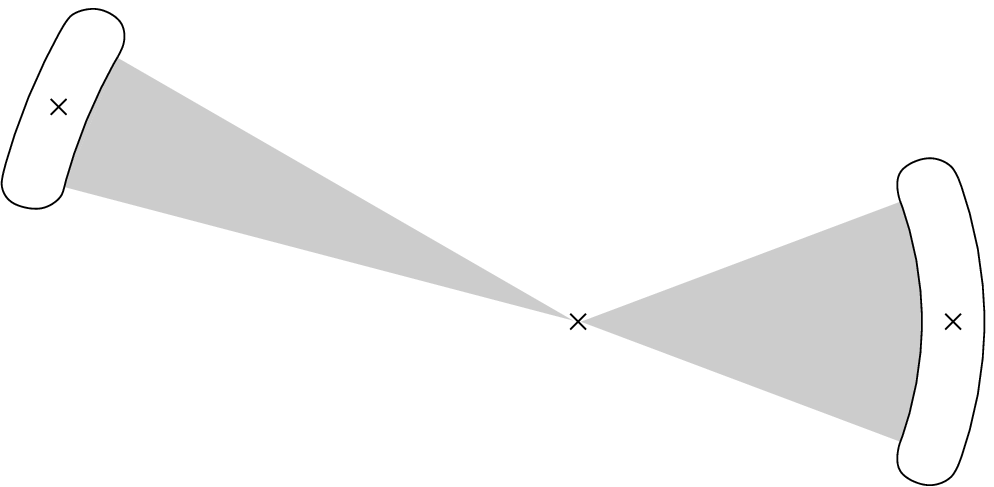}%
\end{picture}%
\setlength{\unitlength}{987sp}%
\begingroup\makeatletter\ifx\SetFigFont\undefined%
\gdef\SetFigFont#1#2#3#4#5{%
  \reset@font\fontsize{#1}{#2pt}%
  \fontfamily{#3}\fontseries{#4}\fontshape{#5}%
  \selectfont}%
\fi\endgroup%
\begin{picture}(18926,9211)(-19500,-8333)
\put(-8399,-4786){\makebox(0,0)[b]{\smash{{\SetFigFont{9}{10.8}{\rmdefault}{\mddefault}{\updefault}$0$}}}}
\put(-2549,-2086){\makebox(0,0)[rb]{\smash{{\SetFigFont{9}{10.8}{\rmdefault}{\mddefault}{\updefault}$\{ V(x) = E_{0} \}$}}}}
\put(-16724,-2311){\makebox(0,0)[lb]{\smash{{\SetFigFont{9}{10.8}{\rmdefault}{\mddefault}{\updefault}$\pi_{x} ( {\mathcal H} )$}}}}
\put(-6899,-5311){\makebox(0,0)[lb]{\smash{{\SetFigFont{9}{10.8}{\rmdefault}{\mddefault}{\updefault}$2 \theta_{1}^{0}$}}}}
\put(-1199,-4786){\makebox(0,0)[b]{\smash{{\SetFigFont{9}{10.8}{\rmdefault}{\mddefault}{\updefault}$a_{1}$}}}}
\put(-18035,-586){\makebox(0,0)[b]{\smash{{\SetFigFont{9}{10.8}{\rmdefault}{\mddefault}{\updefault}$a_{2} e^{i \Theta}$}}}}
\put(-11867,-4186){\makebox(0,0)[lb]{\smash{{\SetFigFont{9}{10.8}{\rmdefault}{\mddefault}{\updefault}$2 \theta_{2}^{0}$}}}}
\put(-2849,-5236){\makebox(0,0)[rb]{\smash{{\SetFigFont{9}{10.8}{\rmdefault}{\mddefault}{\updefault}$\pi_{x} ( {\mathcal H} )$}}}}
\end{picture} \newline
\begin{picture}(28000,4850)(1500,750)
\includegraphics[width=113pt]{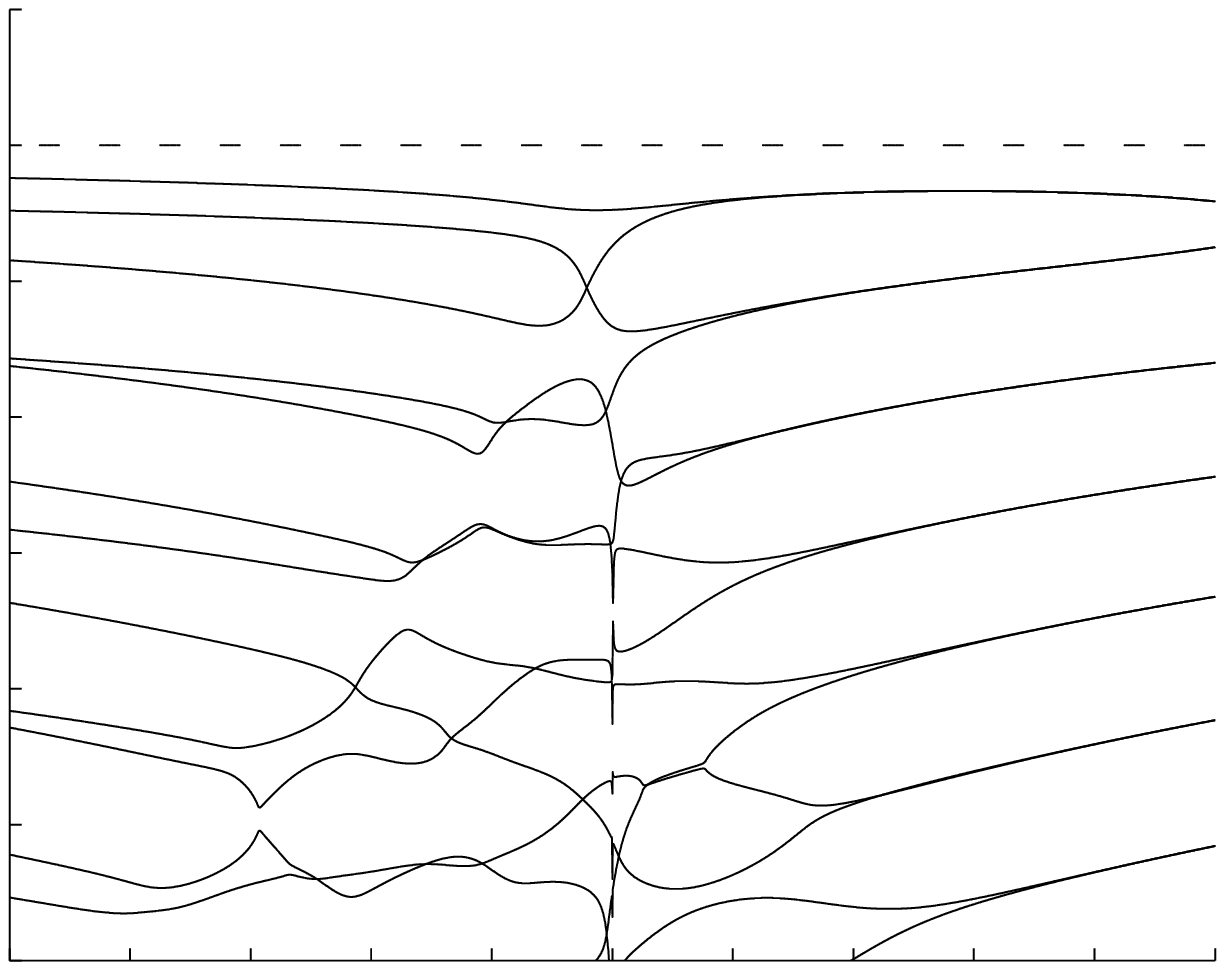}
\includegraphics[width=113pt]{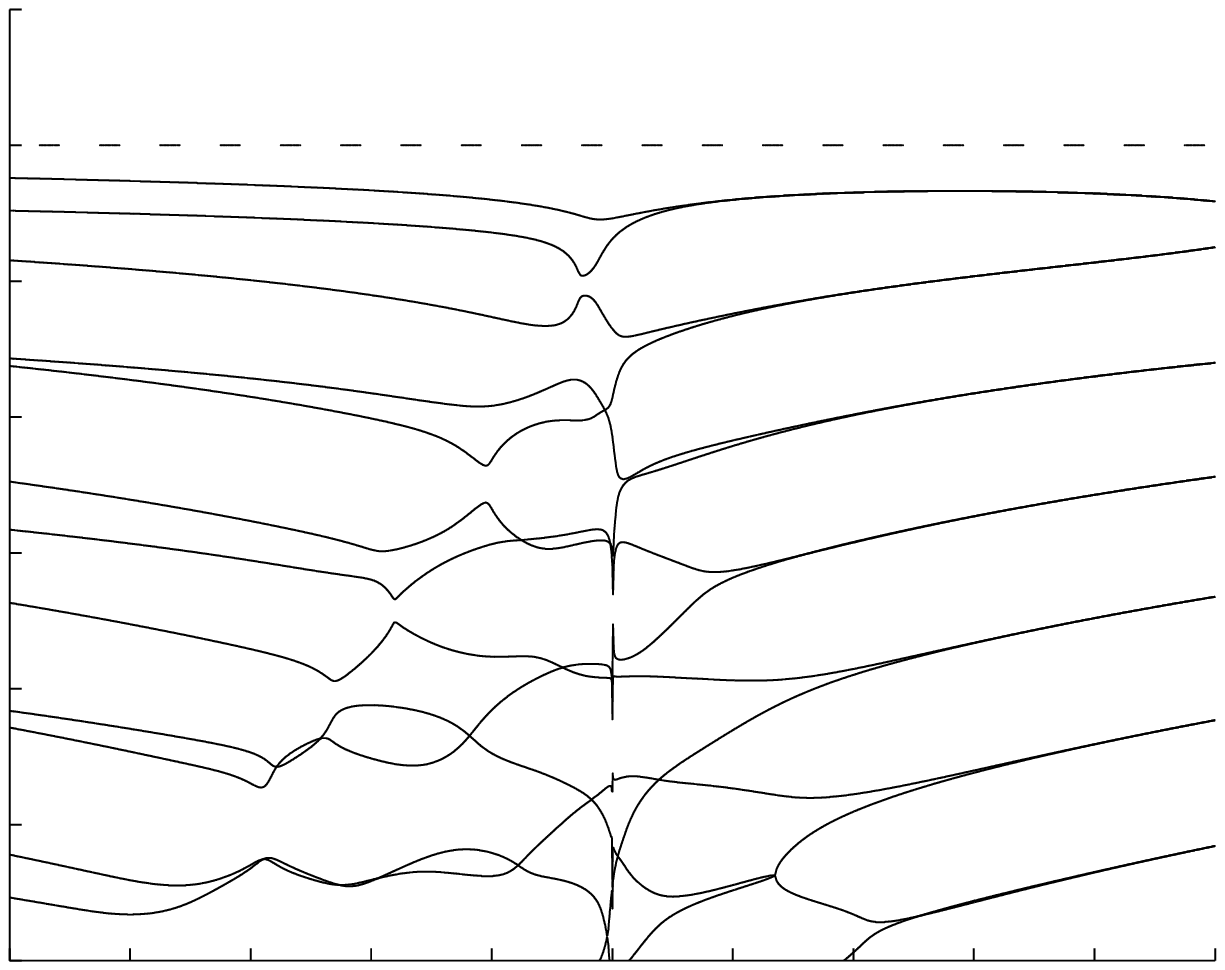}
\includegraphics[width=113pt]{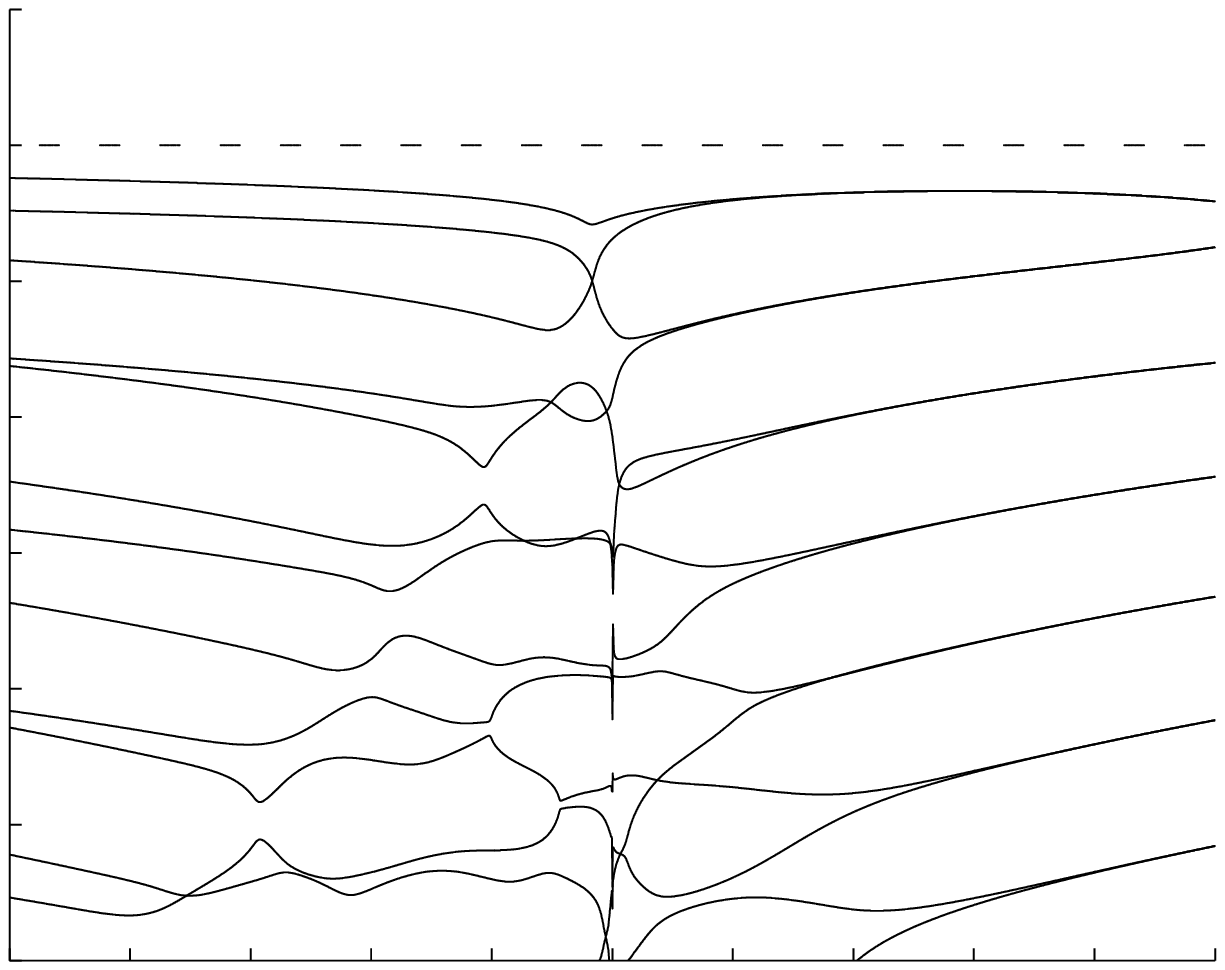}
\includegraphics[width=113pt]{des_3.eps}
\end{picture}
\end{center}
\caption{The geometry of Example \ref{b88} and the corresponding accumulation curves of resonances for different values of $h$.} \label{f17}
\end{figure}

For $\theta_{\bullet}^{1}$ small enough and $\Theta$ close enough to $\pi$, the assumptions \ref{h1}--\ref{h4}, \ref{h6} and \ref{h7} are all satisfied and $\CH^{\infty} = \CH^{\pm \infty}_{\rm tang} = [ - \theta^{0}_{1} , \theta^{0}_{1} ] \cup [ \Theta - \theta^{0}_{2} , \Theta + \theta^{0}_{2} ]$. In fact, to prove these assertions, it is enough to process as in Example \ref{b80} and to remark that, from Lemma \ref{a43} and Proposition \ref{a14}, no Hamiltonian curve coming from a ``croissant'' (say the support of $V_{\infty} ( r - a_{1} ) \psi_{1} ( \theta )$) can touch the other ``croissant'' (the support of $V_{\infty} ( r - a_{2} ) \psi_{2} ( \Theta + \theta )$). Moreover, we have $\alpha ( \omega ) = \omega$, $\CM_{0} ( \alpha ) = 1$ and $\nu ( \alpha ) = - 1$ as in Example \ref{i69}. Eventually, there are two different actions (say $A_{1} , A_{2}$) and times delay (say $T_{1} , T_{2}$) when $a_{1} \neq a_{2}$.

The asymptotic of the resonances is then given by Theorem \ref{i55}. Since there are two different actions, the accumulation curves are periodic functions of $h^{- 1}$. They are drawn in Figure \ref{f17} using a numerical computation of the eigenvalues of $\CT$. On the other hand, with the notations of Remark \ref{b72}, we have
\begin{equation*}
\CM_{0} = 1 \qquad \text{and} \qquad \CJ_{0} ( \sigma ) = 2 \cosh \Big( \frac{\pi}{2} \sigma \Big) \ln \tan \Big( \theta^{0} + \frac{\pi}{4} \Big) ,
\end{equation*}
when $\theta^{0} : = \theta^{0}_{1} = \theta^{0}_{2} $ and $\Theta = \pi$. This computation, Remark \ref{b72} and Theorem \ref{a2} give an explicit resonance free zone which is non-empty but localized in a set where $( \re z - E_{0} ) / h$ is bounded. Note that all the energies in a neighborhood of $E_{0}$ are trapped.

Furthermore, one can also consider the small ``opening angle'' limit as in \eqref{i65}--\eqref{i68}. Assume that $\theta_{1}^{0}$ and $\theta_{2}^{0}$ are small and comparable. In the sequel, we suppose that $\theta_{2}^{0} = \alpha \theta_{1}^{0}$ with $\alpha \neq 0$ to fix the ideas. For $\tau \in \R$, let $\widehat{\CQ} ( \tau , h )$ denote the $2 \times 2$ matrix with coefficients
\begin{equation} \label{i85}
\widehat{\CQ}_{k , \ell} ( \tau , h ) = e^{i A_{k} / h} \frac{\lambda}{\pi} \Gamma \Big( 1 - i \frac{\tau}{\lambda} \Big) e^{i T_{k} \tau} \alpha^{\frac{k + \ell - 2}{2}} \left\{ \begin{aligned}
&( i \lambda )^{- 1 + i \frac{\tau}{\lambda}} && \text{ if } k = \ell , \\
&( i \lambda \cos \Theta )^{- 1 + i \frac{\tau}{\lambda}} && \text{ if } k \neq \ell ,
\end{aligned} \right.
\end{equation}
and $\mu_{1} ( \tau, h ) , \mu_{2} ( \tau, h )$ its two eigenvalues. Equation \eqref{e21} provides an explicit formula for these quantities. The $z_{q , k} ( \tau )$'s of Proposition \ref{i53} verify
\begin{equation} \label{i83}
z_{q , k} ( \tau ) = E_{0} + 2 q \pi \lambda \frac{h}{\vert \ln h \vert} - i \lambda h \frac{\vert \ln \theta_{1}^{0} \vert}{\vert \ln h \vert} + i \ln ( \mu_{k} ( \tau, h ) ) \lambda \frac{h}{\vert \ln h \vert} + o_{\theta_{1}^{0} \to 0} ( 1 ) \frac{h}{\vert \ln h \vert} ,
\end{equation}
for $k = 1 , 2$ and
\begin{equation} \label{i84}
\im z_{q , k} ( \tau ) \leq - \lambda h \frac{\vert \ln \theta_{1}^{0} \vert}{\vert \ln h \vert} - C \frac{h}{\vert \ln h \vert} ,
\end{equation}
for $k \geq 3$ and $\theta_{1}^{0}$ small enough. These relations hold uniformly for $\tau \in [ - C , C ]$ and $h \in ] 0 , 1 ]$ with $C > 0$ fixed. We omit the proof of \eqref{i83}--\eqref{i84} since it is the same as \eqref{i67}--\eqref{i68}, mutatis mutandis.

In some sense, the present setting has some analogies with Example \ref{e6}. The two transversal homoclinic trajectories ``on the opposite size of $0$'' are replaced here by two nappes. Concerning the resonances, similarities also occur: the two first accumulation curves of Figure \ref{f17} look like to the ones of Figure \ref{f25}. In the small ``opening angle'' limit, this resemblance can be proved using \eqref{i83}--\eqref{i84}. Nevertheless, the matrices \eqref{e18} and \eqref{i85} are different.
\end{example}

One can also construct operators with two nappes ``on the same side of $0$'' with different actions.

\begin{figure}
\begin{center}
\begin{picture}(0,0)%
\includegraphics{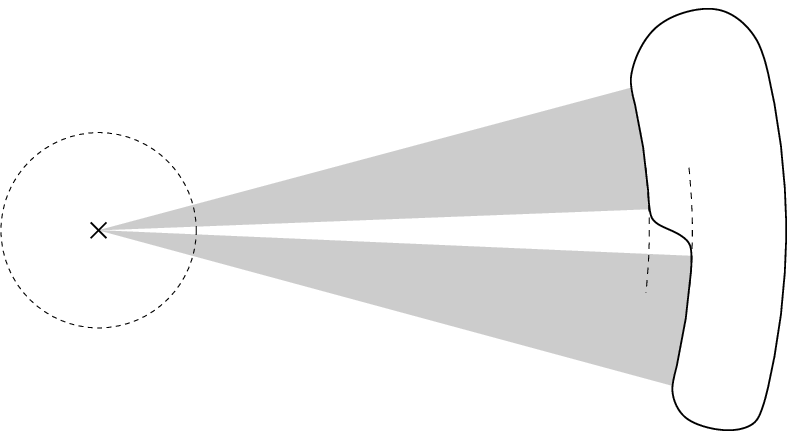}%
\end{picture}%
\setlength{\unitlength}{987sp}%
\begingroup\makeatletter\ifx\SetFigFont\undefined%
\gdef\SetFigFont#1#2#3#4#5{%
  \reset@font\fontsize{#1}{#2pt}%
  \fontfamily{#3}\fontseries{#4}\fontshape{#5}%
  \selectfont}%
\fi\endgroup%
\begin{picture}(15114,8157)(-10289,-9033)
\put(-8774,-5236){\makebox(0,0)[rb]{\smash{{\SetFigFont{9}{10.8}{\rmdefault}{\mddefault}{\updefault}$0$}}}}
\put(-8399,-2986){\makebox(0,0)[b]{\smash{{\SetFigFont{9}{10.8}{\rmdefault}{\mddefault}{\updefault}$\supp V$}}}}
\put(3301,-2836){\makebox(0,0)[b]{\smash{{\SetFigFont{9}{10.8}{\rmdefault}{\mddefault}{\updefault}$\CO$}}}}
\put(-2999,-6136){\makebox(0,0)[b]{\smash{{\SetFigFont{9}{10.8}{\rmdefault}{\mddefault}{\updefault}$I_{2}$}}}}
\put(-2999,-4486){\makebox(0,0)[b]{\smash{{\SetFigFont{9}{10.8}{\rmdefault}{\mddefault}{\updefault}$I_{1}$}}}}
\end{picture} \newline
\begin{picture}(28000,4800)(1500,750)
\includegraphics[width=113pt]{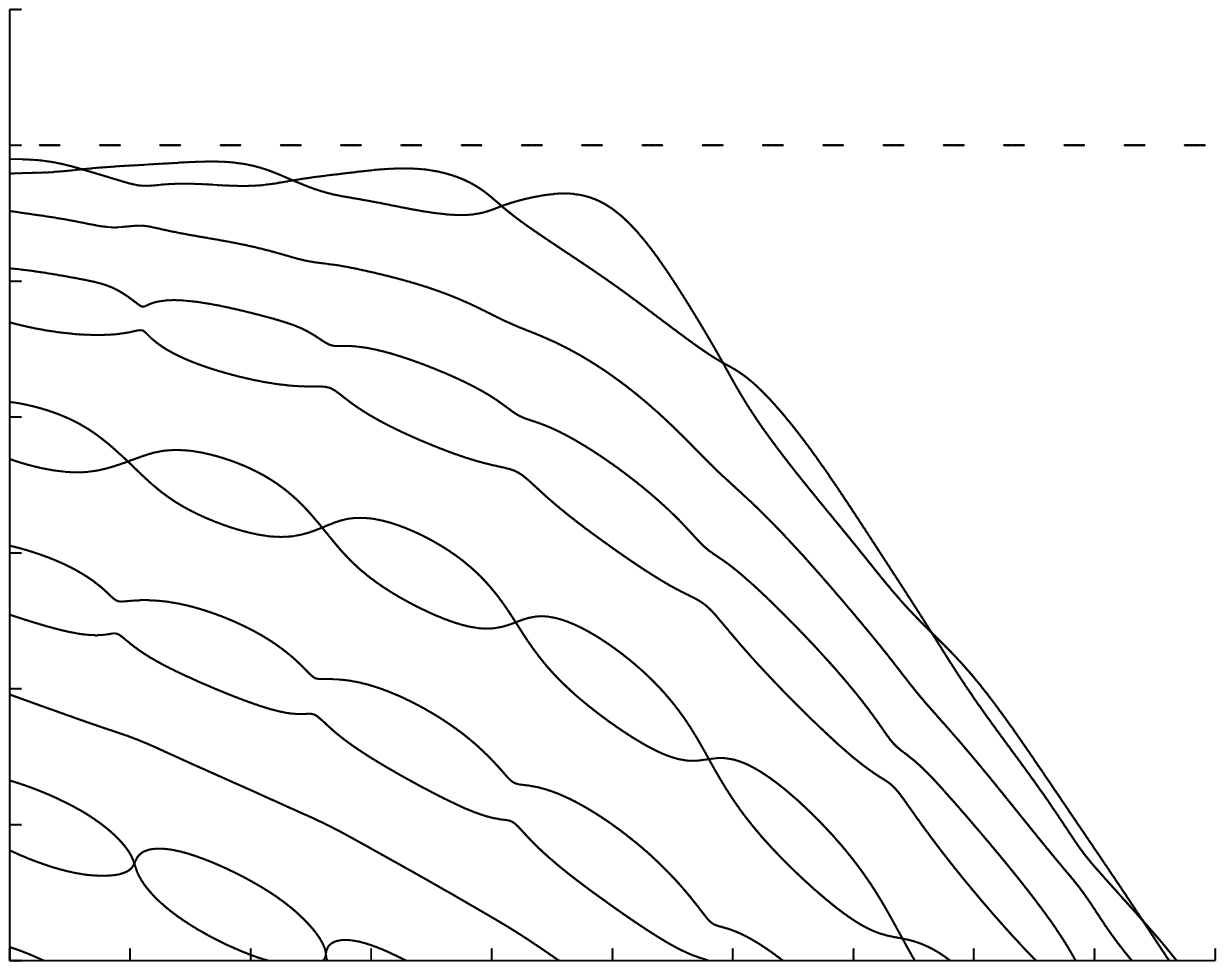}
\includegraphics[width=113pt]{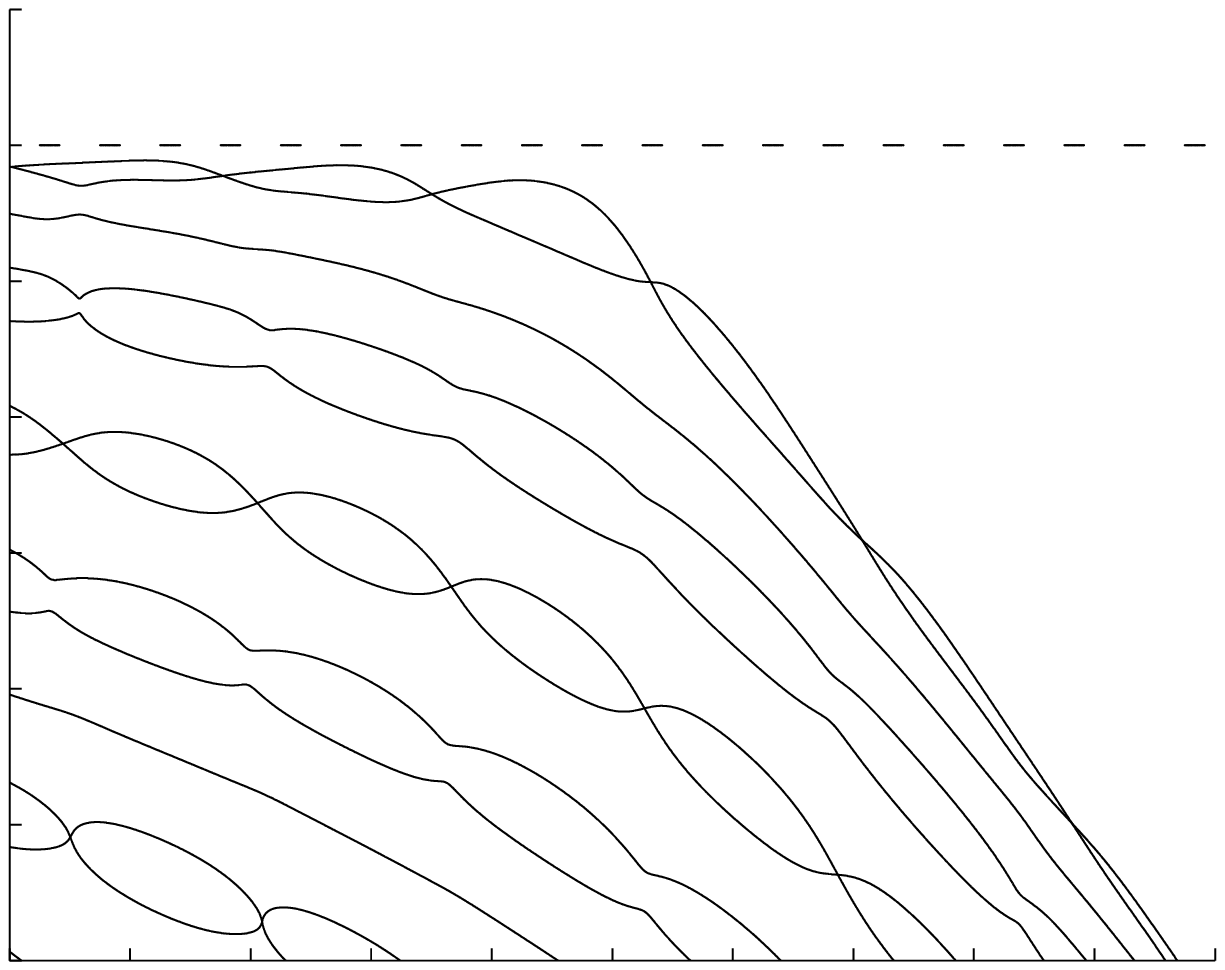}
\includegraphics[width=113pt]{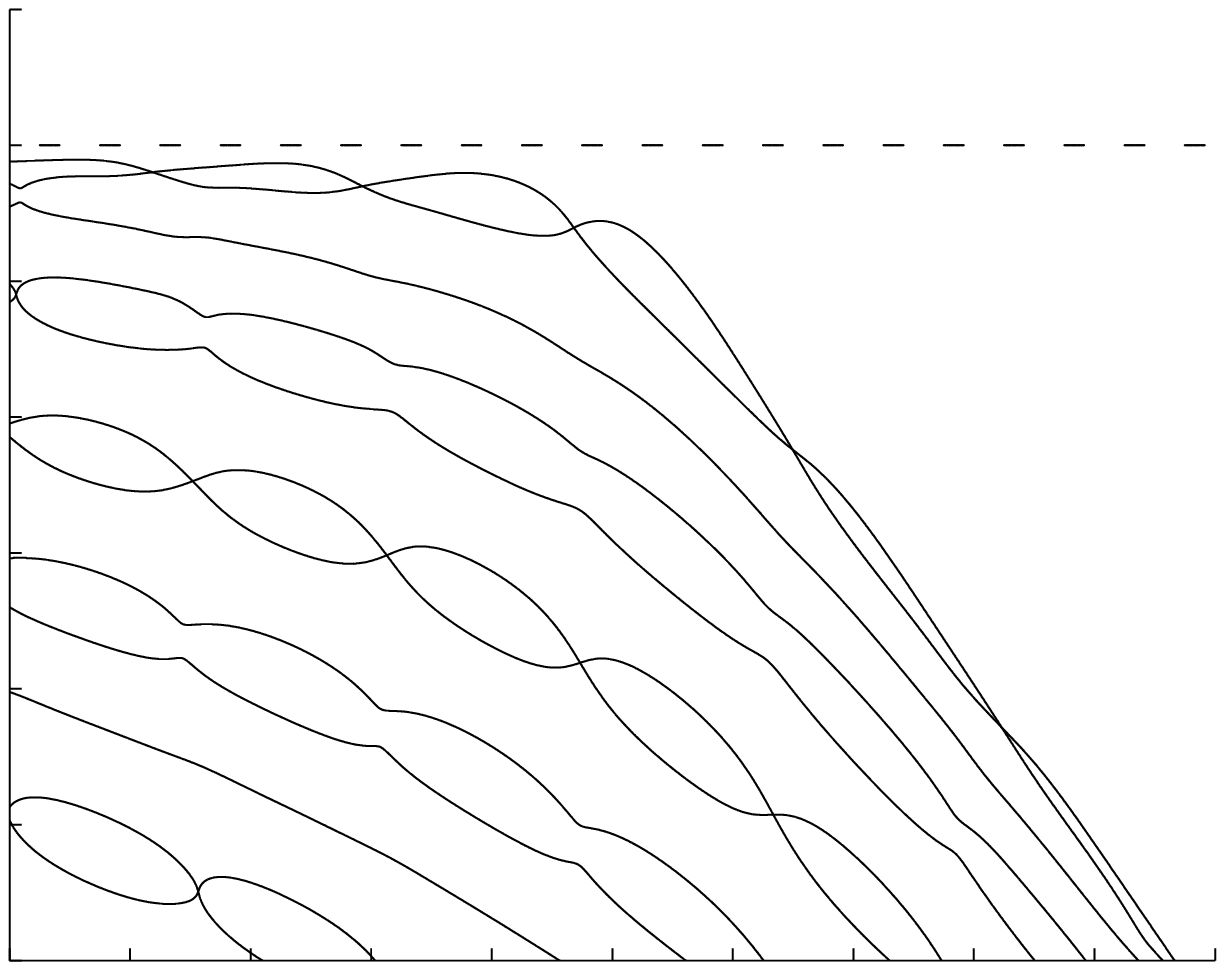}
\includegraphics[width=113pt]{des_7.eps}
\end{picture}
\end{center}
\caption{The geometry of Example \ref{i86} and the corresponding accumulation curves of resonances for different values of $h$.} \label{f42}
\end{figure}

\begin{example}\rm \label{i86}
We adapt Example \ref{i72} and use its notations. Let $P_{0}$ be the operator defined in \eqref{i87}. As previously, let $\CO$ be the obstacle $\CO_{0}$ where the boundary $r = R_{0}$ is replaced by $r = R_{0} + F ( \theta )$. We choose $F \in C^{\infty}$ as in Figure \ref{f42} such that $\CO$ is non-trapping,
\begin{equation*}
F ( \theta ) = \left\{ \begin{aligned}
&a_{1} &&\text{ for } \theta \in I_{1} , \\
&a_{2} &&\text{ for } \theta \in I_{2} ,
\end{aligned} \right.
\end{equation*}
for two intervals $I_{1} , I_{2} \subset K_{0}$ and $F^{\prime} ( \theta ) \neq 0$ for $\theta \in K_{0} \setminus ( I_{1} \cup I_{2} )$ where $K_{0}$ is the set of asymptotic homoclinic directions of $P_{0}$. As usual, we set
\begin{equation*}
P = - h^{2} \Delta_{\R^{2} \setminus \CO} + V (r) .
\end{equation*}

Like in the previous examples, \ref{h1}--\ref{h4}, \ref{h6} are satisfy. As explained in Example \ref{i72}, the homoclinic trajectories correspond to the angles $\theta \in K_{0}$ such that $F^{\prime} ( \theta )= 0$. By construction, we deduce $\CH_{\rm tang}^{\pm \infty} = I_{1} \cup I_{2}$ and $\CH_{\rm trans} = \emptyset$. Moreover, we have $\alpha ( \omega ) = \omega$, $\CM_{0} ( \alpha ) = 1$ and $\nu ( \alpha ) = - 1$. Eventually, there are two different actions (say $A_{1} , A_{2}$) and times delay (say $T_{1} , T_{2}$) when $a_{1} \neq a_{2}$. Thus, \ref{h7} holds true.

The distribution of the resonances is described in Theorem \ref{i55}. The accumulation curves are illustrated in Figure \ref{f42}. They are periodic functions of $h^{- 1}$ since $\CH_{\rm tang}^{\pm \infty}$ consists of two nappes ($I_{1}$ and $I_{2}$) with different actions. The situation is somehow similar to Example \ref{e29}. As in Example \ref{i69} and Example \ref{b88}, this analogy can be specified in the small ``opening angle'' limit. This point is not developed here since it is similar to what has already been made.
\end{example}

\subsubsection{Transition phenomena} \label{s55}

We now examine the behavior of the accumulation curves in the limits $\re \sigma \to \pm \infty$. It can be formally justified by the nature of the trapping below and above $E_{0}$, and illustrates the transitional nature of the homoclinic trappings. In the transversal case, such questions have been treated in details in Section \ref{s14}. By comparison, we have not here the explicit form of the accumulation curves. Thus, we will mainly concentrate on the resonance free domains and use the results of Section \ref{s77} rather than those of Section \ref{s53}. Figure \ref{f15} summarizes the results that will be obtained for large $\pm \re \sigma$. If $\CH^{\pm \infty}_{\rm tang} \neq \emptyset$, three different cases are possible:

\begin{figure}
\begin{center}
\begin{picture}(0,0)%
\includegraphics{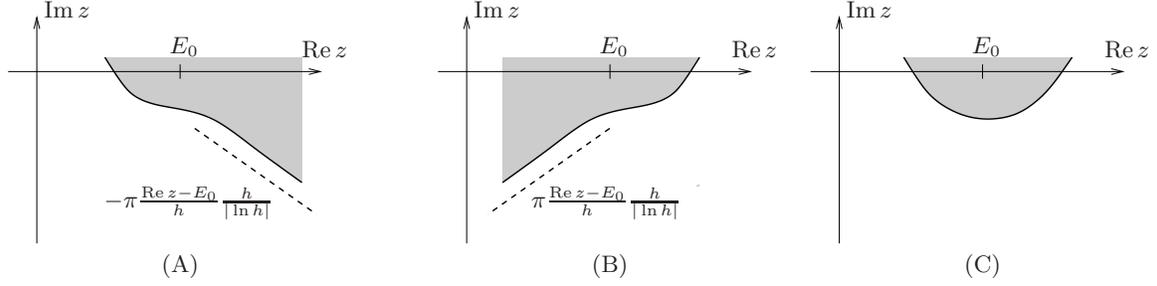}%
\end{picture}%
\setlength{\unitlength}{1184sp}%
\begingroup\makeatletter\ifx\SetFigFont\undefined%
\gdef\SetFigFont#1#2#3#4#5{%
  \reset@font\fontsize{#1}{#2pt}%
  \fontfamily{#3}\fontseries{#4}\fontshape{#5}%
  \selectfont}%
\fi\endgroup%
\begin{picture}(23444,5503)(-621,-6416)
\put(22351,-1861){\makebox(0,0)[lb]{\smash{{\SetFigFont{9}{10.8}{\rmdefault}{\mddefault}{\updefault}$\re z$}}}}
\put(16951,-1036){\makebox(0,0)[lb]{\smash{{\SetFigFont{9}{10.8}{\rmdefault}{\mddefault}{\updefault}$\im z$}}}}
\put(4951,-5011){\makebox(0,0)[rb]{\smash{{\SetFigFont{9}{10.8}{\rmdefault}{\mddefault}{\updefault}$- \pi \frac{\re z - E_{0}}{h} \frac{h}{\vert \ln h \vert}$}}}}
\put(5551,-1861){\makebox(0,0)[lb]{\smash{{\SetFigFont{9}{10.8}{\rmdefault}{\mddefault}{\updefault}$\re z$}}}}
\put(151,-1036){\makebox(0,0)[lb]{\smash{{\SetFigFont{9}{10.8}{\rmdefault}{\mddefault}{\updefault}$\im z$}}}}
\put(14551,-1861){\makebox(0,0)[lb]{\smash{{\SetFigFont{9}{10.8}{\rmdefault}{\mddefault}{\updefault}$\re z$}}}}
\put(9151,-1036){\makebox(0,0)[lb]{\smash{{\SetFigFont{9}{10.8}{\rmdefault}{\mddefault}{\updefault}$\im z$}}}}
\put(10351,-5011){\makebox(0,0)[lb]{\smash{{\SetFigFont{9}{10.8}{\rmdefault}{\mddefault}{\updefault}$\pi \frac{\re z - E_{0}}{h} \frac{h}{\vert \ln h \vert}$}}}}
\put(2776,-1786){\makebox(0,0)[lb]{\smash{{\SetFigFont{9}{10.8}{\rmdefault}{\mddefault}{\updefault}$E_{0}$}}}}
\put(11776,-1786){\makebox(0,0)[lb]{\smash{{\SetFigFont{9}{10.8}{\rmdefault}{\mddefault}{\updefault}$E_{0}$}}}}
\put(19576,-1786){\makebox(0,0)[lb]{\smash{{\SetFigFont{9}{10.8}{\rmdefault}{\mddefault}{\updefault}$E_{0}$}}}}
\put(3001,-6361){\makebox(0,0)[b]{\smash{{\SetFigFont{9}{10.8}{\rmdefault}{\mddefault}{\updefault}{\rm (A)}}}}}
\put(12001,-6361){\makebox(0,0)[b]{\smash{{\SetFigFont{9}{10.8}{\rmdefault}{\mddefault}{\updefault}{\rm (B)}}}}}
\put(19801,-6361){\makebox(0,0)[b]{\smash{{\SetFigFont{9}{10.8}{\rmdefault}{\mddefault}{\updefault}{\rm (C)}}}}}
\end{picture}%
\end{center}
\caption{The resonance free zones for $\re \sigma$ large given by \eqref{b84}.} \label{f15}
\end{figure}

{\rm (A)} All the homoclinic trajectories are ``on the same side of $0$''. That is $g_{-} \cdot g_{+} > 0$ for all $g_{\pm} \in \CH^{\pm \infty}_{\rm tang}$. This condition is satisfied for instance in Example \ref{b86}, Section \ref{s19} {\rm (B)}, Example \ref{i69}, Example \ref{i72}, Example \ref{i80} and Example \ref{i86}. Under the condition {\rm (A)}, Remark \ref{b72} gives
\begin{equation*}
\CA_{0} ( \tau ) \lesssim e^{- \frac{\pi \tau}{2 \lambda}} \Big\vert \Gamma \Big( \frac{n}{2} - i \frac{\tau}{\lambda} \Big) \Big\vert .
\end{equation*}
From Stirling's formula \cite[(6.1.37)]{AbSt64_01}, it implies
\begin{equation} \label{b84}
\CA_{0} ( \tau ) \lesssim \left\{ \begin{aligned}
&e^{- \pi \frac{\tau}{\lambda}} \vert \tau \vert^{\frac{n-1}{2}}  &&\text{ for } \tau \gg 1 ,   \\
&\vert \tau \vert^{\frac{n-1}{2}}  &&\text{ for } \tau \ll - 1 .
\end{aligned} \right.
\end{equation}

Then, Theorem \ref{a2} gives a resonance free zone for the energies greater than $E_{0}$ of the form
\begin{equation} \label{i88}
\Big( - \pi \re \sigma + \lambda \frac{n-1}{2} \ln ( \re \sigma ) + C \Big) \frac{1}{\vert \ln h \vert} \leq \im \sigma \leq 0 ,
\end{equation}
for $\re \sigma \gg 1$ and some $C > 0$. Thus, the resonance free zone increases almost linearly, with a universal rate, as $\re \sigma$ goes to $+ \infty$. In the cited examples, the energies above $E_{0}$ are non-trapped and then at distance $h \vert \ln h \vert$ from the resonances. This can justify the phenomenon: the growth of the resonance free domain is a transition to non-trapping. Remark that \eqref{i88} is analogous to Remark \ref{e23} and coincides exactly with \eqref{e26} in the one dimensional case.

We consider the behavior of the accumulation curves for the energies below $E_{0}$. Since we do not have an explicit formula for the eigenvalues of $\widehat{\CT}$, we only discuss the numerical computations. In Figure \ref{f38}, we see that all the accumulation curves converge to the real axis in the limit $\re \sigma \to - \infty$. In particular, the number of these curves in $[ E_{0} - C h , E_{0} + C h ] + i [ - h \vert \ln h \vert^{- 1} , 0 ]$ seems to diverge as $C \to + \infty$. That the operator of Example \ref{i69} has a lot of degenerate elliptic periodic trajectories for energies below $E_{0}$ could explain this phenomenon. It can be interpreted as a transition to a strong trapping situation. Example \ref{i72} and Example \ref{i86} confirm this description and its interpretation. In the one dimensional case (see Section \ref{s19} {\rm (B)}), the unique accumulation curve is explicitly given by \eqref{e40} and converges to the real axis. Note that \eqref{b84} gives no information in this case since $\vert \tau \vert^{\frac{n-1}{2}} > 1$ for $\tau \ll - 1$ in dimension $n \geq 2$.

{\rm (B)} All the homoclinic trajectories are ``on the opposite side of $0$''. That is $g_{-} \cdot g_{+} < 0$ for all $g_{\pm} \in \CH^{\pm \infty}_{\rm tang}$. As explained after Corollary \ref{e4}, the Schr\"{o}dinger operators can not satisfy this condition when $\CH^{\pm \infty}_{\rm tang} \neq \emptyset$. However, Example \ref{b81} provides an instance of operator satisfying the condition {\rm (B)}. This one dimensional example could probably be generalized in any dimension.

The situation under condition {\rm (B)} is symmetric with respect to the case {\rm (A)}. We get a widening resonance free zone (resp. the convergence of the accumulation curves to the real axis) as $\re \sigma$ goes to $- \infty$ (resp. $+ \infty$). In particular, \eqref{i88} holds true with $\re \sigma$ replaced by $- \re \sigma$. Note that, in Example \ref{b81}, this is in agreement with the geometry of the trapped set at the energies close to $E_{0}$. Indeed, the energies below $E_{0}$ are non-trapped whereas the energies above $E_{0}$ have a well in an island trapping.

{\rm (C)} We have $g_{-} \cdot g_{+} > 0$ for some $g_{\pm} \in \CH^{\pm \infty}_{\rm tang}$ and $g_{-} \cdot g_{+} < 0$ for other $g_{\pm} \in \CH^{\pm \infty}_{\rm tang}$. This assumption is verified in Section \ref{s19} {\rm (C)} and Example \ref{b88}. Watching the numerical computations of Figure \ref{f17}, the set of the resonances for $\re \sigma \ll - 1$ seems to be the gathering of the sets of the resonances generated by the nappes $I_{1}$ and $I_{2}$ separately (see Figure \ref{f38}). As noted in Section \ref{s14}, this illustrates that the potential forms a barrier near $0$ and that the trapped set for energies below $E_{0}$ consists of two separated components similar to Example \ref{i69} on each side of $0$. On the contrary, the accumulation curves seem to combine two by two and to have a regular structure. This can be explained by the nature of the trapped set for energies higher than $E_{0}$. Such facts have already been observed in the transversal situation of Example \ref{e6}.

Under the condition {\rm (C)}, Remark \ref{b72} gives
\begin{equation}
\CA_{0} ( \tau ) \lesssim \left\{ \begin{aligned}
&\vert \tau \vert^{\frac{n-1}{2}}  &&\text{for } \tau \gg 1 ,   \\
&\vert \tau \vert^{\frac{n-1}{2}}  &&\text{for } \tau \ll - 1 .
\end{aligned} \right.
\end{equation}
With these upper bounds, Theorem \ref{a2} gives no resonance free zone for $\re \sigma$ large. This is optimal in Section \ref{s19} {\rm (C)} where no resonance free zone of size $h \vert \ln h \vert^{- 1}$ is possible.

\subsubsection{Stability phenomena} \label{s56}

The first think to note is that the resonances described by Theorem \ref{i55} are very sensitive to the action and the energy. For instance, it follows from \eqref{i61} that a small perturbation of the action $A$ in example \ref{i69} will produce a variation of the characteristic scale $h \vert \ln h \vert^{- 1}$ of the resonances. Formally, one can even consider perturbations of size $h$ of $A$. Numerically, it shows that this quantity must be very carefully computed. We have also seen in the previous Section \ref{s55} that the homoclinic situations are mostly transition regimes in energy. For these reasons, the setting of Theorem \ref{i55} may be seen as unstable.

Nevertheless, as in Section \ref{s24}, the asymptotic of the resonances is stable with respect to the homoclinic set. More precisely, we have

\begin{figure}
\begin{center}
\begin{picture}(0,0)%
\includegraphics{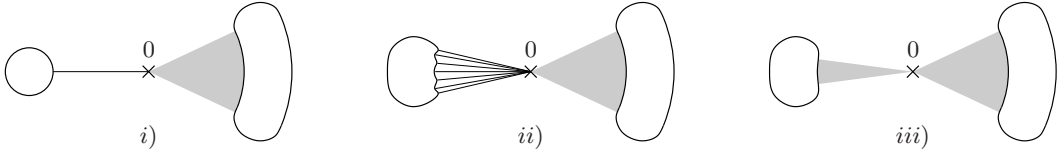}%
\end{picture}%
\setlength{\unitlength}{987sp}%
\begingroup\makeatletter\ifx\SetFigFont\undefined%
\gdef\SetFigFont#1#2#3#4#5{%
  \reset@font\fontsize{#1}{#2pt}%
  \fontfamily{#3}\fontseries{#4}\fontshape{#5}%
  \selectfont}%
\fi\endgroup%
\begin{picture}(26447,3632)(-12022,-7021)
\put(10801,-4786){\makebox(0,0)[b]{\smash{{\SetFigFont{9}{10.8}{\rmdefault}{\mddefault}{\updefault}$0$}}}}
\put(1201,-4786){\makebox(0,0)[b]{\smash{{\SetFigFont{9}{10.8}{\rmdefault}{\mddefault}{\updefault}$0$}}}}
\put(-8399,-4786){\makebox(0,0)[b]{\smash{{\SetFigFont{9}{10.8}{\rmdefault}{\mddefault}{\updefault}$0$}}}}
\put(1201,-6961){\makebox(0,0)[b]{\smash{{\SetFigFont{9}{10.8}{\rmdefault}{\mddefault}{\updefault}$ii)$}}}}
\put(-8399,-6961){\makebox(0,0)[b]{\smash{{\SetFigFont{9}{10.8}{\rmdefault}{\mddefault}{\updefault}$i)$}}}}
\put(10801,-6961){\makebox(0,0)[b]{\smash{{\SetFigFont{9}{10.8}{\rmdefault}{\mddefault}{\updefault}$iii)$}}}}
\end{picture}%
\end{center}
\caption{Three perturbations of Example \ref{i69} corresponding to the three levels of stability described in Remark \ref{i89}.} \label{f43}
\end{figure}

\begin{remark}[Stability]\sl \label{i89}
$i)$ The transversal homoclinic trajectories $\CH_{\rm trans}$ do not appear in the definition of the quantization operator $\CT$. Thus, they play no role in the pseudo-resonances and may contribute to the resonances only through the remainder term $o ( h \vert \ln h \vert^{- 1} )$. In other word, adding, removing or changing $\CH_{\rm trans}$ will have no impact on the (leading term of the) asymptotic of the resonances. Nevertheless, following Section \ref{s79}, it may be possible to see the lower order contribution of $\CH_{\rm trans}$.

$ii)$ The resonances are stable by perturbations of measure zero of the homoclinic tangential set. More precisely, let $Q$ be another operator satisfying the assumptions of Theorem \ref{i55}, whose set of asymptotic directions $\CH^{- \infty}_{\rm tang} ( Q )$ coincides with $\CH^{- \infty}_{\rm tang} ( P )$ modulo a set of measure zero, and such that the dynamical quantities of $P$ and $Q$ are equal on $\CH^{- \infty}_{\rm tang} ( P ) \cap \CH^{- \infty}_{\rm tang} ( Q )$. Then, $L^{2} ( \CH^{- \infty}_{\rm tang} ( P ) ) = L^{2} ( \CH^{- \infty}_{\rm tang} ( Q ) )$, the operators $\CT$ associated to $P$ and $Q$ are the same, and the resonances of $P$ and $Q$ only differ by $o ( h \vert \ln h \vert^{- 1} )$.

$iii)$ Lastly, a small perturbation of the homoclinic set induces a small perturbation of the resonances. In this direction, we state the following result which can be generalized: Let $C , \delta > 0$ and $P$ be a fixed operator. We consider another operator $Q$ and assume, for both $P$ and $Q$, that the assumptions of Theorem \ref{i55} are satisfied, $\alpha ( \omega ) = \omega$, $\vert g_{-} \cdot g_{+} \vert > C^{-1}$ for $g_{\pm} \in \CH^{+ \infty}_{\rm tang}$ and that the geometric quantities $A , T , \CM_{0} , \nu$ are constant on $\CH^{+ \infty}_{\rm tang}$ and do not depend on the operator. Then, if
\begin{equation*}
\mes_{\S^{n - 1}} \big( \CH^{+ \infty}_{\rm tang} ( P ) \Delta \CH^{+ \infty}_{\rm tang} ( Q ) \big) ,
\end{equation*}
is small enough, we have in the domain \eqref{i54}
\begin{equation} \label{i90}
\dist \big( \res (P) , \res (Q) \big) \leq \delta \frac {h}{\vert \ln h \vert} ,
\end{equation}
for $h$ small enough. Here, $M \Delta N = M \setminus N \cup N \setminus M$ denotes the symmetric difference of two sets $M , N$.
\end{remark}

As corollary, Remark \ref{i89} $iii)$ implies that the resonances depend continuously on the opening angle $\theta_{0}$ in Example \ref{i69}. It also shows that one can replace $\CH^{- \infty}_{\rm tang}$ by a nicer compact set in the numerical computations of the eigenvalues of $\CT$. Figure \ref{f43} provides an example for each of the three stabilities of the previous remark.

\section{Generalization to multiple barriers} \label{s42}

In this part, the geometric setting of Section \ref{s2} is generalized. We allow for a finite number of hyperbolic fixed points and heteroclinic orbits between them. A general result is stated in Section \ref{s45} and applied in Section \ref{s89}. The corresponding proofs are given in Section \ref{s46} and Section \ref{s44}.

\Subsection{General result} \label{s45}

As above, let $P = - h^{2} \Delta + V ( x )$ and $E_{0} > 0$. In order to define the resonances through analytic distortion, we still assume \ref{h1}. However, we replace the assumption \ref{h2} by
\begin{hyp} \label{h15}
There exists a finite set $\SV \subset \R^{n}$ such that the potential $V ( x )$ has a non-degenerate maximum at $x = v$ with critical value $E_{0}$ for all $v \in \SV$.
\end{hyp}
In particular, one can find local coordinates such that 
\begin{equation*}
V ( x ) = E_{0} - \sum_{j=1}^{n} \frac{( \lambda_{j}^{v} )^{2}}{4} ( x - v )_{j}^{2} + \CO \left( ( x - v )^{3} \right) ,
\end{equation*}
near $v \in \SV$ with $0 < \lambda_{1}^{v} \leq \lambda_{2}^{v} \leq \dots \leq \lambda_{n}^{v}$. As in Section \ref{s2}, $( v , 0 )$ is a hyperbolic fixed point of the Hamiltonian vector field $H_{p}$. The corresponding incoming/outgoing Lagrangian manifold is denoted by $\Lambda_{\pm}^{v}$.

\begin{figure}
\begin{center}
\begin{picture}(0,0)%
\includegraphics{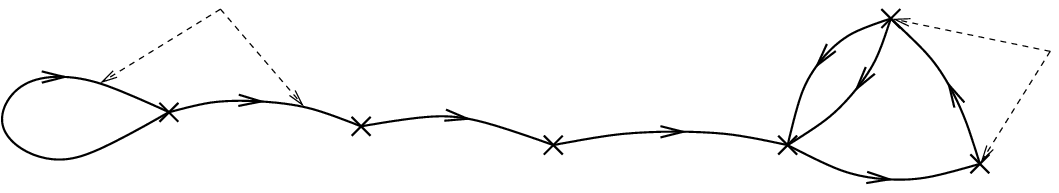}%
\end{picture}%
\setlength{\unitlength}{1184sp}%
\begingroup\makeatletter\ifx\SetFigFont\undefined%
\gdef\SetFigFont#1#2#3#4#5{%
  \reset@font\fontsize{#1}{#2pt}%
  \fontfamily{#3}\fontseries{#4}\fontshape{#5}%
  \selectfont}%
\fi\endgroup%
\begin{picture}(17193,3408)(-3527,-4996)
\put(13651,-2686){\makebox(0,0)[b]{\smash{{\SetFigFont{9}{10.8}{\rmdefault}{\mddefault}{\updefault}$\SV$}}}}
\put(  1,-1711){\makebox(0,0)[b]{\smash{{\SetFigFont{9}{10.8}{\rmdefault}{\mddefault}{\updefault}$\SE$}}}}
\put(7201,-3586){\makebox(0,0)[b]{\smash{{\SetFigFont{9}{10.8}{\rmdefault}{\mddefault}{\updefault}$e$}}}}
\put(3826,-3361){\makebox(0,0)[b]{\smash{{\SetFigFont{9}{10.8}{\rmdefault}{\mddefault}{\updefault}$\widetilde{e}$}}}}
\put(12151,-4936){\makebox(0,0)[b]{\smash{{\SetFigFont{9}{10.8}{\rmdefault}{\mddefault}{\updefault}$\widetilde{v}$}}}}
\put(10726,-1711){\makebox(0,0)[b]{\smash{{\SetFigFont{9}{10.8}{\rmdefault}{\mddefault}{\updefault}$v$}}}}
\put(2326,-4411){\makebox(0,0)[b]{\smash{{\SetFigFont{9}{10.8}{\rmdefault}{\mddefault}{\updefault}$\widetilde{e}^{-}$}}}}
\put(9151,-4711){\makebox(0,0)[b]{\smash{{\SetFigFont{9}{10.8}{\rmdefault}{\mddefault}{\updefault}$e^{+}$}}}}
\put(5326,-4636){\makebox(0,0)[b]{\smash{{\SetFigFont{9}{10.8}{\rmdefault}{\mddefault}{\updefault}$\widetilde{e}^{+} = e^{-}$}}}}
\end{picture}%
\end{center}
\caption{The graph $( \SV , \SE )$ and the notations of Section \ref{s45}.} \label{f49}
\end{figure}

As in \ref{h3}, we assume
\begin{hyp} \label{h16}
The trapped set at energy $E_0$ satisfies
\begin{equation*}
K ( E_{0} ) \subset \bigcup_{v , \widetilde{v} \in \SV} \Lambda_{-}^{v} \cap \Lambda_{+}^{\widetilde{v}} .
\end{equation*}
\end{hyp}
Let $\SE$ denote the set of non-constant Hamiltonian trajectories in $K ( E_{0} )$. This means that a bicharacteristic trajectory $e$ belongs to $\SE$ if and only if $e \subset K ( E_{0} )$ and $e \neq \{ ( v , 0 ) \}$ for all $v \in \SV$. From \ref{h16}, any curve $e \in \SE$ is either a {\it homoclinic trajectory} (if $e \subset \Lambda_{-}^{v} \cap \Lambda_{+}^{v}$ for some $v \in \SV$) or a {\it heteroclinic trajectory} (if $e \subset \Lambda_{-}^{v} \cap \Lambda_{+}^{\widetilde{v}}$ for some $v , \widetilde{v} \in \SV$ with $v \neq \widetilde{v}$). The setting is illustrated in Figure \ref{f49}.

The pair $( \SV , \SE )$ can be regarded as a graph with {\it vertices} $v \in \SV$ and {\it edges} $e \in \SE$. The edges are oriented by the time evolution of the Hamiltonian flow. For an edge $e$, we write $e^{-} \in \SV$ the {\it origin} of $e$ and $e^{+} \in \SV$ the {\it terminus} of $e$. It is allowed to have loops $e$ with $e^{-} = e^{+}$, which correspond to homoclinic trajectories. It is allowed that, for a pair of vertices $( \widetilde{v} , v )$, there are more than one edge $e$ satisfying $e^{-} = \widetilde{v}$ and $e^{+} = v$. Therefore the graph $( \SV , \SE )$ is a directed multigraph or pseudograph. At this stage, there is no assumption on the cardinal of $\SE$ which can be uncountable.

A finite sequence of edges $(e_{1} , \ldots , e_{k} )$ is called a {\it path} when $e_{\ell}^{+} = e_{\ell + 1}^{-}$ for $\ell = 1 , \ldots, k - 1$, and a {\it cycle} when moreover $e_{k}^{+} = e_{1}^{-}$. In the sequel, a cycle $(e_{1} , \ldots , e_{k} )$ is identified with $(e_{2} , \ldots , e_{k} , e_{1} )$ and so on. A cycle is called {\it primitive} if it does not contain any non-trivial sub-cycle. In particular, a primitive cycle is either a single homoclinic trajectory or a sequence of (at least two) heteroclinic trajectories. We shall see that the cycles play the main role in the creation of resonances. To start with, we prove that $P$ has few resonances if there is no cycle. More precisely, let $\Gamma ( h )$ be the union over all the vertices in $\SV$ of the exceptional sets $\Gamma_{v} ( h )$ defined in Section \ref{m81} and assume
\begin{hyp} \label{h17}
There is no cycle in the graph $( \SV , \SE )$.
\end{hyp}
Under the previous assumptions, the following fact holds.

\begin{theorem}\sl \label{j73}
Assume \ref{h1}, \ref{h15}, \ref{h16} and \ref{h17}. For any $C , \delta > 0$, there is no resonance in the domain
\begin{equation} \label{j72}
E_{0} + [ - C h , C  h ] + i [ - C h , 0 ] \setminus \big( \Gamma ( h ) + B ( 0 , \delta h ) \big) ,
\end{equation}
for $h$ small enough. Moreover, for all $\chi \in C^{\infty}_{0} ( \R^{n} )$, there exists $M > 0$ such that
\begin{equation*}
\big\Vert \chi ( P - z )^{-1} \chi \big\Vert \lesssim h^{- M} ,
\end{equation*}
uniformly for $h$ small enough and $z \in \eqref{j72}$.
\end{theorem}

In the Schr\"{o}dinger case, the assumptions of Theorem \ref{j73} imply that $\card \SV \leq 1$ and $\SE = \emptyset$ (see Lemma \ref{k57}). In other words, the trapped set of $P$ at energy $E_{0}$ is either empty or consist of a hyperbolic fixed point. In these situations and for analytic potentials, the asymptotic of resonances has already been obtain by Helffer and Sj\"{o}strand \cite{HeSj86_01}, and Sj\"{o}strand \cite{Sj87_01} respectively. Nevertheless, it is possible to construct more interesting geometries by considering more general operators (see Remark \ref{c13}).

Using the strategy developed in Section \ref{s36} and the test functions constructed in \cite[Section 4]{BoFuRaZe11_01}, one could  show that $P$ has at least one resonance near each element of $\Gamma_{0} ( h )$, the union over the vertices in $\SV$ of the set $\Gamma_{0 , v} ( h )$ defined in \eqref{m80}. In this setting, one might even prove that
\begin{equation*}
\dist \big( \res ( P ) , \Gamma_{0} ( h ) \big) = o ( h ) ,
\end{equation*}
in the domain $B ( E_{0} , C h )$.

To the contrary of \ref{h17}, let us now suppose that
\begin{hyp} \label{h20}
The graph $( \SV , \SE )$ has at least one cycle.
\end{hyp}
We give the asymptotic of resonances under assumptions similar to those of Section \ref{s61}. Adapting the corresponding results from the homoclinic setting, we could have also considered general resonance free domains, tangential intersections of finite order, nappes of heteroclinic trajectories, $\ldots$ As in \ref{h4}, we suppose that
\begin{hyp} \label{h18}
For any $e , \widetilde{e} \in \SE$ with $e^{-} = \widetilde{e}^{+}$, one has $g^{e}_{+} \cdot g^{\widetilde{e}}_{-} \neq 0$.
\end{hyp}
In the previous expression, the asymptotic directions $g^{e}_{\pm} \in \R^{n}$ of the trajectory $e \in \SE$ are defined by
\begin{equation*}
\pi_{x} ( e ( t ) ) = e^{\mp} + g_{\pm}^{e} e^{\pm t \lambda_{1}^{e^{\mp}}} + o \big( e^{\pm t \lambda_{1}^{e^{\mp}}} \big) ,
\end{equation*}
as $t \to \mp \infty$ (see \eqref{d2}). Eventually, we assume that
\begin{hyp} \label{h19}
For any edge $e \in \SE$, the manifolds $\Lambda_{+}^{e^{-}}$ and $\Lambda_{-}^{e^{+}}$ intersect transversally along $e$.
\end{hyp}
In particular, the graph $( \SV , \SE )$ and the set of primitive cycles are finite.

We define a quantity which measures the importance of each cycle. For all $v \in \SV$, we set
\begin{equation*}
\alpha_{v} = \frac{1}{2} \sum_{j = 2}^{n} \frac{\lambda_{j}^{v}}{\lambda_{1}^{v}} \qquad \text{and} \qquad \beta_{v} = \frac{1}{\lambda_{1}^{v}} .
\end{equation*}
For a cycle $\gamma$, let $\SV ( \gamma )$ be the set of vertices belonging to an edge in $\gamma$. Then we denote
\begin{equation} \label{m12}
\alpha  (\gamma ) = \sum_{v \in \SV ( \gamma )} \alpha_{v} , \qquad \beta ( \gamma ) = \sum_{v \in \SV ( \gamma )} \beta_{v} ,
\end{equation}
and finally, we define the {\it damping} of the cycle $\gamma$ to be
\begin{equation} \label{j77}
D ( \gamma ) = \frac{\alpha ( \gamma )}{\beta ( \gamma )} .
\end{equation}
What we call damping truly measures the capacity of a cycle to dissipate the energy. Roughly speaking, $\gamma$ ``tries'' to create resonances satisfying $\im z \approx - D ( \gamma ) h$. Indeed, the quantization rule associated to $\gamma$ would be $h^{\alpha ( \gamma ) - i \beta ( \gamma ) \frac{z - E_{0}}{h}} \approx 1$. Thus, the cycles with the smallest damping play the main role in the asymptotic of the resonances closest to the real axis. We define the real number
\begin{equation} \label{m1}
D_{0} = \min_{\gamma \text{ cycle}} D ( \gamma ) .
\end{equation}
That the infimum is attained is proved in Section \ref{s50}. Moreover, it is enough to minimize over the primitive cycles of the graph $( \SV , \SE )$. A cycle $\gamma$ is called {\it minimal} if $D ( \gamma ) = D_{0}$. Note that $D ( \gamma )= 0$ if and only if $n = 1$. In particular, the notion of damping is relevant only in dimension $n \geq 2$.

In order to state the quantization rule, we define some dynamical quantities. For any $e \in \SE$, let $A_{e} = \int_{e} \xi \cdot d x$ be its action and let $\nu_{e}$ be the Maslov index of $\Lambda_{+}^{e^{-}}$ along $e$. We also recall the definition \eqref{d7} for the Maslov determinant $\CM_{e}^{\pm}$. Then, for all $e , \widetilde{e} \in \SE$ with $e^{-} = \widetilde{e}^{+}$, we set the scalar function
\begin{equation} \label{j83}
\CQ_{e , \widetilde{e}} ( z , h ) = e^{i A_{e} / h} \Gamma \Big( \frac{S_{v}}{\lambda^{v}_{1}} \Big) \sqrt{\frac{\lambda^{v}_{1}}{2 \pi}} \frac{\CM_{e}^{+}}{\CM_{e}^{-}} e^{- \frac{\pi}{2} ( \nu_{e} + \frac{1}{2} ) i} \big\vert g_{-}^{\widetilde{e}} \big\vert \big( i \lambda_{1}^{v} g_{+}^{e} \cdot g_{-}^{\widetilde{e}} \big)^{- S_{v} / \lambda^{v}_{1}} ,
\end{equation}
which encodes the transition from $\widetilde{e}$ to $e$, where
\begin{equation*}
S_{v} = S_{v} ( z , h ) = \sum_{j = 1}^{n} \frac{\lambda_{j}^{v}}{2} - i \frac{z - E_{0}}{h} ,
\end{equation*}
is a new spectral parameter associated with the hyperbolic fixed point $v = e^{-}$. As in \eqref{d10}, this function can be written
\begin{equation} \label{m20}
\CQ_{e , \widetilde{e}} ( z , h ) = \rho_{e} \widetilde{\CQ}_{e , \widetilde{e}} ( \sigma ) ,
\end{equation}
where $\rho_{e} = e^{i A_{e} / h} \in \S^{1}$, the complex variable $\sigma$ is given by \eqref{d92} and $\widetilde{\CQ}_{e , \widetilde{e}}$ is a meromorphic function of $\sigma$ independent of $h$. The asymptotic of the resonances is governed by the $\card \SE \times \card \SE$ matrix whose coefficients are
\begin{equation} \label{m2}
\SQ_{e , \widetilde{e}} ( z , h ) = \left\{
\begin{aligned}
&h^{S_{v} ( z , h ) / \lambda_{1}^{v} - 1 / 2} \CQ_{e , \widetilde{e}} ( z , h ) &&\text{ if } e^{-} = \widetilde{e}^{+} = : v ,  \\
&0 &&\text{ if } e^{-} \neq \widetilde{e}^{+} .
\end{aligned} \right.
\end{equation}
Under the previous assumptions, the set of pseudo-resonances is defined by

\begin{definition}[Quantization rule]\sl \label{j82}
We say that $z$ is a pseudo-resonance if and only if
\begin{equation*}
1 \in \spe ( \SQ ( z , h ) )  .
\end{equation*}
The set of pseudo-resonances is denoted by $\res_{0} ( P )$.
\end{definition}

This definition is similar to Definition \ref{d1}. The main difference is that the quantity $h^{S ( z , h ) / \lambda_{1} - 1 / 2}$ depends now on the pair of edges $( e , \widetilde{e} )$ and can not be factored out. This phenomenon changes drastically the asymptotic of the pseudo-resonances, which is not as simple as in Proposition \ref{d9}. We also define the matrix $\SQ^{\rm p r i n} ( z , h )$ whose entries are
\begin{equation} \label{m3}
\SQ^{\rm p r i n}_{e , \widetilde{e}} ( z , h ) = \left\{
\begin{aligned}
&h^{S_{v} ( z , h ) / \lambda_{1}^{v} - 1 / 2} \CQ_{e , \widetilde{e}} ( z , h ) &&\text{ if } e^{-} = \widetilde{e}^{+} = : v \text{ and } \widetilde{e} , e \\
& &&\text{  belong to minimal cycles} ,  \\
&0 &&\text{ otherwise} .
\end{aligned} \right.
\end{equation}
Let $F , F^{\rm p r i n}$ be the complex functions
\begin{equation*}
F ( z , h ) = \det ( 1 - \SQ ( z , h ) ) \qquad \text{and} \qquad F^{\rm p r i n} ( z , h ) = \det ( 1 - \SQ^{\rm p r i n} ( z , h ) ) .
\end{equation*}
By definition, the pseudo-resonances are the zeros of $F$. For $\tau \in \R$, let $Z = Z ( z , \tau )$ be the new spectral parameter defined by
\begin{equation} \label{m19}
z = E_{0} + \tau h - i D_{0} h + Z \frac{h}{\vert \ln h \vert} .
\end{equation}
Expanding the determinant, the functions $F , F^{\rm p r i n}$ can be written as follows. The proof of these identities and additional informations can be found in Section \ref{s47}.

\begin{remark}\sl \label{m4}
There exists a finite set $\SF \subset [ 0 , + \infty [^{2}$ such that
\begin{equation} \label{m5}
F ( z , h ) = \sum_{( \alpha , \beta ) \in \SF} h^{\alpha} e^{i \beta Z} F_{\alpha , \beta} ( \tau , z , h ) .
\end{equation}
For every $( \alpha , \beta ) \in \SF$, the number $\alpha$ (resp. $\beta$) is a finite sum of $\alpha ( \gamma ) - D_{0} \beta ( \gamma )$ (resp. $\beta ( \gamma )$) for some primitive cycles $\gamma$. In particular, $\alpha = 0$ if and only if all these cycles are minimal.  Moreover, $F_{\alpha , \beta}$ is a finite sum of finite products of $\pm h^{- i \beta_{v} \tau} \CQ_{e , \widetilde{e}} ( z , h )$. Thus, we have
\begin{equation} \label{m23}
F_{\alpha , \beta} ( \tau , z , h ) = \widetilde{F}_{\alpha , \beta} ( \kappa , \rho , \sigma ) ,
\end{equation}
where $\widetilde{F}_{\alpha , \beta}$ is independent of $h$ and holomorphic in $\kappa , \rho , \sigma$ for $\sigma$ outside $( \Gamma_{0} ( h ) - E_{0} ) h^{- 1}$ and $\kappa = ( \kappa_{v} )_{v \in \SV} \in \C^{\card \SV}$ with $\kappa_{v} = h^{- i \beta_{v} \tau} \in \S^{1}$. In particular, $F_{0 , 0} ( \tau , z , h ) = 1$. On the other hand, $F^{\rm p r i n}$ is the sum of the leading terms in \eqref{m5}. More precisely, it means that
\begin{equation} \label{m6}
F^{\rm p r i n} ( z , h ) = \sum_{\beta \in \SB} e^{i \beta Z} F_{0 , \beta} ( \tau , z , h ) ,
\end{equation}
where $\SB = \{ \beta ; \ ( 0 , \beta ) \in \SF \}$.
\end{remark}

Mimicking \eqref{e18}, we define
\begin{equation} \label{m27}
f_{\tau} ( Z , h ) = \sum_{\beta \in \SB} e^{i \beta Z} F_{0 , \beta} ( \tau , E_{0} + \tau h - i D_{0} h , h ) = \sum_{\beta \in \SB} e^{i \beta Z} \widetilde{F}_{0 , \beta} ( \kappa , \rho , \tau - i D_{0} ) ,
\end{equation}
the function $F^{\rm p r i n}$ with the slow variable $\sigma$ fixed at $\tau - i D_{0}$. For $\tau$ and $h$ fixed, $f_{\tau} ( Z , h )$ is an exponential sum which does not vanish identically. The distribution of the zeros of such functions has been studied for a long time (see e.g. Langer \cite{La31_01} or the book of Bellman and Cooke \cite{BeCo63_01}). In particular, it is known that the number of zeros is uniformly bounded on compact sets of given size. In our setting, we have the following estimates proved in Section \ref{s47}. For $C , \delta > 0$, there exists a constant $N > 0$ such that, for all $\tau \in [ - C , C ]$ and $h \in ] 0 , 1 ]$ with $\dist ( E_{0} + \tau h - i D_{0} h , \Gamma_{0} ( h ) ) \geq \delta h$, we have 
\begin{equation} \label{m9}
\big\{ Z  \in \R + i [- C , + \infty [ ; \ f_{\tau} ( Z , h ) = 0 \big\} \subset \R + i [ - C , N ] ,
\end{equation}
and
\begin{equation} \label{m8}
\card \big\{ Z  \in [ - C , C ] + i [ - C , N ] ; \ f_{\tau} ( Z , h ) = 0 \big\} \leq N .
\end{equation}
In the previous expression, the zeros are counted with their multiplicity. Note that these zeros are continuous with respect to $\tau , h$ thanks to Hurwitz's theorem, and analytic except where the multiplicity changes. Lastly, \cite[Theorem 12.7]{BeCo63_01} proves the existence of zeros of $f_{\tau} ( \cdot , h )$ and provides the asymptotic of their number for $\tau , h$ fixed under some natural assumptions.

\begin{proposition}[Asymptotic of the pseudo-resonances]\sl \label{j81}
Assume \ref{h1}, \ref{h15}, \ref{h16}, \ref{h20}--\ref{h19} and let $C , \delta > 0$. Then, uniformly for $\tau \in [ - C , C ]$, the pseudo-resonances $z$ in
\begin{equation} \label{m10}
E_{0} + \tau h + \Big[ - C \frac{h}{\vert \ln h \vert} , C \frac{h}{\vert \ln h \vert} \Big] + i \Big[ - D_{0} h - C \frac{h}{\vert \ln h \vert} , h \Big] \setminus \big( \Gamma_{0} ( h ) + B ( 0 , \delta h ) \big) ,
\end{equation}
satisfy
\begin{equation} \label{m11}
z = E_{0} + \tau h - i D_{0} h + Z \frac{h}{\vert \ln h \vert} + o \Big( \frac{h}{\vert \ln h \vert} \Big) ,
\end{equation}
where $Z$ is a zero of $f_{\tau} ( \cdot , h )$. Conversely, for each $\tau \in [ - C , C ]$ and $Z \in [ - C , C ] + i [ - C , + \infty [$ such that $f_{\tau} ( Z , h ) = 0$, there exists a pseudo-resonance $z$ satisfying \eqref{m11} uniformly with respect to $\tau , Z$.
\end{proposition}

This description is close to Definition \ref{j82} and still implicit. By comparison with the previous asymptotic of the pseudo-resonances (see e.g. Proposition \ref{d9}), the present result corresponds to $\delta ( h ) = C \vert \ln h \vert^{- 1}$. It may perhaps be possible to consider larger zones since the zeros of exponential sums like \eqref{m27} enjoy nice properties (see e.g. \cite[Chapter 12]{BeCo63_01}). Note also that the pseudo-resonances are not expressed in terms of the spectrum of an operator. Using Definition \ref{g80}, the resonances are given by the following theorem in the case of many critical points.

\begin{theorem}[Asymptotic of resonances]\sl \label{j79}
Assume \ref{h1}, \ref{h15}, \ref{h16}, \ref{h20}--\ref{h19} and let $C , \delta > 0$. In the domain
\begin{equation} \label{j80}
E_{0} + [ - C h , C h ] + i \Big[ - D_{0} h - C \frac{h}{\vert \ln h \vert} , h \Big] \setminus \big( \Gamma (h) + B ( 0 , \delta h ) \big),
\end{equation}
we have
\begin{equation*}
\dist \big( \res (P) , \res_{0} (P) \big) = o \Big( \frac {h}{\vert \ln h \vert} \Big) ,
\end{equation*}
as $h$ tends to $0$. Moreover, for all $\chi \in C^{\infty}_{0} ( \R^{n} )$, there exists $M > 0$ such that
\begin{equation*}
\big\Vert \chi ( P -z )^{-1} \chi \big\Vert \lesssim h^{- M} ,
\end{equation*}
uniformly for $h$ small enough and $z \in \eqref{j80}$ with $\dist ( z , \res_{0} ( P ) ) \geq \delta h \vert \ln h \vert^{- 1}$.
\end{theorem}

\begin{figure}
\begin{center}
\begin{picture}(0,0)%
\includegraphics{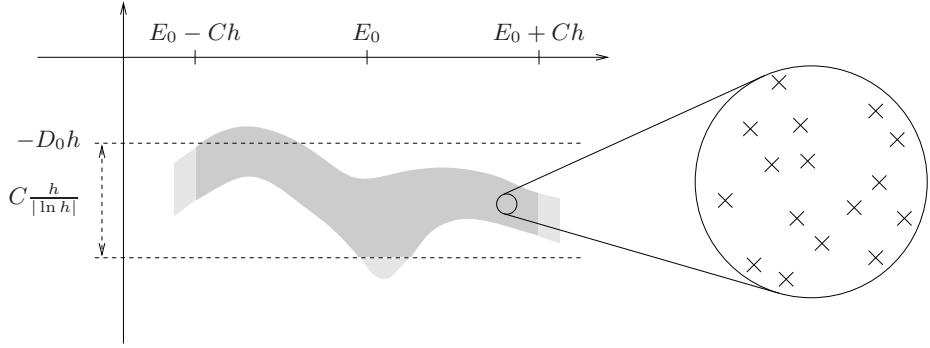}%
\end{picture}%
\setlength{\unitlength}{1184sp}%
\begingroup\makeatletter\ifx\SetFigFont\undefined%
\gdef\SetFigFont#1#2#3#4#5{%
  \reset@font\fontsize{#1}{#2pt}%
  \fontfamily{#3}\fontseries{#4}\fontshape{#5}%
  \selectfont}%
\fi\endgroup%
\begin{picture}(18665,7244)(-621,-6983)
\put(6301,-586){\makebox(0,0)[b]{\smash{{\SetFigFont{9}{10.8}{\rmdefault}{\mddefault}{\updefault}$E_{0}$}}}}
\put(2701,-586){\makebox(0,0)[b]{\smash{{\SetFigFont{9}{10.8}{\rmdefault}{\mddefault}{\updefault}$E_{0} - C h$}}}}
\put(9901,-586){\makebox(0,0)[b]{\smash{{\SetFigFont{9}{10.8}{\rmdefault}{\mddefault}{\updefault}$E_{0} + C h$}}}}
\put(301,-2836){\makebox(0,0)[rb]{\smash{{\SetFigFont{9}{10.8}{\rmdefault}{\mddefault}{\updefault}$- D_{0} h$}}}}
\put(301,-4036){\makebox(0,0)[rb]{\smash{{\SetFigFont{9}{10.8}{\rmdefault}{\mddefault}{\updefault}$C \frac{h}{\vert \ln h \vert}$}}}}
\end{picture}%
\end{center}
\caption{The cloud of resonances of Theorem \ref{j79}.} \label{f50}
\end{figure}

Theorem \ref{j79} implies that the typical size of the imaginary part of the resonances closest to the real axis is $- D_{0} h$. More precisely, combining this result with \eqref{m9} and Proposition \ref{j81}, we deduce that $P$ has no resonance in the domain
\begin{equation} \label{k19}
E_{0} + [ - C h , C  h ] + i \Big[ - D_{0} h + N \frac{h}{\vert \ln h \vert} , h \Big] \setminus \big( \Gamma ( h ) + B ( 0 , \delta h ) \big) ,
\end{equation}
for $h$ small enough. On the other hand, Example \ref{k42}, Example \ref{k61} and Section \ref{s87} below show that there exist situations with (a lot of) resonances satisfying $\im z \approx - D_{0} h$.

Any operator $P$ satisfying the assumptions of Theorem \ref{d8} also verifies those of Theorem \ref{j79}. Thus, this part is a generalization of Section \ref{s61}. Nevertheless, we have here neither the lower bound on the multiplicity of the resonances of Proposition \ref{d91} nor the asymptotic of the pseudo-resonances in large zones of Proposition \ref{d9}. That we used slightly different proofs is the reason for this. Furthermore, the other results of Section \ref{s6} (as the resonances in deeper zones of Section \ref{s43} or the asymptotic of higher order of Section \ref{s79}) will be more technical and harder to understand in the present setting.

Note that $D_{0}$ can be seen as the barycenter over any minimal cycle of the $\alpha_{v} \beta_{v}^{- 1}$'s with the weights $\beta_{v}$. Moreover, we have already shown in Section \ref{s61} that the imaginary part of the resonances is typical like $- \alpha_{0} \beta_{0}^{- 1} h$ when the potential has a unique non-degenerate maximum at $0$ (i.e. $\SV = \{ 0 \}$). Thus, the leading term of the imaginary part of the present resonances is the average of the corresponding quantity associated to each vertex.

\begin{remark}[Structure of the set of resonances]\sl \label{k35}
In some sense, the resonances described in Proposition \ref{j81} and Theorem \ref{j79} still satisfy a two scale asymptotic as in Remark \ref{e2}. Roughly speaking, $\tau$ can be seen as the macroscopic variable, whereas $Z$ plays the role of the microscopic variable. But the situation can now be very different. The resonances do not necessarily concentrate on accumulation curves but can form a ``cloud'' of resonances.

More precisely, if all the $\beta \in \SB$ are $\Z$-related, the function $f_{\tau} ( Z , h )$ defined in \eqref{m27} is polynomial in $e^{i \beta_{0} Z}$ for some $\beta_{0} \in \R \setminus \{ 0 \}$. Let $( \mu_{k} ( \tau , h )^{- 1} )_{k}$ denote the zeros of this polynomial. Then, the pseudo-resonances in \eqref{m10} satisfy
\begin{equation} \label{k39}
z = E_{0} + \tau h + 2 q \pi \beta_{0}^{- 1} \frac{h}{\vert \ln h \vert} - i D_{0} h + i \ln ( \mu_{k} ( \tau , h ) ) \beta_{0}^{- 1} \frac{h}{\vert \ln h \vert} + o \Big( \frac{h}{\vert \ln h \vert} \Big) ,
\end{equation}
for some $k$ and some $q \in \Z$, in the sense of Proposition \ref{d9}. This formula is similar to \eqref{d95}. In particular, the resonances concentrate on accumulation curves. Section \ref{s61} and Section \ref{s71} below provide examples of such situations. There is a priori no link between the number of accumulations curves and the number of cycles or edges.

On the other hand, if the elements of $\SB$ are not $\Z$-related, $f_{\tau} ( Z , h )$ is an exponential sum of general form and its zeros have no periodic structure similar to \eqref{k39}. In this situation, we say that the resonances form a cloud (see Figure \ref{f50}). This is for instance the case in Example \ref{k52} below.
\end{remark}

Even if the resonances do not accumulate on curves, the other phenomena of Section \ref{s18} still hold true in the presence of several barriers. First, the resonances oscillate with $h$. More precisely, assume first that the actions $A_{e}$, with $e \in \SE$ in a minimal cycle, coincide. Then, $f_{\tau} ( \cdot , h )$ and the parameter $Z$ appearing in \eqref{m11} are periodic with respect to $h^{- 1}$. In the general case, these quantities are continuous functions of the $e^{i A_{e} / h}$'s. The same way, the transition phenomena of Section \ref{s14} and the stability properties of Section \ref{s24} can be adapted to the present setting. In fact, Theorem \ref{j79} is already a stability result where the edges of the non-minimal cycles act as negligible perturbations. Indeed, Proposition \ref{j81} implies that these edges do not appear at the main order in the asymptotic of resonances. Nevertheless, if we had a description of the resonances modulo $\CO ( h^{\infty} )$ (by adapting Section \ref{s79} for instance), we should see the influence of the non-minimal cycles. In some sense, the situation is similar to the one of Section \ref{s81} (see the discussion below Remark \ref{k62}).

\Subsection{Applications and examples} \label{s89}

Here, we give the asymptotic of the resonances in various situations using the results of the previous part. The proof of these corollaries can be found in Section \ref{s44}.

\subsubsection{A unique minimal primitive cycle} \label{s71}

In addition to the hypotheses of Theorem \ref{j79}, we assume here that there exists a unique minimal primitive cycle, denoted $\gamma_{0} = ( e_{1} , \ldots , e_{K} )$ in the sequel. This assumption is somehow generic. We have $D_{0} = D ( \gamma_{0} )$ and we define
\begin{equation} \label{k44}
\mu ( \tau , h ) = \prod_{k = 1}^{K} \CQ_{e_{k + 1} , e_{k}} ( E_{0} + \tau h - i D_{0} h , h ) ,
\end{equation}
with the convention $e_{K + 1} = e_{1}$. This quantity can be made explicit using \eqref{j83}. In this case, Theorem \ref{j79} implies the following asymptotic.

\begin{corollary}\sl \label{k40}
In the present setting, let $C , \delta > 0$. The resonances $z$ lying in \eqref{j80} satisfy, in the sense of Proposition \ref{j81},
\begin{equation*}
z = z_{q} ( \tau ) + o ( h \vert \ln h \vert^{- 1} ) ,
\end{equation*}
for some $q \in \Z$ with
\begin{equation} \label{k41}
z_{q} ( \tau ) = E_{0} + \frac{2 q \pi }{\beta ( \gamma_{0} )} \frac{h}{\vert \ln h \vert} - i D_{0} h + i \frac{\ln ( \mu ( \tau , h ) )}{\beta ( \gamma_{0} )}\frac{h}{\vert \ln h \vert} .
\end{equation}
\end{corollary}

This situation generalizes that of a unique homoclinic trajectory treated in Section \ref{s17}. Roughly speaking, the cycle $\gamma_{0}$ counts as a single bicharacteristic curve and the quantization rule becomes scalar (see \eqref{k51}). In particular, \eqref{k41} is similar to \eqref{e8}. Moreover, when $\gamma_{0}$ is reduced to a single homoclinic trajectory, the quantization rule and the asymptotic of the resonances coincide with those of Section \ref{s61}.

This result implies that the resonances accumulate on the unique curve
\begin{equation} \label{k46}
\im \sigma = - i D_{0} + \frac{\ln ( \vert \mu ( \re \sigma , h ) \vert )}{\beta ( \gamma_{0} )}\frac{1}{\vert \ln h \vert} ,
\end{equation}
where $\sigma$ is given by \eqref{d92}. This is coherent with Remark \ref{k35} since $\SB = \{ \beta ( \gamma_{0} ) \}$. As explained in the following remark, there is only few possible geometries for $\gamma_{0}$ in the case of Schr\"{o}dinger operators.

\begin{remark}\sl \label{j75}
In the Schr\"{o}dinger case $P = - h^{2} \Delta + V ( x )$, the assumptions of Corollary \ref{k40} imply that the minimal primitive cycle consists either of \parskip = 0 in

$-$ a single homoclinic trajectory,

$-$ a pair of conjugate heteroclinic trajectories.

\noindent
The last point means that $\gamma_{0} = ( e , J ( e ) )$ where $e \in \SE$ and $J ( x , \xi ) ( t ) = ( x ( - t ) , - \xi ( - t ) )$ is the usual symmetry. This is proved in Section \ref{s44}. Thus, to produce other geometric situations, one should relax the hypothesis that $P$ is a Schr\"{o}dinger operator and consider more general operators as explained in Remark \ref{c13}. \parskip = 0.05 in
\end{remark}

We end this part with a simple operator having a unique primitive cycle.

\begin{figure}
\begin{center}
\begin{picture}(0,0)%
\includegraphics{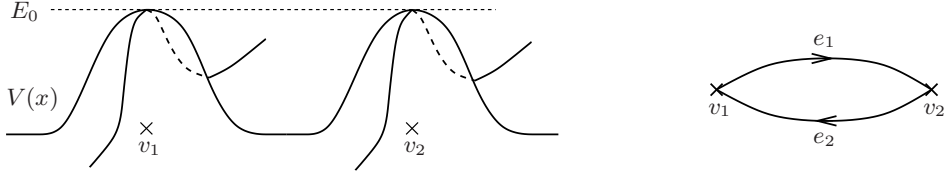}%
\end{picture}%
\setlength{\unitlength}{987sp}%
\begingroup\makeatletter\ifx\SetFigFont\undefined%
\gdef\SetFigFont#1#2#3#4#5{%
  \reset@font\fontsize{#1}{#2pt}%
  \fontfamily{#3}\fontseries{#4}\fontshape{#5}%
  \selectfont}%
\fi\endgroup%
\begin{picture}(23488,4121)(-14368,-4905)
\put(-10724,-4411){\makebox(0,0)[b]{\smash{{\SetFigFont{9}{10.8}{\rmdefault}{\mddefault}{\updefault}$v_{1}$}}}}
\put(9001,-3436){\makebox(0,0)[b]{\smash{{\SetFigFont{9}{10.8}{\rmdefault}{\mddefault}{\updefault}$v_{2}$}}}}
\put(3601,-3436){\makebox(0,0)[b]{\smash{{\SetFigFont{9}{10.8}{\rmdefault}{\mddefault}{\updefault}$v_{1}$}}}}
\put(6226,-1711){\makebox(0,0)[b]{\smash{{\SetFigFont{9}{10.8}{\rmdefault}{\mddefault}{\updefault}$e_{1}$}}}}
\put(6226,-4186){\makebox(0,0)[b]{\smash{{\SetFigFont{9}{10.8}{\rmdefault}{\mddefault}{\updefault}$e_{2}$}}}}
\put(-13499,-1036){\makebox(0,0)[rb]{\smash{{\SetFigFont{9}{10.8}{\rmdefault}{\mddefault}{\updefault}$E_{0}$}}}}
\put(-14324,-3211){\makebox(0,0)[lb]{\smash{{\SetFigFont{9}{10.8}{\rmdefault}{\mddefault}{\updefault}$V ( x )$}}}}
\put(-4049,-4411){\makebox(0,0)[b]{\smash{{\SetFigFont{9}{10.8}{\rmdefault}{\mddefault}{\updefault}$v_{2}$}}}}
\end{picture}%
\end{center}
\caption{The geometric setting of Example \ref{k42}.} \label{f51}
\end{figure}

\begin{example}\rm \label{k42}
Let $V_{1} , V_{2} \in C^{\infty}_{0} ( \R^{n} )$ be two radial barriers of height $E_{0}$. It means that, for $k = 1, 2$, the function $V_{k}$ is radial, satisfies \ref{h2} and $x \cdot \nabla V_{k} ( x ) < 0$ for $x$ in the interior of $\supp V_{k} \setminus \{ 0 \}$. For instance, one can take the potentials constructed in Section \ref{s21}. We then consider
\begin{equation}
P = - h^{2} \Delta + V_{1} ( x - v_{1} ) + V_{2} ( x - v_{2} ) ,
\end{equation}
with $v_{1} , v_{2} \in \R^{n}$ and $\vert v_{1} - v_{2} \vert$ large enough. See Figure \ref{f51} for a description of this potential. This operator satisfies \ref{h1} and \ref{h15}. Moreover, since the $V_{k}$'s are radial, we have
\begin{equation*}
\lambda_{1}^{v_{k}} = \cdots = \lambda_{n}^{v_{k}} = : \lambda_{k} .
\end{equation*}
On the other hand, the trapped set at energy $E_{0}$ is given by
\begin{equation} \label{k43}
K ( E_{0} ) = \big\{ ( v_{1} , 0 ) , ( v_{2} , 0 ) \big\} \cup e_{1} \cup e_{2} ,
\end{equation}
where $e_{1}$ (resp. $e_{2}$) is a heteroclinic trajectory from $v_{1}$ to $v_{2}$ (resp. from $v_{2}$ to $v_{1}$) whose base space projection is the segment $] v_{1} , v_{2} [$ in $\R^{n}$. Besides, one can verify that \ref{h16}, \ref{h18} and \ref{h19} hold true.

From \eqref{k43}, there exists a unique primitive cycle $\gamma_{0} = ( e_{1} , e_{2} )$ which is automatically minimal. Since \ref{h20} is also satisfied, we are in position to apply Corollary \ref{k40}. The damping is given by
\begin{equation*}
D_{0} = D ( \gamma_{0} ) = \frac{\alpha_{v_{1}} + \alpha_{v_{2}}}{\beta_{v_{1}} + \beta_{v_{2}}} = ( n - 1 ) \frac{\lambda_{1} \lambda_{2}}{\lambda_{1} + \lambda_{2}} ,
\end{equation*}
and
\begin{equation*}
\beta ( \gamma_{0} ) = \beta_{v_{1}} + \beta_{v_{2}} = \frac{\lambda_{1} + \lambda_{2}}{\lambda_{1} \lambda_{2}} .
\end{equation*}
Moreover, using that the actions coincide (i.e. $A_{e_{1}} = A_{e_{2}} = : A$) and that the Maslov's indices vanish (i.e. $\nu_{e_{1}} = \nu_{e_{2}} = 0$), the function $\mu$ defined in \eqref{k44} can be written
\begin{equation} \label{k47}
\mu ( \tau , h ) = - i e^{i 2 A / h} \Gamma ( S_{1} ) \Gamma ( S_{2} ) \frac{\sqrt{\lambda_{1} \lambda_{2}}}{2 \pi} \frac{e^{- i \frac{\pi}{2} ( S_{1} + S_{2} )}}{\lambda_{1}^{S_{1}} \lambda_{2}^{S_{2}}} \frac{\CM_{e_{1}}^{+} \CM_{e_{2}}^{+}}{\CM_{e_{1}}^{-} \CM_{e_{2}}^{-}} \frac{\vert g_{-}^{e_{1}} \vert^{1 - S_{2}} \vert g_{-}^{e_{2}} \vert^{1 - S_{1}}}{\vert g_{+}^{e_{1}} \vert^{S_{2}} \vert g_{+}^{e_{2}} \vert^{S_{1}}} ,
\end{equation}
with the notation
\begin{equation*}
S_{k} ( \tau ) : = \frac{n}{2} - \frac{D_{0}}{\lambda_{k}} - i \frac{\tau}{\lambda_{k}} .
\end{equation*}
So, we obtain an explicit asymptotic of the resonances by inserting the previous identities in \eqref{k41}. Using the asymptotic of the Gamma function stated in \eqref{k45}, the (unique) accumulation curve \eqref{k46} satisfies
\begin{equation} \label{n10}
\im \sigma = - i D_{0} + \frac{1}{\vert \ln h \vert}
\left\{ \begin{aligned}
&C_{1 2} + o ( 1 ) &&\text{ as } \re \sigma \to - \infty , \\
&- \pi \re \sigma + C_{1 2} + o ( 1 ) &&\text{ as } \re \sigma \to + \infty , \\
\end{aligned} \right.
\end{equation}
with the constant
\begin{equation*}
C_{1 2} = \frac{\lambda_{1} \lambda_{2}}{\lambda_{1} + \lambda_{2}} \ln \bigg( \frac{\lambda_{1} \lambda_{2}}{\lambda_{1}^{2 S^{0}_{1}} \lambda_{2}^{2 S^{0}_{2}}} \frac{\CM_{e_{1}}^{+} \CM_{e_{2}}^{+}}{\CM_{e_{1}}^{-} \CM_{e_{2}}^{-}} \frac{\vert g_{-}^{e_{1}} \vert^{1 - S^{0}_{2}} \vert g_{-}^{e_{2}} \vert^{1 - S^{0}_{1}}}{\vert g_{+}^{e_{1}} \vert^{S^{0}_{2}} \vert g_{+}^{e_{2}} \vert^{S^{0}_{1}}} \bigg) ,
\end{equation*}
and $S_{k}^{0} = S_{k} ( 0 )$. This equation can be compared with \eqref{e3}. As in Section \ref{s14}, the nature of the trapped set for energies below and above $E_{0}$ justifies the change of behavior between $- \infty$ and $+ \infty$. The accumulation curve is similar to that illustrated in Figure \ref{f24}.

Lastly, we consider the situation $\tau = 0$ in dimension $n = 1$. In this case, we have $D_{0} = 0$, $S_{1} = S_{2} = 1 / 2$, $\CM_{e_{k}}^{-} = \sqrt{\lambda_{3 - k} \vert g_{-}^{e_{k}} \vert}$ and $\CM_{e_{k}}^{+} = \sqrt{\lambda_{k} \vert g_{+}^{e_{k}} \vert}$ (see \eqref{e38} for the proof). Hence, \eqref{k47} becomes
\begin{equation*}
\mu ( 0 , h ) = - 2^{- 1} e^{i 2 A / h} ,
\end{equation*}
since $\Gamma ( 1 / 2 ) = \sqrt{\pi}$. Then, Corollary \ref{k40} implies that the resonances in $B ( E_{0} , C h \vert \ln h \vert^{- 1} )$ satisfy
\begin{equation} \label{k53}
z = E_{0} - \frac{2 A - ( 2 q + 1 ) \pi h + i h \ln 2}{\beta ( \gamma_{0}) \vert \ln h \vert} + o \Big( \frac{h}{\vert \ln h \vert } \Big) ,
\end{equation}
for some $q \in \Z$. Note that this asymptotic was previously obtained by the second and third authors (see \cite[Th\'eor\`eme 0.7]{FuRa98_01}). The remainder term in that paper is better than the present one.
\end{example}

\subsubsection{Identical barriers} \label{s74}

We work under the assumptions of Theorem \ref{j79} and suppose that the $\lambda_{j}^{v}$'s are the same for all the critical points $v \in \SV$. As consequence, the quantities $\alpha_{v}$, $\beta_{v}$ and $S_{v}$ are also independent of $v \in \SV$. In the sequel, we remove the subscript $v$. The simplest way to construct operators satisfying this property is to consider identical bumps localized at different places of $\R^{n}$.

Let $\gamma = ( e_{1} , \ldots , e_{K} )$ be a minimal cycle. Then,
\begin{equation} \label{k55}
D_{0} = D ( \gamma ) = \frac{\alpha ( \gamma )}{\beta ( \gamma )} = \frac{K \alpha}{K \beta} = \frac{1}{2} \sum_{j = 2}^{n} \lambda_{j} .
\end{equation}
The matrix $\SQ$ can be written
\begin{equation} \label{n9}
\SQ ( z , h ) = h^{S ( z , h ) / \lambda_{1} - 1 / 2} \CQ ( z , h ) ,
\end{equation}
where $\CQ$ is the matrix of the $\CQ_{e , \widetilde{e}}$. Mimicking \eqref{e18}, we define
\begin{equation} \label{k56}
\widehat{\CQ} ( \tau , h ) : = \CQ ( E_{0} + \tau h - i D_{0} h , h ) ,
\end{equation}
for $\tau \in \R$, and $\mu_{1} ( \tau ,h ), \ldots , \mu_{\card \SE} ( \tau , h )$ its eigenvalues. As below \eqref{e18}, these functions of $\tau , h$ are continuous and analytic outside of their crossings.

\begin{corollary}\sl \label{k48}
In the present setting, let $C , \delta > 0$. The resonances $z$ lying in \eqref{j80} satisfy, in the sense of Proposition \ref{j81},
\begin{equation*}
z = z_{q , k} ( \tau ) + o ( h \vert \ln h \vert^{- 1} ) ,
\end{equation*}
for some $q \in \Z$ and $k \in \{ 1 , \ldots , \card \SE\}$ with
\begin{equation} \label{k50}
z_{q , k} ( \tau ) = E_{0} + 2 q \pi \lambda_{1} \frac{h}{\vert \ln h \vert} - i h \sum_{j = 2}^{n} \frac{\lambda_{j}}{2} + i \ln ( \mu_{k} ( \tau , h ) ) \lambda_{1} \frac{h}{\vert \ln h \vert} .
\end{equation}
\end{corollary}

The asymptotic distribution of resonances is the same as in Proposition \ref{d9}. This is natural since the quantization rules Definition \ref{d1} and Definition \ref{j82} have the same structure thanks to \eqref{n9}. Furthermore, Corollary \ref{k48} is a natural generalization of Theorem \ref{d8}. Thus, the examples of Section \ref{s18} can be seen as examples of the present situation.

On the other hand, it is possible to directly adapt the proof of Theorem \ref{d8} to the present setting. As consequence, \eqref{k50} holds true in the sense of Proposition \ref{d9} (i.e. we can relax the restriction $\delta = C \vert \ln h \vert^{- 1}$ of Proposition \ref{j81}) and we can give a lower bound for the multiplicity of the resonances by the one of the eigenvalues of $\CQ$ (as in Proposition \ref{d91}).

\subsubsection{Disconnected graphs} \label{s88}

\begin{figure}
\begin{center}
\begin{picture}(0,0)%
\includegraphics{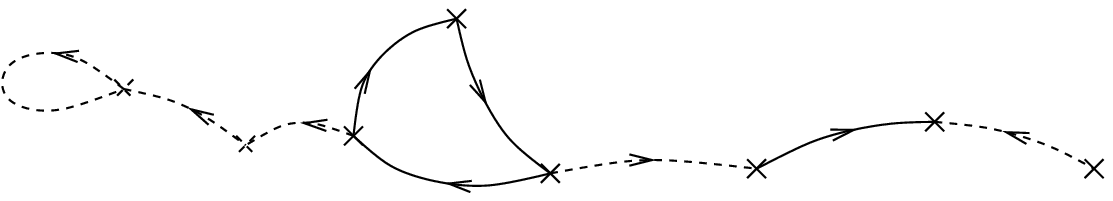}%
\end{picture}%
\setlength{\unitlength}{1184sp}%
\begingroup\makeatletter\ifx\SetFigFont\undefined%
\gdef\SetFigFont#1#2#3#4#5{%
  \reset@font\fontsize{#1}{#2pt}%
  \fontfamily{#3}\fontseries{#4}\fontshape{#5}%
  \selectfont}%
\fi\endgroup%
\begin{picture}(17686,3772)(-6552,-4850)
\put(751,-4786){\makebox(0,0)[b]{\smash{{\SetFigFont{9}{10.8}{\rmdefault}{\mddefault}{\updefault}$( \SV_{1} , \SE_{1} )$}}}}
\put(7201,-4786){\makebox(0,0)[b]{\smash{{\SetFigFont{9}{10.8}{\rmdefault}{\mddefault}{\updefault}$( \SV_{2} , \SE_{2} )$}}}}
\put(10951,-4786){\makebox(0,0)[b]{\smash{{\SetFigFont{9}{10.8}{\rmdefault}{\mddefault}{\updefault}$( \SV_{3} , \SE_{3} )$}}}}
\end{picture}%
\end{center}
\caption{A decomposition of the graph $( \SV , \SE )$ as explained in Section \ref{s88}. The dashed objects belong to $( \SV , \SE )$ but not to the $( \SV_{k} , \SE_{k} )$'s.} \label{f52}
\end{figure}

In this section, we assume that the principal part of the graph $( \SV , \SE )$ can be decomposed into disjoint parts. More precisely, let $\SV_{1} , \ldots , \SV_{K}$ be a finite number of pairwise disjoint subsets of $\SV$. For $k \in  \{ 1 , \ldots , K \}$, we define
\begin{equation} \label{k60}
\SE_{k} = \{ e \in \SE ; \ e^{-} \in \SV_{k} \text{ and } e^{+} \in \SV_{k} \} ,
\end{equation}
the set of edges of the original graph $( \SV , \SE )$ between points of $\SV_{k}$. These sets are necessarily pairwise disjoint. This setting is illustrated in Figure \ref{f52}. We assume that any edge of any minimal cycle of $( \SV , \SE )$ lies in one of the $\SE_{k}$'s. It is equivalent to assume that any minimal cycle of $( \SV , \SE )$ lies in one of the graphs $( \SV_{k} , \SE_{k})$.

Let $P_{k}$ be operators satisfying the assumptions of Theorem \ref{j79} and whose corresponding graph is precisely $( \SV_{k} , \SE_{k} )$. Moreover, we suppose that the symbol of $P_{k}$ coincides with that of $P$ in a neighborhood of $K_{k} ( E_{0} )$, the trapped set of $P_{k}$ at energy $E_{0}$. In this situation, Theorem \ref{j79} leads to the following result.

\begin{corollary}\sl \label{k54}
In the present setting, let $C , \delta > 0$. We have
\begin{equation*}
\dist \bigg( \res ( P ) , \bigcup_{k = 1}^{K} \res ( P_{k} ) \bigg) = o \Big( \frac {h}{\vert \ln h \vert} \Big) ,
\end{equation*}
as $h$ tends to $0$, in the domain \eqref{j80}.
\end{corollary}

Since any cycle in $( \SV_{k} , \SE_{k} )$ is a cycle in $( \SV , \SE )$, the minimal damping of $P_{k}$, denoted $D_{0}^{( k )}$, satisfies
\begin{equation*}
D_{0}^{( k )} \geq D_{0} .
\end{equation*}
Moreover, there exists at least one $k$ such that $D_{0}^{( k )} = D_{0}$. Nevertheless, it is possible that $D_{0}^{( k )} > D_{0}$ for some other $k$. It happens when $( \SV_{k} , \SE_{k} )$ contains no minimal cycle of $( \SV , \SE )$. In this case, $P_{k}$ has no resonance in \eqref{j80}, the set $\SV_{k}$ plays no role in Corollary \ref{k54} and can be removed from the list of the $\SV_{\bullet}$'s.

In this result, the remainder term $o ( h \vert \ln h \vert^{- 1} )$ is coherent with the asymptotic of resonances provided by Theorem \ref{j79}. In general, it is probably not possible to replace it by $\CO ( h^{\infty} )$. Indeed, if the graph is not totally disconnected (that is if there exist edges joining the sub-graphs $( \SV_{k} , \SE_{k} )$), one may imagine that the interaction between the sub-graphs generates some lower order contributions.

It is not always possible to decompose the graph $( \SV , \SE )$. This is for instance the case of all the operators considered in Section \ref{s6}. That the matrix $\CQ$ of \eqref{d4} mixes the contributions of the different edges explains this phenomenon. Note also that we give no lower bound on the multiplicity of the resonance of $P$ in terms of the multiplicity of the resonances of the $P_{k}$'s.

\begin{figure}
\begin{center}
\begin{picture}(0,0)%
\includegraphics{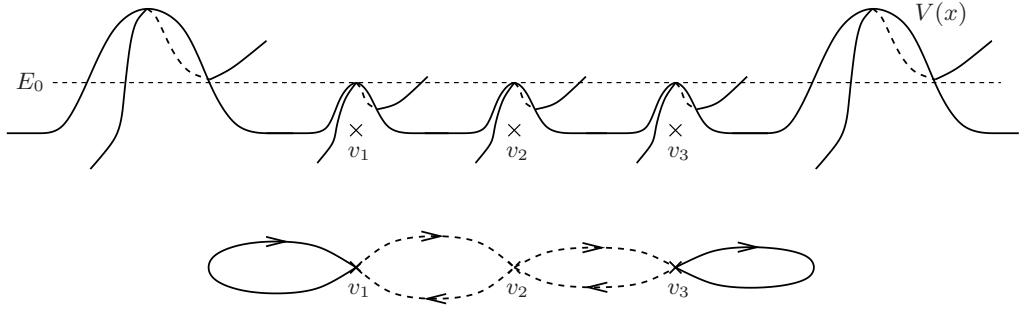}%
\end{picture}%
\setlength{\unitlength}{987sp}%
\begingroup\makeatletter\ifx\SetFigFont\undefined%
\gdef\SetFigFont#1#2#3#4#5{%
  \reset@font\fontsize{#1}{#2pt}%
  \fontfamily{#3}\fontseries{#4}\fontshape{#5}%
  \selectfont}%
\fi\endgroup%
\begin{picture}(25513,7501)(-25918,-8286)
\put(-25649,-2836){\makebox(0,0)[lb]{\smash{{\SetFigFont{9}{10.8}{\rmdefault}{\mddefault}{\updefault}$E_{0}$}}}}
\put(-17024,-4561){\makebox(0,0)[b]{\smash{{\SetFigFont{9}{10.8}{\rmdefault}{\mddefault}{\updefault}$v_{1}$}}}}
\put(-17024,-8011){\makebox(0,0)[b]{\smash{{\SetFigFont{9}{10.8}{\rmdefault}{\mddefault}{\updefault}$v_{1}$}}}}
\put(-8999,-8011){\makebox(0,0)[b]{\smash{{\SetFigFont{9}{10.8}{\rmdefault}{\mddefault}{\updefault}$v_{3}$}}}}
\put(-13049,-8011){\makebox(0,0)[b]{\smash{{\SetFigFont{9}{10.8}{\rmdefault}{\mddefault}{\updefault}$v_{2}$}}}}
\put(-8999,-4561){\makebox(0,0)[b]{\smash{{\SetFigFont{9}{10.8}{\rmdefault}{\mddefault}{\updefault}$v_{3}$}}}}
\put(-13049,-4561){\makebox(0,0)[b]{\smash{{\SetFigFont{9}{10.8}{\rmdefault}{\mddefault}{\updefault}$v_{2}$}}}}
\put(-3074,-1111){\makebox(0,0)[lb]{\smash{{\SetFigFont{9}{10.8}{\rmdefault}{\mddefault}{\updefault}$V ( x )$}}}}
\end{picture}%
\end{center}
\caption{The potential and the graph decomposition in Example \ref{k61}.} \label{f53}
\end{figure}

\begin{example}\rm \label{k61}
In dimension $n \geq 2$, we consider
\begin{equation*}
P = -h^{2} \Delta + V ( x ) ,
\end{equation*}
where the potential $V$ is the sum of five aligned radial bumps. The three central ones have height $E_{0}$ and the two extreme ones have height larger than $E_{0}$. A description of the geometry as well as some notations can be found in Figure \ref{f53}. This operator verifies the hypotheses of Theorem \ref{j79}. We assume that the $\lambda_{1}^{v_{k}} = \cdots = \lambda_{n}^{v_{k}} = : \lambda_{k}$ satisfy
\begin{equation} \label{k59}
\lambda_{1} = \lambda_{3} < \lambda_{2} .
\end{equation}
The primitive cycles of the graph $( \SV , \SE )$ are the two homoclinic trajectories with end point $v_{1}$ and $v_{3}$, and the two pairs of heteroclinic trajectories between $v_{1}$ and $v_{2}$ and between $v_{2}$ and $v_{3}$. The damping of the two homoclinic trajectories is $( n - 1 ) \lambda_{1} /2$, whereas the damping of the two pairs of heteroclinic trajectories is $( n - 1 ) \lambda_{1} \lambda_{2} / ( \lambda_{1} + \lambda_{2} )$. From \eqref{k59}, the minimal primitive cycles are then the two homoclinic trajectories.

Now, we choose $\SV_{-} = \{ v_{1} \}$ and $\SV_{+} = \{ v_{3} \}$. Let $\SE_{\pm}$ denote the set of edges associated to $\SV_{\pm}$ defined in \eqref{k60}. Then, $\SE_{-}$ (resp. $\SE_{+}$) consists of the homoclinic trajectory with end point $v_{1}$ (resp. $v_{3}$). We also define
\begin{equation*}
P_{\pm} = - h^{2} \Delta + V_{\pm} ( x ) ,
\end{equation*}
where $V_{-}$ (resp. $V_{+}$) is the sum of the two bumps on the left (resp. on the right). They verify the hypotheses of Theorem \ref{j79} and their associated graph is precisely $( \SV_{\pm} , \SE_{\pm} )$.

Summing up, we are in position to apply Corollary \ref{k54}. Thus,
\begin{equation*}
\dist \big( \res ( P ) , \res ( P_{-} ) \cup \res ( P_{+} ) \big) = o \Big( \frac {h}{\vert \ln h \vert} \Big) ,
\end{equation*}
in the domain \eqref{j80}. On the other hand, the operators $P_{\pm}$ satisfy the assumptions of Theorem \ref{d8} and the asymptotic of their resonances has been computed in \eqref{e8}--\eqref{e7}. Combining these results, we obtain an explicit asymptotic of the resonances of $P$.
\end{example}

\subsubsection{The one dimensional situation} \label{s85}

In dimension $n = 1$, we have $\alpha_{v} = 0$ for all $v \in \SV$. Thus, the damping of any cycle is null. Under \ref{h20}, it implies that
\begin{equation*}
D_{0} = 0 ,
\end{equation*}
and that any cycle is minimal. In particular, $\SQ = \SQ^{\rm p r i n}$ and $F = F^{\rm p r i n}$. Nevertheless, the $\beta_{v} = ( \lambda_{1}^{v} )^{- 1}$'s can be different and $\Z$-independent. Note also that \eqref{j83}, which defines $\CQ$, can somehow be simplified as in \eqref{k47}--\eqref{k53}.

In compact situations, Colin de Verdi\`ere and Parisse \cite[Section 8--12]{CoPa99_01} have obtained the semiclassical asymptotic of the eigenvalues under similar hypotheses. They have also computed in \cite[Section 13]{CoPa99_01} the scattering matrix in two geometric situations. When $E_{0}$ is the maximum of the potential, the second author \cite{Fu98_01} has given the semiclassical behavior of the scattering matrix in the present geometric setting and interpreted its leading term using Feynman integrals.

\begin{figure}
\begin{center}
\begin{picture}(0,0)%
\includegraphics{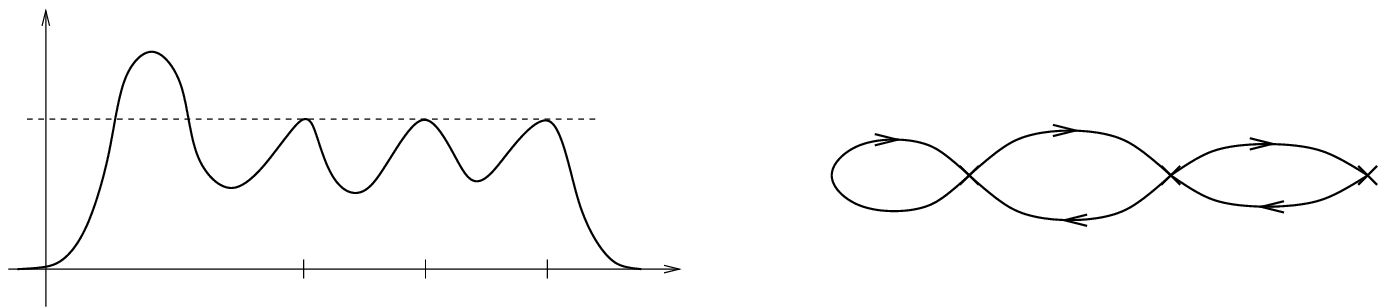}%
\end{picture}%
\setlength{\unitlength}{1184sp}%
\begingroup\makeatletter\ifx\SetFigFont\undefined%
\gdef\SetFigFont#1#2#3#4#5{%
  \reset@font\fontsize{#1}{#2pt}%
  \fontfamily{#3}\fontseries{#4}\fontshape{#5}%
  \selectfont}%
\fi\endgroup%
\begin{picture}(22323,4844)(8011,-5783)
\put(17026,-5611){\makebox(0,0)[b]{\smash{{\SetFigFont{9}{10.8}{\rmdefault}{\mddefault}{\updefault}$v_{3}$}}}}
\put(22501,-2836){\makebox(0,0)[b]{\smash{{\SetFigFont{9}{10.8}{\rmdefault}{\mddefault}{\updefault}$e_{1}$}}}}
\put(25351,-2686){\makebox(0,0)[b]{\smash{{\SetFigFont{9}{10.8}{\rmdefault}{\mddefault}{\updefault}$e_{2}$}}}}
\put(25426,-4861){\makebox(0,0)[b]{\smash{{\SetFigFont{9}{10.8}{\rmdefault}{\mddefault}{\updefault}$e_{3}$}}}}
\put(28501,-2911){\makebox(0,0)[b]{\smash{{\SetFigFont{9}{10.8}{\rmdefault}{\mddefault}{\updefault}$e_{4}$}}}}
\put(28651,-4636){\makebox(0,0)[b]{\smash{{\SetFigFont{9}{10.8}{\rmdefault}{\mddefault}{\updefault}$e_{5}$}}}}
\put(23776,-4186){\makebox(0,0)[b]{\smash{{\SetFigFont{9}{10.8}{\rmdefault}{\mddefault}{\updefault}$v_{1}$}}}}
\put(27001,-4186){\makebox(0,0)[b]{\smash{{\SetFigFont{9}{10.8}{\rmdefault}{\mddefault}{\updefault}$v_{2}$}}}}
\put(30151,-4186){\makebox(0,0)[b]{\smash{{\SetFigFont{9}{10.8}{\rmdefault}{\mddefault}{\updefault}$v_{3}$}}}}
\put(8026,-2911){\makebox(0,0)[lb]{\smash{{\SetFigFont{9}{10.8}{\rmdefault}{\mddefault}{\updefault}$E_{0}$}}}}
\put(11101,-1861){\makebox(0,0)[lb]{\smash{{\SetFigFont{9}{10.8}{\rmdefault}{\mddefault}{\updefault}$V (x)$}}}}
\put(18976,-4936){\makebox(0,0)[lb]{\smash{{\SetFigFont{9}{10.8}{\rmdefault}{\mddefault}{\updefault}$x$}}}}
\put(13126,-5611){\makebox(0,0)[b]{\smash{{\SetFigFont{9}{10.8}{\rmdefault}{\mddefault}{\updefault}$v_{1}$}}}}
\put(15076,-5611){\makebox(0,0)[b]{\smash{{\SetFigFont{9}{10.8}{\rmdefault}{\mddefault}{\updefault}$v_{2}$}}}}
\end{picture} \newline\setlength{\unitlength}{987sp}
\begin{picture}(28000,5200)(1500,750)
\includegraphics[width=113pt]{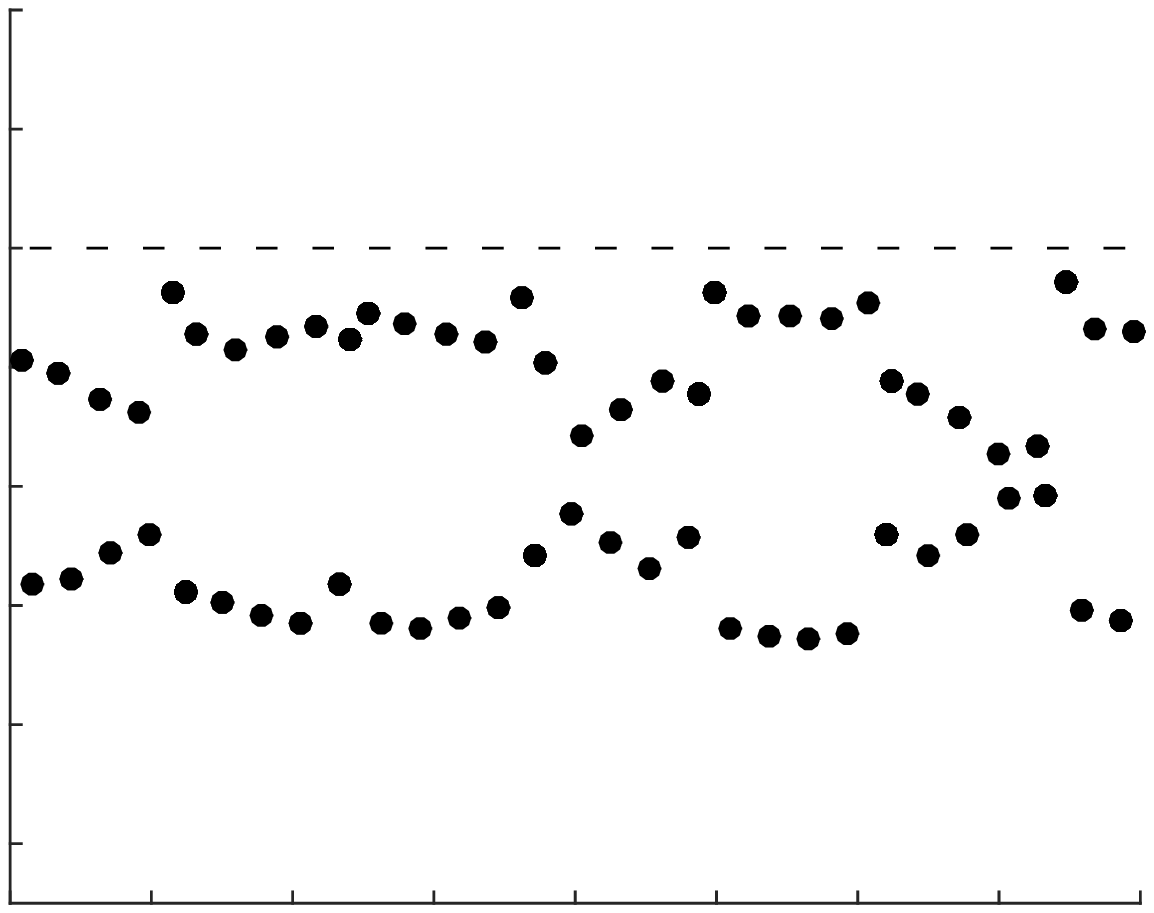}
\includegraphics[width=113pt]{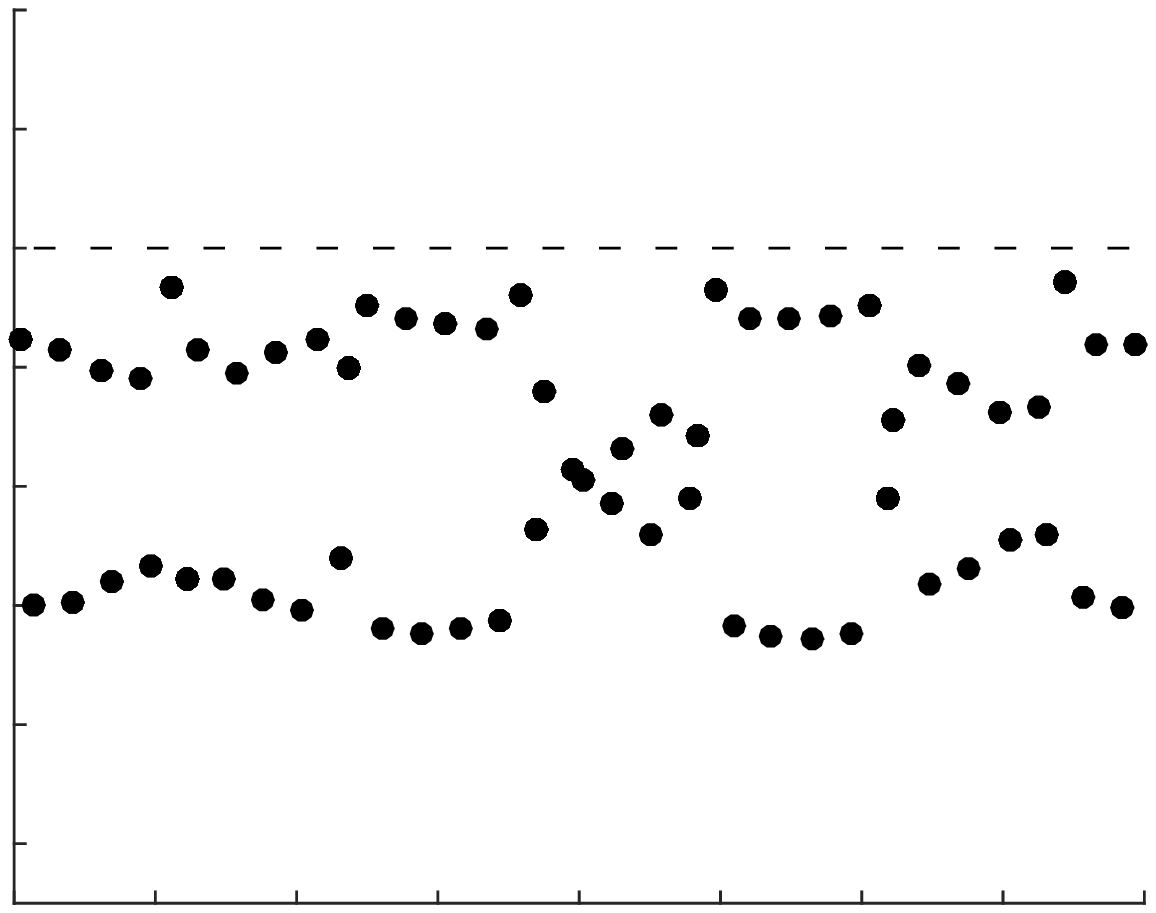}
\includegraphics[width=113pt]{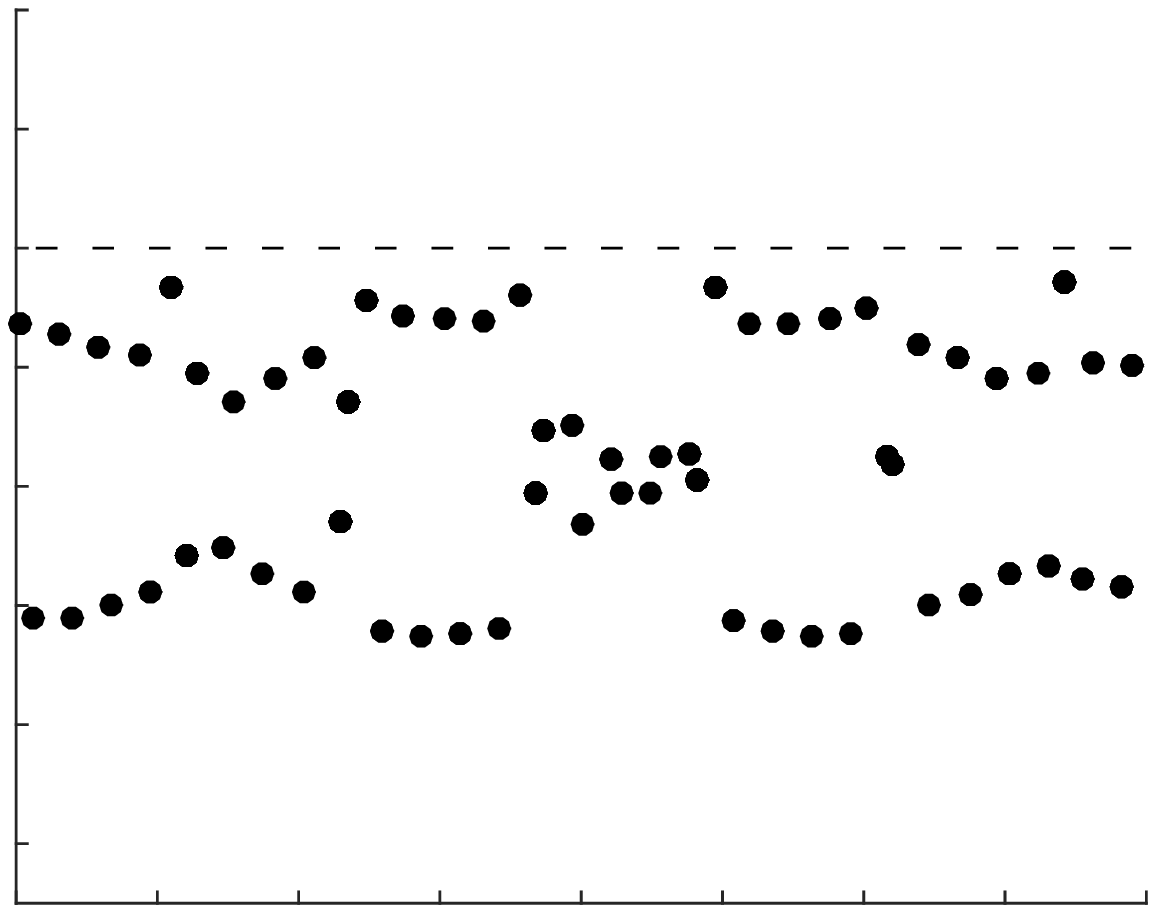}
\includegraphics[width=113pt]{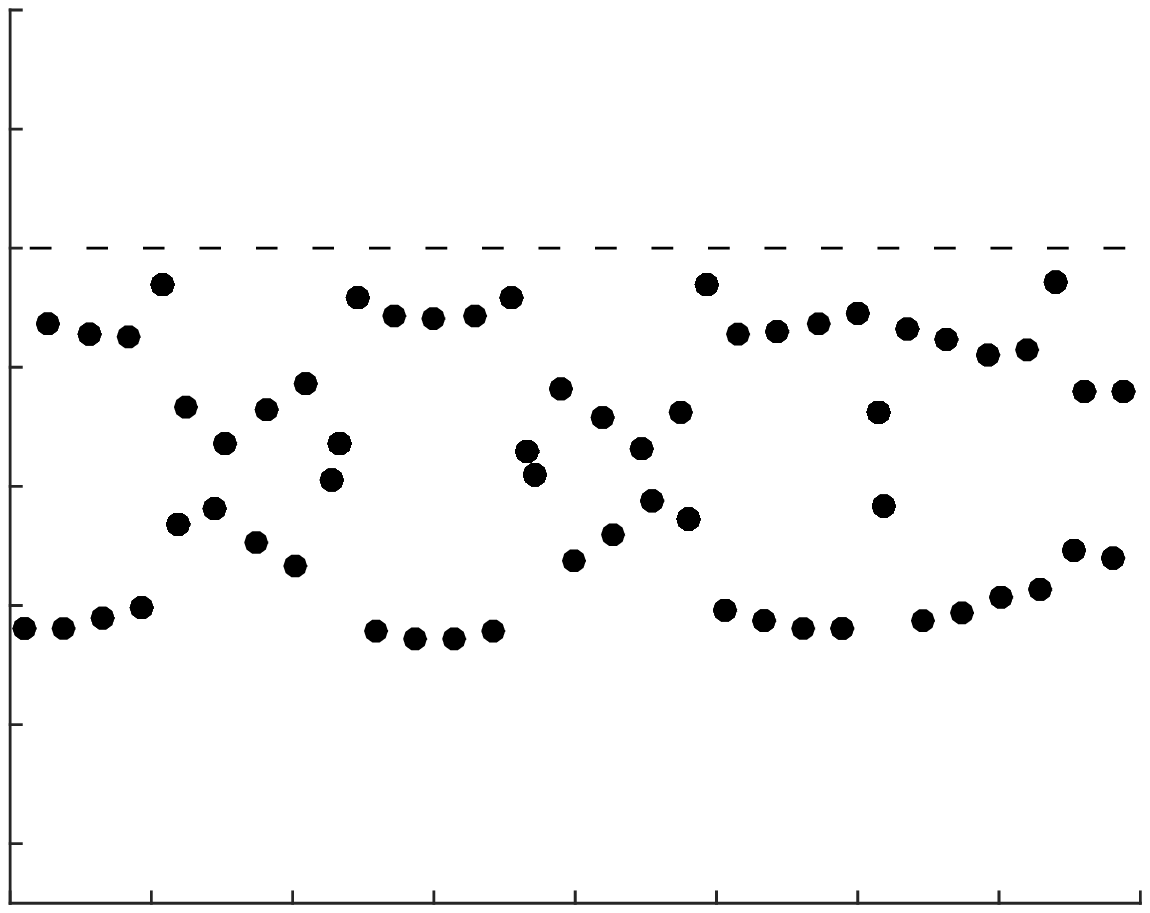}
\end{picture}
\end{center}
\caption{The potential of Example \ref{k52}, the corresponding graph and the resonances in a vicinity of size $h \vert \ln h \vert^{- 1}$ of $E_{0}$ for different values of $h$.} \label{f48}
\end{figure}

\begin{example}\rm \label{k52}
Let $V$ be the potential given by Figure \ref{f48}. Then,
\begin{equation*}
P = - h^{2} \Delta + V ( x ) ,
\end{equation*}
verifies the assumptions of Theorem \ref{j79}. Its trapped set at energy $E_{0}$ consists of three hyperbolic fixed points, one homoclinic trajectory and four heteroclinic trajectories. The associated (minimal) primitive cycles are $( e_{1} )$, $( e_{2} , e_{3} )$ and $( e_{4} , e_{5} )$. Moreover, the matrix $\SQ$ writes
\begin{equation*}
\SQ = \left( \begin{array}{ccccc}
\SQ_{e_{1} , e_{1}} & 0 & \SQ_{e_{1} , e_{3}} & 0 & 0  \\
\SQ_{e_{2} , e_{1}} & 0 & \SQ_{e_{2} , e_{3}} & 0 & 0  \\
0 & \SQ_{e_{3} , e_{2}} & 0 & 0 & \SQ_{e_{3} , e_{5}}  \\
0 & \SQ_{e_{4} , e_{2}} & 0 & 0 & \SQ_{e_{4} , e_{5}}  \\
0 & 0 & 0 & \SQ_{e_{5} , e_{4}} & 0 
\end{array} \right) .
\end{equation*}
Taking the determinant of $1 - \SQ$ and using a numerical computation of the zeros of the exponential sum $F ( z , h )$, we obtain the leading term in the asymptotic of the resonances. The result is drawn in Figure \ref{f48} for some $\Z$-independent $\beta_{v_{1}} , \beta_{v_{2}} , \beta_{v_{3}}$.
\end{example}

\begin{figure}
\begin{center}
\begin{picture}(0,0)%
\includegraphics{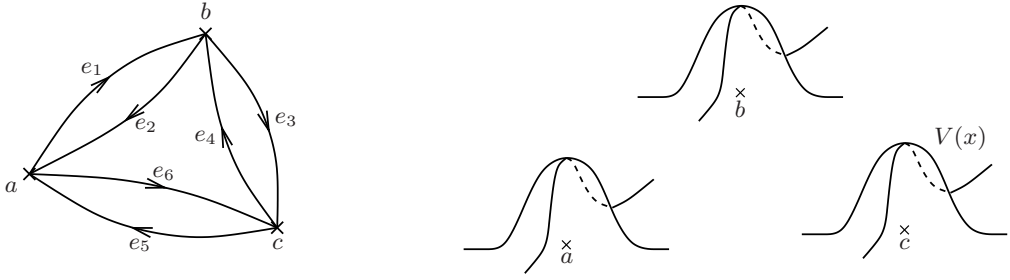}%
\end{picture}%
\setlength{\unitlength}{987sp}%
\begingroup\makeatletter\ifx\SetFigFont\undefined%
\gdef\SetFigFont#1#2#3#4#5{%
  \reset@font\fontsize{#1}{#2pt}%
  \fontfamily{#3}\fontseries{#4}\fontshape{#5}%
  \selectfont}%
\fi\endgroup%
\begin{picture}(25075,6811)(-26639,-6875)
\put(-19949,-6211){\makebox(0,0)[b]{\smash{{\SetFigFont{9}{10.8}{\rmdefault}{\mddefault}{\updefault}$c$}}}}
\put(-4199,-6211){\makebox(0,0)[b]{\smash{{\SetFigFont{9}{10.8}{\rmdefault}{\mddefault}{\updefault}$c$}}}}
\put(-8324,-2911){\makebox(0,0)[b]{\smash{{\SetFigFont{9}{10.8}{\rmdefault}{\mddefault}{\updefault}$b$}}}}
\put(-12674,-6586){\makebox(0,0)[b]{\smash{{\SetFigFont{9}{10.8}{\rmdefault}{\mddefault}{\updefault}$a$}}}}
\put(-3449,-3586){\makebox(0,0)[lb]{\smash{{\SetFigFont{9}{10.8}{\rmdefault}{\mddefault}{\updefault}$V ( x )$}}}}
\put(-24599,-1786){\makebox(0,0)[b]{\smash{{\SetFigFont{9}{10.8}{\rmdefault}{\mddefault}{\updefault}$e_{1}$}}}}
\put(-23249,-3211){\makebox(0,0)[b]{\smash{{\SetFigFont{9}{10.8}{\rmdefault}{\mddefault}{\updefault}$e_{2}$}}}}
\put(-21749,-3511){\makebox(0,0)[b]{\smash{{\SetFigFont{9}{10.8}{\rmdefault}{\mddefault}{\updefault}$e_{4}$}}}}
\put(-19724,-3061){\makebox(0,0)[b]{\smash{{\SetFigFont{9}{10.8}{\rmdefault}{\mddefault}{\updefault}$e_{3}$}}}}
\put(-22799,-4411){\makebox(0,0)[b]{\smash{{\SetFigFont{9}{10.8}{\rmdefault}{\mddefault}{\updefault}$e_{6}$}}}}
\put(-23399,-6211){\makebox(0,0)[b]{\smash{{\SetFigFont{9}{10.8}{\rmdefault}{\mddefault}{\updefault}$e_{5}$}}}}
\put(-26624,-4786){\makebox(0,0)[b]{\smash{{\SetFigFont{9}{10.8}{\rmdefault}{\mddefault}{\updefault}$a$}}}}
\put(-21749,-436){\makebox(0,0)[b]{\smash{{\SetFigFont{9}{10.8}{\rmdefault}{\mddefault}{\updefault}$b$}}}}
\end{picture}%
\end{center}
\caption{The graph in the three bump case and the potential of Example \ref{k74}.} \label{f54}
\end{figure}

\subsubsection{Three bumps} \label{s87}

In this part, we consider operators $P$ satisfying the assumptions of Theorem \ref{j79} and such that $\SV$ consists of three hyperbolic fixed points
\begin{equation*}
\SV = \{ a , b , c \} ,
\end{equation*}
and $\SE$ of six heteroclinic trajectories between these points
\begin{equation*}
\SE = \{ e_{1} , e_{2} , e_{3} , e_{4} , e_{5} , e_{6} \} ,
\end{equation*}
as illustrated in Figure \ref{f54}. In this case, there are five primitive cycles:
\begin{equation*}
( e_{1} , e_{2} ) , \ ( e_{3} , e_{4} ) , \ ( e_{5} , e_{6} ) , \ ( e_{1} , e_{3} , e_{5} ) , \ ( e_{2} , e_{6} , e_{4} ) .
\end{equation*}
Moreover, the matrix $\SQ$ takes the form
\begin{equation*}
\SQ = \left( \begin{array}{cccccc}
0 & \SQ_{e_{1} , e_{2}} & 0 & 0 & \SQ_{e_{1} , e_{5}} & 0  \\
\SQ_{e_{2} , e_{1}} & 0 & 0 & \SQ_{e_{2} , e_{4}} & 0 & 0  \\
\SQ_{e_{3} , e_{1}} & 0 & 0 & \SQ_{e_{3} , e_{4}} & 0 & 0  \\
0 & 0 & \SQ_{e_{4} , e_{3}} & 0 & 0 & \SQ_{e_{4} , e_{6}}  \\
0 & 0 & \SQ_{e_{5} , e_{3}} & 0 & 0 & \SQ_{e_{5} , e_{6}}  \\
0 & \SQ_{e_{6} , e_{2}} & 0 & 0 & \SQ_{e_{6} , e_{5}} & 0
\end{array} \right) .
\end{equation*}
We suppose for simplicity that the dimension is $n = 2$ and that the three vertices are isotropic, that is
\begin{equation*}
\lambda_{1}^{a} = \lambda_{2}^{a} = : \lambda_{a} , \quad \lambda_{1}^{b} = \lambda_{2}^{b} = : \lambda_{b} , 
\quad \lambda_{1}^{c} = \lambda_{2}^{c} = : \lambda_{c} .
\end{equation*}
We first construct operators satisfying these assumptions.

\begin{example}\rm \label{k74}
Consider three radial potentials $V_{1} , V_{2} , V_{3} \in C^{\infty}_{0} ( \R^{2} )$ as in Section \ref{s21} with large scattering angle. We then set
\begin{equation} \label{k75}
V ( x ) = V_{1} ( x - a ) + V_{2} ( x - b ) + V_{3} ( x - c ) ,
\end{equation}
for $a , b , c \in \R^{2}$ (see Figure \ref{f54}). If the points $a , b , c$ are sufficiently spaced out (to avoid that the potentials overlap) and all the angles of the triangle $a b c$ are acute (which ensures that any Hamiltonian curve touches the support of at most two potentials from Proposition \ref{a14}), one can show that the operator 
\begin{equation*}
P = - h ^{2} \Delta + V ( x ) ,
\end{equation*}
satisfies the assumptions of this section.

On the other hand, putting three bumps at the angles of the triangle $a b c$ does not guaranty that the dynamical hypotheses of this part hold. Indeed, consider the potential $W \in C^{\infty}_{0} ( \R^{2} )$ constructed in Example \ref{k76}. It is a radial compactly supported function with a non-degenerate maximum at $x = 0$ satisfying $W ( 0 ) = E_{0}$ and such that $x \cdot \nabla W (x) < 0$ for $x$ in the interior of $\supp W \setminus \{ 0 \}$. As before, we define
\begin{equation}
V ( x ) = W ( x - a ) + W ( x - b ) + W ( x - c ) ,
\end{equation}
where $a b c $ is an equilateral triangle with sides of length $L : = \vert b - a \vert = \vert c - b \vert = \vert a - c \vert  > 0$. This potential looks similar to \eqref{k75}, but the trapped set at energy $E_{0}$ of
\begin{equation*}
P = - h^{2} \Delta + V ( x ) ,
\end{equation*}
has a periodic bicharacteristic curve (touching successively the support of $W ( x - a )$, $W ( x - b )$ and $W ( x - c )$) for all $L$ large enough. This follows from \eqref{k80} and \eqref{k79}. Thus, this operator does not enter into the setting considered here. This shows that an assumption on the shape of the potential (being made of three bumps) does not imply a geometric property of the trapped set (being constituted of three fixed points and six heteroclinic trajectories).

We have seen that any acute triangle can be realized with three bumps. The situation is slightly different if one of the angle (say $\widehat{a b c }$ for instance) is obtuse. It can happen that there is no trajectory from $a$ to $c$ (if $a , b, c$ are almost aligned for instance), or that there are two or more trajectories from $a$ to $c$ (one whose base space projection is the segment $[ a ,  c ]$ and other trajectories passing through the support of $V_{2}$). Nevertheless, in the last case, it is possible to add an absorbing potential supported near these additional trajectories to ``remove'' their contribution in the trapped set (see Remark \ref{c13}). Of course, \ref{h18} forbids to have a right angle.
\end{example}

\begin{figure}
\begin{center}
\begin{picture}(0,0)%
\includegraphics{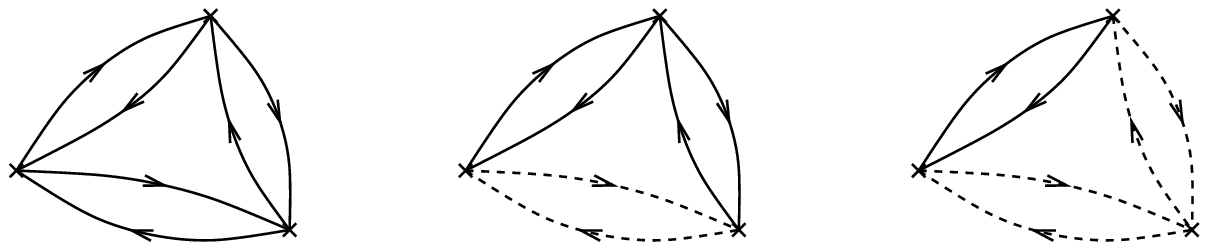}%
\end{picture}%
\setlength{\unitlength}{987sp}%
\begingroup\makeatletter\ifx\SetFigFont\undefined%
\gdef\SetFigFont#1#2#3#4#5{%
  \reset@font\fontsize{#1}{#2pt}%
  \fontfamily{#3}\fontseries{#4}\fontshape{#5}%
  \selectfont}%
\fi\endgroup%
\begin{picture}(23096,6073)(-34814,-6536)
\put(-23210,-6481){\makebox(0,0)[b]{\smash{{\SetFigFont{9}{10.8}{\rmdefault}{\mddefault}{\updefault}(B)}}}}
\put(-20624,-5461){\makebox(0,0)[b]{\smash{{\SetFigFont{9}{10.8}{\rmdefault}{\mddefault}{\updefault}$c$}}}}
\put(-22124,-586){\makebox(0,0)[b]{\smash{{\SetFigFont{9}{10.8}{\rmdefault}{\mddefault}{\updefault}$b$}}}}
\put(-26174,-4186){\makebox(0,0)[b]{\smash{{\SetFigFont{9}{10.8}{\rmdefault}{\mddefault}{\updefault}$a$}}}}
\put(-11924,-5461){\makebox(0,0)[b]{\smash{{\SetFigFont{9}{10.8}{\rmdefault}{\mddefault}{\updefault}$c$}}}}
\put(-13424,-586){\makebox(0,0)[b]{\smash{{\SetFigFont{9}{10.8}{\rmdefault}{\mddefault}{\updefault}$b$}}}}
\put(-17474,-4186){\makebox(0,0)[b]{\smash{{\SetFigFont{9}{10.8}{\rmdefault}{\mddefault}{\updefault}$a$}}}}
\put(-29249,-5461){\makebox(0,0)[b]{\smash{{\SetFigFont{9}{10.8}{\rmdefault}{\mddefault}{\updefault}$c$}}}}
\put(-30749,-586){\makebox(0,0)[b]{\smash{{\SetFigFont{9}{10.8}{\rmdefault}{\mddefault}{\updefault}$b$}}}}
\put(-34799,-4186){\makebox(0,0)[b]{\smash{{\SetFigFont{9}{10.8}{\rmdefault}{\mddefault}{\updefault}$a$}}}}
\put(-14512,-6481){\makebox(0,0)[b]{\smash{{\SetFigFont{9}{10.8}{\rmdefault}{\mddefault}{\updefault}(C)}}}}
\put(-31837,-6481){\makebox(0,0)[b]{\smash{{\SetFigFont{9}{10.8}{\rmdefault}{\mddefault}{\updefault}(A)}}}}
\end{picture}%
\end{center}
\caption{The principal graph in the three different settings. Dashed edges correspond to non-minimal cycles.} \label{f55}
\end{figure}

We now specify the asymptotic of the resonances using Theorem \ref{j79} and distinguishing between the three possible cases (modulo symmetry) for the principal graph. They are illustrated in Figure \ref{f55}.

In case {\rm (A)}, we assume that $\lambda_{a} = \lambda_{b} = \lambda_{c} = : \lambda$. Geometrically, this can achieved by taking identical bumps in Example \ref{k74}. We are then in the framework of Corollary \ref{k48}. All the cycles are minimal and $D_{0} = \lambda / 2$ from \eqref{k55}. Moreover, the resonances are close to
\begin{equation} \label{k65}
z_{q , k} ( \tau ) = E_{0} + 2 q \pi \lambda \frac{h}{\vert \ln h \vert} - i \frac{\lambda}{2} h  + i \ln ( \mu_{k} ( \tau , h ) ) \lambda \frac{h}{\vert \ln h \vert} ,
\end{equation}
where $\mu_{k} ( \tau , h )$ are the eigenvalues of the matrix $\widehat{\CQ} ( \tau , h )$ defined by
\begin{equation} \label{k63}
\widehat{\CQ} = \left( \begin{array}{cccccc}
0 & \widehat{\CQ}_{e_{1} , e_{2}} & 0 & 0 & \widehat{\CQ}_{e_{1} , e_{5}} & 0  \\
\widehat{\CQ}_{e_{2} , e_{1}} & 0 & 0 & \widehat{\CQ}_{e_{2} , e_{4}} & 0 & 0  \\
\widehat{\CQ}_{e_{3} , e_{1}} & 0 & 0 & \widehat{\CQ}_{e_{3} , e_{4}} & 0 & 0  \\
0 & 0 & \widehat{\CQ}_{e_{4} , e_{3}} & 0 & 0 & \widehat{\CQ}_{e_{4} , e_{6}}  \\
0 & 0 & \widehat{\CQ}_{e_{5} , e_{3}} & 0 & 0 & \widehat{\CQ}_{e_{5} , e_{6}}  \\
0 & \widehat{\CQ}_{e_{6} , e_{2}} & 0 & 0 & \widehat{\CQ}_{e_{6} , e_{5}} & 0
\end{array} \right) ,
\end{equation}
with $\widehat{\CQ}_{e , \widetilde{e}} ( \tau , h ) = \CQ_{e , \widetilde{e}} ( E_{0} + \tau h - i \lambda h / 2 , h )$. Unfortunately, an explicit formula for the eigenvalues of a generic matrix of this form seems difficult to obtain.

\begin{figure}
\begin{center}
\begin{picture}(370,170)
\includegraphics[width=350pt]{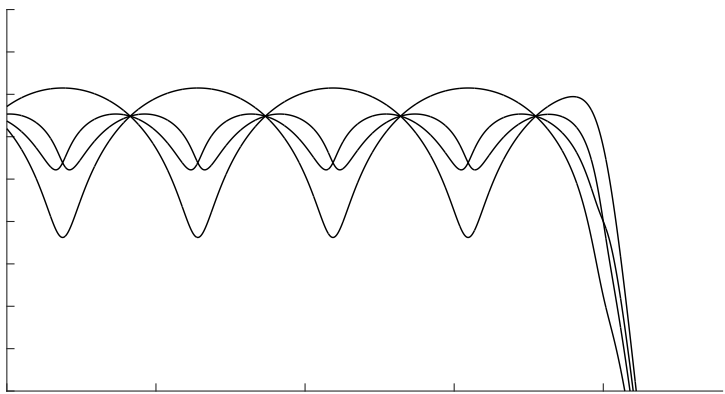}
\end{picture}
\end{center}
\caption{The accumulations curves in the symmetric case.} \label{f56}
\end{figure}

Let us give the asymptotic in the totally symmetric setting: we assume that the triangle $a b c$ is equilateral and that the dynamical quantities are independent of the edge $e \in \SE$. This can be produced by putting the same bump at the three corners of an equilateral triangle in Example \ref{k74}. In this setting, each element of $\widehat{\CQ} ( \tau , h )$ can be written
\begin{equation*}
\widehat{\CQ}_{e , \widetilde{e}} ( \tau , h ) = \CQ_{\infty} ( \tau , h ) \big( \widehat{g}_{+}^{e} \cdot \widehat{g}_{-}^{\widetilde{e}} \big)^{- \frac{1}{2} + i \frac{\tau}{\lambda}} ,
\end{equation*}
if $e^{-} = \widetilde{e}^{+}$ where
\begin{equation} \label{k67}
\CQ_{\infty} ( \tau , h ) = e^{i A / h} \Gamma \Big( \frac{1}{2} - i \frac{\tau}{\lambda} \Big) \frac{1}{\sqrt{2 \pi}} \frac{\CM^{+}}{\CM^{-}} \sqrt{\frac{\vert g_{-} \vert}{\vert g_{+} \vert}} e^{- \frac{\pi}{2} ( \nu + 1 ) i} \big( i \lambda \vert g_{+} \vert \vert g_{-} \vert \big)^{i \frac{\tau}{\lambda}} ,
\end{equation}
is independent of $e , \widetilde{e}$ and $\widehat{g}_{\pm}^{e} = g_{\pm}^{e} / \vert g_{\pm}^{e} \vert$. Moreover, $\widehat{g}_{+}^{e} \cdot \widehat{g}_{-}^{\widetilde{e}} = 1$ when $e$ and $\widetilde{e}$ are conjugate (see Remark \ref{j75}) and $\widehat{g}_{+}^{e} \cdot \widehat{g}_{-}^{\widetilde{e}} = \cos ( \pi / 3 ) = 1 / 2$ otherwise. Hence, the matrix $\widehat{\CQ}$ becomes
\begin{equation} \label{k64}
\widehat{\CQ} = \CQ_{\infty} \left( \begin{array}{cccccc}
0 & 1 & 0 & 0 & \theta & 0  \\
1 & 0 & 0 & \theta & 0 & 0  \\
\theta & 0 & 0 & 1 & 0 & 0  \\
0 & 0 & 1 & 0 & 0 & \theta  \\
0 & 0 & \theta & 0 & 0 & 1  \\
0 & \theta & 0 & 0 & 1 & 0
\end{array} \right) ,
\end{equation}
with
\begin{equation} \label{k68}
\theta ( \tau ) = 2^{\frac{1}{2} - i \frac{\tau}{\lambda}} .
\end{equation}
A direct computation shows that the eigenvalues of the last matrix are
\begin{equation} \label{k66}
\begin{gathered}
\mu_{1} = ( \theta + 1 ) \CQ_{\infty} , \qquad \mu_{2} = ( \theta - 1 ) \CQ_{\infty} ,  \\
\mu_{3} = \mu_{4} = \frac{- \theta + \sqrt{4 - 3 \theta^{2}}}{2} \CQ_{\infty} , \qquad \mu_{5} = \mu_{6} = \frac{- \theta - \sqrt{4 - 3 \theta^{2}}}{2} \CQ_{\infty} .
\end{gathered}
\end{equation}
Inserting \eqref{k66} with \eqref{k67} and \eqref{k68} in \eqref{k65}, we obtain a totally explicit formula for the asymptotic of the resonances. As explained in Section \ref{s74}, Proposition \ref{d91} holds true here. Thus, $P$ has at least two resonances in any neighborhood of size $h \vert \ln h \vert^{- 1}$ of $z_{q , 3} ( \tau ) = z_{q , 4} ( \tau )$ and $z_{q , 5} ( \tau ) = z_{q , 6} ( \tau )$. On the other hand, writing $\mu_{k} = \widehat{\mu}_{k} \CQ_{\infty}$, we get
\begin{equation} \label{k69}
z_{q , k} ( \tau ) = z_{q} ( \tau ) + i \ln ( \widehat{\mu}_{k} ( \tau ) ) \lambda \frac{h}{\vert \ln h \vert} ,
\end{equation}
where $z_{q} ( \tau )$, defined in \eqref{e8}, corresponds to the resonances generated by a unique homoclinic trajectory. The six accumulation curves are described in Figure \ref{f56}. The two central ones, given by $z_{q , 3} = z_{q , 4}$ and $z_{q , 5} = z_{q , 6}$, have multiplicity two. In particular, the number of accumulation curves counted with their multiplicity is equal to the number of edges and bigger than the number of minimal primitive cycles. Since the actions are equal, these curves do not depend on $h$. Lastly, \eqref{k69} shows that they can be seen as the accumulation curve of Corollary \ref{e4} plus periodic functions of $\tau \in \R$.

In case {\rm (B)}, we assume that $\lambda_{a} = \lambda_{c} > \lambda_{b}$. The minimum of the damping is then attained by two primitive cycles $( e_{1} , e_{2} )$ and $( e_{3} , e_{4} )$. Moreover, \eqref{j77} yields
\begin{equation} \label{k70}
D_{0} = \frac{2^{- 1} + 2^{- 1}}{\lambda_{a}^{- 1} + \lambda_{b}^{- 1}} = \frac{\lambda_{a} \lambda_{b}}{\lambda_{a} + \lambda_{b}} .
\end{equation}
In order to give the asymptotic of the resonances, we consider the $2 \times 2$ matrix
\begin{equation} \label{k71}
\widehat{\CQ}^{\rm r e d u} =
\left( \begin{array}{cc}
\widehat{\CQ}_{e_{1} , e_{2}} \widehat{\CQ}_{e_{2} , e_{1}} & \widehat{\CQ}_{e_{1} , e_{2}} \widehat{\CQ}_{e_{2} , e_{4}}  \\
\widehat{\CQ}_{e_{4} , e_{3}} \widehat{\CQ}_{e_{3} , e_{1}} & \widehat{\CQ}_{e_{4} , e_{3}} \widehat{\CQ}_{e_{3} , e_{4}}
\end{array} \right) ,
\end{equation}
where
\begin{equation} \label{k72}
\widehat{\CQ}_{e , \widetilde{e}} ( \tau , h ) = \CQ_{e , \widetilde{e}} \Big( E_{0} + \tau h - i \frac{\lambda_{a} \lambda_{b}}{\lambda_{a} + \lambda_{b}} h , h \Big) .
\end{equation}
Let $\mu_{1} ( \tau , h ) , \mu_{2} ( \tau , h )$ denote the two eigenvalues of $\widehat{\CQ}^{\rm r e d u} ( \tau , h )$. Then, in the sense of Proposition \ref{j81}, the resonances are close to
\begin{equation} \label{k73}
z_{q , k} ( \tau ) = E_{0} + 2 q \pi \frac{\lambda_{a} \lambda_{b}}{\lambda_{a} + \lambda_{b}} \frac{h}{\vert \ln h \vert} - i \frac{\lambda_{a} \lambda_{b}}{\lambda_{a} + \lambda_{b}} h + i \ln ( \mu_{k} ( \tau , h ) ) \frac{\lambda_{a} \lambda_{b}}{\lambda_{a} + \lambda_{b}} \frac{h}{\vert \ln h \vert} ,
\end{equation}
for some $q \in \Z$ and $k \in \{ 1 , 2 \}$. This is proved in Section \ref{s44} using Theorem \ref{j79}. In particular, the resonances concentrate on two accumulation curves, periodic as function of $h^{- 1}$. Roughly speaking, the two primitive cycles can be seen as two interacting homoclinic trajectories around the vertex $b$. This is why the relevant matrix $\widehat{\CQ}^{\rm r e d u}$ has dimension $2$ and not $4$. Thus, the situation is similar to the one of Example \ref{e29} if the angle $\widehat{a b c}$ is acute or of Example \ref{e6} if the angle $\widehat{a b c}$ is obtuse. Nevertheless, the transition phenomenon is different from \eqref{e33} in the last case.

In case {\rm (C)}, we assume that $\lambda_{c} > \lambda_{a} , \lambda_{b}$. The minimum of the damping is attained by the unique primitive cycle $( e_{1} , e_{2} )$ and \eqref{k70} still holds. Hence, we are in the situation of Corollary \ref{k40} and the resonances are close to
\begin{equation}
z_{q} ( \tau ) = E_{0} + 2 q \pi \frac{\lambda_{a} \lambda_{b}}{\lambda_{a} + \lambda_{b}} \frac{h}{\vert \ln h \vert} - i \frac{\lambda_{a} \lambda_{b}}{\lambda_{a} + \lambda_{b}} h + i \ln ( \mu ( \tau , h ) ) \frac{\lambda_{a} \lambda_{b}}{\lambda_{a} + \lambda_{b}} \frac{h}{\vert \ln h \vert} ,
\end{equation}
in the sense of Proposition \ref{j81} with
\begin{equation*}
\mu ( \tau , h ) = \CQ_{e_{1} , e_{2}} \CQ_{e_{2} , e_{1}} \Big( E_{0} + \tau h - i \frac{\lambda_{a} \lambda_{b}}{\lambda_{a} + \lambda_{b}} h , h \Big) .
\end{equation*}
Note that $\mu$ is also given by \eqref{k47} and the unique accumulation curve (independent of $h$) satisfies \eqref{n10}. This is natural since the barrier $c$ plays a negligible role in the asymptotic of these resonances, and that we are exactly in the setting of Example \ref{k42} without this barrier (Corollary \ref{k54} makes this argument rigorous).

\section{Resonant states} \label{s32}

In this section, we are interested in resonant states. Some of their basic properties can be found in the appendix of Petkov and the first author \cite{BoPe13_01}. These generalized eigenvectors play an important role for the resonance expansion of the propagator (see e.g. Lax and Phillips \cite[Theorem III.5.4]{LaPh67_01}) and for the (residue of the) scattering amplitude (see e.g. Lahmar-Benbernou and Martinez \cite{LaMa99_01}). With the notations of Section \ref{s2}, they are defined by

\begin{definition}[Resonant states]\sl \label{j34}
For $z \in \CE_{\theta} \cap \res ( P )$ and $0 \neq u \in H^{2} ( \R^{n} )$, we say that $u$ is a resonant state of $P$ associated to the resonance $z$ if and only if $( P_{\theta} - z ) u = 0$.
\end{definition}

This definition depends in an unessential way on the distortion angle $\theta$ and the function $F$. Indeed, for a resonance $z \in \CE_{\theta}$, the distorted and local resolvents can be written
\begin{equation*}
( P_{\theta} - \lambda )^{- 1} = \sum_{j = 1}^{J_{\theta}} \frac{\Pi^{\theta}_{j}}{( \lambda - z )^{j}} + R_{\theta} ( \lambda ) \qquad \text{and} \qquad ( P - \lambda )^{- 1} = \sum_{j = 1}^{J} \frac{\Pi_{j}}{( \lambda - z )^{j}} + R ( \lambda ),
\end{equation*}
as operators on $L^{2} ( \R^{n} )$ and from $L^{2}_{\rm comp} ( \R^{n} )$ to $L^{2}_{\rm loc} ( \R^{n} )$, respectively. The $\Pi_{j}^{\bullet}$'s are finite rank operators and the operator-valued functions $R_{\bullet} ( \lambda )$ are holomorphic near $z$. By definition, the resonant states are the functions $u_{\theta} \in \im \Pi^{\theta}_{1}$ such that $( P_{\theta} - z ) u_{\theta} = 0$. More generally, a function $u_{\theta}$ belongs to $\im \Pi^{\theta}_{1}$ if and only if $( P_{\theta} - z )^{N} u_{\theta} = 0$ for some $N \geq 1$. These functions are called the generalized resonant states. We will only consider the resonant states in this paper. Let $\chi \in C^{\infty}_{0} ( \R^{n} )$ be supported outside of the complex dilation and such that $\chi = 1$ on a sufficiently large neighborhood of $0$. In \cite[Appendix]{BoPe13_01}, it is proved that for all $\theta , \theta ^{\prime}$ small enough and any resonant state $u_{\theta}$ associated to a resonance $z \in \CE_{\theta} \cap \CE_{\theta^{\prime}}$, there exist unique functions $u \in \im \Pi_{1}$ and $u_{\theta^{\prime}} \in \im \Pi_{1}^{\theta^{\prime}}$ such that $( P - z ) u = ( P_{\theta^{\prime}} - z ) u_{\theta^{\prime}} = 0$ and $\chi u_{\theta} = \chi u = \chi u_{\theta^{\prime}}$. Moreover, the next result shows that it is equivalent to normalize the resonant states on the whole space $\R^{n}$ or on the compact subsets. Its proof can be found in Appendix \ref{s41}. Note also that it is enough to describe the resonant states $u_{\theta}$ near the trapped set $K ( E_{0} )$ since they satisfy $( P_{\theta} - z ) u_{\theta} = 0$ which can be seen as a propagation equation.

\begin{remark}\sl \label{j48}
We assume \ref{h1}. Let $\theta = h \vert \ln h \vert$, $R > 0$ and $\chi \in C^{\infty}_{0} ( \R^{n} ; [ 0 , 1 ])$ be such that $\chi = 1$ on a sufficiently large neighborhood of $0$. We suppose that $P_{\theta} = P$ near the support of $\chi$. Then, there exists $C > 0$ such that, for any resonant state $u$ associated to a resonance $z \in B ( E_{0} , R h )$, we have
\begin{equation*}
\Vert \chi u \Vert_{L^{2} ( \R^{n} )} \leq \Vert u \Vert_{L^{2} ( \R^{n} )} \leq h^{- C} \Vert \chi u \Vert_{L^{2} ( \R^{n} )} ,
\end{equation*}
for $h$ small enough.
\end{remark}

This result is in some sense similar to the resolvent estimates of Proposition \ref{j66}. As in that case, if the distortion function $F$ is well-chosen (that is $H_{p} ( F ( x ) \cdot \xi ) \geq 0$ on the whole energy surface $p^{- 1} ( E_{0} )$) or if the distortion angle satisfies $\theta = \CO ( h )$, it could be possible to replace $h^{- C}$ by $C$ in the remark.

In the semiclassical regime, the asymptotic of the resonant states is already known for shape and barrier-top resonances. In the well in a island situation, Helffer and Sj\"{o}strand \cite[Section 9]{HeSj86_01} have proved that the generalized spectral projection $\Pi_{1}^{\theta}$ associated to an isolated resonance is exponentially close to the spectral projection of the operator restricted to the well. So the resonant states are essentially localized in the well and $\Pi_{1}^{\theta}$ is almost orthogonal. On the other hand, we have obtained in a previous work \cite[Section 4]{BoFuRaZe11_01} the asymptotic of the generalized spectral projection associated to isolated resonances generated by the maximum of the potential. In particular, the resonant states are Lagrangian distributions carried by all the manifold $\Lambda_{+}$ and the norm of $\Pi_{1}^{\theta}$ behaves like a negative power of $h$.

The main result of this part is stated under the assumptions of Section \ref{s61}. The homoclinic set is then the union of $K$ trajectories $\gamma_{1} , \ldots , \gamma_{K}$ along which $\Lambda_{-}$ and $\Lambda_{+}$ intersect transversally. Let us define some geometric quantities. We first fix a point $\rho_{-}^{k} = ( x_{-}^{k} ,\xi_{-}^{k} ) = \gamma_{k} ( t_{-}^{k} )$ on each homoclinic curve $\gamma_{k}$. Let $\Lambda_{+}^{1}$ denote the evolution of $\Lambda_{+}^{0}$ by the Hamiltonian flow after a turn along $\CH$. If $t_{-}^{k}$ is chosen large enough, $\rho_{-}^{k}$ belongs to $\Lambda_{+}^{1}$ and this manifold projects nicely on the $x$-space near $\rho_{-}^{k}$ (see Section \ref{s72}). Then, there exists a unique generating function $\varphi_{+}^{1} \in C^{\infty} ( \R^{n} )$ of $\Lambda_{+}^{1}$ (i.e. $\Lambda_{+}^{1} = \{ ( x , \nabla \varphi_{+}^{1} (x) ) \}$) defined near $x_{-}^{k}$
with the normalization
\begin{equation}
\varphi_{+}^{1} ( x_{-}^{k} ) = \int_{\gamma_{k} ( ] - \infty , t_{-}^{k} ] )} \xi \cdot d x .
\end{equation}
Recall that $A_{k}$ and $\CM_{k}^{-}$ have been defined in \eqref{r5} and \eqref{d7}. Furthermore, with the notations of \eqref{d7}, let
\begin{equation}
\CD_{k} (t) = \sqrt{\Big\vert \det \frac{\partial x_{k} ( s , y )}{\partial ( s , y )} \vert_{s = t , \ y = 0} \Big\vert} ,
\end{equation}
be the Maslov determinant at time $t$. The resonant states satisfy the following asymptotic, whose proof can be found in Section \ref{s40}.

\begin{theorem}[Behavior of the resonant states]\sl \label{j35}
Assume \ref{h1}--\ref{h4}, \ref{h8} and fix $C , \delta > 0$. Let $v = v (h)$ be a family of normalized resonant states associated to a resonance $z = z ( h ) \in \eqref{d90}$. Then, we can find $h^{M} \leq c ( h ) \leq h^{- M}$ for some $M > 0$, such that $u = c v$ has the following properties.

$i)$ The microsupport of $u$ is contained in $\{ ( 0 , 0 ) \} \cup \Lambda_{+}$.

$ii)$ The function $u$ is in $\CI ( \Lambda_{+} , 1 )$ microlocally near any point of $\CH$. We then write
\begin{equation} \label{j36}
u ( x , h ) = e^{- i A_{k} / h} e^{i \frac{z - E_{0}}{h} t_{-}^{k}} \frac{\CM_{k}^{-}}{\CD_{k} ( t_{-}^{k} )} a_{-}^{k} ( x , h ) e^{i \varphi_{+}^{1} (x) / h} ,
\end{equation}
microlocally near $\rho_{-}^{k}$, for some $a_{-}^{k} \in S ( 1 )$.

$iii)$ Let $\SA_{0} ( h ) = ( \SA_{0}^{1} ( h ) , \ldots , \SA_{0}^{K} ( h ) ) \in \C^{K}$ be defined by $\SA_{0}^{k} ( h ) = a_{-}^{k} ( x_{-}^{k} , h )$. This $K$-vector satisfies the equation
\begin{equation} \label{j38}
\big( h^{S ( z , h ) / \lambda_{1} - 1 / 2} \CQ ( z , h ) - 1 \big) \SA_{0} ( h ) = o ( 1 ) ,
\end{equation}
as $h$ goes to $0$, with the normalization $\Vert \SA_{0} ( h ) \Vert_{\ell^{2}} = 1$.
\end{theorem}

In the previous result, $\Vert \cdot \Vert_{\ell^{2}}$ denotes the Euclidean norm in $\C^{K}$. The first point of the theorem implies that $u$ vanishes microlocally near each point $\rho \in T^{*} \R^{n}$ such that $p ( \rho ) \neq E_{0}$, $\rho \in \Lambda_{-} \setminus \Lambda_{+}$ or $\rho$ in the incoming region (see Section \ref{s33}). The normalization factor in \eqref{j36} is natural and can be explained by the asymptotic of $u$ at the barrier-top along the characteristic $\pi_{x} ( \gamma_{k} )$. More precisely, let $b_{-}^{k} \in S ( 1 )$ denote the symbol of $u$ near $\gamma_{k} ( [ t_{-}^{k} , + \infty [ )$ without normalization (i.e. $u = b_{-}^{k} e^{i \varphi_{+}^{1} / h}$ microlocally near each point of $\gamma_{k} ( [ t_{-}^{k} , + \infty [ )$). Then, we have
\begin{equation} \label{j37}
b_{-}^{k} ( x_{k} ( t ) , h ) = e^{- i A_{k} / h} \SA_{0}^{k} ( h ) e^{t \big( \frac{\lambda_{1}}{2} + i \tau \big)} \big( 1 + o_{t \to + \infty} ( 1 ) \big) + o_{h \to 0}^{t} ( 1 ) ,
\end{equation}
where $\tau = ( \re z - E_{0} ) / h$ and $o_{h \to 0}^{t} ( 1 )$ is a function which tends to $0$ as $h$ goes to $0$ for $t$ fixed. The proof of \eqref{j37} can be found in Section \ref{s40}. Without the normalization \eqref{j36}, the vector of the $b_{-}^{k} ( x_{-}^{k} , h )$ would satisfy an equation similar to \eqref{j38} with $\CQ$ replaced by $\CU^{- 1} \CQ \CU$ for some invertible diagonal matrix $\CU$. In some sense, our choice of normalization corresponds formally to the computation of the different quantities at $x_{-}^{k} = 0$.

\begin{remark}\sl \label{j39}
$i)$ From Proposition \ref{d9} and Theorem \ref{d8}, any resonance $z \in \eqref{d90}$ is close to some $z_{q , k} ( \tau )$ given by \eqref{d95} with $q ( h ) \in \Z$, $k ( h ) \in \{ 1 , \ldots , K \}$ and $\re z = E_{0} + \tau h + o ( h )$. In this case, $h^{- S ( z , h ) / \lambda_{1} + 1 / 2} = \mu_{k} ( \tau ,h ) + o ( 1 )$ and \eqref{j38} gives
\begin{equation} \label{l7}
\big( \widehat{\CQ} ( \tau , h ) - \mu_{k} ( \tau , h ) \big) \SA_{0} ( h ) = o ( 1 ) .
\end{equation}

$ii)$ In particular, assume that $\mu_{k} ( \tau , h )$ is a isolated simple eigenvalue in the sense that $\spe ( \widehat{\CQ} ( \tau , h ) ) \cap B ( \mu_{k} , \varepsilon ) = \{ \mu_{k} \}$ for some constant $\varepsilon > 0$. Let $f_{k} ( \tau , h )$ denote a normalized eigenvector of $\widehat{\CQ} ( \tau , h )$ associated to the eigenvalue $\mu_{k} ( \tau , h )$. Then, there exists a complex number $\alpha ( z , h )$ with $\vert \alpha \vert = 1$ such that
\begin{equation} \label{j40}
\SA_{0} = \alpha f_{k} + o( 1 ) ,
\end{equation}
as $h$ goes to $0$. In the general situation, $\SA_{0}$ can be approximated by a linear combination of the generalized eigenvectors associated to the eigenvalues of $\widehat{\CQ}$ near $\mu_{k}$.

$iii)$ Conversely, let $\mu_{k} ( \tau , h )$ be an isolated simple eigenvalue of $\widehat{\CQ} ( \tau , h )$ avoiding a vicinity of $0$ and let $f_{k} ( \tau , h )$ be an associated normalized eigenvector. Then, for all resonance $z$ close to $z_{q , k} ( \tau )$ (in the sense of Proposition \ref{d9}), there exists a resonant state $u$ associated to $z$, normalized as in Theorem \ref{j35} $iii)$, such that
\begin{equation} \label{l6}
\SA_{0} = f_{k} + o( 1 ) ,
\end{equation}
as $h$ goes to $0$.
\end{remark}

\begin{example}\rm \label{l16}
To illustrate the previous remark, we return to Example \ref{e6} and assume that $C_{1} \neq C_{2}$ (see \eqref{e17}). In the limit $\tau \to - \infty$, the two eigenvalues $\mu_{1} , \mu_{2}$ of $\widehat{\CQ} ( \tau , h )$ are isolated and simple since their modulus satisfy
\begin{equation} \label{l14}
\vert \mu_{k} \vert = e^{C_{k} / \lambda} + o_{\tau \to - \infty} ( 1 ) .
\end{equation}
Moreover, there exists a basis of normalized eigenvectors $f_{1} , f_{2}$ of $\widehat{\CQ} ( \tau , h )$ of the form
\begin{equation} \label{j46}
f_{1} = \left( \begin{array}{c}
1 \\
0
\end{array} \right) + o_{\tau \to - \infty} ( 1 ) \qquad \text{and} \qquad f_{2} = \left( \begin{array}{c}
0 \\
1
\end{array} \right) + o_{\tau \to - \infty} ( 1 ) ,
\end{equation}
uniformly with respect to $h$ (see Section \ref{s40}). Then, \eqref{j40} implies that the resonant states associated to the resonances close to $z_{q , k} ( \tau )$ are mainly localized on the characteristic curve $\pi_{x} ( \gamma_{k} )$. This means that their symbol is elliptic near this trajectory and $o_{\tau \to - \infty} ( 1 ) + o^{\tau}_{h \to 0} ( 1 )$ on the other one. This phenomenon can be interpreted as in Section \ref{s14}: for energies below $E_{0}$, the two homoclinic trajectories are ``disconnected'' and then the resonances as well as the resonant states are given by each trajectory separately.

On the contrary, when $\tau$ tends to $+ \infty$, the two eigenvalues $\mu_{1} , \mu_{2}$ satisfy
\begin{equation*}
\mu_{1} = - \mu_{2} + o_{\tau \to + \infty} ( 1 ) \qquad \text{and} \qquad \vert \mu_{k} \vert = e^{C_{\widehat{12}} / \lambda} + o_{\tau \to + \infty} ( 1 ) .
\end{equation*}
where $C_{\widehat{12}}$ is defined in \eqref{e17}. Moreover, there exists a basis of normalized eigenvectors $g_{1} , g_{2}$ of $\widehat{\CQ} ( \tau , h )$ of the form
\begin{equation} \label{j47}
g_{1} = \left( \begin{array}{c}
1 \\
\alpha
\end{array} \right) + o_{\tau \to - \infty} ( 1 ) \qquad \text{and} \qquad g_{2} = \left( \begin{array}{c}
1 \\
- \alpha
\end{array} \right) + o_{\tau \to - \infty} ( 1 ) ,
\end{equation}
for some $\alpha = \alpha ( \tau , h )$ which has a non-zero constant modulus. The explicit expression of $\alpha$ can be found in \eqref{l15}. Thus, any resonant state ``fills'' the two homoclinic trajectories $\gamma_{1}, \gamma_{2}$. This is in agreement with the intuition that these two trajectories recombine for energies higher than $E_{0}$ (see Section \ref{s14}).
\end{example}

As this stage, the results already obtained allow us to recover the symbol $a_{-}^{k}$ of the resonant states in a neighborhood of $\pi_{x} ( \gamma_{k} )$ up to $o ( 1 )$. This follows from the initial condition \eqref{j40} and Lemma \ref{d41} which give $a_{-}^{k}$ in a vicinity of $x_{-}^{k}$, and from the evolution equation $( P_{\theta} - z ) u = 0$. Under some additional hypotheses, it is possible to describe the resonant states modulo $\CO ( h^{\infty} )$. This question is intimately related to the asymptotic of the resonances modulo $\CO ( h^{\infty} )$. Thus, we make the assumptions of Section \ref{s79}. In particular, the homoclinic set $\CH$ is reduced to a single trajectory and we remove the subscript $k$. One can then improve the preceding result as follows. Let $v$ be as in Theorem \ref{j35}. Multiplying $c$ by a complex number of modulus $1$, we can always assume that $\SA_{0} = 1$. Then, the symbol $a_{-}$ satisfies the following asymptotic
\begin{equation} \label{j41}
a_{-} ( x , h ) \simeq \sum_{a = 0}^{+ \infty} \sum_{b = 0}^{B_{a}} \sum_{c = 0}^{C_{a}} a_{-}^{a , b , c} ( x , z , h ) \big( h^{S ( z , h ) / \lambda_{1} - 1 / 2} \big)^{b} ( \ln h )^{c} h^{\widehat{\mu}_{a} / \lambda_{1}} ,
\end{equation}
near $x_{-}$. The symbols $a_{-}^{a , b , c} \in S ( 1 )$ are not explicit but can formally be computed using the values of $p ( x , \xi )$ in a neighborhood of $K ( E_{0} )$ (see Section \ref{s80}). In particular, they can be written $a_{-}^{a , b , c} ( x , z , h ) = e^{i A / h} \widetilde{a}_{-}^{a , b , c} ( x, \sigma )$ with $\widetilde{a}_{-}^{a , b , c}$ holomorphic in $\sigma$ defined in \eqref{d92} (see also \eqref{g14}). The asymptotic \eqref{j41} is a direct consequence of \eqref{g18}. By propagation of singularities, $u$ can be computed modulo $\CO ( h^{\infty} )$ microlocally near any point of $\CH$.

Remark that Theorem \ref{j35} does not give the asymptotic of the resonant states near $( 0 , 0 )$ and near $\Lambda_{+} \setminus \CH$. Nevertheless, if $u$ is known modulo $\CO ( h^{\infty} )$ near $x_{-}^{k}$, \cite[Theorem 5.1]{BoFuRaZe07_01} provides a formula for $u$ in a vicinity of $( 0 , 0 )$. Then, the usual propagation of singularities yields the asymptotic of $u$ on all $\Lambda_{+}$ and hence microlocally near each point of $T^{*} \R^{n}$. Note that \cite{ALBoRa08_01} gives some informations in the orthogonal of the $g_{-}^{k}$'s directions.

As shown in the following example, the resonant states can be unevenly distributed on the homoclinic set.

\begin{example}\rm \label{l17}
We consider again the geometric setting of Example \ref{e29} with three homoclinic trajectories. We assume that $\gamma_{1}$ and $\gamma_{2}$ are symmetric with respect to $\gamma_{3}$ as in Figure \ref{f28} and that all the geometric quantities ($A_{\bullet}$, $g_{\pm}^{\bullet}$, $\nu_{\bullet}$, \ldots) are the same for $\gamma_{1}$ and $\gamma_{2}$. Then, it is proved in Section \ref{s40} that $f = {}^{t} ( 1 , - 1 , 0 ) / \sqrt{2}$ is an eigenvector of $\widehat{\CQ} ( \tau , h )$ for all $\tau \in [ - C , C ]$ and $h$ small enough. Let $\mu_{1} ( \tau , h )$ denote the associated eigenvalue. The geometry can be chosen such that $\mu_{1} ( \tau , h )$ is simple, isolated and far away from $0$ uniformly with respect to $\tau , h$. We consider the resonances close to $z_{q , 1} ( \tau ) \in \eqref{d90}$.

Then, Remark \ref{j39} $iii)$ provides a family of resonant states $u$ normalized as in Theorem \ref{j35} $iii)$ satisfying $\SA_{0} = f + o ( 1 )$. In particular, $u$ is a Lagrangian distribution microlocally near any point of $\CH$ and the leading term of its symbol vanishes on $\gamma_{3}$. In other words, $u$ ``fills'' the trajectories $\gamma_{1}$ and $\gamma_{2}$ but not $\gamma_{3}$. Nevertheless, it may happen that the symbol of $u$ vanishes only at the first order on $\gamma_{3}$ (modulo $o ( 1 )$).
\end{example}

\begin{example}\rm
In the setting of Section \ref{s81}, we assume that the homoclinic set consists of two trajectories $\gamma_{1} , \gamma_{2}$ such that $\Lambda_{-}$ and $\Lambda_{+}$ intersect tangentially of order $m_{1} \geq 2$ along $\gamma_{1}$ and transversally along $\gamma_{2}$ (in other words, $m_{2} = 1$). The geometric situation is illustrated in Figure \ref{f32}. From Theorem \ref{g24}, the resonances close to the real axis are generated by $\gamma_{1}$ (see \eqref{g79}).

On the other hand, one can adapt Theorem \ref{j35} to the case of tangential intersection of finite order. In fact, the role of the total homoclinic set $\CH$ is played only by the homoclinic trajectories of the highest order of tangency (that is $\gamma_{1}$ in this example). In the proof, Lemma \ref{d41} is replaced by Lemma \ref{g46}. Thus, any resonant state $v$ can be polynomially normalized on $\gamma_{1}$. This means that there exists $h^{M} \leq c \leq h^{- M}$ for some $M > 0$ such that $u = c v$ belongs to $\CI ( \Lambda_{+} , 1 )$ microlocally near any point of $\gamma_{1}$ and $\SA_{0}^{1} = 1$. Furthermore, one can show that $u \in \CI ( \Lambda_{+} , 1 )$ with an elliptic symbol microlocally near any point of $\gamma_{2}$. This follows from Lemma \ref{g46}.

Summing up, the resonant states are of the same order on $\gamma_{1}$ and $\gamma_{2}$, whereas the corresponding resonances are only provided by $\gamma_{1}$. Thus, the size of the resonant states near a homoclinic trajectory does not reflect the role of the latter in the asymptotic of the resonances.
\end{example}

We can also generalize the previous results to the setting of heteroclinic trajectories studied in Section \ref{s42}. With the notations of this part, the microsupport of any normalized resonant state $u$ is contained in the set
\begin{equation}
\bigcup_{v \in \SV} \{ ( v , 0 ) \}\cup \bigcup_{e \in \SE} \Lambda_{+}^{e^{-}} .
\end{equation}
Moreover, $u$ is a Lagrangian distribution with Lagrangian manifold $\Lambda_{+}^{e^{-}}$ near each point of $e \in \SE$. Nevertheless, it is more difficult to adapt \eqref{j38} in this situation since the order of these Lagrangian distributions may no longer be uniform and may depend on the edge.

\begin{remark}[Weird distribution of the resonant states]\sl \label{k94}
In the setting of Section \ref{s42}, it may happen that the resonant states are small and large at different edges of the same cycle. More precisely, assume for simplicity that there is a unique minimal primitive cycle in $( \SV , \SE )$ as in Section \ref{s71}. Then, the order of the resonant states as Lagrangian distribution may change by a power of $h$ between the edges of this cycle.
\end{remark}

That the transition through a vertex amplifies or absorbs the resonant states explains this phenomenon. This is why the proof of Theorem \ref{j79} uses weights ($h^{N_{e}}$ and $h^{N_{e \leftarrow \widetilde{e}}}$) which respect the graph structure (see Section \ref{s50}). The setting here is different from Example \ref{l17}. In this example, some resonant states are small along $\gamma_{3}$ but this trajectory does not contribute to the creation of the corresponding resonances (``removing'' $\gamma_{3}$ will only modify these resonances by $o ( h \vert \ln h \vert^{- 1} )$). On the contrary, the resonances are generated by the combination of all the edges of the unique minimal primitive cycle in Remark \ref{k94} (removing any edge will destroy this cycle). The following example illustrates this remark.

\begin{example}\rm \label{k95}
We come back to Example \ref{k42} in dimension $n \geq 2$ and use its notations. The graph $( \SV , \SE )$ consists of two vertices and two heteroclinic trajectories (see Figure \ref{f51}). In particular, the resonances are generated by the unique primitive cycle $( e_{1} , e_{2} )$.

For $k = 1 , 2$, we fix a point $x^{k}_{-} \in \pi_{x} ( e_{k} )$ close to $v_{3 - k}$. Then, for all family of normalized resonant state $v$ associated to a family of resonances $z \in \eqref{j80}$, there exists a (unique) complex number $h^{M} \leq c ( h ) \leq h^{- M}$ for some $M > 0$ such that $u = c v$ satisfies
\begin{equation} \label{k84}
u_{-}^{1} : = u \in \CI ( \Lambda^{e_{1}^{-}}_{+} , 1 ) \text{ microlocally near each point of } e_{1} ,
\end{equation}
and
\begin{equation} \label{k85}
u_{-}^{2} : = u \in \CI \Big( \Lambda^{e_{2}^{-}}_{+} , h^{\frac{( n - 1 ) ( \lambda_{2} - \lambda_{1} )}{2 \lambda_{1} + 2 \lambda_{2}}} \Big) \text{ microlocally near each point of } e_{2} .
\end{equation}
Moreover, if the symbol $a_{-}^{k}$ of $u _{-}^{k}$ is defined as in \eqref{j36}, we have
\begin{equation} \label{k82}
a_{-}^{1} ( x_{-}^{1} , h ) = 1 ,
\end{equation}
and
\begin{equation} \label{k83}
a_{-}^{2} ( x_{-}^{2} , h ) \sim e^{i A / h} \Gamma ( S_{2} ) \frac{\lambda_{2}^{\frac{1}{2} - S_{2}}}{\sqrt{2 \pi}} e^{- i \frac{\pi}{2} ( S_{2} + \frac{1}{2} )} \frac{\CM_{e_{2}}^{+}}{\CM_{e_{2}}^{-}} \frac{\vert g_{-}^{e_{1}} \vert^{1 - S_{2}}}{\vert g_{+}^{e_{2}} \vert^{S_{2}}} h^{\frac{( n - 1 ) ( \lambda_{2} - \lambda_{1} )}{2 \lambda_{1} + 2 \lambda_{2}}} ,
\end{equation}
with the notations of \eqref{k47}. Equation \eqref{k82} can be seen as a normalization relation which guaranties the uniqueness of $c$. The two previous equations show that the symbols $a_{-}^{k}$ are elliptic near $\rho^{k}_{-}$. Thus, the order of Lagrangian distributions in \eqref{k84} and \eqref{k85} is sharp.

Roughly speaking, the above discussion gives
\begin{align*}
\begin{aligned}
&u_{-}^{1} \approx 1 &&\text{and} &&u_{-}^{2} \approx h^{\delta} &&\text{if } \lambda_{1} < \lambda_{2} ,  \\
&u_{-}^{1} \approx 1 &&\text{and} &&u_{-}^{2} \approx 1 &&\text{if } \lambda_{1} = \lambda_{2} , \\
&u_{-}^{1} \approx 1 &&\text{and} &&u_{-}^{2} \approx h^{- \delta} &&\text{if } \lambda_{1} > \lambda_{2} ,
\end{aligned}
\end{align*}
for some $\delta > 0$. Thus, the resonant states can be of different order on the two edges, whereas these conjugate trajectories play a symmetric role in the unique primitive cycle $( e_{1} , e_{2} )$. On the other hand, the edge where all the resonant states are the largest depends on the values of the $\lambda_{k}$'s.
\end{example}

This last example works only in dimension $n \geq 2$. In dimension $n = 1$, \eqref{k84}--\eqref{k83} imply that any resonant state is at the same order on the two edges $e_{1} , e_{2}$. This is a general fact. Indeed, assume that a one dimensional operator satisfying the hypotheses of Theorem \ref{j79} has a unique primitive cycle (recall that all the cycles are minimal in dimension one), then the order of any resonant state associated to any resonance in \eqref{j80} is the same on all the edges of this cycle. Note that the behavior near the vertices is not considered here.

The distribution of the resonant states can even be worse. One may think that they are always of maximal order in the regions of the phase space which generate the corresponding resonances and give their asymptotic (that is neighborhoods of the trajectories providing the dynamical quantities appearing in the quantization rule). This intuition is correct for punctual well in an island situations (see Helffer and Sj\"{o}strand \cite[Section 9--10]{HeSj86_01}) or barrier-top resonances (see our previous work \cite[Section 4]{BoFuRaZe11_01}). Nevertheless, it is wrong in general. There exist situations satisfying the assumptions of Section \ref{s45} where all the resonant states are small near the cycles of $( \SV , \SE )$ and large in irrelevant parts of the trapped set. More precisely, we have the following statement.

\begin{remark}[Delocalization of the resonant states]\sl \label{k81}
For $C_{0} > 0$, there exist operators $P$ satisfying the assumptions of Theorem \ref{j79} (but not of Schr\"{o}dinger form) with resonances in \eqref{j80} and two disjoint compact subsets $K_{1} , K_{2}$ of $K ( E_{0} )$ such that 

$i)$ All the vertices and edges of the cycles in $( \SV , \SE )$ belong to $K_{1}$. In particular, the asymptotic of the resonances (modulo $o ( h \vert \ln h \vert^{- 1} )$) is given by the symbol of $P$ in any neighborhood of $K_{1}$ from Proposition \ref{j81} and Theorem \ref{j79}.

$ii)$ There exists $\varphi_{1} \in C^{\infty}_{0} ( T^{*} \R^{n} )$ with $\varphi_{1} = 1$ near $K_{1}$ such that, for all $\varphi_{2} \in C^{\infty}_{0} ( T^{*} \R^{n} )$ with $\varphi_{2} = 1$ near $K_{2}$, we have
\begin{equation} \label{k86}
\Vert \Op ( \varphi_{1} ) u \Vert \leq h^{C_{0}} \Vert \Op ( \varphi_{2} ) u \Vert ,
\end{equation}
for $h$ small enough and any resonant state $u$ associated to any resonance in \eqref{j80}.
\end{remark}

In this remark, the constant $C_{0}$ depends on the operator $P$. In particular, it is not possible to replace $h^{C_{0}}$ by $\CO ( h^{\infty} )$ in \eqref{k86}. Indeed, in that case, a normalized resonant state $u$ would satisfy $u = 0$ microlocally near each point of $K ( E_{0} )$ thanks to Theorem \ref{a32}, which is impossible from Proposition \ref{a16}. However, there might be other types of trapped geometries whose resonant states are $\CO ( h^{\infty} )$ near the region which generates the corresponding resonances. The same way, there exists a constant $M > 0$ such that $\Vert \Op ( \varphi_{1} ) u \Vert \geq h^{M}$ for $h$ small enough and any normalized resonant state $u$ associated to any resonance in \eqref{j80}.

That a barrier-top may amplify functions in the kernel of $P - z$ by a negative power of $h$ explains the delocalization. More precisely, Corollary \ref{d46} shows that the order of a transversal Lagrangian distribution (i.e. a Lagrangian distribution whose Lagrangian manifold intersects transversally $\Lambda_{-}$ along a trajectory) is multiplied by
\begin{equation*}
h^{( \sum_{j \geq 2} \lambda_{j}^{v} / 2 - D_{0} ) / \lambda_{1}^{v}} ,
\end{equation*}
passing through the vertex $v$ when $\im z \approx - D_{0} h$. Thus, the phenomenon in Remark \ref{k81} is due to a non self-adjoint effect coming from the imaginary part of the spectral parameter $z$. We now give examples of operators justifying Remark \ref{k81}.

\begin{figure}
\begin{center}
\begin{picture}(0,0)%
\includegraphics{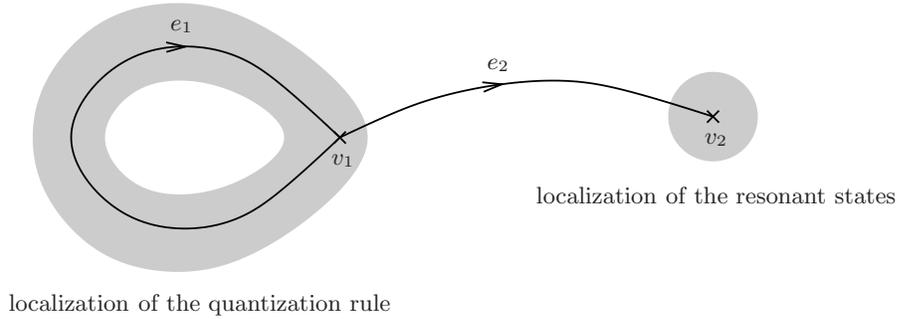}%
\end{picture}%
\setlength{\unitlength}{987sp}%
\begingroup\makeatletter\ifx\SetFigFont\undefined%
\gdef\SetFigFont#1#2#3#4#5{%
  \reset@font\fontsize{#1}{#2pt}%
  \fontfamily{#3}\fontseries{#4}\fontshape{#5}%
  \selectfont}%
\fi\endgroup%
\begin{picture}(18253,7785)(-34117,-8366)
\put(-26324,-4636){\makebox(0,0)[b]{\smash{{\SetFigFont{9}{10.8}{\rmdefault}{\mddefault}{\updefault}$v_{1}$}}}}
\put(-30374,-1261){\makebox(0,0)[b]{\smash{{\SetFigFont{9}{10.8}{\rmdefault}{\mddefault}{\updefault}$e_{1}$}}}}
\put(-16949,-4111){\makebox(0,0)[b]{\smash{{\SetFigFont{9}{10.8}{\rmdefault}{\mddefault}{\updefault}$v_{2}$}}}}
\put(-22424,-2236){\makebox(0,0)[b]{\smash{{\SetFigFont{9}{10.8}{\rmdefault}{\mddefault}{\updefault}$e_{2}$}}}}
\put(-29924,-8311){\makebox(0,0)[b]{\smash{{\SetFigFont{9}{10.8}{\rmdefault}{\mddefault}{\updefault}localization of the quantization rule}}}}
\put(-16949,-5611){\makebox(0,0)[b]{\smash{{\SetFigFont{9}{10.8}{\rmdefault}{\mddefault}{\updefault}localization of the resonant states}}}}
\end{picture}%
\end{center}
\caption{The trapped set in Example \ref{k87}.} \label{f59}
\end{figure}

\begin{example}\rm \label{k87}
In dimension $n \geq 2$, we consider an operator satisfying the assumptions of Theorem \ref{j79} and whose graph $( \SV , \SE )$ is described in Figure \ref{f59}. For simplicity, we assume that both fixed points are isotropic
\begin{equation*}
\lambda_{1}^{1} = \cdots = \lambda_{n}^{1} = : \lambda_{1} \qquad \text{and} \qquad \lambda_{1}^{2} = \cdots = \lambda_{n}^{2} = : \lambda_{2} .
\end{equation*}
In this case, there is a unique primitive cycle (reduced to the homoclinic trajectory $e_{1}$). Then, Corollary \ref{k40} provides the distribution of the resonances. In particular, there are a lot of them in \eqref{j80}. Their asymptotic is given by dynamical quantities defined in any neighborhood of $e_{1} \cup \{ ( v_{1} , 0 ) \}$ modulo $o ( h \vert \ln h \vert^{- 1} )$. Adapting Section \ref{s79}, one could perhaps replace this remainder term by $\CO ( h^{\infty} )$.

We now explain how to realize such a trapped set. Of course, this is not possible for Schr\"{o}dinger operators due to the symmetry of the Hamiltonian curves (see above Example \ref{b81}). The simplest way to realize Figure \ref{f59} is to consider three aligned radial bumps (one higher than $E_{0}$ and two of height $E_{0}$) and then to add an absorbing potential on the heteroclinic trajectory from $v_{2}$ to $v_{1}$ (see Remark \ref{c13}). The obtained operator is dissipative. In Section \ref{s90}, we construct a self-adjoint pseudodifferential operator with such a trapped set. Finally, it is perhaps also possible to realize such a geometric situation for a Schr\"{o}dinger operator with a magnetic potential.

We choose the compact subsets
\begin{equation}
K_{1} = e_{1} \cup \{ ( v_{1} , 0 ) \} \qquad \text{and} \qquad K_{2} = \{ ( v_{2} , 0 ) \} .
\end{equation}
These sets are disjoint and Remark \ref{k81} $i)$ is automatically satisfied. On the other hand, it is proved in Section \ref{s40} that \eqref{k86} holds true for $\lambda_{2} > 0$ small enough (depending on $C_{0}$, $\lambda_{1}$ and \eqref{j80}).

\begin{figure}
\begin{center}
\begin{picture}(0,0)%
\includegraphics{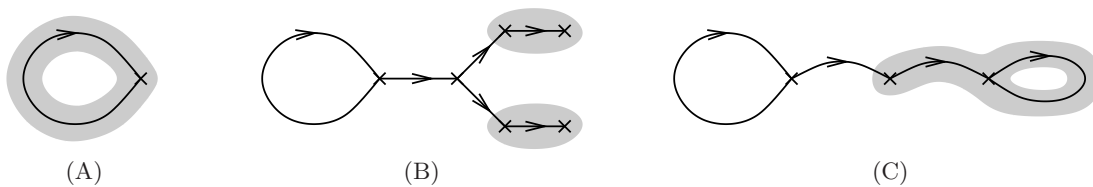}%
\end{picture}%
\setlength{\unitlength}{987sp}%
\begingroup\makeatletter\ifx\SetFigFont\undefined%
\gdef\SetFigFont#1#2#3#4#5{%
  \reset@font\fontsize{#1}{#2pt}%
  \fontfamily{#3}\fontseries{#4}\fontshape{#5}%
  \selectfont}%
\fi\endgroup%
\begin{picture}(27503,4347)(-31575,-6536)
\put(-9374,-6481){\makebox(0,0)[b]{\smash{{\SetFigFont{9}{10.8}{\rmdefault}{\mddefault}{\updefault}(C)}}}}
\put(-21149,-6481){\makebox(0,0)[b]{\smash{{\SetFigFont{9}{10.8}{\rmdefault}{\mddefault}{\updefault}(B)}}}}
\put(-29624,-6481){\makebox(0,0)[b]{\smash{{\SetFigFont{9}{10.8}{\rmdefault}{\mddefault}{\updefault}(A)}}}}
\end{picture}%

\end{center}
\caption{Other geometries with the same pseudo-resonances as in Figure \ref{f59} but with different main localizations of the resonant states (in gray).} \label{f60}
\end{figure}

Moreover, one may change the ``tail'' (i.e. $v_{2}$ and $e_{2}$) and consider other trapped set with the same minimal primitive cycle (i.e. $e_{1}$). For instance, one may consider the geometries illustrated in Figure \ref{f60}. These geometries can be achieved using absorbing potentials. In case {\rm (B)}, we take two symmetric tails. In case {\rm (C)}, we assume that the cycle on the right is not minimal. Since the unique minimal primitive cycle is the same in all these situations, the resonances are the same modulo $o ( h \vert \ln h \vert^{- 1} )$. Nevertheless, the region where the resonant states are of maximal order is totally different (see Figure \ref{f60}). This shows that the asymptotic of resonances and the concentration of the resonant states are not correlated in general. In case {\rm (C)}, we can adjust the geometry so that the homoclinic trajectory on the right is not minimal but that any resonant state is larger here than near $e_{1} \cup \{ ( v_{1} , 0 ) \}$. Thus, the resonances can be generated by one cycle and larger in another cycle.
\end{example}

As for Example \ref{k95}, the conclusions of Remark \ref{k81} do not hold in dimension $n = 1$ without change. Indeed, the imaginary part of the resonances is of order $h \vert \ln h \vert^{- 1}$ in this case, and a barrier-top does not change the order of transversal Lagrangian distributions when the spectral parameter is that close to the real axis (see Corollary \ref{d46}). Nevertheless, if we add a complex valued damping potential of order $h \vert \ln h \vert$ on the homoclinic trajectory $e_{1}$ in Example \ref{k87}, one may artificially produce a delocalization of resonant states in dimension one.

In the setting of Remark \ref{k81}, let $\mu$ be a semiclassical measure associated to a sequence of resonant states normalized in a vicinity of the trapped set (to avoid eventual problems coming from the distortion). The notion of semiclassical measure has been invented by G\'{e}rard \cite{Ge89_01} and Tartar \cite{Ta90_01}. We send back the reader to G\'{e}rard \cite[Section 3]{Ge91_01} for the definition and basic properties. In particular, $\mu$ has mass one and \eqref{k86} implies that $\mu = 0$ near $K_{1}$. This shows that the support of the semiclassical measures is totally disjoint from the trajectories creating the resonances. Thus, it is in general irrelevant to consider semiclassical measures associated to resonant states. This also explains why we work with functions in the general reduction of Section \ref{s3} and not with semiclassical measures like Burq in \cite{Bu02_01}.

Coming back to the setting of Theorem \ref{j35}, one can not distinguish at main order on $\CH$ the resonant states associated to two close resonances on the same accumulation curve. More precisely, we have the following result.

\begin{corollary}\sl \label{j43}
Under the assumptions of Theorem \ref{j35}, let $k_{0} \in \{ 1 , \ldots , K \}$ and $\tau = \tau ( h ) \in [- C , C ]$ for some family of $h$ such that $\mu_{k_{0}} ( \tau , h )$ is an isolated simple eigenvalue of $\widehat{\CQ} ( \tau , h )$ in the sense of Remark \ref{j39} $ii)$. Let $z_{1} , z_{2} \in \eqref{d90}$ with $\re z_{\bullet} = E_{0} + \tau h + o ( h )$ be two resonances such that
\begin{equation}  \label{l23}
z_{\bullet} = z_{q_{\bullet} , k_{0} } ( \tau ) + o \Big( \frac{h}{\vert \ln h \vert} \Big) ,
\end{equation}
where $z_{q , k } ( \tau )$ is given by Proposition \ref{d9}. Then, for any resonant states $u_{1} , u_{2}$ associated to $z_{1} , z_{2}$ normalized as in Theorem \ref{j35}, there exists $\alpha ( h ) \in \C$ with $\vert \alpha \vert = 1$ such that
\begin{equation} \label{j42}
u_{1} = \alpha u_{2} + o ( 1 ) \text{ microlocally near each point of } \CH ,
\end{equation}
as $h$ goes to $0$.
\end{corollary}

This result may seem to be strange since the knowledge of the resonant state $u$ determines the resonance $z$ by the formula $z = P u / u$. Moreover, it is precisely the leading term of $u$ near $\CH$ that is used to compute the asymptotic of the resonances in the proof of Theorem \ref{d8} (see Section \ref{s12}). That resonant states associated to different resonances have the same asymptotic near $\CH$ is a reason why we have not used Grushin problems in this paper. Such phenomenon does not occur for punctual wells in the island situation (see \cite[Section 10]{HeSj86_01}) and also for barrier-top resonances (see \cite[Theorem 4.1]{BoFuRaZe11_01}).

\begin{example}\rm
We come back to Section \ref{s19} {\rm (C)} (see also Figure \ref{f6}). This geometric setting corresponds to a well in the island situation. Thus, for analytic potentials, Helffer and Sj\"{o}strand \cite[Th\'eor\`eme 9.9]{HeSj86_01} have shown that the normalized resonant states are exponentially close to the (sum of) eigenvectors of the operator $P$ restricted to a small neighborhood of the well with Dirichlet boundary condition. In particular, the resonant states are exponentially small outside any neighborhood of the well and almost orthogonal to each other (at least when the corresponding resonances are far away). Thus, two resonant states $u_{1} , u_{2}$ as in Corollary \ref{j43} with $q_{1} \neq q_{2}$ satisfy
\begin{equation} \label{j44}
\< u_{1} , u_{2} \>_{L^{2} ( \R^{n} )} = \CO ( e^{- \varepsilon /h} ) ,
\end{equation}
for some $\varepsilon > 0$ (see Section \ref{s40}). Hence, the resonant states are almost collinear near $\CH$ from \eqref{j42} and almost orthogonal from \eqref{j44}. That $u_{1}$ and $u_{2}$ are different near the barrier-top explains this paradox.

Moreover, from Colin de Verdi{\`e}re and Parisse \cite[Section 5]{CoPa94_01}, the eigenvectors are mainly localized at the top of the barrier, as the resonant states in our case. More precisely, their results imply that
\begin{equation*}
\Vert u \Vert_{L^{2} ( [ - \varepsilon , \varepsilon ] )} \gtrsim \sqrt{\vert \ln h \vert} ,
\end{equation*}
for any $\varepsilon > 0$ and any resonant state $u$ normalized as in Theorem \ref{j35}. It is not clear that this one dimensional behavior still holds true in higher dimensions.
\end{example}

Until now, we have only considered resonant states associated to resonances in \eqref{d90}. But, one can also study the situation in deeper zone (see e.g. Section \ref{s75}). Here, we focus on the resonances of the second kind described in Section \ref{s83}. Recall that the quantization operator $\CQ^{2}$ has been defined in \eqref{g91}.

\begin{proposition}[Resonant states of the second kind]\sl \label{l28}
Assume \ref{h1}--\ref{h4}, \ref{h8}, \ref{h14} and fix $C , \delta > 0$. Let $v = v (h)$ be a family of normalized resonant states associated to a resonance $z = z ( h ) \in \eqref{g87}$. Then, we can find $h^{M} \leq c ( h ) \leq h^{- M}$ for some $M > 0$, such that $u = c v$ has the following properties.

$i)$ The microsupport of $u$ is contained in $\{ ( 0 , 0 ) \} \cup \Lambda_{+}$.

$ii)$ The function $u$ is in $\CI ( \Lambda_{+} , 1 )$ microlocally near any point of $\CH$. We then write
\begin{equation*}
u ( x , h ) = e^{- i A_{k} / h} e^{i \frac{z - E_{0}}{h} t_{-}^{k}} \frac{\CM_{k}^{-}}{\CD_{k} ( t_{-}^{k} )} a_{-}^{k} ( x , h ) e^{i \varphi_{+}^{1} (x) / h} ,
\end{equation*}
microlocally near $\rho_{-}^{k}$, for some $a_{-}^{k} \in S ( 1 )$.

$iii)$ Let $\SA_{0} ( h )$ be the $K$-vector of the $a_{-}^{k} ( x_{-}^{k} , h )$. Let also $\SA_{1} ( h )$ be the $(n - 1 ) K$-vector of the $\partial_{y_{b}} a_{-}^{k} ( x_{-}^{k} , h )$ where $\partial_{y_{b}}$ denotes the derivatives in the $n - 1$ directions of $H_{k}$, the hyperplane of Theorem \ref{a32} associated to $\rho_{-}^{k}$. Then, they satisfy $\SA_{0} ( h ) = \CO ( h )$ and
\begin{equation} \label{l29}
\big( h^{S ( z , h ) / \lambda_{1} + 1 / 2} \CQ^{2} ( z , h ) - 1 \big) \SA_{1} ( h ) = o ( 1 ), 
\end{equation}
as $h$ goes to $0$, with the normalization $\Vert \SA_{1} ( h ) \Vert_{\ell^{2}} = 1$.
\end{proposition}

By comparison with Theorem \ref{j35}, the relevant quantities are now the derivatives of order $1$ of $a_{-}^{k}$ in the transversal directions of $\pi_{x} ( \gamma_{k} )$, and not the value of $a_{-}^{k}$ on this curve. The situation is then similar to that of the harmonic oscillator $- h^{2} \Delta + x^{2}$ (see also \cite{HeSj86_01} for a well in an island trapping and \cite{BoFuRaZe11_01} for barrier-top resonances): the eigenvectors associated to the first eigenvalue are elliptic at $0$ whereas the eigenvectors associated to the second eigenvalue vanish at order $1$ at $0$ modulo $o_{h \to 0} ( 1 )$. Nevertheless, Section \ref{s84} shows that this analogy is only formal.

One can be more precise on the description of resonant states of the second kind. Let $\SA_{2}$ be the $( n - 1 )^{2} K$-vector of the $\partial_{b , c}^{2} a_{-}^{k} ( x_{-}^{k} , h )$. Lemma \ref{g95} and \eqref{i2} give $\SA_{2}$ in terms of $\SA_{1}$ modulo lower order quantities. Then, \eqref{g94} provides an expression for $\SA_{0}$ in terms of $\SA_{1}$. In order to show Proposition \ref{l28}, it is enough to follows the proof of Theorem \ref{j35}, replacing Lemma \ref{d41} by Lemma \ref{g89} and using the reduction \eqref{g93}--\eqref{i37}. We omit the details.

The results of this part can be extended to quasimodes. As usual, they are defined by

\begin{definition}[Quasimodes]\sl
Let $u = u ( h ) \in H^{2} ( \R^{n} )$ be a family of functions with $\Vert u \Vert_{L^{2} ( \R^{n} )} = 1$ and $z = z ( h ) \in \CE_{\theta}$. We say that $u$ is a quasimode of $P$ modulo $\CO ( h^{\infty} )$ associated to the quasiresonance $z$ if and only if $( P_{\theta} - z ) u = \CO ( h^{\infty} )$ as $L^{2} ( \R^{n} )$ function.
\end{definition}

Usually, quasimodes play an important role in the proof of the asymptotic of the resonances (see e.g. G{\'e}rard and Sigal \cite{GeSi92_01}, Tang and Zworski \cite{TaZw98_01} or Popov and Vodev \cite{PoVo99_01} in addition to the papers already cited). This is not the case in the present paper since our strategy use ``test functions'' (see Section \ref{s36}). Nevertheless, the asymptotic of the resonances stated in Theorem \ref{d8} automatically provides resonant states (and then exact quasimodes modulo $\CO ( h^{\infty} )$) which are described in Theorem \ref{j35}. Note that Burq and Christianson \cite{BuCh15_01} have constructed quasimodes modulo $\CO ( h^{2 - \varepsilon} )$ for unidimensional spectral pencils with two homoclinic trajectories. In general, the quasimodes verify the following properties.

\begin{proposition}\sl \label{j45}
Assume \ref{h1}--\ref{h4}, \ref{h8} and fix $C , \delta > 0$. Let $v$ be a quasimode modulo $\CO ( h^{\infty} )$ associated to a quasiresonance $z \in \eqref{d90}$. Then,

$i)$ The quasiresonance $z$ is close to the resonances of $P$. More precisely,
\begin{equation*}
\dist \big( z , \res (P) \big) = o \Big( \frac {h}{\vert \ln h \vert} \Big) .
\end{equation*}

$ii)$ The quasimode $v$ satisfies the conclusions of Theorem \ref{j35}.
\end{proposition}

In dimension $n \geq 2$, this result and Theorem \ref{a1} imply that the quasimodes can not be compactly supported uniformly with respect to $h$. Indeed, adapting \cite[Theorem 1.7]{BoPe13_01}, any quasimode $u$ associated to a quasiresonance $z$ with $\im z \lesssim - h$ satisfies
\begin{equation*}
\Vert \chi u \Vert \leq C \big\Vert \one_{R \leq \vert x \vert \leq R + 1} u \big\Vert + \CO ( h^{\infty} ) ,
\end{equation*}
for $\chi \in C^{\infty}_{0} ( \R^{n} )$ and $R$ large enough. If $u$ is compactly supported, the previous formula becomes $\Vert u \Vert = \CO ( h^{\infty} )$ which gives the contradiction with $\Vert u \Vert = 1$.

On the other hand, under the additional assumptions of Remark \ref{j39} $ii)$, the quasimode $u = c v$ normalized as in Theorem \ref{j35} coincides modulo $o ( 1 )$ with a resonant state microlocally near any point of $\CH$.

\section{General reduction} \label{s3}

\Subsection{The contradiction argument} \label{s31}

We explain here how to obtain a polynomial estimate of the distorted resolvent $( P_{\theta} - z )^{- 1}$ for $z$ in some region $\Omega_{h}$ (see \eqref{a6} for the definition of the distorted operator $P_{\theta}$). The case of the truncated resolvent $\chi ( P - z )^{- 1} \chi$ is reduced to that question from \eqref{c12} with a complex dilation which occurs outside of the support of $\chi \in C^{\infty}_{0} ( \R^{n} )$. We take $\theta = h \vert \ln h \vert$ in the sequel, but we could have chosen $\theta = R h$ with $R > 1$ large enough. For $C > 0$, let $\Omega_{h} \subset B ( E_{0} , C h )$ be a family of compact sets. With quantifiers, we want to prove
\begin{equation} \label{a3}
\exists M > 0 , \quad \exists h_{0} > 0 , \quad \forall h < h_{0} , \quad \forall z \in \Omega_{h} , \qquad \big\Vert ( P_{\theta} - z )^{-1} \big\Vert \leq h^{- M} .
\end{equation}
We shall obtain \eqref{a3} by a contradiction argument. This idea comes from Burq \cite{Bu02_01}. If \eqref{a3} did not hold, there would exist $u = u (h) \in L^{2} ( \R^{n} )$ and $z = z (h) \in \Omega_{h}$ such that
\begin{equation} \label{a4}
\left\{ \begin{aligned}
&( P_{\theta} - z ) u = \CO ( h^{\infty} ) ,    \\
&\Vert u \Vert_{L^{2} ( \R^{n} )} = 1 .
\end{aligned} \right.
\end{equation}
To be precise, \eqref{a4} would occur only for a sequence $( h_{j} )_{j \in \N}$ of positive numbers which converges to $0$. But since this point plays no role in our study, we shall forget to mention the sequence $( h_{j} )_{j}$ in the following and we will do as if \eqref{a4} holds for all $h > 0$ small enough.

In what follows, we suppose that $u \in L^2(\R^n)$ satisfies \eqref{a4}, and we will prove that $\Vert u \Vert = \CO ( h^{\infty} )$, a contradiction. In the rest of this section, we prove that it is sufficient to know that $u$ vanishes microlocally near the trapped set. This reduction has nothing to do with the nature of the trapped set, and relies only on the assumption \ref{h1}. We refer to Appendix \ref{s9} for the notations and some results of microlocal analysis we use.

\Subsection{Energy localization}

We show that $u$ is concentrated in the energy surface $p^{-1} ( E_{0} )$. The main point is a standard elliptic estimate which is a consequence of the ellipticity of $\re (P_{\theta} -z ) \approx p -z$ outside $p^{-1} ( E_{0} )$. More precisely,

\begin{lemma}\sl \label{a12}
For all $\varphi \in C^{\infty}_{0} ( \R )$ with $\varphi = 1$ near $E_{0}$, we have
\begin{equation*}
u = \Op ( \varphi (p) ) u + \CO_{H^{2}} ( h^{\infty} ) .
\end{equation*}
\end{lemma}

\begin{proof}
Let $\varphi_{1} , \varphi_{2} , \ldots$ in $C^{\infty}_{0}( \R )$ with $\one_{\{ E_{0} \}} \prec \cdots \prec \varphi_{2} \prec \varphi_{1} \prec \varphi$. The distorted operator $P_{\theta}$ defined in \eqref{a6} is a differential operator of order $2$ whose symbol $p_{\theta} \in S ( \< \xi \>^{2} )$ satisfies
\begin{align}
p_{\theta} (x, \xi ) &= p \big( x + i \theta F (x) , ( 1 + i \theta {}^{t} ( d F ( x ) ) )^{-1} \xi \big) + S \big( \theta h \< \xi \>^{2} \big)   \nonumber  \\
&= p ( x , \xi ) - i \theta H_{p} ( F(x) \cdot \xi ) + S \big( \theta^{2} \< \xi \>^{2} \big) .  \label{a13}
\end{align}
for $h \lesssim \theta \lesssim 1$. In particular, we deduce $(p_{\theta} -z )^{-1} ( 1 - \varphi ( p) ) \in S ( \< \xi \>^{-2} )$. Then, the pseudodifferential calculus implies
\begin{align*}
\Op \big( (p_{\theta} -z )^{-1} ( 1 - & \varphi (p) ) \big) (P_{\theta} - z)    \\
&= \Op ( 1 - \varphi (p) ) + \Psi \big( h \< \xi \>^{-1} \big) \Op ( 1 - \varphi_{1} (p) ) + \Psi \big( h^{\infty} \< \xi \>^{- \infty} \big) .
\end{align*}
Combining this estimate with \eqref{a4}, we obtain
\begin{equation*}
\Op ( 1 - \varphi (p) ) u = \Psi \big( h \< \xi \>^{-1} \big) \Op ( 1 - \varphi_{1} (p) ) u + \CO_{H^{2}} ( h^{\infty} ).
\end{equation*}
Iterating the previous argument, we get
\begin{equation*}
\Op ( 1 - \varphi (p) ) u = \Psi \big( h^{N} \< \xi \>^{- N} \big) \Op ( 1 - \varphi_{N} (p) ) u + \CO_{H^{2}} ( h^{\infty} ) = \CO_{H^{2}} ( h^{N} ) ,
\end{equation*}
and the lemma follows.
\end{proof}

\Subsection{Spatial localization}

Now, we show that it is enough to estimate $u$ in a compact subset of the position space. This follows from the fact that imaginary part of $P_{\theta} - z$ is elliptic at infinity thanks to the complex distortion.

\begin{lemma}\sl \label{a15}
For $A > 0$ large enough, we have
\begin{equation*}
\Vert u \Vert \lesssim \big\Vert \one_{\{ \vert x \vert \leq A \}} u \big\Vert + \CO ( h^{\infty} ) .
\end{equation*}
\end{lemma}

\begin{proof}
Let $\varphi \in C^{\infty}_{0} ( ] 0 , + \infty [)$ be such that $\varphi = 1$ near $E_{0}$. Combining \eqref{a4} and Lemma \ref{a12}, we obtain
\begin{equation*}
\left\{ \begin{aligned}
&( P_{\theta} - z ) \Op ( \varphi (p) ) u = \CO ( h^{\infty} ) ,   \\
&\Vert \Op ( \varphi (p) ) u \Vert \lesssim 1 .
\end{aligned} \right.
\end{equation*}
Then, using \eqref{a13}, $\vert \im z \vert \lesssim h$ and the pseudodifferential calculus, we get
\begin{align}
\CO ( h^{\infty} ) &= - \im \big\< ( P_{\theta} - z ) \Op ( \varphi (p) ) u , \Op ( \varphi (p) ) u \big\>    \nonumber  \\
&= \big\< \Op \big( \theta \{ p , F (x) \cdot \xi \} + \im z \big) \Op ( \varphi (p) ) u , \Op ( \varphi (p) ) u \big\> + \CO ( \theta^{2} ) \Vert u \Vert^{2}    \nonumber  \\
&= \theta \big\< \Op \big( \{ p , F (x) \cdot \xi \} \varphi^{2} (p) \big) u , u \big\> + \CO ( h ) \Vert u \Vert^{2} . \label{a19}
\end{align}

Since $F (x) = x$ for $x$ large enough, we have, on $p^{-1} ( \supp ( \varphi ) )$,
\begin{equation} \label{a18}
\{ p , F (x) \cdot \xi \} \geq
\left\{ \begin{aligned}
& c > 0 &&\text{ for } \vert x \vert \geq A / 2 ,  \\
& - C &&\text{ for } \vert x \vert \leq A / 2 ,
\end{aligned} \right.
\end{equation}
for $A > 0$ large enough. Let $\chi \in C^{\infty}_{0} ( \R^{n} ; [0 , 1] )$ be such that $\one_{\{ \vert x \vert \leq A /2 \}} \prec \chi \prec \one_{\{ \vert x \vert \leq A \}}$. Thus, \eqref{a18} becomes
\begin{equation*}
\{ p , F (x) \cdot \xi \} \varphi^{2} (p) \geq c (1 - \chi )^{2} \varphi^{2} (p) - C \chi^{2} \varphi^{2} (p) .
\end{equation*}
Now G{\aa}rding's inequality implies
\begin{align*}
\Op \big( \{ p , F (x) \cdot \xi \} \varphi^{2} (p) \big) &\geq c \Op \big( (1 - \chi )^{2} \varphi^{2} (p) \big) - C \Op \big( \chi^{2} \varphi^{2} (p) \big) - \CO (h)    \\
&\geq c \Op \big( (1 - \chi ) \varphi (p) \big)^{2} - C \Op \big( \chi \varphi (p) \big)^{2} - \CO (h) .
\end{align*}
Then, \eqref{a19} together with Lemma \ref{a12} gives
\begin{align*}
\CO ( h^{\infty} ) &\geq c \theta \big\Vert \Op \big( (1 - \chi ) \varphi (p) \big) u \big\Vert^{2} - C \theta \big\Vert \Op \big( \chi \varphi (p) \big) u \big\Vert^{2} + \CO (h) \Vert u \Vert^{2}    \\
&= c \theta \big\Vert(1 - \chi ) \Op ( \varphi (p) ) u \big\Vert^{2} - C \theta \big\Vert \chi \Op ( \varphi (p) ) u \big\Vert^{2} + \CO (h) \Vert u \Vert^{2}   \\
&= c \theta \Vert(1 - \chi ) u \Vert^{2} - C \theta \Vert \chi u \Vert^{2} + \CO (h) \Vert u \Vert^{2} + \CO ( h^{\infty} )    \\
&\geq c \theta /2 \Vert(1 - \chi ) u \Vert^{2} - 2 C \theta \Vert \chi u \Vert^{2} + \CO ( h^{\infty} ) ,
\end{align*}
since $\theta = h \vert \ln h \vert \gg h$. This yields
\begin{equation*}
\Vert u \Vert^{2} \lesssim \Vert \chi u \Vert^{2} + \Vert (1 - \chi ) u \Vert^{2} \lesssim \Vert \chi u \Vert^{2} + \CO ( h^{\infty} ) \leq \big\Vert \one_{\{ \vert x \vert \leq A \}} u \big\Vert^{2} + \CO ( h^{\infty} ) ,
\end{equation*}
and the lemma follows.
\end{proof}

\Subsection{Vanishing in the incoming region} \label{s33}

We now show that $u = 0$ microlocally in the incoming region. This is in agreement with the intuition that the quasimode $u$ must be outgoing. Such a property will be used in the following as an ``initial condition'' at infinity. For $R > 0$, $\varepsilon > 0$ and $\sigma \in [-1 , 1]$, we set
\begin{align*}
\Gamma_{\pm} ( R , \varepsilon , \sigma ) = \Big\{ ( x , \xi ) \in T^{*} \R^{n} ; \ \vert x \vert > R , \ p ( x & , \xi ) \in [E_{0} - \varepsilon , E_{0} + \varepsilon ]   \\
&\text{ and } \pm \cos ( x, \xi ) = \pm \frac{x \cdot \xi}{\vert x \vert \vert \xi \vert} > \pm \sigma \Big\} .
\end{align*}

\begin{lemma}\sl \label{a20}
Let $\varepsilon > 0$, $\sigma < 0$ and then $R > 0$ be large enough. Let $w \in S (1)$ be such that $\supp (w) \subset \Gamma_{-} ( R , \varepsilon , \sigma )$. Then,
\begin{equation*}
\Op (w) u = \CO ( h^{\infty} ) .
\end{equation*}
\end{lemma}

\begin{proof}
This follows from the proof of Theorem 2 of \cite{BoMi04_01}. The unique difference is that the left hand side of \cite[(3.4)]{BoMi04_01} must be replaced by $\CO ( h^{\infty} )$, which does not affect the rest of the proof. As a matter of fact, only resonant states (i.e. function $v$ such that $(P_{\theta} - z )v = 0$) were considered in \cite{BoMi04_01}.
\end{proof}

\Subsection{Propagation of singularities}

In this section, we prove that usual propagation of singularities holds for $P_{\theta} -z$. This result is not completely standard since $P_{\theta}$ is not self-adjoint. It rests crucially on the fact that the imaginary part of $P_{\theta}$ is of size $\theta \lesssim h \vert \ln h \vert$.

\begin{lemma}\sl \label{a21}
For $\rho_{0} \in p^{-1} ( E_{0} )$ and $T \in \R$, denote $\rho_{T} = \exp ( T H_{p} ) ( \rho_{0} )$. If $u =0$ microlocally near $\rho_{0}$, then $u = 0$ microlocally near $\rho_{T}$.
\end{lemma}

\begin{proof}
To show this result, one can follow the standard proof of the propagation of singularities (see e.g. Martinez \cite{Ma02_02}) and check that each step works. Instead, we will reduce the problem $(P_{\theta} -z) u =0$ to an equation of the type $(P - E_{0} ) w = 0$ and then apply the standard propagation of singularities to the operator $P - E_{0}$ which is self-adjoint.

We can assume that $\rho_{0} \neq \rho_{T}$ and then that $p$ is of principal type (i.e. $\nabla p \neq 0$) on the Hamiltonian curve $t \mapsto \exp ( t H_{p} ) ( \rho_{0} )$. Let $\chi \in C^{\infty}_{0} ( T^{*} \R^{n} )$ with $\chi =1$ in a neighborhood of $\{ \exp ( t H_{p} ) ( \rho_{0} ) ; \ t \in [0,T] \}$. For $z \in \Omega_{h}$, we can write $z = E_{0} + h \sigma$ where $\sigma$ is in a compact set of $\C$. From \eqref{a13}, the symbol $p_{\theta} -z \in S ( \< \xi \>^{2} )$ of $P_{\theta} -z$ has the asymptotic in $S ( \< \xi \>^{2} )$
\begin{equation*}
p_{\theta} - z \simeq ( p - E_{0} ) + \sum_{j + k \geq 1} p_{j ,k} ( x , \xi ) \theta^{j} h^{k} ,
\end{equation*}
where the coefficients depend smoothly on $\sigma$ (in fact, $p_{0,1} (x,\xi) = - \sigma$ is the unique coefficient which depend on $\sigma$).

First we construct a weight $G$ such that $e^{G} (P- E_{0}) e^{- G} \approx (P_{\theta} -z)$. Consider $G = \Op (g)$ whose symbol $g \in S ( \theta h^{-1} )$, supported in a neighborhood of $\supp \chi$, satisfies
\begin{equation} \label{a27}
g ( x , \xi , \theta , h) \simeq \sum_{j + k \geq 1} g_{j , k} (x, \xi ) \theta^{j} h^{k-1} ,
\end{equation}
In particular, since $\theta = h \vert \ln h \vert$, we have $g \in S ( \vert \ln h \vert )$. From \cite[Section 4]{DeSjZw04_01}, the operator
\begin{equation} \label{a23}
Q : = e^{G} (P - E_{0}) e^{- G} \simeq \sum_{\ell \geq 0} \frac{1}{\ell !} \ad_{G}^{\ell} (P - E_{0}) ,
\end{equation}
is a pseudodifferential operator of symbol $q \in S ( \< \xi \>^{2} )$ satisfying
\begin{equation*}
q ( x , \xi , \theta , h) \simeq (p - E_{0}) + \sum_{j + k \geq 1} q_{j , k} (x, \xi ) \theta^{j} h^{k} ,
\end{equation*}
with, for all $j+k \geq 1$,
\begin{equation} \label{a24}
q_{j , k} = i \{ p , g_{j,k} \} + r_{j,k} ,
\end{equation}
where $r_{j,k}$ depends only on $p$ and on the $g_{\alpha , \beta}$'s with $\alpha \leq j$, $\beta \leq k$ and $( \alpha , \beta ) \neq ( j , k)$. On the other hand, since $p$ is of principal type on the Hamiltonian curve $t \mapsto \exp ( t H_{p} ) ( \rho_{0} )$, for all $f \in C^{\infty} ( T^{*} \R^{n} )$, one can find a $g \in C^{\infty} ( T^{*} \R^{n} )$ supported in a fixed neighborhood of $\{ \exp ( t H_{p} ) ( \rho_{0} ) ; \ t \in [0,T] \}$ such that
\begin{equation*}
\{ p , g \} = f ,
\end{equation*}
in a fixed neighborhood of $\{ \exp ( t H_{p} ) ( \rho_{0} ) ; \ t \in [0,T] \}$. In fact, it is enough to integrate $f$ along the Hamiltonian curves of $p$ (see \cite[Appendix]{GeSj87_01} for more details). This way, after a possible shrinking of the support of $\chi$, one can construct by induction $g_{j,k}$ supported in a fixed neighborhood of $\supp \chi$ such that
\begin{equation*}
i \{ p , g_{j,k} \} + r_{j,k} = p_{j,k} ,
\end{equation*}
near $\supp \chi$. Using Borel's lemma, we can then find a symbol $g \in S ( \vert \ln h \vert )$ compactly supported and satisfying \eqref{a27}. With such a $g$, \eqref{a23} gives
\begin{equation} \label{a25}
e^{G} (P - E_{0}) e^{- G} \Op ( \chi ) = ( P_{\theta} - E_{0} ) \Op ( \chi ) + \CO ( h^{\infty} ) .
\end{equation}

Now, define $w = e^{- G} \Op ( \chi ) u$. Since $g \in S ( \vert \ln h \vert )$, we have
\begin{equation} \label{a26}
e^{\pm G} \in \Psi ( h^{- C} ) ,
\end{equation}
for some $C > 0$ (see \cite[Section 4]{DeSjZw04_01}). Then, applying \eqref{a25} to $u$ and using \eqref{a4}, \eqref{a26} and the assumptions of the lemma, we get
\begin{equation*}
\left\{ \begin{aligned}
&( P - E_{0} ) w = 0 &&\text{microlocally near } \{ \exp ( t H_{p} ) ( \rho_{0} ) ; \ t \in [0,T] \} ,   \\
&w = 0 &&\text{microlocally near } \rho_{0} ,  \\
&\Vert w \Vert_{H^{2}} \lesssim h^{- C} &&\text{for some } C >0 .
\end{aligned} \right.
\end{equation*}
Then, the standard propagation of singularities for real-valued symbols of principal type implies that $w = 0$ microlocally near $\rho_{T}$. Combining with \eqref{a26}, this yields that $u = 0$ microlocally near $\rho_{T}$.
\end{proof}

\Subsection{Microlocalization near the trapped set}

Combining the previous microlocalizations, we show that, to estimate globally $u$, it is enough to control $u$ microlocally near the trapped set. We begin with a lemma about the classical trajectories which is probably standard.

\begin{lemma}\sl \label{a28}
Let $\rho \in p^{-1} (E_{0})$ and $\rho (t) =\exp ( t H_{p} ) ( \rho )$ be the associated Hamiltonian trajectory. The following alternative holds:

$i)$ either $\rho (t) \to \infty$ as $t \to - \infty$,

$ii)$ or $\dist ( \rho (t) , K (E_{0} ) ) \to 0$ as $t \to - \infty$.
\end{lemma}

\begin{proof}
Assume that $ii)$ does not hold. Then, there exit $\varepsilon > 0$ and a sequence of negative real numbers $(T_{j})_{j \in \N}$ such that $T_{j} < T_{j -1} - 1$ and
\begin{equation*}
\forall j \in \N , \qquad \dist ( \rho ( T_{j} ) , K (E_{0} ) ) > 2 \varepsilon .
\end{equation*}
Since $H_{p}$ is bounded on $p^{-1} ( E_{0} )$, $\rho^{\prime} (t)$ is bounded. In particular, there exits $0 < \delta < 1/2$ such, that for all $j \in \N$ and $t \in [ T_{j} - \delta , T_{j} + \delta ]$, we have
\begin{equation} \label{a29}
\dist ( \rho ( t ) , K (E_{0} ) ) > \varepsilon .
\end{equation}

On the other hand, \cite[Appendix]{GeSj87_01} provides an escape function $G \in C^{\infty} ( T^{*} \R^{n} )$ satisfying
\begin{equation*}
H_p ( G ) \geq \left\{ \begin{aligned}
& 0 &&\text{ on } p^{-1} ( E_{0} ) ,  \\
& 1 &&\text{ on } p^{-1} ( E_{0} ) \setminus ( K (E_{0}) + B (0 , \varepsilon ) ) .
\end{aligned} \right.
\end{equation*}
Therefore, since $\partial_{t} G (\rho (t) ) = H_{p} ( G ) ( \rho (t))$ and $\rho (t) \in p^{-1} (E_{0})$, \eqref{a29} implies that $t \mapsto G ( \rho (t))$ is a non-decreasing function and
\begin{equation*}
G ( \rho ( T_{j} - \delta) ) \leq G ( \rho ( T_{j} + \delta )) - 2 \delta .
\end{equation*}
Then,
\begin{equation*}
G ( \rho ( T_{j} - \delta) ) \leq G ( \rho ( T_{0} )) - 2 j \delta .
\end{equation*}
This implies that $G ( \rho (t)) \to - \infty$ as $t \to - \infty$. Since $G$ is continuous, we must have $\rho (t) \to \infty$ as $t \to - \infty$.
\end{proof}

\begin{proposition}\sl \label{a16}
If $u = 0$ microlocally near each point of $K ( E_{0} )$, then $\Vert u \Vert = \CO ( h^{\infty} )$.
\end{proposition}

\begin{proof}
Let $A > 0$ be given by Lemma \ref{a15} and let $\rho \in p^{-1} ( E_{0} ) \cap \{ \vert x \vert \leq A \}$. As before, we denote $\rho (t) =\exp ( t H_{p} ) ( \rho )$. According to Lemma \ref{a28}, we have two different cases that we treat separately.

Assume first $\rho (t) \to \infty$ as $t \to - \infty$. Then, by a classical argument (see \cite[Appendix]{GeSj87_01}), for $T > 0$ large enough, $\rho ( - T ) \in \Gamma_{-} ( R , \varepsilon , \sigma )$ with $R, \varepsilon , \sigma$ as in the assumptions of Lemma \ref{a20}. Therefore, $u = 0$ microlocally near $\rho (- T )$. Now, the propagation of singularities result given in Lemma \ref{a21} shows that $u =0$ microlocally near $\rho$.

Assume now $\dist ( \rho (t) , K (E_{0} ) ) \to 0$ as $t \to - \infty$. Since $K ( E_{0} )$ is compact, the hypothesis of the proposition yields that $u = 0$ in a neighborhood $U$ of $K ( E_{0} )$, and $\rho (- T) \in U$ for any $T > 0$ large enough. Thus, $u =0$ microlocally near $\rho ( - T)$. Again, $u =0$ microlocally near $\rho$ by propagation of singularities (Lemma \ref{a21}).

By compactness, we have obtained that $u = 0$ microlocally near $p^{-1} ( E_{0} ) \cap \{ \vert x \vert \leq A \}$. The energy localization result stated in Lemma \ref{a12} implies that $\Vert \one_{\{ \vert x \vert \leq A \}} u \Vert = \CO ( h^{\infty} )$. Eventually, the proposition follows from the spatial localization result (Lemma \ref{a15}).
\end{proof}

\begin{remark}\sl \label{a22}
It follows from the contradiction argument explained in Section \ref{s31} and Proposition \ref{a16} that, if $P$ is non-trapping at energy $E_{0}$ (i.e. $K ( E_{0} ) = \emptyset$), then $P$ has no resonance in $\Omega_{h} : = B ( E_{0} , C h )$ for all $C > 0$, and
\begin{equation*}
\big\Vert \chi ( P -z )^{-1} \chi \big\Vert \lesssim h^{- M} ,
\end{equation*}
uniformly for $h$ small enough and $z \in \Omega_{h}$. Note that, with a simpler approach, Martinez \cite{Ma02_01} has proved this result with $\Omega_{h}$ replaced by $[ E_{0} - \varepsilon , E_{0} + \varepsilon] + i [ - C h \vert \ln h \vert , h ]$ for some $\varepsilon > 0$ and all $C > 0$. Working more carefully, one can probably obtain this result in these larger domains.
\end{remark}

\section{Proof of Theorem \ref{a1}} \label{s4}

\begin{figure}
\begin{center}
\begin{picture}(0,0)%
\includegraphics{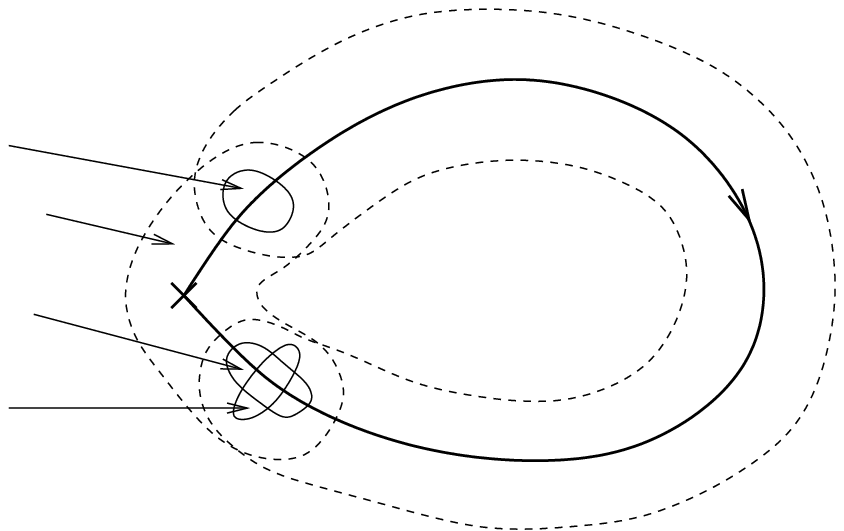}%
\end{picture}%
\setlength{\unitlength}{1579sp}%
\begingroup\makeatletter\ifx\SetFigFont\undefined%
\gdef\SetFigFont#1#2#3#4#5{%
  \reset@font\fontsize{#1}{#2pt}%
  \fontfamily{#3}\fontseries{#4}\fontshape{#5}%
  \selectfont}%
\fi\endgroup%
\begin{picture}(12801,6283)(-1364,-6786)
\put(6976,-961){\makebox(0,0)[lb]{\smash{{\SetFigFont{9}{10.8}{\rmdefault}{\mddefault}{\updefault}$\Omega_{\text{reg}}$}}}}
\put(10726,-4111){\makebox(0,0)[lb]{\smash{{\SetFigFont{9}{10.8}{\rmdefault}{\mddefault}{\updefault}$\CH$}}}}
\put(1051,-3136){\makebox(0,0)[lb]{\smash{{\SetFigFont{9}{10.8}{\rmdefault}{\mddefault}{\updefault}$\Omega_{\text{sing}}$}}}}
\put(1201,-4261){\makebox(0,0)[lb]{\smash{{\SetFigFont{9}{10.8}{\rmdefault}{\mddefault}{\updefault}$V_{+}^{1}$}}}}
\put(-1349,-5386){\makebox(0,0)[lb]{\smash{{\SetFigFont{9}{10.8}{\rmdefault}{\mddefault}{\updefault}$u=u_{-}$ near $V_{-}^{0}$}}}}
\put(-1349,-2236){\makebox(0,0)[lb]{\smash{{\SetFigFont{9}{10.8}{\rmdefault}{\mddefault}{\updefault}$u=u_{+}$ near $V_{+}^{0}$}}}}
\put(3901,-4036){\makebox(0,0)[lb]{\smash{{\SetFigFont{9}{10.8}{\rmdefault}{\mddefault}{\updefault}$0$}}}}
\end{picture}%
\end{center}
\caption{The general setting in the proof of Theorem \ref{a1}.} \label{f5}
\end{figure}

From Proposition \ref{a16}, the theorem will follow if we show that $u$ satisfying \eqref{a4} with $\Omega_{h} = \eqref{m74}$ vanishes microlocally near each point of $K ( E_{0} ) = \{ ( 0 , 0 ) \} \cup \CH$. We will obtain this property by a bootstrap argument using propagation of singularities along the homoclinic set. We can assume that $P_{\theta} = P$ near the base space projection of $\CH$. For $\rho \in T^{*} \R^{n}$ and $t \in \R$, we will use in the sequel the notation $\rho (t) = \exp ( t H_{p} ) ( \rho )$.

\begin{lemma}\sl \label{a33}
We have $u = 0$ microlocally near each point of $\Lambda_{-} \setminus ( \CH \cup \{ ( 0 ,0 ) \} )$.
\end{lemma}

\begin{proof}
Let $\rho \in \Lambda_{-} \setminus ( \CH \cup \{ ( 0 ,0 ) \} )$. Since $\rho \in \Lambda_{-}$, $\rho (t)$ is trapped as $t \to + \infty$. Therefore, $\rho (t)$ is not trapped as $t \to - \infty$ because $\rho \notin K ( E_{0} )$. Thus, $\rho (t) \to \infty$ as $t \to - \infty$. Then, from Lemma \ref{a20}, $u = 0$ microlocally near $\rho ( - T )$ for some $T > 0$ large enough. Eventually, by Lemma \ref{a21}, we get that $u =0$ microlocally near $\rho$.
\end{proof}

\begin{figure}
\begin{center}
\begin{picture}(0,0)%
\includegraphics{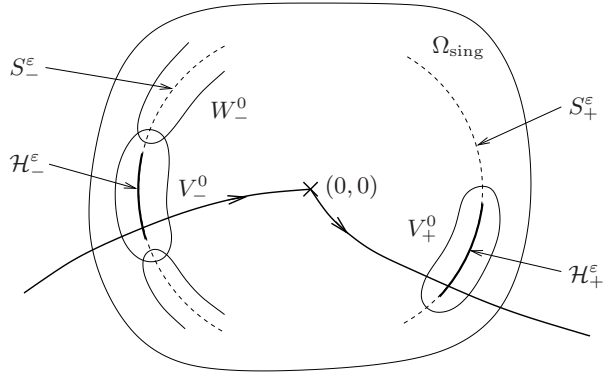}%
\end{picture}%
\setlength{\unitlength}{1184sp}%
\begingroup\makeatletter\ifx\SetFigFont\undefined%
\gdef\SetFigFont#1#2#3#4#5{%
  \reset@font\fontsize{#1}{#2pt}%
  \fontfamily{#3}\fontseries{#4}\fontshape{#5}%
  \selectfont}%
\fi\endgroup%
\begin{picture}(12198,7751)(-914,-7790)
\put(10801,-5911){\makebox(0,0)[lb]{\smash{{\SetFigFont{9}{10.8}{\rmdefault}{\mddefault}{\updefault}$\CH_{+}^{\varepsilon}$}}}}
\put(2626,-4111){\makebox(0,0)[lb]{\smash{{\SetFigFont{9}{10.8}{\rmdefault}{\mddefault}{\updefault}$V_{-}^{0}$}}}}
\put(3301,-2461){\makebox(0,0)[lb]{\smash{{\SetFigFont{9}{10.8}{\rmdefault}{\mddefault}{\updefault}$W_{-}^{0}$}}}}
\put(5701,-4036){\makebox(0,0)[lb]{\smash{{\SetFigFont{9}{10.8}{\rmdefault}{\mddefault}{\updefault}$( 0 , 0 )$}}}}
\put(-899,-3586){\makebox(0,0)[lb]{\smash{{\SetFigFont{9}{10.8}{\rmdefault}{\mddefault}{\updefault}$\CH_{-}^{\varepsilon}$}}}}
\put(-899,-1561){\makebox(0,0)[lb]{\smash{{\SetFigFont{9}{10.8}{\rmdefault}{\mddefault}{\updefault}$S_{-}^{\varepsilon}$}}}}
\put(7426,-4936){\makebox(0,0)[lb]{\smash{{\SetFigFont{9}{10.8}{\rmdefault}{\mddefault}{\updefault}$V_{+}^{0}$}}}}
\put(7951,-1111){\makebox(0,0)[lb]{\smash{{\SetFigFont{9}{10.8}{\rmdefault}{\mddefault}{\updefault}$\Omega_{\text{sing}}$}}}}
\put(10801,-2386){\makebox(0,0)[lb]{\smash{{\SetFigFont{9}{10.8}{\rmdefault}{\mddefault}{\updefault}$S_{+}^{\varepsilon}$}}}}
\end{picture}%
\end{center}
\caption{The geometry near $( 0 , 0 )$.} \label{f3}
\end{figure}

Let $\Omega_{\rm{sing}}$ be a small neighborhood of $(0,0)$ and $\varepsilon > 0$ small enough be given by Theorem \ref{a32} with $\Omega = \Omega_{\rm{sing}}$. Denote $S_{\pm}^{\varepsilon} = \{ ( x, \xi ) \in \Lambda_{\pm}^{0} ; \ \vert x \vert = \varepsilon \}$ and $\CH_{\pm}^{\varepsilon} = S_{\pm}^{\varepsilon} \cap \CH$ (see Figure \ref{f3}). Let also $V_{-}^{0}$ (resp. $W_{-}^{0}$) be a small compact neighborhood of $\CH_{-}^{\varepsilon}$ (resp. $S_{-}^{\varepsilon} \setminus V_{-}^{0}$ with $W_{-}^{0} \cap \CH_{-}^{\varepsilon} = \emptyset$). We define $u_{-}$ as the microlocal restriction of $u$ to a neighborhood of $V_{-}^{0}$. In particular, \eqref{a4} implies
\begin{equation} \label{a34}
\Vert u_{-} \Vert \lesssim 1 .
\end{equation}
Using \eqref{a4} and Lemma \ref{a33}, the setting is the following near $S_{\pm}^{\varepsilon}$:
\begin{equation} \label{a41}
\left\{ \begin{aligned}
&(P -z) u = 0 &&\text{microlocally near } \Omega_{\rm{sing}} ,   \\
&u = u_{-} &&\text{microlocally near } V_{-}^{0} ,   \\
&u = 0 &&\text{microlocally near } W_{-}^{0} ,
\end{aligned} \right.
\end{equation}
and $\Vert u \Vert \lesssim 1$. From \ref{h4}, we know that $g_{-} ( \rho_{-} ) \neq 0$ and $g_{-} ( \rho_{-} ) \cdot g_{+} ( \rho_{+} ) \neq 0$ for all $( \rho_{-} , \rho_{+} ) \in \CH_{-}^{\varepsilon} \times \CH_{+}^{\varepsilon}$. Then, Theorem \ref{a32} gives the existence of a compact neighborhood $V_{+}^{0}$ of $\CH_{+}^{\varepsilon}$ such that
\begin{equation*}
u = \CJ u_{-} \text{ microlocally near } V_{+}^{0} .
\end{equation*}

As above, we define $u_{+}$ as the microlocal restriction of $u$ to a neighborhood of $V_{+}^{0}$. The previous equation and \eqref{m76} say that
\begin{equation} \label{a35}
u_{+} (x) = h^{\sum \frac{\lambda_{j} - \lambda_{1}}{2 \lambda_{1}} - i \frac{z - E_{0}}{\lambda_{1} h}} \int_{\R^{n}} e^{i ( \varphi_{+} ( x ) - \varphi_{-} ( y ) ) / h} \widetilde{d} ( x, y , z , h ) u_{-} ( y ) \, d y ,
\end{equation}
microlocally near $V_{+}^{0}$. Notice that the integrand is microlocalized inside $V_{-}^{0}$ with respect to $y$. Combining this equation with \eqref{a34} and the fact that $- \delta h \leq \im z \leq 0$ for $z \in \eqref{m74}$, we get, microlocally near $V_{+}^{0}$,
\begin{equation} \label{a36}
u_{+} (x) = h^{\sum \frac{\lambda_{j} - \lambda_{1}}{2 \lambda_{1}} - \frac{\delta}{\lambda_{1}}} a (x , h ) e^{i \varphi_{+} (x) / h} \in \CI ( \Lambda_{+}^{0} , h^{- N} ) ,
\end{equation}
for some $a \in S ( 1 )$ and $N \in \R$. See Section \ref{s35} for the definition of $\CI ( \cdot , \cdot )$, the set of semiclassical Lagrangian distributions.

\begin{figure}
\begin{center}
\begin{picture}(0,0)%
\includegraphics{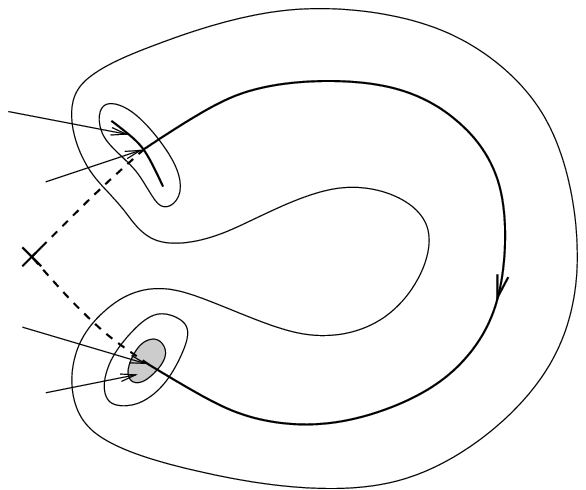}%
\end{picture}%
\setlength{\unitlength}{1184sp}%
\begingroup\makeatletter\ifx\SetFigFont\undefined%
\gdef\SetFigFont#1#2#3#4#5{%
  \reset@font\fontsize{#1}{#2pt}%
  \fontfamily{#3}\fontseries{#4}\fontshape{#5}%
  \selectfont}%
\fi\endgroup%
\begin{picture}(11465,7723)(-1514,-6715)
\put(3751,-1636){\makebox(0,0)[lb]{\smash{{\SetFigFont{9}{10.8}{\rmdefault}{\mddefault}{\updefault}$V_{+}^{0}$}}}}
\put(3826,-4636){\makebox(0,0)[lb]{\smash{{\SetFigFont{9}{10.8}{\rmdefault}{\mddefault}{\updefault}$V_{+}^{1}$}}}}
\put(-149,-3061){\makebox(0,0)[lb]{\smash{{\SetFigFont{9}{10.8}{\rmdefault}{\mddefault}{\updefault}$( 0 , 0 )$}}}}
\put(6601,-1486){\makebox(0,0)[lb]{\smash{{\SetFigFont{9}{10.8}{\rmdefault}{\mddefault}{\updefault}$\Omega_{\text{reg}}$}}}}
\put(  1,-811){\makebox(0,0)[lb]{\smash{{\SetFigFont{9}{10.8}{\rmdefault}{\mddefault}{\updefault}$\CH_{+}^{\varepsilon}$}}}}
\put(301,-4186){\makebox(0,0)[lb]{\smash{{\SetFigFont{9}{10.8}{\rmdefault}{\mddefault}{\updefault}$\rho_{-}$}}}}
\put(-1274,-1936){\makebox(0,0)[lb]{\smash{{\SetFigFont{9}{10.8}{\rmdefault}{\mddefault}{\updefault}$\rho_{-} ( - t ( \rho_{-} ) )$}}}}
\put(-1499,-5311){\makebox(0,0)[lb]{\smash{{\SetFigFont{9}{10.8}{\rmdefault}{\mddefault}{\updefault}$V_{-}^{0} \cap \Lambda_{-}^{0} \cap \CH$}}}}
\end{picture} $\qquad \qquad$ \begin{picture}(0,0)%
\includegraphics{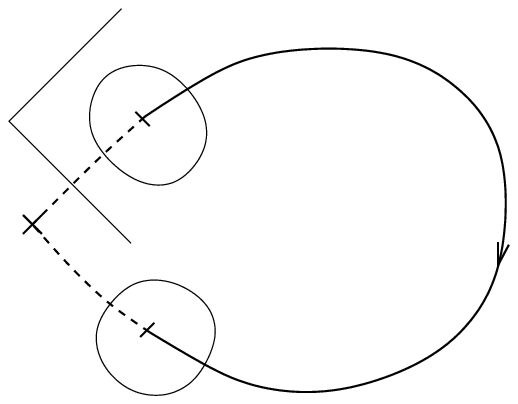}%
\end{picture}%
\setlength{\unitlength}{1184sp}%
\begingroup\makeatletter\ifx\SetFigFont\undefined%
\gdef\SetFigFont#1#2#3#4#5{%
  \reset@font\fontsize{#1}{#2pt}%
  \fontfamily{#3}\fontseries{#4}\fontshape{#5}%
  \selectfont}%
\fi\endgroup%
\begin{picture}(9000,7095)(-164,-6609)
\put(-149,-3061){\makebox(0,0)[lb]{\smash{{\SetFigFont{9}{10.8}{\rmdefault}{\mddefault}{\updefault}$( 0 , 0 )$}}}}
\put(2401,-1036){\makebox(0,0)[lb]{\smash{{\SetFigFont{9}{10.8}{\rmdefault}{\mddefault}{\updefault}$\rho_{+}$}}}}
\put(2926,-2011){\makebox(0,0)[lb]{\smash{{\SetFigFont{9}{10.8}{\rmdefault}{\mddefault}{\updefault}$U$}}}}
\put(3001,-4411){\makebox(0,0)[lb]{\smash{{\SetFigFont{9}{10.8}{\rmdefault}{\mddefault}{\updefault}$\rho_{-}$}}}}
\put(2776,-5386){\makebox(0,0)[lb]{\smash{{\SetFigFont{9}{10.8}{\rmdefault}{\mddefault}{\updefault}$\Lambda_{+}^{1}$}}}}
\put(1126,-1486){\makebox(0,0)[lb]{\smash{{\SetFigFont{9}{10.8}{\rmdefault}{\mddefault}{\updefault}$\Lambda_{+}^{0}$}}}}
\end{picture}%
\end{center}
\caption{The geometric setting near the homoclinic set and $\Lambda_{+}^{1}$.} \label{f4}
\end{figure}

We will now use the standard propagation of singularities along the homoclinic curves $\CH$ to express $u_{-}$ in terms of $u_{+}$. For all $\rho_{-} \in V_{-}^{0} \cap \Lambda_{-}^{0} \cap \CH$, there exists a unique time $t ( \rho_{-} ) > 0$ such that $\rho_{-} ( - t ( \rho_{-} ) ) \in \CH_{+}^{\varepsilon}$. Then, using that $V_{-}^{0} \cap \Lambda_{-}^{0} \cap \CH$ is a compact set and that $\exp ( t ( \rho_{-} ) H_{p} ) ( V_{+}^{0} )$ is a neighborhood of $\rho_{-}$, there exists a finite number of points $\rho_{-}^{j} \in V_{-}^{0} \cap \Lambda_{-}^{0} \cap \CH$, $j = 1 , \ldots , J$, such that
\begin{equation*}
V_{+}^{1} : = \bigcup_{j = 1}^{J} \exp ( t ( \rho_{-}^{j} ) H_{p} ) ( V_{+}^{0} ) ,
\end{equation*}
is a compact neighborhood of $V_{-}^{0} \cap \Lambda_{-}^{0} \cap \CH$. Let also $\Omega_{\text{reg}} \subset T^{*} \R^{n} \setminus \{ ( 0 , 0 ) \}$ be a compact neighborhood of the set
\begin{equation*}
\Big\{ \exp ( t H_{p} ) ( \rho ) ; \ \rho \in V_{+}^{0} \text{ and } 0 \leq t \leq \max_{1 \leq j \leq J} t ( \rho_{-}^{j} ) \Big\} ,
\end{equation*}
(see Figure \ref{f4}). In this region of the homoclinic set, we will use the standard propagation of singularities result since $P$ is of principal type and $\im z = \CO (h)$. Indeed, we know from \eqref{a4} and \eqref{a36} that
\begin{equation} \label{a37}
\left\{ \begin{aligned}
&( P -z ) u = 0 &&\text{microlocally near } \Omega_{\text{reg}} ,  \\
&u = u_{+} \in \CI ( \Lambda_{+}^{0} , h^{- N} ) &&\text{microlocally near } V_{+}^{0} .
\end{aligned} \right.
\end{equation}
By usual propagation of Lagrangian distributions (see \cite{MaFe81_01}), this implies that
\begin{equation} \label{a38}
u = u_{-} \in \CI ( \Lambda_{+}^{1} , h^{- N} ) \text{ microlocally near } V_{+}^{1} ,
\end{equation}
where $\Lambda_{+}^{1}$ is the evolution of $\Lambda_{+}^{0}$ by the Hamiltonian flow when it reaches $V_{+}^{1}$. More precisely, let $\rho_{-} \in V_{+}^{1}$ and recall that $\rho_{+} = \rho_{-} ( - t ( \rho_{-} ) ) \in \CH_{+}^{\varepsilon}$. Then, locally near $\rho_{-}$, $\Lambda_{+}^{1}$ is defined as $\exp ( t ( \rho_{-} ) H_{p} ) ( U )$ where $U$ is a neighborhood of $\rho_{+}$ in $\Lambda_{+}^{0}$ (see Figure \ref{f4}). Note that $V_{+}^{1}$ only depends on $V_{+}^{0}$ and on the classical flow, and may be different from $V_{-}^{0}$.

Now, we insert the information \eqref{a38} into \eqref{a35} and use the assumptions \ref{h5}. But \eqref{a38} does not describe $u_{-}$ in the whole $V_{-}^{0}$. However, the restriction of $u_{-}$ to $V_{-}^{0} \setminus V_{+}^{1}$ does not contribute to $u_{+}$ microlocally near $V_{+}^{0}$. Indeed, $u_{-} =0$ microlocally near each point of $V_{-}^{0} \cap \Lambda_{-}^{0} \setminus V_{+}^{1}$ from Lemma \ref{a33}. Moreover, the values of $u_{-}$ outside $\Lambda_{-}^{0}$ do not contribute to \eqref{a35} by a non stationary phase argument (since $\overline{\widetilde{d} (x , y , z , h ) e^{- i \varphi_{-} ( y ) /h}}$ is microlocalized in $\Lambda_{-}^{0}$ with respect to $y$). In the anisotropic case \ref{h5a}, \eqref{a35} immediately gives, microlocally near $V_{+}^{0}$,
\begin{equation} \label{a39}
u_{+} \in \CI ( \Lambda_{+}^{0} , h^{- N + \alpha} ) \qquad \text{with} \qquad \alpha = \sum_{j=1}^{n} \frac{\lambda_{j} - \lambda_{1}}{2 \lambda_{1}} - \frac{\delta}{\lambda_{1}} > 0 ,
\end{equation}
for $\delta > 0$ small enough. In the transversal case \ref{h5b}, at each point of $\Lambda_{-}^{0} \cap \Lambda_{+}^{1}$, these manifolds intersect transversally along at least one direction. Therefore, using \eqref{a38}, performing a stationary phase expansion in this direction and integrating with respect to the other variables, \eqref{a35} gives, microlocally near $V_{+}^{0}$,
\begin{equation} \label{a40}
u_{+} \in \CI ( \Lambda_{+}^{0} , h^{- N + \alpha} ) \qquad \text{with} \qquad \alpha = \frac{1}{2} - \frac{\delta}{\lambda_{1}} > 0 ,
\end{equation}
for $\delta > 0$ small enough.

Iterating \eqref{a37}--\eqref{a40}, we get
\begin{equation*}
u_{-} \in \CI ( \Lambda_{+}^{1} , h^{- N + \alpha k} ) \text{ microlocally near } V_{+}^{1} ,
\end{equation*}
for all $k \in \N$. This gives $u_{-} = 0$ microlocally near $V_{+}^{1}$ and \eqref{a41} becomes
\begin{equation*}
\left\{ \begin{aligned}
&(P -z) u = 0 &&\text{microlocally near } \Omega_{\rm{sing}} ,   \\
&u = 0 &&\text{microlocally near } S_{-}^{\varepsilon} .
\end{aligned} \right.
\end{equation*}
Then, the uniqueness part in Theorem \ref{a32} yields $u = 0$ microlocally near $(0,0)$. Eventually, since all the trajectories in $\CH$ enter any neighborhood of $(0,0)$,
\begin{equation} \label{a42}
u = 0 \text{ microlocally near } K ( E_{0} ) ,
\end{equation}
by standard propagation of singularities. This, together with Section \ref{s31} and Proposition \ref{a16}, finishes the proof of Theorem \ref{a1}.

\section{Proof of Theorem \ref{a2}} \label{s5}

\Subsection{Reduction to the tangential part} \label{s51}

The proof of Theorem \ref{a2} is similar to that of Theorem \ref{a1}. The difference is that, at each ``turn'' (around the fixed point and the homoclinic set), the considered function $u$ only decays by a small coefficient (typically $\Vert \CT_{0} \Vert < 1$) and not a positive power of $h$ as before. This is due to the fact that the trapping is stronger in the present situation.

\begin{figure}
\begin{center}
\begin{picture}(0,0)%
\includegraphics{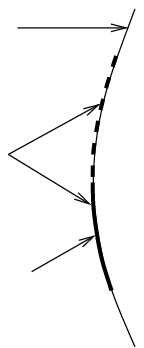}%
\end{picture}%
\setlength{\unitlength}{1184sp}%
\begingroup\makeatletter\ifx\SetFigFont\undefined%
\gdef\SetFigFont#1#2#3#4#5{%
  \reset@font\fontsize{#1}{#2pt}%
  \fontfamily{#3}\fontseries{#4}\fontshape{#5}%
  \selectfont}%
\fi\endgroup%
\begin{picture}(3187,5444)(-764,-6983)
\put(-749,-5986){\makebox(0,0)[lb]{\smash{{\SetFigFont{9}{10.8}{\rmdefault}{\mddefault}{\updefault}$\CH_{{\rm tang} , +}^{\varepsilon}$}}}}
\put(-224,-2011){\makebox(0,0)[lb]{\smash{{\SetFigFont{9}{10.8}{\rmdefault}{\mddefault}{\updefault}$S_{+}^{\varepsilon}$}}}}
\put(-449,-4036){\makebox(0,0)[lb]{\smash{{\SetFigFont{9}{10.8}{\rmdefault}{\mddefault}{\updefault}$\CH_{+}^{\varepsilon}$}}}}
\end{picture} $\qquad \qquad \qquad$ \begin{picture}(0,0)%
\includegraphics{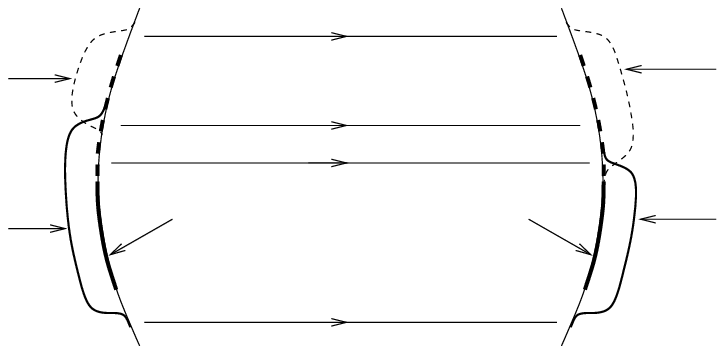}%
\end{picture}%
\setlength{\unitlength}{1184sp}%
\begingroup\makeatletter\ifx\SetFigFont\undefined%
\gdef\SetFigFont#1#2#3#4#5{%
  \reset@font\fontsize{#1}{#2pt}%
  \fontfamily{#3}\fontseries{#4}\fontshape{#5}%
  \selectfont}%
\fi\endgroup%
\begin{picture}(12630,5444)(-1214,-7058)
\put(8101,-5161){\makebox(0,0)[rb]{\smash{{\SetFigFont{9}{10.8}{\rmdefault}{\mddefault}{\updefault}$\CH_{{\rm tang} , -}^{\varepsilon}$}}}}
\put(5401,-2911){\makebox(0,0)[b]{\smash{{\SetFigFont{9}{10.8}{\rmdefault}{\mddefault}{\updefault}$\text{the Hamiltonian flow}$}}}}
\put(-1199,-2761){\makebox(0,0)[lb]{\smash{{\SetFigFont{9}{10.8}{\rmdefault}{\mddefault}{\updefault}$u_{+}^{\rm trans}$}}}}
\put(-1199,-5161){\makebox(0,0)[lb]{\smash{{\SetFigFont{9}{10.8}{\rmdefault}{\mddefault}{\updefault}$u_{+}^{\rm tang}$}}}}
\put(11401,-5161){\makebox(0,0)[lb]{\smash{{\SetFigFont{9}{10.8}{\rmdefault}{\mddefault}{\updefault}$u_{-}^{\rm tang}$}}}}
\put(11401,-2761){\makebox(0,0)[lb]{\smash{{\SetFigFont{9}{10.8}{\rmdefault}{\mddefault}{\updefault}$u_{-}^{\rm trans}$}}}}
\put(2701,-5161){\makebox(0,0)[lb]{\smash{{\SetFigFont{9}{10.8}{\rmdefault}{\mddefault}{\updefault}$\CH_{{\rm tang} , +}^{\varepsilon}$}}}}
\end{picture}%
\end{center}
\caption{The subsets of $S_{\pm}^{\varepsilon}$ and the choice of the supports of $u_{\pm}^{\bullet}$.} \label{f7}
\end{figure}

As in the previous section, to prove Theorem \ref{a2}, it is enough to show that a function $u$ satisfying \eqref{a4} with $\Omega_{h} = \eqref{b68}$ vanishes microlocally near each point of $K ( E_{0} ) = \{ ( 0 , 0 ) \} \cup \CH$ (see Proposition \ref{a16}). Note also that Lemma \ref{a33} is still valid in the present setting. As in Section \ref{s4}, let $\varepsilon > 0$ be small enough, $S_{\pm}^{\varepsilon} = \{ ( x, \xi ) \in \Lambda_{\pm}^{0} ; \ \vert x \vert = \varepsilon \}$ and $\CH_{\pm}^{\varepsilon} = S_{\pm}^{\varepsilon} \cap \CH$. We also define $\CH_{{\rm tang} , \pm}^{\varepsilon} = S_{\pm}^{\varepsilon} \cap \CH_{\rm tang}$ (see Figure \ref{f7}).

We now decompose $u_{\pm}$ into its tangential and transversal trajectories contributions. As above, we define $u_{\pm}$ as the microlocal restriction of $u$ to a neighborhood of $\CH_{\pm}^{\varepsilon}$. Let $V_{\pm}^{\varepsilon} \subset \varepsilon \S^{n-1}$ be small neighborhoods of $\pi_{x} ( \CH_{{\rm tang} , \pm}^{\varepsilon} )$ and let $\chi_{\pm}^{\varepsilon} \in C^{\infty} ( \varepsilon \S^{n-1} ; [ 0 , 1 ] )$ be such that
\begin{equation} \label{j15}
\one_{\pi_{x} ( \CH_{{\rm tang} , \pm}^{\varepsilon} )} \prec \chi_{\pm}^{\varepsilon} \prec \one_{V_{\pm}^{\varepsilon}} .
\end{equation}
Moreover, we can assume that every point of $\CH_{-}^{\varepsilon} \cap V_{-}^{\varepsilon} \times \R^{n}$ comes (for the Hamiltonian flow) from a point of $\CH_{+}^{\varepsilon} \cap \{ \chi_{+}^{\varepsilon} = 1 \} \times \R^{n}$. We then define $u_{\pm}^{\rm tang}$ as the solution of
\begin{equation*}
\left\{ \begin{aligned}
&(P -z) u_{\pm}^{\rm tang} = 0 &&\text{near } \varepsilon \S^{n-1} ,   \\
&u_{\pm}^{\rm tang} = \chi_{\pm}^{\varepsilon} u_{\pm} &&\text{on } \varepsilon \S^{n-1} .
\end{aligned} \right.
\end{equation*}
In particular, $u_{\pm}^{\rm tang}$ is a microlocal restriction of $u$ in a small neighborhood of $\CH_{{\rm tang} , \pm}^{\varepsilon}$. We also define the function $u_{\pm}^{\rm trans} = u_{\pm} - u_{\pm}^{\rm tang}$. So, $u = u_{\pm}^{\rm tang} +  u_{\pm}^{\rm trans}$ microlocally near $\CH_{\pm}^{\varepsilon}$. The setting is illustrated in Figure \ref{f7}. From \eqref{a4}, we get
\begin{equation} \label{a95}
\big\Vert u_{\pm}^{\bullet} \big\Vert \lesssim 1 .
\end{equation}

As in \eqref{a41}, we have
\begin{equation} \label{a96}
\left\{ \begin{aligned}
&(P -z) u = 0 &&\text{microlocally near } ( 0 , 0 ) ,   \\
&u = u_{-}^{\rm tang} + u_{-}^{\rm trans} &&\text{microlocally near } \CH_{-}^{\varepsilon} ,   \\
&u = 0 &&\text{microlocally near } S_{-}^{\varepsilon} \setminus \CH_{-}^{\varepsilon} ,
\end{aligned} \right.
\end{equation}
and $\Vert u \Vert \lesssim 1$. From \ref{h4}, we know that $g_{-} ( \rho_{-} ) \neq 0$ and $g_{-} ( \rho_{-} ) \cdot g_{+} ( \rho_{+} ) \neq 0$ for all $( \rho_{-} , \rho_{+} ) \in \CH_{-}^{\varepsilon} \times \CH_{+}^{\varepsilon}$. Then, Theorem \ref{a32} gives
\begin{align}
u_{+}^{\rm tang} &= \CJ_{{\rm tang} \leftarrow {\rm tang}} u_{-}^{\rm tang} + \CJ_{{\rm tang} \leftarrow {\rm trans}} u_{-}^{\rm trans}  ,   \label{a84} \\
u_{+}^{\rm trans} &= \CJ_{{\rm trans} \leftarrow {\rm tang}} u_{-}^{\rm tang} + \CJ_{{\rm trans} \leftarrow {\rm trans}} u_{-}^{\rm trans} ,    \label{a85}
\end{align}
where the operators $\CJ_{\bullet \leftarrow \star}$ are the restrictions of \eqref{m76} to the microsupport of the corresponding functions:
\begin{equation} \label{a98}
\CJ_{\bullet \leftarrow \star} u_{-}^{\star} (x) = h^{- i \frac{z - E_{0}}{h \lambda}} \int_{\R^{n}} e^{i ( \varphi_{+} (x) - \varphi_{-} (y) )/h} \widetilde{d}_{\bullet \leftarrow \star} ( x, y , z , h ) u_{-}^{\star} ( y ) \, d y ,
\end{equation}
where $\widetilde{d}_{\bullet \leftarrow \star}$ is the restriction of $\widetilde{d}$ near the support of $u_{+}^{\bullet}$ in $x$ and near the support of $u_{-}^{\star}$ in $y$. Note that since $\vert \im z \vert \lesssim h \vert \ln h \vert^{- 1}$, we have
\begin{equation} \label{a99}
\Big\vert h^{- i \frac{z - E_{0}}{h \lambda}} \Big\vert \lesssim 1 .
\end{equation}
Combining \eqref{a95}, \eqref{a98} and \eqref{a99}, we get
\begin{equation} \label{a92}
u_{+}^{\bullet} \in \CI ( \Lambda_{+}^{0} , h^{- N_{+}^{\bullet}} ) ,
\end{equation}
for some $N_{+}^{\bullet} \in \R$.

Now, as in the previous section, we will use the standard propagation of singularities along the homoclinic curves $\CH$ to estimate $u_{-}^{\bullet}$. From \eqref{a4} and \eqref{a92}, we get
\begin{equation} \label{b1}
\left\{ \begin{aligned}
&( P -z ) u = 0 &&\text{microlocally near } \CH_{\rm tang} ,  \\
&u = u_{+}^{\rm tang} \in \CI ( \Lambda_{+}^{0} , h^{- N_{+}^{{\rm tang}}} ) &&\text{microlocally near } \CH_{{\rm tang} , +}^{\varepsilon} ,
\end{aligned} \right.
\end{equation}
with $\Vert u \Vert \lesssim 1$ and $\im z = \CO (h)$. Since the functions $u_{\pm}^{\rm tang}$ have been defined according to the discussion before \eqref{a95} (see also Figure \ref{f7}), the usual propagation of Lagrangian distributions gives
\begin{equation} \label{b2}
u_{-}^{\rm tang} \in \CI ( \Lambda_{+}^{1} , h^{- N_{-}^{\rm tang}} ) ,
\end{equation}
where $\Lambda_{+}^{1}$ is defined below \eqref{a38} and
\begin{equation} \label{b3}
N_{-}^{\rm tang} = N_{+}^{\rm tang} .
\end{equation}
On the other hand, we also have
\begin{equation} \label{a97}
\left\{ \begin{aligned}
&( P -z ) u = 0 &&\text{microlocally near } \CH ,  \\
&u = u_{+}^{\rm tang} + u_{+}^{\rm trans} &&\text{microlocally near } \CH_{+}^{\varepsilon} ,
\end{aligned} \right.
\end{equation}
and $\Vert u \Vert \lesssim 1$. Then, using the definition of $u_{-}^{\rm trans}$ and the proof of \eqref{b2}, we get
\begin{equation} \label{b4}
u_{-}^{\rm trans} \in \CI ( \Lambda_{+}^{1} , h^{- N_{-}^{\rm trans}} ) ,
\end{equation}
where
\begin{equation} \label{b5}
N_{-}^{\rm trans} = \max \big( N_{+}^{\rm tang} , N_{+}^{\rm trans} \big) .
\end{equation}

Note that one can also express the solutions of \eqref{b1} and \eqref{a97} in a more abstract formalism as it as been done in \eqref{a84} and \eqref{a85}. In fact, using an approach similar to the one of Sj\"{o}strand and Zworski \cite{SjZw07_01}, the functions $u_{-}^{\bullet}$ can be written as
\begin{align}
u_{-}^{\rm tang} &= \CM_{{\rm tang} \leftarrow {\rm tang}} u_{+}^{\rm tang} ,     \label{j84} \\
u_{-}^{\rm trans} &= \CM_{{\rm trans} \leftarrow {\rm tang}} u_{+}^{\rm tang} + \CM_{{\rm trans} \leftarrow {\rm trans}} u_{+}^{\rm trans} ,
\end{align}
where $\CM_{\bullet \leftarrow \star}$ are Fourier integral operators whose canonical relation is given by the Hamiltonian flow. Again, the fact that $u_{-}^{\rm tang}$ does not depend on $u_{+}^{\rm trans}$ is due to the choice of the neighborhoods where $u_{\pm}^{\bullet}$ are defined (see Figure \ref{f7}). We will not use this formalism since we will need to compute explicitly the symbol of the Lagrangian distributions in the sequel.

The next result shows that to control all the functions $u_{\pm}^{\bullet}$, it is enough to estimate $N_{-}^{\rm tang}$.

\begin{lemma}\sl \label{a94}
Assume that $u_{-}^{\rm tang} \in \CI ( \Lambda_{+}^{1} , h^{- N} )$ for some constant $N \in \R$. Then, we can take all the constants $N_{\pm}^{\bullet}$ appearing in \eqref{a92}, \eqref{b2} and \eqref{b4} equal to $N$.
\end{lemma}

\begin{proof}
Assume that $u_{-}^{\rm tang} \in \CI ( \Lambda_{+}^{1} , h^{- N} )$. Using the formula \eqref{a98} and the estimate \eqref{a99} which guaranties that the prefactor in \eqref{a98} is of order $\CO ( 1 )$, we get
\begin{equation*}
\CJ_{\bullet \leftarrow {\rm tang}} u_{-}^{\rm tang} \in \CI ( \Lambda_{+}^{0} , h^{- N} ) .
\end{equation*}
On the other hand, $\Lambda_{+}^{1}$ intersects transversally $\Lambda_{-}^{0}$ at each point of the microsupport of $u_{-}^{\rm trans}$. Then, performing a stationary phase expansion in \eqref{a98} as in \eqref{a40} and using also \eqref{a99}, we obtain
\begin{equation} \label{b6}
\CJ_{\bullet \leftarrow {\rm trans}} u_{-}^{\rm trans} \in \CI ( \Lambda_{+}^{0} , h^{- N_{-}^{\rm trans} + \frac{1}{2}} ) .
\end{equation}
From the two last equations, one sees that it is possible to take
\begin{equation} \label{a86}
N_{+}^{\bullet} = \max \Big( N , N_{-}^{\rm trans} - \frac{1}{2} \Big) ,
\end{equation}
in \eqref{a92}.

Furthermore, we have already shown in \eqref{b5} that it is possible to take
\begin{equation} \label{a87}
N_{-}^{\rm trans} = \max \big( N_{+}^{\rm tang} , N_{+}^{\rm trans} \big) .
\end{equation}
Thus, performing a bootstrap argument in \eqref{a86} and \eqref{a87}, one can choose $N_{-}^{\rm trans} = N$ and eventually $N_{+}^{\bullet} = N$.
\end{proof}

\Subsection{Study of the tangential part} \label{s52}

In this part, we will prove that $u_{-}^{\rm tang} = \CO ( h^{\infty} )$. As before, we will show that $u_{-}^{\rm tang}$ is ``smaller'' after a turn around the fixed point and the homoclinic orbits. But, to the contrary of Section \ref{s4}, there is no hope to gain some power of $h$ at each turn since $\CJ_{{\rm tang} \leftarrow {\rm tang}}$ and $\CM_{{\rm tang} \leftarrow {\rm tang}}$ are now of order $0$ on the considered functions. Then, we will have to compute the symbol of the Lagrangian distribution $u_{-}^{\rm tang} \in \CI ( \Lambda_{+}^{1} , h^{- N_{-}^{\rm tang}} )$ (see \eqref{b2}) and prove that it decays after a turn.

We first note that, since $\Lambda_{+}^{1}$ and $\Lambda_{-}^{0}$ have the same tangent space on $\CH_{\rm tang}$, the manifold $\Lambda_{+}^{1}$ projects nicely on the $x$-space near the support of $u_{-}^{\rm tang}$. Thus, we can write
\begin{equation*}
u_{-}^{\rm tang} ( y , h ) = a_{-} ( y , h ) h^{- N_{-}^{\rm tang}} e^{i \varphi^{1}_{+} (y) / h} ,
\end{equation*}
where $\varphi_{+}^{1}$ parametrizes $\Lambda_{+}^{1}$ and $a_{-} \in S (1)$. The same way,
\begin{equation*}
u_{+}^{\rm tang} ( x , h ) = a_{+} ( x , h ) h^{- N_{-}^{\rm tang}} e^{i \varphi_{+} (x) / h} ,
\end{equation*}
for some $a_{+} \in S (1)$.

\begin{figure}
\begin{center}
\begin{picture}(0,0)%
\includegraphics{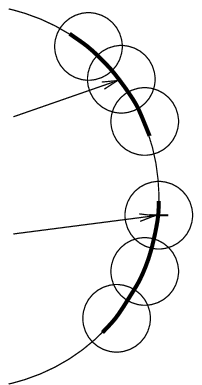}%
\end{picture}%
\setlength{\unitlength}{1184sp}%
\begingroup\makeatletter\ifx\SetFigFont\undefined%
\gdef\SetFigFont#1#2#3#4#5{%
  \reset@font\fontsize{#1}{#2pt}%
  \fontfamily{#3}\fontseries{#4}\fontshape{#5}%
  \selectfont}%
\fi\endgroup%
\begin{picture}(3405,6030)(4486,-8176)
\put(4801,-3961){\makebox(0,0)[rb]{\smash{{\SetFigFont{9}{10.8}{\rmdefault}{\mddefault}{\updefault}$\pi_{x} ( \CH_{{\rm tang} , -}^{\varepsilon} )$}}}}
\put(4501,-7861){\makebox(0,0)[lb]{\smash{{\SetFigFont{9}{10.8}{\rmdefault}{\mddefault}{\updefault}$\varepsilon \S^{n-1}$}}}}
\put(7876,-5611){\makebox(0,0)[lb]{\smash{{\SetFigFont{9}{10.8}{\rmdefault}{\mddefault}{\updefault}$B ( y_{j} , \nu )$}}}}
\put(4801,-5836){\makebox(0,0)[rb]{\smash{{\SetFigFont{9}{10.8}{\rmdefault}{\mddefault}{\updefault}$y_{j}$}}}}
\end{picture} $\qquad \qquad \qquad \qquad$ \begin{picture}(0,0)%
\includegraphics{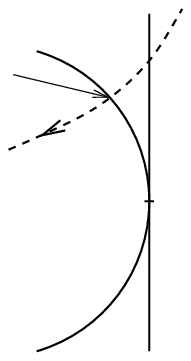}%
\end{picture}%
\setlength{\unitlength}{1184sp}%
\begingroup\makeatletter\ifx\SetFigFont\undefined%
\gdef\SetFigFont#1#2#3#4#5{%
  \reset@font\fontsize{#1}{#2pt}%
  \fontfamily{#3}\fontseries{#4}\fontshape{#5}%
  \selectfont}%
\fi\endgroup%
\begin{picture}(3873,6002)(286,-6662)
\put(3826,-6061){\makebox(0,0)[lb]{\smash{{\SetFigFont{9}{10.8}{\rmdefault}{\mddefault}{\updefault}$H_{j}$}}}}
\put(1051,-6061){\makebox(0,0)[lb]{\smash{{\SetFigFont{9}{10.8}{\rmdefault}{\mddefault}{\updefault}$\varepsilon \S^{n-1}$}}}}
\put(301,-2011){\makebox(0,0)[lb]{\smash{{\SetFigFont{9}{10.8}{\rmdefault}{\mddefault}{\updefault}$Y (y)$}}}}
\put(1726,-3511){\makebox(0,0)[lb]{\smash{{\SetFigFont{9}{10.8}{\rmdefault}{\mddefault}{\updefault}$Y_{y} (t)$}}}}
\put(3826,-1786){\makebox(0,0)[lb]{\smash{{\SetFigFont{9}{10.8}{\rmdefault}{\mddefault}{\updefault}$y$}}}}
\put(3826,-4036){\makebox(0,0)[lb]{\smash{{\SetFigFont{9}{10.8}{\rmdefault}{\mddefault}{\updefault}$y_{j}$}}}}
\end{picture}%
\end{center}
\caption{The cutting of $\pi_{x} ( \CH_{{\rm tang} , -}^{\varepsilon} )$ and the straightening near $y_{j}$.} \label{f8}
\end{figure}

We now estimate $a_{+}$ from $a_{-}$ using the propagation through the fixed point. Since the geometric quantities which govern this transmission appear asymptotically at the fixed point, we will take $\varepsilon > 0$ small enough (see the beginning of Section \ref{s51}). Combining \eqref{a84}, \eqref{b6} and Lemma \ref{a94}, we have
\begin{equation} \label{a91}
u_{+}^{\rm tang} = \CJ_{{\rm tang} \leftarrow {\rm tang}} u_{-}^{\rm tang} + \CI ( \Lambda_{+}^{0} , h^{- N_{-}^{\rm tang} + \frac{1}{2}} ) .
\end{equation}
We will now apply Theorem \ref{a32} to compute $\CJ_{{\rm tang} \leftarrow {\rm tang}} u_{-}^{\rm tang}$. But since we consider $u_{-}^{\rm tang}$ on $\varepsilon \S^{n-1} = \{ \vert y \vert = \varepsilon \}$ whereas the initial conditions in Theorem \ref{a32} live on hyperplanes, we will have to cut $u_{-}^{\rm tang}$ into small pieces and straighten these parts to hypersurfaces. More precisely, let $\nu > 0$ be small enough. There exist $y_{j} \in \pi_{x} ( \CH_{{\rm tang} , -}^{\varepsilon} )$, with $j = 1 , \ldots , 	J$ for some $J = J ( \varepsilon , \nu )$, such that
\begin{equation*}
\pi_{x} \big( \CH_{{\rm tang} , -}^{\varepsilon} \big) \subset \bigcup_{1 \leq j \leq J} B ( y_{j} , \nu ) ,
\end{equation*}
where $\pi_{x} ( y , \eta ) = y$ is the base space projection (see Figure \ref{f8}). Using this partition, one can construct $\varphi_{j} \in C^{\infty}_{0} ( B ( 0 , 2 \nu ) ; [ 0 , 1 ] )$ such that
\begin{equation} \label{b44}
\sum_{1 \leq j \leq J} \varphi_{j} ( y - y_{j} ) = 1 \text{ locally near } \pi_{x} \big( \CH_{{\rm tang} , -}^{\varepsilon} \big) .
\end{equation}
Note that $\varphi_{+}^{1}$ is well-defined and smooth on $\pi_{x} ( \CH_{{\rm tang} , -}^{\varepsilon} ) + B ( 0 , 2 \nu )$ for $\nu > 0$ small enough. For $1 \leq j \leq J$, let
\begin{equation*}
u_{-}^{j} ( y , h ) = a_{-}^{j} ( y , h ) h^{- N_{-}^{\rm tang}} e^{i \varphi^{1}_{+} (y) / h} ,
\end{equation*}
be a WKB solution near $\varepsilon \S^{n-1}$ of
\begin{equation} \label{b7}
\left\{ \begin{aligned}
&( P - z ) u_{-}^{j} = \CO ( h^{\infty} ) ,   \\
&u_{-}^{j} ( y , h ) = \varphi_{j} ( y - y_{j} ) a_{-} ( y , h ) h^{- N_{-}^{\rm tang}} e^{i \varphi_{+}^{1} (y) / h} \text{ for } y \in \varepsilon \S^{n-1} ,
\end{aligned} \right.
\end{equation}
with $a_{-}^{j} \in S (1)$. More precisely, it means that
\begin{equation} \label{b13}
a_{-}^{j} ( y , h ) \simeq \sum_{k \geq 0} a_{- , k}^{j} ( y , h ) h^{k} ,
\end{equation}
where the $a_{- , k}^{j} \in S (1)$'s satisfy the usual transport equations
\begin{align}
2 \nabla \varphi_{+}^{1} \cdot \nabla a_{- , 0}^{j} + \big( \Delta \varphi_{+}^{1} - i z / h \big) a_{- , 0}^{j} &= 0 ,   \label{b8} \\
2 \nabla \varphi_{+}^{1} \cdot \nabla a_{- , k}^{j} + \big( \Delta \varphi_{+}^{1} - i z / h \big) a_{- , k}^{j} &= i \Delta a_{- , k - 1}^{j} \text{ for } k \geq 1 ,  \label{b9}
\end{align}
with initial conditions
\begin{align}
a_{- , 0}^{j} ( y , h ) &= \varphi_{j} ( y - y_{j} ) a_{-} ( y , h ) \text{ for } y \in \varepsilon \S^{n-1} ,   \label{b10}  \\
a_{- , k}^{j} ( y , h ) &= 0 \text{ for } y \in \varepsilon \S^{n-1} \text{ and } k \geq 1 .  \label{b11}
\end{align}
Thus, the functions $u_{-}^{j} \in \CI ( \Lambda_{+}^{1} , h^{- N_{-}^{\rm tang}} )$ verify by construction
\begin{equation} \label{b12}
u_{-}^{\rm tang} = \sum_{1 \leq j \leq J} u_{-}^{j} + \CO ( h ^{\infty} ) ,
\end{equation}
and $\vert a_{-} ( y , h ) \vert = \sum_{j} \vert a_{-}^{j} ( y , h ) \vert$ on $\varepsilon \S^{n-1}$.

Now, we compute $u_{-}^{j}$ on the hypersurface $H_{j} = \{ y \cdot y_{j} = \varepsilon^{2} \}$ which will be the surface of the initial Cauchy data in Theorem \ref{a32} (see Figure \ref{f8}). Using the expansion \eqref{b13}, the transport equation \eqref{b8} and the initial condition \eqref{b10}, we get, for $y \in H_{j}$,
\begin{equation} \label{b14}
a^{j}_{-} ( y , h ) = \varphi_{j} ( Y (y) - y_{j} ) a_{-} ( Y (y) , h ) e^{- \int_{0}^{t (y)} ( \Delta \varphi_{+}^{1} ( Y_{y} (s) ) - i z / h ) d s} + S ( h ) .
\end{equation}
Here, $( Y_{y} (t) , E_{y} (t) ) = \exp ( t H_{p} ) ( y , \nabla \varphi_{+}^{1} ( y ) )$, the time $t (y)$ is the unique (small) time such that $Y_{y} (t) \in \varepsilon \S^{n-1}$ and $Y (y) = Y_{y} ( t (y) )$. Since $y_{j} \in \pi_{x} ( \CH_{{\rm tang} , -}^{\varepsilon} )$, we get that $\nabla \varphi_{+}^{1} ( y_{j} ) = \nabla \varphi_{-} ( y_{j} ) = - \lambda y_{j} / 2 + \CO ( \varepsilon^{2} )$ is not orthogonal to $y_{j}$. Then, applying the implicit function theorem to $( t , y ) \mapsto \vert Y_{y} (t) \vert^{2} - \varepsilon^{2}$, we obtain that $t (y)$ and $Y (y)$ are $C^{\infty}$. In particular,
\begin{equation} \label{b15}
t (y) = \CO_{\varepsilon} ( \nu ) .
\end{equation}
Here, $\CO_{\alpha} ( 1 )$ means a function which is bounded by a constant which may depend on the parameter $\alpha$. On the other hand, since $\Delta \varphi_{+}^{1} ( y_{j} ) = \Delta \varphi_{-} ( y_{j} ) = - n \lambda / 2 + \CO ( \varepsilon )$, we get
\begin{equation} \label{b16}
\Delta \varphi_{+}^{1} ( Y_{y} (s) ) - i z / h = \CO (1) ,
\end{equation}
uniformly for $s \in [ 0 , t (y) ]$. Using \eqref{b15} and \eqref{b16}, \eqref{b14} becomes
\begin{equation} \label{b41}
a^{j}_{-} ( y , h ) = \varphi_{j} ( Y (y) - y_{j} ) a_{-} ( Y (y) , h ) ( 1 + \CO_{\varepsilon} ( \nu ) ) + S ( h )
\end{equation}

From \eqref{a91} and \eqref{b12}, we have
\begin{equation*}
u_{+}^{\rm tang} = \sum_{1 \leq j \leq J} \CJ^{j} u_{-}^{j} + \CI ( \Lambda_{+}^{0} , h^{- N_{-}^{\rm tang} + \frac{1}{2}} ) ,
\end{equation*}
where $\CJ^{j}$ is the operator of Theorem \ref{a32} with initial data on the hypersurface $H_{j}$. Then, modulo $S ( h^{\frac{1}{2}} )$, we have
\begin{equation} \label{b43}
a_{+} ( x , h ) = \sum_{1 \leq j \leq J} a_{+}^{j} ( x , h ) ,
\end{equation}
where
\begin{equation*}
a_{+}^{j} ( x , h ) = h^{- i \frac{z - E_{0}}{h \lambda}} \int_{H_{j}} e^{i (\varphi_{+}^{1} (y) - \varphi_{-} (y) )/h} d^{j} ( x, y , z , h ) a_{-}^{j} ( y , h ) \, d y .
\end{equation*}
From \eqref{a99} and \eqref{b18}, we can replace in the previous equation $d^{j}$ by $d_{0}^{j}$ given in \eqref{b19} modulo $S ( h^{\zeta} )$ for some $0 < \zeta < 1 / 2$. Thus, the previous equation yields
\begin{equation} \label{b42}
\vert a_{+}^{j} ( x , h ) \vert \leq h^{\frac{\im z}{h \lambda}} \int_{H_{j}} \vert d^{j}_{0} ( x, y , z ) \vert \vert a_{-}^{j} ( y , h ) \vert \, d y + \CO_{\varepsilon , \nu} ( h^{\zeta} ) .
\end{equation}
The asymptotic of $d_{0}^{j}$ is given by the following lemma. Its proof, which rests mainly on \eqref{b19} and on properties of the classical dynamic, is postponed to Section \ref{s20}. Equation \eqref{j2} will not be used here but in Section \ref{s27}.

\begin{lemma}\sl \label{b17}
We have
\begin{equation} \label{j2}
d_{0}^{j} ( x , y , z ) = \varepsilon e^{- i n \frac{\pi}{4}} \Big( \frac{\lambda}{2 \pi} \Big)^{\frac{n}{2}} \big( i \lambda x \cdot Y (y) \big)^{- \frac{n}{2} + i \frac{z - E_{0}}{\lambda h}}  \Gamma \Big( \frac{n}{2} - i \frac{z - E_{0}}{\lambda h} \Big)  + \CO ( \varepsilon^{2 - n} ) + \CO_{\varepsilon} ( \nu ) ,
\end{equation}
and
\begin{align}
\vert d_{0}^{j} ( x , y , z ) \vert &= \varepsilon ( 2 \pi )^{- \frac{n}{2}} \big\vert ( x \cdot Y (y) ) \big\vert^{- \frac{n}{2}} e^{- \frac{\pi}{2} \frac{\re z - E_{0} }{\lambda h} \sgn ( x \cdot Y (y) )} \Big\vert \Gamma \Big( \frac{n}{2} - i \frac{\re z - E_{0}}{\lambda h} \Big) \Big\vert   \nonumber \\
&\qquad \qquad \qquad \qquad \qquad \qquad \qquad + \CO ( \varepsilon^{2 - n} ) + \CO_{\varepsilon} ( \nu ) + \CO_{\varepsilon , \nu} \big( \vert \ln h \vert^{-1} \big) ,  \label{j3}
\end{align}
uniformly for $x \in \varepsilon \S^{n-1} \cap \supp a_{+}$ and $y \in H_{j} \cap \supp a_{- , 0}^{j}$.
\end{lemma}

\begin{remark}\sl \label{b40}
A priori the remainder terms in \eqref{b41} and in Lemma \ref{b17} depend on $j$, but following the proof of these results, one can see that they are uniform in $j$.
\end{remark}

Thus, with \eqref{b41} and \eqref{j3} in mind, \eqref{b42} becomes
\begin{align*}
\vert a_{+}^{j} ( x , h ) \vert \leq{}& h^{\frac{\im z}{h \lambda}} \int_{H_{j}} \Big( \varepsilon f_{0} ( x , Y (y) ) + \CO ( \varepsilon^{2 - n} ) + \CO_{\varepsilon} ( \nu ) + \CO_{\varepsilon , \nu} \big( \vert \ln h \vert^{-1} \big) \Big) \\
&\qquad \qquad \qquad \qquad \qquad \qquad \qquad \times \varphi_{j} ( Y (y) - y_{j} ) \vert a_{-} ( Y (y) , h ) \vert \, d y + \CO_{\varepsilon , \nu} ( h^{\zeta} ) ,
\end{align*}
where
\begin{equation} \label{b57}
f_{0}  ( x , Y ) = ( 2 \pi )^{- \frac{n}{2}} \big\vert ( x \cdot Y ) \big\vert^{- \frac{n}{2}} e^{- \frac{\pi}{2} \frac{\re z - E_{0}}{\lambda h} \sgn ( x \cdot Y )} \Big\vert \Gamma \Big( \frac{n}{2} - i \frac{\re z - E_{0}}{\lambda h} \Big) \Big\vert .
\end{equation}
Then, performing the change of variables $H_{j} \ni y \mapsto Y = Y (y) \in \varepsilon \S^{n-1}$,
\begin{align}
\vert a_{+}^{j} ( x , h ) \vert \leq{}& h^{\frac{\im z}{h \lambda}} \int_{\varepsilon \S^{n-1}} \Big( \varepsilon f_{0} ( x , Y ) + \CO ( \varepsilon^{2 - n} ) + \CO_{\varepsilon} ( \nu ) + \CO_{\varepsilon , \nu} \big( \vert \ln h \vert^{-1} \big) \Big) \nonumber \\
&\qquad \qquad \qquad \qquad \qquad \qquad \qquad \times \varphi_{j} ( Y - y_{j} ) \vert a_{-} ( Y , h ) \vert \, d Y + \CO_{\varepsilon , \nu} ( h^{\zeta} ) , \label{n2}
\end{align}
since $d y = ( 1 + \CO_{\varepsilon} ( \nu ) ) d Y$. Summing over $j$, \eqref{b44} and \eqref{b43} imply
\begin{align*}
\vert a_{+} ( x , h ) \vert \leq{}& h^{\frac{\im z}{h \lambda}} \int_{\varepsilon \S^{n-1}} \Big( \varepsilon f_{0} ( x , Y ) + \CO ( \varepsilon^{2 - n} ) + \CO_{\varepsilon} ( \nu ) + \CO_{\varepsilon , \nu} \big( \vert \ln h \vert^{-1} \big) \Big) \vert a_{-} ( Y , h ) \vert \, d Y    \\
& + \CO_{\varepsilon , \nu} ( h^{\zeta} ) ,
\end{align*}
since the remainder terms are uniform with respect to $j$ (see Remark \ref{b40}). Taking $\nu$ small enough depending on $\varepsilon$ in the previous equation, where all the quantities are now on $\varepsilon \S^{n-1}$ and not on small pieces redressed on $H_{j}$, we obtain
\begin{align*}
\vert a_{+} ( x , h ) \vert \leq{}& h^{\frac{\im z}{h \lambda}} \int_{\varepsilon \S^{n-1}} \Big( \varepsilon f_{0} ( x , Y ) + \CO ( \varepsilon^{2 - n} ) + \CO_{\varepsilon} \big( \vert \ln h \vert^{-1} \big) \Big) \vert a_{-} ( Y , h ) \vert \, d Y + \CO_{\varepsilon} ( h^{\zeta} )  \\
={}& h^{\frac{\im z}{h \lambda}} \int_{V_{-}^{\varepsilon}} \varepsilon f_{0} ( x , Y ) \vert a_{-} ( Y , h ) \vert \, d Y  \\
&+ \big( \CO ( \varepsilon ) + \CO_{\varepsilon} \big( \vert \ln h \vert^{-1} \big) \big) \Vert a_{-} ( \cdot , h ) \Vert_{L^{\infty}( V_{-}^{\varepsilon} )} + \CO_{\varepsilon} ( h^{\zeta} ) .
\end{align*}
We now want to replace the integration over $V_{-}^{\varepsilon}$ by an integration over $\CH_{\rm tang}^{- \infty}$ as it appear in the definition of $\CT_{0}$ (see  \eqref{b70}). First, taking $V_{-}^{\varepsilon}$ sufficiently close (depending on $\varepsilon$) to the compact set $\pi_{x} ( \CH_{{\rm tang} , -}^{\varepsilon} )$ (see Figure \ref{f7}) and using \ref{h4} and the regularity of the measure, we obtain
\begin{align}
\vert a_{+} ( x , h ) \vert \leq{}& h^{\frac{\im z}{h \lambda}} \int_{\pi_{x} ( \CH_{{\rm tang} , -}^{\varepsilon} ) / \varepsilon} f_{0} ( \varepsilon^{-1} x , \widetilde{\omega} ) \vert a_{-} ( \varepsilon \widetilde{\omega} , h ) \vert \, d \widetilde{\omega}    \nonumber \\
&+ \big( \CO ( \varepsilon ) + \CO_{\varepsilon} \big( \vert \ln h \vert^{-1} \big) \big) \Vert a_{-} ( \cdot , h ) \Vert_{L^{\infty}( V_{-}^{\varepsilon} )} + \CO_{\varepsilon} ( h^{\zeta} ) .   \label{b50}
\end{align}
Note that in the previous equation, $\pi_{x} ( \CH_{{\rm tang} , -}^{\varepsilon} ) / \varepsilon \subset \S^{n-1}$ and the integrand is of order $\CO (1)$ from \ref{h4} and $x \in V_{+}^{\varepsilon} \subset \varepsilon \S^{n-1}$. We will now use the

\begin{lemma}\sl \label{b46}
We have
\begin{equation*}
\mes_{\S^{n - 1}} \big( \CH_{\rm tang}^{- \infty} \Delta \pi_{x} ( \CH_{{\rm tang} , -}^{\varepsilon} ) / \varepsilon \big) = o_{\varepsilon \to 0} (1) ,
\end{equation*}
where $o_{\varepsilon \to 0} (1)$ notes a function which tends to $0$ as $\varepsilon$ goes to $0$, and $A \Delta B = A \setminus B \cup B \setminus A$ denotes the symmetric difference of two sets $A , B$.
\end{lemma}

\begin{proof}
We first show the following assertion. Let $A \subset \S^{n - 1}$ be a compact set and let $G_{\varepsilon}$ be a family of diffeomorphisms of $\S^{n - 1}$ such that
\begin{equation} \label{j6}
\left\{ \begin{aligned}
&G_{\varepsilon} ( x ) = x + o_{\varepsilon \to 0} ( 1 ) , \\
&\Vert d G_{\varepsilon} ( x ) \Vert \leq C ,
\end{aligned} \right.
\end{equation}
uniformly for $x \in \S^{n - 1}$ and $\varepsilon > 0$ small enough. Then,
\begin{equation} \label{j5}
\mes_{\S^{n - 1}} \big( A \Delta G_{\varepsilon} ( A ) \big) = o_{\varepsilon \to 0} ( 1 ) .
\end{equation}

Indeed, let $\delta > 0$. By the regularity of the Lebesgue measure, there exists an open neighborhood $U \subset \S^{n-1}$ of the compact set $A$ with
\begin{equation} \label{j7}
\mes_{\S^{n - 1}} ( D) \leq \delta ,
\end{equation}
where $B = \S^{n - 1} \setminus U$ and $D = U \setminus A$. Since $A$ and $B$ are disjoint compact sets, there exists $\mu > 0$ such that $A + B ( 0 , \mu ) \subset \S^{n-1} \setminus B = A \cup D$. On the other hand, for $\varepsilon > 0$ small enough and $a \in A$, \eqref{j6} gives $\vert G_{\varepsilon} ( a ) - a \vert < \mu$. Then,
\begin{equation*}
G_{\varepsilon} ( a ) \in A + B ( 0 , \mu ) \subset A \cup D.
\end{equation*}
Thus, for $\varepsilon > 0$ small enough, we have
\begin{equation*}
G_{\varepsilon} ( A ) \subset A \cup D ,
\end{equation*}
and, by the same way,
\begin{equation*}
G_{\varepsilon} ( B ) \subset B \cup D .
\end{equation*}
In particular, $G_{\varepsilon} ( A ) \setminus A \subset D$ and \eqref{j7} implies
\begin{equation} \label{j8}
\mes_{\S^{n - 1}} \big( G_{\varepsilon} ( A ) \setminus A \big) \leq \delta .
\end{equation}
On the other hand, we have
\begin{align*}
A = \S^{n - 1} \setminus ( B \cup D ) \subset \S^{n - 1} \setminus G_{\varepsilon} ( B ) = G_{\varepsilon} ( A ) \cup G_{\varepsilon} ( D ) ,
\end{align*}
since $G_{\varepsilon}$ is a bijection. Combining with \eqref{j6} and \eqref{j7}, we deduce
\begin{align}
\mes_{\S^{n - 1}} \big( A \setminus G_{\varepsilon} ( A ) \big) &\leq \mes_{\S^{n - 1}} ( G_{\varepsilon} ( D ) ) = \int \one_{G_{\varepsilon}^{- 1} ( x ) \in D} d x  \nonumber  \\
&\leq \int \one_{y \in D} \vert d G_{\varepsilon}  ( y ) \vert d y \leq C^{n - 1} \mes_{\S^{n - 1}} ( D ) = C^{n - 1} \delta . \label{j9}
\end{align}
Eventually, \eqref{j5} is a consequence of \eqref{j8} and \eqref{j9}.

Coming back to the proof of Lemma \ref{b46}, we choose $A = \CH_{\rm tang}^{- \infty}$ and $G_{\varepsilon} = F^{-1} ( \varepsilon , \cdot )$ where $F$ is defined and studied in Section \ref{a71}. It is proved in this appendix that $F^{-1} ( \varepsilon , \cdot )$ is a diffeomorphism of $\S^{n - 1}$ (see above Proposition \ref{a65}), that $F^{-1} ( \varepsilon , \CH_{\rm tang}^{- \infty} ) = \pi_{x} ( \CH_{{\rm tang} , -}^{\varepsilon} ) / \varepsilon$ by construction and that \eqref{j6} holds true thanks to \eqref{a68}. Then, Lemma \ref{b46} follows from \eqref{j5}.
\end{proof}

Combining \eqref{b50} with \ref{h4}, \eqref{a99} and Lemma \ref{b46}, we obtain
\begin{align}
\vert a_{+} ( x , h ) \vert \leq{}& h^{\frac{\im z}{h \lambda}} \int_{\CH_{\rm tang}^{- \infty} \cap \pi_{x} ( \CH_{{\rm tang} , -}^{\varepsilon} ) / \varepsilon} f_{0} ( \varepsilon^{-1} x , \widetilde{\omega} ) \vert a_{-} ( \varepsilon \widetilde{\omega} , h ) \vert \, d \widetilde{\omega}    \nonumber \\
&+ \big( o_{\varepsilon \to 0} (1) + \CO_{\varepsilon} \big( \vert \ln h \vert^{-1} \big) \big) \Vert a_{-} ( \cdot , h ) \Vert_{L^{\infty}( V_{-}^{\varepsilon} )} + \CO_{\varepsilon} ( h^{\zeta} ) ,   \label{b55}
\end{align}
for all $x \in V_{+}^{\varepsilon}$. At this point, we can not assume that the integration occurs on the whole $\CH_{\rm tang}^{- \infty}$ since $a_{-} ( \varepsilon \cdot , h )$ is a priori not defined on this set.

\begin{figure}
\begin{center}
\begin{picture}(0,0)%
\includegraphics{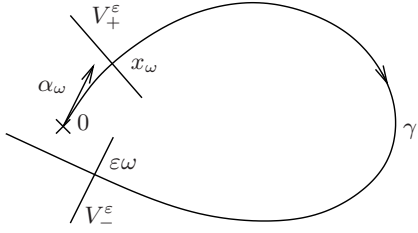}%
\end{picture}%
\setlength{\unitlength}{1184sp}%
\begingroup\makeatletter\ifx\SetFigFont\undefined%
\gdef\SetFigFont#1#2#3#4#5{%
  \reset@font\fontsize{#1}{#2pt}%
  \fontfamily{#3}\fontseries{#4}\fontshape{#5}%
  \selectfont}%
\fi\endgroup%
\begin{picture}(8373,4715)(2368,-6050)
\put(4051,-5986){\makebox(0,0)[lb]{\smash{{\SetFigFont{9}{10.8}{\rmdefault}{\mddefault}{\updefault}$V_{-}^{\varepsilon}$}}}}
\put(5026,-2761){\makebox(0,0)[lb]{\smash{{\SetFigFont{9}{10.8}{\rmdefault}{\mddefault}{\updefault}$x_{\omega}$}}}}
\put(10726,-4111){\makebox(0,0)[lb]{\smash{{\SetFigFont{9}{10.8}{\rmdefault}{\mddefault}{\updefault}$\gamma$}}}}
\put(4201,-1786){\makebox(0,0)[lb]{\smash{{\SetFigFont{9}{10.8}{\rmdefault}{\mddefault}{\updefault}$V_{+}^{\varepsilon}$}}}}
\put(3901,-4036){\makebox(0,0)[lb]{\smash{{\SetFigFont{9}{10.8}{\rmdefault}{\mddefault}{\updefault}$0$}}}}
\put(3076,-3211){\makebox(0,0)[lb]{\smash{{\SetFigFont{9}{10.8}{\rmdefault}{\mddefault}{\updefault}$\alpha_{\omega}$}}}}
\put(4576,-4861){\makebox(0,0)[lb]{\smash{{\SetFigFont{9}{10.8}{\rmdefault}{\mddefault}{\updefault}$\varepsilon \omega$}}}}
\end{picture}%
\end{center}
\caption{The geometric setting in \eqref{b51}.} \label{f10}
\end{figure}

We now estimate $a_{-}$ from $a_{+}$ using the propagation through the homoclinic curves $\CH$. From \eqref{b1}, the usual propagation of the Lagrangian distributions (see \cite{MaFe81_01}) and the choice of the cut-off functions $\chi_{\pm}^{\varepsilon}$ (see Figure \ref{f7}), we get that
\begin{equation*}
u_{-}^{\rm tang} ( \varepsilon \omega , h ) = a_{-} ( \varepsilon \omega , h ) h^{- N_{-}^{\rm tang}} e^{i \varphi^{1}_{+} ( \varepsilon \omega ) / h} ,
\end{equation*}
for $\varepsilon \omega \in V_{-}^{\varepsilon}$, where $a_{-} \in S (1)$ satisfies
\begin{equation} \label{b51}
a_{-} ( \varepsilon \omega , h ) = e^{i ( A_{\omega} - \varphi^{1}_{+} ( \varepsilon \omega ) ) / h} e^{- i \nu_{\omega} \frac{\pi}{2}} e^{i T_{\omega} ( z - E_{0} ) / h} \CM_{\varepsilon} ( \alpha_{\omega} ) \chi_{-}^{\varepsilon} ( \varepsilon \omega ) a_{+} ( x_{\omega} , h ) + S (h) .
\end{equation}
In the previous equation, $x_{\omega}$ is the unique point in $V_{+}^{\varepsilon}$ such that $( x_{\omega} , \nabla \varphi_{+} ( x_{\omega} ) )$ belongs to the Hamiltonian curve $\gamma (t) = ( y (t) , \eta (t) ) = \exp ( t H_{p} ) ( \varepsilon \omega , \nabla \varphi_{+}^{1} ( \varepsilon \omega ) )$ at time $- T_{\omega}$ (see Figure \ref{f10}). The quantity
\begin{equation*}
A_{\omega} = \int_{\gamma ( [ - T_{\omega} , 0 ] )} \eta \cdot d y ,
\end{equation*}
denotes the action along the curve $\gamma ( [ - T_{\omega} , 0 ] )$ and $\nu_{\omega}$ its Maslov's index. Eventually,
\begin{equation} \label{j12}
\alpha_{\omega} = \frac{\Phi ( x_{\omega} )}{\vert \Phi ( x_{\omega} ) \vert} ,
\end{equation}
is the asymptotic direction of the curve $\gamma$ as $t \to - \infty$ and the Maslov determinant $\CM_{\varepsilon}$ is defined in \eqref{a81} (see also Section \ref{a71} for some properties) for $\alpha$ in a neighborhood of $\CH_{\rm tang}^{+ \infty}$. In particular, for $\omega \in \pi_{x} ( \CH_{{\rm tang} , -}^{\varepsilon} ) / \varepsilon$,
\begin{equation} \label{b56}
\alpha_{\omega} = \alpha ( F ( \varepsilon , \omega ) ) \qquad \text{and} \qquad x_{\omega} = \varepsilon F^{-1} ( \varepsilon , \alpha ( F ( \varepsilon , \omega ) ) ) .
\end{equation}
where $F$ is defined in \eqref{a64}. Taking the modulus in \eqref{b51}, we get
\begin{align}
\vert a_{-} ( \varepsilon \omega , h ) \vert &\leq e^{- T_{\omega} \im z / h} \CM_{\varepsilon} ( \alpha_{\omega} ) \vert a_{+} ( x_{\omega} , h ) \vert + \CO_{\varepsilon} (h)   \nonumber \\
&= \CM_{\varepsilon} ( \alpha_{\omega} ) \vert a_{+} ( x_{\omega} , h ) \vert + \CO_{\varepsilon} \big( \vert \ln h \vert^{-1} \big) ,    \label{b52}
\end{align}
for all $\omega \in V_{-}^{\varepsilon} / \varepsilon$. Here, we have used that $\im z = \CO ( h \vert \ln h \vert^{-1} )$ and that the continuous function $\omega \longmapsto T_{\omega}$ is bounded on $V_{-}^{\varepsilon} / \varepsilon$.

Combining \eqref{b55} with \eqref{b52}, we obtain
\begin{align}
\vert a_{-} ( \varepsilon \omega , h ) \vert \leq{}& h^{\frac{\im z}{h \lambda}} \CM_{\varepsilon} ( \alpha_{\omega} ) \int_{\CH_{\rm tang}^{- \infty} \cap \pi_{x} ( \CH_{{\rm tang} , -}^{\varepsilon} ) / \varepsilon} f_{0} ( \varepsilon^{-1} x_{\omega} , \widetilde{\omega} ) \vert a_{-} ( \varepsilon \widetilde{\omega} , h ) \vert \, d \widetilde{\omega}    \nonumber \\
&+ \big( o_{\varepsilon \to 0} (1) + \CO_{\varepsilon} \big( \vert \ln h \vert^{-1} \big) \big) \Vert a_{-} ( \cdot , h ) \Vert_{L^{\infty}( V_{-}^{\varepsilon} )} + \CO_{\varepsilon} ( h^{\zeta} ) .   \label{b54}
\end{align}
Here, we have used that $\CM_{\varepsilon} ( \alpha_{\omega} )$ is uniformly bounded with respect to $\omega \in V_{-}^{\varepsilon} / \varepsilon$ and $\varepsilon$ if $V_{-}^{\varepsilon}$ is taken sufficiently close to $\pi_{x} ( \CH_{{\rm tang} , -}^{\varepsilon} )$. This follows from the fact that $\CM_{\varepsilon}$ is uniformly bounded on $\CH_{\rm tang}^{+ \infty}$ by Proposition \ref{b53}. Let $\widetilde{\CT}_{\varepsilon}$ be the operator on $L^{\infty} ( \S^{n-1} )$ with kernel
\begin{equation} \label{b61}
\widetilde{\CT}_{\varepsilon} ( \omega , \widetilde{\omega} ) = h^{\frac{\im z}{h \lambda}} \one_{V_{-}^{\varepsilon} / \varepsilon} ( \omega ) \CM_{\varepsilon} ( \alpha_{\omega} ) f_{0} ( \varepsilon^{-1} x_{\omega} , \widetilde{\omega} ) \one_{\CH_{\rm tang}^{- \infty} \cap \pi_{x} ( \CH_{{\rm tang} , -}^{\varepsilon} ) / \varepsilon} ( \widetilde{\omega} ) .
\end{equation}
Estimate \eqref{b54} can be written as
\begin{align}
\vert a_{-} ( \varepsilon \omega , h ) \vert \leq{}& \widetilde{\CT}_{\varepsilon} \vert a_{-} ( \varepsilon \cdot , h ) \vert \nonumber  \\
&+ \big( o_{\varepsilon \to 0} (1) + \CO_{\varepsilon} \big( \vert \ln h \vert^{-1} \big) \big) \Vert a_{-} ( \cdot , h ) \Vert_{L^{\infty}( V_{-}^{\varepsilon} )} + \CO_{\varepsilon} ( h^{\zeta} ) . \label{b59}
\end{align}
Moreover, \ref{h4}, \eqref{a99}, \eqref{b57} and the bound on $\CM_{\varepsilon} ( \alpha_{\omega} )$ alluded to after \eqref{b54} imply
\begin{equation} \label{b58}
\big\Vert \widetilde{\CT}_{\varepsilon} \big\Vert_{L^{\infty} ( \S^{n-1} ) \to L^{\infty} ( \S^{n-1} )} = \CO ( 1 ) ,
\end{equation}
uniformly in $\varepsilon$ and $h$. Thus, iterating \eqref{b59}, this implies
\begin{align}
\vert a_{-} ( \varepsilon \omega , h ) \vert \leq{}& \widetilde{\CT}_{\varepsilon}^{k + 1} \vert a_{-} ( \varepsilon \cdot , h ) \vert    \nonumber \\
&+ \big( o_{\varepsilon \to 0}^{k} (1) + \CO_{\varepsilon , k} \big( \vert \ln h \vert^{-1} \big) \big) \Vert a_{-} ( \cdot , h ) \Vert_{L^{\infty}( V_{-}^{\varepsilon} )} + \CO_{\varepsilon , k} ( h^{\zeta} ) ,   \label{b60}
\end{align}
for all $k \in \N$. The notation $o_{\varepsilon \to 0}^{k} (1)$ means that the $o_{\varepsilon \to 0} (1)$ may depend on $k$. We then remark that
\begin{equation} \label{b63}
\widetilde{\CT}_{\varepsilon}^{2} = \widetilde{\CT}_{\varepsilon} \one_{\CH_{\rm tang}^{- \infty} \cap \pi_{x} ( \CH_{{\rm tang} , -}^{\varepsilon} ) / \varepsilon} ( \omega ) \widetilde{\CT}_{\varepsilon} .
\end{equation}
Now, for $\omega \in \CH_{\rm tang}^{- \infty} \cap \pi_{x} ( \CH_{{\rm tang} , -}^{\varepsilon} ) / \varepsilon$, the function $\CM$ is well-defined on $\alpha_{\omega} , \alpha ( \omega ) \in \CH_{\rm tang}^{+ \infty}$. Thus, using \eqref{b56}, \eqref{a64}, the uniform continuity of $\omega \mapsto \alpha ( \omega )$ on the compact set $\CH_{\rm tang}^{- \infty}$ from Lemma \ref{b67} and the uniform continuity of $( \varepsilon , \alpha ) \mapsto \CM_{\varepsilon} ( \alpha ) $ on the compact set $[ 0 , \varepsilon_{0} ] \times \CH_{\rm tang}^{+ \infty}$ from Proposition \ref{b53}, we obtain
\begin{equation} \label{b62}
\CM_{\varepsilon} ( \alpha_{\omega} ) = \CM_{0} ( \alpha ( \omega ) ) + o_{\varepsilon \to 0} (1) ,
\end{equation}
uniformly for $\omega \in \CH_{\rm tang}^{- \infty} \cap \pi_{x} ( \CH_{{\rm tang} , -}^{\varepsilon} ) / \varepsilon$. The same way, from \eqref{b56}, \eqref{a64}, \eqref{a68} and the uniform continuity of $\alpha$ on the compact set $\CH_{\rm tang}^{- \infty}$ (see Lemma \ref{b67})), we have
\begin{equation} \label{b66}
f_{0} ( \varepsilon^{-1} x_{\omega} , \widetilde{\omega} ) = f_{0} ( \alpha ( \omega ) , \widetilde{\omega} ) + o_{\varepsilon \to 0} (1) ,
\end{equation}
uniformly for $\omega , \widetilde{\omega} \in \CH_{\rm tang}^{- \infty} \cap \pi_{x} ( \CH_{{\rm tang} , -}^{\varepsilon} ) / \varepsilon$. Combining \eqref{b61}, \eqref{b58}, \eqref{b62} and \eqref{b66}, \eqref{b63} gives
\begin{equation} \label{b64}
\widetilde{\CT}_{\varepsilon}^{2} = \widetilde{\CT}_{\varepsilon} \CT_{\varepsilon} + o_{\varepsilon \to 0} (1) ,
\end{equation}
where $\CT_{\varepsilon}$ is the operator on $L^{\infty} ( \S^{n-1} )$ with kernel
\begin{equation*}
\CT_{\varepsilon} ( \omega , \widetilde{\omega} ) = h^{\frac{\im z}{h \lambda}} \one_{\CH_{\rm tang}^{- \infty}} ( \omega ) \CM_{0} ( \alpha ( \omega ) ) f_{0} ( \alpha ( \omega ) , \widetilde{\omega} ) \one_{\CH_{\rm tang}^{- \infty} \cap \pi_{x} ( \CH_{{\rm tang} , -}^{\varepsilon} ) / \varepsilon} ( \widetilde{\omega} ) .
\end{equation*}
Thus, \eqref{b60} implies
\begin{align*}
\vert a_{-} ( \varepsilon \omega , h ) \vert \leq{}& \widetilde{\CT}_{\varepsilon} \CT_{\varepsilon}^{k} \one_{\CH_{\rm tang}^{- \infty} \cap \pi_{x} ( \CH_{{\rm tang} , -}^{\varepsilon} ) / \varepsilon} \vert a_{-} ( \varepsilon \cdot , h ) \vert    \\
&+ \big( o_{\varepsilon \to 0}^{k} (1) + \CO_{\varepsilon , k} \big( \vert \ln h \vert^{-1} \big) \big) \Vert a_{-} ( \cdot , h ) \Vert_{L^{\infty}( V_{-}^{\varepsilon} )} + \CO_{\varepsilon , k} ( h^{\zeta} ) .
\end{align*}
Since the function $\vert a_{-} ( \cdot , h ) \vert$ and the kernels of the operators $\widetilde{\CT}_{\varepsilon}$ and $\CT_{\varepsilon}$ are non-negative, the last estimate gives
\begin{align}
\vert a_{-} ( \varepsilon \omega , h ) \vert \leq{}& \widetilde{\CT}_{\varepsilon} \CT^{k} \one_{\CH_{\rm tang}^{- \infty} \cap \pi_{x} ( \CH_{{\rm tang} , -}^{\varepsilon} ) / \varepsilon} \vert a_{-} ( \varepsilon \cdot , h ) \vert    \nonumber  \\
&+ \big( o_{\varepsilon \to 0}^{k} (1) + \CO_{\varepsilon , k} \big( \vert \ln h \vert^{-1} \big) \big) \Vert a_{-} ( \cdot , h ) \Vert_{L^{\infty}( V_{-}^{\varepsilon} )} + \CO_{\varepsilon , k} ( h^{\zeta} ) .  \label{b79}
\end{align}
where $\CT$ is the operator on $L^{\infty} ( \CH_{\rm tang}^{- \infty} )$ with kernel
\begin{equation*}
\CT ( \omega , \widetilde{\omega} ) = h^{\frac{\im z}{h \lambda}} \CM_{0} ( \alpha ( \omega ) ) f_{0} ( \alpha ( \omega ) , \widetilde{\omega} ) .
\end{equation*}
Note that, by \eqref{b57} and \eqref{b70}, we can write
\begin{equation} \label{b78}
\CT = h^{\frac{\im z}{h \lambda}} \CT_{0} \Big( \frac{\re z - E_{0}}{\lambda h} \Big) .
\end{equation}
Using \eqref{b58} and \eqref{b79}, we have proved

\begin{lemma}\sl \label{b65}
There exists a constant $R > 0$ independent of $\varepsilon , k , h$ such that
\begin{align*}
\Vert a_{-} ( \cdot , h ) \Vert_{L^{\infty} ( V_{-}^{\varepsilon} )} \leq{}& R \big\Vert \CT^{k} \big\Vert_{L^{\infty} ( \CH_{\rm tang}^{- \infty} ) \to L^{\infty} ( \CH_{\rm tang}^{- \infty} )} \Vert a_{-} ( \cdot , h ) \Vert_{L^{\infty} ( V_{-}^{\varepsilon} )}   \\
&+ \big( o_{\varepsilon \to 0}^{k} (1) + \CO_{\varepsilon , k} \big( \vert \ln h \vert^{-1} \big) \big) \Vert a_{-} ( \cdot , h ) \Vert_{L^{\infty}( V_{-}^{\varepsilon} )} + \CO_{\varepsilon , k} ( h^{\zeta} ) .
\end{align*}
\end{lemma}

The imaginary part of $z \in \eqref{b68}$ satisfies
\begin{equation*}
\im z \frac{\vert \ln h \vert}{\lambda h} \geq \ln \Big( \CA_{0} \Big( \frac{\re z - E_{0}}{\lambda h} \Big) \Big) + \delta .
\end{equation*}
Combining with \eqref{b71} and \eqref{b78}, we get
\begin{equation*}
\spr ( \CT ) = h^{\frac{\im z}{h \lambda}} \spr \Big( \CT_{0} \Big( \frac{\re z - E_{0}}{\lambda h} \Big) \Big) \leq e^{- \delta} .
\end{equation*}
Then, by Proposition \ref{c8}, we have
\begin{equation*}
\big\Vert \CT^{k} \big\Vert_{L^{\infty} ( \CH_{\rm tang}^{- \infty} ) \to L^{\infty} ( \CH_{\rm tang}^{- \infty} )} \leq e^{- k \frac{\delta}{2}} ,
\end{equation*}
for all $k$ sufficiently large and $z \in \eqref{b68}$. Therefore, taking first $k$ large enough, and then $\varepsilon$ small enough, Lemma \ref{b65} yields
\begin{equation*}
\Vert a_{-} ( \cdot , h ) \Vert_{L^{\infty} ( V_{-}^{\varepsilon} )} \leq \frac{1}{2} \Vert a_{-} ( \cdot , h ) \Vert_{L^{\infty} ( V_{-}^{\varepsilon} )} + \CO ( h^{\zeta} ) ,
\end{equation*}
for $h$ small enough. Therefore,
\begin{equation*}
\Vert a_{-} ( \cdot , h ) \Vert_{L^{\infty} ( V_{-}^{\varepsilon} )} = \CO ( h^{\zeta} ) .
\end{equation*}
Using now the transport equations \eqref{b8}--\eqref{b9} and the initial conditions \eqref{b10}--\eqref{b11}, the previous estimate gives
\begin{equation*}
\Vert a_{-} ( \cdot , h ) \Vert_{L^{\infty}} = \CO ( h^{\zeta} ) .
\end{equation*}
Since $a_{-} \in S (1)$, the Landau--Kolmogorov inequalities (see e.g. Ditzian \cite{Di89_01}) imply $a_{-} \in S ( h^{\frac{\zeta}{2}} )$. In other words, we have proved that
\begin{equation*}
u_{-}^{\rm tang} \in \CI ( \Lambda_{+}^{1} , h^{- N_{-}^{\rm tang}} ) \quad \Longrightarrow \quad u_{-}^{\rm tang} \in \CI ( \Lambda_{+}^{1} , h^{- N_{-}^{\rm tang} + \frac{\zeta}{2}} ) .
\end{equation*}
Then, $u_{-}^{\rm tang} = \CO ( h^{\infty} )$ and, by Lemma \ref{a94},
\begin{equation} \label{d43}
\Vert u_{-} \Vert = \CO ( h^{\infty} ) .
\end{equation}
Now, as in the proof of Theorem \ref{a1} (see the end of Section \ref{s4}), this implies 
\begin{equation*}
u = 0 \text{ microlocally near } K ( E_{0} ) ,
\end{equation*}
and this finishes the proof of Theorem \ref{a2} thanks to Section \ref{s31} and Proposition \ref{a16}.

\section{Proof of the asymptotic of the resonances for a finite number of homoclinic curves} \label{s12}

This section is devoted to the proof of Theorem \ref{d8} and, more generally, of all the results stated in Section \ref{s61}. For that, we follow the general strategy explained in Section \ref{s36}. But before developing this approach, we study the asymptotic behavior of the pseudo-resonances (see Definition \ref{d1}) and the spectral properties of $\CQ$ defined in \eqref{d4}.

\Subsection{From the quantization rule to the asymptotic of the pseudo-resonances} \label{s13}

The aim of this part is to prove Proposition \ref{d9} and to obtain the estimate of the resolvent of $\CQ$ away from the pseudo-resonances stated in Lemma \ref{d13}. We begin with an inequality which strongly relies on the particular structure of the matrix $\CQ ( z , h )$ described in \eqref{d10}.

\begin{lemma}\sl \label{d14}
Let $\alpha > 0$ and $\CK$ be a compact of $\C$ on which
\begin{equation} \label{d15}
\sigma \longmapsto \Gamma \Big( \sum_{j = 1}^{n} \frac{\lambda_{j}}{2 \lambda_{1}} - i \frac{\sigma}{\lambda_{1}} \Big) ,
\end{equation}
is regular. Then, there exists $M > 0$ such that, for all $\rho \in ( \S^{1} )^{K}$, $\sigma \in \CK$ and $\Lambda \in \C$ such that $\dist ( \Lambda , \spe ( \widetilde{\CQ} ( \rho , \sigma ) ) ) \geq \alpha$, we have
\begin{equation} \label{d16}
\big\Vert \big( \widetilde{\CQ} ( \rho , \sigma ) - \Lambda \big)^{- 1} \big\Vert \leq M .
\end{equation}
\end{lemma}

\begin{proof}
First, we remark that, since the function \eqref{d15} is regular on $\CK$, the matrix $\widetilde{\CQ} ( \rho , \sigma )$ is uniformly bounded on $( \S^{1} )^{K} \times \CK$. Thus, \eqref{d16} holds true for $\Lambda$ large enough. In the sequel, we will thus assume that $\vert \Lambda \vert \leq N$ for some $N > 0$.

As before, from the definition of $\widetilde{\CQ}$ and the regularity of the function \eqref{d15} on $\CK$, the application $( \rho , \sigma ) \mapsto \widetilde{\CQ} ( \rho , \sigma )$ is continuous on $( \S^{1} )^{K} \times \CK$. As a consequence, $( \rho , \sigma ) \mapsto \spe ( \widetilde{\CQ} ( \rho , \sigma ) )$ is also continuous by the usual perturbation theory. Thus, the map
\begin{equation}
( \rho , \sigma , \Lambda ) \longmapsto \big( \widetilde{\CQ} ( \rho , \sigma ) - \Lambda \big)^{- 1} ,
\end{equation}
is continuous (as composition of continuous applications) on the compact set
\begin{equation*}
\big\{ ( \rho , \sigma , \Lambda ) \in ( \S^{1} )^{K} \times \CK \times \C ; \ \dist \big( \Lambda , \spe ( \widetilde{\CQ} ( \rho , \sigma ) ) \big) \geq \alpha \text{ and } \vert \Lambda \vert \leq N \big\} .
\end{equation*}
Therefore, this map is bounded and the lemma follows.
\end{proof}

\begin{proof}[Proof of Proposition \ref{d9}]
First, note that $\CQ ( z , h )$ is uniformly bounded on the set defined in \eqref{d12}. Thus, for $z$ in this set satisfying additionally
\begin{equation} \label{d20}
\im z > - \sum_{j = 2}^{n} \frac{\lambda_{j}}{2} h + N \frac{h}{\vert \ln h \vert} ,
\end{equation}
we get, for some $c > 0$,
\begin{equation} \label{d36}
\big\Vert h^{S ( z , h ) / \lambda_{1} - 1 / 2} \CQ ( z , h ) \big\Vert \leq c h^{\sum_{j=1}^{n} \frac{\lambda_{j}}{2 \lambda_{1}} + \frac{\im z}{\lambda_{1} h} - \frac{1}{2}} \leq c e^{- N / \lambda_{1}} \leq 1 / 2 ,
\end{equation}
for $N$ large enough. This implies that no pseudo-resonance verifies \eqref{d20}. On the other hand, using again that $\CQ ( z , h )$ is uniformly bounded, the eigenvalues $\mu_{1} ( \tau , h ) , \ldots , \mu_{K} ( \tau , h )$ are also uniformly bounded. Thus, $z_{q , k} ( \tau )$ never verifies \eqref{d20}. We will then restrict our study to the complement of the domain \eqref{d20}.

We now show that every pseudo-resonance satisfies \eqref{d94}. We fix $C , \delta (h)$ as in Proposition \ref{d9}. Let $\tau \in [ - C , C ]$ and $z$ be a pseudo-resonance in the set
\begin{equation} \label{d21}
E_{0} + [ - C h , C  h ] - i \sum_{j = 2}^{n} \frac{\lambda_{j}}{2} h + i \frac{h}{\vert \ln h \vert} [ - C , N ] ,
\end{equation}
satisfying $\re z \in E_{0} + \tau h + h \delta (h) [ - 1 , 1 ]$. In particular, 
\begin{equation*}
\sigma = \frac{z - E_{0}}{h} \qquad \text{and} \qquad \widetilde{\sigma} = \tau - i \sum_{j = 2}^{n} \frac{\lambda_{j}}{2} ,
\end{equation*}
verify $\sigma = \widetilde{\sigma} + o (1)$ uniformly with respect to $\tau$ and $z$. Eventually, by Definition \ref{d1},
\begin{equation} \label{d18}
h^{-S ( z , h ) / \lambda_{1} + 1 / 2} \in \spe ( \CQ ( z , h ) ) .
\end{equation}

On the other hand, with the notations of \eqref{d10}, we can write $\CQ ( z , h ) = \widetilde{\CQ} ( \rho (h) , \sigma ) = \sum \rho_{k} (h) \widetilde{\CQ}_{k} ( \sigma )$ where $\rho_{k} = e^{i A_{k} / h}$ and $\sigma \mapsto \widetilde{\CQ}_{k} ( \sigma )$ is uniformly continuous (in fact, holomorphic) on the compact set
\begin{equation*}
\sigma \in [ - C , C  ] + i \Big[ - \sum_{j = 2}^{n} \frac{\lambda_{j}}{2} - \nu , 1 \Big] ,
\end{equation*}
for $\nu > 0$ small enough. Combining with $\sigma = \widetilde{\sigma} + o (1)$, it yields
\begin{equation} \label{d29}
\CQ ( z , h ) = \widetilde{\CQ} ( \rho (h) , \widetilde{\sigma} ) + o (1) ,
\end{equation}
uniformly with respect to $\tau$ and $z$. Using in addition that $\CQ ( z , h )$ is uniformly bounded on the set \eqref{d21}, the uniform continuity on bounded regions of the map which, at a $K \times K$ matrix, associates its spectrum (see Theorem II.5.14 of Kato \cite{Ka76_01}) implies that
\begin{equation} \label{d26}
\dist \big( \spe ( \CQ ( z , h ) ) , \spe ( \widetilde{\CQ} ( \rho (h) , \widetilde{\sigma} ) ) \big) = o (1) .
\end{equation}
Then, since the $\mu_{k} ( \tau , h )$ are precisely the eigenvalues of $\widetilde{\CQ} ( \rho (h) , \widetilde{\sigma} )$, \eqref{d18} and the previous estimate prove that there exists $k = k ( \tau , z , h ) \in \{ 1 , \ldots , K \}$ such that
\begin{equation} \label{d19}
h^{-S ( z , h ) / \lambda_{1} + 1 / 2} = \mu_{k} ( \tau , h ) + o (1) ,
\end{equation}
uniformly with respect to $\tau$ and $z$. Moreover, from $\im z > - \sum_{j = 2}^{n} \frac{\lambda_{j}}{2} h - C \frac{h}{\vert \ln h \vert}$, we get
\begin{equation*}
\big\vert h^{-S ( z , h ) / \lambda_{1} + 1 / 2} \big\vert = h^{- \sum_{j=1}^{n} \frac{\lambda_{j}}{2 \lambda_{1}} - \frac{\im z}{\lambda_{1} h} + \frac{1}{2}} \geq e^{- C / \lambda_{1}} > 0 .
\end{equation*}
Then, \eqref{d19} can be written
\begin{equation} \label{d38}
e^{\sum_{j=2}^{n} \frac{\lambda_{j}}{2 \lambda_{1}} \vert \ln h \vert - i \frac{z - E_{0}}{\lambda_{1}} \frac{\vert \ln h \vert}{h}} = h^{-S ( z , h ) / \lambda_{1} + 1 / 2} = \mu_{k} ( \tau , h ) ( 1 + o (1) ) = e^{\ln ( \mu_{k} ( \tau , h ) ) + o (1)} ,
\end{equation}
uniformly with respect to $\tau$ and $z$. This implies \eqref{d94}.

We now show that every complex number given by \eqref{d95} is close to a pseudo-resonance. More precisely, let $C , \delta (h)$ be as in Proposition \ref{d9}. We have then to prove that
\begin{align}
\forall \varepsilon > 0 , \quad \exists h_{0} > 0 , \quad \forall h < h_{0} , \quad & \forall \tau \in [ - C , C ] , \quad \forall z_{q , k} ( \tau ) \in \Omega ( \tau ) , \nonumber \\
&\exists z \in \res_{0} (P) , \qquad \vert z - z_{q , k} ( \tau ) \vert < \varepsilon \frac{h}{\vert \ln h \vert} ,  \label{d22}
\end{align}
with the notation
\begin{equation*}
\Omega ( \tau ) = E_{0} + \tau h + h \delta (h) [ - 1 , 1 ] - i \sum_{j = 2}^{n} \frac{\lambda_{j}}{2} h + i \frac{h}{\vert \ln h \vert} [ - C , N ] .
\end{equation*}
Thus, let $\varepsilon > 0$, $\tau \in [ - C , C ]$ and $z_{q , k} ( \tau ) \in \Omega ( \tau )$. The idea is to show that the operator
\begin{equation} \label{d23}
I : = \int_{\partial \CD} \big( 1 - h^{S ( z , h ) / \lambda_{1} - 1 / 2} \CQ ( z , h ) \big )^{- 1} d z ,
\end{equation}
satisfies $I \neq 0$ for a well-chosen domain $\CD \subset \C$ containing $z_{q , k} ( \tau )$ (uniformly for $h$ small enough, $\tau$ and $z_{q , k} ( \tau )$). Indeed, it will imply that $1 - h^{S ( z , h ) / \lambda_{1} - 1 / 2} \CQ ( z , h )$ is not invertible for all $z$ in $\CD$ and then that $P$ has at least one pseudo-resonance in $\CD$. In the sequel, it will be more convenient to use the rescaled spectral parameter
\begin{equation*}
\Lambda = \Lambda ( z , h ) = h^{-S ( z , h ) / \lambda_{1} + 1 / 2} ,
\end{equation*}
which will be typically of size $1$.

\begin{lemma}\sl \label{d24}
There exists $\alpha > 0$ such that, for all $h$ small enough, $\tau \in [ - C , C]$ and $z_{q , k} ( \tau ) \in \Omega ( \tau )$, there exists $\beta = \beta ( h , \tau , q , k )$ with $\alpha < \beta < \varepsilon$ and
\begin{equation*}
z \in \partial B \Big( z_{q , k} ( \tau ) , \beta \frac{h}{\vert \ln h \vert} \Big) \quad \Longrightarrow \quad \dist \big( \Lambda , \spe ( \widetilde{\CQ} ( \rho (h) , \widetilde{\sigma} ) ) \big) \geq 2 \alpha .
\end{equation*}
\end{lemma}

\begin{proof}
Indeed, for $z = z_{q , k} ( \tau ) + r \frac{h}{\vert \ln h \vert}$, we have
\begin{equation} \label{d25}
\Lambda = \mu_{k} ( \tau , h ) e^{- i r / \lambda_{1}} = \mu_{k} ( \tau , h ) - i r \frac{\mu_{k} ( \tau , h )}{\lambda_{1}} + o_{r \to 0} ( r ) \vert  \mu_{k} ( \tau , h ) \vert .
\end{equation}
On the other hand, $z_{q , k} ( \tau ) \in \Omega ( \tau )$ implies $\vert \mu_{k} ( \tau , h ) \vert \geq e^{- C / \lambda_{1}}$. Moreover, using one more time that $\CQ ( z , h )$ is uniformly bounded on the set \eqref{d12}, there exists $\widetilde{C} > 0$ such that $\vert \mu_{k} ( \tau , h ) \vert \leq \widetilde{C}$. Then, combining \eqref{d25} with the previous estimates on $\mu_{k}$, we can fix $\alpha > 0$ and arrange $\alpha < \beta < \varepsilon$ so that, for $r \in \partial B ( 0 , \beta )$, $\Lambda$ stays at distance at least $2 \alpha$ of the $K$ eigenvalues of $\widetilde{\CQ} ( \rho (h) , \widetilde{\sigma} )$.
\end{proof}

In \eqref{d23}, we choose the domain $\CD = B \big( z_{q , k} ( \tau ) , \beta \frac{h}{\vert \ln h \vert} \big)$ constructed in the previous lemma. Working as in \eqref{d29} and \eqref{d26}, we get
\begin{equation} \label{d30}
\CQ ( z , h ) = \widetilde{\CQ} ( \rho (h) , \widetilde{\sigma} ) + o (1) ,
\end{equation}
and
\begin{equation}
\dist \big( \spe ( \CQ ( z , h ) ) , \spe ( \widetilde{\CQ} ( \rho (h) , \widetilde{\sigma} ) ) \big) = o (1) ,
\end{equation}
uniformly in $h$ small enough, $\tau$, $z_{q , k} ( \tau )$ and $z \in \partial \CD$. Combining with Lemma \ref{d24}, it yields
\begin{equation*}
\dist \big( \Lambda , \spe ( \CQ ( z , h ) ) \cup \spe ( \widetilde{\CQ} ( \rho (h) , \widetilde{\sigma} ) ) \big) \geq \alpha .
\end{equation*}
for all $z \in \partial \CD$. Now Lemma \ref{d14} implies
\begin{equation} \label{d28}
\big\Vert \big( \Lambda - \CQ ( z , h ) \big)^{- 1} \big\Vert + \big\Vert \big( \Lambda - \widetilde{\CQ} ( \rho (h) , \widetilde{\sigma} ) \big)^{- 1} \big\Vert \leq 2 M ,
\end{equation}
for all $z \in \partial \CD$. In particular, combining with \eqref{d30} and the resolvent identity, it gives
\begin{align}
\big( \Lambda - \CQ ( z , h ) \big)^{- 1} ={}& \big( \Lambda  - \widetilde{\CQ} ( \rho (h) , \widetilde{\sigma} ) \big)^{- 1}  \nonumber \\
&+ \big( \Lambda - \CQ ( z , h ) \big)^{- 1} \big( \CQ ( z , h ) - \widetilde{\CQ} ( \rho (h) , \widetilde{\sigma} ) \big) \big( \Lambda  - \widetilde{\CQ} ( \rho (h) , \widetilde{\sigma} ) \big)^{- 1}  \nonumber \\
={}& \big( \Lambda  - \widetilde{\CQ} ( \rho (h) , \widetilde{\sigma} ) \big)^{- 1} + o (1) . \label{d33}
\end{align}
Then, since $\vert \partial \CD \vert$ is of size $h \vert \ln h \vert^{-1}$ and $\Lambda$ is uniformly bounded on $\partial \CD$ (see \eqref{d25}), \eqref{d23} and \eqref{d33} yield
\begin{align}
I &= \int_{\partial \CD} \Lambda \big( \Lambda - \CQ ( z , h ) \big )^{- 1} d z  \nonumber \\
&= \int_{\partial \CD} \Lambda \big( \Lambda  - \widetilde{\CQ} ( \rho (h) , \widetilde{\sigma} ) \big)^{- 1} d z + o \big( h \vert \ln h \vert^{- 1} \big) .   \label{d32}
\end{align}
We now use the change of variable $z \mapsto \Lambda = \Lambda ( z , h )$ which verifies
\begin{equation*}
\Lambda \, d z = i \frac{\lambda_{1} h}{\vert \ln h \vert} d \Lambda .
\end{equation*}
Then, \eqref{d32} becomes
\begin{equation} \label{d34}
I = i \frac{\lambda_{1} h}{\vert \ln h \vert} \int_{\Lambda ( \partial \CD )} \big( \Lambda - \widetilde{\CQ} ( \rho (h) , \widetilde{\sigma} ) \big)^{- 1} d \Lambda + o \big( h \vert \ln h \vert^{- 1} \big) ,
\end{equation}
Since $\Lambda ( \partial \CD )$ is a simple loop around $\mu_{k} ( \tau , h )$ and $\widetilde{\CQ} ( \rho (h) , \widetilde{\sigma} )$ has at least one eigenvalue (namely $\mu_{k} ( \tau , h )$) inside this loop, the Cauchy formula implies
\begin{equation*}
\Vert I \Vert \geq 2 \pi \lambda_{1} \frac{h}{\vert \ln h \vert} ( 1 + o ( 1 ) ) \geq \pi \lambda_{1} \frac{h}{\vert \ln h \vert} \neq 0 .
\end{equation*}
Here also, this estimate holds true uniformly for $h$ small enough, $\tau \in [ - C , C ]$ and $z_{q , k} ( \tau ) \in \Omega ( \tau )$. In other words, there exists at least one pseudo-resonance in $\CD \subset B \big( z_{q , k} ( \tau ) , \varepsilon \frac{h}{\vert \ln h \vert} \big)$ and \eqref{d22} is verified.
\end{proof}

We finish this part with an estimate of the resolvent of $\CQ$ away from the pseudo-resonances, which will be used in the proof of the asymptotic of the resonances.

\begin{lemma}\sl \label{d13}
Let $\beta , C > 0$. Then, there exists $M > 0$ such that, for all $z$ in \eqref{d12}, we have
\begin{equation*}
\dist \big( z , \res_{0} (P) \big) > \beta \frac{h}{\vert \ln h \vert} \quad \Longrightarrow \quad \big\Vert \big( 1 - h^{S ( z , h ) / \lambda_{1} - 1 / 2} \CQ ( z , h ) \big)^{- 1} \big\Vert \leq M .
\end{equation*}
\end{lemma}

\begin{proof}
As already remarked in \eqref{d20}--\eqref{d36}, we can assume that $z$ belongs to \eqref{d21} since the required estimate is clear for $z$ satisfying \eqref{d20}. Mimicking the previous proofs, we define
\begin{equation*}
\sigma = \frac{z - E_{0}}{h} , \qquad \widetilde{\sigma} = \re \sigma - i \sum_{j = 2}^{n} \frac{\lambda_{j}}{2} \qquad \text{and} \qquad \Lambda = h^{-S ( z , h ) / \lambda_{1} + 1 / 2} .
\end{equation*}
We claim that there exists $\alpha > 0$ such that, for all $z \in \eqref{d21}$, we have
\begin{equation} \label{d37}
\dist \big( z , \res_{0} (P) \big) > \beta \frac{h}{\vert \ln h \vert} \quad \Longrightarrow \quad \dist \big( \Lambda , \spe ( \widetilde{\CQ} ( \rho (h) , \sigma ) ) \big) > \alpha .
\end{equation}
Indeed, assume that \eqref{d37} does not hold true. Then, there exists a sequence of $z = z (h) \in \eqref{d21}$ with $\dist ( z , \res_{0} (P) ) > \beta h \vert \ln h \vert^{-1}$ and such that $\dist ( \Lambda , \spe ( \widetilde{\CQ} ( \rho (h) , \sigma ) ) ) = o (1)$. Using $\sigma = \widetilde{\sigma} + o (1)$, we conclude as in \eqref{d29}--\eqref{d38} that there exists $k \in \{ 1 , \ldots , K \}$ and $q \in \Z$ (which may depend on $h$) such that
\begin{equation*}
z = z_{q , k} ( \re \sigma ) + o \Big( \frac{h}{\vert \ln h \vert} \Big) .
\end{equation*}
Applying the second part of Proposition \ref{d9}, this yields $\dist ( z , \res_{0} (P) )= o ( h \vert \ln h \vert^{-1} )$ which is clearly in contradiction with the left hand side of \eqref{d37}. This proves \eqref{d37} and the lemma follows from Lemma \ref{d14}.
\end{proof}

\Subsection{Resonance free zone and resolvent estimate}  \label{s72}

Following the general strategy for the asymptotic of resonances explained in Section \ref{s36}, we first show that $P$ has no resonance away from the pseudo-resonances. In this way, we begin with the following proposition which is the half of Theorem \ref{d8}. The rest of this part is devoted to the proof of this result.

\begin{proposition}\sl \label{d39}
Let $C , \delta > 0$. For $h$ small enough, $P$ has no resonance in the domain
\begin{align}
E_{0} + [ - C h , C  h ] + i \Big[ & - \sum_{j = 2}^{n} \frac{\lambda_{j}}{2} h - C \frac{h}{\vert \ln h \vert} , h \Big]   \nonumber \\
&\setminus \big( \Gamma (h) + B ( 0 , \delta h ) \big) \bigcup \Big( \res_{0} (P) + B \Big( 0 , \delta \frac{h}{\vert \ln h \vert} \Big) \Big) .  \label{d40}
\end{align}
Moreover, there exists $M > 0$ such that
\begin{equation*}
\big\Vert( P_{\theta} -z )^{-1} \big\Vert \lesssim h^{- M} ,
\end{equation*}
uniformly for $h$ small enough and $z \in \eqref{d40}$.
\end{proposition}

In order to prove this result, we apply the general contradiction argument developed in Section \ref{s3} with the domain $\Omega_{h} = \eqref{d40}$. Thus, to prove Proposition \ref{d39}, it is enough to show that every $u$ satisfying \eqref{a4} vanishes microlocally near each point of $K ( E_{0} ) = \{ ( 0 , 0 ) \} \cup \CH$. As before, we can assume that $P_{\theta} = P$ near the base space projection of $\CH$. Since the assumptions of Theorem \ref{a1} and Theorem \ref{d8} are quite similar, we will use the notations and some results of Section \ref{s4}.

\begin{figure}
\begin{center}
\begin{picture}(0,0)%
\includegraphics{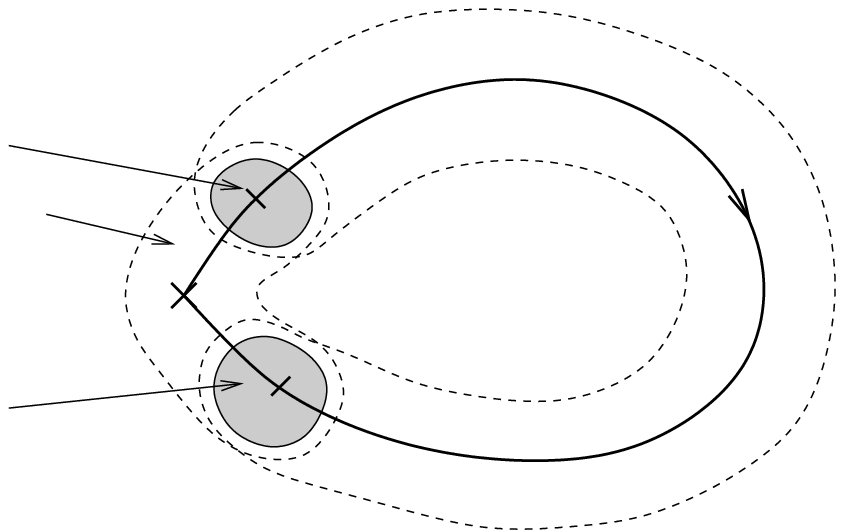}%
\end{picture}%
\setlength{\unitlength}{1579sp}%
\begingroup\makeatletter\ifx\SetFigFont\undefined%
\gdef\SetFigFont#1#2#3#4#5{%
  \reset@font\fontsize{#1}{#2pt}%
  \fontfamily{#3}\fontseries{#4}\fontshape{#5}%
  \selectfont}%
\fi\endgroup%
\begin{picture}(12801,6283)(-1364,-6786)
\put(6976,-961){\makebox(0,0)[lb]{\smash{{\SetFigFont{9}{10.8}{\rmdefault}{\mddefault}{\updefault}$\Omega_{\text{reg}}$}}}}
\put(1051,-3136){\makebox(0,0)[lb]{\smash{{\SetFigFont{9}{10.8}{\rmdefault}{\mddefault}{\updefault}$\Omega_{\text{sing}}$}}}}
\put(3901,-4036){\makebox(0,0)[lb]{\smash{{\SetFigFont{9}{10.8}{\rmdefault}{\mddefault}{\updefault}$0$}}}}
\put(-1349,-2236){\makebox(0,0)[lb]{\smash{{\SetFigFont{9}{10.8}{\rmdefault}{\mddefault}{\updefault}$u=u_{+}^{k}$ near $U_{+}^{k}$}}}}
\put(-1349,-5386){\makebox(0,0)[lb]{\smash{{\SetFigFont{9}{10.8}{\rmdefault}{\mddefault}{\updefault}$u=u_{-}^{k}$ near $U_{-}^{k}$}}}}
\put(4576,-3061){\makebox(0,0)[lb]{\smash{{\SetFigFont{9}{10.8}{\rmdefault}{\mddefault}{\updefault}$\rho_{+}^{k}$}}}}
\put(4201,-5386){\makebox(0,0)[lb]{\smash{{\SetFigFont{9}{10.8}{\rmdefault}{\mddefault}{\updefault}$\rho_{-}^{k}$}}}}
\put(10501,-2761){\makebox(0,0)[lb]{\smash{{\SetFigFont{9}{10.8}{\rmdefault}{\mddefault}{\updefault}$\gamma_{k}$}}}}
\end{picture}%
\end{center}
\caption{The geometric setting in the proof of Proposition \ref{d39}.} \label{f21}
\end{figure}

Recall that $\Omega_{\rm{sing}} = \Omega$ is a small neighborhood of $(0,0)$ and $\varepsilon > 0$ is a small enough constant given by Theorem \ref{a32}. As before, we define $S_{\pm}^{\varepsilon} = \{ ( x, \xi ) \in \Lambda_{\pm}^{0} ; \ \vert x \vert = \varepsilon \}$. For $k \in \{ 1 , \ldots , K \}$, we denote by $\rho_{\pm}^{k} = ( x_{\pm}^{k} , \xi_{\pm}^{k} )$ the unique point in the intersection $\gamma_{k} \cap S_{\pm}^{\varepsilon}$. In particular, there exists a unique time $t_{\pm}^{k} \in \R$ such that $\rho_{\pm}^{k} = \gamma_{k} ( t_{\pm}^{k} )$. Let $U_{+}^{k}$ be a restriction near $\rho_{+}^{k}$ of the set $V_{+}^{0}$ constructed in Section \ref{s4}. We then define $U_{-}^{k} = \exp ( ( t_{-}^{k} - t_{+}^{k} ) H_{p} ) ( U_{+}^{k} )$, the neighborhood of $\rho_{-}^{k}$ which corresponds to $V_{+}^{1}$ in the notations of Section \ref{s4}. The geometric setting and the previous definitions are illustrated in Figure \ref{f21}.

Mimicking the previous sections, we define $u_{\pm}^{k}$ as the microlocal restriction of $u$ to a neighborhood of $U_{\pm}^{k}$. Since $z \in \eqref{d40}$ avoids $\Gamma (h)$, we can apply Theorem \ref{a32} and we obtain, as in \eqref{a36} and \eqref{a38}, that
\begin{equation} \label{d44}
u_{-}^{k} \in \CI ( \Lambda_{+}^{1} , h^{-N} ) \qquad \text{and} \qquad u_{+}^{k} \in \CI ( \Lambda_{+}^{0} , h^{-N} ) ,
\end{equation}
for some $N \in \R$. From \ref{h8}, $\Lambda_{-}$ and $\Lambda_{+}$ intersect transversally along the Hamiltonian curve $\gamma_{k}$. Then, Proposition C.1 of \cite{ALBoRa08_01} states that $\Lambda_{+}^{1}$ projects nicely on the $x$-space near $U_{-}^{k}$ (after a possible shrinking of $U_{-}^{k}$ around $\rho_{-}^{k}$). Thus, let $\varphi_{+}^{1} (x)$ be the unique generating function of $\Lambda_{+}^{1}$ (i.e. $\Lambda_{+}^{1} = \Lambda_{\varphi_{+}^{1}} : = \{ ( x , \nabla \varphi_{+}^{1} (x) ) \}$) defined near $\pi_{x} ( U_{-}^{k} \cap \Lambda_{+}^{1} )$ with the normalization
\begin{equation} \label{d45}
\varphi_{+}^{1} ( x_{-}^{k} ) = \int_{\gamma_{k} ( ] - \infty , t_{-}^{k} ] )} \xi \cdot d x .
\end{equation}
Then, there exist symbols $a_{\pm}^{k} \in S ( h^{- N} )$ defined near $\pi_{x} ( U_{-}^{k} \cap \Lambda_{+}^{1} )$ and $\pi_{x} ( U_{+}^{k} \cap \Lambda_{+}^{0} )$ such that
\begin{equation} \label{d48}
\left\{ \begin{aligned}
u_{-}^{k} (x) &= e^{- i A_{k} / h} e^{i \frac{z - E_{0}}{h} t_{-}^{k}} \frac{\CM_{k}^{-}}{\CD_{k} ( t_{-}^{k} )} a_{-}^{k} ( x , h ) e^{i \varphi_{+}^{1} (x) / h} ,  \\
u_{+}^{k} (x) &= a_{+}^{k} ( x , h ) e^{i \varphi_{+} (x) / h} ,
\end{aligned} \right.
\end{equation}
where, using the notations of \eqref{d7},
\begin{equation} \label{g73}
\CD_{k} (t) = \sqrt{\Big\vert \det \frac{\partial x_{k} ( s , y )}{\partial ( s , y )} \vert_{s = t , \ y = 0} \Big\vert} .
\end{equation}
We have added the renormalization factor $e^{- i A_{k} / h} e^{i \frac{z - E_{0}}{h} t_{-}^{k}} \CM_{k}^{-} ( \CD_{k} ( t_{-}^{k} ) )^{- 1}$ so that the quantities calculated below will be expressed simply in terms of the matrix $\CQ$. We now compute the symbol $a_{-}^{k} ( x , h )$ after a turn around the critical point and $\gamma_{k}$.

\begin{lemma}\sl \label{d41}
There exist $\zeta > 0$ and symbols $\CP_{k , \ell} \in S (1)$ independent of $u$ such that
\begin{equation} \label{d42}
a_{-}^{k} ( x , h ) = h^{S ( z , h ) / \lambda_{1} - 1 / 2} \sum_{\ell = 1}^{K} \CP_{k , \ell} ( x , h ) a_{-}^{\ell} ( x_{-}^{\ell} , h ) + S ( h^{- N + \zeta} ) ,
\end{equation}
for all $x$ near $\pi_{x} ( U_{-}^{k} \cap \Lambda_{+}^{1} )$. Moreover, $\CP_{k , \ell} ( x_{-}^{k} , h ) = \CQ_{k , \ell} ( z , h )$ (see \eqref{d4}).
\end{lemma}

\begin{proof}
From \eqref{a4} and Lemma \ref{a33}, the function $u$ satisfies
\begin{equation} \label{d62}
\left\{ \begin{aligned}
&(P -z) u = 0 &&\text{microlocally near } \Omega_{\rm{sing}} ,   \\
&u = u_{-}^{k} &&\text{microlocally near } U_{-}^{k} ,   \\
&u = 0 &&\text{microlocally near } S_{-}^{\varepsilon} \setminus \big( U_{-}^{1} \cup \cdots \cup U_{-}^{K} \big) ,
\end{aligned} \right.
\end{equation}
and $\Vert u \Vert \leq 1$. By linearity, $u$ is then the sum over $\ell$ of the solutions of the following microlocal Cauchy problems
\begin{equation*}
\left\{ \begin{aligned}
&(P -z) u^{\ell} = 0 &&\text{microlocally near } \Omega_{\rm{sing}} ,   \\
&u^{\ell} = u_{-}^{\ell} &&\text{microlocally near } U_{-}^{\ell} ,   \\
&u^{\ell} = 0 &&\text{microlocally near } S_{-}^{\varepsilon} \setminus U_{-}^{\ell} ,
\end{aligned} \right.
\end{equation*}
where $u_{-}^{\ell}$ is the semiclassical Lagrangian distribution described in \eqref{d48}. Moreover, \ref{h4} implies that $g_{+} ( \rho_{+}^{k} ) \cdot g_{-} ( \rho_{-}^{\ell} ) \neq 0$ for all $k \in \{ 1 , \ldots , K \}$. Thus, we can apply Corollary \ref{d46} to compute $u^{\ell}$ in the outgoing region and eventually obtain
\begin{equation*}
a_{+}^{k} ( x , h ) = h^{S ( z , h ) / \lambda_{1} - 1 / 2} \sum_{\ell = 1}^{K} \CR_{k , \ell} ( x , h ) e^{- i A_{\ell} / h} e^{i \frac{z - E_{0}}{h} t_{-}^{\ell}} \frac{\CM_{\ell}^{-}}{\CD_{\ell} ( t_{-}^{\ell} )} a_{-}^{\ell} ( x_{-}^{\ell} , h ) + S ( h^{-N + \zeta} ) ,
\end{equation*}
for all $x$ near $\pi_{x} ( U_{+}^{k} \cap \Lambda_{+}^{0} )$ and some $0 < \zeta < 1$. Moreover, the symbols $\CR_{k , \ell} \in S (1)$ are given by \eqref{m79} and satisfy, in particular,
\begin{align*}
\CR_{k , \ell} ( x_{+}^{k} , h ) &= e^{i A_{\ell} / h} e^{- i \frac{\pi}{4}} \sqrt{\frac{\lambda_{1}}{2 \pi}} \Gamma \big( S ( z , h ) / \lambda_{1} \big) \big\vert g_{-} ( \rho_{-}^{\ell} ) \big\vert \big( i \lambda_{1} g_{+} ( \rho_{+}^{k} ) \cdot g_{-} ( \rho_{-}^{\ell} ) \big)^{- S ( z , h ) / \lambda_{1}}     \\
&\qquad \qquad \qquad \qquad \qquad \times \frac{\CD_{\ell} ( t_{-}^{\ell} )}{\CD_{k} ( t_{+}^{k} )} \lim_{s \to - \infty} \frac{\CD_{k} ( s + t_{+}^{k} )}{e^{s \sum_{j} \lambda_{j} / 2}} \lim_{s \to + \infty} \frac{e^{s ( \sum_{j} \lambda_{j} - 2 \lambda_{1} ) / 2}}{\CD_{\ell} ( s + t_{-}^{\ell} )}    \\
&= e^{i A_{\ell} / h} e^{- i \frac{\pi}{4}} \sqrt{\frac{\lambda_{1}}{2 \pi}} \Gamma \big( S ( z , h ) / \lambda_{1} \big) \big\vert g_{-}^{\ell} \big\vert \big( i \lambda_{1} g_{+}^{k} \cdot g_{-}^{\ell} \big)^{- S ( z , h ) / \lambda_{1}}     \\
&\qquad \qquad \qquad \qquad \qquad \qquad \qquad \qquad \qquad \times e^{i \frac{z - E_{0}}{h} ( t_{+}^{k} - t_{-}^{\ell} )} \frac{\CD_{\ell} ( t_{-}^{\ell} )}{\CD_{k} ( t_{+}^{k} )} \frac{\CM_{k}^{+}}{\CM_{\ell}^{-}} ,
\end{align*}
see also \eqref{d49}. Then, it yields
\begin{equation} \label{d51}
a_{+}^{k} ( x , h ) = h^{S ( z , h ) / \lambda_{1} - 1 / 2} \sum_{\ell = 1}^{K} \widetilde{\CR}_{k , \ell} ( x , h ) a_{-}^{\ell} ( x_{-}^{\ell} , h ) + S ( h^{-N + \zeta} ) ,
\end{equation}
where the symbols $\widetilde{\CR}_{k , \ell} \in S (1)$ satisfy
\begin{equation*}
\widetilde{\CR}_{k , \ell} ( x_{+}^{k} , h ) = e^{- i \frac{\pi}{4}} \sqrt{\frac{\lambda_{1}}{2 \pi}} \Gamma \big( S ( z , h ) / \lambda_{1} \big) \big\vert g_{-}^{\ell} \big\vert \big( i \lambda_{1} g_{+}^{k} \cdot g_{-}^{\ell} \big)^{- S ( z , h ) / \lambda_{1}} e^{i \frac{z - E_{0}}{h} t_{+}^{k}} \frac{\CM_{k}^{+}}{\CD_{k} ( t_{+}^{k} )} .
\end{equation*}

We now compute $u$ along the curve $\gamma_{k}$. From \eqref{a4}, we have
\begin{equation} \label{d65}
\left\{ \begin{aligned}
&( P -z ) u = 0 &&\text{microlocally near } \Omega_{\text{reg}} ,  \\
&u = u_{+}^{k} &&\text{microlocally near } U_{+}^{k} ,
\end{aligned} \right.
\end{equation}
where $\Vert u \Vert \leq 1$ and $u_{+}^{k}$ is a semiclassical Lagrangian distribution given by \eqref{d48} and \eqref{d51}. The compact set $\Omega_{\text{reg}}$ is as in Section \ref{s4}. By propagation of Lagrangian distributions (see \cite{MaFe81_01}), the usual transport equations give
\begin{equation} \label{d52}
e^{- i A_{k} / h} e^{i \frac{z - E_{0}}{h} t_{-}^{k}} \frac{\CM_{k}^{-}}{\CD_{k} ( t_{-}^{k} )} a_{-}^{k} ( x , h ) = \CS_{k} ( x , h ) a_{+}^{k} ( \widetilde{x} (x)  , h ) + S ( h^{-N + 1} ) ,
\end{equation}
for all $x$ near $\pi_{x} ( U_{-}^{k} \cap \Lambda_{+}^{1} )$. Here, $\widetilde{x} (x) = \pi_{x} ( \exp ( ( t_{-}^{k} - t_{+}^{k} ) H_{p} ) ( x , \nabla \varphi_{+}^{1} (x) ) ) \in \pi_{x} ( U_{+}^{k} \cap \Lambda_{+}^{0} )$ and $\CS_{k} \in S (1)$ is a symbol which satisfies
\begin{equation*}
\CS_{k} ( x_{-}^{k} , h ) = e^{- i \frac{\pi}{2} \nu_{k}} e^{i \frac{z - E_{0}}{h} ( t_{-}^{k} - t_{+}^{k} )} \frac{\CD_{k} ( t_{+}^{k} )}{\CD_{k} ( t_{-}^{k} )} .
\end{equation*}

Combining \eqref{d51} and \eqref{d52}, we finally obtain
\begin{equation} \label{i7}
a_{-}^{k} ( x , h ) = h^{S ( z , h ) / \lambda_{1} - 1 / 2} \sum_{\ell = 1}^{K} \CP_{k , \ell} ( x , h ) a_{-}^{\ell} ( x_{-}^{\ell} , h ) + S ( h^{-N + \zeta} ) ,
\end{equation}
for all $x$ near $\pi_{x} ( U_{-}^{k} \cap \Lambda_{+}^{1} )$, where the symbols $\CP_{k , \ell} \in S (1)$ satisfy
\begin{align}
\CP_{k , \ell} ( x_{-}^{k} , h ) &= e^{i A_{k} / h} e^{- i \frac{z - E_{0}}{h} t_{-}^{k}} \frac{\CD_{k} ( t_{-}^{k} )}{\CM_{k}^{-}} \CS_{k} ( x_{-}^{k} , h ) \widetilde{\CR}_{k , \ell} ( x_{+}^{k} , h )    \nonumber \\
&= e^{i A_{k} / h}\Gamma \big( S ( z , h ) / \lambda_{1} \big) \sqrt{\frac{\lambda_{1}}{2 \pi}} \frac{\CM_{k}^{+}}{\CM_{k}^{-}} e^{- i \frac{\pi}{2} ( \nu_{k} + \frac{1}{2} )} \big\vert g_{-}^{\ell} \big\vert \big( i \lambda_{1} g_{+}^{k} \cdot g_{-}^{\ell} \big)^{- S ( z , h ) / \lambda_{1}}    \nonumber \\
&= \CQ_{k , \ell} ( z , h ) .  \label{d76}
\end{align}
This shows the lemma.
\end{proof}

We come back to the proof of Proposition \ref{d39}. Applying \eqref{d42} with $x = x_{-}^{k}$, we get
\begin{equation} \label{e72}
\big( 1 - h^{S ( z , h ) / \lambda_{1} - 1 / 2} \CQ ( z , h ) \big) a_{-} ( x_{-} , h ) = \CO ( h^{-N + \zeta} ) ,
\end{equation}
where $a_{-} ( x_{-} , h )$ designs the $K$-column of the $a_{-}^{k} ( x_{-}^{k} , h )$. Using now that $z \in \eqref{d40}$ is at distance $h \vert \ln h \vert^{-1}$ of the pseudo-resonances, the previous estimate and Lemma \ref{d13} give
\begin{equation*}
\forall k \in \{ 1 , \ldots , K \} , \qquad a_{-}^{k} ( x_{-}^{k} , h ) = \CO ( h^{-N + \zeta} ) .
\end{equation*}
Applying one more time \eqref{d42} and using that $\vert h^{S ( z , h ) / \lambda_{1} - 1 / 2} \vert \leq e^{C / \lambda_{1}} \lesssim 1$ for $z \in \eqref{d40}$, we then deduce $a_{-}^{k} ( x , h ) \in S ( h^{-N + \zeta} )$ on $\pi_{x} ( U_{-}^{k} \cap \Lambda_{+}^{1} )$. In other words, starting from \eqref{d44}, we have proved that
\begin{equation*}
u_{-}^{k} \in \CI ( \Lambda_{+}^{1} , h^{- N + \zeta} ) \qquad \text{and} \qquad u_{+}^{k} \in \CI ( \Lambda_{+}^{0} , h^{- N + \zeta} ) .
\end{equation*}
Now, the standard bootstrap argument (see the end of Section \ref{s4} or \eqref{d43}) implies that $u_{\pm}^{k} = \CO ( h^{\infty} )$ and then
\begin{equation*}
u = 0 \text{ microlocally near } K ( E_{0} ) ,
\end{equation*}
which finishes the proof of Proposition \ref{d39} thanks to Section \ref{s31} and Proposition \ref{a16}.

\Subsection{Existence of resonances near the pseudo-resonances} \label{s73}

As explained in Section \ref{s36}, we now show that $P$ has at least one resonance near each pseudo-resonance. This is the aim of the following proposition. Note that Theorem \ref{d8} is a direct consequence of Proposition \ref{d39} and of this result.

\begin{proposition}\sl \label{d50}
Let $C , \delta , r > 0$ and assume that $h$ is small enough. For any pseudo-resonance $z$ in the set
\begin{equation} \label{d53}
E_{0} + [ - C h , C  h ] + i \Big[ - \sum_{j = 2}^{n} \frac{\lambda_{j}}{2} h - C \frac{h}{\vert \ln h \vert} , h \Big] \setminus \big( \Gamma (h) + B ( 0 , \delta h ) \big) ,
\end{equation}
the operator $P$ has at least one resonance in $ B ( z , r h \vert \ln h \vert^{- 1} )$.
\end{proposition}

\begin{remark}\sl
This result will be proved by a contradiction argument. Alternatively, we could have given a direct proof using methods similar to Lemma \ref{d24}. We do not follow this path since the indirect proof is much simpler, especially in the continuous setting (see Section \ref{s29}).
\end{remark}

Writing Proposition \ref{d50} with quantifiers, we have to prove
\begin{align}
\forall r > 0 , \quad \exists h_{0} & > 0 , \quad \forall h < h_{0} ,  \nonumber \\
&\forall z \in \res_{0} (P) \cap \eqref{d53} , \qquad \card \Big( \res (P) \cap B \Big( z , 2 r \frac{h}{\vert \ln h \vert} \Big) \Big) \geq 1 . \label{d54}
\end{align}
If \eqref{d54} does not hold true, then there exist $r > 0$, a (decreasing) sequence of $h$ which goes to $0$ and a sequence of $z = z (h) \in \res_{0} (P) \cap \eqref{d53}$ such that
\begin{equation} \label{d55}
P \text{ has no resonance in } B \Big( z , 2 r \frac{h}{\vert \ln h \vert} \Big) .
\end{equation}

We first reduce the problem using some compactness arguments. As in \eqref{d92}, we denote
\begin{equation*}
\sigma = \frac{z - E_{0}}{h}.
\end{equation*}
Remark that Proposition \ref{d9} and the fact that $Q ( z , h )$ is uniformly bounded imply that the pseudo-resonances belong to 
\begin{equation*}
E_{0} + [ - C h , C  h ] - i \sum_{j = 2}^{n} \frac{\lambda_{j}}{2} h + i \frac{h}{\vert \ln h \vert} [ - C , N ] ,
\end{equation*}
for some fixed $N$. This argument was already used in \eqref{d20}--\eqref{d36}. Then, up to the extraction of a subsequence, we can assume that
\begin{equation*}
\sigma \longrightarrow \sigma_{0} = \tau_{0} -  i \sum_{j = 2}^{n} \frac{\lambda_{j}}{2} \qquad \text{as } h \to 0 ,
\end{equation*}
for some fixed $\tau_{0} \in [ - C , C ]$. The same way, since $\rho (h) = ( e^{i A_{1} / h} , \ldots , e^{i A_{K} / h} )$ belongs to the compact set $( \S^{1} )^{K}$, we can assume that $\rho (h) = \rho_{0} + o (1)$ where $\rho_{0}$ is independent of $h$. In particular, \eqref{d10} gives
\begin{equation} \label{d57}
\widehat{\CQ} ( \tau_{0} , h ) = \CQ_{0} + o (1) ,
\end{equation}
where $\CQ_{0} = \widetilde{\CQ} ( \rho_{0} , \sigma_{0} )$. From the continuity of the map which, at a $K \times K$ matrix, associates its spectrum (see \eqref{d26} for a similar argument), we obtain
\begin{equation} \label{m82}
\dist \big( \spe ( \widehat{\CQ} ( \tau_{0} , h ) ) , \spe ( \CQ_{0}) \big) = o (1) .
\end{equation}
Then, Proposition \ref{d9} implies the existence of $k = k (h) \in \{ 1 , \ldots , K \}$ and $q = q (h) \in \Z$ such that $z = z_{q , k}^{0} + o ( h \vert \ln h \vert^{- 1} )$ where
\begin{equation*}
z_{q , k}^{0} = E_{0} + 2 q \pi \lambda_{1} \frac{h}{\vert \ln h \vert} - i h \sum_{j = 2}^{n} \frac{\lambda_{j}}{2} + i \ln ( \mu_{k} ) \lambda_{1} \frac{h}{\vert \ln h \vert} ,
\end{equation*}
and the $\mu_{k}$'s are the $K$ eigenvalues of $\CQ_{0}$. Since $k$ is in a finite set, we can assume that $k = k_{0}$ is a constant and we denote $\mu_{0} = \mu_{k_{0}}$. Summing up, there exists an eigenvalue $\mu_{0}$ of $\CQ_{0}$ such that $z = z_{q}^{0} + o ( h \vert \ln h \vert^{-1} )$ with
\begin{equation} \label{d58}
z_{q}^{0} = E_{0} + 2 q \pi \lambda_{1} \frac{h}{\vert \ln h \vert} - i h \sum_{j = 2}^{n} \frac{\lambda_{j}}{2} + i \ln ( \mu_{0} ) \lambda_{1} \frac{h}{\vert \ln h \vert} ,
\end{equation}
and $\re z = \tau_{0} + o (h)$. Moreover, since $z \in \eqref{d53}$, we necessarily have
\begin{equation} \label{d81}
\mu_{0} \neq 0 .
\end{equation}
In other words, we have fixed all the parameters of the sequence of pseudo-resonances $z = z (h)$ excepted $q \in \Z$ (which is however limited by the condition $\re z_{q}^{0} = \tau_{0} h + o (h)$). Of course, \eqref{d55} implies
\begin{equation} \label{d59}
P \text{ has no resonance in } B \Big( z_{q}^{0} , r \frac{h}{\vert \ln h \vert} \Big) .
\end{equation}

\begin{figure}
\begin{center}
\begin{picture}(0,0)%
\includegraphics{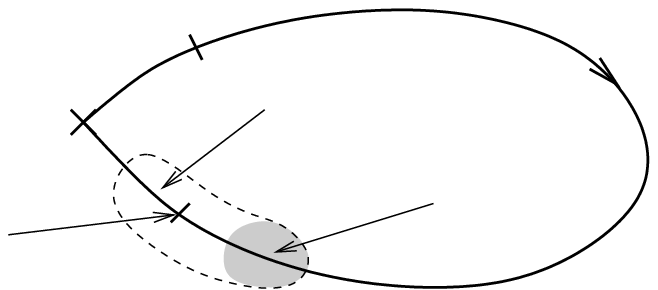}%
\end{picture}%
\setlength{\unitlength}{1579sp}%
\begingroup\makeatletter\ifx\SetFigFont\undefined%
\gdef\SetFigFont#1#2#3#4#5{%
  \reset@font\fontsize{#1}{#2pt}%
  \fontfamily{#3}\fontseries{#4}\fontshape{#5}%
  \selectfont}%
\fi\endgroup%
\begin{picture}(8322,3396)(2086,-5975)
\put(3901,-4036){\makebox(0,0)[lb]{\smash{{\SetFigFont{9}{10.8}{\rmdefault}{\mddefault}{\updefault}$0$}}}}
\put(7951,-5011){\makebox(0,0)[lb]{\smash{{\SetFigFont{9}{10.8}{\rmdefault}{\mddefault}{\updefault}$MS ( v_{k} )$}}}}
\put(6001,-3886){\makebox(0,0)[lb]{\smash{{\SetFigFont{9}{10.8}{\rmdefault}{\mddefault}{\updefault}$MS ( \Op ( \chi_{k} ) \widehat{v}_{k} )$}}}}
\put(10126,-3136){\makebox(0,0)[lb]{\smash{{\SetFigFont{9}{10.8}{\rmdefault}{\mddefault}{\updefault}$\gamma_{k}$}}}}
\put(2101,-5386){\makebox(0,0)[lb]{\smash{{\SetFigFont{9}{10.8}{\rmdefault}{\mddefault}{\updefault}$\rho_{-}^{k}$}}}}
\put(4426,-2836){\makebox(0,0)[lb]{\smash{{\SetFigFont{9}{10.8}{\rmdefault}{\mddefault}{\updefault}$\rho_{+}^{k}$}}}}
\end{picture}%
\end{center}
\caption{The microsupports of $\Op ( \chi_{k} ) \widehat{v}_{k}$ and $v_{k}$.} \label{f22}
\end{figure}

We now construct a test function related to the previous reductions. For that, we use the same notations $\Omega_{\rm{sing}}$, $\Omega_{\text{reg}}$, $\varepsilon$, $S_{\pm}^{\varepsilon}$, $\rho_{\pm}^{k} = ( x_{\pm}^{k} , \xi_{\pm}^{k} )$, $t_{\pm}^{k}$ and $U_{\pm}^{k}$ as in Section \ref{s72}. Let $w_{0} = ( w_{0}^{1} , \ldots , w_{0}^{K} ) \in \C^{K}$ be a fixed (non-zero) eigenvector of $\CQ_{0}$ associated to the eigenvalue $\mu_{0}$. Let $\widetilde{v}_{k}$ be a usual WKB solution of
\begin{equation} \label{d89}
\left\{ \begin{aligned}
&( P - \widetilde{z} ) \widetilde{v}_{k} = 0 &&\text{near } x_{-}^{k} ,   \\
&\widetilde{v}_{k} (x) = w_{0}^{k} e^{i \varphi_{+}^{1} (x) / h} &&\text{on } \vert x \vert = \varepsilon \text{ near } x_{-}^{k} ,
\end{aligned} \right.
\end{equation}
holomorphic in $\widetilde{z} \in B ( E_{0} , R h )$ for some large $R > 1$. Then, $\widetilde{v}_{k} (x) = \widetilde{a}_{k} ( x , h ) e^{i \varphi_{+}^{1} (x) / h}$ is a Lagrangian distribution where the symbol $\widetilde{a}_{k} \in S (1)$ has an asymptotic expansions $\widetilde{a}_{k} ( x , h ) \simeq \widetilde{a}_{k}^{0} (x) + \widetilde{a}_{k}^{1} (x) h + \cdots$ in $S (1)$. Consider now cut-off functions $\chi_{k} , \psi_{k} \in C^{\infty}_{0} ( T^{*} \R^{n} )$ such that $\chi_{k} = 1$ near $U_{-}^{k}$, $\psi_{k} = 1$ near $\supp ( \nabla \chi_{k} ) \cap \gamma_{k} ( ] - \infty , t_{-}^{k} [ )$ and such that $\chi_{k}$ and $\psi_{k}$ are supported in the region where $\widetilde{v}_{k}$ is defined (see Figure \ref{f22}). Then we set
\begin{equation}
\widehat{v}_{k} (x) = e^{- i A_{k} / h} e^{i \frac{\widetilde{z} - E_{0}}{h} t_{-}^{k}} \frac{\CM_{k}^{-}}{\CD_{k} ( t_{-}^{k} )} \widetilde{v}_{k} (x) ,
\end{equation}
and eventually
\begin{equation} \label{d60}
v = \sum_{k = 1}^{K} v_{k} \qquad \text{with} \qquad v_{k} = \Op ( \psi_{k} ) \big[ P , \Op ( \chi_{k} ) \big] \widehat{v}_{k} .
\end{equation}
The function $v$ will be our ``test function''. Note that $v$ is holomorphic with respect to $z$ in the set \eqref{d53}.

Let $0 < s < r$ be a small enough number which will be fixed in the sequel and denote $\CD = B ( z_{q}^{0} , s h \vert \ln h \vert^{-1} )$. Thus, $\partial \CD$ is a positively oriented path around $z_{q}^{0}$. Let $u \in L^{2} ( \R^{n} )$ be the solution of
\begin{equation} \label{d61}
( P_{\theta} - \widetilde{z} ) u = v ,
\end{equation}
for $\widetilde{z} \in \partial \CD$. From \eqref{m82}, \eqref{d58} and Proposition \ref{d9}, $\partial \CD$ is at distance $h \vert \ln h \vert^{- 1}$ from the pseudo-resonances for $s > 0$ fixed small enough. Then, Proposition \ref{d39} shows that $u = ( P_{\theta} - \widetilde{z} )^{-1} v$ is well-defined, holomorphic near $\partial \CD$ and $\Vert u \Vert \lesssim h^{- M}$ uniformly for $\widetilde{z} \in \partial \CD$. We now compute $u$ following Section \ref{s3}, Section \ref{s4} and Section \ref{s72}. By construction, \eqref{d60} implies that $v \in \CI ( \Lambda_{+}^{1} , 1 )$ with $\Lambda_{+}^{1} \subset \Lambda_{+} \subset p^{-1} ( E_{0} )$. Then, the energy localization (Lemma \ref{a12}), the vanishing in the incoming region (Lemma \ref{a20}), the propagation of singularities (Lemma \ref{a21} with the additional assumption $\exp ( [ 0 , T ] H_{p} ) ( \rho_{0} ) \cap \supp ( v ) \neq \emptyset$) and the vanishing outside $\CH$ (Lemma \ref{a33}) still hold true uniformly for $\widetilde{z} \in \partial \CD$.

As usual, we define $u_{\pm}^{k}$ as the microlocal restriction of $u$ to a neighborhood of $U_{\pm}^{k}$. Since the microsupport of $v$ does not meet $\Omega_{\rm{sing}}$, \eqref{d61} gives
\begin{equation} \label{d63}
\left\{ \begin{aligned}
&( P - \widetilde{z} ) u = 0 &&\text{microlocally near } \Omega_{\rm{sing}} ,   \\
&u = u_{-}^{k} &&\text{microlocally near } U_{-}^{k} ,   \\
&u = 0 &&\text{microlocally near } S_{-}^{\varepsilon} \setminus \big( U_{-}^{1} \cup \cdots \cup U_{-}^{K} \big) ,
\end{aligned} \right.
\end{equation}
and $\Vert u \Vert \lesssim h^{- M}$. Note that this microlocal Cauchy problem is similar to \eqref{a41} or \eqref{d62}. Moreover, we remark that $\widetilde{z} \in \partial \CD$ is at distance $h$ of $\Gamma (h)$ since $\widetilde{z} = z + \CO ( h \vert \ln h \vert^{-1} )$ with $z \in \eqref{d53}$. Thus, we can apply Theorem \ref{a32} to solve \eqref{d63} and we obtain, as in \eqref{a36} or \eqref{d44},
\begin{equation} \label{d66}
u_{k}^{+} \in \CI ( \Lambda_{+}^{0} , h^{-N} ) ,
\end{equation}
for some $N \in \R$ uniformly for $\widetilde{z} \in \partial \CD$.

On the other hand, along the regular part $\gamma_{k}$, $u$ satisfies the microlocal Cauchy problem
\begin{equation} \label{d64}
\left\{ \begin{aligned}
&( P - \widetilde{z} ) u = v_{k} &&\text{microlocally near } \Omega_{\text{reg}} ,  \\
&u = u_{+}^{k} &&\text{microlocally near } U_{+}^{k} ,
\end{aligned} \right.
\end{equation}
with $\Vert u \Vert \lesssim h^{- M}$. By linearity, $u$ is the sum of the solutions of
\begin{equation} \label{d67}
\left\{ \begin{aligned}
&( P - \widetilde{z} ) u_{\rm{usual}}^{k} = 0 &&\text{microlocally near } \Omega_{\text{reg}} ,  \\
&u_{\rm{usual}}^{k} = u_{+}^{k} &&\text{microlocally near } U_{+}^{k} ,
\end{aligned} \right.
\end{equation}
and
\begin{equation} \label{d68}
\left\{ \begin{aligned}
&( P - \widetilde{z} ) u_{\rm{test}}^{k} = v_{k} &&\text{microlocally near } \Omega_{\text{reg}} ,  \\
&u_{\rm{test}}^{k} = 0 &&\text{microlocally near } U_{+}^{k} .
\end{aligned} \right.
\end{equation}
Remark that \eqref{d67} is similar to \eqref{a37} and \eqref{d65}. In particular, as in \eqref{a38} or \eqref{d44}, the usual propagation of Lagrangian distribution yields
\begin{equation} \label{d69}
u_{\rm{usual}}^{k} \in \CI ( \Lambda_{+}^{1} , h^{- N} ) \text{ microlocally near } U_{-}^{k} ,
\end{equation}
uniformly for $\widetilde{z} \in \partial \CD$. Furthermore, one can check that
\begin{equation} \label{d70}
u_{\rm{test}}^{k} = \Op ( \chi_{k} ) \widehat{v}_{k} \in \CI ( \Lambda_{+}^{1} , 1 ) ,
\end{equation}
is the explicit solution of \eqref{d68} microlocally near $\Omega_{\text{reg}}$. Combining the two last equations, we deduce
\begin{equation} \label{d71}
u_{-}^{k} = u_{\rm{usual}}^{k} + u_{\rm{test}}^{k} \in \CI \big( \Lambda_{+}^{1} , \max ( 1 , h^{- N} ) \big) ,
\end{equation}
uniformly for $\widetilde{z} \in \partial \CD$. Then, as in \eqref{d48}, there exist symbols $a_{\pm}^{k} \in S ( \max ( 1 , h^{- N} ) )$ defined near $\pi_{x} ( U_{-}^{k} \cap \Lambda_{+}^{1} )$ and $\pi_{x} ( U_{+}^{k} \cap \Lambda_{+}^{0} )$ such that
\begin{equation} \label{d87}
\left\{ \begin{aligned}
u_{-}^{k} (x) &= e^{- i A_{k} / h} e^{i \frac{\widetilde{z} - E_{0}}{h} t_{-}^{k}} \frac{\CM_{k}^{-}}{\CD_{k} ( t_{-}^{k} )} a_{-}^{k} ( x , h ) e^{i \varphi_{+}^{1} (x) / h} ,  \\
u_{+}^{k} (x) &= a_{+}^{k} ( x , h ) e^{i \varphi_{+} (x) / h} .
\end{aligned} \right.
\end{equation}
Making a turn along $\CH$, these symbols satisfies the following relation.

\begin{lemma}\sl \label{d72}
Uniformly for $\widetilde{z} \in \partial \CD$, we have
\begin{equation} \label{d73}
a_{-}^{k} ( x , h ) = h^{S ( \widetilde{z} / h ) / \lambda_{1} - 1 / 2} \sum_{\ell = 1}^{K} \CP_{k , \ell} ( x , h ) a_{-}^{\ell} ( x_{-}^{\ell} , h ) + \widetilde{a}_{k} ( x , h ) + S ( h^{-N + \zeta} ) ,
\end{equation}
where $\zeta > 0$ and the symbols $\CP_{k , \ell}$ are those of Lemma \ref{d41}.
\end{lemma}

\begin{proof}
This proof is similar to the one of Lemma \ref{d41} and we will use some of its estimates. We first remark that the microlocal Cauchy problem \eqref{d63} is exactly the same as \eqref{d62}. Thus, we obtain
\begin{equation} \label{d74}
a_{+}^{k} ( x , h ) = h^{S ( \widetilde{z} / h ) / \lambda_{1} - 1 / 2} \sum_{\ell = 1}^{K} \widetilde{\CR}_{k , \ell} ( x , h ) a_{-}^{\ell} ( x_{-}^{\ell} , h ) + S ( h^{-N + \zeta} ) ,
\end{equation}
as in \eqref{d51}.

As explained previously, the microlocal Cauchy problem \eqref{d64} can be decomposed in \eqref{d67} and \eqref{d68}. The microlocal Cauchy problem \eqref{d67} is exactly \eqref{d65}. Thus, \eqref{d52} holds true for the symbol of the Lagrangian distribution $u_{\rm{usual}}^{k}$. Moreover, $u_{\rm{test}}^{k}$ is given by \eqref{d70}. Summing up and using that $\chi_{k} = 1$ near $U_{-}^{k}$, it yields
\begin{align}
e^{- i A_{k} / h} e^{i \frac{\widetilde{z} - E_{0}}{h} t_{-}^{k}} \frac{\CM_{k}^{-}}{\CD_{k} ( t_{-}^{k} )} a_{-}^{k} ( x , h ) ={}& \CS_{k} ( x , h ) a_{+}^{k} ( \widetilde{x} (x) , h ) \nonumber  \\
&+ e^{- i A_{k} / h} e^{i \frac{\widetilde{z} - E_{0}}{h} t_{-}^{k}} \frac{\CM_{k}^{-}}{\CD_{k} ( t_{-}^{k} )} \widetilde{a}_{k} ( x , h ) + S ( h^{-N + 1} ) .   \label{d75}
\end{align}
The lemma follows from \eqref{d74}, \eqref{d75} and the computation made in \eqref{d76}.
\end{proof}

As in Section \ref{s13}, we define
\begin{equation*}
\Lambda = \Lambda ( \widetilde{z} , h ) = h^{-S ( \widetilde{z} / h ) / \lambda_{1} + 1 / 2} ,
\end{equation*}
which satisfies $1 \lesssim \vert \Lambda \vert = h^{- \sum_{j=2}^{n} \frac{\lambda_{j}}{2 \lambda_{1}} - \frac{\im \widetilde{z}}{\lambda_{1} h}} \lesssim 1$ uniformly for $\widetilde{z} \in \partial \CD$. We now fix $s$ such that the new spectral parameter $\Lambda$ avoids $\spe ( \CQ_{0} ) \setminus \{ \mu_{0} \}$ on $\partial \CD$.

\begin{lemma}\sl \label{d80}
For $0 < s < r$ (fixed) small enough, $\Lambda ( \partial \CD )$ is a simple loop around $\mu_{0}$, the rest of the spectrum of $\CQ_{0}$ is outside of $\Lambda ( \CD )$ and all the spectrum of $\CQ_{0}$ is at distance larger than some constant $\alpha > 0$ of $\Lambda ( \partial \CD )$.
\end{lemma}

\begin{proof}
This result is close to Lemma \ref{d24}. For $\widetilde{z} = z_{q}^{0} + \rho s \frac{h}{\vert \ln h \vert}$ with $\vert \rho \vert = 1$, we have
\begin{equation*}
\Lambda = \mu_{0} e^{- i \rho s / \lambda_{1}} = \mu_{0} - i \rho s \frac{\mu_{0}}{\lambda_{1}} + o_{s \to 0} (s) ,
\end{equation*}
and the lemma follows from \eqref{d81} since the matrix $\CQ_{0}$ has a discrete spectrum.
\end{proof}

We then prove that $N = 0$.

\begin{lemma}\sl \label{d77}
We have
\begin{equation*}
u_{-}^{k} \in \CI ( \Lambda_{+}^{1} , 1) \qquad \text{and} \qquad u_{+}^{k} \in \CI ( \Lambda_{+}^{0} , 1 ) ,
\end{equation*}
uniformly for $\widetilde{z} \in \partial \CD$.
\end{lemma}

\begin{proof}
Indeed, applying \eqref{d73} with $x = x_{-}^{k}$, we get
\begin{equation} \label{d82}
\big( 1 - h^{S ( \widetilde{z} / h ) / \lambda_{1} - 1 / 2} \CQ ( \widetilde{z} , h ) \big) a_{-} ( x_{-} , h ) = w_{0} + \CO ( h^{-N + \zeta} ) ,
\end{equation}
where $a_{-} ( x_{-} , h )$ designs the $K$-column of the $a_{-}^{k} ( x_{-}^{k} , h )$. Note here that $\widetilde{a}_{k} ( x_{-}^{k} , h ) = w_{0}^{k}$ by construction (see \eqref{d89}). Using $\Lambda = \CO (1)$, this yields
\begin{equation} \label{d78}
\Vert a_{-} ( x_{-} , h ) \Vert \lesssim \big\Vert \big( \Lambda - \CQ ( \widetilde{z} , h ) \big)^{-1} \big\Vert \big( \CO (1) + \CO ( h^{-N + \zeta} ) \big) .
\end{equation}
Since $\widetilde{\sigma} : = ( \widetilde{z} - E_{0} ) / h= \sigma_{0} + o (1)$ and $\rho = \rho_{0} + o (1)$, \eqref{d10} gives
\begin{equation} \label{d86}
\CQ ( \widetilde{z} , h ) = \CQ_{0} + o (1) .
\end{equation}
Then, the continuity of the spectrum with respect to the $K \times K$ matrix implies
\begin{equation*}
\dist \big( \spe ( \CQ ( \widetilde{z} , h ) ) , \spe ( \CQ_{0}) \big) = o (1) .
\end{equation*}
Combining with Lemma \ref{d80}, we deduce $\dist ( \Lambda , \CQ ( \widetilde{z} , h ) ) \geq \alpha / 2$ for all $\widetilde{z} \in \partial \CD$ and $h$ small enough. Then, Lemma \ref{d14} yields
\begin{equation} \label{d83}
\big\Vert \big( \Lambda - \CQ ( \widetilde{z} , h ) \big)^{-1} \big\Vert \leq M ,
\end{equation}
and \eqref{d78}  becomes
\begin{equation*}
\forall k \in \{ 1 , \ldots , K \} , \qquad a_{-}^{k} ( x_{-}^{k} , h ) = \CO (1) + \CO ( h^{-N + \zeta} ) ,
\end{equation*}
uniformly for $\widetilde{z} \in \partial \CD$. Applying one more time \eqref{d73}, this estimate implies $a_{-}^{k} \in S ( h^{-N + \zeta} ) + S (1)$ starting from $a_{-}^{k} \in S ( h^{-N} )$. The lemma follows then by the usual bootstrap argument, limited by the class of symbols $S (1)$.
\end{proof}

\begin{lemma}\sl \label{d79}
We have
\begin{equation} \label{d84}
\frac{1}{2 i \pi} \int_{\partial \CD} u ( \widetilde{z} ) \, d \widetilde{z} = b ( x , h ) e^{i \varphi_{+}^{1} (x) / h} \text{ microlocally near } U_{-}^{k} ,
\end{equation}
where $b \in S ( h \vert \ln h \vert^{-1} )$ satisfies
\begin{equation*}
b ( x_{-}^{k} , h ) = i e^{- i A_{k} / h} e^{i \sigma_{0} t_{-}^{k}} \frac{\lambda_{1} \CM_{k}^{-}}{\CD_{k} ( t_{-}^{k} )} w_{0}^{k} \frac{h}{\vert \ln h \vert} + o \Big( \frac{h}{\vert \ln h \vert} \Big) .
\end{equation*}
\end{lemma}

\begin{proof}
Formula \eqref{d87} implies \eqref{d84} with
\begin{equation} \label{d88}
b ( x , h ) = e^{- i A_{k} / h} \frac{\CM_{k}^{-}}{\CD_{k} ( t_{-}^{k} )} \frac{1}{2 i \pi} \int_{\partial \CD} e^{i \frac{\widetilde{z} - E_{0}}{h} t_{-}^{k}} a_{-}^{k} ( x , h ) \, d \widetilde{z} ,
\end{equation}
near $\pi_{x} ( U_{-}^{k} \cap \Lambda_{+}^{1} )$. Since $a_{-}^{k} \in S (1)$ uniformly for $\widetilde{z} \in\partial \CD$ from Lemma \ref{d77} and $\vert \partial \CD \vert$ is of size $h \vert \ln h \vert^{-1}$, we deduce $b \in S ( h \vert \ln h \vert^{-1} )$. It remains to compute $b ( x_{-}^{k} , h )$.

Using Lemma \ref{d77} and \eqref{d83}, \eqref{d82} gives
\begin{equation} \label{d85}
a_{-} ( x_{-} , h ) = \Lambda \big( \Lambda - \CQ ( \widetilde{z} , h ) \big)^{- 1} w_{0} + \CO ( h^{\zeta} ) ,
\end{equation}
uniformly for $\widetilde{z} \in \partial \CD$. On the other hand, Lemma \ref{d80} also implies
\begin{equation*}
\Vert ( \Lambda - \CQ_{0} )^{- 1} \Vert \lesssim 1 ,
\end{equation*}
uniformly for $\widetilde{z} \in \partial \CD$. Then, the resolvent identity and \eqref{d86} give
\begin{align*}
\big( \Lambda - \CQ ( \widetilde{z} , h ) \big)^{- 1} ={}& ( \Lambda  - \CQ_{0} )^{- 1} + \big( \Lambda - \CQ ( \widetilde{z} , h ) \big)^{- 1} \big( \CQ ( \widetilde{z} , h ) - \CQ_{0} \big) ( \Lambda  - \CQ_{0} )^{- 1}   \\
={}& ( \Lambda  - \CQ_{0} )^{- 1} + o (1) ,
\end{align*}
and \eqref{d85} becomes
\begin{equation}
a_{-} ( x_{-} , h ) = \Lambda ( \Lambda - \CQ_{0} )^{- 1} w_{0} + o (1) .
\end{equation}
Moreover, $\widetilde{\sigma} = \sigma_{0} + o (1)$ shows
\begin{equation*}
e^{i \frac{\widetilde{z} - E_{0}}{h} t_{-}^{k}} = e^{i \sigma_{0} t_{-}^{k}} + o (1) .
\end{equation*}
Thus, using the change of variable $\widetilde{z} \mapsto \Lambda$ which verifies $d \widetilde{z} = \frac{i \lambda_{1} h}{\vert \ln h \vert \Lambda} d \Lambda$, \eqref{d88} yields
\begin{align}
b ( x_{-}^{k} , h ) &= e^{- i A_{k} / h} \frac{\CM_{k}^{-}}{\CD_{k} ( t_{-}^{k} )} e^{i \sigma_{0} t_{-}^{k}} \frac{1}{2 i \pi} \int_{\partial \CD} \Big[ \Lambda ( \Lambda - \CQ_{0} )^{- 1} w_{0} \Big]_{k} \, d \widetilde{z} + o \Big( \frac{h}{\vert \ln h \vert} \Big)   \nonumber \\
&= e^{- i A_{k} / h} \frac{\CM_{k}^{-}}{\CD_{k} ( t_{-}^{k} )} e^{i \sigma_{0} t_{-}^{k}} \frac{i \lambda_{1} h}{\vert \ln h \vert} \frac{1}{2 i \pi} \bigg[ \int_{\Lambda ( \partial \CD )} ( \Lambda - \CQ_{0} )^{- 1} w_{0} \, d \Lambda \bigg]_{k} + o \Big( \frac{h}{\vert \ln h \vert} \Big)  \nonumber \\
&= e^{- i A_{k} / h} \frac{\CM_{k}^{-}}{\CD_{k} ( t_{-}^{k} )} e^{i \sigma_{0} t_{-}^{k}} \frac{i \lambda_{1} h}{\vert \ln h \vert} w_{0}^{k} + o \Big( \frac{h}{\vert \ln h \vert} \Big) ,
\end{align}
since $w_{0}$ is an eigenvector of $\CQ_{0}$ associated to the eigenvalue $\mu_{0}$. Here, $[ y ]_{k}$ denotes the $k$-th coordinate of $y \in \R^{K}$. This proves the lemma.
\end{proof}

Eventually, the last result allows us to conclude. Since $v ( \widetilde{z} )$ is holomorphic in \eqref{d53}, $u = ( P_{\theta} - \widetilde{z} )^{-1} v$ and \eqref{d59} imply that $u ( \widetilde{z} )$ is holomorphic near $\CD$. In particular,
\begin{equation*}
\frac{1}{2 i \pi} \int_{\partial \CD} u ( \widetilde{z} ) \, d \widetilde{z} = 0 ,
\end{equation*}
which is in contradiction with Lemma \ref{d79}. This implies Proposition \ref{d50}.

\Subsection{Proof of Proposition \ref{d91}}

In this part, we explain how to adapt the proof of Proposition \ref{d50} in order to obtain this result. Since Proposition \ref{d91} implies Proposition \ref{d50}, we could have directly shown the first proposition. But, we have chosen this presentation for the sake of clarity.

As in \eqref{d55}, if Proposition \ref{d91} does not hold true, there exist $s > 0$, a sequence of $h$ which goes to $0$ and a sequence of $z = z (h) \in \res_{0} ( P) \cap \eqref{d90}$ such that
\begin{equation*}
\card \Big( \res (P) \cap B \Big( z , 2 s \frac{h}{\vert \ln h \vert} \Big) \Big) < \card \Big\{ ( q , k ) \in \Z \times \{ 1 , \ldots , K \} ; \ z_{q , k} ( \re \sigma ) \in B \Big( z , s \frac{h}{\vert \ln h \vert} \Big) \Big\} .
\end{equation*}
Then we can fix the various parameters as in \eqref{d59} and assume, up to the extraction of a subsequence, that
\begin{align}
\card \Big( \res (P) \cap B &\Big( z_{q}^{0} , 2 s \frac{h}{\vert \ln h \vert} \Big) \Big)  \nonumber \\
&< \card \Big\{ ( \widehat{q} , k ) \in \Z \times \{ 1 , \ldots , K \} ; \ z_{\widehat{q} , k}^{0} \in B \Big( z_{q}^{0} , r s \frac{h}{\vert \ln h \vert} \Big) \Big\} .  \label{d96}
\end{align}
where $z_{q}^{0}$ is defined in \eqref{d58},
\begin{equation*}
z_{\widehat{q} , k}^{0} = E_{0} + 2 \widehat{q} \pi \lambda_{1} \frac{h}{\vert \ln h \vert} - i h \sum_{j = 2}^{n} \frac{\lambda_{j}}{2} + i \ln ( \mu_{k} ) \lambda_{1} \frac{h}{\vert \ln h \vert} ,
\end{equation*}
with $\mu_{k} \in \spe ( \CQ_{0} )$, and $1 < r < 2$ is fixed such that $\partial B ( z_{q}^{0} , r s h \vert \ln h \vert^{- 1} )$ avoids the $z_{\widehat{q} , k}^{0}$'s. In particular, the right hand side of \eqref{d96} is a constant, denoted $N_{0} > 1$ in the sequel. Moreover, if $s$ is small enough, there exists, for each $k \in \{ 1 , \ldots , K \} $, at most one $\widehat{q} \in \Z$ such that $z_{\widehat{q} , k}^{0}$ belongs to $\CD : = B ( z_{q}^{0} , r s h \vert \ln h \vert^{- 1} )$. Now, let
\begin{equation*}
\Pi_{k} = - \frac{1}{2 i \pi} \oint_{\gamma_{k}} ( \CQ_{0} - z )^{- 1} d z ,
\end{equation*}
where $\gamma_{k}$ is a simple loop around $\mu_{k}$, the spectral projection of $\CQ_{0}$ associated to its eigenvalue $\mu_{k}$. We defined $\Pi_{0} = \sum_{k} \Pi_{k}$ where the sum is restricted to the $k$'s which contribute to the right hand side of \eqref{d96}. Thus, $\rank \Pi_{0} = N_{0}$. Let $w_{0 , \ell} \in \R^{K}$, $\ell = 1 , \ldots , N_{0}$, be a basis of $\im \Pi_{0}$.

Following \eqref{d89}--\eqref{d60}, we construct a test function $v_{\ell}$ for all $\ell = 1 , \ldots , N_{0}$. Mimicking \eqref{d61}, we then consider the problem
\begin{equation} \label{d97}
( P_{\theta} - \widetilde{z} ) u_{\ell} = v_{\ell} ,
\end{equation}
for $\widetilde{z} \in \partial \CD$. As in Lemma \ref{d79}, we have
\begin{equation} \label{d98}
\frac{1}{2 i \pi} \int_{\partial \CD} u_{\ell} ( \widetilde{z} ) \, d \widetilde{z} = b_{\ell} ( x , h ) e^{i \varphi_{+}^{1} (x) / h} \text{ microlocally near } U_{-}^{k} ,
\end{equation}
where $b_{\ell} \in S ( h \vert \ln h \vert^{-1} )$ satisfies
\begin{equation*}
b_{\ell} ( x_{-}^{k} , h ) = i e^{- i A_{k} / h} e^{i \sigma_{0} t_{-}^{k}} \frac{\lambda_{1} \CM_{k}^{-}}{\CD_{k} ( t_{-}^{k} )} w_{0 , \ell}^{k} \frac{h}{\vert \ln h \vert} + o \Big( \frac{h}{\vert \ln h \vert} \Big) .
\end{equation*}
The proof of \eqref{d98} is the same as the one of \eqref{d84} except that $w_{0 , \ell}$ is now a generalized eigenvector of $\CQ_{0}$.

From \eqref{d98} and the independence of the $w_{0 , \ell}$, we deduce that the $N_{0}$ functions
\begin{equation*}
\frac{1}{2 i \pi} \int_{\partial \CD} u_{\ell} ( \widetilde{z} ) \, d \widetilde{z} ,
\end{equation*}
are independent in $L^{2} ( \R^{n} )$. On the other hand, since the $v_{\ell}$'s are holomorphic, \eqref{d97} implies that the above functions belong to the image of the sum of the spectral projections of $P_{\theta}$ associated to resonances lying in $\CD$. These two remarks yield
\begin{equation}
\card \Big( \res (P) \cap B \Big( z_{q}^{0} , r s \frac{h}{\vert \ln h \vert} \Big) \Big) \geq N_{0} ,
\end{equation}
which is in contradiction with \eqref{d96}.

\section{Proof of the other results of Section \ref{s6}} \label{s23}

This part is devoted to the demonstration of the statements and results of Section \ref{s6} which are not proved in Section \ref{s12}. We begin with the

\begin{proof}[Proof of Remark \ref{d3} $i)$]
Let $\widehat{\gamma}_{k} (t) = \gamma_{k} ( t + t_{k} )$ be a new time parametrization of $\gamma_{k}$. From \eqref{d6} and \eqref{d7}, we deduce $\widehat{g}^{k}_{\pm} = g^{k}_{\pm} e^{\pm \lambda_{1} t_{k}}$,
\begin{equation} \label{d49}
\widehat{\CM}_{k}^{+} = \CM_{k}^{+} e^{t_{k} \sum_{j} \lambda_{j} / 2} \qquad \text{and} \qquad \widehat{\CM}_{k}^{-} = \CM_{k}^{-} e^{t_{k} ( \sum_{j} \lambda_{j} - 2 \lambda_{1} ) / 2} .
\end{equation}
Then, \eqref{d4} implies
\begin{equation*}
\widehat{\CQ}_{k , \ell} = e^{( \lambda_{1} - S ( z , h ) ) t_{k}} \CQ_{k , \ell} e^{- ( \lambda_{1} - S ( z , h ) ) t_{\ell}} .
\end{equation*}
In other words, using the $K \times K$ matrix $U = \diag ( e^{( \lambda_{1} - S ( z , h ) ) t_{k}} )$, we have $\widehat{\CQ} = U \CQ U^{- 1}$. Therefore, $\CQ$ and $\widehat{\CQ}$ have the same eigenvalues.
\end{proof}

\Subsection{Proof of the stability results} \label{s25}

We prove here the assertions of Section \ref{s24}.

\begin{proof}[Proof of Proposition \ref{i38}]
If $K = K_{0}$, this result follows directly from \eqref{d4} (which shows that $\CQ_{\bullet}$ only depends on the values of the symbol of $P_{\bullet}$ near $\CH_{\bullet}$), Proposition \ref{d9} and Theorem \ref{d8}. We now assume $K > K_{0}$. Using the $K \times K$ matrix  $U = \diag ( \vert g_{-}^{k} \vert^{1 / 2} )$, we set
\begin{equation} \label{i43}
\overline{\CQ} ( \tau , h ) : = U \CQ \Big( E_{0} + h \tau - i h \sum_{j = 2}^{n} \frac{\lambda_{j}}{2} , h \Big) U^{- 1} .
\end{equation}
The same way, we define the $K_{0} \times K_{0}$ matrix $\overline{\CQ}_{0} ( \tau , h )$. The spectrum of $\CQ_{\bullet}$ coincides with the spectrum of $\overline{\CQ}_{\bullet}$. From Proposition \ref{d9} and Theorem \ref{d8}, it is enough to show that, in the domain $\{ \lambda \in \C ; \ \vert \lambda \vert \geq e^{- C / 2 \lambda_{1}} \}$, we have
\begin{equation} \label{i39}
\dist \big( \spe \big( \overline{\CQ} ( \tau , h ) \big) , \spe \big( \overline{\CQ}_{0} ( \tau , h ) \big) \big) = o_{\varepsilon \to 0} ( 1 ) ,
\end{equation}
uniformly for $\tau \in [ - C , C ]$ and $h$ small enough. Using the explicit form of the $\CQ_{\bullet}$'s (see \eqref{d4}) and the ``smallness of the perturbation'' (i.e. \eqref{i40}), we can write
\begin{equation*}
\overline{\CQ} ( \tau , h ) = \left( \begin{array}{cc} 
\overline{\CQ}_{0} ( \tau , h ) & \CO ( 1 ) \\
\CO ( \varepsilon ) & \CO ( \varepsilon )
\end{array} \right) = \left( \begin{array}{cc} 
\overline{\CQ}_{0} ( \tau , h ) & \CO ( 1 ) \\
0 & 0
\end{array} \right) + \CO ( \varepsilon ) .
\end{equation*}
Moreover, we have
\begin{equation*}
\spe \left( \begin{array}{cc} 
\overline{\CQ}_{0} ( \tau , h ) & \CO ( 1 ) \\
0 & 0
\end{array} \right) = \spe \big( \overline{\CQ}_{0} ( \tau , h ) \big) \cup \{ 0 \} .
\end{equation*}
Then, \eqref{i39} follows from the previous relations and the uniform continuity on bounded regions of the map which, at a $K \times K$ matrix, associates its spectrum.
\end{proof}

\begin{proof}[Proof of Remark \ref{i42}]
We come back to Example \ref{i41} with $A_{2} \geq A_{1} + 1$ and
\begin{equation} \label{i46}
\frac{\CM^{+}_{2}}{\CM^{-}_{2}} \sqrt{\frac{\vert g_{-}^{2} \vert}{\vert g_{+}^{2} \vert}} = \varepsilon .
\end{equation}
Thus, $P_{0}$ and $P$ satisfy the assumptions of Proposition \ref{i38} for $\varepsilon > 0$ small enough. We define $\overline{\CQ}_{0}$ and $\overline{\CQ}$ as in \eqref{i43} and study their eigenvalues. The $2 \times 2$ matrix $\overline{\CQ}$ takes the form \eqref{i44} with $a , b = \CO ( 1 )$ and $c , d = \CO ( \varepsilon )$. Using \eqref{e21}, the eigenvalue $\mu_{1} ( \tau , h )$ of $\overline{\CQ} ( \tau , h )$ with the largest modulus (which will provide the resonances closest to the real axis thanks to \eqref{d95}) is given by
\begin{align}
\mu_{1} ( \tau , h ) &= \frac{a + d + a \sqrt{1 - 2 d a^{- 1} + d^{2} a^{- 2} + 4 b c a^{- 2}}}{2}  \nonumber \\
&= \frac{a + d + a - d + 2 b c a^{- 1} + \CO ( \varepsilon^{2} )}{2} \nonumber \\
&= a + b c a^{- 1} + \CO ( \varepsilon^{2} ) , \label{i45}
\end{align}
since $a^{- 1} = \CO ( 1 )$. Of course, $\mu_{1}^{0} ( \tau , h ) = a$ is the (unique) eigenvalue of $\overline{\CQ}_{0}$. We deduce from \eqref{d4}, \eqref{d10}, \eqref{i46} and \eqref{i45} that
\begin{align*}
\mu_{1}^{0} ( \tau , h ) &= e^{i A_{1} / h} M ( \tau ) \\
\mu_{1} ( \tau , h ) &= e^{i A_{1} / h} M ( \tau ) + e^{i A_{2} / h} \varepsilon N ( \tau ) + \CO ( \varepsilon^{2} ) ,
\end{align*}
where $M$ and $N$ do not vanish on $\R$.

The previous equations yield
\begin{align}
\vert \mu_{1} ( \tau , h ) \vert &\leq  \big\vert M ( \tau ) + e^{i ( A_{2} - A_{1} ) / h} \varepsilon N ( \tau ) \big\vert + \CO ( \varepsilon^{2} )  \nonumber \\
&\leq \vert M ( \tau ) \vert + \varepsilon \vert N ( \tau ) \vert \cos \big( \theta ( \tau ) + ( A_{2} - A_{1} ) / h \big) + \CO ( \varepsilon^{2} ) , \label{i47}
\end{align}
where $\theta ( \tau )$ is a continuous argument of $N ( \tau ) \overline{M ( \tau )}$. We now take any interval $[ A , B ] \subset \R$ and $\theta_{0} \in \R$ such that $\vert \theta ( \tau ) - \theta_{0} \vert \leq 1/2$ for all $\tau \in [ A , B ]$. We now restrict $h$ to a sequence of positive numbers converging to $0$ such that
\begin{equation*}
\dist \big( \theta ( \tau ) + ( A_{2} - A_{1} ) / h  , \pi + 2 \pi \Z \big) \leq 1/2 ,
\end{equation*}
for all $\tau \in [ A , B]$ and $h$ in this sequence. Then, \eqref{i47} becomes
\begin{equation*}
\vert \mu_{1} ( \tau , h ) \vert \leq \vert \mu_{1}^{0} ( \tau , h ) \vert- \varepsilon \mu + \CO ( \varepsilon^{2} ) ,
\end{equation*}
for some $\mu > 0$ and all $\tau \in [ A , B ]$ and $h$ in this sequence. Applying \eqref{d95} and Theorem \ref{d8}, the conclusion of Remark \ref{i42} holds true in
\begin{equation*}
E_{0} + [ A h , B  h ] + i \Big[ - \sum_{j = 2}^{n} \frac{\lambda_{j}}{2} h - C \frac{h}{\vert \ln h \vert} , h \Big] \setminus \big( \Gamma (h) + B ( 0 , \eta h ) \big) ,
\end{equation*}
for any $\eta > 0$.

Finally, we explain how to avoid $\Gamma (h) + B ( 0 , \eta h )$. Since $\Gamma (h)$ is uniformly bounded in compact regions of size $h$, there exists $R \in \N$ such that
\begin{equation*}
\card \bigg( \Gamma ( h ) \cap E_{0} + [ A h , B  h ] + i \Big[ - \sum_{j = 2}^{n} \frac{\lambda_{j}}{2} h - C \frac{h}{\vert \ln h \vert} , h \Big] \bigg) \leq R ,
\end{equation*}
for all $h$. We now decompose $[ A , B ]$ in $2 R + 1$ intervals
\begin{equation*}
[ A , B ] = \bigcup_{j = 1}^{2 R + 1} I_{j} \qquad \text{with} \qquad I_{j} = \Big[ A + ( j - 1 ) \frac{B - A}{2 R + 1} , A + j \frac{B - A}{2 R + 1} \Big] ,
\end{equation*}
and fix $\eta < ( B - A ) / ( 4 R + 2 )$. For all $h$ in the sequence previously constructed, there exists at least a $j ( h ) \in \{ 1 , \ldots , 2 R + 1 \}$ such that
\begin{equation*}
E_{0} + h I_{j ( h )} + i \Big[ - \sum_{j = 2}^{n} \frac{\lambda_{j}}{2} h - C \frac{h}{\vert \ln h \vert} , h \Big] \bigcap \Gamma (h) + B ( 0 , \eta h ) = \emptyset .
\end{equation*}
Then, up to the extraction of a subsequence, we can always assume that this $j ( h )$ is constant and the conclusion of Remark \ref{i42} (without $\Gamma ( h )$) holds true in the interval $E_{0} + h I_{j}$.
\end{proof}

\Subsection{Resonance free domains below the accumulation curves} \label{s78}

The aim of this part is to demonstrate the results of Section \ref{s15} and \ref{s16}, more precisely Theorem \ref{e11}, Proposition \ref{e15} and Proposition \ref{e19}. For that, we use the general strategy explained in Section \ref{s36} and developed in Section \ref{s72}. But since we have already proved Theorem \ref{d8}, we only have to show that $P$ has no resonance and that its truncated resolvent has a polynomial estimate in some appropriate domains (as explained in \eqref{e84}). Thus, we follow Section \ref{s72}. We do not recall here all the steps and detail only the new arguments. The main difference with the proof of Proposition \ref{d39} is that the ``spectral parameter'' $h^{S ( z / h ) / \lambda_{1} - 1 / 2}$ is no longer bounded in the domain of study but rather verifies
\begin{equation} \label{e85}
\big\vert h^{S ( z / h ) / \lambda_{1} - 1 / 2} \big\vert \lesssim h^{- \alpha / \lambda_{1}} ,
\end{equation}
for $\im z \geq - h \sum_{j = 2}^{n} \lambda_{j} / 2 - h \alpha$.

\begin{proof}[Proof of Theorem \ref{e11}]
We show that the resolvent of $P$ has a polynomial estimate in \eqref{e84}. The beginning of the Section \ref{s72} can be applied without changes. In particular, we have $u_{-}^{k} \in \CI ( \Lambda_{+}^{1} , h^{-N} )$ as in \eqref{d44}. Its symbol is defined as in \eqref{d48}. But, since \eqref{e85} holds true for $z \in \eqref{e84}$, \eqref{d42} must be replaced by
\begin{equation} \label{e89}
a_{-}^{k} ( x , h ) = h^{S ( z , h ) / \lambda_{1} - 1 / 2} \sum_{\ell = 1}^{K} \CP_{k , \ell} ( x , h ) a_{-}^{\ell} ( x_{-}^{\ell} , h ) + S ( h^{-N + \zeta - \frac{\alpha}{\lambda_{1}}} ) ,
\end{equation}
for all $x$ near $\pi_{x} ( U_{-}^{k} \cap \Lambda_{+}^{1} )$. The additional factor $h^{- \frac{\alpha}{\lambda_{1}}}$ in the remainder term comes from the factor $h^{S ( z , h ) / \lambda_{1} - 1 / 2}$ in Corollary \ref{d46}. Then, with the notations of \eqref{e72} and assuming that $\alpha$ is small enough, we deduce
\begin{equation} \label{e86}
\big( 1 - h^{S ( z , h ) / \lambda_{1} - 1 / 2} \CQ ( z , h ) \big) a_{-} ( x_{-} , h ) = \CO ( h^{-N + \frac{\zeta}{2}} ) .
\end{equation}
We write
\begin{align}
\big( 1 - h^{S ( z , h ) / \lambda_{1} - 1 / 2} & \CQ ( z , h ) \big)^{- 1} \nonumber \\
&= - h^{- S ( z , h ) / \lambda_{1} + 1 / 2} \CQ ( z , h )^{- 1} \big( 1 - h^{- S ( z , h ) / \lambda_{1} + 1 / 2} \CQ ( z , h )^{- 1} \big)^{- 1} .    \label{e87}
\end{align}
For $z \in \eqref{e84}$, we have
\begin{equation} \label{e90}
\big\vert h^{- S ( z , h ) / \lambda_{1} + 1 / 2} \big\vert = h^{- \sum_{j=1}^{n} \frac{\lambda_{j}}{2 \lambda_{1}} - \frac{\im z}{\lambda_{1} h} + \frac{1}{2}} \leq e^{- C / \lambda_{1}} .
\end{equation}
Moreover, \ref{h10} and \eqref{e10} assure that $\Vert \CQ ( z , h )^{- 1} \Vert \lesssim 1$ uniformly for $z \in \eqref{e84}$. Then, if $C$ is chosen large enough, \eqref{e87} implies
\begin{equation} \label{e93}
\big\Vert \big( 1 - h^{S ( z , h ) / \lambda_{1} - 1 / 2} \CQ ( z , h ) \big)^{- 1} \big\Vert \lesssim \big\vert h^{- S ( z , h ) / \lambda_{1} + 1 / 2} \big\vert .
\end{equation}
and then
\begin{equation*}
a_{-} ( x_{-} , h ) = \CO \big( h^{-N + \frac{\zeta}{2}} \big\vert h^{- S ( z , h ) / \lambda_{1} + 1 / 2} \big\vert \big) ,
\end{equation*}
from \eqref{e86}. Eventually, inserting this estimate in \eqref{e89} gives
\begin{equation} \label{e91}
u_{-}^{k} \in \CI ( \Lambda_{+}^{1} , h^{- N + \frac{\zeta}{2}} ) .
\end{equation}
In other words, we have gained $h^{\frac{\zeta}{2}}$ in the order of $u_{-}$. Thus, the proof follows from the bootstrap argument at the end of Section \ref{s4}.
\end{proof}

\begin{proof}[Proof of Proposition \ref{e15}]
We have to prove here that the resolvent of $P$ satisfies a polynomial estimate in the domain $E_{0} + h \Omega \setminus ( \Gamma (h) + B ( 0 , \delta h ) )$. We follow the strategy of the proof of Theorem \ref{e11}. Nevertheless, we can not use \eqref{e89} here since this estimate is only valid for $\im z \geq - h \sum_{j = 2}^{n} \lambda_{j} / 2 - h \alpha$ with $\alpha > 0$ small enough and not deeper in the complex plane. Instead, we will used \eqref{m78} which shows that
\begin{equation} \label{e92}
a_{-}^{k} ( x , h ) = h^{S ( z , h ) / \lambda_{1} - 1 / 2} \sum_{\ell = 1}^{K} \widetilde{\CP}_{k , \ell} ( x , h ) a_{-}^{\ell} ( x_{-}^{\ell} , h ) + S \big( h^{- N + 1} \big\vert h^{S ( z , h ) / \lambda_{1} - 1 / 2} \big\vert \big) ,
\end{equation}
where $\widetilde{\CP}_{k , \ell} = \CP_{k , \ell} + S ( h^{\zeta} )$. Taking $x = x_{-}^{k}$ in the previous equation, we deduce
\begin{equation*}
\big( 1 - h^{S ( z , h ) / \lambda_{1} - 1 / 2} \big( \CQ ( z , h ) + \CO ( h^{\zeta} ) \big) \big) a_{-} ( x_{-} , h ) = \CO \big( h^{- N + 1} \big\vert h^{S ( z , h ) / \lambda_{1} - 1 / 2} \big\vert \big) .
\end{equation*}
As in \eqref{e93}, we have
\begin{equation*}
\big\Vert \big( 1 - h^{S ( z , h ) / \lambda_{1} - 1 / 2} \big( \CQ ( z , h ) + \CO ( h^{\zeta} ) \big) \big)^{- 1} \big\Vert \lesssim \big\vert h^{- S ( z , h ) / \lambda_{1} + 1 / 2} \big\vert ,
\end{equation*}
and then
\begin{equation*}
a_{-} ( x_{-} , h ) = \CO ( h^{- N + 1} ) .
\end{equation*}
Combining this estimate with \eqref{e92}, we obtain
\begin{equation} \label{e94}
a_{-}^{k} \in S \big( h^{- N + 1} \big\vert h^{S ( z , h ) / \lambda_{1} - 1 / 2} \big\vert \big) .
\end{equation}
By assumption, there exists $\nu > 0$ such that $\im \sigma \geq - \sum_{j = 2}^{n} \lambda_{j} / 2 - \lambda_{1} + \nu$ for all $\sigma \in \Omega$. Then, we get
\begin{equation*}
\big\vert h^{S ( z , h ) / \lambda_{1} - 1 / 2} \big\vert = h^{\sum_{j=1}^{n} \frac{\lambda_{j}}{2 \lambda_{1}} + \frac{\im z}{\lambda_{1} h} - \frac{1}{2}} \leq h^{- 1 + \frac{\nu}{\lambda_{1}}} ,
\end{equation*}
and \eqref{e94} gives $a_{-}^{k} \in S ( h^{- N + \frac{\nu}{\lambda_{1}}} )$. Thus, we have gained $h^{\frac{\nu}{\lambda_{1}}}$ in the order of $u_{-}$. The usual bootstrap argument implies the proposition.
\end{proof}

\begin{proof}[Proof of Proposition \ref{e19}]
From Theorem \ref{d8}, we have to show again that the resolvent of $P$ admits a polynomial estimate in the set \eqref{e84}. Thus, this proof is similar to the one of Theorem \ref{e11}. The unique difference is that we can not write \eqref{e87} since $\CQ$ is no longer invertible. In order to estimate the resolvent of $\CQ$, we write
\begin{equation} \label{e88}
\big( 1 - h^{S ( z , h ) / \lambda_{1} - 1 / 2} \CQ ( z , h ) \big)^{- 1} = h^{- S ( z , h ) / \lambda_{1} + 1 / 2} \big( h^{- S ( z , h ) / \lambda_{1} + 1 / 2} - \CQ ( z , h ) \big)^{- 1} .
\end{equation}
From the hypotheses of Proposition \ref{e19}, the eigenvalues of $\CQ ( z , h)$ are $0$ with multiplicity $K - 1$ and 
\begin{equation*}
\mu \Big( \frac{z - E_{0}}{h} + i \sum_{j = 2}^{n} \frac{\lambda_{j}}{2} , h \Big) ,
\end{equation*}
with multiplicity $1$. Moreover, the coefficients of $\CQ ( z ,h)$ are uniformly bounded in \eqref{e84}. Using that the inverse of a matrix is the transpose of its cofactor matrix divided by its determinant, \eqref{e88} implies
\begin{equation*}
\big\Vert \big( 1 - h^{S ( z , h ) / \lambda_{1} - 1 / 2} \CQ ( z , h ) \big)^{- 1} \big\Vert \lesssim \frac{\big\vert h^{- S ( z , h ) / \lambda_{1} + 1 / 2} \big\vert}{\big\vert h^{- S ( z , h ) / \lambda_{1} + 1 / 2} \big\vert^{K - 1} \big\vert h^{- S ( z , h ) / \lambda_{1} + 1 / 2} - \mu \big\vert} .
\end{equation*}
Since $\mu ( \tau , h )$ avoids a neighborhood of $0$ and $h^{S ( z , h ) / \lambda_{1} - 1 / 2}$ satisfies \eqref{e85} and \eqref{e90}, the last estimate becomes
\begin{equation*}
\big\Vert \big( 1 - h^{S ( z , h ) / \lambda_{1} - 1 / 2} \CQ ( z , h ) \big)^{- 1} \big\Vert \lesssim h^{- \alpha K / \lambda_{1}} \big\vert h^{- S ( z , h ) / \lambda_{1} + 1 / 2} \big\vert \leq h^{- \frac{\zeta}{4}} \big\vert h^{- S ( z , h ) / \lambda_{1} + 1 / 2} \big\vert ,
\end{equation*}
for $\alpha$ small enough and $C$ large enough. Then, similarly to \eqref{e91}, we get $u_{-}^{k} \in \CI ( \Lambda_{+}^{1} , h^{- N + \frac{\zeta}{4}} )$ and the proposition follows.
\end{proof}

\begin{lemma}\sl \label{j59}
Consider the situation of Example \ref{e48} {\rm (B)}. There exists $\alpha > 0$ such that, for all $C , \delta > 0$, $P$ has no resonance in
\begin{equation} \label{j60}
E_{0} + [ - C h , C  h ] + i \Big[ - \sum_{j = 2}^{n} \frac{\lambda_{j}}{2} h - \alpha h , h \Big] \setminus \big( \Gamma (h) + B ( 0 , \delta h ) \big) .
\end{equation}
for $h$ small enough in the sequence $\big( \frac{A_{2} - A_{1}}{( 2 j + 1 ) \pi} \big)_{j \in \N}$.
\end{lemma}

\begin{proof}
We show that the resolvent of $P$ has a polynomial estimate in \eqref{j60} following the proof of Theorem \ref{e11}. Formulas \eqref{e89}--\eqref{e86} are still valid in the present case, but the resolvent estimate of $\CQ$ stated in \eqref{e93} can not be used here. Thanks to \eqref{j61} and the assumptions of the lemma, the function $\mu ( \tau , h )$ vanishes for all $\tau \in \C$ (outside the singularities $- i \lambda_{1} ( \N + 1 / 2 )$) and $h$ in the sequence $\big( \frac{A_{2} - A_{1}}{( 2 j + 1 ) \pi} \big)_{j \in \N}$. It implies that the $K$ eigenvalues of $\CQ ( z , h )$ are $0$ for all $z \in \eqref{j60}$. In particular, the Cayley--Hamilton theorem implies that $\CQ ( z , h )$ is nilpotent of order at most $K$ (i.e. $\CQ ( z , h )^{K} = 0$). Then, its resolvent satisfies
\begin{equation} \label{j62}
\big( 1 - h^{S ( z , h ) / \lambda_{1} - 1 / 2} \CQ ( z , h ) \big)^{- 1} = \sum_{k = 0}^{K - 1} \big( h^{S ( z , h ) / \lambda_{1} - 1 / 2} \big)^{k} \CQ ( z , h )^{k} .
\end{equation}
On the other hand, we have
\begin{equation*}
\big\vert h^{S ( z , h ) / \lambda_{1} - 1 / 2} \big\vert \leq h^{- \frac{\alpha}{\lambda_{1}}} ,
\end{equation*}
for $z \in \eqref{j60}$. Thus, \eqref{j62} gives
\begin{equation}
\big\Vert \big( 1 - h^{S ( z , h ) / \lambda_{1} - 1 / 2} \CQ ( z , h ) \big)^{- 1} \big\Vert \lesssim h^{- \frac{\alpha ( K - 1 )}{\lambda_{1}}} \leq h^{- \frac{\zeta}{4}} ,
\end{equation}
for $\alpha > 0$ small enough. Combining this estimate with \eqref{e86}, we deduce
\begin{equation*}
a_{-} ( x_{-} , h ) = \CO \big( h^{-N + \frac{\zeta}{4}} \big) .
\end{equation*}
Eventually, as in \eqref{e91}, we get $u_{-}^{k} \in \CI ( \Lambda_{+}^{1} , h^{- N + \frac{\zeta}{4} - \frac{\alpha}{\lambda_{1}}} ) \subset \CI ( \Lambda_{+}^{1} , h^{- N + \frac{\zeta}{8}} )$, for $\alpha > 0 $ small enough. It implies the lemma as in the previous proofs.
\end{proof}

\Subsection{Proof when a change of multiplicity occurs} \label{s76}

The aim of this part is to prove Proposition \ref{e56}. As in Section \ref{s13}, we start with the asymptotic of the pseudo-resonances stated in Lemma \ref{e52}.

Following \eqref{d20}--\eqref{d36}, we first exclude some parts of the set \eqref{e53} which give no contribution to the distribution of resonances. Let
\begin{equation} \label{e63}
\rho = \frac{z - E_{0}}{h} + i \sum_{j = 2}^{n} \frac{\lambda_{j}}{2} - \tau_{0} ,
\end{equation}
be the relevant rescaled spectral parameter and define the set
\begin{equation} \label{e57}
E_{0} + h \tau_{0} + h \varepsilon (h) [ - 1 , 1 ] - i h \sum_{j = 2}^{n} \frac{\lambda_{j}}{2} + i h \Big[  - N \frac{\ln \vert \ln h \vert}{\vert \ln h \vert} , \frac{N}{\vert \ln h \vert} \Big] ,
\end{equation}
with $N > 0$ large enough. In particular,
\begin{equation} \label{e58}
h^{S ( z , h ) / \lambda_{1} - 1 / 2} \mu \Big( \frac{z - E_{0}}{h} + i \sum_{j = 2}^{n} \frac{\lambda_{j}}{2} \Big) = e^{i \tau_{0} \vert \ln h \vert / \lambda_{1}} e^{i \rho \vert \ln h \vert / \lambda_{1}} \mu ( \tau_{0} + \rho ) .
\end{equation}
For $z \in \eqref{e53} \setminus \eqref{e57}$, we have
\begin{equation*}
\im \rho > \frac{N}{\vert \ln h \vert} \qquad \text{or} \qquad \im \rho < - N \frac{\ln \vert \ln h \vert}{\vert \ln h \vert} .
\end{equation*}
If $\im \rho > N / \vert \ln h \vert$, we get, for some $c > 0$,
\begin{equation} \label{e59}
\big\vert e^{i \rho \vert \ln h \vert / \lambda_{1}} \mu ( \tau_{0} + \rho ) \big\vert  \leq c e^{- \im \rho \vert \ln h \vert / \lambda_{1}} \leq c e^{- N / \lambda_{1}} \leq 1 / 2 ,
\end{equation}
for $N$ large enough. If now $\im \rho < - N \ln \vert \ln h \vert / \vert \ln h \vert$, we obtain the same way using also \ref{h12}, for some new constant $c > 0$,
\begin{equation} \label{e60}
\big\vert e^{i \rho \vert \ln h \vert / \lambda_{1}} \mu ( \tau_{0} + \rho ) \big\vert \geq c \frac{\ln^{\ell} \vert \ln h \vert}{\vert \ln h \vert^{\ell}} e^{- \im \rho \vert \ln h \vert / \lambda_{1}} \geq \frac{\ln^{\ell} \vert \ln h \vert}{\vert \ln h \vert^{\ell}} \vert \ln h \vert^{N / \lambda_{1}} \geq 2 ,
\end{equation}
for $N$ large enough. Therefore, combining \eqref{e58} with the estimates \eqref{e59} and \eqref{e60}, we obtain
\begin{equation} \label{e61}
\bigg\vert \bigg( h^{S ( z , h ) / \lambda_{1} - 1 / 2} \mu \Big( \frac{z - E_{0}}{h} + i \sum_{j = 2}^{n} \frac{\lambda_{j}}{2} \Big) - 1 \bigg)^{-1} \bigg\vert \leq 2,
\end{equation}
for all $z \in \eqref{e53} \setminus \eqref{e57}$. In particular, this set has no pseudo-resonance. In the sequel, we will thus assume that $z$ belongs to the set \eqref{e57}.

Before proving Lemma \ref{e52}, we collect some properties on the Lambert function $W$, the multivalued inverse of the complex function $x \mapsto x e^{x}$. We follow the presentation of \cite{CoGoHaJeKn96_01}. In particular, $W_{q}$ with $q \in \Z$ denotes the $q$-th branch of the Lambert function. We first give estimates on $W_{q} ( A )$ for large $A$.

\begin{lemma}\sl \label{e55}
There exists $C > 1$ such that, for all $x , A \in \C$ with $\vert A \vert > C$ and
\begin{equation*}
x e^{x} = A ,
\end{equation*}
we have
\begin{equation*}
\vert x \vert > \frac{\ln \vert A \vert}{2} \qquad \text{and} \qquad \ln \vert A \vert - \ln \vert x \vert = \re x < \ln \vert A \vert .
\end{equation*}
\end{lemma}

\begin{proof}
We first prove that $\vert x \vert > \ln \vert A \vert / 2$. If it was not the case, then
\begin{equation}
\vert A \vert = \vert x e^{x} \vert \leq \frac{\ln \vert A \vert}{2} e^{\ln \vert A \vert / 2} = \sqrt{\vert A \vert} \ln \vert A \vert / 2 ,
\end{equation}
which is impossible for $\vert A \vert > C \gg 1$. In particular, taking $C$ large enough, we can always assume that $\vert x \vert > \ln C / 2 > 1$.

Now, taking again the modulus of $x e^{x} = A$, we deduce
\begin{equation*}
\ln \vert x \vert + \re x = \ln \vert A \vert .
\end{equation*}
Using $\ln \vert x \vert > 0$, it yields $\ln \vert A \vert - \ln \vert x \vert = \re x < \ln \vert A \vert$.
\end{proof}

We also give an upper bound on the derivative of the Lambert function.

\begin{lemma}\sl \label{e62}
There exists $C > 1$ such that, for all $A \in \C$ with $\vert A \vert > C$, we have
\begin{equation*}
\vert W^{\prime} ( A ) \vert \lesssim \frac{1}{\vert A \vert} .
\end{equation*}
\end{lemma}

\begin{proof}
Differentiating the formula $W (A) e^{W ( A )} = A$, we obtain
\begin{equation*}
A W^{\prime} ( A ) \Big( 1 + \frac{1}{W ( A )} \Big) = 1 .
\end{equation*}
Since Lemma \ref{e55} implies that $\vert W ( A ) \vert > 2$ for $A$ large enough, the lemma follows.
\end{proof}

\begin{proof}[Proof of Lemma \ref{e52}]
As in the proof of Proposition \ref{d9}, we first show that every pseudo-resonance satisfies \eqref{e54}. So, we consider a pseudo-resonance $z \in \eqref{e53}$. From \eqref{e61}, we have in fact $z \in \eqref{e57}$. Then, the rescaled spectral variable $\rho$, defined in \eqref{e63}, belongs to
\begin{equation} \label{e64}
\varepsilon (h) [ - 1 , 1 ] + i \Big[  - N \frac{\ln \vert \ln h \vert}{\vert \ln h \vert} , \frac{N}{\vert \ln h \vert} \Big] .
\end{equation}
In particular, $\rho = o (1)$. Using the assumption \ref{h12} and \eqref{e58}, we deduce
\begin{equation*}
\alpha \rho^{\ell} \big( 1 + \CO ( \rho ) \big) e^{i \rho \vert \ln h \vert / \lambda_{1}} = e^{- i \tau_{0} \vert \ln h \vert / \lambda_{1}} .
\end{equation*}
Then, we can write
\begin{equation*}
i \frac{\rho \vert \ln h \vert}{\ell \lambda_{1}} e^{i \frac{\rho \vert \ln h \vert}{\ell \lambda_{1}}} = i \frac{\vert \ln h \vert}{\ell \lambda_{1} \alpha^{1 / \ell}} e^{- i \frac{\tau_{0}}{\ell \lambda_{1}} \vert \ln h \vert} e^{i 2 \pi \beta / \ell} ( 1 + o (1) ) ,
\end{equation*}
for some $\beta \in \{ 0 , \ldots , \ell - 1\}$. By definition, it yields
\begin{equation*}
\rho = - i \ell \lambda_{1} W_{q^{\prime}} \Big( i \frac{\vert \ln h \vert}{\ell \lambda_{1} \alpha^{1 / \ell}} e^{- i \frac{\tau_{0}}{\ell \lambda_{1}} \vert \ln h \vert} e^{i 2 \pi \beta / \ell} ( 1 + o (1) ) \Big) \frac{1}{\vert \ln h \vert} ,
\end{equation*}
for some $q^{\prime} \in \Z$. Using now Lemma \ref{e62} and the mean value theorem, we obtain
\begin{equation} \label{e69}
\rho = - i \ell \lambda_{1} W_{q} \Big( i \frac{\vert \ln h \vert}{\ell \lambda_{1} \alpha^{1 / \ell}} e^{- i \frac{\tau_{0}}{\ell \lambda_{1}} \vert \ln h \vert} e^{i 2 \pi \beta / \ell} \Big) \frac{1}{\vert \ln h \vert} + o \Big( \frac{1}{\vert \ln h \vert} \Big) ,
\end{equation}
for some $q \in q^{\prime} + \{ - 1 , 0 , 1 \}$ (the branch of the Lambert function may change). Coming back to the original spectral variable $z$, we obtain exactly \eqref{e54}.

We now prove that every complex number given by \eqref{e65} is close to a pseudo-resonance. Let $q \in \Z$ and $\beta \in \{ 0 , \ldots , \ell - 1 \}$ be such that $z_{q , \beta} \in \eqref{e53}$. As in \eqref{e61}, we first show that $z_{q , \beta} \in \eqref{e57}$ for $N$ large enough. If it was not the case, we have
\begin{equation} \label{n8}
\re w < - N \qquad \text{or} \qquad \re w > N \ln \vert \ln h \vert ,
\end{equation}
with
\begin{equation*}
w = W_{q} \Big( i \frac{\vert \ln h \vert}{\ell \lambda_{1} \alpha^{1 / \ell}} e^{- i \frac{\tau_{0}}{\ell \lambda_{1}} \vert \ln h \vert} e^{i 2 \pi \beta / \ell} \Big) .
\end{equation*}
If $\re w < - N$, then
\begin{equation*}
\frac{\vert \ln h \vert}{\ell \lambda_{1} \alpha^{1 / \ell}} = \big\vert w e^{w} \big\vert = \vert w \vert e^{\re w} \lesssim e^{- N} \vert \ln h \vert ,
\end{equation*}
since $\vert w \vert \lesssim \vert \ln h \vert$ for $z_{q , \beta} \in \eqref{e53}$. This inequality is impossible for $N$ large enough. If $\re w > N \ln \vert \ln h \vert$, we have
\begin{equation*}
\frac{\vert \ln h \vert}{\ell \lambda_{1} \alpha^{1 / \ell}} = \big\vert w e^{w} \big\vert = \vert w \vert e^{\re w} \gtrsim e^{N \ln \vert \ln h \vert} \ln \vert \ln h \vert = \vert \ln h \vert^{N} \ln \vert \ln h \vert ,
\end{equation*}
from Lemma \ref{e55}. We get again a contradiction for $N$ large enough. Thus, \eqref{n8} does not holds and $z_{q , \beta} \in \eqref{e57}$. In order to find a pseudo-resonance $z$ satisfying \eqref{e54}, it is equivalent to find $r = o (1)$ such that
\begin{equation} \label{e66}
\mu \Big( \tau_{0} + \rho_{q , \beta} + \frac{r}{\vert \ln h \vert} \Big) e^{i \tau_{0}\vert \ln h \vert / \lambda_{1}} e^{i \rho_{q , \beta} \vert \ln h \vert / \lambda_{1}} e^{i r / \lambda_{1}} = 1 ,
\end{equation}
where
\begin{equation} \label{e67}
\rho_{q , \beta} = - i \ell \lambda_{1} W_{q} \Big( i \frac{\vert \ln h \vert}{\ell \lambda_{1} \alpha^{1 / \ell}} e^{- i \frac{\tau_{0}}{\ell \lambda_{1}} \vert \ln h \vert} e^{i 2 \pi \beta / \ell} \Big) \frac{1}{\vert \ln h \vert} .
\end{equation}
Using \eqref{e67} and $e^{W_{q} (x)} = x / W_{q} (x)$, we have
\begin{equation*}
\Big( \rho_{q , \beta} + \frac{r}{\vert \ln h \vert} \Big)^{\ell} e^{i \rho_{q , \beta} \vert \ln h \vert / \lambda_{1}} = \alpha^{- 1} e^{- i \tau_{0} \vert \ln h \vert / \lambda_{1}} \Big( 1 + \frac{r}{\rho_{q , \beta} \vert \ln h \vert} \Big)^{\ell} .
\end{equation*}
Combining with \ref{h12}, \eqref{e66} becomes
\begin{equation} \label{n7}
\Big( 1 + \frac{r}{\rho_{q , \beta} \vert \ln h \vert} \Big)^{\ell} \Big( 1 + g \Big( \rho_{q , \beta} + \frac{r}{\vert \ln h \vert} \Big) \Big) e^{i r / \lambda_{1}} - 1 = 0 ,
\end{equation}
for some function $g$ holomorphic near $0$ with $g ( 0 ) = 0$. The left hand side of \eqref{n7} can be written $F ( r , y )$ where $y \in \C^{3}$ is a shortcut for
\begin{equation} \label{e68}
y = \Big( \rho_{q , \beta} , \frac{1}{\vert \ln h \vert} , \frac{1}{\rho_{q , \beta} \vert \ln h \vert} \Big) ,
\end{equation}
$F$ is holomorphic near $0$, $F ( 0 , 0 ) = 0$ and $\partial_{r} F ( 0 , 0 ) = i / \lambda_{1}$. Note that $y = o (1)$ thanks to Lemma \ref{e55} and $\rho_{q , \beta} \in \eqref{e64}$. Then, the implicit function theorem provides a solution $r = r (y) = o (1)$ of $F ( r , y ) = 0$, and Lemma \ref{e52} follows.
\end{proof}

As in Lemma \ref{d13}, we now show that the ``inverse of the quantization rule'' is uniformly bounded away from the pseudo-resonances. Since the proof of this result is similar to the first part of the proof of Lemma \ref{e52} (see \eqref{e64}--\eqref{e69}), we omit it.

\begin{lemma}\sl \label{e71}
Let $\delta > 0$. Then, there exists $M > 0$ such that, for all $z$ in \eqref{e53}, we have
\begin{equation*}
\dist \big( z , \res_{0} (P) \big) > \delta \frac{h}{\vert \ln h \vert} \quad \Longrightarrow \quad \bigg\vert \bigg( 1 - h^{S ( z , h ) / \lambda_{1} - 1 / 2} \mu \Big( \frac{z - E_{0}}{h} + i \sum_{j = 2}^{n} \frac{\lambda_{j}}{2} \Big) \bigg)^{-1} \bigg\vert \leq M .
\end{equation*}
\end{lemma}

\begin{proof}[Proof of Proposition \ref{e56}]
We use here the same strategy as for the proof of Theorem \ref{d8} and we only explain the necessary changes. We first show that $P$ has no resonance and that its resolvent satisfies a polynomial estimate in the set
\begin{align}
E_{0} + h \tau_{0} + h \varepsilon (h) [ - 1 , 1 ] &+ i h \Big[ - \sum_{j = 2}^{n} \frac{\lambda_{j}}{2} - \alpha , 1 \Big] \nonumber \\
&\setminus \big( \Gamma (h) + B ( 0 , \delta h ) \big) \bigcup \Big( \res_{0} (P) + B \Big( 0 , \delta \frac{h}{\vert \ln h \vert} \Big) \Big) .   \label{e73}
\end{align}
For that, we follow Section \ref{s72} and Section \ref{s78}. As in \eqref{e86} (see also \eqref{e72}), 
\begin{equation} \label{e74}
\big( 1 - h^{S ( z , h ) / \lambda_{1} - 1 / 2} \CQ ( z , h ) \big) a_{-} ( x_{-} , h ) = \CO ( h^{-N + \frac{\zeta}{2}} ) ,
\end{equation}
if $\alpha > 0$ is small enough. For $z \in \eqref{e73}$, we have $\vert h^{S ( z , h ) / \lambda_{1} - 1 / 2} \vert \lesssim h^{- \alpha / \lambda_{1}}$ (see \eqref{e85}) and the eigenvalues of $Q$ are $0$ with multiplicity $K - 1$ and $\mu$ with multiplicity $1$. Using that the inverse a matrix is given by the transpose of its cofactor matrix divided by its determinant, we deduce
\begin{align}
\big\Vert \big( 1 - h^{S ( z , h ) / \lambda_{1} - 1 / 2} \CQ ( z , h ) \big)^{- 1} \big\Vert &\lesssim \frac{h^{- \alpha ( K - 1 ) / \lambda_{1}}}{\vert 1 - h^{S ( z , h ) / \lambda_{1} - 1 / 2} \mu \vert}  \nonumber \\
&\lesssim h^{- \frac{\zeta}{4}} \bigg\vert \bigg( 1 - h^{S ( z , h ) / \lambda_{1} - 1 / 2} \mu \Big( \frac{z - E_{0}}{h} + i \sum_{j = 2}^{n} \frac{\lambda_{j}}{2} \Big) \bigg)^{-1} \bigg\vert , \label{e77}
\end{align}
for $\alpha > 0$ is small enough. Combining this estimate with Lemma \ref{e71}, \eqref{e74} yields
\begin{equation*}
a_{-} ( x_{-} , h ) = \CO ( h^{-N + \frac{\zeta}{4}} ) ,
\end{equation*}
and thus $a_{-}^{k} ( x , h ) \in S ( h^{-N + \frac{\zeta}{8}} )$ on $\pi_{x} ( U_{-}^{k} \cap \Lambda_{+}^{1} )$. Then, we conclude that $u = \CO ( h^{\infty} )$ and that $P$ has no resonance, together with a polynomial estimate of its resolvent in \eqref{e73} by the usual bootstrap argument developed at the end of Section \ref{s4}.

We now prove the existence of resonances near the pseudo-resonances following the same approach as in Section \ref{s73}. We can no longer define appropriate test functions with $\tau$ replaced by $\tau_{0}$ (see \eqref{d57}) since $\mu ( \tau_{0} ) = 0$ from \ref{h12}. Instead, we will use an explicit eigenvector of $\CQ ( z , h )$. Let $w_{0} ( z , h ) = ( w_{0}^{1} , \ldots , w_{0}^{K} ) \in \C^{K}$ be the $K$-vector with entries
\begin{equation*}
w_{0}^{k} = e^{i A_{k} / h} \Gamma \big( S ( z , h ) / \lambda_{1} \big) \sqrt{\frac{\lambda_{1}}{2 \pi}} \frac{\CM_{k}^{+}}{\CM_{k}^{-}} e^{- \frac{\pi}{2} ( \nu_{k} + \frac{1}{2} ) i} \big\vert g_{-}^{\ell} \big\vert \big( i \lambda_{1} \vert g_{+}^{k} \vert \big)^{- S ( z , h ) / \lambda_{1}} .
\end{equation*}
Since $\lambda_{1} < \lambda_{2}$ and all the asymptotic directions are at the same side of $0$, the coefficients of the matrix $\CQ ( z , h )$ have the product structure $\CQ_{k ,\ell} ( z , h ) = w_{0}^{k} \vert g_{-}^{\ell} \vert^{- S ( z , h ) / \lambda_{1}}$. Thus, with the notation of \eqref{e13}, we get
\begin{equation} \label{e75}
\CQ ( z , h ) w_{0} = \mu \Big( \frac{z - E_{0}}{h} + i \sum_{j = 2}^{n} \frac{\lambda_{j}}{2} \Big) w_{0} .
\end{equation}
Moreover, the vector $w_{0}$ does not cancel (more precisely, $1 \lesssim \Vert w_{0} \Vert \lesssim 1$ where $\Vert \cdot \Vert$ is any norm on $\C^{K}$) and $z \mapsto w_{0} ( z, h )$ is holomorphic by construction. We now build a test function $v$ as in \eqref{d60} and consider the solution $u$ of
\begin{equation}
( P_{\theta} - \widetilde{z} ) u = v ,
\end{equation}
where $\widetilde{z} \in \partial \CD$ and $\CD = B ( z_{q , \beta} , s h \vert \ln h \vert^{-1} )$ with $s > 0$ small enough and $z_{q , \beta} \in \eqref{e70}$. Note here that $v$ is holomorphic in $\CD$. We now follows Section \ref{s73} and obtain, as in Lemma \ref{d72} and Lemma \ref{d77}, that $u_{-}^{k} \in \CI ( \Lambda_{+}^{1} ,  1 )$,
\begin{equation*}
u_{-}^{k} (x) = e^{- i A_{k} / h} e^{i \frac{\widetilde{z} - E_{0}}{h} t_{-}^{k}} \frac{\CM_{k}^{-}}{\CD_{k} ( t_{-}^{k} )} a_{-}^{k} ( x , h ) e^{i \varphi_{+}^{1} (x) / h} ,
\end{equation*}
and
\begin{equation} \label{e76}
a_{-}^{k} ( x , h ) = h^{S ( \widetilde{z} / h ) / \lambda_{1} - 1 / 2} \sum_{\ell = 1}^{K} \CP_{k , \ell} ( x , h ) a_{-}^{\ell} ( x_{-}^{\ell} , h ) + \widetilde{a}_{k} ( x , h ) + S ( h^{\frac{\zeta}{2}} ) ,
\end{equation}
where $\widetilde{a}_{k} \in S ( 1 )$ is the symbol of the WKB solution $\widetilde{v}_{k}$ constructed as in \eqref{d89}. As usual, we apply \eqref{e76} with $x = x_{-}^{k}$ and obtain
\begin{equation*}
\big( 1 - h^{S ( \widetilde{z} / h ) / \lambda_{1} - 1 / 2} \CQ ( \widetilde{z} , h ) \big) a_{-} ( x_{-} , h ) = w_{0} + \CO ( h^{\frac{\zeta}{2}} ) ,
\end{equation*}
Since $\widetilde{z} \in \partial \CD$ is at distance $h \vert \ln h \vert^{- 1}$ from the pseudo-resonances, \eqref{e77} and Lemma \ref{e71} imply
\begin{align}
a_{-} ( x_{-} , h ) &= \big( 1 - h^{S ( \widetilde{z} / h ) / \lambda_{1} - 1 / 2} \CQ ( \widetilde{z} , h ) \big)^{- 1} w_{0} + \CO ( h^{\frac{\zeta}{4}} ) \nonumber \\
&= \bigg( 1 - h^{S ( \widetilde{z} , h ) / \lambda_{1} - 1 / 2} \mu \Big( \frac{\widetilde{z} - E_{0}}{h} + i \sum_{j = 2}^{n} \frac{\lambda_{j}}{2} \Big) \bigg)^{-1} w_{0} + \CO ( h^{\frac{\zeta}{4}} ) , \label{e78}
\end{align}
thanks to \eqref{e75}. On the other hand, using \eqref{e66}--\eqref{e68}, one can see that the holomorphic function in $\CD$
\begin{equation*}
1 - h^{S ( \widetilde{z} , h ) / \lambda_{1} - 1 / 2} \mu \Big( \frac{\widetilde{z} - E_{0}}{h} + i \sum_{j = 2}^{n} \frac{\lambda_{j}}{2} \Big) = - F \Big( ( \widetilde{z} - z_{q , \beta} ) \frac{\vert \ln h \vert}{h} , y \Big) ,
\end{equation*}
makes a simple loop around $0$ when $\widetilde{z}$ goes around $\partial \CD$. Moreover, recall that $\partial_{r} F ( 0 , 0 ) = i / \lambda_{1}$. Then, \eqref{e78} and the Cauchy formula give, as in Lemma \ref{d79},
\begin{equation} \label{e79}
\frac{1}{2 i \pi} \int_{\partial \CD} u ( \widetilde{z} ) \, d \widetilde{z} = b ( x , h ) e^{i \varphi_{+}^{1} (x) / h} \text{ microlocally near } U_{-}^{k} ,
\end{equation}
where $b \in S ( h \vert \ln h \vert^{-1} )$ satisfies
\begin{equation*}
b ( x_{-}^{k} , h ) = i e^{- i A_{k} / h} e^{- i \sigma_{0} t_{-}^{k}} \frac{\lambda_{1} \CM_{k}^{-}}{\CD_{k} ( t_{-}^{k} )} w_{0}^{k} \frac{h}{\vert \ln h \vert} + o \Big( \frac{h}{\vert \ln h \vert} \Big) .
\end{equation*}
As in the end of Section \ref{s73}, this implies that $P$ has a resonance in $\CD$. Indeed, otherwise, since $u = ( P_{\theta} - \widetilde{z} )^{- 1} v$ and $v$ is holomorphic in $\CD$, we deduce
\begin{equation*}
\frac{1}{2 i \pi} \int_{\partial \CD} u ( \widetilde{z} ) \, d \widetilde{z} = 0 ,
\end{equation*}
which is in contradiction with \eqref{e79}.
\end{proof}

\Subsection{Proof of the other subsets of resonances}\label{s84}

The goal of this section is to prove the results stated in Section \ref{s83}. We first explain how the matrix $\CQ^{2}$ is constructed. For that, we make Lemma \ref{d41} more precise. Recall that this lemma gives a closed relation on the symbol of the restrictions $u_{-}^{k} \in \CI ( \Lambda_{+}^{1} , h^{-N} )$ of a solution $u$ of \eqref{a4}. We send back the reader to Section \ref{s4} and Section \ref{s72} for notations and some preliminary results.

\begin{lemma}\sl \label{g89}
There exist $\zeta > 0$ and matrices $\CP^{( \bullet )}_{k} \in S ( 1 )$ independent of $u$ such that
\begin{align}
a_{-}^{k} ( x , h ) &= h^{S ( z , h ) / \lambda_{1} - 1 / 2} \big( ( \CP^{(0)}_{k} ( x , h ) + S ( h^{\zeta} ) ) \SA_{0} \nonumber \\
&\qquad \qquad + h \CP^{(1)}_{k} ( x , h ) \SA_{1} +  h \CP^{(2)}_{k} ( x , h ) \SA_{2} + S ( h^{-N + 1 + \zeta} ) \big) ,   \label{g90}
\end{align}
for all $x$ near $\pi_{x} ( U_{-}^{k} \cap \Lambda_{+}^{1} )$. Here, $\SA_{0}$ (resp. $\SA_{1}$ and $\SA_{2}$) is the $K$-vector (resp. $( n - 1 ) K$-vector and $( n - 1 )^{2} K$-vector) of the $a_{-}^{\ell} ( x_{-}^{\ell} , h )$ (resp. $\partial_{y_{a}} a_{-}^{\ell} ( x_{-}^{\ell} , h )$ and $\partial_{y_{a} , y_{b}}^{2} a_{-}^{\ell} ( x_{-}^{\ell} , h )$) where $\partial_{y_{\bullet}}$ denotes the derivatives in the $n - 1$ directions of $H_{\ell}$, the hyperplane of Theorem \ref{a32} associated to $\rho_{-}^{\ell}$. Moreover, the following properties are satisfied.

$i)$ For all $\alpha \in \N^{n}$, the matrix $\partial_{x}^{\alpha} \CP^{( \bullet )}_{k} ( x_{-}^{k} , h )$ can be written $\rho_{k} \widetilde{\CP}^{\alpha , ( \bullet )}_{k} ( \sigma )$ with $\widetilde{\CP}^{\alpha , ( \bullet )}_{k}$ independent of $h$ and holomorhic in $\sigma$ near $\im \sigma = - \sum_{j = 2}^{n} \lambda_{j} / 2 - \lambda_{1}$.

$ii)$ $\CP^{(0)}_{k} ( x_{-}^{k} , h ) = \CQ_{k} ( z , h )$ where $\CQ_{k}$ is the $k$-th line of the matrix $\CQ$ defined in \eqref{d4}.

$iii)$ There exists a $K \times ( n - 1 )^{2} K$ matrix $\CL = \CL ( z , h )$ independent of $x$ and $k$ such that $\CP^{(2)}_{k} ( x , h ) = \CP^{(0)}_{k} ( x , h ) \CL$.
\end{lemma}

\begin{proof}
Let us follow the proof of Lemma \ref{d41}. To simplify the exposition, we only study the contribution of $u_{\ell}$ in the computation of $u_{k}$. We start with the propagation through the fixed point $( 0 , 0 )$. To solve the microlocal Cauchy problem \eqref{d62}, we apply Theorem \ref{a32} instead of Corollary \ref{d46}. Thus, \eqref{a31} gives, modulo $S ( h^{\infty} )$,
\begin{align}
a_{+}^{k} ( x , h ) &= h^{S ( z , h ) / \lambda_{1} - n / 2} e^{i \frac{z - E_{0}}{h} t_{-}^{\ell}}  \nonumber \\
&\qquad \qquad \quad \times \frac{\CM_{\ell}^{-}}{\CD_{\ell} ( t_{-}^{\ell} )} \int_{H_{\ell}} e^{i ( \varphi_{+}^{1} - \varphi_{-} - A_{\ell} ) ( y ) / h} d ( x, y , z , h ) a_{-}^{\ell} ( y , h ) \, d y . \label{g98}
\end{align}
From \ref{h8}, the phase function $\varphi_{+}^{1} - \varphi_{-}$ has a unique non-degenerate critical point at $x_{-}^{\ell}$. Thus, performing a stationary phase expansion in \eqref{g98} (see Proposition 5.2 of \cite{DiSj99_01} for instance), \eqref{d51} is replaced by
\begin{equation} \label{g97}
a_{+}^{k} ( x , h ) =  h^{S ( z , h ) / \lambda_{1} - 1 / 2} ( \SB_{0} + h \SB_{1} + S ( h^{- N + 2} ) ) ,
\end{equation}
where the symbols $\SB_{j} \in S ( h^{- N} )$, $j = 0 , 1$, have the form
\begin{equation*}
\SB_{j} ( x , h ) = \big( B_{j}^{\ell} ( \partial_{y} ) d ( x, y , z , h ) a_{-}^{\ell} ( y , h ) \big) \vert_{y = x_{-}^{\ell}} ,
\end{equation*}
and $B_{j}^{\ell}$ are differential operator of order $\leq 2 j$ independent of $x$ and $k$. In particular, $B_{0}^{\ell}$ is the usual non-zero constant. On one hand, we have (by comparison with \eqref{d51} or by an application of Corollary \ref{d46})
\begin{equation} \label{g99}
\SB_{0} ( x , h ) = \big( \widetilde{\CR}_{k , \ell} ( x , h ) + S ( h^{\zeta} ) \big) a_{-}^{\ell} ( x_{-}^{\ell} , h ) ,
\end{equation}
and
\begin{equation} \label{i9}
\widetilde{\CR}_{k , \ell} ( x , h ) = B_{0}^{\ell} d_{0} ( x , x_{-}^{\ell} , z ) ,
\end{equation}
with $\widetilde{\CR}_{k , \ell}$ defined in \eqref{d51} and some $\zeta > 0$. On the other hand, a direct computation shows
\begin{align}
\SB_{1} ( x , h ) &= \big( B_{j}^{\ell} ( \partial_{y} ) d_{0} ( x, y , z ) a_{-}^{\ell} ( y , h ) \big) \vert_{y = x_{-}^{\ell}} + S ( h^{- N + \zeta} )   \nonumber \\
&= \sum_{a , b = 1}^{n - 1} \CC_{a , b} ( x , h ) \partial_{y_{a} , y_{b}}^{2} a_{-}^{\ell} ( x_{-}^{\ell} , h ) + \sum_{a = 1}^{n - 1} \CC_{a} ( x , h ) \partial_{y_{a}} a_{-}^{\ell} ( x_{-}^{\ell} , h )   \nonumber \\
&\qquad \qquad \qquad \qquad \qquad \qquad \qquad \qquad + \CC ( x , h ) a_{-}^{\ell} ( x_{-}^{\ell} , h ) + S ( h^{- N + \zeta} ), \label{i3}
\end{align}
for some $\CC_{\bullet} \in S (1)$. Moreover, \eqref{b19} implies that these symbols, as well as their derivatives, computed at $x = x_{+}^{k}$ can be written as functions of the rescaled parameter $\sigma$ that are independent of $h$, namely
\begin{equation} \label{i4}
\partial_{x}^{\alpha} \CC_{\bullet} ( x_{+}^{k} , h ) = \widetilde{\CC}^{\alpha}_{\bullet} ( \sigma ) .
\end{equation}
Eventually, since $B_{1}^{\ell}$ is a differential operator of order at most $2$, we have
\begin{equation} \label{i5}
\CC_{a , b} ( x , h ) = \widehat{\CC}_{a , b}^{\ell} d_{0} ( x , x_{-}^{\ell} , z ) ,
\end{equation}
for some constants $\widehat{\CC}_{a , b}^{\ell} = \widehat{\CC}_{a , b}^{\ell} ( z , h )$ independent of $k$. Combining \eqref{g97} with \eqref{g99} and \eqref{i3}, we deduce
\begin{align}
a_{+}^{k} ( x , h ) ={}& h^{S ( z , h ) / \lambda_{1} - 1 / 2} \Big( \big( \widetilde{\CR}_{k , \ell} ( x , h ) + S ( h^{\zeta} ) \big) a_{-}^{\ell} ( x_{-}^{\ell} , h )   \nonumber \\
&+ h \sum_{a = 1}^{n - 1} \CC_{a} ( x , h ) \partial_{y_{a}} a_{-}^{\ell} ( x_{-}^{\ell} , h ) + h \sum_{a , b = 1}^{n - 1} \CC_{a , b} ( x , h ) \partial_{y_{a} , y_{b}}^{2} a_{-}^{\ell} ( x_{-}^{\ell} , h ) + S ( h^{- N + 1 + \zeta} ) \Big) .    \label{i6}
\end{align}

The evolution equation satisfied by $u$ along the Hamiltonian curve $\gamma_{k}$ is given by the microlocal Cauchy problem \eqref{d65}. Taking into account that the initial condition verifies \eqref{i6} and the form of the usual transport equations, the propagation of Lagrangian distributions implies
\begin{align}
e^{- i A_{k} / h} e^{i \frac{z - E_{0}}{h} t_{-}^{k}} \frac{\CM_{k}^{-}}{\CD_{k} ( t_{-}^{k} )} a_{-}^{k} ( x , & h ) = h^{S ( z , h ) / \lambda_{1} - 1 / 2} \Big( \big( \CS_{k} ( x , h ) \widetilde{\CR}_{k , \ell} ( \widetilde{x} (x) , h ) + S ( h^{\zeta} ) \big) a_{-}^{\ell} ( x_{-}^{\ell} , h )   \nonumber \\
&+ h \sum_{a = 1}^{n - 1} \CS_{k} ( x , h ) \CC_{a} ( \widetilde{x} (x) , h ) \partial_{y_{a}} a_{-}^{\ell} ( x_{-}^{\ell} , h )   \label{i8}  \\
&+ h \sum_{a , b = 1}^{n - 1} \CS_{k} ( x , h ) \CC_{a , b} ( \widetilde{x} (x) , h ) \partial_{y_{a} , y_{b}}^{2} a_{-}^{\ell} ( x_{-}^{\ell} , h ) + S ( h^{- N + 1 + \zeta} ) \Big) , \nonumber
\end{align}
for all $x$ near $\pi_{x} ( U_{-}^{k} \cap \Lambda_{+}^{1} )$. In the previous expression, $\widetilde{x} (x)$ and $\CS_{k}$ are as in \eqref{d52}. The previous equation yields \eqref{g90} with
\begin{equation} \label{i10}
\left\{ \begin{aligned}
&\CP^{(0)}_{k} ( x , h ) : = \CS_{k} ( x , h ) \widetilde{\CR}_{k , \ell} ( \widetilde{x} (x) , h ) e^{i A_{k} / h} e^{- i \frac{z - E_{0}}{h} t_{-}^{k}} \CD_{k} ( t_{-}^{k} ) / \CM_{k}^{-} ,   \\
&\big( \CP^{(1)}_{k} ( x , h ) \big)_{a} : = \CS_{k} ( x , h ) \CC_{a} ( \widetilde{x} (x) , h ) e^{i A_{k} / h} e^{- i \frac{z - E_{0}}{h} t_{-}^{k}} \CD_{k} ( t_{-}^{k} ) / \CM_{k}^{-} ,  \\
&\big( \CP^{(2)}_{k} ( x , h ) \big)_{a , b} : = \CS_{k} ( x , h ) \CC_{a , b} ( \widetilde{x} (x) , h ) e^{i A_{k} / h} e^{- i \frac{z - E_{0}}{h} t_{-}^{k}} \CD_{k} ( t_{-}^{k} ) / \CM_{k}^{-} .
\end{aligned} \right.
\end{equation}

We now prove the properties on the matrices $\CP^{(\bullet)}_{k}$ stated in Lemma \ref{g89}. The point $i)$ follows from \eqref{i10}, the form of $\widetilde{\CR}_{k , \ell}$ and $\CS_{k}$ defined the proof of Lemma \ref{d41} and the definition \eqref{i3} of $\CC_{a}$ and $\CC_{a , b}$. On the other hand, \eqref{d76} implies $\CP^{(0)}_{k} ( x_{-}^{k} , h ) = \CQ_{k} ( z , h )$, that is $ii)$. Eventually, using \eqref{i9}, \eqref{i5} and \eqref{i10}, we deduce
\begin{equation*}
\big( \CP_{k}^{( 2 )} (x) \big)_{\ell , a , b} = \big( \CP_{k}^{( 0 )} (x) \big)_{\ell} \widehat{\CC}_{a , b}^{\ell} / B_{0}^{\ell} ,
\end{equation*}
and then $\CP^{(2)}_{k} (x) = \CP^{(0)}_{k} (x) \CL$ where the $K \times ( n - 1 )^{2} K$ matrix $\CL$, defined by $\CL_{( \ell^{\prime} ) , ( \ell , a , b)} = \widehat{\CC}_{a , b}^{\ell} / B_{0}^{\ell} \delta_{\ell^{\prime} , \ell}$, is independent of $x$ and $k$. This shows $iii)$.
\end{proof}

We assume that $z \in \eqref{g87}$. In particular, we have
\begin{equation} \label{g93}
1 \lesssim \vert h^{S ( z , h ) / \lambda_{1} + 1 / 2} \vert \lesssim 1 ,
\end{equation}
uniformly for $z \in \eqref{g87}$. Computing \eqref{g90} at $x = x_{-}^{k}$, we deduce
\begin{align}
\SA_{0} - h^{S ( z , h ) / \lambda_{1} - 1 / 2} &( \CP_{0} + \CO ( h^{\zeta} ) ) \SA_{0}  \nonumber \\
&= h^{S ( z , h ) / \lambda_{1} + 1 / 2} \big( \CP_{1} \SA_{1} + \CP_{2} \SA_{2} + \CO ( h^{- N + \zeta} ) \big) .  \label{g92}
\end{align}
In the previous expression, $\CP_{j} = \CP_{j} ( z , h )$ is the $K \times ( n - 1 )^{j} K$ matrix of the $\CP^{( j )}_{k} ( x_{-}^{k} , h )$'s for $j = 0 , 1, 2$. On the other hand, from \ref{h14} and (an expression similar to) \eqref{e10}, $\CP_{0} ( z , h ) = \CQ (z , h)$ is uniformly invertible in \eqref{g87}. Then, \eqref{g93} and \eqref{g92} give
\begin{align}
\SA_{0} &= \big( 1 - h^{S ( z , h ) / \lambda_{1} - 1 / 2} \CP_{0} + \CO ( h^{- 1 + \zeta } ) \big)^{-1} h^{S ( z , h ) / \lambda_{1} + 1 / 2}   \nonumber  \\
&\qquad \qquad \qquad \qquad \qquad \qquad \qquad \big( \CP_{1} \SA_{1} + \CP_{2} \SA_{2} + \CO ( h^{- N + \zeta} ) \big)  \nonumber  \\
&= - h \CP_{0}^{-1} \CP_{1} \SA_{1} - h \CP_{0}^{-1} \CP_{2} \SA_{2} + \CO ( h^{- N + 1 + \zeta} ) . \label{g94}
\end{align}
In particular, this implies $\SA_{0} = \CO ( h^{- N + 1} )$.

Let $\CP_{j}^{\prime}$ (resp. $\CP_{j}^{\prime \prime}$) denote the $( n - 1 ) K \times ( n - 1 )^{j} K$ (resp. $( n - 1 )^{2} K \times ( n - 1 )^{j} K$) matrix of the $\partial_{y_{a}} \CP^{( j )}_{k} ( x_{-}^{k} , h )$ (resp. $\partial_{y_{a} , y_{b}}^{2} \CP^{( j )}_{k} ( x_{-}^{k} , h )$) with $k \in \{ 1 , \ldots , K \}$ and $a , b \in \{ 1 , \ldots , n - 1 \}$. Differentiating \eqref{g90} with respect to $y$ and using \eqref{g94}, we deduce
\begin{align}
\SA_{1} &= h^{S ( z , h ) / \lambda_{1} + 1 / 2} ( \CP_{1}^{\prime} - \CP_{0}^{\prime} \CP_{0}^{-1} \CP_{1} ) \SA_{1}  \nonumber  \\
&\qquad \qquad \qquad \qquad + h^{S ( z , h ) / \lambda_{1} + 1 / 2} ( \CP_{2}^{\prime} - \CP_{0}^{\prime} \CP_{0}^{-1} \CP_{2} ) \SA_{2} + \CO ( h^{- N + \zeta} ) ,  \label{i1} \\
\SA_{2} &= h^{S ( z , h ) / \lambda_{1} + 1 / 2} ( \CP_{1}^{\prime \prime} - \CP_{0}^{\prime \prime} \CP_{0}^{-1} \CP_{1} ) \SA_{1}  \nonumber \\
&\qquad \qquad \qquad \qquad + h^{S ( z , h ) / \lambda_{1} + 1 / 2} ( \CP_{2}^{\prime \prime} - \CP_{0}^{\prime \prime} \CP_{0}^{-1} \CP_{2} ) \SA_{2} + \CO ( h^{- N + \zeta} ) .   \label{i2}
\end{align}
At this point of the proof, it seems that both derivatives of order $1$ and $2$ contribute to the resonances in the domain \eqref{g87}. But, the following result shows that the derivatives of order $2$ play no role.

\begin{lemma}\sl \label{g95}
We have
\begin{equation*}
\CP_{2}^{\prime} - \CP_{0}^{\prime} \CP_{0}^{-1} \CP_{2} = 0 \qquad \text{and} \qquad \CP_{2}^{\prime \prime} - \CP_{0}^{\prime \prime} \CP_{0}^{-1} \CP_{2} = 0 .
\end{equation*}
\end{lemma}

\begin{proof}
From Lemma \ref{g89} $iii)$, we can write $\partial_{x}^{\alpha} \CP^{(2)}_{k} ( x_{-}^{k} , h ) = \partial_{x}^{\alpha} \CP^{(0)}_{k} ( x_{-}^{k} , h ) \CL$ for all $\alpha \in \N^{n}$. It yields the expressions
\begin{equation*}
\CP_{2} = \CP_{0} \CL , \qquad \CP_{2}^{\prime} = \CP_{0}^{\prime} \CL \qquad \text{and} \qquad \CP_{2}^{\prime \prime} = \CP_{0}^{\prime \prime} \CL .
\end{equation*}
They imply
\begin{equation*}
\CP_{2}^{\prime} - \CP_{0}^{\prime} \CP_{0}^{-1} \CP_{2} = \CP_{0}^{\prime} \CL - \CP_{0}^{\prime} \CP_{0}^{-1} \CP_{0} \CL = 0 .
\end{equation*}
The relation $\CP_{2}^{\prime \prime} - \CP_{0}^{\prime \prime} \CP_{0}^{-1} \CP_{2} = 0$ can be proved the same way.
\end{proof}

Equation \eqref{i2} and Lemma \ref{g95} imply that $\SA_{2}$ is an explicit function of $\SA_{1}$.  Moreover, $\SA_{2}$ disappears in \eqref{i1} which becomes a closed relation on $\SA_{1}$. The quantization matrix is then defined as
\begin{equation} \label{i37}
\CQ^{2} ( z , h ) : = \CP_{1}^{\prime} - \CP_{0}^{\prime} \CP_{0}^{-1} \CP_{1} .
\end{equation}
It verifies

\begin{lemma}\sl \label{g96}
The matrix $\CQ^{2} ( z , h )$ takes the form \eqref{g91}.
\end{lemma}

\begin{proof}
Let $R = \diag ( e^{i A_{k} / h} )$ and let $R^{\prime}$ be the $( n - 1 ) K \times ( n - 1 ) K$ diagonal matrix whose $( k , a )$-th diagonal coefficient (with $k \in \{ 1 , \ldots , K \}$ and $a \in \{ 1 , \ldots n - 1 \}$) is $e^{i A_{k} / h}$. From Lemma \ref{g89} $i)$, we can write
\begin{equation*}
\CP_{\bullet}^{\star} = R^{\star} \CN_{\bullet}^{\star} ,
\end{equation*}
where $\CN_{\bullet}^{\star} = \CN_{\bullet}^{\star} ( \sigma )$ is independent of $h$ and holomorhic in $\sigma$ near $\im \sigma = - \sum_{j = 2}^{n} \lambda_{j} / 2 - \lambda_{1}$. Moreover, $\CN_{0}^{-1} ( \sigma )$ is well-defined and holomorphic thanks to \ref{h14}. Then,
\begin{equation*}
\CQ^{2} ( z , h ) = R^{\prime} \CN_{1}^{\prime} - R^{\prime} \CN_{0}^{\prime} \CN_{0}^{-1} R^{- 1} R \CN_{1} = R^{\prime} \big( \CN_{1}^{\prime} - \CN_{0}^{\prime} \CN_{0}^{-1} \CN_{1} \big) ,
\end{equation*}
which gives the lemma.
\end{proof}

As stated in Section \ref{s83}, the particular form \eqref{g91} of $\CQ^{2} ( z , h )$ allows us to adapt the proof of Proposition \ref{d9} in order to demonstrate Proposition \ref{g88}. It remains to give the

\begin{proof}[Proof of Theorem \ref{g86}]
We follow the strategy of the proof of Theorem \ref{d8}: we first show that $P$ has no resonance and a polynomial estimate of its resolvent away from the pseudo-resonances of the second kind and then prove that $P$ has a resonance near each pseudo-resonance of the second kind.

Thus, we first prove that there exists $M > 0$ such that
\begin{equation} \label{i11}
\big\Vert( P_{\theta} -z )^{-1} \big\Vert \lesssim h^{- M} ,
\end{equation}
uniformly for $z \in \eqref{g87}$ with $\dist ( z , \res_{0}^{2} ( P ) ) \geq \delta h \vert \ln h \vert^{- 1}$. From the general argument of Section \ref{s3} and the strategy of Section \ref{s72}, this problem reduces to the study of
\begin{equation} \label{i26}
\left\{ \begin{aligned}
&( P - z ) u = 0 &&\text{microlocally near } K ( E_{0} ) ,  \\
&u = 0 &&\text{microlocally near } \Lambda_{-} \setminus ( \{ 0 \} \cup \CH ) ,
\end{aligned} \right.
\end{equation}
with $z \in \eqref{g87}$, $\dist ( z , \res_{0}^{2} ( P ) ) \geq \delta h \vert \ln h \vert^{- 1}$ and $\Vert u \Vert_{L^{2} ( \R^{n} )} \leq 1$. Mimicking \eqref{d44}, one can prove that
\begin{equation*}
u_{-}^{k} \in \CI ( \Lambda_{+}^{1} , h^{-N} ) ,
\end{equation*}
for some $N \in \R$. As usual, the symbol $a_{-}^{k} \in S ( h^{- N} )$ of $u_{-}^{k}$ is defined by
\begin{equation*}
u_{-}^{k} ( x , h ) = e^{- i A_{k} / h} e^{i \frac{z - E_{0}}{h} t_{-}^{k}} \frac{\CM_{k}^{-}}{\CD_{k} ( t_{-}^{k} )} a_{-}^{k} ( x , h ) e^{i \varphi_{+}^{1} (x) / h} .
\end{equation*}
We now use the notations and the results of the beginning of this section. From \eqref{i1}, Lemma \ref{g95} and the definition of $\CQ^{2}$, we get
\begin{equation} \label{i28}
\SA_{1} = h^{S ( z , h ) / \lambda_{1} + 1 / 2} \CQ^{2} ( z , h) \SA_{1} + \CO ( h^{- N + \zeta} ) .
\end{equation}
Since $z$ is supposed at distance $h \vert \ln h \vert^{- 1}$ from the pseudo-resonances of the second kind (see Definition \ref{g85}), the previous estimate and the adaptation to the present setting of lemma \ref{d13} imply
\begin{equation} \label{i13}
\SA_{1} = \big( 1 - h^{S ( z , h ) / \lambda_{1} + 1 / 2} \CQ^{2} ( z , h) \big)^{-1} \CO ( h^{- N + \zeta} ) = \CO ( h^{- N + \zeta} ) .
\end{equation}
Combining \eqref{i2}, \eqref{i13} with Lemma \ref{g95}, we deduce $\SA_{2} = \CO ( h^{- N + \zeta} )$. Then, \eqref{g94} gives $\SA_{0} = \CO ( h^{- N + 1 + \zeta} )$. Eventually, the previous estimates on the $\SA_{\bullet}$, \eqref{g90} and \eqref{g93} yield $a_{-}^{k} \in S ( h^{- N + \zeta} )$ for some universal constant $\zeta > 0$. Consequently, $u = \CO ( h^{\infty} )$ and \eqref{i11} holds true by the bootstrap argument explained in Section \ref{s4}.

It remains to show the existence of a resonance near each pseudo-resonance of the second kind. More precisely, we have to prove
\begin{align}
\forall r > 0 , \quad \exists & h_{0} > 0 , \quad \forall h < h_{0} ,  \nonumber \\
&\forall z \in \res_{0}^{2} ( P ) \cap \eqref{g87} , \qquad \card \Big( \res (P) \cap B \Big( z , 2 r \frac{h}{\vert \ln h \vert} \Big) \Big) \geq 1 . \label{i21}
\end{align}
For that, we follow Section \ref{s73}. The main difference is the choice of the ``test function''. We first reduce \eqref{i21} using compactness arguments as in \eqref{d54}--\eqref{d59}. If \eqref{i21} does not hold true, then there exist $\tau_{0} \in [ a , b ]$, $\rho_{0} \in ( \S^{1} )^{K}$ and a non-zero eigenvalue $\mu_{0}$ of 
\begin{equation*}
\CQ_{0}^{2} = \sum_{k = 1}^{K} \rho_{0 , k} \widetilde{\CQ}^{2}_{k} \Big( \tau_{0} - i \sum_{j = 2}^{n} \frac{\lambda_{j}}{2} - i \lambda_{1} \Big) ,
\end{equation*}
such that, for some sequence of $h$ going to $0$,
\begin{equation} \label{i22}
P \text{ has no resonance in } B \Big( z_{q}^{2 , 0} , r \frac{h}{\vert \ln h \vert} \Big) ,
\end{equation}
where $\rho ( h ) = \rho_{0} + o (1)$,
\begin{equation*}
z_{q}^{2 , 0} = E_{0} + 2 q \pi \lambda_{1} \frac{h}{\vert \ln h \vert} - i h \sum_{j = 2}^{n} \frac{\lambda_{j}}{2} - i h \lambda_{1} + i \ln ( \mu_{0} ) \lambda_{1} \frac{h}{\vert \ln h \vert} ,
\end{equation*}
and $q = q ( h ) \in \Z$ is not fixed but satisfies $\re z_{q}^{2 , 0} = \tau_{0} h + o ( h )$. Note that $\CQ_{0}^{2}$ corresponds to the matrix $\CQ^{2}$ with all the parameters fixed (see \eqref{g91}).

We now construct the ``test function''. Let $w_{0} = ( w_{0}^{k , a} )_{k , a} \in \C^{( n - 1 ) K}$, for $k \in \{ 1 , \ldots , K \}$ and $a \in \{ 1 , \ldots , n - 1 \}$, be a fixed (non-zero) eigenvector of $\CQ_{0}^{2}$ associated to the eigenvalue $\mu_{0}$. We define
\begin{equation} \label{i17}
w_{0}^{k} ( y ) = \sum_{a = 1}^{n - 1} w_{0}^{k , a} y_{a} ,
\end{equation}
where $y$ denotes, as before, the $n - 1$ variables on $H_{k}$ centered at $x_{-}^{k}$. Mimicking \eqref{d89}, let $\widetilde{v}_{k} = \widetilde{a}_{k} e^{i \varphi_{+}^{1} / h}$ be a WKB solution of
\begin{equation} \label{i15}
\left\{ \begin{aligned}
&( P - \widetilde{z} ) \widetilde{v}_{k} = 0 &&\text{near } x_{-}^{k} ,   \\
&\widetilde{v}_{k} ( x ) = w_{0}^{k} ( y ( x ) ) e^{i \varphi_{+}^{1} (x) / h} &&\text{for } x \in H_{k} \text{ near } x_{-}^{k} ,
\end{aligned} \right.
\end{equation}
for $\widetilde{z} \in B ( E_{0} , R h )$ for some large $R > 1$. In particular, $\widetilde{a}_{k} ( x , h ) = \widetilde{a}_{k}^{0} ( x ) + \widetilde{a}_{k}^{1} ( x ) h + \cdots$ in $S ( 1 )$ and
\begin{equation} \label{i16}
\widetilde{a}_{k} ( x , h ) = w_{0}^{k} ( y ( x ) ) ,
\end{equation}
for $x \in H_{k}$. Consider $\CD = B ( z_{q}^{2 , 0} , s h \vert \ln h \vert^{-1} )$ with $0 < s < r$ small enough. As in \eqref{d61}, we define $u \in L^{2} ( \R^{n} )$ by
\begin{equation} \label{i12}
( P_{\theta} - \widetilde{z} ) u = v ,
\end{equation}
for $\widetilde{z} \in \partial \CD$ and $v$ constructed as in \eqref{d60}. Thanks to \eqref{i11}, $u$ is well-defined, holomorphic and polynomially bounded for $z$ in a neighborhood of $\partial \CD$. As in \eqref{d71} and \eqref{d87}, we get
\begin{equation*}
u_{-}^{k} (x) = e^{- i A_{k} / h} e^{i \frac{\widetilde{z} - E_{0}}{h} t_{-}^{k}} \frac{\CM_{k}^{-}}{\CD_{k} ( t_{-}^{k} )} a_{-}^{k} ( x , h ) e^{i \varphi_{+}^{1} (x) / h} ,
\end{equation*}
for some $a_{-}^{k} \in S ( h^{- N} )$ with $N \geq 0$. Adapting Lemma \ref{d72} and Lemma \ref{g89}, the solution $u$ of \eqref{i12} satisfies
\begin{align}
a_{-}^{k} ( x , h ) &= h^{S ( z , h ) / \lambda_{1} - 1 / 2} \big( ( \CP^{(0)}_{k} ( x , h ) + S ( h^{\zeta} ) ) \SA_{0} \nonumber \\
&\qquad \qquad + h \CP^{(1)}_{k} ( x , h ) \SA_{1} +  h \CP^{(2)}_{k} ( x , h ) \SA_{2} + S ( h^{-N + 1 + \zeta} ) \big) + \widetilde{a}_{k} ( x , h ) ,   \label{i14}
\end{align}
for all $x$ near $\pi_{x} ( U_{-}^{k} \cap \Lambda_{+}^{1} )$. Since $\widetilde{a}_{k} ( x_{-}^{k} , h ) = 0$ from \eqref{i17} and \eqref{i16}, the proof of \eqref{g94} together with \eqref{i14} give
\begin{equation} \label{i20}
\SA_{0} = - h \CP_{0}^{-1} \CP_{1} \SA_{1} - h \CP_{0}^{-1} \CP_{2} \SA_{2} + \CO ( h^{- N + 1 + \zeta} ) .
\end{equation}
On the other hand, using Lemma \ref{g95} and the particular form of $\widetilde{a}_{k}$ on $H_{k}$ (see \eqref{i17} and \eqref{i16}), the relations \eqref{i1} and \eqref{i2} are replaced by
\begin{align}
\SA_{1} &= h^{S ( \widetilde{z} , h ) / \lambda_{1} + 1 / 2} \CQ^{2} ( \widetilde{z} , h ) \SA_{1} + w_{0} + \CO ( h^{- N + \zeta} ) ,  \label{i18} \\
\SA_{2} &= h^{S ( \widetilde{z} , h ) / \lambda_{1} + 1 / 2} ( \CP_{1}^{\prime \prime} - \CP_{0}^{\prime \prime} \CP_{0}^{-1} \CP_{1} ) \SA_{1} + \CO ( h^{- N + \zeta} ) .   \label{i19}
\end{align}

Now the rest of the proof follows closely Section \ref{s73}. The expressions \eqref{i14}--\eqref{i19} play the role of \eqref{d73} and \eqref{d82} in the demonstration of Proposition \ref{d50}. Since $\widetilde{z} \in \partial \CD$ is at distance $h \vert \ln h \vert^{- 1}$ from $\res_{0}^{2} ( P )$, we can perform the bootstrap argument of \eqref{i13} and prove, as in Lemma \ref{d77}, that $N = 0$. In other words, $u_{-}^{k} \in \CI ( \Lambda_{+}^{1} , 1 )$ uniformly for $\widetilde{z} \in \partial \CD$. Then, as in the proof of Lemma \ref{d79}, \eqref{i18} yields
\begin{align*}
\SA_{1} &= \big( 1 - h^{S ( \widetilde{z} , h ) / \lambda_{1} + 1 / 2} \CQ_{0}^{2} \big)^{- 1} w_{0} + o ( 1 )   \\
&= \big( 1 - h^{S ( \widetilde{z} , h ) / \lambda_{1} + 1 / 2} \mu_{0} \big)^{- 1} w_{0} + o ( 1 ) .
\end{align*}
Thus, as in \eqref{d84}, we eventually obtain
\begin{equation*}
\frac{1}{2 i \pi} \int_{\partial \CD} u ( \widetilde{z} ) \, d \widetilde{z} = b ( x , h ) e^{i \varphi_{+}^{1} ( x ) / h},
\end{equation*}
with $b \in S ( h \vert \ln h \vert^{- 1} )$ and
\begin{align}
\partial_{y} b ( x_{-}^{k} , h ) &= e^{- i A_{k} / h} \frac{\CM_{k}^{-}}{\CD_{k} ( t_{-}^{k} )} \frac{1}{2 i \pi} \int_{\partial \CD} e^{- i \frac{\widetilde{z} - E_{0}}{h} t_{-}^{k}} \big[ \SA_{1} ( \widetilde{z} ) \big]_{k}  \, d \widetilde{z}  \nonumber \\
&= i e^{- i A_{k} / h} e^{i \sigma_{0} t_{-}^{k}} \frac{\lambda_{1} \CM_{k}^{-}}{\CD_{k} ( t_{-}^{k} )} w_{0}^{k} \frac{h}{\vert \ln h \vert} + o \Big( \frac{h}{\vert \ln h \vert} \Big) .
\end{align}
This is in contradiction with the holomorphy of the resolvent $\widetilde{z} \mapsto( P_{\theta} - \widetilde{z} )^{- 1}$ inside $\CD$. That is, there exists at least one resonance in $\CD$ and \eqref{i22} can not be satisfied.
\end{proof}

\begin{proof}[Proof of Remark \ref{i23}]
In the one dimensional case, the hypersurface $H_{k}$ of Theorem \ref{a32} is a point and the integration in \eqref{g98} is simply the value of integrand at $x_{-}^{\ell}$. Then, \eqref{g90} can be replaced by
\begin{equation} \label{i24}
a_{-}^{k} ( x , h ) = h^{- i \frac{z - E_{0}}{\lambda_{1} h}} \CP_{k} ( x , h ) \SA_{0}  + S ( h^{\infty} ) ,
\end{equation}
for all $x$ near $\pi_{x} ( U_{-}^{k} \cap \Lambda_{+}^{1} )$. Here, $\CP_{k} \in S ( 1 )$ and
\begin{equation*}
\CP_{k} ( x_{-}^{k} , h ) = \CQ_{k} ( z , h ) + \CO ( h^{\zeta} ) ,
\end{equation*}
where $\zeta > 0$ and $\CQ_{k}$ is the $k$-th line of the matrix $\CQ$ defined in \eqref{d4}.

We now apply the contradiction argument of Section \ref{s3}. Let $u$ be a solution of \eqref{i26}. Applying \eqref{i24} at $x = x_{-}^{k}$, we get
\begin{equation} \label{i25}
\Big( 1 - h^{- i \frac{z - E_{0}}{\lambda_{1} h}} \big( \CQ ( z , h ) + \CO ( h^{\zeta} ) \big) \Big) \SA_{0} = \CO ( h^{\infty} ) .
\end{equation}
Since $\Omega$ is a compact subset of $\R + i ] - \infty , 0 [$, we have $\vert h^{- i \frac{z - E_{0}}{\lambda_{1} h}} \vert \geq h^{- \nu}$ for some $\nu > 0$. Moreover, $\Vert \CQ ( z , h )^{- 1} \Vert \lesssim 1$ uniformly for $z \in E_{0} + h \Omega$ because $\det \widetilde{\CZ} \neq 0$ on $\Omega$. Then, we have
\begin{equation*}
\Big\Vert \Big( 1 - h^{- i \frac{z - E_{0}}{\lambda_{1} h}} \big( \CQ ( z , h ) + \CO ( h^{\zeta} ) \big) \Big)^{- 1} \Big\Vert \lesssim 1 .
\end{equation*}
This estimate is similar to \eqref{e93}. Combining with \eqref{i25}, it yields $\SA_{0} = \CO ( h^{\infty} )$ and then $a_{-}^{k} \in S ( h^{\infty} )$ by \eqref{i24}. Making one more turn along $K ( E_{0} )$, we eventually obtain that $u = 0$ microlocally near $K ( E_{0} )$ which gives the remark thanks to Section \ref{s3}.
\end{proof}

\begin{proof}[Proof of \eqref{i30}]
We compute $\CQ^{2} ( z , h )$ under the symmetry assumptions of Example \ref{i29}. For that, we follows the proof of Lemma \ref{g89} and give some formulas for the $\CP_{\bullet}^{\star}$'s. We decompose $x = ( x_{1} , x^{\prime} ) \in \R^{n}$. Along all this proof, we will use that a smooth function of $x$ which is stable by the rotations in the $x^{\prime}$ variables is necessarily a smooth function of $x_{1}$ and $x^{\prime} {}^{2}$. By assumption, $V$ is an example of such a function.

We first consider the stationary phase expansion in \eqref{g98}. The symmetry of $V$ and the form of the Hamiltonian vector field imply, by the previous argument, that $\varphi_{-} ( y )$ and $\varphi_{+}^{1} ( y )$ are smooth functions of $y_{1}$ and $y^{\prime} {}^{2}$. Restricting on the hypersurface $H = \{ x \in \R^{n} ; \ x_{1} = c_{1} \}$ of Theorem \ref{a32}, it yields
\begin{equation*}
\big( \varphi_{+}^{1} - \varphi_{-} \big) ( c_{1} , y^{\prime} ) = A_{1} + \alpha y^{\prime} {}^{2} + \CO ( y^{\prime} {}^{4} ) ,
\end{equation*}
for some $\alpha \neq 0$ since \ref{h8} holds true. Following the proof of \cite[Lemma 5.1]{DiSj99_01}, there exists a local diffeomorphism $\kappa$ defined near $0$ with
\begin{equation} \label{i32}
\kappa ( y^{\prime} ) = y^{\prime} + \CO ( y^{\prime} {}^{3} ) ,
\end{equation}
and such that
\begin{equation*}
( \varphi_{+}^{1} - \varphi_{-} ) \circ \kappa ( y^{\prime} ) = A_{1} + \alpha y^{\prime} {}^{2} .
\end{equation*}
Moreover, since $V_{1} ( x ) = E_{0} - \lambda_{1}^{2} x^{2} / 4$ near $0$, we deduce from \eqref{b19} that
\begin{equation} \label{i33}
d_{0} ( x , y^{\prime} , z ) = C f ( x^{2} ) g ( y^{\prime}{}^{2} ) ( x \cdot y )^{- S ( z, h ) / \lambda_{1}} ,
\end{equation}
for some smooth functions $f , g$ such that $f ( 0 ) \neq 0$ and $g ( 0 ) \neq 0$. Here and in the sequel, $C = C ( z , h )$ designs a constant which may change from line to line such that $1 \lesssim \vert C ( z , h ) \vert \lesssim 1$ uniformly for $z \in \eqref{i27}$ and $h$ small enough. Using the diffeomorphism $\kappa$, \eqref{g98} can be written
\begin{equation*}
a_{+}^{1} ( x , h ) = h^{S ( z , h ) / \lambda_{1} - n / 2} C \int_{\R^{n - 1}} e^{i \alpha y^{\prime}{}^{2} / h} d ( x , \kappa ( y^{\prime} ) , z , h ) a_{-}^{1} ( c_{1} , \kappa ( y^{\prime} ) , h ) \bigg\vert \frac{\partial \kappa ( y^{\prime} )}{\partial y^{\prime}} \bigg\vert \, d y^{\prime} .
\end{equation*}
Since the phase is now a quadratic form, we can use the explicit formula for the stationary phase method (see \cite[Page 45]{DiSj99_01}). In particular, the differential operator $B_{1}^{1}$ is equal to $\Delta_{y^{\prime}}$ modulo some non-zero multiplicative constant. The previous expression becomes
\begin{align*}
a_{+}^{1} ( x , h ) &= h^{S ( z , h ) / \lambda_{1} - 1 / 2} \bigg( S ( 1 ) a_{-}^{1} ( x_{-}^{1} , h )  \\
&\qquad \qquad + h C \Delta_{y^{\prime}} \bigg( d ( x , \kappa ( y^{\prime} ) , z , h ) a_{-}^{1} ( c_{1} , \kappa ( y^{\prime} ) , h ) \bigg\vert \frac{\partial \kappa ( y^{\prime} )}{\partial y^{\prime}} \bigg\vert \bigg) \Big\vert_{y^{\prime} = 0} + S ( h^{2} ) \bigg) .
\end{align*}
Then, using \eqref{i32}, \eqref{i33}, \eqref{b18} and the notations of \eqref{i3}, we obtain
\begin{equation} \label{i35}
\CC_{a} ( x , h ) = C x_{a} f ( x^{2} ) x_{1}^{- 1 - S ( z, h ) / \lambda_{1}} ,
\end{equation}
for $a = 2 , \ldots , n$. On the other hand, using again the symmetry of the potential $V$, the functions $\widetilde{x} ( x )$ and $\CS_{1} ( x , h )$ defined in \eqref{d52} can be written, for all $x \in H$,
\begin{equation} \label{i36}
\widetilde{x} ( x ) = x_{+}^{1} + k ( x^{\prime} {}^{2} ) x^{\prime}  \qquad \text{and} \qquad \CS_{1} ( x , h ) = \CS_{1} ( x_{+}^{1} , h ) + x^{\prime} {}^{2} \ell ( x^{\prime} {}^{2} , h ) ,
\end{equation}
where $k$ is a smooth $n \times ( n - 1 )$ matrix such that $k ( 0 ) = \beta {}^{t} ( 0 , Id )$ with $\beta \neq 0$ and $\ell$ is a smooth function.

Combining \eqref{i10} with \eqref{i35} and \eqref{i36}, the quantities $\CP_{1}^{\star}$ defined at the beginning of the section (see \eqref{g92}--\eqref{i2}) verify
\begin{equation*}
\CP_{1} = 0  \qquad \text{and} \qquad \CP_{1}^{\prime} = C Id .
\end{equation*}
Then, \eqref{i37} implies that $\CQ^{2} ( z , h ) = q^{2} ( z , h ) Id$ for some $q^{2} ( z , h )$ satisfying $1 \lesssim \vert q^{2} ( z , h ) \vert \lesssim 1$ for $z \in \eqref{i27}$ and $h$ small enough. Since $\CQ^{2}$ also satisfies \eqref{g91}, the function $q^{2} ( z , h )$ must be as in \eqref{i30}.
\end{proof}

\Subsection{Proof of the asymptotic of higher order} \label{s80}

This part is devoted to the proof of the results stated in Section \ref{s79}. Since there is only one homoclinic trajectory from \ref{h9}, we remove the subscripts $k , \ell =1$ which were used in the equations to indicate the number of the trajectory.

We explain how $\CQ_{\rm tot}$ is constructed. The first step is to obtain the complete asymptotic (i.e. up to order $\CO ( h^{\infty} )$) in Lemma \ref{d41}. Recall that this lemma gives a closed relation on the symbol of the restriction $u_{-} \in \CI ( \Lambda_{+}^{1} , h^{-N} )$ of a solution $u$ of \eqref{a4}. We send back the reader to Section \ref{s4} and Section \ref{s72} for notations and some preliminary results.

\begin{lemma}\sl \label{g11}
There exist symbols $\CP^{( \alpha )} \in S ( h^{\vert \alpha \vert / 2} )$ independent of $u$ such that
\begin{equation} \label{g12}
a_{-} ( x , h ) \simeq h^{S ( z , h ) / \lambda_{1} - 1 / 2} \sum_{\alpha \in \N^{n - 1}} \CP^{(\alpha)} ( x , h ) \partial_{y}^{\alpha} a_{-} ( x_{-} , h ) ,
\end{equation}
for all $x$ near $\pi_{x} ( U_{-} \cap \Lambda_{+}^{1} )$. Here, $\partial_{y}$ denote the derivatives in the $n - 1$ directions of $H$, the hyperplane of Theorem \ref{a32}. Moreover, the symbols $\CP^{(\alpha)}$ have the form
\begin{equation*}
\CP^{(\alpha)} ( x , h ) \simeq \sum_{\widehat{\mu}_{a} \geq \vert \alpha \vert \lambda_{1} / 2} \sum_{c = 0}^{C_{a}} \CP^{(\alpha)}_{a , c} ( x , z , h ) ( \ln h )^{c} h^{\widehat{\mu}_{a} / \lambda_{1}} ,
\end{equation*}
with $\CP^{(\alpha)}_{a , c} ( x , z , h ) = e^{i A / h} \widetilde{\CP}^{(\alpha)}_{a , c} ( x , \sigma )$, $\widetilde{\CP}^{(\alpha)}_{a , c} \in S (1)$ holomorphic with respect to $\sigma$ and $\CP^{( 0 )} = \CP + S ( ( \ln h )^{C_{1}} h^{\widehat{\mu}_{1} / \lambda_{1}} )$ where $\CP$ is given in Lemma \ref{d41}.
\end{lemma}

For shorter expressions, the dependence of $\CP^{(\alpha)}$, $a_{-}$, $\ldots$ in $z$ is not made explicit in the notations. Note that the leading term in \eqref{g12} coincides with \eqref{d42}. As usual, \eqref{g12} means that the difference between $a_{-} ( x , h )$ and the sum on the right hand side up to $\alpha$ belongs to $S ( h^{\sum_{j=2}^{n} \frac{\lambda_{j}}{2 \lambda_{1}} + \frac{\im z}{\lambda_{1} h} + \vert \alpha \vert / 2 - N} )$.

\begin{proof}
The proof of this result is similar to the one of Lemma \ref{d41} or Lemma \ref{g89}. We start with the propagation through the fixed point $( 0 , 0 )$. To solve the microlocal Cauchy problem \eqref{d62}, we apply Theorem \ref{a32} instead of Corollary \ref{d46}. Performing a stationary phase expansion in \eqref{a31}, \eqref{d51} is replaced by
\begin{equation} \label{g9}
a_{+} ( x , h ) \simeq  h^{S ( z , h ) / \lambda_{1} - 1 / 2} \sum_{\alpha \in \N^{n - 1}} \widetilde{\CR}^{( \alpha )} ( x , h ) \partial_{y}^{\alpha} a_{-} ( x_{-} , h ) ,
\end{equation}
in the same sense than \eqref{g12}. Here, $\widetilde{\CR}^{( \alpha )} \in S ( h^{\vert \alpha \vert / 2} )$ has the form
\begin{equation*}
\widetilde{\CR}^{( \alpha )} ( x , h ) \simeq \sum_{\widehat{\mu}_{a} \geq \vert \alpha \vert \lambda_{1} / 2} \sum_{b = 0}^{B^{(\alpha)}_{a}} \widetilde{\CR}_{a , b}^{( \alpha )} ( x , z , h ) ( \ln h )^{b} h^{\widehat{\mu}_{a} / \lambda_{1}} ,
\end{equation*}
and $\widetilde{\CR}_{a , b}^{( \alpha )}$ can be written $\widetilde{\CR}_{a , b}^{( \alpha )} ( x , z , h ) = e^{i A / h} \widehat{\CR}_{a , b}^{( \alpha )} ( x , \sigma )$ with $\widehat{\CR}_{a , b}^{( \alpha )} \in S ( 1 )$. Eventually, Corollary \ref{d46} gives $B^{(0)}_{0} = 0$ and $\widetilde{\CR}^{( 0 )}_{0 , 0} = \widetilde{\CR}$ where $\widetilde{\CR}$ is defined in \eqref{d51}.

We now consider the microlocal Cauchy problem \eqref{d65} which describes the propagation through the Hamiltonian trajectories $\CH$. As in \eqref{d52}, the usual propagation of Lagrangian distributions gives 
\begin{equation} \label{g10}
e^{- i A / h} e^{i \frac{z - E_{0}}{h} t_{-}} \frac{\CM^{-}}{\CD ( t_{-} )} a_{-} ( x , h ) \simeq \sum_{\alpha \in \N^{n}} \CS^{( \alpha )} ( x , h ) \partial_{x}^{\alpha} a_{+} ( \widetilde{x} (x)  , h ) ,
\end{equation}
in the same sense than \eqref{g12}. The symbols $\CS^{( \alpha )} \in S ( h^{\vert \alpha \vert / 2} )$ satisfy the asymptotic 
\begin{equation*}
\CS^{( \alpha )} ( x , h ) \simeq \sum_{a \geq \vert \alpha \vert / 2} \CS_{a}^{( \alpha )} ( x , z , h ) h^{a} ,
\end{equation*}
with $\CS_{a}^{( \alpha )} ( x , z , h ) = \widetilde{\CS}_{a}^{( \alpha )} ( x , \sigma )$ and $\widetilde{\CS}_{a}^{( \alpha )} \in S (1)$. Eventually, we have $\CS_{0}^{( 0 )} = \CS$ where $\CS$ is defined in \eqref{d52}.

Combining \eqref{g9} and \eqref{g10} with the computation \eqref{d76}, we deduce Lemma \ref{g11}.
\end{proof}

We now construct the coefficients $\CQ_{a , b , c}$ of $\CQ_{\rm tot}$. The idea is to consider formally the relation \eqref{g12} and to transform it into a closed relation on $a_{-} ( x_{-} , h )$ only. Let $\CA \in \N \setminus \{ 0 \}$. Keeping only the derivatives up to order $\CA$ in \eqref{g12} leads to
\begin{equation} \label{g15}
a_{-} ( x , h ) = h^{S ( z , h ) / \lambda_{1} - 1 / 2} \sum_{0 \leq \vert \alpha \vert \leq \CA} \CP^{(\alpha)} ( x , h ) \partial_{y}^{\alpha} a_{-} ( x_{-} , h ) + S ( h^{\CA / 2 - N} ) .
\end{equation}
After derivation, this formula gives
\begin{equation*}
\partial_{y}^{\beta} a_{-} ( x_{-} , h ) = h^{S ( z , h ) / \lambda_{1} - 1 / 2} \sum_{0 \leq \vert \alpha \vert \leq \CA} \partial_{y}^{\beta} \CP^{(\alpha)} ( x_{-} , h ) \partial_{y}^{\alpha} a_{-} ( x_{-} , h ) + \CO ( h^{\CA / 2 - N} ) ,
\end{equation*}
for all $\beta \in \N^{n - 1}$. In particular, we have
\begin{align}
\partial_{y}^{\beta} a_{-} ( x_{-} , h ) -{}& h^{S ( z , h ) / \lambda_{1} - 1 / 2} \sum_{\vert \alpha \vert = \CA} \partial_{y}^{\beta} \CP^{(\alpha)} ( x_{-} , h ) \partial_{y}^{\alpha} a_{-} ( x_{-} , h )  \nonumber \\
&= h^{S ( z , h ) / \lambda_{1} - 1 / 2} \sum_{0 \leq \vert \alpha \vert \leq \CA - 1} \partial_{y}^{\beta} \CP^{(\alpha)} ( x_{-} , h ) \partial_{y}^{\alpha} a_{-} ( x_{-} , h ) + \CO ( h^{\CA / 2 - N} ) , \label{e96}
\end{align}
for all $\vert \beta \vert = \CA$. Since $h^{S ( z , h ) / \lambda_{1} - 1 / 2} \partial_{y}^{\beta} \CP^{(\alpha)} ( x_{-} , h ) = \CO ( h^{\CA / 2} )$ for $\vert \alpha \vert = \CA$, we can inverse the previous equations and obtain
\begin{equation} \label{g45}
\partial_{y}^{\beta} a_{-} ( x_{-} , h ) = \sum_{0 \leq \vert \alpha \vert \leq \CA - 1} \CR^{(\alpha)}_{\CA , \beta} \partial_{y}^{\alpha} a_{-} ( x_{-} , h ) + \CO ( h^{\CA / 2 - N} ) ,
\end{equation}
for all $\vert \beta \vert = \CA$. As in \eqref{e95}, $\CR^{(\alpha)}_{\CA , \beta} \in S ( h^{\CA / 2 - \varepsilon} )$, for all $\varepsilon > 0$, satisfy
\begin{equation} \label{g16}
\CR^{(\alpha)}_{\CA , \beta} ( z , h ) \simeq \sum_{\widehat{\mu}_{a} \geq \CA \lambda_{1} / 2} \sum_{b = 0}^{B_{a}} \sum_{c = 0}^{C_{a}} ( \CR^{(\alpha)}_{\CA , \beta} )_{a , b , c} ( z , h ) \big( h^{S ( z , h ) / \lambda_{1} - 1 / 2} \big)^{b} ( \ln h )^{c} h^{\widehat{\mu}_{a} / \lambda_{1}} ,
\end{equation}
where $( \CR^{(\alpha)}_{\CA , \beta} )_{a , b , c}$ has the form \eqref{g14}. Inserting \eqref{g45} in \eqref{g15}, we have
\begin{equation} \label{g13}
a_{-} ( x , h ) = h^{S ( z , h ) / \lambda_{1} - 1 / 2} \sum_{0 \leq \vert \alpha \vert \leq \CA - 1} \widehat{\CP}^{(\alpha)}_{\CA} ( x , h ) \partial_{y}^{\alpha} a_{-} ( x_{-} , h ) + S ( h^{\CA / 2 - N} ) ,
\end{equation}
where $\widehat{\CP}^{(\alpha)}_{\CA} ( x , h )$ is like \eqref{e95} except that the coefficients $\widetilde{\CQ}_{a , b , c}$ are symbols in $x$ and that the sum over $a \in \N$ is restricted to $\widehat{\mu}_{a} \geq \vert \alpha \vert \lambda_{1} / 2$. Iterating the process \eqref{g15}--\eqref{g13}, we eventually obtain
\begin{equation} \label{g18}
a_{-} ( x , h ) = h^{S ( z , h ) / \lambda_{1} - 1 / 2} \widehat{\CP}^{(0)}_{\CA} ( x , h ) a_{-} ( x_{-} , h ) + S ( h^{\CA / 2 - N} ) ,
\end{equation}
for some new $\widehat{\CP}^{(0)}_{\CA}$ satisfying the same properties as before. At this stage, we can no longer inverse $1 - h^{S ( z , h ) / \lambda_{1} - 1 / 2} \widehat{\CP}^{(0)}_{\CA} ( x_{-} , h )$ as in \eqref{e96}, since $\widehat{\CP}^{(0)}_{\CA}$ is of order $1$.

Thus, we have constructed a function $\widehat{\CP}^{(0)}_{\CA} ( x_{-} , h )$ which satisfies \eqref{e95}. A priori, it depends on $\CA$. Nevertheless, thanks to \eqref{g16}, we have
\begin{equation*}
\widehat{\CP}^{(0)}_{\CA} ( x_{-} , h ) = \widehat{\CP}^{(0)}_{\CB} ( x_{-} , h ) + \CO \big( h^{\min ( \CA , \CB ) / 2} \big) .
\end{equation*}
In other words, each coefficient in the asymptotic of $\widehat{\CP}^{(0)}_{\CA} ( x_{-} , h )$ is independent of $\CA$ for $\CA$ large enough. We call $\CQ_{a , b , c} ( z , h )$ this coefficient. Using Borel's lemma, we can finally find a symbol $\CQ_{\rm tot}$ satisfying the properties stated below \eqref{e95}. Note that it is not unique.

\begin{proof}[Proof of Proposition \ref{g1}]
First, the pseudo-resonances at infinite order satisfy the asymptotic stated in Proposition \ref{d9}. On can indeed adapt the proof of this result since $\CQ_{\rm tot} ( z , h )$ is holomorphic in \eqref{d12} and satisfies \eqref{e99}. We omit the details.

It remains to show that the pseudo-resonances at infinite order near any $z_{q} ( \tau )$ is unique. Thus, taking into account the expression of the $z_{q} ( \tau )$'s given in \eqref{e8}, it is enough to show that
\begin{equation}
\exists \alpha > 0 , \quad \forall z_{1} \neq z_{2} \in \res_{\infty} (P) \cap \eqref{d12} , \qquad \vert z_{1} - z_{2} \vert \geq \alpha \frac{h}{\vert \ln h \vert} . \label{e98}
\end{equation}
Assume that \eqref{e98} does not hold true. Then, there exists a sequence of $z_{1} \neq z_{2} \in \res_{\infty} (P) \cap \eqref{d12}$ such that $z_{1} - z_{2} = o ( h \vert \ln h \vert^{- 1} )$. As explained in the beginning of the proof of Proposition \ref{d9}, we have in fact $z_{1} , z_{2} \in \eqref{d21}$. We write
\begin{equation*}
\partial_{z} \big( h^{S ( z , h ) / \lambda_{1} - 1 / 2} \CQ_{\rm tot} ( z , h ) \big) = i \frac{\vert \ln h \vert}{\lambda_{1} h} h^{S ( z , h ) / \lambda_{1} - 1 / 2} \CQ_{\rm tot} ( z , h ) + h^{S ( z , h ) / \lambda_{1} - 1 / 2} \partial_{z} \CQ_{\rm tot} ( z , h ) .
\end{equation*}
From \eqref{e95}--\eqref{g14}, we deduce $\partial_{z}^{k} \CQ_{\rm tot} ( z , h ) = \CO ( h^{- k} )$ for all $k \in \N$. Moreover, \eqref{g4} and \eqref{e99} imply $\vert \CQ_{\rm tot} ( z , h ) \vert \gtrsim 1$. Eventually, as in Section \ref{s13}, we have
\begin{equation*}
1 \lesssim \big\vert h^{S ( z , h ) / \lambda_{1} - 1 / 2} \big\vert = h^{\sum_{j=2}^{n} \frac{\lambda_{j}}{2 \lambda_{1}} + \frac{\im z}{\lambda_{1} h}} \lesssim 1 ,
\end{equation*}
for all $z \in \eqref{d21}$. Then,
\begin{equation} \label{g5}
\big\vert \partial_{z} \big( h^{S ( z , h ) / \lambda_{1} - 1 / 2} \CQ_{\rm tot} ( z , h ) \big) \big\vert \geq \beta \frac{\vert \ln h \vert}{h} ,
\end{equation}
for some $\beta > 0$, $h$ small enough and $z \in \eqref{d21}$. On the other hand,
\begin{equation} \label{g6}
\partial_{z}^{k} \big( h^{S ( z , h ) / \lambda_{1} - 1 / 2} \CQ_{\rm tot} ( z , h ) \big) = \CO \Big( \frac{\vert \ln h \vert^{k}}{h^{k}} \Big) ,
\end{equation}
for all $k \in \N$ and $z \in \eqref{d21}$. Using \eqref{g5} and \eqref{g6}, the Taylor formula yields
\begin{align}
0 &= \big\vert h^{S ( z_{1} , h ) / \lambda_{1} - 1 / 2} \CQ_{\rm tot} ( z_{1} , h ) - h^{S ( z_{2} , h ) / \lambda_{1} - 1 / 2} \CQ_{\rm tot} ( z_{2} , h ) \big\vert   \nonumber \\
&\geq \beta \frac{\vert \ln h \vert}{h} \vert z_{1} - z_{2} \vert + \CO \Big( \frac{\vert \ln h \vert^{2}}{h^{2}} \Big) \vert z_{1} - z_{2} \vert^{2} , \label{g19}
\end{align}
and then
\begin{equation*}
\beta = \CO \Big( \frac{\vert \ln h \vert}{h} \Big) \vert z_{1} - z_{2} \vert = o \Big( \frac{\vert \ln h \vert}{h} \frac{h}{\vert \ln h \vert} \Big) = o ( 1 ) ,
\end{equation*}
which is absurd. Thus, we get \eqref{e98} by contradiction and the proposition follows.
\end{proof}

\begin{proof}[Proof of Theorem \ref{g3}]
From Theorem \ref{d8} and \eqref{g7}, we already know that
\begin{equation*}
\dist \big( \res (P) , \res_{\infty} (P) \big) = o \Big( \frac {h}{\vert \ln h \vert} \Big) ,
\end{equation*}
in the domain \eqref{d90}. Moreover, Proposition \ref{g1} implies that the distance between two pseudo-resonances at infinite order is at least $\pi \lambda_{1} h \vert \ln h \vert^{-1}$. Thus, to complete the proof of Theorem \ref{g3}, it is enough to show \eqref{g8}. Since a part of this estimate has already been obtained in Theorem \ref{d8}, it remains to show that \eqref{g8} holds true if
\begin{equation} \label{g17}
h^{C} \leq \dist ( z , \res_{\infty} ( P ) ) \leq \nu \frac{h}{\vert \ln h \vert} ,
\end{equation}
for some $\nu > 0$.

To prove this polynomial estimate, we follow the strategy of Section \ref{s72}. If it does not hold true, there exists a solution $u$ of \eqref{a4}. Following the general reduction of Section \ref{s3} and the first part of Section \ref{s72}, the restriction $u_{-} \in \CI ( \Lambda_{+}^{1} , h^{-N} )$ of $u$ satisfies Lemma \ref{g11}.

Let $\CA \geq 2 C + 2$. Using the computations made in \eqref{g18}, we deduce
\begin{equation*}
a_{-} ( x_{-} , h ) = h^{S ( z , h ) / \lambda_{1} - 1 / 2} \widehat{\CP}^{(0)}_{\CA} ( x_{-} , h ) a_{-} ( x_{-} , h ) + \CO ( h^{\CA / 2 - N} ) .
\end{equation*}
By construction of $\CQ_{\rm tot}$, we have $\widehat{\CP}^{(0)}_{\CA} ( x_{-} , h ) = \CQ_{\rm tot} ( z , h ) + \CO ( h^{\CA / 2} )$. The last equation becomes
\begin{equation} \label{g20}
\big( 1 - h^{S ( z , h ) / \lambda_{1} - 1 / 2} \CQ_{\rm tot} ( z , h ) \big) a_{-} ( x_{-} , h ) = \CO ( h^{\CA / 2 - N} ) .
\end{equation}
From \eqref{g17}, there exists  $\widetilde{z} \in \res_{\infty} ( P )$ such that $h^{C} \leq \vert \widetilde{z} - z \vert \leq \nu h \vert \ln h \vert^{- 1}$. As in \eqref{g19}, we have
\begin{align*}
\big\vert 1 - h^{S ( z , h ) / \lambda_{1} - 1 / 2} \CQ_{\rm tot} ( z , h ) \big\vert &= \big\vert h^{S ( \widetilde{z} , h ) / \lambda_{1} - 1 / 2} \CQ_{\rm tot} ( \widetilde{z} , h ) - h^{S ( z , h ) / \lambda_{1} - 1 / 2} \CQ_{\rm tot} ( z , h ) \big\vert    \\
&\geq \beta \frac{\vert \ln h \vert}{h} \vert \widetilde{z} - z \vert + \CO \Big( \frac{\vert \ln h \vert^{2}}{h^{2}} \Big) \vert \widetilde{z} - z \vert^{2}     \\
&\geq \frac{\vert \ln h \vert}{h} \vert \widetilde{z} - z \vert \big( \beta + \CO ( \nu ) \big) \geq \frac{\beta}{2} \frac{\vert \ln h \vert}{h} \vert \widetilde{z} - z \vert \gtrsim h^{C} ,
\end{align*}
for $\nu$ small enough. Then, \eqref{g20} gives $a_{-} ( x_{-} , h ) = \CO ( h^{\CA / 2 - N - C} ) = \CO ( h^{- N + 1} )$. Applying again \eqref{g18}, we obtain
\begin{equation*}
u_{-} \in \CI ( \Lambda_{+}^{1} , h^{- N + 1} ) .
\end{equation*}
Thus, we have gained a factor $h$ in the order of $u_{-}$. The rest of the proof follows from the usual bootstrap argument at the end of Section \ref{s4}.
\end{proof}

\Subsection{Proof for tangential intersection of finite order} \label{s82}

This part is devoted to the proofs of the results of Section \ref{s81}. We first demonstrate the geometric statements of that section.

\subsubsection{Proof of the geometric assertions} The manifold $\Lambda_{+}^{1}$ (the part of $\Lambda_{+}$ after a turn which is defined in \eqref{a38}) projects nicely on the $x$-space near $\pi_{x} ( \gamma_{k} \cap \Lambda_{-}^{0} )$. This follows from Proposition C.1 of \cite{ALBoRa08_01} when $\Lambda_{-}^{0}$ and $\Lambda_{+}^{1}$ intersect transversally along $\gamma_{k}$ and from \eqref{a57} when $\Lambda_{-}^{0}$ and $\Lambda_{+}^{1}$ are tangent along $\gamma_{k}$. Thus, the phase functions $\varphi_{+}^{k}$ are well-defined. Moreover, the assumption \ref{h13} implies \eqref{g21}. In the tangent case (i.e. when $m_{k} \geq 2$), we have $\Delta \varphi_{+}^{k} = \Delta \varphi_{-}$ on $\pi_{x} ( \gamma_{k} )$ which implies that the limit \eqref{g44} belongs in $] 0 , + \infty [$. We now show the comparison formula \eqref{g72}.

\begin{proof}[Proof of \eqref{g72}] We use simultaneously the notations of Section \ref{s77}, Section \ref{s81} and Section \ref{s53}. Let $T_{k} \in \R$ be such that $x_{+} ( t , \widehat{g}_{+}^{k} ) = x_{k} ( t + T_{k} )$. By definition of $t_{+}^{\varepsilon} ( \widehat{g}_{+}^{k} )$, we have $\varepsilon = \vert x_{+} ( t_{+}^{\varepsilon} ( \widehat{g}_{+}^{k} ) , \widehat{g}_{+}^{k}) \vert = \vert x_{k} ( t_{+}^{\varepsilon} ( \widehat{g}_{+}^{k} ) + T_{k} ) \vert$. Then, \eqref{d6} implies
\begin{equation} \label{g74}
\vert g_{+}^{k} \vert e^{\lambda_{1} ( t_{+}^{\varepsilon} ( \widehat{g}_{+}^{k} ) + T_{k} )} = \varepsilon ( 1 + o_{\varepsilon \to 0} ( 1 ) ) .
\end{equation}
The same way,
\begin{equation} \label{g75}
\vert g_{-}^{k} \vert e^{- \lambda_{1} ( t_{-}^{\varepsilon} ( \widehat{g}_{+}^{k} ) + T_{k} )} = \varepsilon ( 1 + o_{\varepsilon \to 0} ( 1 ) ) .
\end{equation}
On the other hand, from \eqref{a81} and \eqref{g73}, we get
\begin{equation*}
\CM_{\varepsilon} ( \widehat{g}_{+}^{k} ) = \frac{\CD_{k} ( t_{+}^{\varepsilon} ( \widehat{g}_{+}^{k} ) + T_{k} ) e^{- ( t_{+}^{\varepsilon} ( \widehat{g}_{+}^{k} ) + T_{k} ) \sum_{j} \lambda_{j} / 2}}{\CD_{k} ( t_{-}^{\varepsilon} ( \widehat{g}_{+}^{k} ) + T_{k} ) e^{( t_{-}^{\varepsilon} ( \widehat{g}_{+}^{k} ) + T_{k} ) \sum_{j} \lambda_{j} / 2}} e^{( t_{+}^{\varepsilon} ( \widehat{g}_{+}^{k} ) + t_{-}^{\varepsilon} ( \widehat{g}_{+}^{k} ) + 2 T_{k}) \sum_{j} \lambda_{j} / 2} .
\end{equation*}
Then, the definitions of $\CM_{k}^{\pm}$ and $\CM_{0}$ give
\begin{equation} \label{g76}
\frac{\CM_{k}^{+}}{\CM_{k}^{-}} = \CM_{0} ( \widehat{g}_{+}^{k} ) e^{- ( t_{+}^{\varepsilon} ( \widehat{g}_{+}^{k} ) + t_{-}^{\varepsilon} ( \widehat{g}_{+}^{k} ) + 2 T_{k} ) \sum_{j} \lambda_{j} / 2} ( 1 + o_{\varepsilon \to 0} ( 1 ) ) .
\end{equation}
Combining \eqref{g74}, \eqref{g75} and \eqref{g76} leads to
\begin{align*}
&\frac{\CM_{k}^{+}}{\CM_{k}^{-}} \big\vert g_{-}^{\ell} \big\vert^{\frac{\lambda_{1} + \lambda_{2}}{\lambda_{1}}} \big( i \lambda_{1} g_{+}^{k} \cdot g_{-}^{\ell} \big)^{- S ( z , h ) / \lambda_{1}}     \\
&\ = \CM_{0} ( \widehat{g}_{+}^{k} ) \big( i \lambda_{1} \widehat{g}_{+}^{k} \cdot \widehat{g}_{-}^{\ell} \big)^{- S ( z , h ) / \lambda_{1}} \frac{\vert g_{-}^{\ell} \vert^{\frac{\lambda_{1} + \lambda_{2}}{2 \lambda_{1}} + i \frac{z - E_{0}}{\lambda_{1} h}}}{\vert g_{-}^{k} \vert^{\frac{\lambda_{1} + \lambda_{2}}{2 \lambda_{1}} + i \frac{z - E_{0}}{\lambda_{1} h}}} e^{i ( t_{-}^{\varepsilon} ( \widehat{g}_{+}^{k} ) - t_{+}^{\varepsilon} ( \widehat{g}_{+}^{k} ) ) \frac{z - E_{0}}{h}} \varepsilon^{2 i \frac{z - E_{0}}{\lambda_{1} h}} ( 1 + o_{\varepsilon \to 0} ( 1 ) ) .
\end{align*}
Remark that this quantity is independent of $\varepsilon$. Thus, taking the limit $\varepsilon \to 0$ and using the definition \eqref{i50} of the time delay, we eventually obtain \eqref{g72}.
\end{proof}

It remains to give the

\begin{proof}[Proof of \eqref{g22}]
Consider $k \in \{ 1 , \ldots , K \}$ such that $m_{k} \geq 2$. As in Section \ref{a71}, we project the Hamiltonian flow $H_{p}$ restricted to $\Lambda_{-}^{0}$ onto $\R^{2}$. Let
\begin{equation*}
H_{p}^{-} = 2 \nabla \varphi_{-} ( x ) \cdot \partial_{x} ,
\end{equation*}
be the Hamiltonian vector field restricted to $\Lambda_{-}^{0}$. This means that $( x (t) , \xi (t) )$ is a Hamiltonian curve in $\Lambda_{-}^{0}$ if and only if $x (t)$ is an integral curve of $H_{p}^{-}$ and $\xi (t) = \nabla \varphi_{-} ( x (t) )$. The flow of $H_{p}^{-}$ is contracting since $2 \nabla \varphi_{-} ( x ) = ( - \lambda_{1} x_{1} , - \lambda_{2} x_{2} ) + \CO ( x^{2} )$ from \eqref{a57}. Then, the Hartman Theorem (see \cite[Page 127]{Pe01_01}) provides a $C^{1}$-diffeomorphism $\Phi : \R^{2} \to \R^{2}$ defined in a  neighborhood of $0$ such that $\Phi ( 0 ) = 0$, $d \Phi ( 0 ) = Id_{\R^{2}}$ and
\begin{equation} \label{g31}
\Phi \big( \exp ( t H_{p}^{-} ) ( x ) \big) = e^{- L t} \Phi ( x ) ,
\end{equation}
with $L = \diag ( \lambda_{1} , \lambda_{2} )$.

For $t_{0} > 0$ large enough, let $\delta^{0} = ( \delta_{x}^{0} , \delta_{\xi}^{0} ) \in T_{\gamma_{k} ( t_{0} )} \Lambda_{-}^{0} = T_{\gamma_{k} ( t_{0} )} \Lambda_{+}^{1}$ be the tangent vector defined by
\begin{equation} \label{g30}
\delta_{x}^{0} = \big( d \Phi ( x_{k} ( t_{0} ) ) \big)^{- 1} ( 0 , 1 ) \qquad \text{and} \qquad \delta_{\xi}^{0} = \Hess \varphi_{-} ( x_{k} ( t_{0} ) ) \delta_{x}^{0} .
\end{equation}
Let $\rho_{0} (s) : \R \to \Lambda_{+}^{1}$ be a smooth curve defined near $0$ such that $\rho_{0} ( 0 ) = \gamma_{k} ( t_{0} )$ and $\partial_{s} \rho_{0} ( 0 ) = \delta^{0}$. Of course, this curve can not be a Hamiltonian trajectory since $\delta^{0}$ is almost orthogonal to the Hamiltonian vector field. We now set, for $t \geq 0$,
\begin{equation} \label{g29}
\rho ( t , s ) = ( x ( t , s ) , \xi ( t , s ) ) : = \exp ( t H_{p} ) ( \rho_{0} ( s ) ) .
\end{equation}
The idea is that the behavior of $\alpha_{k} ( x_{1} )$ as $x_{1}$ goes to $0$ is given by the derivatives of $\rho$.

Since $\rho ( t , 0 ) = \gamma_{k} ( t )$ is expandible from \cite[Remark 3.10]{HeSj85_01}, the choice of coordinates and the notations of \eqref{d6} imply
\begin{equation} \label{g39}
x ( t , 0 ) = \big( \vert g_{-}^{k} \vert e^{- \lambda_{1} t_{0}} e^{- \lambda_{1} t} , 0 \big) + \CO (  e^{- ( \lambda_{1} + \varepsilon ) t} ) ,
\end{equation}
for some $\varepsilon > 0$, and
\begin{equation} \label{g32}
\partial_{t} x ( t , 0 ) = 2 \nabla \varphi_{-} ( x ( t , 0 ) ) = \big( - \lambda_{1} \vert g_{-}^{k} \vert e^{- \lambda_{1} t_{0}} e^{- \lambda_{1} t} , 0 \big) + o (  e^{- \lambda_{1} t} ) ,
\end{equation}
as $t \to + \infty$. We now study the large time behavior of the tangent vector $\partial_{s} \rho ( t , 0 ) \in T_{\gamma_{k} ( t + t_{0} )} \Lambda_{+}^{1} = T_{\gamma_{k} ( t + t_{0} )} \Lambda_{-}^{0}$. Using \eqref{g31}--\eqref{g29}, we get
\begin{align}
\partial_{s} x ( t , 0 ) &= d \pi_{x} d \exp ( t H_{p} ) ( \delta^{0} ) = d \Phi^{-1} e^{- L t} d \Phi ( \delta^{0}_{x} ) \nonumber   \\
&= d \Phi^{-1} e^{- L t} ( 0 , 1 ) = ( 0 , e^{- \lambda_{2} t} ) + o ( e^{- \lambda_{2} t} ) . \label{g33}
\end{align}
In other words, \eqref{g32} and \eqref{g33} show that the derivatives $\partial_{t} x ( t , 0 )$ and $\partial_{s} x ( t , 0 )$ correspond roughly speaking to the two base vectors $( 1 , 0 )$ and $( 0 , 1 )$ respectively.

We now compute the behavior of $\varphi^{k}_{+} ( x ( t , s ) ) - \varphi_{-} ( x ( t , s ) )$. For that, we construct new symplectic coordinates as in \cite[Lemma 2.1]{HeSj85_01} or \cite[Proposition C.1]{ALBoRa08_01}. After the first symplectic diffeomorphism $F_{1} : ( x , \xi ) \longmapsto ( x , \xi - \nabla \varphi_{-} ( x ) )$, the manifold $\Lambda_{+}^{0}$ becomes $\{ ( x , \nabla \varphi_{+} ( x ) - \nabla \varphi_{-} ( x ) ) \} = \{ ( \nabla \psi_{+} ( \xi ) , \xi ) \}$ for some $\psi_{+} \in C^{\infty} ( \R^{n} )$. We then define the symplectic diffeomorphism $F_{2} : ( x , \xi ) \longmapsto ( x - \nabla \psi_{+} ( \xi ) , \xi )$ and introduce the new symplectic local coordinates $( y , \eta ) = F_{2} \circ F_{1} ( x , \xi )$. In particular,
\begin{equation} \label{g35}
y_{j} = \frac{x_{j}}{2} - \frac{\xi_{j}}{\lambda_{j}} + \CO ( ( x , \xi )^{2} ) \qquad \text{and} \qquad \eta_{j} = \xi_{j} + \frac{\lambda_{j}}{2} x_{j} + \CO ( ( x , \xi )^{2} ) .
\end{equation}
In these coordinates, $\Lambda_{-}^{0}$ is given by $\xi = 0$ and $\Lambda_{+}^{0}$ is given by $x = 0$. Then, $p ( x , \xi ) = E_{0} + A ( y , \eta ) y \cdot \eta$ where $A \in C^{\infty} ( \R^{4} )$ satisfies $A ( 0 , 0 ) = - \diag ( \lambda_{1} , \lambda_{2} )$. By construction, we also have
\begin{equation} \label{g34}
\eta ( t , s ) = \xi ( t , s ) - \nabla \varphi_{-} ( x ( t , s ) ) = \nabla \varphi_{+}^{k} ( x ( t , s ) ) - \nabla \varphi_{-} ( x ( t , s ) ) .
\end{equation}
From \ref{h13}, the manifolds $\Lambda_{-}^{0}$ and $\Lambda_{+}^{1}$ have an intersection of order $m_{k}$ along the trajectory $\gamma_{k}$. From \eqref{g34}, it implies that $\eta ( t , s )$ vanishes at order $m_{k}$ on $s = 0$. It means that
\begin{equation} \label{g36}
\partial_{s}^{\alpha} \eta ( t , 0 ) = 0 ,
\end{equation}
for all $0 \leq \alpha \leq m_{k} - 1$ and $\partial_{s}^{m_{k}} \eta ( t , 0 ) \neq 0$. Combining with the particular form of $p$, it yields
\begin{align}
\partial_{t} \partial_{s}^{m_{k}} \eta ( t , 0 ) &= \partial_{s}^{m_{k}} \partial_{t} \eta ( t , 0 )  \nonumber \\
&= - \partial_{s}^{m_{k}} ( \partial_{y} p ) ( y ( t , s ) , \eta ( t , s ) ) \vert_{s=0}  \nonumber \\
&= - \partial_{s}^{m_{k}} \big( A ( y , \eta ) \eta + ( \partial_{y} A ) ( y , \eta ) y \cdot \eta \big) \vert_{s=0}   \nonumber \\
&= - A ( \gamma_{k} (t) ) \partial_{s}^{m_{k}} \eta ( t , 0 ) - ( \partial_{y} A ) ( \gamma_{k} (t) ) y ( t , 0 ) \cdot \partial_{s}^{m_{k}} \eta ( t , 0 ) , \label{g37}
\end{align}
since $\gamma_{k} (t) = ( y ( t , 0 ) , 0 )$ in the variables $( y , \eta )$. Using $\rho ( t , s ) \in \Lambda_{+}^{1} \subset p^{-1} ( E_{0} )$, we deduce
\begin{equation*}
\big( A ( y , \eta ) y \cdot \eta \big) ( t , s ) = 0 ,
\end{equation*}
for all $t , s \in \R$. Differentiating $m_{k}$ times with respect to $s$ and using \eqref{g36}, it yields
\begin{equation} \label{g38}
A ( \gamma_{k} (t) ) y ( t , 0 ) \cdot \partial_{s}^{m_{k}} \eta ( t , 0 ) = 0 .
\end{equation}
Let $( c_{1} ( t ) , c_{2} ( t ) ) : = \partial_{s}^{m_{k}} \eta ( t , 0 )$ denote the two coordinates of $\partial_{s}^{m_{k}} \eta ( t , 0 ) \in \R^{2}$. From \eqref{g39}, \eqref{g35} and $\xi ( t , 0 ) = \nabla \varphi_{-} ( x ( t , 0 ) )$, we get
\begin{equation*}
y ( t , 0 ) = \big( \vert g_{-}^{k} \vert e^{- \lambda_{1} t_{0}} e^{- \lambda_{1} t} , 0 \big) + \CO (  e^{- ( \lambda_{1} + \varepsilon ) t} ) ,
\end{equation*}
for some $\varepsilon > 0$. Then, \eqref{g38} gives
\begin{equation} \label{g40}
c_{1} ( t ) = \CO ( e^{- \varepsilon t} ) c_{2} ( t ) .
\end{equation}
Combining this equation with \eqref{g37}, we deduce
\begin{equation*}
\partial_{t} c_{2} ( t ) = \big( \lambda_{2} + r ( t ) \big) c_{2} ( t ) .
\end{equation*}
with $r (t) = \CO ( e^{- \varepsilon t} )$. It implies
\begin{equation} \label{g43}
c_{2} ( t ) = e^{\lambda_{2} t} e^{\int_{0}^{t} r ( u ) \, d u} c_{2} ( 0 ) = e^{\lambda_{2} t} \beta + \CO ( e^{( \lambda_{2} - \varepsilon ) t} ) .
\end{equation}
Since the contact between $\Lambda_{-}^{0}$ and $\Lambda_{+}^{1}$ is of order $m_{k}$, we have $\partial_{s}^{m_{k}} \eta ( t , 0 ) \neq 0$ and then $c_{2} ( t ) \neq 0$ from \eqref{g40}. It implies that $\beta = e^{\int_{0}^{+ \infty} r ( u ) \, d u} c_{2} ( 0 ) \neq 0$.

Define $G ( x_{1} , x_{2} ) : = \varphi^{k}_{+} ( x_{1} , x_{2} ) - \varphi_{-} ( x_{1} , x_{2} )$. From the assumption \ref{h13}, the manifolds $\Lambda_{-}^{0}$ and $\Lambda_{+}^{1}$ have an intersection of order $m_{k}$ along $\gamma_{k}$. It implies that $G$ vanishes of order $1 + m_{k}$ on the curve $\pi_{x} ( \gamma_{k} )$. In particular,
\begin{equation} \label{g41}
\partial_{x_{1}}^{\alpha} \partial_{x_{2}}^{\beta} G ( x_{k} (t) ) = 0 ,
\end{equation}
for all $\alpha , \beta \in \N$ with $\alpha + \beta \leq m_{k}$. Derivating with respect to $t$ yields
\begin{equation*}
d \partial_{x_{1}}^{\alpha} \partial_{x_{2}}^{\beta} G ( x_{k} (t) ) ( \partial_{t} x ( t , 0 ) ) = 0 .
\end{equation*}
From \eqref{g32}, it implies
\begin{equation*}
\partial_{x_{1}}^{1 + \alpha} \partial_{x_{2}}^{\beta} G ( x_{k} (t) ) = o \big( \partial_{x_{1}}^{\alpha} \partial_{x_{2}}^{1 + \beta} G ( x_{k} (t) ) \big) .
\end{equation*}
By iteration, we deduce
\begin{equation} \label{g42}
\partial_{x_{1}}^{1 + \alpha} \partial_{x_{2}}^{\beta} G ( x_{k} (t) ) = o \big( \partial_{x_{2}}^{1 + m_{k}} G ( x_{k} (t) ) \big) ,
\end{equation}
for all $\alpha , \beta \in \N$ with $\alpha + \beta = m_{k}$.

Using \eqref{g34}, we can write $\eta ( t , s ) = \nabla G ( x ( t , s ) )$. Differentiating $m_{k}$ times with respect to $s$ and using \eqref{g41}, we get
\begin{equation*}
\partial_{s}^{m_{k}} \eta ( t , 0 ) = \nabla d^{m_{k}} G ( x_{k} (t) ) \big( \partial_{s} x ( t , 0 ) , \ldots , \partial_{s} x ( t , 0 ) \big) .
\end{equation*}
Considering the second coordinate in this equality leads to
\begin{equation*}
c_{2} ( t ) = e^{- m_{k} \lambda_{2} t } ( 1 + o ( 1 ) ) \partial_{x_{2}}^{1 + m_{k}} G ( x_{k} (t) ) ,
\end{equation*}
from \eqref{g33} and \eqref{g42}. Combining with \eqref{g39} and \eqref{g43}, it gives
\begin{equation*}
\partial_{x_{2}}^{1 + m_{k}} G ( x_{k} (t) ) = \beta e^{( 1 + m_{k} ) \lambda_{2} t} ( 1 + o ( 1 ) ) = \alpha_{k}^{\infty} ( 1 + m_{k} ) ! x_{1}^{- ( 1 + m_{k} ) \lambda_{2} / \lambda_{1}} ( 1 + o ( 1 ) ) ,
\end{equation*}
with
\begin{equation*}
\alpha_{k}^{\infty} : = \frac{\beta}{( 1 + m_{k} ) !} \big( \vert g_{-}^{k} \vert e^{- \lambda_{1} t_{0}} \big)^{( 1 + m_{k} ) \lambda_{2} / \lambda_{1}} \neq 0 .
\end{equation*}
Since $G = \varphi^{k}_{+} - \varphi_{-}$, we have proved \eqref{g22}.
\end{proof}

\subsubsection{Proof of Theorem \ref{g24}} The proof of this result is similar to the one of Theorem \ref{d8}. The only thing which has to be changed is the computation of the quantum monodromy along the trapped set $K ( E_{0} )$. In other words, we must adapt Lemma \ref{d41} to the present setting. With a new version of this lemma, the resonance free zone and the polynomial resolvent estimate will follow as in Section \ref{s72}. Moreover, considering ``test functions'' microlocalized on the $K_{1}$ first homoclinic curves (corresponding to the most tangential and thus the most trapping trajectories), one can prove the existence of resonances near the pseudo-resonances as in Section \ref{s73}.

We use the notations of Section \ref{s72}. Let $u$ be a solution of
\begin{equation} \label{g50}
\left\{ \begin{aligned}
&( P - z ) u = 0 &&\text{microlocally near } K ( E_{0} ) ,  \\
&u = 0 &&\text{microlocally near } \Lambda_{-} \setminus ( \{ 0 \} \cup \CH ) ,
\end{aligned} \right.
\end{equation}
with $z \in \eqref{g25}$ and $\Vert u \Vert_{L^{2} ( \R^{2} )} \leq 1$. As in \eqref{d44}, one can show that
\begin{equation*}
u_{-}^{k} \in \CI ( \Lambda_{+}^{1} , h^{-N} ) \qquad \text{and} \qquad u_{+}^{k} \in \CI ( \Lambda_{+}^{0} , h^{-N} ) ,
\end{equation*}
for some $N \in \R$. Following \eqref{d48}, the symbols $a_{\pm}^{k} \in S ( h^{- N} )$ are defined by
\begin{equation*}
\left\{ \begin{aligned}
u_{-}^{k} (x) &= e^{i \frac{z - E_{0}}{h} t_{-}^{k}} \frac{\CM_{k}^{-}}{\CD_{k} ( t_{-}^{k} )} a_{-}^{k} ( x , h ) e^{i \varphi_{+}^{k} (x) / h} ,  \\
u_{+}^{k} (x) &= a_{+}^{k} ( x , h ) e^{i \varphi_{+} (x) / h} .
\end{aligned} \right.
\end{equation*}
If $m_{1} \geq 2$, Lemma \ref{d41} is replaced by

\begin{lemma}\sl \label{g46}
There exist $\zeta > 0$ and symbols $\CP_{k , \ell} \in S (1)$ independent of $u$ such that
\begin{equation} \label{g58}
a_{-}^{k} ( x , h ) = h^{S ( z , h ) / \lambda_{1} - \frac{m_{1}}{1 + m_{1}}} \sum_{\ell = 1}^{K_{1}} \CP_{k , \ell} ( x , h ) a_{-}^{\ell} ( x_{-}^{\ell} , h ) + S ( h^{-N + \zeta} ) ,
\end{equation}
for all $k \in \{ 1 , \ldots , K \}$ and $x$ near $\pi_{x} ( U_{-}^{k} \cap \Lambda_{+}^{1} )$. Moreover, $\CP_{k , \ell} ( x_{-}^{k} , h ) = \CQ_{k , \ell} ( z , h )$ where $\CQ$ is given by \eqref{g23} for $k \in \{ 1 , \ldots , K_{1} \}$.
\end{lemma}

Note that only the contributions of the $K_{1}$ ``most tangential'' trajectories appear in the leading term. The other contributions enter into the remainder term $S ( h^{-N + \zeta} )$.

\begin{proof}
As in the proof of Lemma \ref{d41}, we begin with the propagation through the hyperbolic fixed point $( 0 , 0 )$. The function $u$ satisfies \eqref{d62} and we can then apply Theorem \ref{a32} to compute $u_{+}$. More precisely, $u_{+}$ can be written as the sum over $\ell \in \{ 1 , \ldots , K \}$ of the contributions provided by \eqref{a31}. Thus, we have
\begin{equation*}
a_{+}^{k} ( x , h ) = h^{\frac{\lambda_{2} - \lambda_{1}}{2 \lambda_{1}} - i \frac{z - E_{0}}{\lambda_{1} h}} \sum_{\ell = 1}^{K} \int_{H_{\ell}} e^{i ( \varphi_{+}^{\ell} - \varphi_{-} ) ( y ) / h} e^{i \frac{z - E_{0}}{h} t_{-}^{\ell}} \frac{\CM_{\ell}^{-}}{\CD_{\ell} ( t_{-}^{\ell} )} d ( x , y , z , h ) a_{-}^{\ell} ( y , h ) \, d y ,
\end{equation*}
where $H_{\ell} = \{ x \in \R^{n} ; \ x_{1} = \varepsilon_{0} \}$ and $\varepsilon_{0} > 0$ denotes the first coordinate of $x_{-}^{\ell}$. Recall that, in the previous expression, we have made a linear changes of coordinates adapted to each trajectory $\gamma_{\ell}$ (i.e. the $y$ variables depend on $\ell$). We now compute these integrals by the method of the (degenerate) stationary phase. We note that the hyperplane $H_{\ell}$ is unidimensional and that the phase function $\varphi_{+}^{\ell} - \varphi_{-}$ has a zero of finite order $1 + m_{\ell}$ at $x_{2}^{k} ( y_{1} )$ from \eqref{g21}. Then, applying (7.7.30) and (7.7.31) of H\"{o}rmander \cite{Ho90_01} to compute these previous oscillatory integrals, we obtain
\begin{align*}
a_{+}^{k} ( x , h ) = h^{\frac{\lambda_{2} - \lambda_{1}}{2 \lambda_{1}} - i \frac{z - E_{0}}{\lambda_{1} h}} \sum_{\ell = 1}^{K} & h^{\frac{1}{1 + m_{\ell}}} \vert \alpha_{\ell} ( x_{-}^{\ell} ) \vert^{-\frac{1}{1 + m_{\ell}}} \frac{2 b_{\ell}}{1 + m_{\ell}} \Gamma \Big( \frac{1}{1 + m_{\ell}} \Big) \frac{\CM_{\ell}^{-}}{\CD_{\ell} ( t_{-}^{\ell} )}    \\
&\times e^{i \frac{z - E_{0}}{h} t_{-}^{\ell}} d ( x , x_{-}^{\ell} , z , h ) a_{-}^{\ell} ( x_{-}^{\ell} , h ) + S \big( h^{\frac{\lambda_{2} - \lambda_{1}}{2 \lambda_{1}} - i \frac{z}{\lambda_{1} h} + \frac{2}{1 + m_{\ell}} - N} \big),
\end{align*}
with the notations $\alpha_{\ell} ( x_{1} , x_{2} ) : = \alpha_{\ell} ( x_{1} )$ and
\begin{equation*}
b_{\ell} = \left\{ \begin{aligned}
&e^{i \frac{\pi}{2 + 2 m_{\ell}} \sgn ( \alpha_{\ell}^{\infty} )}  &&\text{ for odd } m_{\ell} ,  \\
&\cos \Big( \frac{\pi}{2 + 2 m_{\ell}} \Big) &&\text{ for even } m_{\ell} .
\end{aligned} \right.
\end{equation*}
Taking only the leading terms, it gives
\begin{equation} \label{g47}
a_{+}^{k} ( x , h ) = h^{S ( z , h ) / \lambda_{1} - \frac{m_{1}}{1 + m_{1}}} \sum_{\ell = 1}^{K_{1}} \widetilde{\CR}_{k , \ell} ( x , h ) a_{-}^{\ell} ( x_{-}^{\ell} , h ) + S ( h^{- N + \zeta} ) ,
\end{equation}
with $\zeta > 0$ and $\widetilde{\CR}_{k , \ell} \in S ( 1 )$ defined by
\begin{equation*}
\widetilde{\CR}_{k , \ell} ( x , h ) = \vert \alpha_{\ell} ( x_{-}^{\ell} ) \vert^{-\frac{1}{1 + m_{1}}} \frac{2 b_{\ell}}{1 + m_{1}} \Gamma \Big( \frac{1}{1 + m_{1}} \Big) \frac{\CM_{\ell}^{-}}{\CD_{\ell} ( t_{-}^{\ell} )} e^{i \frac{z - E_{0}}{h} t_{-}^{\ell}} d_{0} ( x , x_{-}^{\ell} , z ) .
\end{equation*}

We now explicit $\widetilde{\CR}_{k , \ell}$ in terms of geometric quantities. For that, we use the usual trick: we obtain a formula of $\widetilde{\CR}_{k , \ell}$ in terms of a free parameter (concretely, the distance to $0$ of the hyperplane of integration), and we compute its asymptotic behavior when this parameter diverges. Such a method has been used in \cite[(6.26)]{BoFuRaZe07_01} for instance. For $0 < \varepsilon \leq \varepsilon_{0}$, let $t_{\varepsilon} > 0$ be such that the first coordinate of $x_{\ell} ( t_{\varepsilon} )$ is $\varepsilon$. In particular, $t_{\varepsilon_{0}} = t_{-}^{\ell}$, $x_{\ell} ( t_{\varepsilon_{0}} ) = x_{-}^{\ell}$ and $t_{\varepsilon} \to + \infty$ as $\varepsilon \to 0$. We remark that the symbol $\widetilde{\CR}_{k , \ell}$ is independent of the initial condition $a_{-}^{\ell} ( x_{-}^{\ell} , h )$. Thus, to compute $\widetilde{\CR}_{k , \ell}$, we can always assume that $N = 0$, $a_{-}^{\ell} ( x_{-}^{\ell} , h ) = 1$ and $a_{-}^{\ell^{\prime}} ( x , h ) = 0$ for $\ell \neq \ell^{\prime}$. The idea is to compute $a_{-}$ near $x_{\ell} ( t_{\varepsilon} )$ using the usual propagation of Lagrangian distributions and then to apply Theorem \ref{a32} with initial condition on the hyperplane $\{ y_{1} = \varepsilon \}$. Since $( P - z ) u_{-} = 0$ microlocally near $\pi_{x} ( \gamma_{\ell} \cap \Lambda_{-}^{0} )$, the usual propagation of WKB solutions gives
\begin{equation} \label{g49}
a_{-} ( x_{\ell} ( t_{\varepsilon} ) , h ) = a_{-} ( x_{\ell} ( t_{\varepsilon_{0}} ) , h ) \frac{\CD_{\ell} ( t_{\varepsilon_{0}} )}{\CD_{\ell} ( t_{\varepsilon} )} e^{i \frac{z - E_{0}}{h} ( t_{\varepsilon} - t_{\varepsilon_{0}} )} + S ( h^{- N + 1} ) .
\end{equation}
We now compute $a_{+} ( x , h )$ using Theorem \ref{a32} with initial condition on $\{ x_{1} = \varepsilon \}$. In other words, applying \eqref{g47} with $\varepsilon_{0}$ replaced by $\varepsilon$, we deduce
\begin{align}
a_{+}^{k} ( x , h ) = h^{S ( z , h ) / \lambda_{1} - \frac{m_{1}}{1 + m_{1}}} \vert \alpha_{\ell} ( x_{\ell} ( t_{\varepsilon} ) ) & \vert^{-\frac{1}{1 + m_{1}}} \frac{2 b_{\ell}}{1 + m_{1}} \Gamma \Big( \frac{1}{1 + m_{1}} \Big) \frac{\CM_{\ell}^{-}}{\CD_{\ell} ( t_{\varepsilon_{0}} )}   \nonumber \\
&\times e^{i \frac{z - E_{0}}{h} t_{\varepsilon_{0}}} d_{0} ( x , x_{\ell} ( t_{\varepsilon} ) , z ) a_{-}^{\ell} ( x_{\ell} ( t_{\varepsilon} ) , h ) + S ( h^{\zeta} ) . \label{g48}
\end{align}
Combining \eqref{g49} and \eqref{g48} with the uniqueness of the solution of the microlocal Cauchy problem \eqref{g50}, we get
\begin{equation} \label{g54}
\widetilde{\CR}_{k , \ell} ( x , h ) = \vert \alpha_{\ell} ( x_{\ell} ( t_{\varepsilon} ) ) \vert^{-\frac{1}{1 + m_{1}}} \frac{2 b_{\ell}}{1 + m_{1}} \Gamma \Big( \frac{1}{1 + m_{1}} \Big) \frac{\CM_{\ell}^{-}}{\CD_{\ell} ( t_{\varepsilon} )} e^{i \frac{z - E_{0}}{h} t_{\varepsilon}} d_{0} ( x , x_{\ell} ( t_{\varepsilon} ) , z ) ,
\end{equation}
for all $0 < \varepsilon \leq \varepsilon_{0}$. We now compute the asymptotic as $\varepsilon \to 0$ of this quantity which is independent of $\varepsilon$. As in Lemma \ref{b17}, the symbol $d_{0}$ satisfies the following estimate. Its proof is postponed at the end of the section.

\begin{lemma}\sl \label{g57}
We have
\begin{align}
d_{0} ( x_{+}^{k} , x_{\ell} ( t_{\varepsilon} ) , z ) = \frac{\sqrt{\lambda_{1} \lambda_{2}}}{2 \pi} e^{- i \frac{\pi}{2}} \Gamma \big( S ( z , h ) / \lambda_{1} \big) & \frac{{\CM_{k}^{+}}}{\CD_{k} ( t_{+}^{k} )} \big\vert g_{-}^{\ell} \big\vert \big( i \lambda_{1} g_{+}^{k} \cdot g_{-}^{\ell} \big)^{- S ( z , h ) / \lambda_{1}}    \nonumber \\
&\quad \times e^{i \frac{z - E_{0}}{h} t_{+}^{k}} e^{t_{\varepsilon} S ( z , h )} e^{- \lambda_{1} t_{\varepsilon}} ( 1 + o_{\varepsilon \to 0} ( 1 ) ) .  \label{g52}
\end{align}
\end{lemma}

On the other hand, \eqref{g44} gives
\begin{equation} \label{g51}
\CD_{\ell} ( t_{\varepsilon} ) = \CM_{\ell}^{-} e^{- t_{\varepsilon} \sum_{j} \lambda_{j} / 2} ( 1 + o_{\varepsilon \to 0} ( 1 ) ) .
\end{equation}
Eventually, \eqref{g22} implies
\begin{equation} \label{g53}
\vert \alpha_{\ell} ( x_{\ell} ( t_{\varepsilon} ) ) \vert^{-\frac{1}{1 + m_{1}}} = \vert \alpha_{k}^{\infty} \vert^{-\frac{1}{1 + m_{1}}} \vert g_{-}^{\ell} \vert^{\lambda_{2} / \lambda_{1}} e^{- \lambda_{2} t_{\varepsilon}} ( 1 + o_{\varepsilon \to 0} ( 1 ) ) .
\end{equation}
since $x_{\ell} ( t_{\varepsilon} ) = g_{-}^{\ell} e^{- \lambda_{1} t_{\varepsilon}} ( 1 + o_{\varepsilon \to 0} ( 1 ) )$ from \eqref{d2}. We now insert the asymptotics \eqref{g52}, \eqref{g51} and \eqref{g53} in the formula \eqref{g54}. Using that $\widetilde{\CR}_{k , \ell}$ is equal to its limit as $\varepsilon \to 0$ since it is independent of $\varepsilon$, it yields
\begin{align}
\widetilde{\CR}_{k , \ell} ( x_{+}^{k} , h ) = \vert \alpha_{k}^{\infty} \vert^{-\frac{1}{1 + m_{1}}} b_{\ell} \frac{\sqrt{\lambda_{1} \lambda_{2}}}{\pi ( 1 + m_{1} )} \Gamma \Big( & \frac{1}{1 + m_{1}} \Big) \Gamma \big( S ( z , h ) / \lambda_{1} \big) \frac{{\CM_{k}^{+}}}{\CD_{k} ( t_{+}^{k} )} e^{- i \frac{\pi}{2}}    \nonumber \\
& \times \vert g_{-}^{\ell} \vert^{\frac{\lambda_{1} + \lambda_{2}}{\lambda_{1}}} \big( i \lambda_{1} g_{+}^{k} \cdot g_{-}^{\ell} \big)^{- S ( z , h ) / \lambda_{1}} e^{i \frac{z - E_{0}}{h} t_{+}^{k}} .  \label{g55}
\end{align}

We now express $a_{-}$ in terms of $a_{+}$. From \eqref{g50}, the function $u$ verifies the microlocal Cauchy problem \eqref{d65}. Then, exactly as in \eqref{d52}, we have
\begin{equation} \label{g56}
e^{i \frac{z - E_{0}}{h} t_{-}^{k}} \frac{\CM_{k}^{-}}{\CD_{k} ( t_{-}^{k} )} a_{-}^{k} ( x , h ) = \CS_{k} ( x , h ) a_{+}^{k} ( \widetilde{x} (x)  , h ) + S ( h^{-N + 1} ) ,
\end{equation}
where $\widetilde{x} (x) = \pi_{x} ( \exp ( ( t_{-}^{k} - t_{+}^{k} ) H_{p} ) ( x , \nabla \varphi_{+}^{1} (x) ) )$ and $\CS_{k} \in S (1)$ is a symbol satisfying
\begin{equation} \label{g59}
\CS_{k} ( x_{-}^{k} , h ) = e^{i A_{k} / h} e^{- i \frac{\pi}{2} \nu_{k}} e^{i \frac{z - E_{0}}{h} ( t_{-}^{k} - t_{+}^{k} )} \frac{\CD_{k} ( t_{+}^{k} )}{\CD_{k} ( t_{-}^{k} )} .
\end{equation}
Combining \eqref{g47} and \eqref{g56}, we obtain \eqref{g58}. Eventually, \eqref{g55} and \eqref{g59} give
\begin{align}
\CP_{k , \ell} ( x_{-}^{k} , h ) &= e^{- i \frac{z - E_{0}}{h} t_{-}^{k}} \frac{\CD_{k} ( t_{-}^{k} )}{\CM_{k}^{-}} \CS_{k} ( x_{-}^{k} , h ) \widetilde{\CR}_{k , \ell} ( x_{+}^{k} , h )    \nonumber \\
&= e^{i A_{k} / h} \vert \alpha_{k}^{\infty} \vert^{-\frac{1}{1 + m_{1}}} b_{\ell} \frac{\sqrt{\lambda_{1} \lambda_{2}}}{\pi ( 1 + m_{1} )} \Gamma \Big( \frac{1}{1 + m_{1}} \Big) \Gamma \big( S ( z , h ) / \lambda_{1} \big) \frac{\CM_{k}^{+}}{\CM_{k}^{-}}    \nonumber \\
&\qquad \qquad \qquad \qquad \qquad \quad \times e^{- \frac{\pi}{2} ( \nu_{k} + 1 ) i} \big\vert g_{-}^{\ell} \big\vert^{\frac{\lambda_{1} + \lambda_{2}}{\lambda_{1}}} \big( i \lambda_{1} g_{+}^{k} \cdot g_{-}^{\ell} \big)^{- S ( z , h ) / \lambda_{1}}   \nonumber \\
&= \CQ_{k , \ell} ( z , h ) ,
\end{align}
and Lemma \ref{g46} follows.
\end{proof}

\begin{proof}[Proof of Lemma \ref{g57}]
This result is similar to Lemma \ref{b17} and we used here some estimates of its proof (see Section \ref{s20}). First, coming back to the definitions of $t_{\varepsilon}$ and $g_{-}^{\ell}$, we obtain
\begin{equation*}
\varepsilon = \vert g_{-}^{\ell} \vert e^{- \lambda_{1} t_{\varepsilon}} ( 1 + o_{\varepsilon \to 0} ( 1 ) ) .
\end{equation*}
Moreover, with the notation of \eqref{d2}, we have
\begin{equation*}
g_{-} ( \gamma_{\ell} ( t_{\varepsilon} ) ) = g_{-}^{\ell} e^{- \lambda_{1} t_{\varepsilon}} \qquad \text{and} \qquad g_{+} ( \rho_{+}^{k} ) = g_{+}^{k} e^{\lambda_{1} t_{+}^{k}} .
\end{equation*}
Otherwise, as in \eqref{b23} and \eqref{b24}, one can show
\begin{equation*}
\big\vert \det \nabla^{2}_{y^{\prime} , y^{\prime} } \varphi_{-} ( x_{\ell} ( t_{\varepsilon} ) ) \big\vert^{\frac{1}{2}} = \sqrt{\frac{\lambda_{2}}{2}} + o_{\varepsilon \to 0} ( 1 ) \qquad \text{and} \qquad \big\vert \partial_{\xi_{1}} p ( \gamma_{\ell} ( t_{\varepsilon} ) ) \big\vert^{\frac{1}{2}} = \sqrt{\lambda_{1} \varepsilon} ( 1 + o_{\varepsilon \to 0} ( 1 ) ) .
\end{equation*}
The Liouville formula $\partial_{t} \ln \CD_{k} ( t ) = \Delta \varphi_{+} ( x_{k} ( t ) )$ for $t \ll - 1$ and \eqref{d7} give
\begin{equation*}
e^{\int_{0}^{- \infty} ( \Delta \varphi_{+} ( x (s) ) - \sum_{j} \lambda_{j} / 2 ) \, d s} =  \frac{{\CM_{k}^{+}}}{\CD_{k} ( t_{+}^{k} )} e^{t_{+}^{k} \sum_{j} \lambda_{j} / 2} ,
\end{equation*}
with the notations of \eqref{b19} in the left hand side. Eventually, we will prove that
\begin{equation} \label{g70}
\lim_{t\to + \infty} \frac{e^{(\sum_{j} \lambda_{j} /2 -\lambda_{1})t}}{{\sqrt{\Big\vert \det \frac{\partial y ( t , y^{\prime} , \eta^{\prime} )}{\partial ( t , y^{\prime} )} \vert_{\eta^{\prime} = \partial_{y^{\prime}} \varphi_{-} ( y )} \Big\vert}}} = \sqrt{\frac{2}{\lambda_{1} \varepsilon}} ( 1 + o_{\varepsilon \to 0} ( 1 ) ) .
\end{equation}
The previous estimates together with \eqref{b19} imply the lemma.

It remains to show \eqref{g70}. Computing the derivative with respect to $t$ as in \eqref{b26} yields
\begin{equation} \label{g68}
\frac{\partial y ( t , y^{\prime} , \eta^{\prime} )}{\partial t} \vert_{\eta^{\prime} = \partial_{y^{\prime}} \varphi_{-} (y)} = - \lambda_{1} \varepsilon e^{- \lambda_{1} t}
\bigg( \begin{gathered}
1  \\
0
\end{gathered} \bigg)
+ \CO ( e^{- ( \lambda_{1} + \nu ) t } ) ,
\end{equation}
for some $\nu > 0$. We now compute the derivative of $y ( t , y^{\prime} , \eta^{\prime} )$ with respect to $y^{\prime}$ at $\eta^{\prime} = \partial_{y^{\prime}} \varphi_{-} (y)$ as in \eqref{b35}. This estimate can not be directly applied here since we have not assumed that $\lambda_{1} = \lambda_{2}$. Note that $y^{\prime} = y_{2}$ since $n = 2$ in the present setting. As in \eqref{b29}, we have to follow the evolution along the Hamiltonian flow of a tangent vector $( \delta_{y} , \delta_{\eta} )$ of $\Lambda_{\eta^{\prime}}$. In particular,
\begin{equation} \label{g60}
\delta_{y} (0) = ( 0 , \delta_{y^{\prime}} (0) ) \qquad \text{and} \qquad \delta_{\eta} (0) = \big( \partial_{y^{\prime}} f_{-} ( \varepsilon , y^{\prime} , \eta^{\prime} ) \cdot \delta_{y^{\prime}} (0) , 0 \big) .
\end{equation}
Let symplectic local coordinates $( k , \kappa ) \in T^{*} \R^{2}$ centered at $( 0 , 0 )$ such that $\Lambda_{-}^{0}$ (resp. $\Lambda_{+}^{0}$) is given by $k = 0$ (resp. $\kappa =0$) with
\begin{equation} \label{g64}
k_{j} = \frac{1}{\sqrt{\lambda_{j}}} \Big( \eta_{j} + \frac{\lambda_{j}}{2} y_{j} \Big) + \CO \big( ( y , \eta )^{2} \big)   \quad \text{and} \quad \kappa_{j} = \frac{1}{\sqrt{\lambda_{j}}} \Big( \eta_{j} - \frac{\lambda_{j}}{2} y_{j} \Big) + \CO \big( ( y , \eta )^{2} \big) .
\end{equation}
Then, $p ( y , \eta ) = E_{0} + A ( k , \kappa ) k \cdot \kappa$ with $A (0,0) = \diag ( \lambda_{1} , \lambda_{2} )$. In these coordinates, the initial condition \eqref{g60} becomes
\begin{align}
\delta_{k} (0) &= \frac{1}{2} \Big( 2 \lambda_{1}^{- 1 / 2} \partial_{y^{\prime}} f_{-} ( \varepsilon , y^{\prime} , \eta^{\prime} ) \cdot \delta_{y^{\prime}} (0) , \lambda_{2}^{1 / 2} \delta_{y^{\prime}} (0) \Big) + \CO ( \varepsilon ) \vert \delta_{y^{\prime}} (0) \vert ,   \label{g66} \\
\delta_{\kappa} (0) &= \frac{1}{2} \Big( 2 \lambda_{1}^{- 1 / 2} \partial_{y^{\prime}} f_{-} ( \varepsilon , y^{\prime} , \eta^{\prime} ) \cdot \delta_{y^{\prime}} (0) , - \lambda_{2}^{1 / 2} \delta_{y^{\prime}} (0) \Big) + \CO ( \varepsilon ) \vert \delta_{y^{\prime}} (0) \vert .
\end{align}
Moreover, as in \eqref{b27}, the evolution of the tangent vector $( \delta_{k} , \delta_{\kappa} )$ verifies
\begin{equation}
\frac{d}{d t}
\left( \begin{array}{c}
\delta_{k} \\
\delta_{\kappa}
\end{array} \right)
=
\left( \begin{array}{cc}
\diag ( \lambda_{1} , \lambda_{2} ) + \CO \big( \varepsilon e^{-\lambda_{1} t} \big) & 0 \\
\CO \big( \varepsilon e^{- \lambda_{1} t} \big) & - \diag ( \lambda_{1} , \lambda_{2} ) + \CO \big( \varepsilon e^{- \lambda_{1} t} \big)
\end{array} \right) \left( \begin{array}{c}
\delta_{k} \\
\delta_{\kappa}
\end{array} \right) .
\end{equation}
From the first equation
\begin{equation} \label{g61}
\partial_{t} \delta_{k} = \big( \diag ( \lambda_{1} , \lambda_{2} ) + \CO \big( \varepsilon e^{-\lambda_{1} t} \big) \big) \delta_{k} ,
\end{equation}
we deduce
\begin{equation} \label{g65}
\vert \delta_{k} ( t ) \vert = \CO ( e^{\lambda_{2} t} ) \vert \delta_{k} (0) \vert .
\end{equation}
Thus, the second line of \eqref{g61} can be written
\begin{equation*}
\partial_{t} \delta_{k_{2}} = \lambda_{2} \delta_{k_{2}} + \CO \big( \varepsilon e^{( \lambda_{2} - \lambda_{1} ) t} \big) \vert \delta_{k} ( 0 ) \vert ,
\end{equation*}
which implies $\partial_{t} ( e^{- \lambda_{2} t} \delta_{k_{2}} ) = \CO ( \varepsilon e^{- \lambda_{1} t} ) \vert \delta_{k} ( 0 ) \vert$. As consequence,
\begin{equation} \label{g62}
\delta_{k_{2}} ( t ) = e^{\lambda_{2} t} \delta_{k_{2}} ( 0 ) + \CO \big( \varepsilon  e^{\lambda_{2} t} \big) \vert \delta_{k} ( 0 ) \vert .
\end{equation}
On the other hand, as in \eqref{b32}, \cite[Lemma 5.6]{HeSj85_01} gives $\delta_{\kappa} ( t ) = \CO_{\varepsilon} ( e^{- \lambda_{1} t} ) \delta_{k} (t)$ and then
\begin{equation} \label{g63}
\delta_{\kappa} ( t ) = \CO_{\varepsilon} \big( e^{( \lambda_{2} - \lambda_{1} ) t} \big) \vert \delta_{k} ( 0 ) \vert ,
\end{equation}
from \eqref{g65}. Inverting the change of variables \eqref{g64}, the relations \eqref{g65}--\eqref{g63} imply
\begin{equation} \label{g67}
\delta_{y} (t) = \frac{e^{\lambda_{2} t}}{\sqrt{\lambda_{2}}} \big( \CO ( 1 ) \vert \delta_{k} (0) \vert , \delta_{k_{2}} (0) \big) + \CO_{\varepsilon} \big( e^{( \lambda_{2} - \lambda_{1} ) t} \big) \vert \delta_{k} ( 0 ) \vert + \CO \big( \varepsilon  e^{\lambda_{2} t} \big) \vert \delta_{k} ( 0 ) \vert .
\end{equation}
As in \eqref{b37}, we have $\partial_{y^{\prime}} f_{-} ( \varepsilon , y^{\prime} , \partial_{y^{\prime}} \varphi_{-} (y) ) = \CO (1)$ and \eqref{g66} becomes
\begin{equation*}
\delta_{k} (0) = \frac{\sqrt{\lambda_{2}}}{2} \big( \CO ( 1 ) \vert \delta_{y^{\prime}} (0) \vert, \delta_{y^{\prime}} (0) \big) + \CO ( \varepsilon ) \vert \delta_{y^{\prime}} (0) \vert .
\end{equation*}
Thus, \eqref{g67} can be written
\begin{equation*}
\delta_{y} (t) = \frac{e^{\lambda_{2} t}}{2} \big( \CO ( 1 ) \vert \delta_{y^{\prime}} (0) \vert , \delta_{y^{\prime}} (0) \big) + \CO_{\varepsilon} \big( e^{( \lambda_{2} - \lambda_{1} ) t} \big) \vert \delta_{y^{\prime}} (0) \vert + \CO \big( \varepsilon  e^{\lambda_{2} t} \big) \vert \delta_{y^{\prime}} (0) \vert .
\end{equation*}
In other words, we have proved
\begin{equation}
\frac{\partial y ( t , y^{\prime} , \eta^{\prime} )}{\partial y^{\prime}} \vert_{\eta^{\prime} = \partial_{y^{\prime}} \varphi_{-} (y)} = \frac{1}{2} e^{\lambda_{2} t} 
\left( \begin{gathered}
\CO ( 1 ) \\
1
\end{gathered} \right)
+ \CO_{\varepsilon} \big( e^{( \lambda_{2} - \lambda_{1} ) t} \big) + \CO \big( \varepsilon  e^{\lambda_{2} t} \big) .  \label{g69}
\end{equation}
Combining \eqref{g68} and \eqref{g69}, we get
\begin{equation*}
\Big\vert \det \frac{\partial y ( t , y^{\prime} , \eta^{\prime} )}{\partial ( t , y^{\prime} )}\vert_{\eta^{\prime} = \partial_{y^{\prime}} \varphi_{-} ( y )} \Big\vert = \frac{\lambda_{1} \varepsilon}{2} e^{( \lambda_{2} - \lambda_{1} ) t} + \CO_{\varepsilon} \big( e^{( \lambda_{2} - \lambda_{1} - \nu ) t} \big) + \CO \big( \varepsilon^{2} e^{( \lambda_{2} - \lambda_{1} ) t} \big) ,
\end{equation*}
for some new $\nu > 0$. Eventually, this gives \eqref{g70}.
\end{proof}

\section{Proof of the asymptotic of the resonances for a nappe of homoclinic curves} \label{s27}

The proof of Theorem \ref{i55} is a combination of that of resonance free domains for strong trapping (see Section \ref{s5}) and that of the asymptotic of the resonances generated by transversal homoclinic trajectories (see Section \ref{s12}). This is done in Section \ref{s28} and Section \ref{s29}. Lastly, the proof of the other results of Section \ref{s26} can be found in Section \ref{s30}.

The first step is to give the asymptotic of the pseudo-resonances and to obtain some useful estimates on the classical quantization operator $\CT$. It is crucial to recall that, thanks to \ref{h7}, the operator $\CT$ can be written $\CT ( z , h ) = \widetilde{\CT} ( \rho , \sigma )$ (see \eqref{i81}). Moreover, $( \rho , \sigma) \mapsto \widetilde{\CT} ( \rho , \sigma)$ is an analytic function on $( \S^{1} )^{K} \times \C$ (for $\im \sigma > - n \lambda / 2$) with values in the set of compact operators.

Proposition \ref{i53} is similar to Proposition \ref{d9}. The difference is that $\CT$ may be a compact operator on an infinite dimensional space whereas $\CQ$ is a matrix. Nevertheless, the proof of Proposition \ref{d9} can be followed here since the number of eigenvalues of $\CT$ playing a role in the set \eqref{i54} is uniformly bounded (see \eqref{i52}) and since the perturbation theory for finite systems of eigenvalues works as in the case of finite dimension (see Chapter IV.5 of Kato \cite{Ka76_01}). We omit the details. Moreover, one can obtain the following result adapting the proof of Lemma \ref{d13}.

\begin{lemma}\sl \label{i92}
Let $\beta , C > 0$. There exists $M > 0$ such that, for all $z$ in \eqref{i54},
\begin{equation*}
\dist \big( z , \res_{0} (P) \big) > \beta \frac{h}{\vert \ln h \vert} \quad \Longrightarrow \quad \big\Vert \big( 1 - h^{- i \frac{z - E_{0}}{\lambda h}} \CT ( z , h ) \big)^{- 1} \big\Vert \leq M .
\end{equation*}
\end{lemma}

\Subsection{Resonance free zone and resolvent estimate} \label{s28}

As in Section \ref{s72}, we show that $P$ has no resonance away from the pseudo-resonances and that its resolvent satisfies a polynomial estimate. More precisely, we have

\begin{proposition}\sl \label{i93}
Let $C , \delta > 0$. For $h$ small enough, $P$ has no resonance in the domain
\begin{equation} \label{i94}
E_{0} + [ - C h , C  h ] + i \Big[ - C \frac{h}{\vert \ln h \vert} , \frac{h}{\vert \ln h \vert} \Big] \setminus \Big( \res_{0} (P) + B \Big( 0 , \delta \frac{h}{\vert \ln h \vert} \Big) \Big) .
\end{equation}
Moreover, there exists $M > 0$ such that
\begin{equation*}
\big\Vert( P_{\theta} -z )^{-1} \big\Vert \lesssim h^{- M} ,
\end{equation*}
uniformly for $h$ small enough and $z \in \eqref{i94}$.
\end{proposition}

This result is proved in the rest of this part following the general contradiction argument of Section \ref{s3}, the strategy developed in Section \ref{s72} to deal with homoclinic trajectories and the constructions of Section \ref{s5} made for the present strong trapping case.

Following Section \ref{s3}, we consider a solution $u$ of \eqref{a4} with $\Omega_{h} = \eqref{i94}$ and show that it vanishes microlocally near each point of $K ( E_{0} ) = \{ ( 0 , 0 ) \} \cup \CH$. This will prove Proposition \ref{i93}. We assume that $P_{\theta} = P$ microlocally in a neighborhood of $K ( E_{0} )$ and use the notations, constructions and results of Section \ref{s5} (the hypotheses of Section \ref{s77} hold true here). Thus, for $\varepsilon > 0$ be small enough, we set $S_{\pm}^{\varepsilon} = \{ ( x, \xi ) \in \Lambda_{\pm}^{0} ; \ \vert x \vert = \varepsilon \}$, $\CH_{\pm}^{\varepsilon} = S_{\pm}^{\varepsilon} \cap \CH$ and $\CH_{{\rm tang} , \pm}^{\varepsilon} = S_{\pm}^{\varepsilon} \cap \CH_{\rm tang}$.

For $\beta \in \S^{n - 1}$, let $F ( \varepsilon , \beta ) \in \S^{n - 1}$ denote the normalized asymptotic direction of the Hamiltonian trajectory of $\Lambda_{-}^{0}$ passing through $\varepsilon \beta$. This function is studied in Section \ref{a71}. To be precise, $F$ concerns initially the Hamiltonian trajectories of $\Lambda_{+}^{0}$ but can also be used on $\Lambda_{-}^{0}$ (see the paragraph after \eqref{a70}). For $\varepsilon$ small enough, $F ( \varepsilon , \cdot )$ is a diffeomorphism of $\S^{n - 1}$ and the image of $\CH_{{\rm tang} , -}^{\varepsilon}$ is precisely $\CH_{\rm tang}^{- \infty}$. From \ref{h7} and \eqref{a68}, there exists $K$ disjoint open subsets $V_{1} , \ldots , V_{K}$ of $\S^{n - 1}$
with
\begin{equation*}
\pi_{x} \big( \CH_{\rm tang}^{- \infty} \cup \CH_{{\rm tang} , -}^{\varepsilon} \big) \subset \CV : = \bigcup_{1 \leq k \leq K} V_{k} ,
\end{equation*}
for $\varepsilon$ small enough and such that the action $A$ and the Maslov index $\nu$ are constant on each $V_{k} \cap \pi_{x} ( \CH_{\rm tang}^{- \infty} \cup \CH_{{\rm tang} , -}^{\varepsilon} )$. We also consider $V_{\pm}^{\varepsilon} \subset \varepsilon \S^{n - 1}$ small neighborhoods of $\pi_{x} ( \CH_{{\rm tang} , \pm}^{\varepsilon} )$ such that $V_{-}^{\varepsilon} \Subset \varepsilon \CV$

We construct the functions $u_{\pm}^{\rm tang} , u_{\pm}^{\rm trans} \in L^{2} ( \R^{n} )$ as in the beginning of Section \ref{s51}. In particular, the microsupport of $u_{-}^{\rm tang}$ intersected with $\varepsilon \S^{n - 1}$ is a subset of $V_{-}^{\varepsilon}$. It has already been proved in \eqref{a92}--\eqref{b4} that
\begin{equation} \label{i95}
u_{-}^{\bullet} \in \CI ( \Lambda_{+}^{1} , h^{- N_{-}^{\bullet}} ) \qquad \text{and} \qquad u_{+}^{\rm tang} \in \CI ( \Lambda_{+}^{0} , h^{- N_{-}^{\rm tang}} ) ,
\end{equation}
for some $N_{-}^{\bullet} \in \R$. Moreover, Lemma \ref{a94} implies that we always have
\begin{equation} \label{i96}
N_{-}^{\rm trans} \leq N_{-}^{\rm tang} .
\end{equation}
In other words, the transversal part $u_{-}^{\rm trans}$ is controlled by the tangential part $u_{-}^{\rm tang}$.

We now deal with the tangential part following Section \ref{s52}. We set
\begin{equation*}
u_{-}^{\rm tang} ( y , h ) = a_{-} ( y , h ) e^{- i A^{y} / h} e^{i \varphi^{1}_{+} (y) / h} \qquad \text{and} \qquad
u_{+}^{\rm tang} ( x , h ) = a_{+} ( x , h ) e^{i \varphi_{+} (x) / h} ,
\end{equation*}
with $a_{\pm} \in S ( h^{- N_{-}^{\rm tang}} )$. The phase function $\varphi_{+}^{1}$ parametrizes $\Lambda_{+}^{1}$ with the same convention mutatis mutandis as in \eqref{d45}. For $y \in \supp ( a_{-} ) \subset V_{-}^{\varepsilon} \subset \varepsilon \CV$, there is a unique $k \in \{ 1 , \ldots , K \}$ such that $y \in V_{k}$ and the action term $A^{y}$ is chosen equal to $A_{k}$. Since the $V_{k}$'s are disjoint, this choice does not affect the regularity of $a_{-}$. Note that these notations and the renormalization factor $e^{- i A^{y} / h}$ are slightly different from those of Section \ref{s52} and \eqref{d48}.

We now obtain a closed equation on $a_{-}$ making a turn along the trapped set. The following result corresponds to Lemma \ref{d41} (see also Lemma \ref{g89} and Lemma \ref{g11}) in the strong trapping case.

\begin{lemma}\sl \label{i97}
There exist $R , \zeta > 0$ independent of $u$ such that, for $\varepsilon$ small enough and $V_{-}^{\varepsilon}$ sufficiently close to $\pi_{x} ( \CH_{{\rm tang} , -}^{\varepsilon} )$, we have
\begin{equation} \label{i98}
\Vert a_{-} ( \varepsilon \cdot , h ) \Vert_{L^{\infty} ( \S^{n - 1} )} \leq R \Vert a_{-} ( \varepsilon \cdot , h ) \Vert_{L^{2} ( \CH_{\rm tang}^{- \infty} )} + \CO_{\varepsilon} ( h^{- N_{-}^{\rm tang} + \zeta} ) ,
\end{equation}
and
\begin{align}
a_{-} ( \varepsilon \cdot, h ) &= h^{- i \frac{z - E_{0}}{h \lambda}} \CT ( z , h ) a_{-} ( \varepsilon \cdot , h ) + o_{\varepsilon  \to 0} ( 1 ) \Vert a_{-} ( \varepsilon \cdot , h ) \Vert_{L^{2} ( \CH_{\rm tang}^{- \infty} )} + \CO_{\varepsilon} ( h^{- N_{-}^{\rm tang} + \zeta} ) ,   \label{i99}
\end{align}
as functions in $L^{2} ( \CH_{\rm tang}^{- \infty} )$.
\end{lemma}

\begin{proof}
We first express $a_{+}$ in terms of $a_{-}$. For that, we follow the beginning of Section \ref{s52}. Equation \eqref{a91} states that
\begin{equation*}
u_{+}^{\rm tang} = \CJ_{{\rm tang} \leftarrow {\rm tang}} u_{-}^{\rm tang} + \CI ( \Lambda_{+}^{0} , h^{- N_{-}^{\rm tang} + \frac{1}{2}} ) .
\end{equation*}
To compute $\CJ_{{\rm tang} \leftarrow {\rm tang}} u_{-}^{\rm tang}$, we decompose $u_{-}^{\rm tang}$ into small parts as in Section \ref{s52}. Thus, for $\nu > 0$ be small enough, we consider the partition of the unity constructed in \eqref{b44}. It induces the decompositions \eqref{b12} and  \eqref{b43}--\eqref{b42} of $a_{-}$ and $a_{+}$. That is
\begin{equation} \label{j31}
a_{-} ( y , h ) = \sum_{1 \leq j \leq J} a_{-}^{j} ( y , h ) + S ( h^{\infty} ) ,
\end{equation}
where $a_{-}^{j} ( y , h ) = \varphi_{j} ( y - y_{j} ) a_{-} ( y , h )$ for $y \in \varepsilon \S^{n - 1}$ and the $\varphi_{j} \in C^{\infty}_{0} ( \R^{n} ; [ 0 , 1 ] )$ satisfy
\begin{equation*}
\sum_{1 \leq j \leq J} \varphi_{j} ( y - y_{j} ) = 1 \text{ locally near } \pi_{x} \big( \CH_{{\rm tang} , -}^{\varepsilon} \big) .
\end{equation*}
Moreover, we have
\begin{equation} \label{j30}
a_{+} ( x , h ) = \sum_{1 \leq j \leq J} a_{+}^{j} ( x , h ) + S ( h^{\infty} ) ,
\end{equation}
with
\begin{equation*}
a_{+}^{j} ( x , h ) = h^{- i \frac{z - E_{0}}{h \lambda}} \int_{H_{j}} e^{i (\varphi_{+}^{1} (y) - \varphi_{-} (y) ) / h} d^{j}_{0} ( x, y , z ) e^{- i A^{y} / h} a_{-}^{j} ( y , h ) \, d y + S ( h^{- N_{-}^{\rm tang} + \zeta} ) ,
\end{equation*}
for some $0 < \zeta < 1 / 2$. Applying \eqref{a99}, \eqref{j2} which provides the asymptotic of $d_{0}^{j}$ and \eqref{b41} which express $a_{-}^{j}$ on $H_{j}$ in terms of $a_{-}$ on $\varepsilon \S^{n - 1}$, the previous equation becomes
\begin{align*}
a_{+}^{j} ( x , h ) &= \int_{H_{j}} \big( \varepsilon f ( x , Y (y) ) + \CO ( \varepsilon^{2 - n} ) + \CO_{\varepsilon} ( \nu ) \big) \varphi_{j} ( Y (y) - y_{j} )    \\
&\qquad \qquad \qquad \qquad \times e^{i (\varphi_{+}^{1} (y) - \varphi_{-} (y) - i A^{y} ) / h} a_{-} ( Y (y) , h ) \, d y + \CO_{\varepsilon , \nu} ( h^{- N_{-}^{\rm tang} + \zeta} ) ,
\end{align*}
where
\begin{equation*}
f  ( x , Y ) = h^{- i \frac{z - E_{0}}{h \lambda}} e^{- i n \frac{\pi}{4}} \Big( \frac{\lambda}{2 \pi} \Big)^{\frac{n}{2}} \big( i \lambda x \cdot Y \big)^{- \frac{n}{2} + i \frac{z - E_{0}}{\lambda h}}  \Gamma \Big( \frac{n}{2} - i \frac{z - E_{0}}{\lambda h} \Big) .
\end{equation*}
Performing the change of variables $H_{j} \ni y \mapsto Y = Y (y) \in \varepsilon \S^{n-1}$ as in \eqref{n2}, it yields
\begin{align}
a_{+}^{j} ( x , h ) &= \int_{\varepsilon \S^{n - 1}} \big( \varepsilon f ( x , Y ) + \CO ( \varepsilon^{2 - n} ) + \CO_{\varepsilon} ( \nu ) \big) \varphi_{j} ( Y - y_{j} )    \nonumber \\
&\qquad \qquad \qquad \times e^{i (\varphi_{+}^{1} ( y ( Y ) ) - \varphi_{-} ( y ( Y ) ) - A^{y ( Y )} ) / h} a_{-} ( Y , h ) \, d Y + \CO_{\varepsilon , \nu} ( h^{- N_{-}^{\rm tang} + \zeta} ) . \label{j4}
\end{align}
Note that the remainder terms are uniform with respect to $j$. Thus, summing over $j$, \eqref{j31}--\eqref{j30} and \eqref{j4} give
\begin{align*}
a_{+} ( x , h ) &= \int_{\varepsilon \S^{n - 1}} \big( \varepsilon f ( x , Y ) + \CO ( \varepsilon^{2 - n} ) + \CO_{\varepsilon} ( \nu ) \big) e^{i (\varphi_{+}^{1} ( y ( Y ) ) - \varphi_{-} ( y ( Y ) ) - A^{y ( Y )} ) / h}    \\
&\qquad \qquad \qquad \qquad \qquad \qquad \qquad \qquad \qquad \times  a_{-} ( Y , h ) \, d Y + \CO_{\varepsilon , \nu} ( h^{- N_{-}^{\rm tang} + \zeta} ) .
\end{align*}
Choosing now $\nu$ small enough depending on $\varepsilon$ and using that $a_{-} ( Y , h ) = 0$ for $Y \in \varepsilon \S^{n - 1} \setminus V_{-}^{\varepsilon}$ by construction, it yields
\begin{align*}
a_{+} ( x , h ) &= \int_{V_{-}^{\varepsilon}} \big( \varepsilon f ( x , Y ) + \CO ( \varepsilon^{2 - n} ) \big) e^{i (\varphi_{+}^{1} ( y ( Y ) ) - \varphi_{-} ( y ( Y ) ) - A^{y ( Y )} ) / h}  \nonumber \\
&\qquad \qquad \qquad \qquad \qquad \qquad \qquad \qquad \times  a_{-} ( Y , h ) \, d Y + \CO_{\varepsilon} ( h^{- N_{-}^{\rm tang} + \zeta} )   \\
&= \int_{V_{-}^{\varepsilon}} \varepsilon f ( x , Y ) e^{i (\varphi_{+}^{1} ( y ( Y ) ) - \varphi_{-} ( y ( Y ) ) - A^{y ( Y )} ) / h} a_{-} ( Y , h ) \, d Y    \\
&\qquad \qquad \qquad \qquad \qquad \qquad + \CO ( \varepsilon ) \Vert a_{-} ( \varepsilon \cdot , h ) \Vert_{L^{2} ( \S^{n - 1} )} + \CO_{\varepsilon} ( h^{- N_{-}^{\rm tang} + \zeta} ) .
\end{align*}
Thanks to the regularity of the measure, $V_{-}^{\varepsilon}$ can be taken arbitrarily close to the compact set $\pi_{x} ( \CH_{{\rm tang} , -}^{\varepsilon} )$. Moreover, for $Y \in \pi_{x} ( \CH_{{\rm tang} , -}^{\varepsilon} )$, the Hamiltonian curve $\gamma ( t ) = ( y ( t ) , \eta ( t ) ) =  \exp ( t H_{p} ) ( y ( Y ) , \nabla \varphi_{+}^{1} ( y ( Y ) ) )$ lies both in $\Lambda_{-}^{0}$ and $\Lambda_{+}^{1}$ and then
\begin{equation*}
\varphi_{+}^{1} ( y ( Y ) ) - \varphi_{-} ( y ( Y ) ) = \int_{\gamma ( ] - \infty , 0 ] )} \xi \cdot d x - \varphi_{-} ( 0 ) + \int_{\gamma ( ] 0 , + \infty ] )} \xi \cdot d x = A^{y ( Y )} .
\end{equation*}
Combining these arguments, the previous equation becomes
\begin{align*}
a_{+} ( x , h ) &= \int_{\pi_{x} ( \CH_{{\rm tang} , -}^{\varepsilon} ) / \varepsilon} \varepsilon^{2 i \frac{z - E_{0}}{\lambda h}} f ( \varepsilon^{- 1} x , \widetilde{\omega} ) a_{-} ( \varepsilon \widetilde{\omega} , h ) \, d \widetilde{\omega}      \\
&\qquad \qquad \qquad \qquad \qquad + \CO ( \varepsilon ) \Vert a_{-} ( \varepsilon \cdot , h ) \Vert_{L^{2} ( \S^{n - 1} )} + \CO_{\varepsilon} ( h^{- N_{-}^{\rm tang} + \zeta} ) .
\end{align*}
We now remark that $\vert f ( \varepsilon^{- 1} x , \widetilde{\omega} ) \vert \lesssim 1$ uniformly for $x \in \varepsilon \S^{n - 1} \cap \supp a_{+}$ and $\widetilde{\omega} \in \S^{n - 1} \cap \supp a_{-} ( \varepsilon \cdot , h )$. Thus, Lemma \ref{b46} and this estimate allow to replace $\pi_{x} ( \CH_{{\rm tang} , -}^{\varepsilon} ) / \varepsilon$ by $\CH_{\rm tang}^{- \infty}$ in the domain of integration. It yields
\begin{align}
a_{+} ( x , h ) &= \int_{\CH_{\rm tang}^{- \infty}} \varepsilon^{2 i \frac{z - E_{0}}{\lambda h}} f ( \varepsilon^{- 1} x , \widetilde{\omega} ) a_{-} ( \varepsilon \widetilde{\omega} , h ) \, d \widetilde{\omega}     \nonumber \\
&\qquad \qquad \qquad \qquad \qquad + o_{\varepsilon  \to 0} ( 1 ) \Vert a_{-} ( \varepsilon \cdot , h ) \Vert_{L^{2} ( \S^{n - 1} )} + \CO_{\varepsilon} ( h^{- N_{-}^{\rm tang} + \zeta} ) . \label{j10}
\end{align}

We now express $a_{-}$ in terms of $a_{+}$ following the evolution of $u$ along the homoclinic set $\CH$. Here, we can directly apply \eqref{b51} which says
\begin{equation} \label{j11}
a_{-} ( \varepsilon \omega , h ) = e^{i A^{\varepsilon \omega} / h} e^{- i \nu_{\omega} \frac{\pi}{2}} e^{i T_{\omega} \frac{z - E_{0}}{h}} \CM_{\varepsilon} ( \alpha_{\omega} ) \chi_{-}^{\varepsilon} ( \varepsilon \omega ) a_{+} ( x_{\omega} , h ) + S ( h^{- N_{-}^{\rm tang} + 1} ) .
\end{equation}
As explained below \eqref{b54}, the coefficient $\CM_{\varepsilon} ( \alpha_{\omega} )$ is uniformly bounded with respect to $\omega \in \supp a_{-} ( \varepsilon \cdot , h )$ and $\varepsilon$. Thus, \eqref{j10} and \eqref{j11} imply
\begin{equation} \label{j14}
\vert a_{-} ( \varepsilon \omega , h ) \vert \leq R \Vert a_{-} ( \varepsilon \cdot , h ) \Vert_{L^{2} ( \CH_{\rm tang}^{- \infty} )} + o_{\varepsilon  \to 0} ( 1 ) \Vert a_{-} ( \varepsilon \cdot , h ) \Vert_{L^{2} ( \S^{n - 1} )} + \CO_{\varepsilon} ( h^{- N_{-}^{\rm tang} + \zeta} ) ,
\end{equation}
uniformly for $\omega \in \S^{n - 1}$. Taking the $L^{2}$ norm leads to
\begin{equation} \label{j16}
\Vert a_{-} ( \varepsilon \cdot , h ) \Vert_{L^{2} ( \S^{n - 1} )} \leq R \Vert a_{-} ( \varepsilon \cdot , h ) \Vert_{L^{2} ( \CH_{\rm tang}^{- \infty} )} + \CO_{\varepsilon} ( h^{- N_{-}^{\rm tang} + \zeta} ) ,
\end{equation}
for $\varepsilon$ small enough. Using again \eqref{j14}, we finally obtain \eqref{i98}.

We pass to the proof of \eqref{i99}. Combining \eqref{j10}, \eqref{j11} and \eqref{j16}, we obtain
\begin{align} 
a_{-} ( \varepsilon \omega , & h ) = \one_{\pi_{x} ( \CH_{{\rm tang} , -}^{\varepsilon} ) / \varepsilon} e^{i A^{\varepsilon \omega} / h} e^{- i \nu_{\omega} \frac{\pi}{2}} \varepsilon^{2 i \frac{z - E_{0}}{\lambda h}} e^{i T_{\omega} \frac{z - E_{0}}{h}} \CM_{\varepsilon} ( \alpha_{\omega} )     \nonumber \\
&\times \int_{\CH_{\rm tang}^{- \infty}} f ( \varepsilon^{- 1} x_{\omega} , \widetilde{\omega} ) a_{-} ( \varepsilon \widetilde{\omega} , h ) \, d \widetilde{\omega} + o_{\varepsilon  \to 0} ( 1 ) \Vert a_{-} ( \varepsilon \cdot , h ) \Vert_{L^{2} ( \CH_{\rm tang}^{- \infty} )} + \CO_{\varepsilon} ( h^{- N_{-}^{\rm tang} + \zeta} ) , \label{n3}
\end{align}
as functions in $L^{2} ( \CH_{\rm tang}^{- \infty} )$. In the previous equation, we have used \eqref{j15} and Lemma \ref{b46} in order to replace $\chi_{-}^{\varepsilon} ( \varepsilon \cdot )$ by $\one_{\pi_{x} ( \CH_{{\rm tang} , -}^{\varepsilon} ) / \varepsilon}$. Assume that $\omega \in \CH_{\rm tang}^{- \infty} \cap \pi_{x} ( \CH_{{\rm tang} , -}^{\varepsilon} ) / \varepsilon$. In this case, the Hamiltonian trajectory $\gamma ( t ) = \exp ( t H_{p} ) ( \varepsilon \omega , \nabla \varphi_{-} ( \varepsilon \omega ) )$ is homoclinic and its asymptotic direction $\alpha_{\omega} \in \CH_{\rm tang}^{+ \infty}$ is given by \eqref{j12}. There is no reason for $\alpha_{\omega}$ to be equal to $\alpha ( \omega )$ even if \eqref{b56}, \eqref{a64} and Lemma \ref{b67} imply that
\begin{equation} \label{j13}
\alpha_{\omega} = \alpha ( \omega ) + o_{\varepsilon  \to 0} ( 1 ) ,
\end{equation}
uniformly with respect to $\omega$. In particular, we have
\begin{equation} \label{j17}
e^{i A^{\varepsilon \omega} / h} = e^{i A ( \alpha_{\omega} ) / h} = e^{i A ( \alpha ( \omega ) ) / h} ,
\end{equation}
for $\varepsilon$ small enough, since the action is locally constant from \ref{h7}. This is why we make this hypothesis. For the same reason, we have $e^{- i \nu_{\omega} \frac{\pi}{2}} = e^{- i \nu ( \alpha ( \omega ) ) \frac{\pi}{2}}$. On the other hand, the uniform convergence of $\CM_{\varepsilon} ( \alpha )$ to $\CM_{0} ( \alpha )$ (see \eqref{i49}) and \eqref{j13} yield
\begin{equation*}
\CM_{\varepsilon} ( \alpha_{\omega} ) = \CM_{0} ( \alpha ( \omega ) ) + o_{\varepsilon  \to 0} ( 1 ) .
\end{equation*}
It remains to deal with $T_{\omega}$. By definition, we have $T_{\omega} = t_{-}^{\varepsilon} ( \alpha_{\omega} ) - t_{+}^{\varepsilon} ( \alpha_{\omega} )$. By Proposition \ref{a65} and Proposition \ref{a74}, we get
\begin{equation*}
\varepsilon^{2 i \frac{z - E_{0}}{\lambda h}} e^{i T_{\omega} \frac{z - E_{0}}{h}} = e^{i q ( \alpha_{\omega} ) \frac{z - E_{0}}{h} + i \varepsilon ( r - m ) ( \varepsilon , \alpha_{\omega} ) \frac{z - E_{0}}{h}} ,
\end{equation*}
with $m , r \in C^{0} ( [ - \varepsilon_{0} , \varepsilon_{0} ] \times \CH_{\rm tang}^{+ \infty} )$. Moreover, $q ( \alpha ) = T ( \alpha ) \in C^{0} ( \CH_{\rm tang}^{+ \infty} )$ from \eqref{i50}. Combining with \eqref{j13}, we obtain
\begin{equation*}
\varepsilon^{2 i \frac{z - E_{0}}{\lambda h}} e^{i T_{\omega} \frac{z - E_{0}}{h}} = e^{i T ( \alpha ( \omega ) ) \frac{z - E_{0}}{h}} + o_{\varepsilon  \to 0} ( 1 ) ,
\end{equation*}
uniformly with respect to $\omega$. Lastly, using $x_{\omega} = \varepsilon F^{-1} ( \varepsilon , \alpha ( F ( \varepsilon , \omega ) ) )$ (see \eqref{b56}), \eqref{a64} and \eqref{a68}, we deduce $x_{\omega} = \varepsilon \alpha ( \omega ) + o_{\varepsilon  \to 0} ( \varepsilon )$ uniformly with respect to $\omega$. From the form of $f$, it implies
\begin{equation*}
f ( \varepsilon^{- 1} x_{\omega} , \widetilde{\omega} ) = f ( \alpha ( \omega ) , \widetilde{\omega} ) + o_{\varepsilon  \to 0} ( 1 ) ,
\end{equation*}
uniformly with respect to $\omega , \widetilde{\omega}$. Using the previous asymptotics together with \eqref{n3}, we obtain
\begin{align}
a_{-} ( \varepsilon \omega , h ) ={}& \one_{\pi_{x} ( \CH_{{\rm tang} , -}^{\varepsilon} ) / \varepsilon} e^{i A ( \alpha ( \omega ) ) / h} e^{- i \nu ( \alpha ( \omega ) ) \frac{\pi}{2}} e^{i T ( \alpha ( \omega ) ) \frac{z - E_{0}}{\lambda h}} \CM_{0} ( \alpha ( \omega ) )  \nonumber \\
&\times \int_{\CH_{\rm tang}^{- \infty}} f ( \alpha ( \omega ) , \widetilde{\omega} ) a_{-} ( \varepsilon \widetilde{\omega} , h ) \, d \widetilde{\omega} + o_{\varepsilon  \to 0} ( 1 ) \Vert a_{-} ( \varepsilon \cdot , h ) \Vert_{L^{2} ( \CH_{\rm tang}^{- \infty} )} + \CO_{\varepsilon} ( h^{- N_{-}^{\rm tang} + \zeta} ) ,
\end{align}
as functions in $L^{2} ( \CH_{\rm tang}^{- \infty} )$. In this formula, we recognize the kernel of the operator $\CT ( z , h )$ given in \eqref{c19} and we get
\begin{align*}
a_{-} ( \varepsilon \omega , h ) = \one_{\pi_{x} ( \CH_{{\rm tang} , -}^{\varepsilon} ) / \varepsilon} & h^{- i \frac{z - E_{0}}{h \lambda}} \CT ( z , h ) a_{-} ( \varepsilon \cdot , h ) \\
&\qquad + o_{\varepsilon  \to 0} ( 1 ) \Vert a_{-} ( \varepsilon \cdot , h ) \Vert_{L^{2} ( \CH_{\rm tang}^{- \infty} )} + \CO_{\varepsilon} ( h^{- N_{-}^{\rm tang} + \zeta} ) .
\end{align*}
Using \eqref{c19}, \eqref{j15} and Lemma \ref{b46}, we have
\begin{equation*}
\big\Vert \big( 1 - \one_{\pi_{x} ( \CH_{{\rm tang} , -}^{\varepsilon} ) / \varepsilon} \big) \CT ( z , h ) ( \omega , \widetilde{\omega} ) \big\Vert_{L^{2} ( \CH_{\rm tang}^{- \infty} \times \CH_{\rm tang}^{- \infty} )} = o_{\varepsilon  \to 0} ( 1 ) .
\end{equation*}
Then, we finally obtain
\begin{equation*}
a_{-} ( \varepsilon \cdot, h ) = h^{- i \frac{z - E_{0}}{h \lambda}} \CT ( z , h ) a_{-} ( \varepsilon \cdot , h ) + o_{\varepsilon  \to 0} ( 1 ) \Vert a_{-} ( \varepsilon \cdot , h ) \Vert_{L^{2} ( \CH_{\rm tang}^{- \infty} )} + \CO_{\varepsilon} ( h^{- N_{-}^{\rm tang} + \zeta} ) ,
\end{equation*}
as functions in $L^{2} ( \CH_{\rm tang}^{- \infty} )$. This is precisely \eqref{i99}.
\end{proof}

From \eqref{i99}, we can write
\begin{equation*}
\big( 1 - h^{- i \frac{z - E_{0}}{h \lambda}} \CT ( z , h ) \big) a_{-} ( \varepsilon \cdot, h ) = o_{\varepsilon  \to 0} ( 1 ) \Vert a_{-} ( \varepsilon \cdot , h ) \Vert_{L^{2} ( \CH_{\rm tang}^{- \infty} )} + \CO_{\varepsilon} ( h^{- N_{-}^{\rm tang} + \zeta} ) .
\end{equation*}
Using that $z \in \eqref{i94}$ is at distance $h \vert \ln h \vert^{-1}$ of the pseudo-resonances, Lemma \ref{i92} gives
\begin{align*}
\Vert a_{-} ( \varepsilon \cdot, h ) \Vert_{L^{2} ( \CH_{\rm tang}^{- \infty} )} &\leq o_{\varepsilon  \to 0} ( 1 ) \Vert a_{-} ( \varepsilon \cdot , h ) \Vert_{L^{2} ( \CH_{\rm tang}^{- \infty} )} + \CO_{\varepsilon} ( h^{- N_{-}^{\rm tang} + \zeta} )    \\
&\leq \CO ( h^{- N_{-}^{\rm tang} + \zeta} ) ,
\end{align*}
for $\varepsilon$ small enough. Applying \eqref{i98}, we have
\begin{equation}
\Vert a_{-} ( \varepsilon \cdot , h ) \Vert_{L^{\infty} ( \S^{n - 1} )} = \CO ( h^{- N_{-}^{\rm tang} + \zeta} ) .
\end{equation}
Now the rest of the proof follows the end of Section \ref{s52}. We successively show that the $L^{\infty}$ norm of $a_{-} ( \cdot , h )$ is a $\CO ( h^{- N_{-}^{\rm tang} + \zeta} )$ in a vicinity of $\varepsilon \S^{n - 1}$, that $a_{-} \in S ( h^{- N_{-}^{\rm tang} + \frac{\zeta}{2}} )$ by the Landau--Kolmogorov inequalities, that $u_{-} = \CO ( h^{\infty} )$ by the standard bootstrap argument and finally that $u = \CO ( h^{\infty} )$ thanks to Proposition \ref{a16}. Thus, Proposition \ref{i93} follows from the general strategy of Section \ref{s31}.

\Subsection{Existence of resonances near the pseudo-resonances} \label{s29}

To prove Theorem \ref{i55}, it remains to show that $P$ has at least one resonance near each pseudo-resonance. This is done in the following proposition.

\begin{proposition}\sl \label{j18}
Let $C , r > 0$ and assume that $h$ is small enough. For any pseudo-resonance $z$ in the set \eqref{i54}, the operator $P$ has at least one resonance in $ B ( z , r h \vert \ln h \vert^{- 1} )$.
\end{proposition}

This result is the analogous to Proposition \ref{d50} in the present strong trapping geometry. It will be proved the same way following Section \ref{s73}. Note that Theorem \ref{i55} is a direct consequence of Proposition \ref{i93} and Proposition \ref{j18}.

As in the beginning of Section \ref{s73}, we assume that Proposition \ref{j18} does not hold true. The compactness argument of \eqref{d59} can be adapted here. Thus, there exists a (decreasing) sequence of $h$ which goes to $0$ such that
\begin{equation*}
\tau ( h ) = \re \frac{z - E_{0}}{h} \longrightarrow \tau_{0} \in [ - C , C ] ,
\end{equation*}
$\rho ( h ) \longrightarrow \rho_{0} \in ( \S^{1} )^{K}$ and
\begin{equation} \label{j19}
P \text{ has no resonance in } B \Big( z_{q}^{0} , r \frac{h}{\vert \ln h \vert} \Big) ,
\end{equation}
where
\begin{equation*}
z_{q}^{0} = E_{0} + 2 q \pi \lambda \frac{h}{\vert \ln h \vert}+ i \ln ( \mu_{0} ) \lambda \frac{h}{\vert \ln h \vert} \in\eqref{i54} ,
\end{equation*}
the complex number $\mu_{0} \neq 0$ is an eigenvalue of $\CT_{0} : = \widetilde{\CT} ( \rho_{0} , \tau_{0} )$ and $q = q ( h ) \in \Z$ is such that $\re z_{q}^{0} = \tau_{0} + o_{h \to 0} ( 1 )$.

\begin{figure}
\begin{center}
\begin{picture}(0,0)%
\includegraphics{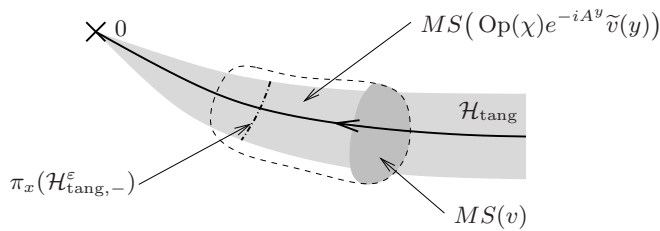}%
\end{picture}%
\setlength{\unitlength}{1579sp}%
\begingroup\makeatletter\ifx\SetFigFont\undefined%
\gdef\SetFigFont#1#2#3#4#5{%
  \reset@font\fontsize{#1}{#2pt}%
  \fontfamily{#3}\fontseries{#4}\fontshape{#5}%
  \selectfont}%
\fi\endgroup%
\begin{picture}(8866,3243)(2236,-7021)
\put(8701,-3961){\makebox(0,0)[lb]{\smash{{\SetFigFont{9}{10.8}{\rmdefault}{\mddefault}{\updefault}$MS \big( \Op ( \chi ) e^{- i A^{y}} \widetilde{v} ( y ) \big)$}}}}
\put(3901,-4036){\makebox(0,0)[lb]{\smash{{\SetFigFont{9}{10.8}{\rmdefault}{\mddefault}{\updefault}$0$}}}}
\put(9301,-5311){\makebox(0,0)[lb]{\smash{{\SetFigFont{9}{10.8}{\rmdefault}{\mddefault}{\updefault}$\CH_{\rm tang}$}}}}
\put(9226,-6961){\makebox(0,0)[lb]{\smash{{\SetFigFont{9}{10.8}{\rmdefault}{\mddefault}{\updefault}$MS ( v )$}}}}
\put(2251,-6436){\makebox(0,0)[lb]{\smash{{\SetFigFont{9}{10.8}{\rmdefault}{\mddefault}{\updefault}$\pi_{x} ( \CH^{\varepsilon}_{{\rm tang} , -} )$}}}}
\end{picture}%
\end{center}
\caption{The microsupport of $v$.} \label{f45}
\end{figure}

We now construct a test function as in \eqref{d60}. Let $w_{0} \in L^{2} ( \CH_{\rm tang}^{- \infty} )$ be a normalized eigenvector of $\CT_{0}$ associated to the eigenvalue $\mu_{0}$. We extend $w_{0}$ by $0$ on the set $\S^{n - 1} \setminus \CH_{\rm tang}^{- \infty}$. Since $w_{0}$ is a priori non-regular, we consider $w \in C^{\infty} ( \S^{n - 1} )$ satisfying
\begin{equation} \label{j20}
\Vert w - w_{0} \Vert_{L^{2} ( \S^{n - 1} )} \leq \eta ,
\end{equation}
for some small $\eta > 0$ which will be fixed in the sequel. Let $\widetilde{v}$ be a WKB solution of
\begin{equation*}
\left\{ \begin{aligned}
&( P - \widetilde{z} ) \widetilde{v} = 0 &&\text{near } \varepsilon \S^{n - 1} ,   \\
&\widetilde{v} (x) = \chi_{-}^{\varepsilon} ( x ) w ( x / \varepsilon ) e^{i \varphi_{+}^{1} (x) / h} &&\text{on } \varepsilon \S^{n - 1} ,
\end{aligned} \right.
\end{equation*}
holomorphic in $\widetilde{z} \in B ( E_{0} , ( C + 1 ) h )$. In particular, $\widetilde{v} ( x ) = \widetilde{a} ( x , h ) e^{i \varphi_{+}^{1} (x) / h}$ is a Lagrangian distribution whose symbol $\widetilde{a} \in S (1)$ vanishes outside a small neighborhood of $\pi_{x} ( \CH_{{\rm tang} , -}^{\varepsilon} )$ (near $\varepsilon \S^{n - 1}$) and has a classical expansion $\widetilde{a} ( x , h ) = \widetilde{a}^{0} (x) + \widetilde{a}^{1} (x) h + \cdots$ in $S (1)$. Consider now cut-off functions $\chi , \psi \in C^{\infty}_{0} ( T^{*} \R^{n} )$ such that $\chi = 1$ near $V_{-}^{\varepsilon}$, $\psi = 1$ near $\supp ( \nabla \chi ) \cap \{ \exp ( t H_{p} ) ( x , \nabla \varphi_{+}^{1} ( x ) ) ; \ t < 0 \text{ and } x \in V_{-}^{\varepsilon} \}$ and such that $\chi$ and $\psi$ are supported in the region where $\widetilde{v}$ is defined (see Figure \ref{f45}). Then, our test function is
\begin{equation} \label{j21}
v = \Op ( \psi ) \big[ P , \Op ( \chi ) \big] e^{- i A^{y} / h} \widetilde{v} ( y ) ,
\end{equation}
which is holomorphic in \eqref{i54}. Note that, the function $A^{y}$ is smooth on the support of $\widetilde{v}$ thanks to \eqref{j20}.

As in Section \ref{s73}, we denote $\CD = B ( z_{q}^{0} , s h \vert \ln h \vert^{-1} )$ for some $0 < s < r$ fixed below. Let $u \in L^{2} ( \R^{n} )$ be the solution of
\begin{equation} \label{j25}
( P_{\theta} - \widetilde{z} ) u = v ,
\end{equation}
for $\widetilde{z} \in \partial \CD$. The set $\partial \CD$ is at distance $h \vert \ln h \vert^{- 1}$ from the pseudo-resonances. From Proposition \ref{i93}, $u = ( P_{\theta} - \widetilde{z} )^{-1} v$ is well-defined, holomorphic near $\partial \CD$ and $\Vert u \Vert \lesssim h^{- M}$ uniformly for $\widetilde{z} \in \partial \CD$. We define $u_{-}^{\rm tang}$ and $u_{-}^{\rm trans}$ as in Section \ref{s51}. Following Section \ref{s52}, we get
\begin{equation*}
u_{-}^{\bullet} \in \CI ( \Lambda_{+}^{1} , h^{- N_{-}^{\bullet}} ) ,
\end{equation*}
for some $N_{-}^{\bullet} \in \R$ uniformly for $\widetilde{z} \in \partial \CD$. Furthermore, $N_{-}^{\rm trans} \leq N_{-}^{\rm tang}$. As in the previous part, we write
\begin{equation*}
u_{-}^{\rm tang} ( y , h ) = a_{-} ( y , h ) e^{- i A^{y} / h} e^{i \varphi^{1}_{+} (y) / h} ,
\end{equation*}
for some $a_{-} \in S ( h^{- N_{-}^{\rm tang}} )$. Following Section \ref{s73} and Section \ref{s28}, one can prove the following result.

\begin{lemma}\sl \label{j22}
There exist $R , \zeta > 0$ such that, for $\varepsilon$ small enough and $V_{-}^{\varepsilon}$ sufficiently close to $\pi_{x} ( \CH_{{\rm tang} , -}^{\varepsilon} )$, we have, uniformly for $\widetilde{z} \in \partial \CD$,
\begin{equation} \label{n4}
\Vert a_{-} ( \varepsilon \cdot , h ) \Vert_{L^{\infty} ( \S^{n - 1} )} \leq R \Vert a_{-} ( \varepsilon \cdot , h ) \Vert_{L^{2} ( \CH_{\rm tang}^{- \infty} )} + \CO_{\eta} ( 1 ) + \CO_{\varepsilon , \eta} ( h^{- N_{-}^{\rm tang} + \zeta} ) ,
\end{equation}
and
\begin{align}
a_{-} ( \varepsilon \cdot, h ) &= h^{- i \frac{\widetilde{z} - E_{0}}{h \lambda}} \CT ( \widetilde{z} , h ) a_{-} ( \varepsilon \cdot , h ) + \widetilde{a} ( \varepsilon \cdot , h )   \nonumber \\
&\qquad \qquad \quad + o_{\varepsilon  \to 0} ( 1 ) \Vert a_{-} ( \varepsilon \cdot , h ) \Vert_{L^{2} ( \CH_{\rm tang}^{- \infty} )} + o_{\varepsilon  \to 0}^{\eta} ( 1 ) + \CO_{\varepsilon , \eta} ( h^{-N _{-}^{\rm tang} + \zeta} ) ,   \label{j24}
\end{align}
as functions in $L^{2} ( \CH_{\rm tang}^{- \infty} )$. The notation $o_{a \to 0}^{b} ( 1 )$ designs a function which tends to $0$ as $a$ goes to $0$ with $b$ fixed.
\end{lemma}

\begin{proof}
This result is similar to Lemma \ref{i97} and we follow its proof. By construction, the test function $v$ does not affect the propagation of singularities through the fixed point. Thus, \eqref{j10} still holds true. Moreover, the propagation along $\CH_{\rm tang}$ yields that \eqref{j11} is replaced by the equation
\begin{equation*}
a_{-} ( \varepsilon \omega , h ) = e^{i A^{\varepsilon \omega} / h} e^{- i \nu_{\omega} \frac{\pi}{2}} e^{i T_{\omega} \frac{z - E_{0}}{h}} \CM_{\varepsilon} ( \alpha_{\omega} ) \chi_{-}^{\varepsilon} ( \varepsilon \omega ) a_{+} ( x_{\omega} , h ) + \widetilde{a} ( \varepsilon \omega , h ) + S ( h^{- N_{-}^{\rm tang} + 1} ) ,
\end{equation*}
for $\omega \in \S^{n - 1}$. This expression is analogous to \eqref{d71}. Since
\begin{equation*}
\Vert \widetilde{a} ( \varepsilon \omega , h ) \Vert_{L^{\infty} ( \S^{n - 1} )} = \big\Vert \chi_{-}^{\varepsilon} ( \varepsilon \omega ) w ( \omega ) \big\Vert_{L^{\infty} ( \S^{n - 1} )} = \CO_{\eta} ( 1 ) ,
\end{equation*}
we obtain \eqref{n4} as \eqref{i98}. On the other hand, \eqref{j24} is proved as \eqref{i99}. The additional term $o_{\varepsilon  \to 0}^{\eta} ( 1 )$ comes from the remainder term $\CO_{\eta} ( 1 )$ of \eqref{n4}.
\end{proof}

We define the new spectral parameter
\begin{equation*}
\Lambda = \Lambda ( \widetilde{z} , h ) = h^{i \frac{\widetilde{z} - E_{0}}{h \lambda}} ,
\end{equation*}
which satisfies $1 \lesssim \vert \Lambda \vert \lesssim 1$ uniformly for $\widetilde{z} \in \partial \CD$. We now fix $s$ small enough such that $\Lambda ( \partial \CD )$ does not meet $\spe ( \CT_{0} ) \setminus \{ \mu_{0} \}$. Then, $\spe ( \CT_{0} ) \cap \CD = \{ \mu_{0} \}$, $\Lambda ( \partial \CD )$ is a simple loop around $\mu_{0}$ and the spectrum of $\CT_{0}$ is at distance $\beta > 0$ of $\Lambda ( \partial \CD )$. This can be proved as in Lemma \ref{d80}. We then obtain the

\begin{lemma}\sl \label{n5}
For $\varepsilon$ small enough and $V_{-}^{\varepsilon}$ sufficiently close to $\pi_{x} ( \CH_{{\rm tang} , -}^{\varepsilon} )$, we have $u_{-}^{\rm tang} \in \CI ( \Lambda_{+}^{1} , h^{- c} )$ for all $c > 0$ and
\begin{equation} \label{n6}
\Vert a_{-} ( \varepsilon \cdot , h ) \Vert_{L^{\infty} ( \S^{n - 1} )} = \CO_{\varepsilon , \eta} ( 1 ) ,
\end{equation}
uniformly for $\widetilde{z} \in \partial \CD$.
\end{lemma}

In other words, $N_{-}^{\rm tang} = c$ for any $c > 0$. In the sequel, we take $c = \zeta / 2$. This result follows from Lemma \ref{j22} and the bootstrap argument at the end of Proposition \ref{i93}. It is limited by $\widetilde{a} ( \varepsilon \cdot , h ) = \CO_{\eta} ( 1 )$. The difference with Lemma \ref{d77} is that the Landau--Kolmogorov inequalities do not allow to reach $u_{-}^{\rm tang} \in \CI ( \Lambda_{+}^{1} , 1 )$ but only $u_{-}^{\rm tang} \in \CI ( \Lambda_{+}^{1} , h^{- c} )$ for all $c > 0$. Applying again Lemma \ref{j22}, we deduce \eqref{n6}.

\begin{lemma}\sl \label{j23}
For $\varepsilon$ small enough and $V_{-}^{\varepsilon}$ sufficiently close to $\pi_{x} ( \CH_{{\rm tang} , -}^{\varepsilon} )$, we have
\begin{equation} \label{j28}
\frac{1}{2 i \pi} \int_{\partial \CD} u ( \widetilde{z} ) \, d \widetilde{z} = b ( x , h ) e^{i \varphi_{+}^{1} (x) / h} \text{ microlocally near } \CH_{-}^{\varepsilon} ,
\end{equation}
where $b \in S ( h^{1 - \zeta / 2} \vert \ln h \vert^{-1} )$. Moreover,
\begin{equation} \label{j29}
b ( \varepsilon \cdot , h ) = \frac{i \lambda}{2 \pi} \frac{h}{\vert \ln h \vert} w_{0} + \frac{h}{\vert \ln h \vert} \big( \CO ( \eta ) + o_{\varepsilon  \to 0}^{\eta} ( 1 ) + o_{h \to 0}^{\varepsilon , \eta} ( 1 ) \big) ,
\end{equation}
as function in $L^{2} ( \CH_{\rm tang}^{- \infty} )$.
\end{lemma}

\begin{proof}
That $u$ is a holomorphic function in $\CI ( \Lambda_{+}^{1} , h^{- \zeta / 2} )$ microlocally near $\CH_{-}^{\varepsilon}$ implies \eqref{j28}. Since $\CT ( \widetilde{z} , h ) = \CT_{0} + o_{h \to 0} ( 1 )$ uniformly for $\widetilde{z} \in \partial \CD$, we deduce
\begin{equation*}
\dist \big( \spe ( \CT ( \widetilde{z} , h ) ) , \spe ( \CT_{0}) \big) = o_{h \to 0} (1) .
\end{equation*}
Then, Lemma \ref{i92} and the paragraph above Lemma \ref{n5} imply
\begin{equation} \label{j26}
\big\Vert \big( \Lambda - \CT ( \widetilde{z} , h ) \big)^{-1} \big\Vert + \Vert ( \Lambda - \CT_{0} )^{-1} \Vert \lesssim 1 ,
\end{equation}
uniformly for $\widetilde{z} \in \partial \CD$. Combining with the resolvent identity, we get
\begin{align*}
\big( \Lambda - \CT ( \widetilde{z} , h ) \big)^{- 1} ={}& ( \Lambda  - \CT_{0} )^{- 1} + \big( \Lambda - \CT ( \widetilde{z} , h ) \big)^{- 1} \big( \CT ( \widetilde{z} , h ) - \CT_{0} \big) ( \Lambda  - \CT_{0} )^{- 1}   \\
={}& ( \Lambda  - \CT_{0} )^{- 1} + o_{h \to 0} (1) ,
\end{align*}
uniformly for $\widetilde{z} \in \partial \CD$. Then, \eqref{j24} can be written
\begin{align}
a_{-} ( \varepsilon \cdot, h ) &= \Lambda \big( \Lambda - \CT ( \widetilde{z} , h ) \big)^{- 1} \widetilde{a} ( \varepsilon \cdot , h )  \nonumber \\
&\qquad \qquad \ + o_{\varepsilon  \to 0} ( 1 ) \Vert a_{-} ( \varepsilon \cdot , h ) \Vert_{L^{2} ( \CH_{\rm tang}^{- \infty} )} + o_{\varepsilon  \to 0}^{\eta} ( 1 ) + \CO_{\varepsilon , \eta} ( h^{\zeta / 2} ) \nonumber \\
&= \Lambda \big( \Lambda - \CT_{0} \big)^{- 1} \widetilde{a} ( \varepsilon \cdot , h ) + o_{h \to 0}^{\eta} ( 1 )   \nonumber \\
&\qquad \qquad \ + o_{\varepsilon  \to 0} ( 1 ) \Vert a_{-} ( \varepsilon \cdot , h ) \Vert_{L^{2} ( \CH_{\rm tang}^{- \infty} )} + o_{\varepsilon  \to 0}^{\eta} ( 1 ) + \CO_{\varepsilon , \eta} ( h^{\zeta / 2} ) ,  \label{j27}
\end{align}
as function in $L^{2} ( \CH_{\rm tang}^{- \infty} )$. Taking the $L^{2}$ norm and using \eqref{j20}, we get
\begin{align*}
\Vert a_{-} ( \varepsilon \cdot, h ) \Vert_{L^{2} ( \CH_{\rm tang}^{- \infty} )} &\leq M + o_{\varepsilon  \to 0} ( 1 ) \Vert a_{-} ( \varepsilon \cdot , h ) \Vert_{L^{2} ( \CH_{\rm tang}^{- \infty} )} + o_{\varepsilon  \to 0}^{\eta} ( 1 ) + o_{h \to 0}^{\varepsilon , \eta} (1) \\
&\leq M + o_{\varepsilon  \to 0}^{\eta} ( 1 ) + o_{h \to 0}^{\varepsilon , \eta} ( 1 ) ,
\end{align*}
for some $M > 0$ and $\varepsilon$ small enough. Combining with Lemma \ref{b46}, \eqref{j27} becomes
\begin{align}
a_{-} ( \varepsilon \cdot, h ) &= \Lambda \big( \Lambda - \CT_{0} \big)^{- 1} \widetilde{a} ( \varepsilon \cdot , h ) + o_{\varepsilon  \to 0}^{\eta} ( 1 ) + o_{h \to 0}^{\varepsilon , \eta} ( 1 )  \nonumber \\
&= \Lambda \big( \Lambda - \CT_{0} \big)^{- 1} \chi_{-}^{\varepsilon} ( \varepsilon \cdot ) w ( \cdot ) + o_{\varepsilon  \to 0}^{\eta} ( 1 ) + o_{h \to 0}^{\varepsilon , \eta} ( 1 )   \nonumber \\
&= \Lambda \big( \Lambda - \CT_{0} \big)^{- 1} w + o_{\varepsilon  \to 0}^{\eta} ( 1 ) + o_{h \to 0}^{\varepsilon , \eta} ( 1 ) ,
\end{align}
uniformly for $\widetilde{z} \in \partial \CD$. Integrating with respect to $\widetilde{z}$, we get
\begin{align*}
b ( \varepsilon \cdot , h ) &= \frac{1}{2 i \pi} \int_{\partial \CD} \Lambda \big( \Lambda - \CT_{0} \big)^{- 1} w \, d \widetilde{z} + \frac{h}{\vert \ln h \vert} \big( o_{\varepsilon  \to 0}^{\eta} ( 1 ) + o_{h \to 0}^{\varepsilon , \eta} ( 1 ) \big)   \\
&= \frac{\lambda}{2 \pi} \frac{h}{\vert \ln h \vert} \int_{\Lambda ( \partial \CD )} \big( \Lambda - \CT_{0} \big)^{- 1} w \, d \Lambda + \frac{h}{\vert \ln h \vert} \big( o_{\varepsilon  \to 0}^{\eta} ( 1 ) + o_{h \to 0}^{\varepsilon , \eta} ( 1 ) \big)   \\
&= \frac{i \lambda}{2 \pi} \frac{h}{\vert \ln h \vert} \Pi_{0} w + \frac{h}{\vert \ln h \vert} \big( o_{\varepsilon  \to 0}^{\eta} ( 1 ) + o_{h \to 0}^{\varepsilon , \eta} ( 1 ) \big)  \\
&= \frac{i \lambda}{2 \pi} \frac{h}{\vert \ln h \vert} w_{0} + \frac{h}{\vert \ln h \vert} \big( \CO ( \eta ) + o_{\varepsilon  \to 0}^{\eta} ( 1 ) + o_{h \to 0}^{\varepsilon , \eta} ( 1 ) \big) ,
\end{align*}
where $\Pi_{0}$ is the spectral projection of $\CT_{0}$ associated to the eigenvalue $\mu_{0}$.
\end{proof}

Eventually, we can end the proof of Proposition \ref{j18}. The function $v$ and the operator $( P_{\theta} - \widetilde{z} )^{- 1}$ are holomorphic in $z \in \CD$ (see \eqref{j19}). Then, 
\begin{equation*}
\frac{1}{2 i \pi} \int_{\partial \CD} u ( \widetilde{z} ) \, d \widetilde{z} = 0 ,
\end{equation*}
which is in contradiction with Lemma \ref{j23}, choosing $\eta$ small enough, then $\varepsilon$ small enough and eventually $h$ small enough.

\Subsection{Proof of the additional results of Section \ref{s26}} \label{s30}

We first compare the quantization operators $\CQ$ (Section \ref{s6}) and $\CT$ (Section \ref{s26}) in dimension $n = 1$. In this case, these two operators are square matrices whose size is the number of homoclinic trajectories.

\begin{lemma}\sl \label{i82}
In dimension $n = 1$, we have
\begin{equation*}
\CQ ( z , h ) = U^{- 1} \CT ( z , h ) U \qquad \text{with} \qquad  U = \diag \Big( \vert g_{-}^{k} \vert^{\frac{1}{2} + i \frac{z - E_{0}}{\lambda h}} \Big) .
\end{equation*}
In particular, the pseudo-resonances of Definition \ref{d1} and Definition \ref{i51} coincide.
\end{lemma}

\begin{proof}
We use the indices $k , \ell \in \{ - 1 , 1 \}$ and not the parameters $\omega , \alpha$. We have already seen in \eqref{e38} and above \eqref{b75} that $\CM_{0} = 1$ and $\CM_{k}^{\pm} = \sqrt{\lambda \vert g_{\pm}^{k} \vert}$. Moreover, from \eqref{d2} and the definition of $t_{\pm}^{\varepsilon}$, we have
\begin{equation*}
\big( \vert g_{\pm}^{k} \vert + o_{\varepsilon \to 0} ( 1 ) \big) e^{\pm \lambda t_{\pm}^{\varepsilon}} = \varepsilon .
\end{equation*}
Using \eqref{i50} and taking the limit $\varepsilon \to 0$, it gives
\begin{equation}
e^{\lambda T_{k}} = \vert g_{-}^{k} \vert \vert g_{+}^{k} \vert .
\end{equation}
Thus, \eqref{d4} becomes
\begin{align*}
\CQ_{k , \ell} &= e^{i A_{k} / h} \Gamma \Big( \frac{1}{2} - i \frac{z - E_{0}}{\lambda h} \Big) \Big( \frac{\lambda}{2 \pi} \Big)^{\frac{1}{2}} e^{- i ( \nu_{k} \frac{\pi}{2} + \frac{\pi}{4} )} \big( i \lambda \widehat{g}_{+}^{k} \cdot \widehat{g}_{-}^{\ell} \big)^{- \frac{1}{2} + i \frac{z - E_{0}}{\lambda h}} \frac{\vert g_{+}^{k} \vert^{i \frac{z - E_{0}}{\lambda h}}}{\vert g_{-}^{k} \vert^{\frac{1}{2}}} \vert g_{-}^{\ell} \vert^{\frac{1}{2} + i \frac{z - E_{0}}{\lambda h}} \\
&= e^{i A_{k} / h} \Gamma \Big( \frac{1}{2} - i \frac{z - E_{0}}{\lambda h} \Big) \Big( \frac{\lambda}{2 \pi} \Big)^{\frac{1}{2}} e^{- i ( \nu_{k} \frac{\pi}{2} + \frac{\pi}{4} - T_{k} \frac{z - E_{0}}{h} )} \big( i \lambda \widehat{g}_{+}^{k} \cdot \widehat{g}_{-}^{\ell} \big)^{- \frac{1}{2} + i \frac{z - E_{0}}{\lambda h}} \frac{\vert g_{-}^{\ell} \vert^{\frac{1}{2} + i \frac{z - E_{0}}{\lambda h}}}{\vert g_{-}^{k} \vert^{\frac{1}{2} + i \frac{z - E_{0}}{\lambda h}}} ,
\end{align*}
where $\widehat{g}_{\pm}^{\bullet} : = g_{\pm}^{\bullet} / \vert g_{\pm}^{\bullet} \vert$. Comparing with \eqref{c19}, we eventually get Lemma \ref{i82}.
\end{proof}

We then study the asymptotic behavior of the accumulation curves when the opening angle is small.

\begin{proof}[Proof of \eqref{i65}]
In this proof, we note $\widehat{\CT}_{\theta_{0}}$ the operator $\widehat{\CT} ( \tau , h )$ with opening angle $2 \theta_{0}$ and
\begin{equation} \label{i62}
k ( x , y ) = e^{i A / h} \frac{\lambda}{2 \pi} \Gamma \Big( 1 - i \frac{\tau}{\lambda} \Big) e^{i T \tau} \big( i \lambda \cos ( x - y ) \big)^{- 1 + i \frac{\tau}{\lambda}} ,
\end{equation}
the kernel defined for all $x , y \in \R$ with $\cos (x - y ) \neq 0$. Thus, the distribution kernel of $\widehat{\CT}_{\theta_{0}}$ is $k$ restricted to $[ - \theta_{0} , \theta_{0} ]^{2}$. Let also define
\begin{equation*}
\begin{aligned}
V_{\theta_{0}} : \\
{}^{}
\end{aligned}
\left\{ \begin{aligned}
& L^{2} ( [ - \theta_{0} , \theta_{0} ] ) && \longrightarrow && L^{2} ( [ - 1 , 1 ] )  \\
& f ( x ) && && \sqrt{\theta_{0}} f ( \theta_{0} x )
\end{aligned} \right.
\end{equation*}
which is a unitary transform. Then, the eigenvalues of $\widehat{\CT}_{\theta_{0}}$ are the eigenvalues of the operator $\theta_{0} K_{\theta_{0}} : = V_{\theta_{0}} \widehat{\CT}_{\theta_{0}} V_{\theta_{0}}^{- 1}$ acting on $L^{2} ( [ - 1 , 1 ] )$. A direct computation shows that the kernel of $K_{\theta_{0}}$ is given by
\begin{equation} \label{i63}
k_{\theta_{0}} ( x , y ) = k ( \theta_{0} x , \theta_{0} y ) ,
\end{equation}
for $x , y \in [ - 1 , 1 ]$. From \eqref{i62} (or more generally, from the properties stated in Section \ref{s53} on the quantities appearing in $\CT$), $k$ is a continuous function. Then, \eqref{i63} gives
\begin{equation} \label{i64}
k_{\theta_{0}} ( x , y ) = k ( 0 , 0 ) + o_{\theta_{0} \to 0} ( 1 ) = \frac{e^{i A / h} q_{\tau}}{2} + o_{\theta_{0} \to 0} ( 1 ) ,
\end{equation}
uniformly for $x , y \in [ - 1 , 1 ]$, $\tau \in [ - C , C ]$ and $h \in ] 0 , 1 ]$. Let $\CL$ be the operator on $L^{2} ( [ - 1 , 1 ] )$ with kernel equal to $1$. We now use that, for any operator $M$ on $L^{2} ( [ - 1 , 1 ] )$ with kernel $m ( x , y )$, we have
\begin{equation*}
\Vert M \Vert \leq \Vert M \Vert_{\rm H S} = \Vert m \Vert_{L^{2} ( [ - 1 , 1 ]^{2} )} .
\end{equation*}
Combining with \eqref{i64}, it gives
\begin{equation}
K_{\theta_{0}} = \frac{e^{i A / h} q_{\tau}}{2} \CL + o_{\theta_{0} \to 0} ( 1 ) = \frac{e^{i A / h} q_{\tau}}{2} \big( \CL + o_{\theta_{0} \to 0} ( 1 ) \big) .
\end{equation}
Since $\CL$ can be written $\CL = \one_{[ - 1 , 1 ]} ( \one_{[ - 1 , 1 ]} , \cdot )$, it is a self-adjoint bounded rank one operator. Moreover, its non-zero eigenvalue is $2$ and $\one_{[ - 1 , 1 ]}$ is an associated eigenfunction. From the last equation and the perturbation theory, the spectrum of $K_{\theta_{0}}$ can be decomposed as
\begin{equation} \label{i66}
\sigma ( K_{\theta_{0}} ) = \{ \widetilde{\mu}_{\tau , h , \theta_{0}}\} \cup \widetilde{R}_{\tau , h , \theta_{0}} ,
\end{equation}
where $\widetilde{\mu}_{\tau , h , \theta_{0}}$ is a simple eigenvalue satisfying
\begin{equation*}
\widetilde{\mu}_{\tau , h , \theta_{0}} = e^{i A / h} q_{\tau} + o_{\theta_{0} \to 0} ( 1 ) ,
\end{equation*}
and $\widetilde{R}_{\tau , h , \theta_{0}} \subset B ( 0 , o_{\theta_{0} \to 0} ( 1 ) )$. Note that the $o_{\theta_{0} \to 0} ( 1 )$'s appearing in the previous expressions are uniform with respect to $\tau \in [ - C , C ]$ and $h \in ] 0 , 1 ]$. Since the eigenvalues of $\widehat{\CT} ( \tau , h )$ are those of $\theta_{0} K_{\theta_{0}}$, \eqref{i66} implies \eqref{i65}.
\end{proof}

We now demonstrate some properties of the accumulation curves near $E_{0}$ under the stronger assumptions of Remark \ref{i60}.

\begin{proof}[Proof of Remark \ref{i60}]
As in \eqref{i73}, we have
\begin{equation} \label{i77}
\widehat{\CT} ( 0 , h ) ( \omega , \widetilde{\omega} ) = - i e^{i A / h} \frac{\lambda}{2 \pi} \CM_{0} e^{- i \nu \frac{\pi}{2}} ( i \lambda \omega \cdot \widetilde{\omega} )^{- 1} .
\end{equation}
Then, in order to prove the remark, it is enough to show that the operator $A$ defined on $L^{2} ( X )$, with kernel $a ( x , y ) = \cos ( x - y )^{- 1}$, has an infinity of non-zero eigenvalues. Here, $X = \CH^{+ \infty}_{\rm tang}$ is seen as a real subset of length less than $\pi / 2$. Since $A$ is a compact and self-adjoint operator, it is equivalent to show that its range, $\im A$, is of infinite dimension.

Assume that $\im A$ is a finite dimensional vector space. In particular, $\im A$ is closed. Let
\begin{equation*}
Y = \big\{ x \in X ; \ \mes_{\S^{n - 1}} ( X \cap x + [ - \varepsilon , \varepsilon ] ) \neq 0 \text{ for all } \varepsilon > 0 \big\} ,
\end{equation*}
be the subset of charged points of $X$. Since $\mes_{\S^{n - 1}} ( X ) > 0$, $Y$ is infinite. For $y_{0} \in Y$ and $\varepsilon > 0$, we have
\begin{equation*}
A ( \one_{y_{0} + [ - \varepsilon , \varepsilon ]} ) = \big( \cos ( x - y_{0} )^{-1} + o_{\varepsilon \to 0} ( 1 ) \big) \mes_{\S^{n - 1}} ( X \cap y_{0} + [ - \varepsilon , \varepsilon ] ) .
\end{equation*}
Since $\im A$ is closed, we deduce that $\cos ( x - y_{0} )^{-1} \in \im A$ for all $y_{0} \in Y$. We now claim that $( \cos ( x - y_{0} )^{-1} )_{y_{0} \in Y}$ is a free family in $L^{2} ( X )$. If it was not the case, there would exist $y_{1} , \ldots , y_{K} \in Y$ and $a_{1} , \ldots , a_{K} \in \C \setminus \{ 0 \}$ such that
\begin{equation} \label{i74}
\sum_{k = 1}^{K} a_{k} \cos ( x - y_{k} )^{-1} = 0 ,
\end{equation}
in $L^{2} ( X )$. In particular, \eqref{i74} holds punctually almost everywhere on $X$. Using the meromorphy of $\cos ( x - y_{k} )^{-1}$ and $\mes_{\S^{n - 1}} ( X ) > 0$, \eqref{i74} holds true on all $\C$. On the other hand, since $Y$ is a subset of length less than $\pi / 2$, the singularities of $\cos ( x - y_{k} )^{-1}$ and $\cos ( x - y_{\ell} )^{-1}$ are different for $k \neq \ell$. It implies that \eqref{i74} is not possible on $\C$. Thus, the family $( \cos ( x - y_{0} )^{-1} )_{y_{0} \in Y}$ is free, $\im A$ can not be a finite dimensional vector space and the remark follows.
\end{proof}

\begin{proof}[Proof of Remark \ref{i76}]
By assumption, the kernel of $\widehat{\CT} ( P ) ( 0 , h )$ (resp. $\widehat{\CT} ( Q ) ( 0 , h )$) is given by \eqref{i77} on $B_{1} : = \CH^{+ \infty}_{\rm tang} ( P )$ (resp. $B_{2} : = \CH^{+ \infty}_{\rm tang} ( Q )$). Then, as in the proof of Remark \ref{i60}, it is enough to show that
\begin{equation} \label{i78}
\spr ( A_{1} ) \leq \spr ( A_{2} ) ,
\end{equation}
where $A_{\bullet}$ is defined on $L^{2} ( B_{\bullet})$ with kernel $a ( x , y ) = \cos ( x - y )^{- 1}$.

For $u \in L^{2} ( B_{1} )$, let $\widetilde{u} \in L^{2} ( B_{2} )$ denote its extension by $0$. Since $A_{1}$ and $A_{2}$ have same kernel on $B_{1}$, we get $A_{2} \widetilde{u} ( x ) = A_{1} u ( x )$ for all $x \in B_{1}$. It yields
\begin{equation*}
\Vert A_{1} u \Vert_{L^{2} ( B_{1} )} \leq \Vert A_{2} \widetilde{u} \Vert_{L^{2} ( B_{2} )} \leq \Vert A_{2} \Vert \Vert \widetilde{u} \Vert_{L^{2} ( B_{2} )} = \Vert A_{2} \Vert \Vert u \Vert_{L^{2} ( B_{1} )} .
\end{equation*}
We deduce $\Vert A_{1} \Vert \leq \Vert A_{2} \Vert$. Since $A_{\bullet}$ is self-adjoint, $\spr ( A_{\bullet} ) = \Vert A_{\bullet} \Vert$. Then, \eqref{i78} is verified and the remark follows.
\end{proof}

Lastly, we prove the stability phenomenon of Remark \ref{i89} $iii)$.

\begin{proof}[Proof of \eqref{i90}]
We follow the proof of Proposition \ref{i38} and set
\begin{equation*}
\varepsilon = \mes_{\S^{n - 1}} \big( \CH^{+ \infty}_{\rm tang} ( P ) \Delta \CH^{+ \infty}_{\rm tang} ( Q ) \big) .
\end{equation*}
From Proposition \ref{i53} and Theorem \ref{i55} and since the action is the same for $P$ and $Q$, we only have to show that, in the domain $\{ z \in \C ; \ \vert z \vert \geq e^{- 2 C / \lambda} \}$, we have
\begin{equation} \label{i91}
\dist \big( \spe \big( \widetilde{\CT} ( P ) ( \tau ) \big) , \spe \big( \widetilde{\CT} ( Q ) ( \tau ) \big) \big) = o_{\varepsilon \to 0} ( 1 ) ,
\end{equation}
uniformly for $Q$ and $\tau \in [ - C , C ]$. Here, $\widetilde{\CT} ( \tau ) = e^{- i A / h} \CT ( \tau , h )$ as in \eqref{i81}. Let $\overline{\CT} ( P , \tau )$ be the operator defined on $L^{2} ( \S^{n - 1} )$ with kernel
\begin{equation*}
\overline{t}_{P , \tau} ( \omega , \widetilde{\omega} ) = \left\{ \begin{aligned}
&\Gamma \Big( \frac{n}{2} - i \frac{\tau}{\lambda} \Big) \Big( \frac{\lambda}{2 \pi} \Big)^{\frac{n}{2}} \CM_{0} e^{- i ( \nu \frac{\pi}{2} + n \frac{\pi}{4} - T \tau)} \big( i \lambda \omega \cdot \widetilde{\omega} \big)^{- \frac{n}{2} + i \frac{\tau}{\lambda}} && \text{ if } \omega , \widetilde{\omega} \in \CH^{- \infty}_{\rm tang} ( P ) ,  \\
&0 && \text{ otherwise}.
\end{aligned} \right.
\end{equation*}
The operator $\overline{\CT} ( Q , \tau )$ is defined the same way. In particular, the non-zero eigenvalues  of $\overline{\CT} ( \tau )$ are those of $\widetilde{\CT} ( \tau )$. By hypothesis, $\overline{t}_{P , \tau} ( \omega , \widetilde{\omega} )$ and $\overline{t}_{Q , \tau} ( \omega , \widetilde{\omega} )$ are uniformly bounded and coincide for $\omega , \widetilde{\omega} \in \CH^{- \infty}_{\rm tang} ( P ) \cap \CH^{- \infty}_{\rm tang} ( Q )$. Then, using $\Vert K \Vert \leq \Vert K \Vert_{\rm H S} = \Vert k \Vert_{L^{2}}$ for any operator $K$ with kernel $k$, we deduce
\begin{equation}
\overline{\CT} ( Q , \tau ) = \overline{\CT} ( P , \tau ) + o_{\varepsilon \to 0} ( 1 ) .
\end{equation}
uniformly for $Q$ and $\tau \in [ - C , C ]$. This estimate implies \eqref{i91} by the perturbation theory for finite systems of eigenvalues (see Chapter IV.5 of Kato \cite{Ka76_01}).
\end{proof}

\section{Proof of the main results of Section \ref{s42}} \label{s46}

This part is devoted to the proof of Proposition \ref{j81} and Theorem \ref{j79}, the proof of the other results of Section \ref{s42} being postponed to Section \ref{s44}. We follow the general strategy explained in Section \ref{s36} and carried out in the previous sections. We begin with some notions specific to the graph structure.

\Subsection{Notations and tools} \label{s50}

\begin{figure}
\begin{center}
\begin{picture}(0,0)%
\includegraphics{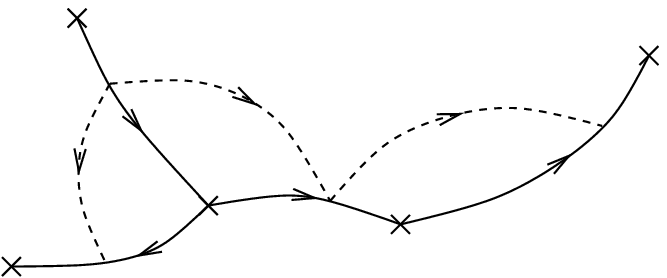}%
\end{picture}%
\setlength{\unitlength}{1184sp}%
\begingroup\makeatletter\ifx\SetFigFont\undefined%
\gdef\SetFigFont#1#2#3#4#5{%
  \reset@font\fontsize{#1}{#2pt}%
  \fontfamily{#3}\fontseries{#4}\fontshape{#5}%
  \selectfont}%
\fi\endgroup%
\begin{picture}(10566,4341)(-1082,-4894)
\put(3901,-4261){\makebox(0,0)[b]{\smash{{\SetFigFont{9}{10.8}{\rmdefault}{\mddefault}{\updefault}$\widetilde{e} = c^{-}$}}}}
\put(8326,-3586){\makebox(0,0)[b]{\smash{{\SetFigFont{9}{10.8}{\rmdefault}{\mddefault}{\updefault}$e = c^{+}$}}}}
\put(6076,-1861){\makebox(0,0)[b]{\smash{{\SetFigFont{9}{10.8}{\rmdefault}{\mddefault}{\updefault}$c = ( e , \widetilde{e} )$}}}}
\end{picture}%
\end{center}
\caption{The graph $\SG = ( \SV , \SE )$ in plain curves and the adjoint graph $\SG^{*} = ( \SE , \SC )$ in dashed lines.} \label{f57}
\end{figure}

In Section \ref{s45}, we have seen the trapped set $K ( E_{0} )$ as the graph $\SG = ( \SV , \SE )$. But, the pertinent quantum quantities live on the edges $\SE$ and are transported to other edges (as the function $u_{-}^{k}$ of Section \ref{s72}, for instance). Thus, we rather use the {\it adjoint graph} $\SG^{*} = ( \SE , \SC )$ of $\SG$. The vertex of $\SG^{*}$ are the edges of $\SG$ and the edges of $\SG^{*}$ are the {\it ways} $c \in \SC$, that is the pairs $( e , \widetilde{e} ) \in \SE^{2}$ with $e^{-} = \widetilde{e}^{+}$. For a way $c = ( e , \widetilde{e} )$, we define the {\it starting edge} $c^{-} = \widetilde{e}$, the {\it ending edge} $c^{+} = e$ and the {\it intermediate vertex} $v ( c ) = e^{-} = \widetilde{e}^{+}$ (see Figure \ref{f57}). A finite sequence of ways $p = ( c_{1} , \ldots , c_{J} )$ is a path in the adjoint graph $\SG^{*}$ when $c_{j}^{+} = c_{j + 1}^{-}$ for all $1 \leq j \leq J - 1$, and a cycle when moreover $c_{J}^{+} = c_{1}^{-}$. As usual, a cycle $( c_{1} , \ldots , c_{J} )$ is identified with the cycle $( c_{2} , \ldots , c_{J} , c_{1} )$. Note that
\begin{equation} \label{m84}
\left\{ \begin{aligned}
&\text{cycles of } \SG^{*} && \longrightarrow && \text{cycles of } \SG  \\
&( c_{1} , \ldots , c_{J} ) && && ( c_{1}^{-} , \ldots , c_{J}^{-} )
\end{aligned} \right.
\end{equation}
provides a natural bijection between the cycles of $\SG^{*}$ and $\SG$. We also import the notions of primitive cycle and minimal cycle to $\SG^{*}$. In the present section, we will only consider paths and cycles in the adjoint graph $\SG^{*}$ and not in $\SG$.

For a path $p = ( c_{1} , \ldots , c_{J} )$ of the adjoint graph $\SG^{*}$, we define its {\it starting edge} $p^{-} = c_{1}^{-}$, its {\it ending edge} $p^{+} = c_{J}^{+}$, its {\it length} $\ell ( p ) = J$ and its {\it strength} by
\begin{equation} \label{m13}
G ( p ) = \alpha ( p ) - D_{0} \beta ( p ) ,
\end{equation}
with
\begin{equation*}
\alpha ( p ) = \sum_{1 \leq j \leq J} \alpha_{v ( c_{j} )} \qquad \text{and} \qquad \beta ( p ) = \sum_{c \in p} \beta_{v ( c )} .
\end{equation*}
The previous definitions of $\alpha , \beta$ on $\SG^{*}$ are coherent with that on $\SG$ (see \eqref{m12}) in the sense that they are preserved by the bijection \eqref{m84} for cycles. For any cycle $p$, we have $G ( p ) \geq 0$ and $G ( p ) = 0$ if and only if $p$ is minimal (see \eqref{m1}).

We explain how to decompose a path into cycles. Let $p$ be a path of the adjoint graph $\SG^{*}$. If $\widetilde{p}$ is a cycle inside $p$, we write $p = \widetilde{p} \oplus ( p \setminus \widetilde{p} )$ where $p \setminus \widetilde{p}$ is the path $p$ without (a copy of) the cycle $\widetilde{p}$. Moreover, if there is a cycle inside $p$, then there is a primitive cycle inside $p$. Then, iterating this procedure, we see that any path $p$ can be decomposed as
\begin{equation} \label{k3}
p = p_{\rm resi} \oplus p_{1} \oplus \cdots \oplus p_{K} ,
\end{equation}
where $p_{1} , \ldots , p_{K}$ are primitive cycles, and the residual path $p_{\rm resi}$ has no cycle and same starting and ending edges as $p$. This decomposition is not necessarily unique. From the definition of the strength \eqref{m13}, we deduce
\begin{equation} \label{k4}
G ( p ) = G ( p_{\rm resi} ) + G ( p_{1} ) + \cdots + G ( p_{K} ) .
\end{equation}
Moreover, since $p_{\rm resi}$ has no cycle, $\ell ( p_{\rm resi} ) < \card \SE$. Then, there exists a constant $G_{0} \in \R$ independent of $p_{\rm resi}$ such that
\begin{equation} \label{k2}
G ( p_{\rm resi} ) \geq - G_{0} ,
\end{equation}
since the number of paths without cycle is finite. We now show that the infimum in \eqref{m1} is attained. Let $p$ be a cycle of the adjoint graph $\SG^{*}$. From \eqref{k3}, we can write
\begin{equation*}
p = p_{1} \oplus \cdots \oplus p_{K} ,
\end{equation*}
where $p_{1} , \ldots , p_{K}$ are primitive cycles. There is no residual path here since $p$ is a cycle. Let $k_{0} \in \{ 1 , \ldots , K \}$ be such that
\begin{equation*}
D ( p_{k_{0}} ) = \min_{k \in \{ 1 , \ldots , K \}} D ( p_{k} ) .
\end{equation*}
Then, we have
\begin{equation}
D ( p ) = \frac{\sum \alpha ( p_{k} )}{\sum \beta ( p_{k} )} \geq \frac{\sum D ( p_{k_{0}} ) \beta ( p_{k} )}{\sum \beta ( p_{k} )} = D ( p_{k_{0}} ) .
\end{equation}
Thus, it is enough to consider primitive cycles in \eqref{m1}. Since there is a finite number of such cycles, the infimum is attained.

Since the vertices in $\SV$ are not identical, we introduce weights that respect the dynamic on $\SG^{*}$. So, for $e \in \SE$, we define
\begin{equation} \label{m14}
N_{e} = \min_{p \text{ path} , \ p^{+} = e} G ( p ) .
\end{equation}
Here we allow $p = \emptyset$ with the convention $G ( \emptyset ) = 0$. In the previous formula, the infimum is attained. Indeed, let $p$ be a path in $\SG^{*}$ with $p^{+} = e$. Using a decomposition \eqref{k3} of $p$, the formula \eqref{k4} and $G ( \widetilde{p} ) \geq 0$ for all cycle $\widetilde{p}$, we obtain
\begin{equation*}
G ( p ) = G ( p_{\rm resi} ) + G ( p_{1} ) + \cdots + G ( p_{K} ) \geq G ( p_{\rm resi} ) ,
\end{equation*}
with $p_{\rm resi}^{+} = e$. Thus, it is enough to minimize in \eqref{m14} over the finite set of the paths without cycle. It implies that the infimum is actually a minimum. Moreover, since $\SE$ is finite, there exists $N_{\infty} \geq 0$ such that
\begin{equation} \label{m52}
- N_{\infty} \leq N_{e} \leq 0 ,
\end{equation}
for all $e \in \SE$. The $N_{e}$ are ``coherent'' in the following sense.

\begin{lemma}\sl \label{m15}
For all $e \in \SE$, we have
\begin{equation*}
N_{e} = \min_{\fract{\widetilde{e} \in \SE , \ p \text{ path}}{p^{-} = \widetilde{e} , \ p^{+} = e}} \big( G ( p ) + N_{\widetilde{e}} \big) .
\end{equation*}
\end{lemma}

\begin{proof}
First, consider $\widetilde{e} \in \SE$ and a path $p$ with $p^{-} = \widetilde{e}$ and $p^{+} = e$. From \eqref{m14}, there exists a path $\widetilde{p}$ with $\widetilde{p}^{+} = \widetilde{e}$ and $N_{\widetilde{e}} = G ( \widetilde{p} )$. Then, $\widetilde{p} \cup p$ is a path with ending edge $e$. Using again \eqref{m14} and the additivity of $G$, we get
\begin{equation*}
G ( p ) + N_{\widetilde{e}} = G ( p ) +  G ( \widetilde{p} ) = G ( \widetilde{p} \cup p ) \geq N_{e} .
\end{equation*}
Taking the minimum of the left hand side, we obtain
\begin{equation} \label{m16}
\min_{\fract{\widetilde{e} \in \SE , \ p \text{ path}}{p^{-} = \widetilde{e} , \ p^{+} = e}} \big( G ( p ) + N_{\widetilde{e}} \big) \geq N_{e} .
\end{equation}
On the other hand, taking the path $p = \emptyset$ from $e$ to $e$ leads to
\begin{equation} \label{m17}
\min_{\fract{\widetilde{e} \in \SE , \ p \text{ path}}{p^{-} = \widetilde{e} , \ p^{+} = e}} \big( G ( p ) + N_{\widetilde{e}} \big) \leq G ( \emptyset ) + N_{e} = N_{e} .
\end{equation}
The lemma follows from \eqref{m16} and \eqref{m17}.
\end{proof}

Only half of Lemma \ref{m15} (more precisely \eqref{m16}) is used in the sequel. Remark that $N_{e} = 0$ for all $e \in \SE$ under the assumptions of Section \ref{s61}. Thus, this construction is useful only when the vertices are different.

We also define a quantity which measures the distance between two edges on the adjoint graph $\SG^{*}$. For $e ,\widetilde{e} \in \SE$, we set
\begin{equation} \label{m58}
N_{e \leftarrow \widetilde{e}} = \min_{\fract{p \text{ path}}{p^{-} = \widetilde{e} , \ p^{+} = e}} G ( p ) ,
\end{equation}
with the convention that $N_{e \leftarrow \widetilde{e}} = + \infty$ if there is no path from $\widetilde{e}$ to $e$. Here again, we allow the path $p = \emptyset$ from $e$ to $e$ for all $e \in \SE$. As for \eqref{m14}, it is enough to minimize on the paths without cycle. Then, this infimum is attained and
\begin{equation} \label{m65}
N_{e \leftarrow e} = 0 ,
\end{equation}
for all $e\in \SE$ (the unique path without cycle being $\emptyset$ in this case). Moreover, following the proof of Lemma \ref{m15}, we obtain the result below showing that these weights are coherent with the graph structure.

\begin{lemma}\sl \label{m59}
For all $e , \widetilde{e} \in \SE$, we have
\begin{equation*}
N_{e \leftarrow \widetilde{e}} = \min_{\fract{\widehat{e} \in \SE , \ p \text{ path}}{p^{-} = \widehat{e} , \ p^{+} = e}} \big( G ( p ) + N_{\widehat{e} \leftarrow \widetilde{e}} \big) .
\end{equation*}
\end{lemma}

\Subsection{Study of the classical system} \label{s47}

We prove here the assertions of Section \ref{s45} concerning the classical system, and we estimate the resolvent of $\SQ ( z , h )$ away from the pseudo-resonances. This part is somehow similar to Section \ref{s13}. We first recall that the determinant of a $N \times N$ matrix $M = ( M_{j , k} )_{j , k}$ can be computed using cycles. More precisely, we have
\begin{equation} \label{m18}
\det ( M ) = \sum_{C} \prod_{p \in C} ( - 1 )^{\ell ( p ) - 1} \prod_{c \in p} M_{c} ,
\end{equation}
where the sum is taken over all the sets of primitive cycles in $\{ 1 , \ldots , N \}$ whose supports are two by two disjoint and fill $\{ 1 , \ldots , N \}$. In this expression, we use the same definition of cycles (i.e. finite sequence of consecutive ways) than at the beginning of Section \ref{s50}. For $c = ( j , k )$, $M_{c}$ is a shortcut for $M_{j , k}$. Using this identity, we can give the

\begin{proof}[Proof of Remark \ref{m4}]
We compute $F ( z , h ) = \det ( 1 - \SQ ( z , h ) )$ applying \eqref{m18}. Consider a primitive cycle $p$ of $\{ 1 , \ldots , \card \SE \}$. Then, $\prod_{c \in p} ( 1 - \SQ )_{c} \neq 0$ only if $p$ is a cycle of the adjoint graph $\SG^{*}$ or if $p = ( ( e , e ) )$ for some $e \in \SE$. In the first case and if $p \neq ( ( e , e ) )$, one can write
\begin{align}
\prod_{c \in p} ( 1 - \SQ ( z , h ) )_{c} &= ( - 1 )^{\ell ( p )} \prod_{c \in p} \SQ_{c} ( z , h )  \nonumber \\
&= h^{\sum_{c \in p} \alpha_{v ( c )} - i \beta_{v ( c )} \sigma} ( - 1 )^{\ell ( p )} \prod_{c \in p} \CQ_{c} ( z , h )    \nonumber \\
&= h^{G ( p )} e^{i \beta ( p ) Z} h^{- i \beta ( p ) \tau} ( - 1 )^{\ell ( p )} \prod_{c \in p} \CQ_{c} ( z , h ) .  \label{m21}
\end{align}
Here, we have used \eqref{m3}, \eqref{m19}, \eqref{m13} and
\begin{equation*}
S_{v} ( z , h ) / \lambda_{1}^{v} - 1 / 2 = \alpha_{v} - i \beta_{v} \sigma .
\end{equation*}
In the second case, we have
\begin{align}
\prod_{c \in p} ( 1 - \SQ ( z , h ) )_{c} &= 1 + ( - 1 )^{\ell ( p )} \prod_{c \in p} \SQ_{c} ( z , h )  \nonumber \\
&=1 + h^{G ( p )} e^{i \beta ( p ) Z} h^{- i \beta ( p ) \tau} ( - 1 )^{\ell ( p )} \prod_{c \in p} \CQ_{c} ( z , h ) .   \label{m22}
\end{align}
Combining \eqref{m18} with \eqref{m21} and \eqref{m22}, we obtain \eqref{m5} with
\begin{equation*}
\SF = \Big\{ ( \alpha , \beta ) \in \R^{2} ; \ \alpha = \sum_{p \in C^{\prime}} G ( p ) \text{ and } \beta = \sum_{p \in C^{\prime}} \beta ( p ) \Big\} ,
\end{equation*}
where $C^{\prime}$ goes through the sets of primitive cycles in $\SG^{*}$ whose supports are two by two disjoint. Let $( \alpha , \beta ) \in \SF$. The discussion below \eqref{m13} implies $\alpha , \beta \geq 0$. Moreover, $\alpha = 0$ if and only if all the primitive cycles of the corresponding $C^{\prime}$ are minimal. Besides, \eqref{m23} follows from \eqref{m20} and the previous equations. Lastly, the coefficient corresponding to $\alpha = \beta = 0$ in \eqref{m5} is simply $1$ from \eqref{m18}. In other words, $F_{0 , 0} ( \tau , z  , h ) = 1$ for all $\tau , z  , h$.

In order to deal with $F^{\rm p r i n}$, we define $\SC^{\rm p r i n}$ as the set of ways $c = ( e , \widetilde{e} ) \in \SC$ appearing in a minimal cycle in $\SG^{*} = ( \SE , \SC )$. We claim that
\begin{equation} \label{m25}
\text{All the cycles in } \SG^{*}_{\rm p r i n}  : = ( \SE , \SC^{\rm p r i n} ) \text{ are minimal.}
\end{equation}
If it is not the case, there exists a cycle $p = ( c_{1} , \ldots , c_{J} )$ in $( \SE , \SC^{\rm p r i n} )$ with $G ( p ) > 0$. Since $c_{j} \in \SC^{\rm p r i n}$, there exists a minimal cycle $p_{j}$ containing $c_{j}$. We can write $p_{j} = \{ c_{j} \} \cup p_{j}^{c}$ where $p_{j}^{c}$ is a path with starting edge $c_{j}^{+}$ and ending edge $c_{j}^{-}$. Hence, $\widetilde{p} = p_{J}^{c} \cup \cdots \cup p_{1}^{c}$ is a cycle. The additivity of $G$ gives
\begin{equation*}
0 = G ( p_{1} ) + \cdots + G ( p_{J} ) = G ( p ) + G ( \widetilde{p} ) .
\end{equation*}
Thus, $G ( \widetilde{p} ) = - G ( p ) < 0$ which is impossible. This shows \eqref{m25}.

Following the same strategy, one can show that
\begin{equation} \label{m83}
( e , \widetilde{e} ) \in \SC^{\rm p r i n} \text{ if and only if } e^{-} = \widetilde{e}^{+} \text{ and } \widetilde{e} , e \text{ belong to minimal cycles of } \SG .
\end{equation}
The direct sense follows from the definition of $\SC^{\rm p r i n}$ and \eqref{m84}. Conversely, if $e$ (resp. $\widetilde{e}$) belongs to a minimal cycle of the graph $\SG$, there exists a minimal cycle of the form $\gamma = \{ e \} \cup \gamma^{c}$ (resp. $\widetilde{\gamma} = \{ \widetilde{e} \} \cup \widetilde{\gamma}^{c}$) in $\SG$. If moreover $e^{-} = \widetilde{e}^{+}$, then $\widehat{\gamma} = \{ \widetilde{e} \} \cup \{ e \} \cup \gamma^{c} \cup \widetilde{\gamma}^{c}$ is a cycle in $\SG$ and
\begin{equation*}
G ( \widehat{\gamma} ) = G ( \gamma ) + G ( \widetilde{\gamma} ) = 0 ,
\end{equation*}
by additivity of $G$. Thus, $\widehat{\gamma}$ is minimal. Eventually, the inverse of the application \eqref{m84} provides a minimal cycle in $\SG^{*}$ containing the way $( e , \widetilde{e} )$. This implies $( e , \widetilde{e} ) \in \SC^{\rm p r i n}$ and proves \eqref{m83}. Using \eqref{m3} and \eqref{m83}, the matrix $\SQ^{\rm p r i n}$ can be written
\begin{equation*}
\SQ^{\rm p r i n}_{e , \widetilde{e}} ( z , h ) = \left\{
\begin{aligned}
&h^{S_{v} ( z , h ) / \lambda_{1}^{v} - 1 / 2} \CQ_{e , \widetilde{e}} ( z , h ) &&\text{ if } ( e , \widetilde{e} ) \in \SC^{\rm p r i n} , \\
&0 &&\text{ if } ( e , \widetilde{e} ) \notin \SC^{\rm p r i n} .
\end{aligned} \right.
\end{equation*}

We now compute $F^{\rm p r i n} ( z , h ) = \det ( 1 - \SQ^{\rm p r i n} ( z , h ) )$ using \eqref{m18} and compare with the corresponding expression for $F ( z , h )$. Let $p$ be a primitive cycle in $\SG^{*}$ with $p \neq ( ( e , e ) )$ for all $e \in \SE$. If $p$ is minimal (i.e. $G ( p ) = 0$), all its ways are in $\SC^{\rm p r i n}$ and
\begin{equation} \label{m66}
\prod_{c \in p} ( 1 - \SQ^{\rm p r i n}( z , h ) )_{c} = \eqref{m21} .
\end{equation}
Whereas if $p$ is not minimal (i.e. $G ( p ) \neq 0$), at least one of its ways is not in $\SC^{\rm p r i n}$ from \eqref{m25} and
\begin{equation} \label{m67}
\prod_{c \in p} ( 1 - \SQ^{\rm p r i n} ( z , h ) )_{c} = 0 .
\end{equation}
When $p = ( ( e , e ) )$, we have similar equations. From \eqref{m18}, it implies that $F^{\rm p r i n} ( z , h )$ satisfies also \eqref{m5} except that only the terms associated to $( \alpha , \beta ) \in \SF$ with $\alpha = 0$ appear. This is exactly \eqref{m6}.
\end{proof}

We obtain some estimates on $f_{\tau} ( \cdot , h )$ and study its zeros. The results proved here are already known when $\tau , h$ are fixed (see Theorem 3 of Langer \cite{La31_01}). We begin with the proof of \eqref{m9} and define
\begin{equation*}
\Gamma_{0} : = ( \Gamma_{0} ( h ) - E_{0} ) h^{- 1} ,
\end{equation*}
which is independent of $h$. From Remark \ref{m4}, we have $F_{0 , 0} ( \tau , z  , h ) = 1$ for all $\tau , z  , h$. Let
\begin{equation*}
\beta_{\rm min} = \min_{\SB \setminus \{ 0 \}} \beta > 0 .
\end{equation*}
For all $\beta \in \SB$, the function $\widetilde{F}_{0 , \beta} ( \kappa , \rho , \tau - i D_{0} )$ is continuous and then bounded on the compact set $( \S^{1} )^{\card \SV + \card \SE} \times [ - C , C ] \setminus ( \Gamma_{0} + i D_{0} + B ( 0 , \delta ) )$. Summing up, \eqref{m6} gives
\begin{equation} \label{m68}
f_{\tau} ( Z , h ) = 1 + \CO ( e^{- \beta_{\rm min} \im Z} ) ,
\end{equation}
uniformly for $\tau \in [ - C , C ] \setminus ( \Gamma_{0} + i D_{0} + B ( 0 , \delta ) )$, $h \in ] 0 , 1 ]$ and $\im Z > 0$. Thus, there exists $N > 0$ such that $f_{\tau} ( \cdot , h )$ does not vanish in $\R + i [ N , +\infty [$ and
\begin{equation} \label{m32}
\big\vert f_{\tau} ( Z , h ) \big\vert^{- 1} \leq 2 ,
\end{equation}
for all $\tau \in [ - C , C ] \setminus ( \Gamma_{0} + i D_{0} + B ( 0 , \delta ) )$, $h \in ] 0 , 1 ]$ and $Z \in \R + i [ N , +\infty [$. This implies \eqref{m9}. Working the same way, we can also prove that there exists (a new) $N > 0$ such that
\begin{equation} \label{m29}
\vert F ( z , h ) \vert^{- 1} \leq 2 \qquad \text{and} \qquad \big\vert {F^{\rm p r i n}} ( z , h ) \big\vert^{- 1} \leq 2 ,
\end{equation}
for all $h \in ] 0 , 1 ]$ and
\begin{equation} \label{m35}
z \in [ - C h , C h ] + i \Big[ - D_{0} h + N \frac{h}{\vert \ln h \vert} , h \Big] \setminus \big( \Gamma_{0} ( h ) + B ( 0 , \delta h ) \big) .
\end{equation}
Equations \eqref{m32} and \eqref{m29} will be used to exclude the presence of (pseudo-)resonances near the real axis.

Let us now prove \eqref{m8}. From \eqref{m27}, we can write
\begin{equation} \label{m31}
f_{\tau} ( Z , h ) = \widetilde{f} ( Z , \lambda ) ,
\end{equation}
where $\widetilde{f}$ is independent of $h$ and holomorphic in $( Z , \lambda )$ for $\lambda = ( \kappa , \rho , \tau )$ near the set
\begin{equation} \label{m30}
( \S^{1} )^{\card \SV + \card \SE} \times [ - C , C ] \setminus \big( \Gamma_{0} + i D_{0} + B ( 0 , \delta ) \big) .
\end{equation}
For $\lambda_{0} \in \eqref{m30}$ fixed, the holomorphic function $Z \mapsto \widetilde{f} ( Z , \lambda_{0} )$ has a finite number of zeros in $[ - C , C ] + i [ - C , N ]$. Using the continuity of $\widetilde{f}$ with respect to $\lambda$ and the Rouch\'e theorem, we deduce
\begin{align*}
\card \big\{ Z  \in [ - C , C ] + i [ &- C , N ] ; \ \widetilde{f} ( Z , \lambda ) = 0 \big\}  \\
&\leq \card \big\{ Z  \in [ - C , C ] + i [ - C , N ] ; \ \widetilde{f} ( Z , \lambda_{0} ) = 0 \big\} ,
\end{align*}
for $\lambda$ in some neighborhood of $\lambda_{0}$. By compactness of \eqref{m30}, there exists $M > 0$ such that
\begin{equation*}
\card \big\{ Z  \in [ - C , C ] + i [ - C , N ] ; \ \widetilde{f} ( Z , \lambda ) = 0 \big\} \leq M ,
\end{equation*}
for all $\lambda \in \eqref{m30}$. Combining with \eqref{m31}, it yields \eqref{m8}.

Let
\begin{equation} \label{m34}
\Lambda ( \tau , h ) = \{ Z  \in \C ; \ f_{\tau} ( Z , h ) = 0 \} ,
\end{equation}
denote the set of zeros of $f_{\tau} ( \cdot , h )$. We have the following estimate outside $\Lambda ( \tau , h )$ which is equivalent to Lemma \ref{d14}.

\begin{lemma}\sl \label{m28}
Let $C , \delta , \alpha > 0$ and $\CK$ be a compact of $\C$. Then, there exists $M > 0$ such that
\begin{equation*}
\big\vert f_{\tau} ( Z , h ) \big\vert^{- 1} \leq M ,
\end{equation*}
for all $Z \in \CK$, $\tau \in [ - C , C ] \setminus ( \Gamma_{0} + i D_{0} + B ( 0 , \delta ) )$ and $h \in ] 0 , 1]$ with $\dist ( Z , \Lambda ( \tau , h ) ) \geq \alpha$.
\end{lemma}

\begin{proof}
Mimicking \eqref{m34} and using the notations of \eqref{m31}, we define
\begin{equation*}
\widetilde{\Lambda} ( \lambda ) = \{ Z  \in \C ; \ \widetilde{f} ( Z , \lambda ) = 0 \} ,
\end{equation*}
and
\begin{equation}
\big\{ ( Z , \lambda ) \in \CK \times \eqref{m30} ; \ \dist ( Z , \widetilde{\Lambda} ( \lambda ) ) \geq \alpha \big\} . \label{m33}
\end{equation}
Since $\widetilde{f} ( Z , \lambda )$ is continuous with respect to $\lambda$ and does not vanish identically, the set of its zeros $\widetilde{\Lambda} ( \lambda )$ depends also continuously on $\lambda$. This implies that \eqref{m33} is a compact set. Then, the continuous function $( Z , \lambda) \mapsto \widetilde{f} ( Z , \lambda )^{- 1}$ is bounded on this compact set. The lemma follows from \eqref{m31}.
\end{proof}

\begin{proof}[Proof of Proposition \ref{j81}]
Thanks to \eqref{m29}, there is no pseudo-resonances in \eqref{m35}. The same way, the zeros of $f_{\tau} ( \cdot , h )$ belong to $\R + i [ - C , N ]$ from \eqref{m9}. Thus, it is enough to work in the set
\begin{equation} \label{m36}
E_{0} + \tau h - i D_{0} h + \big( [ - C , C ] + i [ - C , N ] \big) \frac{h}{\vert \ln h \vert} \setminus \big( \Gamma_{0} ( h ) + B ( 0 , \delta h ) \big) ,
\end{equation}
instead of \eqref{m10} and to consider $Z \in [ - C , C ] + i [ - C , N ]$.

We now compare $F ( z , h )$ and $f_{\tau} ( Z , h )$. For $z \in \eqref{m36}$, we have $\sigma = \tau - i D_{0} + \CO ( \vert \ln h \vert^{- 1} )$. Using the analyticity of $\widetilde{F}_{0 , \beta}$, we get
\begin{equation*}
\widetilde{F}_{0 , \beta} ( \kappa , \rho , \sigma ) = \widetilde{F}_{0 , \beta} ( \kappa , \rho , \tau - i D_{0} ) + \CO \big( \vert \ln h \vert^{- 1} \big) .
\end{equation*}
Then, \eqref{m6} and \eqref{m27} imply
\begin{equation*}
F^{\rm p r i n} ( z , h ) = f_{\tau} ( Z , h ) + \CO \big( \vert \ln h \vert^{- 1} \big) ,
\end{equation*}
uniformly for $z \in \eqref{m36}$, $\tau \in [ - C , C ]$ with $Z$ defined by \eqref{m19}. On the other hand, Remark \ref{m4} gives $F ( z , h ) = F^{\rm p r i n} ( z , h ) + \CO ( h^{\zeta} )$ in \eqref{m36}, for some $\zeta > 0$. Thus, we have
\begin{equation} \label{m37}
F ( z , h ) = f_{\tau} ( Z , h ) + \CO \big( \vert \ln h \vert^{- 1} \big) ,
\end{equation}
uniformly for $z \in \eqref{m36}$, $\tau \in [ - C , C ]$ with $Z$ defined by \eqref{m19}.

Let $\alpha > 0$ be small enough and $Z \in [ - C , C ] + i [ - C , N ]$ with $\dist ( Z , \Lambda ( \tau , h ) ) \geq \alpha$. Combining Lemma \ref{m28} and \eqref{m37}, we get
\begin{equation} \label{m41}
\vert F ( z , h ) \vert \geq \frac{1}{M} + \CO \Big( \frac{1}{\vert \ln h \vert} \Big) \geq \frac{1}{2 M} ,
\end{equation}
for $h$ small enough. Thus, $z$ is not a pseudo-resonance. Taking the contrapositive, we have just proved that, for every pseudo-resonance $z \in \eqref{m10}$, there exists a zero $Z$ of $f_{\tau} ( \cdot , h )$ with
\begin{equation} \label{m38}
\Big\vert z - \Big( E_{0} + \tau h - i D_{0} h + Z \frac{h}{\vert \ln h \vert} \Big) \Big\vert \leq \alpha \frac{h}{\vert \ln h \vert} .
\end{equation}

On the other hand, let $Z \in \Lambda ( \tau , h ) \cap [ - C , C ] + i [ - C , N ]$. Using \eqref{m8} (with $C$ replaced by $C + 1$), there exists $\alpha < r < 2 \alpha N$ such that $\dist ( \widetilde{Z} , \Lambda ( \tau , h ) ) \geq \alpha$ for all $\widetilde{Z} \in \partial B ( Z , r )$. Then Lemma \ref{m28} shows that
\begin{equation*}
\big\vert f_{\tau} ( \widetilde{Z} , h ) \big\vert \geq M^{- 1} ,
\end{equation*}
for all $\widetilde{Z} \in \partial B ( Z , r )$. Together with \eqref{m37} and the Rouch\'e theorem, we deduce that $F ( z , h )$ has a zero in 
\begin{equation} \label{m39}
E_{0} + \tau h - i D_{0} h + B ( Z , r ) \frac{h}{\vert \ln h \vert} .
\end{equation}
Proposition \ref{j81} follows from \eqref{m38} and \eqref{m39} since $\alpha > 0$ can be taken arbitrarily small.
\end{proof}

Adapting Lemma \ref{m28} to the function $F ( z , h )$, we obtain

\begin{lemma}\sl \label{m40}
Let $C , \delta , \alpha > 0$. Then, there exists $M > 0$ such that
\begin{equation*}
\vert F ( z , h ) \vert^{- 1} \leq M ,
\end{equation*}
for all $z \in \eqref{m10}$ with $\dist ( z , \res_{0} ( P ) ) \geq \alpha h \vert \ln h \vert^{- 1}$.
\end{lemma}

\begin{proof}
In the region \eqref{m35}, this estimate follows from  \eqref{m29}. Assume now that $z \in \eqref{m36}$. Let $\tau \in [ - C , C ] \setminus ( \Gamma_{0} + i D_{0} + B ( 0 , \delta ) )$ and $Z \in [ - C , C ] + i [ - C , N ]$ such that \eqref{m19} holds true. Since $\dist ( z , \res_{0} ( P ) ) \geq \alpha h \vert \ln h \vert^{- 1}$, Proposition \ref{j81} yields $\dist ( Z , \Lambda ( \tau , h ) ) \geq \alpha / 2$ for $h$ small enough. Then, \eqref{m41} provides the required inequality.
\end{proof}

This estimate controls the determinant of $1 - \SQ ( z , h )$, but not the inverse of this operator. Indeed, the coefficients of its adjugate matrix may be polynomial in $h^{- 1}$ and then unbounded. The situation was different in Lemma \ref{d13} where the adjugate matrix is always bounded. Let $M = ( M_{j , k} )_{j , k}$ be a general $N \times N$ matrix. The $( j , k )$ entry of the adjugate matrix, denoted $\adj ( M )$, is defined by
\begin{equation} \label{m45}
\adj ( M ) _{j , k} = \det \big( \widehat{M}_{k , j} \big) ,
\end{equation}
where $\widehat{M}_{k , j}$ is the matrix $M$ except that the entries on the line $k$ and on the column $j$ are $0$ but the coefficient $( k , j )$ is $1$. As in Remark \ref{m4}, the adjugate matrix of $1 - \SQ$ satisfies

\begin{lemma}\sl \label{m44}
 For all $e , \widetilde{e} \in \SE$, there exists a finite set $\SF_{e , \widetilde{e}} \subset [ 0 , + \infty [^{2}$ such that
\begin{equation} \label{m62}
\adj \big( 1 - \SQ ( z , h ) \big)_{e , \widetilde{e}} = h^{N_{e \leftarrow \widetilde{e}}} \sum_{( \alpha , \beta ) \in \SF_{e , \widetilde{e}}} h^{\alpha} e^{i \beta Z} F_{\alpha , \beta}^{e , \widetilde{e}} ( \tau , z , h ) ,
\end{equation}
with the notations of \eqref{m19}. If there is no path from $\widetilde{e}$ to $e$ in the adjoint graph $\SG^{*}$ and $e \neq \widetilde{e}$ (i.e. if $N_{e \leftarrow \widetilde{e}} = + \infty$), this coefficient is null. The functions $F_{\alpha , \beta}^{e , \widetilde{e}}$ can be written
\begin{equation} \label{m63}
F_{\alpha , \beta}^{e , \widetilde{e}} ( \tau , z , h ) = \widetilde{F}_{\alpha , \beta}^{e , \widetilde{e}} ( \kappa , \rho , \sigma ) ,
\end{equation}
where $\widetilde{F}_{\alpha , \beta}$ is independent of $h$ and holomorphic in $\kappa , \rho , \sigma$ for $\sigma$ outside $\Gamma_{0}$. In particular, we have
\begin{equation} \label{m42}
\adj \big( 1 - \SQ ( z , h ) \big)_{e , \widetilde{e}} = \CO ( h^{N_{e \leftarrow \widetilde{e}}} ) ,
\end{equation}
uniformly for $z \in \eqref{m10}$.

If there is a path from $\widetilde{e}$ to $e$ in the graph $\SG^{*}_{\rm p r i n}$ or if $e = \widetilde{e}$, we have
\begin{equation} \label{m60}
\adj \big( 1 - \SQ^{\rm p r i n} ( z , h ) \big)_{e , \widetilde{e}} = h^{N_{e \leftarrow \widetilde{e}}} \sum_{( 0 , \beta ) \in \SF_{e , \widetilde{e}}} e^{i \beta Z} F_{0 , \beta}^{e , \widetilde{e}} ( \tau , z , h ) .
\end{equation}
In particular, there exists $\zeta > 0$ independent of $e , \widetilde{e}$ such that
\begin{equation} \label{m61}
\adj \big( 1 - \SQ^{\rm p r i n} ( z , h ) \big)_{e , \widetilde{e}} = \adj \big( 1 - \SQ ( z , h ) \big)_{e , \widetilde{e}} + \CO ( h^{N_{e \leftarrow \widetilde{e}} + \zeta} ) ,
\end{equation}
uniformly for $z \in \eqref{m10}$.
\end{lemma}

\begin{proof}
We begin with the first part of the lemma. From \eqref{m45}, we have to compute the cofactor
\begin{equation*}
\det \big( \widehat{( 1 - \SQ )}_{\widetilde{e} , e} \big) .
\end{equation*}
For that, we use \eqref{m18}. Thus, let $C$ be a set of primitive cycles in $\{ 1 , \ldots , \card \SE \}$ whose supports are two by two disjoint and fill $\{ 1 , \ldots , \card \SE \}$. If the quantity
\begin{equation} \label{m46}
\prod_{p \in C} ( - 1 )^{\ell ( p ) - 1} \prod_{c \in p} \big( \widehat{( 1 - \SQ )}_{\widetilde{e} , e} \big)_{c} ,
\end{equation}
is not zero (i.e. if the contribution of $C$ in $\adj ( 1 - \SQ )_{e , \widetilde{e}}$ is not null), there exists one (and only one) cycle $p_{0}$ of $C$ containing the way $c_{0} = ( \widetilde{e} , e )$. Moreover, all the other ways appearing in the cycles of $C$ must be in the adjoint graph $\SG^{*}$ or must be of the form $( \widehat{e} , \widehat{e} )$ for some $\widehat{e} \in \SE$. Denote $\widetilde{p}_{0} = p_{0} \setminus \{ c_{0} \}$ which is necessarily a path of the adjoint graph with starting edge $\widetilde{e}$ and ending edge $e$ since $p_{0} = \{ ( \widetilde{e} , e ) \} \cup \widetilde{p}_{0}$ is a cycle. Again, we use the convention $\widetilde{p}_{0} = \emptyset$ as a formal path from $\widetilde{e}$ to $e$ if $e = \widetilde{e}$. In particular, if there is no path from $\widetilde{e}$ to $e$ in the adjoint graph $\SG^{*}$ and $e \neq \widetilde{e}$, $\eqref{m46} = 0$ for all $C$ and eventually
\begin{equation}  \label{m49}
\adj \big( 1 - \SQ ( z , h ) \big)_{e , \widetilde{e}} = 0 .
\end{equation}

Assume now that there exists at least one path from $\widetilde{e}$ to $e$ in the adjoint graph $\SG^{*}$ or that $e = \widetilde{e}$. Let $p \neq p_{0}$ be a cycle of $C$. We have
\begin{equation*}
\prod_{c \in p} \big( \widehat{( 1 - \SQ )}_{\widetilde{e} , e} \big)_{c} = \prod_{c \in p} ( 1 - \SQ )_{c} .
\end{equation*}
Using \eqref{m21} or \eqref{m22} to compute the right hand side, we obtain
\begin{equation} \label{m47}
\prod_{c \in p} \big( \widehat{( 1 - \SQ )}_{\widetilde{e} , e} \big)_{c} = h^{G ( p )} e^{i \beta ( p ) Z} h^{- i \beta ( p ) \tau} ( - 1 )^{\ell ( p )} \prod_{c \in p} \CQ_{c} ( z , h ) + \left\{ \begin{aligned}
&0 &&\text{ if } p \neq ( ( \widehat{e} , \widehat{e} ) ) ,  \\
&1 &&\text{ if } p = ( ( \widehat{e} , \widehat{e} ) ) ,
\end{aligned} \right.
\end{equation}
with $G ( p ) \geq 0$ since $p$ is a cycle. On the other hand, for $p = p_{0}$, we have as in \eqref{m21}
\begin{align}
\prod_{c \in p_{0}} \big( \widehat{( 1 - \SQ )}_{\widetilde{e} , e} \big)_{c} &= ( - 1 )^{\ell ( \widetilde{p}_{0} )} \prod_{c \in \widetilde{p}_{0}} \SQ_{c}  \nonumber \\
&= h^{G ( \widetilde{p}_{0} )} e^{i \beta ( \widetilde{p}_{0} ) Z} h^{- i \beta ( \widetilde{p}_{0} ) \tau} ( - 1 )^{\ell ( \widetilde{p}_{0} )} \prod_{c \in \widetilde{p}_{0}} \CQ_{c} ( z , h )  \nonumber \\
&= h^{N_{e \leftarrow \widetilde{e}}} h^{G ( \widetilde{p}_{0} ) - N_{e \leftarrow \widetilde{e}}} e^{i \beta ( \widetilde{p}_{0} ) Z} h^{- i \beta ( \widetilde{p}_{0} ) \tau} ( - 1 )^{\ell ( \widetilde{p}_{0} )} \prod_{c \in \widetilde{p}_{0}} \CQ_{c} ( z , h ) ,  \label{m48}
\end{align}
with $G ( \widetilde{p}_{0} ) - N_{e \leftarrow \widetilde{e}} \geq 0$ from \eqref{m58}. Combining \eqref{m46} with \eqref{m47} and \eqref{m48}, we get \eqref{m62}. The particular form for $F_{\alpha , \beta}^{e , \widetilde{e}}$ stated in \eqref{m63} follows from \eqref{m20}. Lastly, \eqref{m42} is a consequence of \eqref{m62} and $\im Z \geq - C$ for $z \in \eqref{m10}$.

We now prove the second part of the lemma and assume that there is a path from $\widetilde{e}$ to $e$ in the graph $\SG^{*}_{\rm p r i n}$ or that $e = \widetilde{e}$. We first show that
\begin{equation} \label{m64}
\begin{gathered}
\text{The paths } p \text{ from } \widetilde{e} \text{ to } e \text{ in } \SG^{*} \text{ such that }  \\
G ( p ) = N_{e \leftarrow \widetilde{e}} \text{ are the paths from } \widetilde{e} \text{ to } e \text{ in } \SG^{*}_{\rm p r i n} .
\end{gathered}
\end{equation}
In the case $e = \widetilde{e}$, the paths $p$ from $e$ to $e$ (i.e. the cycles passing through $e$) in $\SG^{*}$ such that $G ( p ) = N_{e \leftarrow e} = 0$ (see \eqref{m65}) are by definition the minimal cycles. Thus, \eqref{m64} follows from the definition of $\SG^{*}_{\rm p r i n}$ in this case. Note also  that the formal path $p = \emptyset$ from $e$ to $e$ satisfies $G ( p ) = N_{e \leftarrow e}$ but is not in $\SG^{*}$ and $\SG^{*}_{\rm p r i n}$. Suppose now that $e \neq \widetilde{e}$ and consider $p$ (resp. $\widehat{p}$) a path in $\SG^{*}_{\rm p r i n}$ (resp. $\SG^{*}$) from $\widetilde{e}$ to $e$. Such paths exist by assumption. Working as below \eqref{m25}, we can construct a path $\widetilde{p}$ in $\SG^{*}_{\rm p r i n}$ from $e$ to $\widetilde{e}$. Then, $p \cup \widetilde{p}$ and $\widehat{p} \cup \widetilde{p}$ are two cycles in $\SG^{*}$ and $p \cup \widetilde{p}$ is minimal thanks to \eqref{m25}. The additivity of $G$ gives
\begin{equation*}
G ( p ) = G ( p \cup \widetilde{p} ) - G ( \widetilde{p} ) \leq G ( \widehat{p} \cup \widetilde{p} ) - G ( \widetilde{p} ) = G ( \widehat{p} ) .
\end{equation*}
This implies that $G ( p ) = N_{e \leftarrow \widetilde{e}}$ for all paths $p$ in $\SG^{*}_{\rm p r i n}$ from $\widetilde{e}$ to $e$. On the other hand, let $\widehat{p}$ be a path from $\widetilde{e}$ to $e$ in $\SG^{*}$ such that $G ( \widehat{p} ) = N_{e \leftarrow \widetilde{e}}$. Consider a path $p$ from $\widetilde{e}$ to $e$ in $\SG^{*}_{\rm p r i n}$ and let $\widetilde{p}$ be as before. Then, we have
\begin{equation*}
G ( \widehat{p} \cup \widetilde{p} ) = G ( \widehat{p} ) + G ( \widetilde{p} ) = N_{e \leftarrow \widetilde{e}} + G ( \widetilde{p} ) = G ( p ) + G ( \widetilde{p} ) = G ( p \cup \widetilde{p} ) = 0 .
\end{equation*}
Thus, $\widehat{p} \cup \widetilde{p}$ is minimal and $\widehat{p}$ is in $\SG^{*}_{\rm p r i n}$ by definition. This finishes the proof of \eqref{m64}.

In order to show \eqref{m60}, we use again \eqref{m18} and \eqref{m45} which allows us to write $\adj ( 1 - \SQ^{\rm p r i n} )_{e , \widetilde{e}}$ as a finite sum of terms like \eqref{m46} mutatis mutandis. We then compare the resulting expression with \eqref{m62}. As explained in \eqref{m66}--\eqref{m67}, only the contributions from the minimal cycles (i.e. the terms of order $h^{0}$) of \eqref{m47} appear in the expansion of $\adj ( 1 - \SQ^{\rm p r i n} )_{e , \widetilde{e}}$. The same way, \eqref{m64} shows that the terms \eqref{m48} appearing in the expansion of $\adj ( 1 - \SQ^{\rm p r i n} )_{e , \widetilde{e}}$ are those such that $G ( \widetilde{p}_{0} ) - N_{e \leftarrow \widetilde{e}} = 0$. Then, $\adj ( 1 - \SQ^{\rm p r i n} )_{e , \widetilde{e}}$ is given by the sum of the terms in the right hand side of \eqref{m62} corresponding to $\alpha = 0$. This proves \eqref{m60}. Finally, \eqref{m61} follows from \eqref{m62}, \eqref{m60} and $\im Z \geq - C$ for $z \in \eqref{m10}$.
\end{proof}

\Subsection{Resonance free zone and resolvent estimate}  \label{s48}

Since the present setting has some similarities with that of Section \ref{s61}, we use some constructions and intermediate results of Section \ref{s12}. First we prove a resolvent estimate away from pseudo-resonances. More precisely,

\begin{proposition}\sl \label{j86}
Let $C , \delta >0$. For $h$ small enough, $P$ has no resonance in the domain
\begin{equation} \label{j87}
\begin{aligned}
E_{0} + [ - C h , C  h ] & + i \Big[ - D_{0} h - C \frac{h}{\vert \ln h \vert} , h \Big]    \\
&\setminus \big( \Gamma (h) + B ( 0 , \delta h ) \big) \bigcup \Big( \res_{0} (P) + B \Big( 0 , \delta \frac{h}{\vert \ln h \vert} \Big) \Big) .
\end{aligned}
\end{equation}
Moreover, there exists $M > 0$ such that
\begin{equation*}
\big\Vert( P_{\theta} -z )^{-1} \big\Vert \lesssim h^{- M} ,
\end{equation*}
for $h$ small enough and $z \in \eqref{j87}$.
\end{proposition}

In order to do this, we use the general contradiction argument of Section \ref{s3}. So, it is enough to show that every $u$ satisfying $( P_{\theta} - z ) u = \CO ( h^{\infty} )$ and $\Vert u \Vert = 1$ vanishes microlocally near each point of $K ( E_{0} )$.

For each edge $e \in \SE$, we choose a point $\rho_{\pm}^{e} = ( x_{\pm}^{e} , \xi_{\pm}^{e} ) \in e$ close to $( e^{\mp} , 0 )$. We define $u_{\pm}^{e}$ as the microlocal restriction of $u$ to a neighborhood of $\rho_{\pm}^{e}$. Contrary to Section \ref{s4} and Section \ref{s72}, we do not specify the neighborhoods (i.e. $U_{\pm}^{\bullet} , V_{\pm}^{\bullet}$) where these functions are localized since the present geometric setting is similar mutatis mutandis to that of these sections.

The first step is to express $u_{-}^{e}$ in terms of the $u_{-}^{\bullet}$'s using the propagation of singularities through the vertex $e^{-}$ and along the trajectory $e$. We use an abstract formulation similar to \eqref{a84} and \eqref{j84}. From Lemma \ref{a33}, $u$ is a solution of the microlocal Cauchy problem
\begin{equation} \label{j88}
\left\{ \begin{aligned}
&(P -z) u = 0 &&\text{microlocally near } e^{-} ,   \\
&u = u_{-}^{\widetilde{e}} &&\text{microlocally near } \rho_{-}^{\widetilde{e}} \text{ if } e^{-} = \widetilde{e}^{+} ,   \\
&u = 0 &&\text{microlocally near each point of } \Lambda_{-}^{e^{-}} \setminus K ( E_{0} ) ,
\end{aligned} \right.
\end{equation}
with $\Vert u \Vert \leq 1$. As in \eqref{d62}, the unique solution of \eqref{j88} can be seen as the sum over $\widetilde{e}$ (with $e^{-} = \widetilde{e}^{+} = : v$) of the solutions of the microlocal Cauchy problems with non-zero initial data only on the curve $\widetilde{e}$. Moreover, \ref{h18} guaranties that $g_{+} ( \rho_{+}^{e} ) \cdot g_{-} ( \rho_{-}^{\widetilde{e}} ) \neq 0$. Then, applying Theorem \ref{a32}, we obtain
\begin{equation} \label{j85}
u_{+}^{e} = \sum_{e^{-} = \widetilde{e}^{+}} \SJ_{e \leftarrow \widetilde{e}} u_{-}^{\widetilde{e}} \text{ microlocally near } \rho_{+}^{e} ,
\end{equation}
where the operator $\SJ_{e \leftarrow \widetilde{e}}$ is given by \eqref{m76} for the vertex $e^{-}$ and restricted to microlocal neighborhoods of $\rho_{+}^{e}$ and $\rho_{-}^{\widetilde{e}}$ (see \eqref{a98} for a similar situation). This means that
\begin{equation*}
\SJ_{e \leftarrow \widetilde{e}} u_{-}^{\widetilde{e}} ( x ) = h^{S_{v} ( z , h ) / \lambda_{1}^{v} - n / 2} \int_{\R^{n}} e^{i ( \varphi_{+}^{v} ( x ) - \varphi_{-}^{v} ( y ) ) / h} \widetilde{d}_{e \leftarrow \widetilde{e}} ( x, y , z , h ) u_{-}^{\widetilde{e}} ( y ) \, d y ,
\end{equation*}
where $\varphi_{\pm}^{v}$ are generating function near $v$ of the manifold $\Lambda_{\pm}^{v}$ with $\varphi_{\pm}^{v} ( v ) = 0$, and $\widetilde{d}_{e \leftarrow \widetilde{e}}$ is the symbol $\widetilde{d}$ of \eqref{m76} restricted to neighborhoods of $x = x_{+}^{e}$ and $y = x_{-}^{\widetilde{e}}$. Using $\Vert u \Vert \leq 1$ and $z \in \eqref{j87}$, we deduce
\begin{equation} \label{j89}
u_{+}^{e} \in \CI ( \Lambda_{+}^{e^{-}} , h^{- N} ) ,
\end{equation}
for some $N \in \R$ independent of $e$. We now compute $u_{-}^{e}$ from $u_{+}^{e}$ using the propagation of singularities along $e$. Like in \eqref{a37}, $u$ satisfies
\begin{equation}
\left\{ \begin{aligned}
&( P -z ) u = 0 &&\text{microlocally near } e ,  \\
&u = u_{+}^{e} &&\text{microlocally near } \rho_{+}^{e} ,
\end{aligned} \right.
\end{equation}
and $\Vert u \Vert \leq 1$. Microlocally near $\rho_{-}^{e}$, the solution of this microlocal Cauchy problem is given by $e^{- i t_{e} ( P - z ) / h} u_{+}^{e}$, where $t_{e} \in \R$ is such that $\exp ( t_{e} H_{p} ) ( \rho_{+}^{e} ) = \rho_{-}^{e}$. Thus, as explained below \eqref{b5}, we have
\begin{equation} \label{j91}
u_{-}^{e} = \SM_{e} u_{+}^{e} \text{ microlocally near } \rho_{-}^{e} ,
\end{equation}
where $\SM_{e}$ is a Fourier integral operator of order $0$ with canonical transformation $\exp ( t_{e} H_{p} )$. In particular, combining with \eqref{j89}, we obtain
\begin{equation} \label{k90}
u_{-}^{e} \in \CI ( \Lambda_{+}^{e^{-}} , h^{- N} ) .
\end{equation}
Contrary to the previous sections, we do not specify the part (i.e. $\Lambda^{0}_{+} , \Lambda^{1}_{+}$) of the Lagrangian manifold $\Lambda_{+}^{e^{-}}$ appearing in \eqref{j89} and \eqref{j90}. Using the coherent weights $N_{e} \geq 0$ defined in \eqref{m14} and satisfying \eqref{m52}, we have a fortiori
\begin{equation} \label{j90}
u_{-}^{e} \in \CI ( \Lambda_{+}^{e^{-}} , h^{N_{e} - N} ) .
\end{equation}

Combining \eqref{j85} and \eqref{j91}, we can write
\begin{equation} \label{j92}
u_{-}^{e} = \sum_{e^{-} = \widetilde{e}^{+}} \SP_{e , \widetilde{e}} u_{-}^{\widetilde{e}} + S ( h^{\infty} ) ,
\end{equation}
with $\SP_{e , \widetilde{e}} = \SM_{e} \SJ_{e \leftarrow \widetilde{e}}$. We consider this equation at the level of the symbols of the $u_{-}^{e}$'s. As in \eqref{d48}, we define the symbol $a_{-}^{e} ( x , h ) \in S ( h^{N_{e} - N} )$ for $e \in \SE$ by
\begin{equation} \label{k15}
u_{-}^{e} ( x , h ) = e^{- i A_{e} / h} e^{i \frac{z - E_{0}}{h} t_{-}^{e}} \frac{\CM_{e}^{-}}{\CD_{e} ( t_{-}^{e} )} a_{-}^{e} ( x , h ) e^{i \varphi_{+}^{1 , e} ( x ) / h} .
\end{equation}
Here, $A_{e} , \CM_{e}^{-}$ are defined in Section \ref{s42}, the time $t_{-}^{e} \in \R$ is such that $e ( t_{-}^{e} ) = \rho_{-}^{e}$,
\begin{equation*}
\CD_{e} (t) = \sqrt{\Big\vert \det \frac{\partial x_{e} ( s , y )}{\partial ( s , y )} \vert_{s = t , \ y = 0} \Big\vert} ,
\end{equation*}
with the notation of \eqref{d7} adapted to the present setting and $\varphi_{+}^{1 , e}$ is the generating function near $x_{-}^{e}$ of the extension of the outgoing stable manifold associated with $e^{-}$ and normalized by
\begin{equation*}
\varphi_{+}^{1 , e} ( x_{-}^{e} ) = \int_{e ( ] - \infty , t_{-}^{e}  ] )} \xi \cdot d x .
\end{equation*}
The symbols $a_{-}^{e}$ satisfy the following identity.

\begin{lemma}\sl \label{k10}
There exist $\zeta > 0$ and symbols $\CP_{e , \widetilde{e}} \in S( h^{\alpha_{v} - D_{0} \beta_{v}} )$ independent of $u$, where $v$ is the intermediate vertex $v ( ( e , \widetilde{e} ) )$, such that
\begin{equation} \label{k13}
a_{-}^{e} ( x , h ) = \sum_{e^{-} = \widetilde{e}^{+}} \CP_{e , \widetilde{e}} ( x , h ) a_{-}^{\widetilde{e}} ( x_{-}^{\widetilde{e}} , h ) + S ( h^{N_{e} - N + \zeta} ) ,
\end{equation}
for $e \in \SE$ and $x$ near $x_{-}^{e}$. Moreover, $\CP_{e , \widetilde{e}} ( x_{-}^{e} , h ) = \SQ_{e , \widetilde{e}} ( z , h )$.
\end{lemma}

\begin{proof}
From \eqref{j92}, we only have to compute the symbol of $\SP_{e , \widetilde{e}} u_{-}^{\widetilde{e}}$, that is the contribution of the way $( e , \widetilde{e} )$. This can be achieved as in Lemma \ref{d41} (see also Lemma \ref{g89}, Lemma \ref{g46}, \ldots) and we omit the details. We only mention that we apply Corollary \ref{d46} instead of Theorem \ref{a32} in \eqref{j85}, that $\CP_{e , \widetilde{e}} \in S( h^{\alpha_{v} - D_{0} \beta_{v}} )$ follows from
\begin{equation*}
\big\vert h^{S_{v} ( z , h ) / \lambda_{1}^{v} - 1 / 2} \big\vert \leq h^{\alpha_{v} - D_{0} \beta_{v}} ,
\end{equation*}
uniformly for $z \in \eqref{j80}$, and that the remainder terms are estimated using
\begin{equation} \label{m53}
h^{\alpha_{v} - D_{0} \beta_{v}} h^{N_{\widetilde{e}} - N + \zeta} = h^{G ( p )} h^{N_{\widetilde{e}} - N + \zeta} \leq h^{N_{e} - N + \zeta} ,
\end{equation}
where $p = ( ( e , \widetilde{e} ) )$ is a path from $\widetilde{e}$ to $e$. The last inequality follows from Lemma \ref{m15}.
\end{proof}

We can now finish the proof of Proposition \ref{j86}. Let $\SA$ be the $\card \SE$-vector with entries $\SA_{e} = a_{-}^{e} ( x_{-}^{e} , h )$. Applying \eqref{k13} with $x = x_{-}^{e}$, we get
\begin{equation*}
\SA = \SQ ( z , h ) \SA + \SR ,
\end{equation*}
where $\SR$ is a $\card \SE$-vector satisfying $\SR_{e} = \CO ( h^{N_{e} - N + \zeta} )$. This identity can be written
\begin{equation*}
\SA = \big( 1 - \SQ ( z , h ) \big)^{- 1} \SR = \frac{1}{\det \big( 1 - \SQ ( z , h ) \big)} \adj \big( 1 - \SQ ( z , h ) \big) \SR .
\end{equation*}
Using Lemma \ref{m40} (to estimate the determinant) and \eqref{m42} (to estimate the adjugate matrix), it yields, for all $e \in \SE$,
\begin{equation} \label{k29}
\SA_{e} = \sum_{\widetilde{e} \in \SE} \CO ( h^{N_{e \leftarrow \widetilde{e}}} ) \SR_{\widetilde{e}} = \sum_{\widetilde{e} \in \SE} \CO ( h^{N_{e \leftarrow \widetilde{e}} + N_{\widetilde{e}} - N + \zeta} ) = \CO ( h^{N_{e} - N + \zeta} ) ,
\end{equation}
uniformly for $z \in \eqref{j87}$. Here, we have used that Lemma \ref{m15} implies
\begin{equation*}
G ( p ) + N_{\widetilde{e}} \geq N_{e} ,
\end{equation*}
for all path $p$ from $\widetilde{e}$ to $e$, and then
\begin{equation*}
N_{e \leftarrow \widetilde{e}} + N_{\widetilde{e}} \geq N_{e} ,
\end{equation*}
from \eqref{m58}. A new application of \eqref{k13} gives
\begin{equation} \label{j99}
a_{-}^{e} \in \sum_{e^{-} = \widetilde{e}^{+}} S ( h^{\alpha_{v} - D_{0} \beta_{v} + N_{\widetilde{e}} - N + \zeta} ) + S ( h^{N_{e} - N + \zeta} ) \subset S ( h^{N_{e} - N + \zeta} ) ,
\end{equation}
from \eqref{m53}. In other words, starting with $a_{-}^{e} \in S ( h^{N_{e} - N} )$ for all $e \in \SE$, we have obtained $a_{-}^{e} \in S ( h^{N_{e} - N + \zeta} )$ for all $e \in \SE$. Then, a bootstrap argument yields
\begin{equation}
u = 0 \text{ microlocally near each point of } e ,
\end{equation}
for all $e \in \SE$. The uniqueness part of Theorem \ref{a32} implies
\begin{equation*}
u = 0 \text{ microlocally near } K ( E_{0} ) ,
\end{equation*}
and Proposition \ref{j86} follows thanks to Section \ref{s3}.

\Subsection{Existence of resonances near the pseudo-resonances}  \label{s49}

In order to prove Theorem \ref{j79}, it remains to show the existence of a resonance near each pseudo-resonance. This is done in the following proposition which is analogous to Proposition \ref{d50}.

\begin{proposition}\sl \label{k14}
Let $C , \delta , r > 0$ and assume that $h$ is small enough. For any pseudo-resonance $z \in \eqref{j80}$, the operator $P$ has at least one resonance in $ B ( z , r h \vert \ln h \vert^{- 1} )$.
\end{proposition}

We follow closely Section \ref{s73} and use the notations of the proof of Proposition \ref{j86}. As usual, this result is proved by a contradiction argument. If Proposition \ref{k14} is not satisfied, there exist $r > 0$, a sequence of $h$ which goes to $0$ and a sequence $z = z ( h ) \in \res_{0} ( P ) \cap \eqref{j80}$ such that
\begin{equation} \label{m54}
P \text{ has no resonance in } B \Big( z , 2 r \frac{h}{\vert \ln h \vert} \Big) .
\end{equation}

As for \eqref{d55}, we use compactness arguments to reduce this problem. As in \eqref{d92}, we set $\sigma ( h ) = ( z ( h ) - E_{0} ) / h$ which belongs to the compact set
\begin{equation*}
[ - C , C ] + i [ - D_{0} - 1 , 1 ] \setminus ( \Gamma_{0} + B ( 0 , \delta ) ) .
\end{equation*}
Moreover, from \eqref{m9} and Proposition \ref{j81}, we have
\begin{equation} \label{m56}
\im z ( h ) \subset - D_{0} h + [ - C , N ] \frac{h}{\vert \ln h \vert} ,
\end{equation}
and then $\im \sigma ( h ) = - D_{0} + o ( 1 )$ as $h \to 0$. Then, up to the extraction of a subsequence, we can assume that
\begin{equation*}
\sigma ( h ) \longrightarrow \tau_{0} - i D_{0} \qquad \text{as } h \to 0 ,
\end{equation*}
for some $\tau_{0} \in [ - C , C ] \setminus ( \Gamma_{0} + i D_{0} + B ( 0 , \delta ) )$. We set
\begin{equation*}
\tau ( h ) = \re \sigma ( h ) = \re \frac{z ( h ) - E_{0}}{h} ,
\end{equation*}
which verifies $\tau ( h ) = \tau_{0} + o ( 1 )$. Moreover, since
\begin{equation*} 
Z ( h ) = \vert \ln h \vert \Big( \frac{z ( h ) - E_{0}}{h} - \tau ( h ) + i D_{0} \Big) \in i [ - C , N ] ,
\end{equation*}
from \eqref{m56}, we can extract a subsequence of $h$ such that
\begin{equation*}
Z ( h ) \longrightarrow Z_{0} \qquad \text{as } h \to 0 ,
\end{equation*}
for some constant $Z_{0} \in [ - C , C ] + i [ - C , N ]$. On the other hand, using that the parameter $\kappa ( h ) = ( h^{- i \beta_{v} \tau ( h )} )_{v \in \SV}$ belongs to the compact set $( \S^{1} )^{\card \SV}$, we can assume that $\kappa ( h ) = \kappa_{0} + o ( 1 )$ with $\kappa_{0} \in ( \S^{1} )^{\card \SV}$ independent of $h$. The same way, since $\rho ( h ) = ( e^{i A_{e} / h} )_{e \in \SE}$ belongs to the compact set $( \S^{1} )^{\card \SE}$, we have, up to extraction of a subsequence, $\rho ( h ) = \rho_{0} + o ( 1 )$ for some $\rho_{0} \in ( \S^{1} )^{\card \SE}$ independent of $h$. For $Z \in \C$, we define
\begin{equation} \label{m69}
F^{\rm f i x} ( Z ) = \sum_{\beta \in \SB} e^{i \beta Z} \widetilde{F}_{0 , \beta} ( \kappa_{0} , \rho_{0} , \tau_{0} - i D_{0} ) ,
\end{equation}
with the notations of Remark \ref{m4}. This function is holomorphic and independent of $h$. Since $\widetilde{F}_{0 , \beta}$ is holomorphic, the previous limits and \eqref{m27} yield
\begin{equation} \label{m70}
F^{\rm f i x} ( Z ) = f_{\tau ( h )} ( Z , h ) + o ( 1 ) ,
\end{equation}
as $h \to 0$, uniformly for $Z$ in a compact set of $\C$. In particular,
\begin{equation*}
F^{\rm f i x} ( Z_{0} ) = f_{\tau ( h )} ( Z ( h ) , h ) + o ( 1 ) .
\end{equation*}
Moreover, we have
\begin{equation*}
f_{\tau ( h )} ( Z ( h ) , h ) = F ( z ( h ) , h ) + \CO \big( \vert \ln h \vert^{- 1} \big) ,
\end{equation*}
from \eqref{m37}. On the other hand, $F ( z ( h ) , h ) = 0$ because $z ( h )$ is a pseudo-resonance. Summing up, $F^{\rm f i x} ( Z_{0} ) = 0 + o ( 1 )$ and then
\begin{equation} \label{m57}
F^{\rm f i x} ( Z_{0} ) = 0 ,
\end{equation}
because this quantity is independent of $h$. Moreover, adapting the proof of \eqref{m68}, one can show that $F^{\rm f i x} ( Z ) \to 1$ as $\im Z \to + \infty$. Thus, $F^{\rm f i x}$ does not vanish identically. Lastly, \eqref{m54} gives
\begin{equation} \label{k16}
P \text{ has no resonance in } E_{0} + \tau ( h ) h - i D_{0} h + Z_{0} \frac{h}{\vert \ln h \vert} + B ( 0 , r ) \frac{h}{\vert \ln h \vert} .
\end{equation}

Let $\SQ^{\rm f i x} ( Z )$ be the $\card \SE \times \card \SE$ matrix with entries
\begin{equation*}
\SQ^{\rm f i x}_{e , \widetilde{e}} ( Z ) = \left\{
\begin{aligned}
&h^{\alpha_{v} - D_{0} \beta_{v}} e^{i \beta_{v} Z} \kappa_{0 , v} \rho_{0 , e} \widetilde{\CQ}_{e , \widetilde{e}} ( \tau_{0} - i D_{0} ) &&\text{ if } e^{-} = \widetilde{e}^{+} = : v  \\
& &&\text{ and } ( \widetilde{e} , e ) \in \SC^{\rm p r i n} , \\
&0 &&\text{ otherwise} ,
\end{aligned} \right.
\end{equation*}
where $\widetilde{\CQ}$ is defined in \eqref{m20}. In some sense, $\SQ^{\rm f i x} ( Z )$ is the matrix $\SQ^{\rm p r i n} ( z , h )$ with $z$ given by \eqref{m19} and with $\kappa , \rho , \sigma$ fixed at $\kappa_{0} , \rho_{0} , \tau_{0} - i D_{0}$ (see \eqref{m3}). In particular, following the proof of \eqref{m6} and using \eqref{m23} and \eqref{m69}, we deduce
\begin{equation*}
F^{\rm f i x} ( Z ) = \det \big( 1 - \SQ^{\rm f i x} ( Z ) \big) .
\end{equation*}
The same way, since each coefficient of the adjugate matrix $\adj ( M )$ is a finite sum of a finite product of coefficients of the matrix $M$ (see \eqref{m18} and \eqref{m45}), the adjugate matrix of $1 - \SQ^{\rm f i x}$ can be written
\begin{equation*}
\adj \big( 1 - \SQ^{\rm f i x} ( Z ) \big)_{e , \widetilde{e}} = h^{N_{e \leftarrow \widetilde{e}}} \sum_{( 0 , \beta ) \in \SF_{e , \widetilde{e}}} e^{i \beta Z} \widetilde{F}_{0 , \beta}^{e , \widetilde{e}} ( \kappa_{0} , \rho_{0} , \tau_{0} - i D_{0} ) ,
\end{equation*}
with $\widetilde{F}_{0 , \beta}^{e , \widetilde{e}}$ given by \eqref{m63}. Since $F^{\rm f i x}$ is not identically zero, the previous equations yield
\begin{equation} \label{k37}
\big( 1 - \SQ^{\rm f i x} \big)^{- 1}_{e , \widetilde{e}} ( Z ) = \frac{1}{\det \big( 1 - \SQ^{\rm f i x} ( Z ) \big)} \adj \big( 1 - \SQ^{\rm f i x} ( Z ) \big)_{e , \widetilde{e}} = h^{N_{e \leftarrow \widetilde{e}}} M_{e , \widetilde{e}} ( Z ) ,
\end{equation}
where the $M_{e , \widetilde{e}}$'s are meromorphic functions. From \eqref{m57}, $1 - \SQ^{\rm f i x} ( Z_{0} )$ is not invertible. Thus, there exists a pair of edges $e_{1} , e_{2}$, such that $M_{e_{2} , e_{1}} ( Z )$ has a pole at $Z = Z_{0}$. We can write
\begin{equation} \label{k38}
M_{e_{2} , e_{1}} ( Z ) = \frac{c_{0}}{( Z - Z_{0} )^{m_{0}}} + \frac{R ( Z )}{( Z - Z_{0} )^{m_{0} - 1}} ,
\end{equation}
where $m_{0} \geq 1$, $c_{0} \neq 0$ and $R ( Z )$ is holomorphic near $Z_{0}$.

For $0 < s < r$, we define the domain
\begin{equation*}
\CD = E_{0} + \tau ( h ) h - i D_{0} h + B ( Z_{0} , s ) \frac{h}{\vert \ln h \vert} ,
\end{equation*}
If $s$ is fixed small enough, $Z_{0}$ is the unique zero of $F^{\rm f i x}$ in $B ( Z_{0} , s )$ and $F^{\rm f i x} ( Z )$ is away from zero for $Z$ near $\partial B ( Z_{0} , s )$. Then, Proposition \ref{j81} and \eqref{m70} imply that, for $h$ small enough, $\partial \CD$ is at distance at least $h \vert \ln h \vert^{- 1}$ from the pseudo-resonances of $P$. Applying Proposition \ref{j86}, we deduce
\begin{equation} \label{m71}
\big\Vert( P_{\theta} - \widetilde{z} )^{-1} \big\Vert \lesssim h^{- M} ,
\end{equation}
uniformly for $\widetilde{z} \in \partial \CD$.

We construct a test function microsupported near the edge $e_{1}$. Let $\widetilde{v}$ be a WKB solution of
\begin{equation}
\left\{ \begin{aligned}
&( P - \widetilde{z} ) \widetilde{v} = 0 &&\text{near } x_{-}^{e_{1}} ,   \\
&\widetilde{v} (x) = e^{i \varphi_{+}^{1 , e_{1}} ( x ) / h} &&\text{on } H \text{ near } x_{-}^{e_{1}} ,
\end{aligned} \right.
\end{equation}
holomorphic in $\widetilde{z} \in B ( E_{0} , R h )$ with $R > 1$ large enough. Here, $H$ is a hyperplane transversal to $\pi_{x} ( e_{1} )$ at $x_{-}^{e_{1}}$. In particular, this Lagrangian distribution can be written
\begin{equation} \label{k34}
\widetilde{v} ( x , h ) = \widetilde{a} ( x , h ) e^{i \varphi_{+}^{1 , e_{1}} ( x ) / h} ,
\end{equation}
for some $\widetilde{a} \in S ( 1 )$ holomorphic in $\widetilde{z}$ such that $\widetilde{a} ( x_{-}^{e_{1}}  , h ) = 1$. Considering cut-off functions $\chi , \psi \in C^{\infty}_{0} ( T^{*} \R^{n} )$ as in \eqref{d60}, we define successively
\begin{equation*}
\widehat{v} ( x , h ) = e^{- i A_{e_{1}} / h} e^{i \frac{\widetilde{z} - E_{0}}{h} t_{-}^{e_{1}}} \frac{\CM_{e_{1}}^{-}}{\CD_{e_{1}} ( t_{-}^{e_{1}} )} \widetilde{v} ( x , h ) ,
\end{equation*}
and
\begin{equation}
v = \Op ( \psi ) \big[ P , \Op ( \chi ) \big] \widehat{v} .
\end{equation}
Thus, the microsupport of $v$ is close to $e_{1} ( ] - \infty , t_{-}^{e_{1}} - \nu ] )$ for some $\nu > 0$ (see Figure \ref{f22}). Eventually, let $u \in L^{2} ( \R^{n} )$ be the solution of
\begin{equation}
( P_{\theta} - \widetilde{z} ) u = v ,
\end{equation}
for $\widetilde{z} \in \partial \CD$. Thanks to \eqref{m71}, $u = ( P_{\theta} - \widetilde{z} )^{- 1 } v$ is a holomorphic function of $\widetilde{z} \in \partial \CD$ and satisfies $\Vert u \Vert \lesssim h^{- M}$ uniformly on $\partial \CD$. Let $u_{\pm}^{e}$ be the microlocal restriction of $u$ to a neighborhood of $\rho_{\pm}^{e}$. Following the proofs of Lemma \ref{d72} and Proposition \ref{j86}, one can show
\begin{equation} \label{j96}
u_{-}^{e} \in \CI ( \Lambda_{+}^{e^{-}} , h^{- N} ) ,
\end{equation}
for some $N \in \R$ independent of $e \in \SE$. Moreover, \eqref{j92} still holds true for $e \neq e_{1}$, whereas this equation is replaced by
\begin{equation} \label{j95}
u_{-}^{e_{1}} = \sum_{e_{1}^{-} = \widetilde{e}^{+}} \SP_{e_{1} , \widetilde{e}} u_{-}^{\widetilde{e}} + \widetilde{v} + S ( h^{\infty} ) ,
\end{equation}
for $e = e_{1}$ because the microsupport of $v$ is ``before'' $\rho_{-}^{e_{1}}$ on $e_{1}$. Since we want to use the weights $N_{e \leftarrow e_{1}}$ constructed in \eqref{m58}, we have to show the following result.

\begin{lemma}\sl \label{k26}
For all $e \in \SE$ such that $N_{e \leftarrow e_{1}} = + \infty$, we have
\begin{equation} \label{k25}
u_{-}^{e} = \CO ( h^{\infty} ) ,
\end{equation}
uniformly for $\widetilde{z} \in \partial \CD$.
\end{lemma}

\begin{proof}
Let
\begin{equation*}
\SE^{\rm r e s i} = \{ e \in \SE ; \ N_{e \leftarrow e_{1}} = + \infty \} ,
\end{equation*}
be the set of edges $e \neq e_{1}$ such that there is no path in $\SG^{*}$ from $e_{1}$ to $e$. We set
\begin{equation*}
u_{-}^{{\rm r e s i} , e} = u_{-}^{e} ,
\end{equation*}
for $e \in \SE^{\rm r e s i}$. We remark that, for all $e \in \SE^{\rm r e s i}$ and $\widetilde{e} \in \SE$ such that $e^{-} = \widetilde{e}^{+}$, we have $\widetilde{e} \in \SE^{\rm r e s i}$. Then, \eqref{j92} gives
\begin{equation} \label{j97}
u_{-}^{{\rm r e s i} , e} = \sum_{e^{-} = \widetilde{e}^{+}} \SP_{e , \widetilde{e}} u_{-}^{{\rm r e s i} , \widetilde{e}} + S ( h^{\infty} ) ,
\end{equation}
for all $e \in \SE^{\rm r e s i}$. Roughly speaking, this means that $\SE^{\rm r e s i}$ is a closed system. From \eqref{m52} and \eqref{j96}, we have $u_{-}^{{\rm r e s i} , e} \in \CI ( \Lambda_{+}^{e^{-}} , h^{N_{e} - N^{\rm r e s i}} )$ for all $e \in \SE^{\rm r e s i}$ with $N^{\rm r e s i} = N$. Let $a_{-}^{{\rm r e s i} , e} \in S ( h^{N_{e} - N^{\rm r e s i}} )$ be the symbol of $u_{-}^{{\rm r e s i} ,e}$ defined as in \eqref{k15}. From \eqref{j97}, they satisfy \eqref{k13} for $e \in \SE^{\rm r e s i}$. Thus, the $\card \SE$-vector $\SA^{\rm r e s i}$ defined by
\begin{equation*}
\SA^{\rm r e s i}_{e} = \left\{ \begin{aligned}
& a_{-}^{{\rm r e s i} , e} ( x_{-}^{e} , h ) &&\text{ if } e \in \SE^{\rm r e s i} , \\
& 0 &&\text{ if } e \in \SE \setminus \SE^{\rm r e s i} ,
\end{aligned} \right.
\end{equation*}
verifies the equation
\begin{equation*}
\SA^{\rm r e s i} = \SQ ( \widetilde{z} , h ) \SA^{\rm r e s i} + \SR^{\rm r e s i} ,
\end{equation*}
and then
\begin{equation} \label{j98}
\SA^{\rm r e s i} = \big( 1 - \SQ ( \widetilde{z} , h ) \big)^{- 1} \SR^{\rm r e s i} ,
\end{equation}
where the $\card \SE$-vector $\SR^{\rm r e s i}$ is of the form
\begin{equation*}
\SR^{\rm r e s i}_{e} = \left\{ \begin{aligned}
& \CO ( h^{N_{e} - N^{\rm r e s i} + \zeta} ) &&\text{ if } e \in \SE^{\rm r e s i} , \\
& \CO ( h^{\widetilde{N}} ) &&\text{ if } e \in \SE \setminus \SE^{\rm r e s i} ,
\end{aligned} \right.
\end{equation*}
for some $\widetilde{N} \in \R$. By definition of $1 - \SQ$ and $\SE^{\rm r e s i}$, the matrix $\SQ$ can be written, separating the edges in $\SE^{\rm r e s i}$ and those in $\SE \setminus \SE^{\rm r e s i}$,
\begin{equation*}
1 - \SQ = \left( \begin{array}{cc}
1 - \SQ^{\rm r e s i} & 0 \\
A & B
\end{array} \right) ,
\end{equation*}
for some matrices $A , B$. When this matrix is invertible, its inverse takes the form
\begin{equation*}
( 1 - \SQ )^{- 1} = \left( \begin{array}{cc}
( 1 - \SQ^{\rm r e s i} )^{- 1} & 0 \\
- B^{- 1} A ( 1 - \SQ^{\rm r e s i} )^{- 1} & B^{- 1}
\end{array} \right) .
\end{equation*}
From \eqref{j98}, it implies that only the entries of $\SR^{\rm r e s i}$ corresponding to edges in $\SE^{\rm r e s i}$ are useful to estimate $\SA^{\rm r e s i}$. Since $\partial \CD$ is at distance at least $h \vert \ln h \vert^{- 1}$ of the pseudo-resonances, we can work as for \eqref{j99} and prove that $\SA^{\rm r e s i}_{e} = \CO ( h^{N_{e} - N^{\rm r e s i} + \zeta} )$ for all $e \in \SE^{\rm r e s i}$. Another application of \eqref{k13} show that
\begin{equation}
a_{-}^{{\rm r e s i} , e} \in S ( h^{N_{e} - N^{\rm r e s i} + \zeta} ) ,
\end{equation}
for all $e \in \SE^{\rm r e s i}$. Then, a bootstrap argument yields \eqref{k25}. Note that the previous estimates hold true uniformly for $\widetilde{z} \in \partial \CD$.
\end{proof}

From \eqref{j96} and Lemma \ref{k26}, there exists a new $N \in \R$ independent of $e \in \SE$ such that
\begin{equation}
u_{-}^{e} \in \CI ( \Lambda_{+}^{e^{-}} , h^{N_{e \leftarrow e_{1}} - N} ) ,
\end{equation}
for all $e \in \SE$. Let $a_{-}^{e} \in S ( h^{N_{e \leftarrow e_{1}} - N} )$ be the symbol of $u_{-}^{e}$ defined as in \eqref{k15}. As in Lemma \ref{d72}, they satisfy

\begin{lemma}\sl \label{k27}
Let $\widetilde{a}_{e_{1}} = \widetilde{a}$ and $\widetilde{a}_{e} = 0$ for $e \neq e_{1}$. Uniformly for $\widetilde{z} \in \partial \CD$, we have
\begin{equation} \label{k28}
a_{-}^{e} ( x , h ) = \sum_{e^{-} = \widetilde{e}^{+}} \CP_{e , \widetilde{e}} ( x , h ) a_{-}^{\widetilde{e}} ( x_{-}^{\widetilde{e}} , h ) + \widetilde{a}_{e} ( x, h ) + S ( h^{N_{e \leftarrow e_{1}} - N + \zeta} ) ,
\end{equation}
where $\zeta > 0$ and the symbols $\CP_{e , \widetilde{e}}$ are the ones of Lemma \ref{k10}.
\end{lemma}

\begin{proof}
This result can be proved as Lemma \ref{k10} and only few changes have to be made. The $u_{-}^{e}$'s are linked by \eqref{j95} instead of \eqref{j92}. Moreover, Lemma \ref{m15} is replaced by Lemma \ref{m59} which guaranties that the weights $N_{e \leftarrow e_{1}}$ are coherent.
\end{proof}

As in Section \ref{s48}, let $\SA$ (resp. $\widetilde{\SA}$) be the $\card \SE$-vector of the $a_{-}^{e} ( x_{-}^{e} , h )$ (resp. $\widetilde{a}_{e} ( x_{-}^{e} , h )$). Taking $x = x_{-}^{e}$ in \eqref{k28}, we get
\begin{equation*}
\SA = \SQ ( \widetilde{z} , h ) \SA + \widetilde{\SA} + \SR ,
\end{equation*}
with $\SR_{e} = \CO ( h^{N_{e \leftarrow e_{1}} - N + \zeta} )$. Hence, it yields
\begin{equation} \label{k32}
\SA = \big( 1 -  \SQ ( \widetilde{z} , h ) \big)^{- 1} \widetilde{\SA} + \big( 1 -  \SQ ( \widetilde{z} , h ) \big)^{- 1} \SR .
\end{equation}
Using that $\partial \CD$ is at distance at least $h \vert \ln h \vert^{- 1}$ of $\res_{0} ( P )$, we obtain as in \eqref{k29}
\begin{align}
\SA_{e} &= \sum_{\widetilde{e} \in \SE} \CO ( h^{N_{e \leftarrow \widetilde{e}}} ) \widetilde{\SA}_{\widetilde{e}} + \sum_{\widetilde{e} \in \SE} \CO ( h^{N_{e \leftarrow \widetilde{e}}} ) \SR_{\widetilde{e}}  \nonumber \\
&= \sum_{\widetilde{e} \in \SE} \CO ( h^{N_{e \leftarrow \widetilde{e}} + N_{\widetilde{e} \leftarrow e_{1}}} ) + \sum_{\widetilde{e} \in \SE} \CO ( h^{N_{e \leftarrow \widetilde{e}} + N_{\widetilde{e} \leftarrow e_{1}} - N + \zeta} ) ,  \label{k30}
\end{align}
since $\widetilde{\SA}_{\widetilde{e}} = \CO ( h^{N_{\widetilde{e} \leftarrow e_{1}}} )$. From Lemma \ref{m59}, we have $G ( p ) + N_{\widetilde{e} \leftarrow e_{1}} \geq N_{e \leftarrow e_{1}}$ for all path $p$ from $\widetilde{e}$ to $e$. Minimizing over $p$ leads to
\begin{equation*}
N_{e \leftarrow \widetilde{e}} + N_{\widetilde{e} \leftarrow e_{1}} \geq N_{e \leftarrow e_{1}} ,
\end{equation*}
from \eqref{m58}. Hence, \eqref{k30} becomes
\begin{equation} \label{k33}
\SA_{e} = \CO ( h^{N_{e \leftarrow e_{1}}} ) + \CO ( h^{N_{e \leftarrow e_{1}} - N + \zeta} ) .
\end{equation}
Applying one more time \eqref{k28} and Lemma \ref{m59}, we deduce
\begin{equation*}
a_{-}^{e} \in S ( h^{N_{e \leftarrow e_{1}}} ) + S ( h^{N_{e \leftarrow e_{1}} - N + \zeta} ) .
\end{equation*}
Then, the usual bootstrap argument gives
\begin{equation} \label{k31}
a_{-}^{e} \in S ( h^{N_{e \leftarrow e_{1}}} ) ,
\end{equation}
uniformly for $\widetilde{z} \in \partial \CD$ (i.e. $N = 0$).

Computing the component on the edge $e_{2}$ in \eqref{k32} and estimating the remainder term as in \eqref{k33}, we get
\begin{align*}
a_{-}^{e_{2}} ( x_{-}^{e_{2}} , h ) &= \big( 1 -  \SQ ( \widetilde{z} , h ) \big)^{- 1}_{e_{2} , e_{1}} \widetilde{a}_{e_{1}} ( x_{-}^{e_{1}} , h ) + \CO ( h^{N_{e_{2} \leftarrow e_{1}} + \zeta} )  \\
&= \big( 1 -  \SQ ( \widetilde{z} , h ) \big)^{- 1}_{e_{2} , e_{1}} + \CO ( h^{N_{e_{2} \leftarrow e_{1}} + \zeta} ) ,
\end{align*}
from \eqref{k34}. From \eqref{m5}, \eqref{m6} and Lemma \ref{m40}, we have
\begin{equation*}
\vert F ( \widetilde{z} , h ) \vert^{- 1} + \big\vert F^{\rm p r i n} ( \widetilde{z} , h ) \big\vert^{- 1} \leq 3 M ,
\end{equation*}
for all $\widetilde{z} \in \partial \CD$. Together with \eqref{m5}, \eqref{m6} and \eqref{m61}, it shows that one can replace $\SQ$ by $\SQ^{\rm p r i n}$ in the previous equation. That is
\begin{equation*}
a_{-}^{e_{2}} ( x_{-}^{e_{2}} , h ) = \big( 1 -  \SQ^{\rm p r i n} ( \widetilde{z} , h ) \big)^{- 1}_{e_{2} , e_{1}} + \CO ( h^{N_{e_{2} \leftarrow e_{1}} + \zeta} ) .
\end{equation*}
The same way, since the quantities $\sigma (  h ) , \tau ( h ) , \kappa ( h ) , \rho ( h )$ converge to $\tau_{0} - i D_{0} , \tau_{0} , \kappa_{0} , \rho_{0}$ as $h \to 0$ (see \eqref{m56}--\eqref{m69}), one can replace $\SQ^{\rm p r i n}$ by $\SQ^{\rm f i x}$. More precisely, the previous equation becomes
\begin{align}
a_{-}^{e_{2}} ( x_{-}^{e_{2}} , h ) &= \big( 1 -  \SQ^{\rm f i x} ( \widetilde{Z} ) \big)^{- 1}_{e_{2} , e_{1}} + o ( h^{N_{e_{2} \leftarrow e_{1}}} )  \nonumber \\
&= h^{N_{e_{2} \leftarrow e_{1}}} M_{e_{2} , e_{1}} ( \widetilde{Z} ) + o ( h^{N_{e_{2} \leftarrow e_{1}}} ) ,  \label{k36}
\end{align}
with the new variable
\begin{equation*}
\widetilde{Z} = \vert \ln h \vert \Big( \frac{\widetilde{z} - E_{0}}{h} - \tau ( h ) + i D_{0} \Big) \in \partial B ( Z_{0} , s ) .
\end{equation*}
The last equality of \eqref{k36} follows from \eqref{k37}. We set
\begin{equation*}
z_{0} ( h ) = E_{0} + \tau ( h )  h - i D_{0} h + Z_{0} \frac{h}{\vert \ln h \vert} .
\end{equation*}
From the choice of $s$, $Z_{0}$ is the unique singularity of $M_{e_{2} , e_{1}}$ in $B ( Z_{0} , s )$. Moreover, the behavior of $M_{e_{2} , e_{1}}$ at $Z_{0}$ is described in \eqref{k38}. Combining with \eqref{k36}, the Cauchy formula and the previous change of variables, we obtain
\begin{align}
\frac{1}{2 i \pi} \int_{\partial \CD} &( \widetilde{z} - z_{0} ( h ) )^{m_{0} - 1} a_{-}^{e_{2}} ( x_{-}^{e_{2}} , h ) \, d \widetilde{z}   \nonumber  \\
&= \frac{1}{2 i \pi} \frac{h^{m_{0} + N_{e_{2} \leftarrow e_{1}}}}{\vert \ln h \vert^{m_{0}}} \int_{\partial B ( Z_{0} , s )} ( \widetilde{Z} - Z_{0} )^{m_{0} - 1} M_{e_{2} , e_{1}} ( \widetilde{Z} ) \, d \widetilde{Z} + o \Big( \frac{h^{m_{0} + N_{e_{2} \leftarrow e_{1}}}}{\vert \ln h \vert^{m_{0}}} \Big)  \nonumber \\
&= \frac{1}{2 i \pi} \frac{h^{m_{0} + N_{e_{2} \leftarrow e_{1}}}}{\vert \ln h \vert^{m_{0}}} \int_{\partial B ( Z_{0} , s )} \frac{c_{0}}{\widetilde{Z} - Z_{0}} \, d \widetilde{Z} + o \Big( \frac{h^{m_{0} + N_{e_{2} \leftarrow e_{1}}}}{\vert \ln h \vert^{m_{0}}} \Big)   \nonumber \\
&= c_{0} \frac{h^{m_{0} + N_{e_{2} \leftarrow e_{1}}}}{\vert \ln h \vert^{m_{0}}} + o \Big( \frac{h^{m_{0} + N_{e_{2} \leftarrow e_{1}}}}{\vert \ln h \vert^{m_{0}}} \Big) .
\end{align}
Thus, the Lagrangian distribution
\begin{equation*}
\frac{1}{2 i \pi} \int_{\partial \CD} ( \widetilde{z} - z_{0} ( h ) )^{m_{0} - 1} u_{-}^{e_{2}} ( \widetilde{z} ) \, d \widetilde{z} \in \CI \Big( \Lambda_{+}^{e_{2}^{-}} , \frac{h^{m_{0} + N_{e_{2} \leftarrow e_{1}}}}{\vert \ln h \vert^{m_{0}}} \Big) ,
\end{equation*}
is elliptic at $x = x_{-}^{e_{2}}$. On the other hand, \eqref{k16} and the choice of $v$ guaranty that $u_{-}^{e_{2}} ( \widetilde{z} )$ is holomorphic in $\CD$. In particular,
\begin{equation*}
\frac{1}{2 i \pi} \int_{\partial \CD} ( \widetilde{z} - z_{0} ( h ) )^{m_{0} - 1} u_{-}^{e_{2}} ( \widetilde{z} ) \, d \widetilde{z} = 0 .
\end{equation*}
It contradicts the ellipticity of this Lagrangian distribution. Then, Proposition \ref{k14} follows.

\section{Proofs of the other results of Section \ref{s42}} \label{s44}

Here, we show the results of Section \ref{s42} which have not been proved in Section \ref{s46}. We first deal with graphs without cycle.

\begin{proof}[Proof of Theorem \ref{j73}]
We follow the general strategy developed in Section \ref{s3}. From Proposition \ref{a16}, it is enough to show that if $u ( x , h ) \in L^{2} ( \R^{n} )$ with $\Vert u \Vert \leq 1$ satisfies $( P_{\theta} - z ) u = \CO ( h^{\infty} )$ for $z$ in the domain $\Omega_{h} = \eqref{j72}$, then it is microlocally $0$ near $K ( E_{0} )$.

Using a contradiction argument, we now prove that
\begin{equation} \label{j74}
u = 0 \text{ microlocally near each point of } e ,
\end{equation}
for all $e \in \SE$. If it does not hold true, there exists at least one edge $e_{1}$ such that $u$ is not $0$ microlocally near each point of $e_{1}$. Let us examine the situation on the incoming manifold $\Lambda_{-}^{e_{1}^{-}}$ associated to $e_{1}^{-}$.  As in Lemma \ref{a33}, $u = 0$ microlocally near each point of $\Lambda_{-}^{e_{1}^{-}} \setminus K ( E_{0} )$. If in addiction $u = 0$ microlocally near each point of each edge $e$ with $e^{+} = e_{1}^{-}$, then the uniqueness part in Theorem \ref{a32} gives $u = 0$ microlocally near $( e_{1}^{-} , 0)$  and then near each point of $e_{1}$. This is impossible. Thus, there should be at least one edge $e_{2}$ with $e_{2}^{+} = e_{1}^{-}$ such that $u$ is not $0$ microlocally near each point of $e_{2}$. Repeating this procedure, we construct an infinite sequence $e_{1} , e_{2} , e_{3} , \ldots$ of successive edges (in the sense $e_{k + 1}^{+} = e_{k}^{-}$ for  all $k \geq 1$) on which $u$ is not microlocally $0$ near each point. Since the number of vertices is finite, there should be a cycle in this infinite sequence and then in the graph $( \SV , \SE )$. This contradicts \ref{h17} and \eqref{j74} is verified.

Let $v \in \SV$. From Lemma \ref{a33} and \eqref{j74}, we get $u = 0$ microlocally near each point of $\Lambda_{-}^{v}$. Then, the uniqueness part in Theorem \ref{a32} implies $u = 0$ microlocally near $( v , 0 )$. Summing up, we have proved $u = 0$ microlocally near $K ( E_{0} )$ and the theorem follows.
\end{proof}

We now show two geometric results specific to the Schr\"{o}dinger operators.

\begin{lemma}\sl \label{k57}
Let $P$ be a Schr\"{o}dinger operator satisfying the assumptions of Theorem \ref{j73}. Then, the following alternative holds:

$i)$ $P$ is non-trapping at energy $E_{0}$,

$ii)$ $E_{0}$ is the global maximum of $V$ attained at a single point $v_{0} \in \R$.
\end{lemma}

\begin{proof}
Let $C : = \{ x \in \R^{n} ; \ V ( x ) \geq E_{0} \}$. We suppose that $K ( E_{0} ) \neq \emptyset$ and that
\begin{equation} \label{k58}
\card C \geq 2 .
\end{equation}
From \ref{h16}, there exists at least one vertex $x_{0} \in \SV$, which is automatically an isolated point of $C$. Let $g : [ t_{-} , t_{+} ] \longrightarrow \R^{n}$ be a geodesic from $x_{0}$ to $C \setminus \{ x_{0} \}$ for the Jacobi metric $( E_{0} - V ( x ) )_{+} d x^{2}$. From Jacobi's theorem (see  Abraham and Marsden \cite[Theorem 3.7.7]{AbMa78_01}), $g$ is the base space projection of a Hamiltonian trajectory $\gamma : ] - \infty , t_{1} [ \longrightarrow p^{- 1} ( E_{0} )$, up to a new time parametrization. We define
\begin{equation*}
x_{1} = g ( t_{+} ) = \lim_{t \to t_{1}} \pi_{x} ( \gamma ( t ) ) \in C .
\end{equation*}

If $t_{1} = + \infty$, we have $x_{1} \in \SV$ and $\gamma$ is a heteroclinic trajectory from $x_{0}$ to $x_{1}$. Using the symmetry of the Hamiltonian trajectories $( x (t) , \xi (t) ) \longmapsto ( x ( - t ) ,  - \xi ( - t ) )$ for the Schr\"{o}dinger operators, we obtain an heteroclinic trajectory $\widetilde{\gamma}$ from $x_{1}$ to $x_{0}$. Then, $( \gamma , \widetilde{\gamma} )$ is a cycle in the graph $( \SV , \SE )$.

Assume now that $t_{1} < + \infty$ and let $\gamma ( t ) = ( x ( t ) , \xi ( t ) )$ be the components of $\gamma$. A direct computation shows that
\begin{equation*}
\widehat{\gamma} ( t ) = \left\{ \begin{aligned}
& \gamma ( t ) &&\text{ for } t \in ] - \infty , t_{1} [ , \\
& ( x_{1} , 0 ) &&\text{ for } t = t_{1} ,  \\
& ( x ( 2 t_{1} - t) , - \xi ( 2 t_{1} - t) ) &&\text{ for } t \in ] t_{1} , + \infty [ ,
\end{aligned} \right.
\end{equation*}
defines a Hamiltonian trajectory because $P$ is a Schr\"{o}dinger operator. Since this curve is a homoclinic trajectory (from $x_{0}$ to $x_{0}$), $( \widehat{\gamma} )$ is a cycle in the graph $( \SV , \SE )$.

In both cases, we have proved that $( \SV , \SE )$ has a cycle. This is in contradiction with \ref{h17}. Then, \eqref{k58} does not hold and the lemma follows.
\end{proof}

\begin{proof}[Proof of Remark \ref{j75}]
Assume that $\gamma_{0}$ is not reduced to a homoclinic trajectory. Then, as for any primitive cycle, $\gamma_{0}$ is a finite sequence of heteroclinic trajectories without sub-cycle. Let $e_{0} \in \gamma_{0}$ be such that
\begin{equation} \label{j78}
\frac{\alpha_{e_{0}^{-}} + \alpha_{e_{0}^{+}}}{\beta_{e_{0}^{-}} + \beta_{e_{0}^{+}}} = \min_{e \in \gamma_{0}} \frac{\alpha_{e^{-}} + \alpha_{e^{+}}}{\beta_{e^{-}} + \beta_{e^{+}}} .
\end{equation}
This quantity is the damping of the primitive cycle $\gamma_{1} = ( e_{0} , J ( e_{0} ) )$. Here, we have used that $J ( e_{0} ) \in \SE$ since $P$ is a Schr\"{o}dinger operator. From \eqref{j77} and \eqref{j78}, we get
\begin{equation*}
D ( \gamma_{0} ) = \frac{2 \alpha ( \gamma_{0} )}{2 \beta ( \gamma_{0} )} = \frac{\sum_{e \in \gamma_{0}} \alpha_{e^{-}} + \alpha_{e^{+}}}{\sum_{e \in \gamma_{0}} \beta_{e^{-}} + \beta_{e^{+}}} \geq \frac{\sum_{e \in \gamma_{0}} D ( \gamma_{1} ) ( \beta_{e^{-}} + \beta_{e^{+}} )}{\sum_{e \in \gamma_{0}} \beta_{e^{-}} + \beta_{e^{+}}} = D ( \gamma_{1} ) .
\end{equation*}
Then $\gamma_{1}$ is also a minimal primitive cycle. By uniqueness, we finally deduce $\gamma_{0} = \gamma_{1}$.
\end{proof}

It remains to prove the miscellaneous corollaries of Section \ref{s89}.

\begin{proof}[Proof of Corollary \ref{k40}]
To prove this result, it is enough to remark that
\begin{equation*}
\SQ^{\rm p r i n} = \left( \begin{array}{ccccc}
0 & & & \SQ_{e_{1} , e_{K}} & \\
\SQ_{e_{2} , e_{1}} & \ddots & & & \\
& \ddots & \ddots & & \\
& & \SQ_{e_{K} , e_{K - 1}} & 0 & \\
& & & & 0
\end{array} \right) .
\end{equation*}
Taking the determinant of $1 - \SQ^{\rm p r i n}$, it yields
\begin{align*}
F^{\rm p r i n} ( z , h ) &= 1 - \prod_{k = 1}^{K} \SQ_{e_{k + 1} , e_{k}} ( z , h )  \\
&= 1 - e^{i \beta ( \gamma_{0} ) Z} e^{i \beta ( \gamma_{0} ) \tau \vert \ln h \vert} \prod_{k = 1}^{K} \CQ_{e_{k + 1} , e_{k}} ( z , h ) .
\end{align*}
Then, \eqref{m27} gives
\begin{equation} \label{k51}
f_{\tau} ( Z , h ) = 1 - e^{i \beta ( \gamma_{0} ) Z} e^{i \beta ( \gamma_{0} ) \tau \vert \ln h \vert} \mu ( \tau , h ) .
\end{equation}
Its zeros are given by
\begin{equation*}
Z = - \tau \vert \ln h \vert + \frac{2 q \pi }{\beta ( \gamma_{0} )} + i \frac{\ln ( \mu ( \tau , h ) )}{\beta ( \gamma_{0} )} ,
\end{equation*}
for some $q \in \Z$. Eventually, Corollary \ref{k40} is a consequence of \eqref{m11}, the previous formula and Theorem \ref{j79}.
\end{proof}

\begin{proof}[Proof of Corollary \ref{k48}]
As in \eqref{k55}, any cycle $\gamma$ on $( \SV , \SE )$ of length $K$ verifies
\begin{equation*}
\frac{\alpha ( \gamma )}{\beta ( \gamma )} = \frac{K \alpha}{K \beta} = D_{0}.
\end{equation*}
Hence, $\alpha ( \gamma ) - D_{0} \beta ( \gamma ) = 0$ (i.e. all the cycles are minimal). Thus, all the coefficients $\alpha$ appearing in \eqref{m5} satisfy $\alpha = 0$. Comparing \eqref{m5} and \eqref{m6}, we deduce $F = F^{\rm p r i n}$. Then, \eqref{m27}, \eqref{k56} and $h^{S ( z , h ) / \lambda_{1} - 1 / 2} = e^{i \beta Z} h^{- i \beta \tau}$ imply
\begin{equation}
f_{\tau} ( Z , h ) = \det \Big( 1 - e^{i \beta Z} h^{- i \beta \tau} \widehat{\CQ} ( \tau , h ) \Big) .
\end{equation}
In particular, $Z \in [ - C , C ] + i [ - C , + \infty [$ is a zeros of $f_{\tau} ( \cdot , h )$ if and only if $e^{-i \beta Z} h^{i \beta \tau} \in \spe ( \widehat{\CQ} ( \tau , h ) )$ if and only if
\begin{equation*}
Z = - \tau \vert \ln h \vert + 2 q \pi \lambda_{1} + i \ln ( \mu_{k} ( \tau , h ) ) \lambda_{1} ,
\end{equation*}
for some $q \in \Z$ and $k \in \{ 1 , \ldots , \card \SE \}$. Finally, \eqref{k50} is a direct consequence of Proposition \ref{j81} and Theorem \ref{j79}.
\end{proof}

\begin{proof}[Proof of Corollary \ref{k54}]
As explained below Corollary \ref{k54}, we can always assume that each graph $( \SV_{k} , \SE_{k} )$ contains a minimal cycle of $( \SV , \SE )$. Using that the $\SV_{k}$'s are disjoint and \eqref{k60}, we get $\SQ^{\rm p r i n}_{e , \widetilde{e}} = 0$ for all $e \in \SE_{k}$ and $\widetilde{e} \in \SE_{\widetilde{k}}$ with $k \neq \widetilde{k}$. Then, the matrix $\SQ^{\rm p r i n}$ can be written
\begin{equation*}
\SQ^{\rm p r i n} = \left( \begin{array}{cccc}
\SQ^{\rm p r i n}_{1} & & &  \\
& \ddots & &  \\
& & \SQ^{\rm p r i n}_{K} &  \\
& & & 0
\end{array} \right) ,
\end{equation*}
where $\SQ^{\rm p r i n}_{k}$ is the principal matrix $\SQ^{\rm p r i n}$ associated to the operator $P_{k}$. Taking the determinant, we deduce
\begin{equation*}
F^{\rm p r i n} ( z , h ) = \prod_{k = 1}^{K} F^{\rm p r i n}_{k} ( z , h ) .
\end{equation*}
From \eqref{m6} and \eqref{m27}, it implies
\begin{equation*}
f_{\tau} ( z , h ) = \prod_{k = 1}^{K} f_{\tau}^{k} ( z , h ) ,
\end{equation*}
where $f_{\tau}^{k}$ is the function $f_{\tau}$ associated to $P_{k}$. Eventually, the corollary is a consequence of  Proposition \ref{j81} and Theorem \ref{j79}.
\end{proof}

\begin{proof}[Proof of \eqref{k73}]
Since the minimal primitive cycles are $( e_{1} , e_{2} )$ and $( e_{3} , e_{4} )$, the principal matrix $\SQ^{\rm p r i n}$ defined in \eqref{m3} can be written
\begin{equation}
\SQ^{\rm p r i n} = \left( \begin{array}{cccccc}
0 & \SQ_{e_{1} , e_{2}} & 0 & 0 & 0 & 0  \\
\SQ_{e_{2} , e_{1}} & 0 & 0 & \SQ_{e_{2} , e_{4}} & 0 & 0  \\
\SQ_{e_{3} , e_{1}} & 0 & 0 & \SQ_{e_{3} , e_{4}} & 0 & 0  \\
0 & 0 & \SQ_{e_{4} , e_{3}} & 0 & 0 & 0  \\
0 & 0 & 0 & 0 & 0 & 0  \\
0 & 0 & 0 & 0 & 0 & 0
\end{array} \right) .
\end{equation}
Taking the determinant of $1 - \SQ^{\rm p r i n}$, we get
\begin{equation*}
F^{\rm p r i n} = 1 - \SQ_{e_{1} , e_{2}} \SQ_{e_{2} , e_{1}} - \SQ_{e_{3} , e_{4}} \SQ_{e_{4} , e_{3}} + \SQ_{e_{1} , e_{2}} \SQ_{e_{2} , e_{1}} \SQ_{e_{3} , e_{4}} \SQ_{e_{4} , e_{3}} - \SQ_{e_{1} , e_{2}} \SQ_{e_{2} , e_{4}} \SQ_{e_{4} , e_{3}} \SQ_{e_{3} , e_{1}} .
\end{equation*}
Now, \eqref{m19}, \eqref{m27} and \eqref{k72} give
\begin{align*}
f_{\tau} ( Z , h ) ={}& 1 - e^{i Z \frac{\lambda_{a} + \lambda_{b}}{\lambda_{a} \lambda_{b}}} h^{- i \tau \frac{\lambda_{a} + \lambda_{b}}{\lambda_{a} \lambda_{b}}} \widehat{\CQ}_{e_{1} , e_{2}} \widehat{\CQ}_{e_{2} , e_{1}} - e^{i Z \frac{\lambda_{a} + \lambda_{b}}{\lambda_{a} \lambda_{b}}} h^{- i \tau \frac{\lambda_{a} + \lambda_{b}}{\lambda_{a} \lambda_{b}}} \widehat{\CQ}_{e_{3} , e_{4}} \widehat{\CQ}_{e_{4} , e_{3}}   \\
&+ e^{2 i Z \frac{\lambda_{a} + \lambda_{b}}{\lambda_{a} \lambda_{b}}} h^{- 2 i \tau \frac{\lambda_{a} + \lambda_{b}}{\lambda_{a} \lambda_{b}}} \widehat{\CQ}_{e_{1} , e_{2}} \widehat{\CQ}_{e_{2} , e_{1}} \widehat{\CQ}_{e_{3} , e_{4}} \widehat{\CQ}_{e_{4} , e_{3}}    \\
&- e^{2 i Z \frac{\lambda_{a} + \lambda_{b}}{\lambda_{a} \lambda_{b}}} h^{- 2 i \tau \frac{\lambda_{a} + \lambda_{b}}{\lambda_{a} \lambda_{b}}} \widehat{\CQ}_{e_{1} , e_{2}} \widehat{\CQ}_{e_{2} , e_{4}} \widehat{\CQ}_{e_{4} , e_{3}} \widehat{\CQ}_{e_{3} , e_{1}} ,
\end{align*}
since
\begin{equation*}
h^{S_{a} ( z , h ) / \lambda_{a} - 1 / 2} h^{S_{b} ( z , h ) / \lambda_{b} - 1 / 2} = h^{S_{c} ( z , h ) / \lambda_{c} - 1 / 2} h^{S_{b} ( z , h ) / \lambda_{b} - 1 / 2} = e^{i Z \frac{\lambda_{a} + \lambda_{b}}{\lambda_{a} \lambda_{b}}} h^{- i \tau \frac{\lambda_{a} + \lambda_{b}}{\lambda_{a} \lambda_{b}}} .
\end{equation*}
Using this identity and \eqref{k71}, we deduce
\begin{equation}
f_{\tau} ( Z , h ) = \det \Big( 1 - e^{i Z \frac{\lambda_{a} + \lambda_{b}}{\lambda_{a} \lambda_{b}}} h^{- i \tau \frac{\lambda_{a} + \lambda_{b}}{\lambda_{a} \lambda_{b}}} \widehat{\CQ}^{\rm r e d u} ( \tau , h ) \Big) .
\end{equation}

Then, $f_{\tau} ( Z , h ) = 0$ if and only if
\begin{equation*}
e^{- i Z \frac{\lambda_{a} + \lambda_{b}}{\lambda_{a} \lambda_{b}}} h^{i \tau \frac{\lambda_{a} + \lambda_{b}}{\lambda_{a} \lambda_{b}}} = \mu_{k} ( \tau , h ) ,
\end{equation*}
for some $k \in \{ 1 , 2 \}$. This equation is equivalent to
\begin{equation}
Z = - \tau \vert \ln h \vert + 2 q \pi \frac{\lambda_{a} \lambda_{b}}{\lambda_{a} + \lambda_{b}} + i \ln( \mu_{k} ( \tau , h ) ) \frac{\lambda_{a} \lambda_{b}}{\lambda_{a} + \lambda_{b}} ,
\end{equation}
for some $q \in \Z$. Eventually, combining this formula with Proposition \ref{j81} and Theorem \ref{j79}, we obtain \eqref{k73}.
\end{proof}

\section{Proofs of the asymptotic of the resonant states} \label{s40}

This part is devoted to the proof of the results and assertions of Section \ref{s32}. We begin with the main result and its corollaries collecting arguments developed for the asymptotic of resonances.

\begin{proof}[Proof of Theorem \ref{j35}]
Since the resonant states satisfies \eqref{a4}, we can apply the general results of Section \ref{s3}. Moreover, we are in the same geometric setting as in Section \ref{s61}. Thus, we will use the notations and some technical lemmas of Section \ref{s72}.

Let $v (h)$ be a family of normalized resonant states associated to a resonance $z ( h ) \in \eqref{d90}$. As in Section \ref{s72}, we define $v_{\pm}^{k}$ as the microlocal restriction of $v$ to a neighborhood of $U_{\pm}^{k}$ (see Figure \ref{f21}). These restrictions satisfy \eqref{d44}. In particular, there exist a constant $N \in \R$ such that $v_{-}^{k} \in \CI ( \Lambda_{+}^{1} , h^{- N} )$, for $k = 1 , \ldots , K$. Hence, there exist symbols $\widetilde{a}_{-}^{k} \in S ( h^{- N} )$ defined near $\pi_{x} ( U_{-}^{k} \cap \Lambda_{+}^{1} )$ and satisfying
\begin{equation} \label{r1}
v_{-}^{k} (x,h) = e^{- i A_{k} / h} e^{ i \frac{z - E_{0}}{h} t_{-}^{k}} \frac{\CM_{k}^{-}}{\CD_{k} ( t_{-}^{k} )} 
\widetilde{a}_{-}^{k} ( x , h ) e^{i \varphi_{+}^{1} (x) / h} ,
\end{equation}
where the action $A_{k}$, the dynamical quantities $\CD_{k}, \CM_{k}^{-}$ and the phase function $\varphi_{+}^{1}$ are defined in Section \ref{s61} and Section \ref{s32}. After a turn along the critical point and the homoclinic trajectories $\CH$, Lemma \ref{d41} gives the following formula for the symbol $\widetilde{a}_{-}^{k} ( x , h )$: 
\begin{equation} \label{r2}
\widetilde{a}_{-}^{k} ( x , h ) = h^{S ( z , h ) / \lambda_{1} - 1 / 2} \sum_{\ell = 1}^{K} \CP_{k , \ell} ( x , h ) \widetilde{a}_{-}^{\ell} ( x_{-}^{\ell} , h ) + S ( h^{- N + \zeta} ) ,
\end{equation}
for all $x$ near $\pi_{x} ( U_{-}^{k} \cap \Lambda_{+}^{1} )$. Here the constant $\zeta > 0$ and the symbols $\CP_{k , \ell} \in S ( 1 )$ are independent of $v$ and satisfy $\CP_{k , \ell} ( x_{-}^{k} , h ) = \CQ_{k , \ell} ( z , h )$ (see \eqref{d4}).

We set $\widetilde{\SA}_{0} ( h ) = ( \widetilde{a}_{-}^{1} ( x_{-}^{1} , h ) , \ldots , \widetilde{a}_{-}^{K} ( x_{-}^{K} , h ) )$, the $K$-vector associated to the resonant state $v$. Now, we are ready to prove the following estimate:
\begin{equation} \label{r3}
\text{There exists } M > 0 \text{ such that } \big\Vert \widetilde{\SA}_{0} ( h ) \big\Vert_{\ell^{2}} \in \big[ h^{M} , h^{-M} \big] \text{ for } h \text{ small enough.}
\end{equation}
First, remembering that $\widetilde{a}_{-}^{k} \in S ( h^{- N} )$ for $k = 1 , \ldots , K$, we get $\Vert \widetilde{a}_{-} ( h ) \Vert_{\ell^2} \leq  h^{- M}$ for $h$ small enough with $M = N + 1$. For the second inequality, we use a contradiction argument. If it is not satisfied, then $\Vert \widetilde{\SA}_{0} ( h ) \Vert_{\ell^{2}} = \CO ( h^{\infty} )$ for a sequence of positive numbers $h$ which converges to $0$. Then, \eqref{r2} yields $\widetilde{a}_{-}^{k} \in S ( h^{- N + \zeta} )$. Now, the standard bootstrap argument (see the end of Section \ref{s4}) implies that $\widetilde{a}_{-}^{k} \in S ( h^{\infty} )$ and hence $\Vert v \Vert = \CO ( h^{\infty} )$. This is in contradiction with $\Vert v \Vert = 1$. Thus, there exist $h_{0} , M > 0$ such that for all $0 < h \leq h_{0}$, $\Vert \widetilde{\SA}_{0} ( h ) \Vert_{\ell^{2}} \geq h^{M}$. Thus, \eqref{r3} holds true.

Taking $c ( h ) = \Vert \widetilde{\SA}_{0} ( h ) \Vert_{\ell^{2}}^{-1}$, we define $u = c v$. By transcribing the qualitative/quantitative properties of $v$, we will obtain the required properties for $u$. Lemma \ref{a33} and \eqref{r3} imply that the microsupport of $u$ is contained in $\{(0,0)\}\cup \Lambda_{+}$. From \eqref{r1} and \eqref{r3}, the symbol $a_{-}^{k} = c \widetilde{a}_{-}^{k} \in S ( h^{- N - M} )$ is such that
\begin{equation*}
u ( x , h ) = e^{- i A_{k} / h} e^{ i \frac{z - E_{0}}{h} t_{-}^{k}} \frac{\CM_{k}^{-}}{\CD_{k} ( t_{-}^{k} )} a_{-}^{k} ( x , h ) e^{i \varphi_{+}^{1} ( x ) / h} ,
\end{equation*}
microlocally near $\rho_{-}^{k}$. Moreover, the $K$-vector $\SA_{0} ( h )$ defined in Theorem \ref{j35} $iii)$ verifies $\SA_{0} = c \widetilde{\SA}_{0}$. Thus, our choice of $c$ guaranties that $\Vert \SA_{0} ( h ) \Vert_{\ell^{2}} = 1$. On the other hand, \eqref{r2} becomes
\begin{equation} \label{r4} 
a_{-}^{k} ( x , h ) = h^{S ( z , h ) / \lambda_{1} - 1 / 2} \sum_{\ell = 1}^{K} \CP_{k , \ell} ( x , h ) \SA_{0}^{\ell} ( h ) + S ( h^{- R + \zeta} ) ,
\end{equation}
if $a_{-}^{\ell} \in S ( h^{- R} )$ for all $\ell$. Combining the previous equation with $\vert h^{S ( z , h ) / \lambda_{1} - 1 / 2} \vert \lesssim 1$ for $z \in \eqref{d90}$ and $\Vert \SA_{0} ( h ) \Vert_{\ell^{2}} = 1$, we deduce
\begin{equation*} 
a_{-}^{k} \in S ( 1 ) + S ( h^{- R + \zeta} ) ,
\end{equation*}
if $a_{-}^{\ell} \in S ( h^{- R} )$ for all $\ell$. Then, a bootstrap argument implies $a_{-}^{k} \in S ( 1 )$. Lastly, taking $x = x_{-}^{k}$ in \eqref{r4} with $R = 0$, we get 
\begin{equation*} 
\big( h^{S ( z , h ) / \lambda_{1} - 1 / 2} \CQ ( z , h ) - 1 \big) \SA_{0} ( h ) = o ( 1 ) ,
\end{equation*}
as $h$ goes to $0$. This ends the proof of Theorem \ref{j35}.
\end{proof}

\begin{proof}[Proof of \eqref{j37}]
This formula relies on the propagation of Lagrangian distributions. Since $u$ satisfies $( P - z ) u =0$ and $u ( x , h ) = b_{-}^{k} ( x , h ) e^{i \varphi_{+}^{1} (x) / h}$ microlocally near $\gamma_{k} ( [ t_{-}^{k} , + \infty [ )$, the usual transport equations yield
\begin{equation*}
b_{-}^{k} ( x_{k} ( t ) , h ) = e^{i \frac{z - E_{0}}{h} ( t - t_{-}^{k} )} \frac{\CD_{k} ( t_{-}^{k} )}{\CD_{k} ( t )} b_{-}^{k} ( x_{k} ( t_{-}^{k} ) , h ) + o_{h \to 0}^{t} ( 1 ) ,
\end{equation*}
for all $t > t_{-}^{k}$. Using the normalization \eqref{j36}, the previous equation becomes
\begin{equation} \label{l5}
b_{-}^{k} ( x_{k} ( t ) , h ) = e^{- i A_{k} / h} e^{i \frac{z - E_{0}}{h} t} \frac{\CM_{k}^{-}}{\CD_{k} ( t )} \SA_{0}^{k} ( h ) + o_{h \to 0}^{t} ( 1 ) .
\end{equation}
Since \eqref{d7} can be written
\begin{equation*}
\CD_{k} ( t ) =  e^{t \frac{\sum \lambda_{j} - 2 \lambda_{1}}{2}} \CM_{k}^{-} ( 1 + o_{t \to + \infty} ( 1 ) ) ,
\end{equation*}
\eqref{l5} implies
\begin{equation} \label{l10}
b_{-}^{k} ( x_{k} ( t ) , h ) = e^{- i A_{k} / h} \SA_{0}^{k} ( h ) e^{t \big( i \frac{z - E_{0}}{h} - \frac{\sum \lambda_{j} - 2 \lambda_{1}}{2} \big)} \big( 1 + o_{t \to + \infty} ( 1 ) \big) + o_{h \to 0}^{t} ( 1 ) .
\end{equation}
We now use the asymptotic of the resonances. Combining Proposition \ref{d9} and Theorem \ref{d8}, we have
\begin{equation*}
\frac{z - E_{0}}{h} = \tau - i \Big( \frac{\sum \lambda_{j} - \lambda_{1}}{2} \Big) + o ( 1 ).
\end{equation*}
Then,
\begin{equation*}
e^{t \big( i \frac{z - E_{0}}{h} - \frac{\sum \lambda_{j} - 2 \lambda_{1}}{2} \big)} = e^{t \big( \frac{\lambda_{1}}{2} + i \tau + o ( 1 ) \big)} = e^{t \big( \frac{\lambda_{1}}{2} + i \tau \big)} + o_{h \to 0}^{t} ( 1 ) .
\end{equation*}
Using this equation and $\Vert \SA_{0} ( h ) \Vert = 1$, \eqref{l10} becomes
\begin{equation*}
b_{-}^{k} ( x_{k} ( t ) , h ) = e^{- i A_{k} / h} \SA_{0}^{k} ( h ) e^{t \big( \frac{\lambda_{1}}{2} + i \tau \big)} \big( 1 + o_{t \to + \infty} ( 1 ) \big) + o_{h \to 0}^{t} ( 1 ) .
\end{equation*}
which is \eqref{j37}.
\end{proof}

\begin{proof}[Proof of Remark \ref{j39}]
$i)$ The regularity of $\CQ$ with respect to the spectral parameter, the asymptotic \eqref{d95} of the resonance $z$ and the definition \eqref{e18} of $\widehat{\CQ}$ imply $\CQ ( z , h ) = \widehat{\CQ} ( \tau , h ) + o ( 1 )$ as $h$ goes to $0$. The same way, \eqref{d5} and \eqref{d95} lead to $h^{- S ( z , h ) / \lambda_{1} + 1 / 2} = \mu_{k} ( \tau ,h ) + o ( 1 )$. Moreover, $1 \lesssim \vert \mu_{k} ( \tau , h ) \vert \lesssim 1$ since $z \in \eqref{d90}$ (see the discussion after Proposition \ref{d9}). Eventually, \eqref{l7} follows from the previous arguments and \eqref{j38}.

$ii)$ We prove \eqref{j40} using a perturbation argument standard in spectral theory (see e.g. Helffer and Sj{\"o}strand \cite[Proposition 2.5]{HeSj84_01}). Let $\Pi ( \tau , h )$ the spectral projection of $\widehat{\CQ} ( \tau , h )$ associated to the eigenvalue $\mu_{k} ( \tau , h )$. Since $\mu_{k}$ is simple and isolated, this operator can be written
\begin{equation} \label{l9}
\Pi = - \frac{1}{2 i \pi} \oint_{\partial \CD} \big( \widehat{\CQ} - s \big)^{- 1} d s ,
\end{equation}
with $\CD = B ( \mu_{k} , \varepsilon / 2 )$. Combining with the resolvent identity, \eqref{l7} and Lemma \ref{d14}, we get
\begin{align}
\Pi \SA_{0} ( h ) - \SA_{0} ( h ) &= - \frac{1}{2 i \pi} \oint_{\partial \CD} \Big( \big( \widehat{\CQ} - s \big)^{- 1} - ( \mu_{k} - s \big)^{- 1} \Big) d s \, \SA_{0} ( h )  \nonumber \\
&= \frac{1}{2 i \pi} \oint_{\partial \CD} \big( \widehat{\CQ} - s \big)^{- 1} ( \mu_{k} - s )^{- 1} d s \big( \widehat{\CQ} - \mu_{k} \big) \SA_{0} ( h )  \nonumber \\
&= o ( 1 ) , \label{l8}
\end{align}
as $h$ goes to $0$.

On the other hand, we have $\Pi \SA_{0} \in \im \Pi = \C f_{k}$ since $\mu_{k}$ is simple. Thus, there exists $\beta ( z , h ) \in \C$ such that
\begin{equation*}
\SA_{0} = \beta f_{k} + o ( 1 ) .
\end{equation*}
Using that $\SA_{0}$ and $f_{k}$ are both normalized, we deduce $\vert \beta \vert = 1 + o ( 1 )$. Then, if we set $\alpha = \beta / \vert \beta \vert$, the previous equation yields
\begin{equation*}
\SA_{0} = \alpha f_{k} + \alpha ( \vert \beta \vert - 1 ) f_{k} + o ( 1 ) = \alpha f_{k} + o ( 1 ) ,
\end{equation*}
and \eqref{j40} follows.

If $\mu_{k}$ is not isolated, the previous argument can be adapted. More precisely, we fix $\delta > 0$ and remark that $\widehat{\CQ}$ has at most $K$ different eigenvalues. Then, for all $( \tau , h )$, there exists $\varepsilon ( \tau , h )$ such that $\delta / 2 < \varepsilon < \delta$ and
\begin{equation*}
\dist \big( \partial \CD , \spe ( \widehat{\CQ} ( \tau , h ) ) \big) \geq \frac{\delta}{4 ( K + 1 )} ,
\end{equation*}
with $\CD = B ( \mu_{k} , \varepsilon )$. We define again $\Pi$ using \eqref{l9}. This operator is now the sum of the spectral projections of $\widehat{\CQ}$ associated to the eigenvalues in $\CD$. Following the proof of \eqref{l8}, we obtain $\SA_{0} = \Pi \SA_{0} + o ( 1 )$ which shows that $\SA_{0}$ can be approximated by a linear combination of generalized eigenvectors associated to the eigenvalues of $\widehat{\CQ}$ in $B ( \mu_{k} , \delta )$.

$iii)$ Let $v$ be a resonant state of $P$ associated to $z$ and normalized as in Theorem \ref{j35} $iii)$. Thanks to \eqref{j40}, there exists $\alpha$ of modulus $1$ such that $\SA_{0} = \alpha f_{k} + o ( 1 )$. Then, $u = v / \alpha$ satisfies the required properties.
\end{proof}

We now treat the different examples using techniques of linear algebra.

\begin{proof}[Proof of Example \ref{l16}]
\eqref{j46} and \eqref{j47} could have been shown applying perturbations arguments. Instead, we will use explicit formulas for the eigenvectors of a generic $2 \times 2$ matrix \eqref{i44}. When $\vert d - a \vert \gtrsim 1$ and $\vert b \vert , \vert c \vert \ll 1$, a direct calculus proves that $( f_{1} , f_{2} )$ forms a basis of eigenvectors of \eqref{i44} where
\begin{equation} \label{l12}
f_{1} = \left( \begin{array}{c}
1 \\
\frac{- 2 c}{d - a + \sqrt{( d - a )^{2} + 4 b c}}
\end{array} \right) , \qquad
f_{2} = \left( \begin{array}{c}
\frac{2 b}{d - a + \sqrt{( d - a )^{2} + 4 b c}} \\
1
\end{array} \right) ,
\end{equation}
and $\sqrt{( d - a )^{2} + 4 b c}$ is the square root of $( d - a )^{2} + 4 b c$ close to $d - a$. On the other hand, when $\vert b \vert , \vert c \vert \gtrsim 1$ and $\vert a \vert , \vert d \vert \ll 1$, $( g_{1} , g_{2} )$ forms a basis of eigenvectors of \eqref{i44} where
\begin{equation} \label{l13}
g_{1} = \left( \begin{array}{c}
1 \\
\frac{- 2 c}{d - a + \sqrt{( d - a )^{2} + 4 b c}}
\end{array} \right)
\qquad \text{and} \qquad
g_{2} = \left( \begin{array}{c}
1 \\
\frac{- 2 c}{d - a - \sqrt{( d - a )^{2} + 4 b c}}
\end{array} \right) .
\end{equation}

Using the explicit expression for $\CQ$ given in \eqref{d4}, we get
\begin{equation*}
\widehat{\CQ} ( \tau , h ) = \left( \begin{array}{cc}
e^{i A_{1} / h} \mu_{1} ( \tau ) & 0 \\
0 & e^{i A_{2} / h} \mu_{2} ( \tau )
\end{array} \right) + o_{\tau \to - \infty} ( 1 ) ,
\end{equation*}
where the $\mu_{\bullet} ( \tau )$ are given in \eqref{l11}. From \eqref{l14}, we can apply \eqref{l12} and \eqref{j46} follows. On the other hand, in the limit $\tau \to + \infty$, the definition \eqref{d4} of $\CQ$ leads to
\begin{equation*}
\widehat{\CQ} ( \tau , h ) = \left( \begin{array}{cc}
0 & e^{i A_{1} / h} \widetilde{\mu}_{1} ( \tau ) \\
e^{i A_{2} / h} \widetilde{\mu}_{2} ( \tau ) & 0
\end{array} \right) + o_{\tau \to - \infty} ( 1 ) ,
\end{equation*}
with
\begin{equation*}
\widetilde{\mu}_{k} ( \tau ) = \Gamma \Big( \frac{1}{2} - i \frac{\tau}{\lambda} \Big) \sqrt{\frac{\lambda}{2 \pi}} \frac{\CM^{+}_{k}}{\CM^{-}_{k}} e^{- \frac{\pi}{2} \nu_{k} i} \vert g_{-}^{3 - k} \vert \big( \lambda \vert g_{+}^{k} \cdot g_{-}^{3 - k} \vert \big)^{- \frac{1}{2} + i \frac{\tau}{\lambda}} e^{\frac{\pi \tau}{2 \lambda}} .
\end{equation*}
Using \eqref{e36}, we get $\vert \widetilde{\mu}_{k} \vert \gtrsim 1$. Then, \eqref{l13} implies \eqref{j47} with
\begin{align}
\alpha ( \tau , h ) &= - \sqrt{\frac{e^{i A_{2} / h} \widetilde{\mu}_{2} ( \tau )}{e^{i A_{1} / h} \widetilde{\mu}_{1} ( \tau )}}   \nonumber \\
&= - e^{i \frac{A_{2} - A_{1}}{2 h}} \sqrt{\frac{\CM^{+}_{2} \CM^{-}_{1}}{\CM^{-}_{2} \CM^{+}_{1}}} e^{- \frac{\pi}{4} ( \nu_{2} - \nu_{1} ) i} \sqrt{\frac{\vert g_{-}^{1} \vert}{\vert g_{-}^{2} \vert}} \bigg( \frac{\vert g_{+}^{2} \cdot g_{-}^{1} \vert}{\vert g_{+}^{1} \cdot g_{-}^{2} \vert} \bigg)^{- \frac{1}{4} + i \frac{\tau}{2 \lambda}} . \label{l15}
\end{align}
In particular, $\alpha ( \tau , h )$ has a non-zero constant modulus.
\end{proof}

\begin{proof}[Proof of Example \ref{l17}]
Since $\gamma_{1}$ and $\gamma_{2}$ are symmetric with respect to $\gamma_{3}$, the matrix $\widehat{\CQ} ( \tau , h )$ can be written
\begin{equation} \label{l18}
\widehat{\CQ} = \left( \begin{array}{ccc}
a & b & d \\
b & a & d \\
c & c & e
\end{array} \right) ,
\end{equation}
where the coefficients are given by
\begin{align*}
a &= e^{i A_{1} / h} \Gamma \Big( \frac{1}{2} - i \frac{\tau}{\lambda} \Big) \sqrt{\frac{\lambda}{2 \pi}} \frac{\CM_{1}^{+}}{\CM_{1}^{-}} e^{- \frac{\pi}{2} ( \nu_{1} + 1 ) i} \big\vert g_{-}^{1} \big\vert \big( \lambda \vert g_{+}^{1} \vert \vert g_{-}^{1} \vert \big)^{{- \frac{1}{2} + i \frac{\tau}{\lambda}}} e^{- \frac{\pi \tau}{2 \lambda}} ,  \\
b &= e^{i A_{1} / h} \Gamma \Big( \frac{1}{2} - i \frac{\tau}{\lambda} \Big) \sqrt{\frac{\lambda}{2 \pi}} \frac{\CM_{1}^{+}}{\CM_{1}^{-}} e^{- \frac{\pi}{2} ( \nu_{1} + 1 ) i} \big\vert g_{-}^{1} \big\vert \big( \lambda g_{+}^{1} \cdot g_{-}^{2} \big)^{{- \frac{1}{2} + i \frac{\tau}{\lambda}}} e^{- \frac{\pi \tau}{2 \lambda}} ,   \\
c &= e^{i A_{3} / h} \Gamma \Big( \frac{1}{2} - i \frac{\tau}{\lambda} \Big) \sqrt{\frac{\lambda}{2 \pi}} \frac{\CM_{3}^{+}}{\CM_{3}^{-}} e^{- \frac{\pi}{2} ( \nu_{3} + 1 ) i} \big\vert g_{-}^{3} \big\vert \big( \lambda g_{+}^{3} \cdot g_{-}^{1} \big)^{{- \frac{1}{2} + i \frac{\tau}{\lambda}}} e^{- \frac{\pi \tau}{2 \lambda}} , \\
d &= e^{i A_{1} / h} \Gamma \Big( \frac{1}{2} - i \frac{\tau}{\lambda} \Big) \sqrt{\frac{\lambda}{2 \pi}} \frac{\CM_{1}^{+}}{\CM_{1}^{-}} e^{- \frac{\pi}{2} ( \nu_{1} + 1 ) i} \big\vert g_{-}^{1} \big\vert \big( \lambda g_{+}^{1} \cdot g_{-}^{3} \big)^{{- \frac{1}{2} + i \frac{\tau}{\lambda}}} e^{- \frac{\pi \tau}{2 \lambda}} , \\
e &= e^{i A_{3} / h} \Gamma \Big( \frac{1}{2} - i \frac{\tau}{\lambda} \Big) \sqrt{\frac{\lambda}{2 \pi}} \frac{\CM_{3}^{+}}{\CM_{3}^{-}} e^{- \frac{\pi}{2} ( \nu_{3} + 1 ) i} \big\vert g_{-}^{3} \big\vert \big( \lambda \vert g_{+}^{3} \vert \vert g_{-}^{3} \vert \big)^{{- \frac{1}{2} + i \frac{\tau}{\lambda}}} e^{- \frac{\pi \tau}{2 \lambda}} .
\end{align*}
Computing the zeros of the characteristic polynomial of \eqref{l18}, one can verify that the eigenvalues of $\widehat{\CQ}$ are given by
\begin{align*}
\mu_{1} &= a - b , \\
\mu_{2} &= \frac{a + b + e}{2} + \frac{\sqrt{( a + b - e )^{2} + 8 c d}}{2} , \\
\mu_{3} &= \frac{a + b + e}{2} - \frac{\sqrt{( a + b - e )^{2} + 8 c d}}{2} .
\end{align*}
Moreover, $f = {}^{t} ( 1 , - 1 , 0 ) / \sqrt{2}$ is always a normalized eigenvector of $\widehat{\CQ}$ associated to the eigenvalue $\mu_{1}$.

Since $0 < g_{+}^{1} \cdot g_{-}^{2} < \vert g_{+}^{1} \vert \vert g_{-}^{1} \vert$, there exists $\varepsilon > 0$ such that $\vert \mu_{1} ( \tau , h ) \vert \geq \varepsilon$ for all $\tau \in [ - C , C]$ and $h \in ] 0 , 1 ]$. Thus, $\mu_{1}$ is far away from $0$. We now choose the geometry in a such way that the quantity
\begin{equation} \label{l19}
\frac{\CM^{+}_{3}}{\CM^{-}_{3}} \sqrt{\frac{\vert g_{-}^{3} \vert}{\vert g_{+}^{3} \vert}} ,
\end{equation}
is small. Example \ref{i41} provides different methods to realize this property. For instance, one can take the obstacle $\CO$ in a such way that the radius of curvature of $\partial \CO$ near $x_{3}$ is small. Moreover, we can always parametrize $\gamma_{3}$ such that $\vert g_{-}^{3} \vert = 1$. In this setting, $a , b , d$ are independent of \eqref{l19}, and $c , e$ go to $0$ with \eqref{l19}. Thus, $\mu_{2}$ and $\mu_{3}$ are close to $0$ and $a + b$ for \eqref{l19} small. Since $a - b$ avoid these values, $\mu_{1}$ is simple and isolated.
\end{proof}

We now study the delocalization phenomenon.

\begin{proof}[Proof of Example \ref{k87}]
Here, we show that \eqref{k86} holds true in Example \ref{k87} for $\lambda_{2}$ small enough. If it is not the case, there exists a (decreasing) sequence of $h$ which goes to $0$ and a sequence of normalized resonant states $u = u ( h )$ associated to some resonances $z = z ( h ) \in \eqref{j80}$ such that
\begin{equation} \label{k88}
\Vert \Op ( \varphi_{1} ) u \Vert > h^{C_{0}} \Vert \Op ( \varphi_{2} ) u \Vert ,
\end{equation}
for $h$ small enough. In the previous inequality, $\one_{K_{j}} \prec \varphi_{j} \in C^{\infty}_{0} ( T^{*} \R^{n} )$ is supported in a small neighborhood of $K_{j}$.

We use notations and results of Section \ref{s48}. For $k = 1 , 2$, let $u_{-}^{k}$ be the microlocal restriction of $u$ near $\rho_{-}^{k} = ( x_{-}^{k} , \xi_{-}^{k} ) \in e_{k}$ close to $( v_{k} , 0 )$. As in \eqref{k90}, we have
\begin{equation} \label{k92}
u_{-}^{1} \in \CI ( \Lambda_{+}^{v_{1}} , h^{- N} ) \qquad \text{and} \qquad u_{-}^{2} \in \CI ( \Lambda_{+}^{v_{1}} , h^{- N} ) ,
\end{equation}
for some $N \in \R$. Their symbols $a_{-}^{k} ( x , h ) \in S ( h^{- N} )$ are defined by
\begin{equation*}
u_{-}^{k} ( x , h ) = e^{- i A_{k} / h} e^{i \frac{z - E_{0}}{h} t_{-}^{k}} \frac{\CM_{k}^{-}}{\CD_{k} ( t_{-}^{k} )} a_{-}^{k} ( x , h ) e^{i \varphi_{+}^{1 , k} ( x ) / h} .
\end{equation*}
(see \eqref{k15}) and they satisfy
\begin{equation} \label{k91}
\left\{ \begin{aligned}
a_{-}^{1} ( x , h ) = \CP_{1 , 1} ( x , h ) a_{-}^{1} ( x_{-}^{1} , h ) + S ( h^{- N + \zeta} ) ,  \\
a_{-}^{2} ( x , h ) = \CP_{2 , 1} ( x , h ) a_{-}^{1} ( x_{-}^{1} , h ) + S ( h^{- N + \zeta} ) ,
\end{aligned} \right.
\end{equation}
from Lemma \ref{k10}.

Mimicking \eqref{r3}, one can show that there exists $M > 0$ such that 
\begin{equation*}
\big\vert a_{-}^{1} ( x_{-}^{1} , h ) \big\vert \in \big[ h^{M} , h^{- M} \big]
\end{equation*}
for $h$ small enough. Then, after a renormalization, we can always assume that
\begin{equation}
a_{-}^{1} ( x_{-}^{1} , h ) = 1 ,
\end{equation}
$N = 0$ in \eqref{k92}--\eqref{k91} and $\Vert u \Vert \leq h^{- R}$ for some $R \in \R$. Theorem \ref{a32} yields
\begin{equation*}
u = \CO ( h^{- N_{0}} ) \text{ microlocally near } ( v_{1} , 0 ) ,
\end{equation*}
where $N_{0}$ depends only of $\lambda_{1}$ and \eqref{j80}. Combining with \eqref{k92}, this gives
\begin{equation} \label{k93}
\Vert \Op ( \varphi_{1} ) u \Vert \lesssim h^{- N_{0}} .
\end{equation}

On the other hand, \eqref{k91} and the definition \eqref{j83} of $\CP_{2 , 1} ( x_{-}^{2} , h ) = \CQ_{2 , 1} ( z , h )$ imply that the symbol $a_{-}^{2} (x , h )$ is elliptic near $x_{-}^{2}$. Let $\rho \in \Lambda_{+}^{v_{2}}$ close to $( v_{2} , 0 )$ be such that
\begin{equation*}
g_{+} ( \rho ) \cdot g_{-}^{2} \neq 0 .
\end{equation*}
From Corollary \ref{d46} and the ellipticity of $a_{-}^{2}$, we deduce that
\begin{equation*}
u \in \CI \Big( \Lambda_{+}^{v_{1}} , h^{\frac{n - 1}{2} - \frac{D_{0}}{\lambda_{2}}} \Big) \text{ microlocally near } \rho ,
\end{equation*}
with an elliptic symbol. Since $D_{0} = ( n - 1 ) \lambda_{1} / 2$, this equation becomes
\begin{equation*}
u \in \CI \Big( \Lambda_{+}^{v_{1}} , h^{- \frac{n - 1}{2} \frac{\lambda_{1} - \lambda_{2}}{\lambda_{2}}} \Big) \text{ microlocally near } \rho ,
\end{equation*}
with an elliptic symbol. In particular,
\begin{equation} \label{k89}
\Vert \Op ( \varphi_{2} ) u \Vert \gtrsim h^{- \frac{n - 1}{2} \frac{\lambda_{1} - \lambda_{2}}{\lambda_{2}}} \geq h^{- N_{0} - C_{0} - 1} ,
\end{equation}
for $\lambda_{2}$ small enough. Eventually, \eqref{k93} and \eqref{k89} are in contradiction with \eqref{k88}. This ends the proof of \eqref{k86}.
\end{proof}

Let us show that the resonant states are almost identical on $\CH$.

\begin{proof}[Proof of Corollary \ref{j43}]
Let $f_{k_{0}} \in \C^{K}$ be a normalized eigenvector of $\widehat{\CQ} ( \tau , h )$ associated to the simple eigenvalue $\mu_{k_{0}} ( \tau , h )$. Let $a_{j} , \SA_{j}$ denote the symbols provided by Theorem \ref{j35} $ii)$ for the resonant state $u_{j}$, $j = 1 , 2$. That is
\begin{equation} \label{j64}
u_{j} ( x , h ) = e^{- i A_{k} / h} e^{i \frac{z_{j} - E_{0}}{h} t_{-}^{k}} \frac{\CM_{k}^{-}}{\CD_{k} ( t_{-}^{k} )} a_{j}^{k} ( x , h ) e^{i \varphi_{+}^{1} (x) / h} ,
\end{equation}
microlocally near $\rho_{-}^{k}$ and $\SA_{j}^{k} ( h ) = a_{j}^{k} ( x_{-}^{k} , h )$. From Remark \ref{j39} $ii)$, there exists $\alpha_{j} (  h ) \in \C$ with $\vert \alpha_{j} \vert = 1$ such that
\begin{equation*}
\SA_{j}  = \alpha_{j} f_{k_{0}} + o ( 1 ) ,
\end{equation*}
as $h$ goes to $0$. We then define $\alpha = \alpha_{1} / \alpha_{2}$ which satisfies $\vert \alpha \vert = 1$. The last equation gives
\begin{equation} \label{l22}
\SA_{1}  = \alpha \SA_{2} + o ( 1 ) .
\end{equation}

Now, let $\CP_{k , \ell}^{j} \in S ( 1 )$ denote the symbols of Lemma \ref{d41} computed at $z = z_{j}$. From \eqref{d42}, we have
\begin{equation} \label{l26}
a_{j}^{k} ( x , h ) = h^{S ( z_{j} , h ) / \lambda_{1} - 1 / 2} \sum_{\ell = 1}^{K} \CP_{k , \ell}^{j} ( x , h ) \SA_{j}^{k} ( h ) + S ( h^{\zeta} ) .
\end{equation}
Following the proof of Lemma \ref{d41} and using that $z_{1} = z_{2} + o ( h )$, one can verify that
\begin{equation} \label{l24}
\CP_{k , \ell}^{1} ( x , h ) = \CP_{k , \ell}^{2} ( x , h ) + o ( 1 ) ,
\end{equation}
uniformly for $x$ near $\pi_{x} ( U_{-}^{k} \cap \Lambda_{+}^{1} )$. Moreover, \eqref{d95} and \eqref{l23} imply
\begin{equation*}
z_{1} = z_{2} + 2 \pi ( q_{1} - q_{2} ) \lambda_{1} \frac{h}{\vert \ln h \vert} + o \Big( \frac{h}{\vert \ln h \vert} \Big) ,
\end{equation*}
and then
\begin{equation} \label{l25}
h^{S ( z_{1} , h ) / \lambda_{1} - 1 / 2} = h^{S ( z_{2} , h ) / \lambda_{1} - 1 / 2} e^{i ( z_{1} - z_{2} ) \frac{\vert \ln h \vert}{\lambda_{1} h}} = h^{S ( z_{2} , h ) / \lambda_{1} - 1 / 2} + o ( 1 ) ,
\end{equation}
since $1 \lesssim \vert h^{S ( z_{j} , h ) / \lambda_{1} - 1 / 2} \vert \lesssim 1$. Combining \eqref{l26} with the estimates \eqref{l22}, \eqref{l24} and \eqref{l25}, we obtain
\begin{equation} \label{j65}
a_{1}^{k} ( x , h ) = \alpha a_{2}^{k} ( x , h ) + o ( 1 ) ,
\end{equation}
uniformly for $x$ near $\pi_{x} ( U_{-}^{k} \cap \Lambda_{+}^{1} )$. Using again $z_{1} = z_{2} + o ( h )$, we get $e^{i \frac{z_{1} - E_{0}}{h} t_{-}^{k}} = e^{i \frac{z_{2} - E_{0}}{h} t_{-}^{k}} + o ( 1 )$. Then, \eqref{j64} and \eqref{j65} imply
\begin{equation} \label{l27}
u_{1} = \alpha u_{2} + o ( 1 ) \text{ microlocally near } \rho_{-}^{k}.
\end{equation}
Finally, the usual propagation of singularities along $\CH$ and $( P - z_{1} ) ( u_{1} - \alpha u_{2} ) = o ( 1 )$ imply that \eqref{l27} holds true microlocally near each point of $\CH$.
\end{proof}

We now show that resonant states are almost orthogonal in the well in an island situation. This result is essentially contained in \cite{HeSj86_01}. We provide some details since our distortion is slightly different from that of Helffer and Sj\"{o}strand.

\begin{figure}
\begin{center}
\begin{picture}(0,0)%
\includegraphics{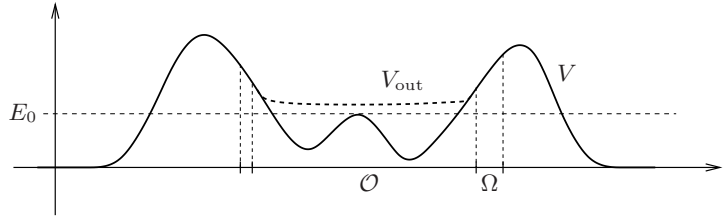}%
\end{picture}%
\setlength{\unitlength}{987sp}%
\begingroup\makeatletter\ifx\SetFigFont\undefined%
\gdef\SetFigFont#1#2#3#4#5{%
  \reset@font\fontsize{#1}{#2pt}%
  \fontfamily{#3}\fontseries{#4}\fontshape{#5}%
  \selectfont}%
\fi\endgroup%
\begin{picture}(17916,5466)(-11882,-5194)
\put(1801,-1786){\makebox(0,0)[lb]{\smash{{\SetFigFont{9}{10.8}{\rmdefault}{\mddefault}{\updefault}$V$}}}}
\put(-11249,-2761){\makebox(0,0)[rb]{\smash{{\SetFigFont{9}{10.8}{\rmdefault}{\mddefault}{\updefault}$E_{0}$}}}}
\put(-2924,-4561){\makebox(0,0)[b]{\smash{{\SetFigFont{9}{10.8}{\rmdefault}{\mddefault}{\updefault}$\CO$}}}}
\put(113,-4561){\makebox(0,0)[b]{\smash{{\SetFigFont{9}{10.8}{\rmdefault}{\mddefault}{\updefault}$\Omega$}}}}
\put(-2624,-1936){\makebox(0,0)[lb]{\smash{{\SetFigFont{9}{10.8}{\rmdefault}{\mddefault}{\updefault}$V_{\rm out}$}}}}
\end{picture}%
\end{center}
\caption{The sets $\CO , \Omega$ and the potential $V_{\rm out}$ in the proof of \eqref{j44}.} \label{f47}
\end{figure}

\begin{proof}[Proof of \eqref{j44}]
Let $\CO \Subset \Omega$ be two small smooth neighborhoods of the well (the compact component of $\{ x ; \ V ( x ) \leq E_{0} \}$). We construct two operator $P_{\rm in}$ and $P_{\rm out}$ as follows. First $P_{\rm in}$ is the restriction of $P$ in $\Omega$ with Dirichlet boundary condition. Secondly $P_{\rm out} = - h^{2} \Delta + V_{\rm out} ( x )$ is an operator $P$ with the well filled up inside $\CO$ (see Figure \ref{f47}). Eventually, consider $\chi \in C^{\infty}_{0} ( \Omega )$ such that $\chi = 1$ near $\CO$.

From \cite[Th\'eor\`eme 9.9]{HeSj86_01}, there exist two normalized functions $w_{1} , w_{2}$, which belong to the sum of the eigenspaces of $P_{\rm in}$ associated to its eigenvalues close to $z_{1} , z_{2}$ respectively, such that
\begin{equation} \label{r6}
u_{j} = w_{j} + \CO ( e^{- \varepsilon /h} ) ,
\end{equation}
inside $\Omega$ for some $\varepsilon > 0$ which may change from line to line. Moreover, if $P_{\rm out , \theta}$ designs the complex distortion of $P_{\rm out}$, we can write
\begin{equation*}
( P_{\rm out , \theta} - z ) ( 1 - \chi ) u_{j} = ( P_{\theta} - z ) ( 1 - \chi ) u_{j} = ( 1 - \chi ) ( P_{\theta} - z ) u_{j} - [ P , \chi ] u_{j} = - [ P , \chi ] u_{j} .
\end{equation*}
Since $[ P , \chi ]$ is localized outside the well, \eqref{r6} and the standard properties of the eigenvectors of $P_{\rm in}$ imply that $[ P , \chi ] u_{j} = \CO ( e^{- \varepsilon /h} )$. Moreover, since $P_{\rm out}$ is non-trapping at energy $E_{0}$, \cite{Ma02_01} gives a polynomial estimate of the resolvent of $P_{\rm out , \theta}$. Thus, the last equation becomes
\begin{equation} \label{r7}
( 1 - \chi ) u_{j} = - ( P_{\rm out , \theta} - z )^{- 1} [ P , \chi ] u_{j} = \CO ( e^{- \varepsilon /h} ) .
\end{equation}

Summing up, \eqref{r6} and \eqref{r7} give
\begin{equation} \label{r8}
u_{j} = w_{j} + \CO ( e^{- \varepsilon /h} ) ,
\end{equation}
in $L^{2} ( \R^{n} )$. On the other hand, the functions $w_{1}$ and $w_{2}$ are orthogonal in $L^{2} ( \Omega )$ since $\vert z_{1} - z_{2} \vert \gtrsim h \vert \ln h \vert^{- 1}$. Then, \eqref{r8} yields
\begin{equation*}
\< u_{1} , u_{2} \>_{L^{2} ( \R^{n} )} = \< w_{1} , w_{2} \>_{L^{2} ( \Omega )} + \CO ( e^{- \varepsilon /h} ) = \CO ( e^{- \varepsilon /h} ) ,
\end{equation*}
and \eqref{j44} follows.
\end{proof}

Lastly, we consider the quasimodes and the quasiresonances.

\begin{proof}[Proof of Proposition \ref{j45}]
Assume that $i)$ is not satisfied. Then, there exist $\delta > 0$ and a (decreasing) sequence of $h$ which goes to $0$ such that
\begin{equation*}
\dist \big( z , \res (P) \big) \geq \delta \frac {h}{\vert \ln h \vert} ,
\end{equation*}
for $h$ in this sequence. In particular, Theorem \ref{d8} and Proposition \ref{j66} give $\Vert ( P_{\theta} - z )^{- 1} \Vert \lesssim h^{- M}$. Applying this estimate to $( P_{\theta} - z ) v $, we deduce
\begin{equation*}
1 = \Vert v \Vert \lesssim h^{- M} \Vert ( P_{\theta} - z ) v \Vert = \CO ( h^{\infty} ) .
\end{equation*}
Thus, we get a contradiction and $i)$ holds true. To obtain $ii)$, it is enough to follow the proof of Theorem \ref{j35}.
\end{proof}

\appendix

\section{Review of semiclassical analysis} \label{s9}

\Subsection{Microlocal terminology} \label{s35}

Basic notions of semiclassical microlocal analysis used throughout the paper are recalled in this part. We send back the reader to the books of Dimassi and Sj\"{o}strand \cite{DiSj99_01}, Guillemin and Sternberg \cite{GuSt13_01}, Martinez \cite{Ma02_02}, Robert \cite{Ro87_01} and Zworski \cite{Zw12_01} for more details on this theory.

We say that $m : \R^{d} \longrightarrow [ 0 , + \infty [$ is an {\it order function} if there exists $C > 0$ such that
\begin{equation*}
m ( x ) \leq C \< x - y \>^{C} m ( y ) ,
\end{equation*}
for all $x , y \in \R^{d}$. If $m ( x )$ is an order function, we say that a smooth function $a ( x , h )$ is a {\it symbol of class} $S ( m )$ when
\begin{equation*}
\forall \alpha \in \N^{d} , \quad \exists C > 0 , \quad \forall h \in ] 0 , 1 ] , \qquad \vert \partial_{x}^{\alpha} a ( x , h ) \vert \leq C m ( x ) .
\end{equation*}
If $f ( h )$ is a function of $h$ only, we write  $S ( f (  h ) m )$ instead of $f ( h ) S ( m )$.

Let $m$ be an order function on $T^{*} \R^{n}$. If $a ( x , \xi , h )$ is a symbol of class $S ( m )$, the {\it semiclassical pseudodifferential operator} $\Op ( a )$ with symbol $a$ is defined by
\begin{equation}
( \Op ( a ) u ) ( x ) = \frac{1}{(2 \pi h)^{n}} \iint e^{i ( x - y ) \cdot \xi / h } a \Big( \frac{x + y}{2} , \xi \Big) u(y) \, d y  \, d \xi ,
\end{equation}
for $u \in \CS^{\prime} ( \R^{n} )$ in the sense of tempered distributions. Since we are interested in spectral analysis and deal with self-adjoint operators, it is natural to use the Weyl quantization. We denote by $\Psi ( m )$ the space of operators $\Op ( S ( m ) )$.

A function $u \in L^{2} ( \R^{n} )$ is said to be {\it polynomially bounded} when $\Vert u \Vert \lesssim h^{- C}$ for some $C > 0$. For such a function $u$ and $V \subset T^{*} \R^{n}$ bounded, we say that $u = 0$ {\it microlocally near} $V$ if and only if
\begin{equation*}
\Vert \Op ( \psi ) u \Vert_{L^{2}} = \CO ( h^{\infty} ) ,
\end{equation*}
for some $\psi \in C^{\infty}_{0} ( T^{*} \R^{n} )$ with $\psi = 1$ near $V$.

We now recall the definition of the semiclassical Lagrangian distributions. For the general theory in the classical setting, we refer to H\"{o}rmander \cite{Ho71_01} or \cite[Section 25.1]{Ho94_01}. In the present semiclassical setting, this notion is developed in the book of Ivrii \cite[Section 1.2]{Iv98_01}, in the PhD thesis of Dozias \cite{Do94_01} and in the lecture notes of Colin de Verdi{\`e}re \cite{Co16_01}.

A manifold $\Lambda \subset T^{*} \R^{n}$ is called a {\it Lagrangian manifold} when $d \xi \wedge d x \vert_{\Lambda} = 0$ and $\dim \Lambda = n$. Consider $\varphi (x , \theta ) \in C^{\infty} ( \Omega )$ where $\Omega$ is an open set of $\R^{n + r}$. We say that $\varphi$ is a {\it non-degenerate phase function} if and only if, for all $( x , \theta ) \in C_{\varphi}$ with
\begin{equation*}
C_{\varphi} = \{ ( x , \theta ) \in \Omega ; \ \partial_{\theta} \varphi = 0 \} ,
\end{equation*}
the $r$ differentials $d \partial_{\theta_{1}} \varphi , \ldots , d \partial_{\theta_{r}} \varphi$ are linearly independent. If $\varphi$ is a non-degenerate phase function, $C_{\varphi}$ is a manifold of dimension $n$ and
\begin{equation*}
\begin{aligned}
j_{\varphi} : \\
{}
\end{aligned}
\left\{ \begin{aligned}
&C_{\varphi} &&\longrightarrow &&T^{*} \R^{n} \\
&( x , \theta ) && &&( x , \partial_{x} \varphi )
\end{aligned} \right. 
\end{equation*}
is locally a diffeomorphism whose image $\Lambda_{\varphi} : = j_{\varphi} ( C_{\varphi} )$ is a Lagrangian manifold.

\begin{definition}\sl \label{a30}
Let $\Lambda$ be a Lagrangian manifold and let $0 \leq m ( h ) \leq h^{- N}$ for some $N > 0$. We say that $u ( x, h ) \in L^{2} ( \R^{n} )$ is a Lagrangian distribution of class $\CI ( \Lambda , m )$ if and only if $u$ is polynomially bounded and, for all $\rho \in T^{*} \R^{n}$, we have
\begin{equation}
u ( x ) = h^{- \frac{r}{2}} \int_{\R^{r}} e^{i \varphi ( x , \theta ) / h} a ( x , \theta , h ) \, d \theta.
\end{equation}
microlocally near $\rho$. Here, the symbol $a \in S ( m )$ is compactly supported in $x , \theta$ (uniformly with respect to $h$), and the function $\varphi$ is a non-degenerate phase function defined near the support of $a$ such that $\Lambda_{\varphi} \subset \Lambda$.
\end{definition}

In particular, $u \in \CI ( \Lambda , m )$ vanishes microlocally near each point $\rho$ not in $\Lambda$. When $\Lambda$ projects nicely on the $x$-space near $\rho$, the Lagrangian distributions can be written more simply. In this case, there exists a function $\psi ( x  ) \in C^{\infty} ( \R^{n} )$ defined near $\pi_{x} ( \rho )$ such that $\Lambda_{\psi} = \{ ( x , \xi ) ; \ \xi = \nabla \psi ( x ) \}$ coincides with $\Lambda$ near $\rho$. Then, for all $u \in \CI ( \Lambda , m )$, there exists $a ( x , h ) \in S ( m )$ supported near $\pi_{x} ( \rho )$ such that
\begin{equation}
u ( x ) = a ( x , h ) e^{i \psi ( x ) / h} ,
\end{equation}
microlocally near $\rho$. The Lagrangian distributions considered in this paper are often of this form.

Eventually, let $( \mu_{k} )_{k \geq 0}$ be the increasing sequence of linear combinations over $\N$ of the $\lambda_{j}$'s given by \ref{h2}. Following Helffer and Sj\"{o}strand \cite{HeSj86_01}, we say that a smooth function $u ( t , x ) : [ 0 , + \infty [ \times \R^{d} \longrightarrow \C$ is an {\it expandible function} when for all $K \in \N$, $\varepsilon > 0$ and $( \alpha , \beta ) \in \N^{1 + d}$, we have
\begin{equation} \label{m92}
\partial^{\alpha}_{t} \partial_{x}^{\beta} \Big( u ( t , x ) - \sum_{k = 1}^{K} u_{k} ( t , x ) e^{- \mu_{k} t} \Big) = \CO \big( e^{- ( \mu_{K + 1} - \varepsilon ) t} \big) ,
\end{equation}
for a sequence $( u_{k} )_{k \in \N}$ of smooth functions which are polynomials in $t$.

\Subsection{Propagation through a hyperbolic fixed point} \label{m81}

In this part, we recall the solution of the homogeneous microlocal Cauchy problem near the hyperbolic fixed point $( 0 , 0 )$. In other words, we give the asymptotic of the function $u$ satisfying
\begin{equation*}
\left\{ \begin{aligned}
&( P - z ) u = 0 &&\text{ microlocally near } ( 0 , 0 ) , \\
&u = u_{-} &&\text{ microlocally near } \Lambda^{0}_{-} ,
\end{aligned} \right.
\end{equation*}
where $u_{-}$ is a given initial data. The results stated here come from our previous work \cite{BoFuRaZe07_01}. In dimension $n = 1$, this problem has also been considered by Helffer and Sj\"{o}strand \cite[Appendix B]{HeSj89_01} in the analytic category and by Colin de Verdi\`ere and Parisse \cite{CoPa94_01,CoPa94_02} in the present $C^{\infty}$ category. Since we work microlocally near $( 0 , 0 )$, we only assume \ref{h2} in this section. The results are given for the Schr\"{o}dinger operators but hold true for general pseudodifferential operators mutatis mutandis.

The classical dynamic near the fixed point is described in Section \ref{s2}. We use here the notations and the statements of this section. Furthermore, let $( \widehat{\mu}_{k} )_{k \geq 0}$ denote the increasing sequence of linear combinations over $\N$ of the $( \mu_{j} - \mu_{1} )$'s with $j \geq 1$. In particular, $\mu_{0} = \widehat{\mu}_{0} = 0$, $\mu_{1} = \lambda_{1}$ and $\widehat{\mu}_{1} = \mu_{2} - \lambda_{1} > 0$.

Summarizing Theorems 2.1, 2.5 and 2.6 of \cite{BoFuRaZe07_01}, we get the result below. We send back the reader to this paper for more details, some remarks and the proofs. First, we define the exceptional set
\begin{equation} \label{m80}
\Gamma_{0} ( h ) = \bigg\{ E_{0} - i h \sum_{j = 1}^{n} \lambda_{j} \Big( \frac{1}{2} + \alpha_{j} \Big) ; \ \alpha \in \N^{n} \bigg\} ,
\end{equation}
which approximates the set of resonances generated by a barrier-top. There exists a discrete set $\Gamma ( h ) \subset \C$ containing $\Gamma_{0} ( h )$ with $\card \Gamma ( h ) \cap B ( E_{0} , C h )$ uniformly bounded for all $C > 0$ and $\Gamma ( h ) \subset \{ \im z < - \delta_{0} h \}$ for some $\delta_{0} > 0$, such that the following theorem holds.

\begin{figure}
\begin{center}
\begin{picture}(0,0)%
\includegraphics{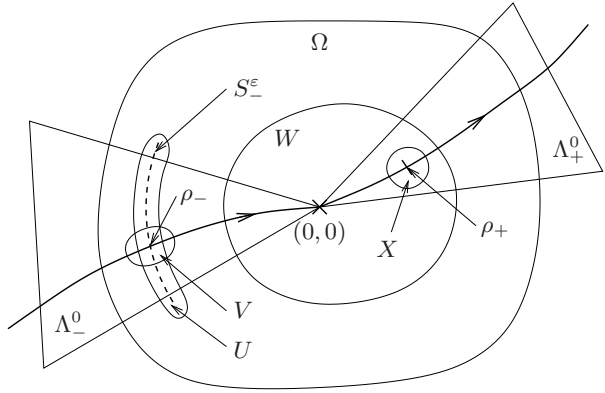}%
\end{picture}%
\setlength{\unitlength}{1184sp}%
\begingroup\makeatletter\ifx\SetFigFont\undefined%
\gdef\SetFigFont#1#2#3#4#5{%
  \reset@font\fontsize{#1}{#2pt}%
  \fontfamily{#3}\fontseries{#4}\fontshape{#5}%
  \selectfont}%
\fi\endgroup%
\begin{picture}(12505,8126)(-1157,-7790)
\put(10276,-2911){\makebox(0,0)[lb]{\smash{{\SetFigFont{9}{10.8}{\rmdefault}{\mddefault}{\updefault}$\Lambda_{+}^{0}$}}}}
\put(5401,-4636){\makebox(0,0)[b]{\smash{{\SetFigFont{9}{10.8}{\rmdefault}{\mddefault}{\updefault}$( 0 , 0 )$}}}}
\put(5401,-661){\makebox(0,0)[b]{\smash{{\SetFigFont{9}{10.8}{\rmdefault}{\mddefault}{\updefault}$\Omega$}}}}
\put(4726,-2686){\makebox(0,0)[b]{\smash{{\SetFigFont{9}{10.8}{\rmdefault}{\mddefault}{\updefault}$W$}}}}
\put(3601,-6286){\makebox(0,0)[lb]{\smash{{\SetFigFont{9}{10.8}{\rmdefault}{\mddefault}{\updefault}$V$}}}}
\put(3601,-7111){\makebox(0,0)[lb]{\smash{{\SetFigFont{9}{10.8}{\rmdefault}{\mddefault}{\updefault}$U$}}}}
\put(3601,-1636){\makebox(0,0)[lb]{\smash{{\SetFigFont{9}{10.8}{\rmdefault}{\mddefault}{\updefault}$S_{-}^{\varepsilon}$}}}}
\put(2476,-3811){\makebox(0,0)[lb]{\smash{{\SetFigFont{9}{10.8}{\rmdefault}{\mddefault}{\updefault}$\rho_{-}$}}}}
\put(8776,-4486){\makebox(0,0)[lb]{\smash{{\SetFigFont{9}{10.8}{\rmdefault}{\mddefault}{\updefault}$\rho_{+}$}}}}
\put(6826,-5011){\makebox(0,0)[b]{\smash{{\SetFigFont{9}{10.8}{\rmdefault}{\mddefault}{\updefault}$X$}}}}
\put(-149,-6586){\makebox(0,0)[lb]{\smash{{\SetFigFont{9}{10.8}{\rmdefault}{\mddefault}{\updefault}$\Lambda_{-}^{0}$}}}}
\end{picture}%
\end{center}
\caption{The geometric setting of Theorem \ref{a32}.} \label{f63}
\end{figure}

\begin{theorem}\sl \label{a32}
Let $\Omega \subset T^{*} \R^{n}$ be a small enough neighborhood of $( 0 ,0 )$ and
\begin{equation*}
S^{\varepsilon}_{\pm} = \{ ( x, \xi ) \in \Lambda^{0}_{\pm} ; \ \vert x \vert = \varepsilon \} \subset \Omega ,
\end{equation*}
with $\varepsilon > 0$ sufficiently small. Let also $C , \delta > 0$, $\rho_{-} = ( x_{-} , \xi_{-} ) \in S^{\varepsilon}_{-}$ with $g_{-} ( \rho_{-} ) \neq 0$ and $U$ (resp. $V$) be a small enough neighborhood of $S^{\varepsilon}_{-}$ (resp. $\rho_{-}$). Then, there exists an open neighborhood $W$ of $( 0 , 0 )$ containing $S^{\varepsilon}_{+}$ with the following property. For all $z \in B ( E_{0} ,  C h )$ with $\dist ( z , \Gamma ( h ) ) \geq \delta h$ and all $u_{-} \in L^{2} ( \R^{n} )$ with $\Vert u_{-} \Vert_{L^{2}} \leq 1$ such that
\begin{equation*}
\left\{ \begin{aligned}
&( P - z ) u_{-} =0 &&\text{ microlocally near } U , \\
&u_{-} = 0 &&\text{ microlocally near } U \setminus V ,
\end{aligned} \right.
\end{equation*}
the microlocal Cauchy problem
\begin{equation}
\left\{ \begin{aligned}
&( P - z ) u = 0 &&\text{ microlocally near } \Omega , \\
&u = u_{-} &&\text{ microlocally near } U ,
\end{aligned} \right.
\end{equation}
has a polynomially bounded solution $u$ which is unique microlocally near $W$ and satisfies
\begin{equation} \label{m75}
\Vert u \Vert_{L^{2}} \leq h^{- \frac{C}{\lambda_{1}} + \frac{\sum \lambda_{j}}{\lambda_{1}} - \frac{n}{2} - 2} .
\end{equation}

Moreover, if $\rho_{+} = ( x_{+} , \xi_{+} ) \in \Lambda^{0}_{+} \cap W$ satisfies $g_{-} ( \rho_{-} ) \cdot g_{+} ( \rho_{+} ) \neq 0$, then
\begin{equation} \label{a31}
u ( x ) = \CJ u_{-} ( x ) = h^{S ( z , h ) / \lambda_{1} - n / 2} \int_{H} e^{i ( \varphi_{+} ( x ) - \varphi_{-} ( y ) ) / h} d ( x, y , z , h ) u_{-} ( y ) \, d y ,
\end{equation}
microlocally near $X$, a neighborhood of $\rho_{+}$ independent of $u_{-}$. In this formula,
\begin{equation} \label{n13}
S ( z , h ) = \sum_{j = 1}^{n} \frac{\lambda_{j}}{2} - i \frac{z - E_{0}}{h} ,
\end{equation}
$H \subset \R^{n}$ is any hyperplane transversal to the base space projection of the Hamiltonian vector field restricted to $\Lambda^{0}_{-}$ at $x_{-}$. In the sequel, we assume that $H = \{ x \in \R^{n} ; \ x_{1} = c_{1} \}$ for some $c_{1} > 0$. This is always the case modulo a rotation of the coordinates. The symbol $d$ satisfies the asymptotic
\begin{equation} \label{b18}
d ( x , y , z , h ) \simeq \sum_{k \geq 0} d_{k} ( x , y , z , \ln h ) h^{\widehat{\mu}_{k} / \lambda_{1}} ,
\end{equation}
in $S ( 1 )$ where $d_{j}$ is polynomial in $\ln h$. The functions $d$ and $d_{j}$ are holomorphic functions of $\sigma = ( z - E_{0} ) /h$ for $z \in B ( E_{0} , C h )$ with $\dist ( z , \Gamma_{0} ( h ) ) \geq \delta h$. Furthermore, the principal symbol $d_{0}$ of $d$ is independent of $\ln h$ and given by
\begin{align}
d_{0} ( x , y , z ) ={}& ( 2 \pi )^{- \frac{n}{2}} \sqrt{\lambda_{1}} e^{-i n \frac{\pi}{4}} \Gamma \big( S ( z , h ) / \lambda_{1} \big) \big( i \lambda_{1} g_{-} ( \rho_{y}^{-} ) \cdot g_{+} ( \rho_{x}^{+} ) \big)^{- S ( z , h ) / \lambda_{1}}  \nonumber   \\
&\times \big\vert g_{-} ( \rho_{y}^{-} ) \big\vert \big\vert \det \nabla^{2}_{y^{\prime} , y^{\prime}} \varphi_{-} ( y ) \big\vert^{1 / 2} \big\vert \partial_{\xi_{1}} p ( \rho_{y}^{-} ) \big\vert^{1 / 2}        \label{b19} \\
&\times e^{\int_{0}^{- \infty} ( \Delta \varphi_{+} ( x ( s ) ) - \sum_{j} \lambda_{j} / 2 ) \, d s}
\lim_{t\to + \infty} \frac{e^{t ( \sum_{j} \lambda_{j} / 2 -\lambda_{1} )}}{{\sqrt{\Big\vert \det \frac{\partial y ( t , y^{\prime} , \eta^{\prime} )}{\partial ( t , y^{\prime} )}\vert_{\eta^{\prime} = \partial_{y^{\prime}} \varphi_{-} ( y )} \Big\vert}}} ,    \nonumber
\end{align}
where $\rho_{x}^{\pm} = ( x , \nabla \varphi_{\pm} ( x ) )$ and $x ( t )$ (resp. $y ( t , y^{\prime} , \eta^{\prime} )$) denotes the base space projection of the trajectory $\exp ( t H_{p} ) ( \rho_{x}^{+} )$ (resp. $\exp ( t H_{p} ) ( c_{1} , y^{\prime} , f_{-} ( c_{1} , y^{\prime} , \eta^{\prime} ) , \eta^{\prime} )$) with $y = ( c_{1} , y^{\prime} ) \in H$ and $f_{-} ( y , \eta^{\prime} ) = - \sqrt{E_{0} - \eta^{\prime} {}^{2} - V ( y )}$. In particular, $\rho^{\pm}_{x_{\pm}} = \rho_{\pm}$. The limits of the two last terms of \eqref{b19} exist and are real positive numbers. Eventually, we use the convention $( i a )^{b} = \vert a \vert^{b} e^{i \frac{\pi}{2} \sgn ( a ) b}$ for $a \in \R$, $b \in \C$.
\end{theorem}

The geometric setting is illustrated in Figure \ref{f63}. As explained in \cite[Remark 2.7]{BoFuRaZe07_01}, it is possible to replace the hyperplane $H$ by a small neighborhood of $x_{-}$ in \eqref{a31}. More precisely, we have
\begin{equation} \label{m76}
u ( x ) = h^{S ( z , h ) / \lambda_{1} - n / 2} \int_{\R^{n}} e^{i ( \varphi_{+} ( x ) - \varphi_{-} ( y ) ) / h} \widetilde{d} ( x , y , z , h ) u_{-} ( y ) \, d y ,
\end{equation}
microlocally near $X$. The symbol $\widetilde{d}$ verifies properties similar to those of $d$. Moreover, it is supported in the $y$ variables near $x_{-}$ (in any prescribed vicinity of $\pi_{x} ( V )$).

When the initial data $u_{-}$ is a Lagrangian distribution whose Lagrangian manifold is transverse to $\Lambda^{0}_{-}$, Theorem \ref{a32} takes the following form. This result follows from \cite[Section 6]{BoFuRaZe07_01} (one can also make a stationary phase expansion in \eqref{a31}).

\begin{corollary}\sl \label{d46}
In addition to the assumptions of Theorem \ref{a32}, let $\Lambda$ be a Lagrangian manifold that intersects $\Lambda^{0}_{-}$ along the Hamiltonian curve $\exp ( t H_{p} ) ( \rho_{-} )$ and such that $\Lambda$ projects nicely on the base space near $\exp ( t H_{p} ) ( \rho_{-} )$ for all $t \geq 0$. We suppose that
\begin{equation*}
u_{-} ( x ) = a_{-} ( x , h ) e^{i \psi ( x ) / h} ,
\end{equation*}
microlocally near $V$, with $\Lambda = \{ (x , \xi ) ; \ \xi = \nabla \psi (x) \}$ and $a_{-} \in S ( 1 )$. Then, 
\begin{equation} \label{m77}
u ( x ) = h^{S ( z , h ) / \lambda_{1} - 1 / 2} a_{+} ( x , h ) e^{i \varphi_{+} ( x ) / h} ,
\end{equation}
microlocally near $X$. The symbol $a_{+} ( x , h ) \in S ( 1 )$ can be written
\begin{equation} \label{m78}
a_{+} ( x , h ) = J ( x , h ) a_{-} ( x_{-} , h ) + S ( h ) ,
\end{equation}
with $J ( x , h ) = J_{0} ( x , h ) + S ( h^{\widehat{\mu}_{1} / 2 \lambda_{1}} )$ and
\begin{align}
J_{0} ( x , h ) ={}& e^{i A_{-} / h} e^{- i \frac{\pi}{4}} \sqrt{\frac{\lambda_{1}}{2 \pi}} \Gamma \big( S ( z , h ) / \lambda_{1} \big) \big\vert g_{-} ( \rho_{y}^{-} ) \big\vert \big( i \lambda_{1} g_{-} ( \rho_{y}^{-} ) \cdot g_{+} ( \rho_{x}^{+} ) \big)^{- S ( z , h ) / \lambda_{1}}  \nonumber \\
&\times \frac{D_{-} ( 0 )}{D_{+} ( 0 , x )} \lim_{t \to - \infty} \frac{D_{+} ( t , x )}{e^{t \sum_{j} \lambda_{j} / 2}} \lim_{t \to + \infty} \frac{e^{t ( \sum_{j} \lambda_{j} - 2 \lambda_{1} ) / 2}}{D_{-} ( t )} . \label{m79}
\end{align}
Here, $A_{-} = \psi ( x_{-} ) - \varphi_{-} ( x_{-} )$ and the Maslov determinants are given by
\begin{equation*}
D_{-} ( t ) = \sqrt{\Big\vert \det \frac{\partial x_{-} ( s , y )}{\partial ( s , y )} \vert_{s = t , \ y = 0} \Big\vert} \quad \text{and} \quad D_{+} ( t , x ) = \sqrt{\Big\vert \det \frac{\partial x_{+} ( s , x , y )}{\partial ( s , y )} \vert_{s = t , \ y = 0} \Big\vert} .
\end{equation*}
The function $( x_{-} ( t , y ) , \xi_{-} ( t , y ) ) : \R \times \R^{n - 1} \longrightarrow T^{*} \R^{n}$, defined near $\R \times \{ 0 \}$, is a smooth parametrization of $\Lambda$ by Hamiltonian curves such that $x_{-} ( 0 , 0 ) = x_{-}$ and $\partial_{( t , y )} x_{-} ( 0 , 0 )$ is invertible. The same way, $( x_{+} ( t , x , y ) , \xi_{+} ( t , x , y ) ) : \R \times \R^{n} \times \R^{n - 1} \longmapsto T^{*} \R^{n}$, defined near $\R \times \{ x_{+} \} \times \{ 0 \}$, is a smooth parametrization of $\Lambda^{0}_{+}$ by Hamiltonian curves such that $x_{+} ( 0 , x , 0 ) = x$ and $\partial_{( t , y )} x_{+} ( 0 , x , 0 )$ is invertible. The two limits in \eqref{m79} exist and provide real positive numbers.
\end{corollary}

That $\Lambda$ projects nicely on the base space near $\exp ( t H_{p} ) ( \rho_{-} )$ for $t$ large enough is in fact a consequence of the transversality of the intersection and of Proposition C.1 of \cite{ALBoRa08_01}. In the previous results, we have not described the solution $u$ near $( 0 , 0 )$ or in the orthogonal directions of $\rho_{-}$ (that is microlocally near $\rho_{+} \in \Lambda^{0}_{+}$ such that $g_{-} ( \rho_{-} ) \cdot g_{+} ( \rho_{+} ) = 0$). Some informations concerning these questions can be found in \cite{ALBoRa08_01,BoFuRaZe07_01}.

\section{Some properties of the Hamiltonian flow} \label{s10}

\Subsection{Consequences of the assumption \ref{h4}} \label{s11}

In this section, we give some properties of the classical flow. We begin with a perhaps standard result which, roughly speaking, says that the Hamiltonian curves in the energy surface $p^{-1} ( E_{0} )$ close to $(0,0)$ have an outgoing direction orthogonal to their incoming direction. For $\rho \in T^{*} \R^{n}$ and $t \in \R$, we will use in the sequel the notation $\rho (t) = \exp ( t H_{p} ) ( \rho )$.

\begin{lemma}\sl \label{a43}
Assume \ref{h2}. Let $U$ be a small open neighborhood of $(0,0)$ and $\rho_{\pm} \in \Lambda_{\pm}^{0} \cap U$ be such that
\begin{equation*}
g_{-} ( \rho_{-} ) \cdot g_{+} ( \rho_{+}) \neq 0 .
\end{equation*}
Then, there exist $U_{\pm} \subset U$, neighborhood of $\rho_{\pm}$, such that no piece of Hamiltonian curve in $U \cap p^{-1} ( E_{0} )$ starts in $U_{-}$ and ends in $U_{+}$.
\end{lemma}

In other words, there exits no Hamiltonian curve $\rho (t) \in p^{-1} ( E_{0} )$ such that $\rho ( t_{\pm} ) \in U_{\pm}$ and $\rho (t) \in U$ for all $t \in [ t_{-} , t_{+} ]$. The following remark gives a quantum proof of this lemma. Below, we give a proof at the classical level using some constructions and results of \cite{BoFuRaZe07_01}.

\begin{remark}\sl \label{a53}
Lemma \ref{a43} is in fact a corollary of Theorem \ref{a32}. Indeed, if this result did not hold, there would exist a piece of Hamiltonian curve in $U \cap p^{-1} ( E_{0} )$ which starts in $U_{-}$ and ends in $U_{+}$. Let $u$ be a solution of $( P - E_{0} ) u= 0$ and assume that $u$ only charges this piece of Hamiltonian curve in $U_{-}$. Then, by usual propagation of singularities along this curve, $u \neq 0$ microlocally near $U_{+}$. However, since $u =0$ microlocally near $\Lambda_{-} \cap U_{-}$, $u = 0$ microlocally near $U_{+}$ by Theorem \ref{a32}.
\end{remark}

\begin{proof}
Assume that Lemma \ref{a43} does not hold. Then, there exist $\rho_{\pm} \in \Lambda_{\pm}^{0} \cap U$ with
\begin{equation} \label{a54}
g_{-} ( \rho_{-} ) \cdot g_{+} ( \rho_{+}) \neq 0 ,
\end{equation}
a sequence of points $( \rho_{n} )_{n \in \N}$ and a sequence of positive times $(t_{n} )_{n \in \N}$ such that
\begin{equation} \label{a48}
\rho_{n} \longrightarrow \rho_{-} \qquad \text{and} \qquad \rho_{n} ( t_{n} ) \longrightarrow \rho_{+} ,
\end{equation}
as $n \to + \infty$ (see Figure \ref{f2}). Since the Hamiltonian flow is continuous and $\rho_{-} (t) \to 0$ as $t \to + \infty$, we necessarily have $t_{n} \to + \infty$ as $n \to + \infty$.

\begin{figure}
\begin{center}
\begin{picture}(0,0)%
\includegraphics{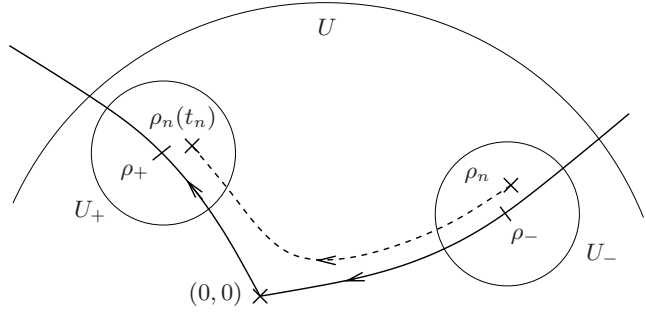}%
\end{picture}%
\setlength{\unitlength}{1184sp}%
\begingroup\makeatletter\ifx\SetFigFont\undefined%
\gdef\SetFigFont#1#2#3#4#5{%
  \reset@font\fontsize{#1}{#2pt}%
  \fontfamily{#3}\fontseries{#4}\fontshape{#5}%
  \selectfont}%
\fi\endgroup%
\begin{picture}(13323,6350)(718,-6544)
\put(7201,-886){\makebox(0,0)[lb]{\smash{{\SetFigFont{9}{10.8}{\rmdefault}{\mddefault}{\updefault}$U$}}}}
\put(4501,-6436){\makebox(0,0)[lb]{\smash{{\SetFigFont{9}{10.8}{\rmdefault}{\mddefault}{\updefault}$( 0 , 0 )$}}}}
\put(11251,-5161){\makebox(0,0)[lb]{\smash{{\SetFigFont{9}{10.8}{\rmdefault}{\mddefault}{\updefault}$\rho_{-}$}}}}
\put(10276,-3886){\makebox(0,0)[lb]{\smash{{\SetFigFont{9}{10.8}{\rmdefault}{\mddefault}{\updefault}$\rho_{n}$}}}}
\put(3676,-2761){\makebox(0,0)[lb]{\smash{{\SetFigFont{9}{10.8}{\rmdefault}{\mddefault}{\updefault}$\rho_{n} ( t_{n} )$}}}}
\put(12826,-5611){\makebox(0,0)[lb]{\smash{{\SetFigFont{9}{10.8}{\rmdefault}{\mddefault}{\updefault}$U_{-}$}}}}
\put(3076,-3811){\makebox(0,0)[lb]{\smash{{\SetFigFont{9}{10.8}{\rmdefault}{\mddefault}{\updefault}$\rho_{+}$}}}}
\put(2101,-4711){\makebox(0,0)[lb]{\smash{{\SetFigFont{9}{10.8}{\rmdefault}{\mddefault}{\updefault}$U_{+}$
}}}}
\end{picture}%
\end{center}
\caption{The geometric setting of Lemma \ref{a43} and its proof.} \label{f2}
\end{figure}

We now use the notations of \cite[Section 5]{BoFuRaZe07_01} and use that $g_{-} ( \rho_{-} ) \neq 0$. Setting $\rho_{\pm} = ( x_{\pm} , \xi_{\pm} )$, $\rho_{n} = ( x_{n} , \xi_{n} )$ and using the notation $y = (y_{1} , y^{\prime} )$ with $y_{1} \in \R$ and $y^{\prime} \in \R^{n-1}$ for $y \in \R^{n}$, we can suppose that $\rho_{-} \in H_{-} = \{ x_{1} = \varepsilon \}$ and that the vector field $H_{p}$ is transverse to $H_{-} \times \R^{n}$ at $\rho_{-}$ (see \cite[Page 98]{BoFuRaZe07_01}). Then, there exists $\eta_{n}^{\prime} \in \R^{n-1}$ with $\eta_{n}^{\prime} \to \xi_{-}^{\prime}$ as $n \to + \infty$ such that (see \cite[Lemma 5.5]{BoFuRaZe07_01})
\begin{equation*}
\rho_{n} \in \Lambda_{\psi_{\eta_{n}^{\prime}}} = \big\{ ( x , \xi ) \in T^{*} \R^{n} ; \ x \text{ near } x_{-} \text{ and } \xi = \nabla \psi_{\eta_{n}^{\prime}} (x) \big\} ,
\end{equation*}
where
\begin{equation} \label{a49}
\left\{ \begin{aligned}
&p ( x , \nabla \psi_{\eta^{\prime}} ) = E_{0} ,  \\
&\psi_{\eta^{\prime}} (x) = x^{\prime} \cdot \eta^{\prime} \text{ for } x \in H_{-} .
\end{aligned} \right.
\end{equation}
As in \cite[(5.22)]{BoFuRaZe07_01}, we note
\begin{equation*}
\Gamma_{0}^{\eta^{\prime}} = \big\{ ( x , \xi ) \in \Lambda_{\psi_{\eta^{\prime}}} ; \ \psi_{\eta^{\prime}} (x) = \psi_{\eta^{\prime}} ( x ( \eta^{\prime} ) ) \big \} ,
\end{equation*}
where $\rho_{\eta^{\prime}} = (x ( \eta^{\prime} ) , \nabla \psi_{\eta^{\prime}} ( x ( \eta^{\prime} ) ) )$ is the unique point in the intersection $\Lambda_{-} \cap \Lambda_{\psi_{\eta^{\prime}}} \cap H_{-}$. Using \eqref{a49} and the implicit function theorem, there exists $\varepsilon_{n} \to 0$ as $n \to + \infty$ such that
\begin{equation} \label{a55}
\rho_{n} ( \varepsilon_{n} ) \in \Gamma_{0}^{\eta^{\prime}_{n}} .
\end{equation}
More precisely, let
\begin{equation*}
f (t , \rho ) = \psi_{\eta^{\prime} ( \rho )} ( x ( t , \rho ) ) - \psi_{\eta^{\prime} ( \rho )} ( x ( \eta^{\prime} ( \rho ))) ,
\end{equation*}
where $x (t , \rho )$ is the base space projection of $\rho (t)$ and $\eta^{\prime} ( \rho )$ is such that $\rho \in \Lambda_{\psi_{\eta^{\prime} ( \rho )}}$. Then, $f (t , \rho )$ is $C^{\infty}$ near $(0 , \rho_{-} )$. Moreover, we have $f ( 0 , \rho_{-} ) = 0$ and
\begin{equation*}
\partial_{t} f (0 , \rho_{-} ) = \partial_{t} x ( 0 , \rho_{-} ) \cdot \nabla \psi_{\eta^{\prime} ( \rho_{-} )} ( x ( 0 , \rho_{-} ) ) = 2 \xi ( 0 , \rho_{-} ) \cdot \xi ( 0 , \rho_{-} ) = 2 \vert \xi_{-} \vert^{2} \neq 0 .
\end{equation*}
Thus, by the implicit function theorem, there exists a $C^{\infty}$ function $t ( \rho )$ with $t ( \rho_{-} ) = 0$ such that $f ( t ( \rho ) , \rho ) = 0$. In particular, \eqref{a55} holds with $\varepsilon_{n} = t ( \rho_{n} )$.

Then, by \cite[(5.23) and (5.38)]{BoFuRaZe07_01}, $\rho_{n} (t + \varepsilon_{n} ) \in \Lambda_{t}^{\eta^{\prime}_{n}}$ where, from \cite[(5.44)--(5.45)]{BoFuRaZe07_01},
\begin{equation} \label{a50}
\Lambda_{t}^{\eta^{\prime}} = \big\{ (x , \xi ) ; \ \xi = \nabla_{x} \varphi (t ,x ,\eta^{\prime} ) \big\} \qquad \text{with} \qquad \partial_{t} \varphi + p ( x , \nabla_{x} \varphi ) = E_{0} ,
\end{equation}
and, as expandible functions, (see \cite[Lemma 5.10]{BoFuRaZe07_01})
\begin{equation} \label{a51}
\varphi ( t , x , \eta^{\prime} ) \simeq \varphi_{+} ( x ) + \widetilde{\psi} ( \eta^{\prime} ) + \sum_{j \geq 1} e^{- \mu_{j} t} \varphi_{j} (t ,x ,\eta^{\prime} ) ,
\end{equation}
where the $\varphi_{j} (t ,x ,\eta^{\prime} )$ are polynomials in $t$. Moreover, $\varphi_{1}$ does not depend on $t$ and satisfies (see \cite[(7.29)]{ALBoRa08_01})
\begin{equation} \label{a52}
\varphi_{1} (x , \eta^{\prime} ) = - \lambda_{1} g_{-} ( \rho_{\eta^{\prime}} ) \cdot g_{+} \big( (x , \nabla \varphi_{+} (x) ) \big) .
\end{equation}
Using $p ( \rho_{n} ( t_{n} ) ) = E_{0}$ and \eqref{a50}, we get $\partial_{t} \varphi ( t_{n} - \varepsilon_{n} , x_{n} ( t_{n} ) , \eta_{n}^{\prime} ) =0$. Since $\varepsilon_{n} \to 0$, \eqref{a51} and \eqref{a52} imply
\begin{equation*}
- \lambda_{1} g_{-} ( \rho_{\eta^{\prime}_{n}} ) \cdot g_{+} \big( ( x_{n} ( t_{n} ) , \nabla \varphi_{+} ( x_{n} ( t_{n} ) ) ) \big) = \varphi_{1} ( x_{n} ( t_{n} ) , \eta_{n}^{\prime} ) = \CO \big( e^{- \delta t_{n}} \big) ,
\end{equation*}
for some $\delta > 0$. Taking the limit $n \to + \infty$ in the previous equation and using the continuity of $\rho \mapsto g_{\pm} ( \rho )$, and the limits $\rho_{\eta^{\prime}_{n}} \to \rho_{-}$ and $\rho_{n} ( t_{n} ) \to \rho_{+}$, we eventually obtain
\begin{equation*}
g_{-} ( \rho_{-} ) \cdot g_{+} ( \rho_{+} ) = 0 ,
\end{equation*}
which is in contradiction with \eqref{a54}.
\end{proof}

We now recall the Hartman--Grobman theorem (see \cite[Page 120]{Pe01_01}). There exists a homeomorphism $F$ from $U$ onto $V$, two open neighborhoods of $(0,0)$, with the following property. For all $\rho \in U$, there is an open interval $0 \in I \subset \R$ containing zero such that
\begin{equation} \label{a45}
\forall t \in I , \qquad F \big( \exp ( t H_{p} ) ( \rho ) \big) = e^{t L} F ( \rho ) ,
\end{equation}
where $L = \diag ( \lambda_{1} , \ldots , \lambda_{n} , - \lambda_{1} , \ldots , - \lambda_{n} )$.

\begin{proposition}\sl \label{a44}
Assume \ref{h1}--\ref{h4}. There exists $\varepsilon > 0$ such that
\begin{equation*}
\rho \in \CH \text{ and } \vert \rho \vert \leq \varepsilon \quad \Longrightarrow \quad \rho \in \Lambda_{-}^{0} \cup \Lambda_{+}^{0} .
\end{equation*}
\end{proposition}

\begin{proof}
We prove this result using a contradiction argument. Let $\delta > 0$ be small enough such that $\overline{B ( 0 , \delta )} \subset V$ and such that Lemma \ref{a43} is true in a vicinity of $W = F^{-1} ( \overline{B ( 0 , \delta )} )$. Assume that Proposition \ref{a44} does not hold. Then, there exists $( \rho_{n} )_{n \in \N} \in \CH \setminus ( \Lambda_{-}^{0} \cup \Lambda_{+}^{0} )$ with $\rho_{n} \to 0$ as $n \to + \infty$. For $n$ large enough, we have $\rho_{n} \in W$ and we define
\begin{equation*}
\pi_{n} (t) = ( y_{n} (t) , \eta_{n} (t) ) = F ( \rho_{n} (t) ) = e^{t L} F ( \rho_{n} ) = e^{t L} \pi_{n} = e^{t L} ( y_{n} , \eta_{n} ) .
\end{equation*}
In particular, $\pi_{n} \to 0$ as $n \to + \infty$.

\begin{figure}
\begin{center}
\begin{picture}(0,0)%
\includegraphics{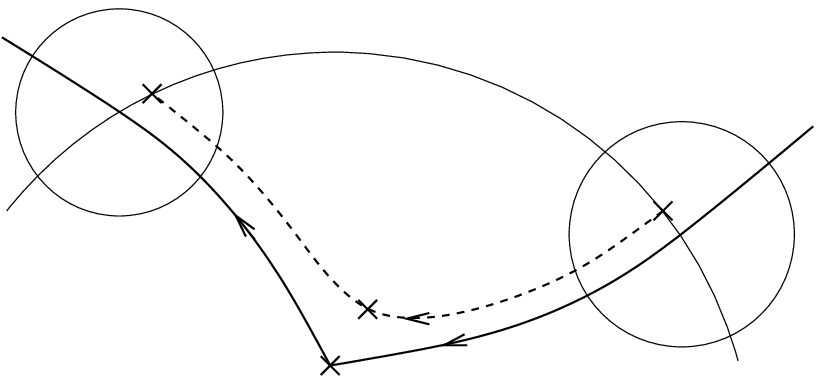}%
\end{picture}%
\setlength{\unitlength}{1184sp}%
\begingroup\makeatletter\ifx\SetFigFont\undefined%
\gdef\SetFigFont#1#2#3#4#5{%
  \reset@font\fontsize{#1}{#2pt}%
  \fontfamily{#3}\fontseries{#4}\fontshape{#5}%
  \selectfont}%
\fi\endgroup%
\begin{picture}(13041,5905)(718,-6544)
\put(2326,-4486){\makebox(0,0)[lb]{\smash{{\SetFigFont{9}{10.8}{\rmdefault}{\mddefault}{\updefault}$U_{+}$
}}}}
\put(11926,-4336){\makebox(0,0)[lb]{\smash{{\SetFigFont{9}{10.8}{\rmdefault}{\mddefault}{\updefault}$\rho_{-}$}}}}
\put(2401,-2761){\makebox(0,0)[lb]{\smash{{\SetFigFont{9}{10.8}{\rmdefault}{\mddefault}{\updefault}$\rho_{+}$}}}}
\put(6676,-5086){\makebox(0,0)[lb]{\smash{{\SetFigFont{9}{10.8}{\rmdefault}{\mddefault}{\updefault}$\rho_{n}$}}}}
\put(11176,-3286){\makebox(0,0)[lb]{\smash{{\SetFigFont{9}{10.8}{\rmdefault}{\mddefault}{\updefault}$\rho_{n} ( t_{n}^{-} )$}}}}
\put(2026,-1636){\makebox(0,0)[lb]{\smash{{\SetFigFont{9}{10.8}{\rmdefault}{\mddefault}{\updefault}$\rho_{n} ( t_{n}^{+} )$}}}}
\put(6751,-2161){\makebox(0,0)[lb]{\smash{{\SetFigFont{9}{10.8}{\rmdefault}{\mddefault}{\updefault}$\partial W$}}}}
\put(4501,-6436){\makebox(0,0)[lb]{\smash{{\SetFigFont{9}{10.8}{\rmdefault}{\mddefault}{\updefault}$( 0 , 0 )$}}}}
\put(11476,-2236){\makebox(0,0)[lb]{\smash{{\SetFigFont{9}{10.8}{\rmdefault}{\mddefault}{\updefault}$U_{-}$}}}}
\end{picture}%
\end{center}
\caption{The geometric setting of the proof of Proposition \ref{a44}.} \label{f1}
\end{figure}

Since $\pi_{n} \notin F ( \Lambda_{-}^{0} )$ (resp. $\pi_{n} \notin F ( \Lambda_{+}^{0} )$), $\pi_{n} (t)$ leaves $\overline{B ( 0 , \delta )}$ at some time $t^{+}_{n} > 0$ (resp. $t^{-}_{n} < 0$). At this positive time, we have
\begin{equation*}
\vert \pi_{n} ( t_{n}^{+} ) \vert = \delta \qquad \text{and} \qquad \vert \eta_{n} ( t_{n}^{+} ) \vert = \big\vert e^{- \lambda t_{n}^{+}} \eta_{n} \big\vert \leq \vert \eta_{n} \vert ,
\end{equation*}
from \eqref{a45}. Then, $\eta_{n} ( t_{n}^{+} ) \to 0$ as $n \to + \infty$. The same way, $\vert \pi_{n} ( t_{n}^{-} ) \vert = \delta$ and $\vert y_{n} ( t_{n}^{-} ) \vert \to 0$ as $n \to + \infty$. Coming back to the original variables, this gives
\begin{equation} \label{a46}
\rho_{n} ( t_{n}^{\pm} ) \in \CH \cap \partial W \qquad \text{and} \qquad \dist ( \rho_{n} ( t_{n}^{\pm} ) , \Lambda_{\pm}^{0} ) \longrightarrow 0 ,
\end{equation}
as $n \to + \infty$ (see Figure \ref{f1}).

Using that $W$ is bounded, we can always assume that, up to the extraction of a subsequence, $\rho_{n} ( t_{n}^{\pm} ) \to \rho_{\pm}$ as $n \to + \infty$. Since $\partial W \subset T^{*} \R^{n} \setminus \{ ( 0 , 0 ) \}$, $K ( E_{0} )$ and $\Lambda_{\pm}^{0} \cap W$ are closed sets, \eqref{a46} yields
\begin{equation} \label{a47}
\rho_{\pm} \in \Lambda_{\pm}^{0} \cap \CH \cap W \qquad \text{and} \qquad \rho_{\pm} \neq ( 0 , 0 ).
\end{equation}
Then, \ref{h4} and Lemma \ref{a43} provide $U_{\pm}$, neighborhoods of $\rho_{\pm}$, such that no piece of Hamiltonian curve in $W \cap p^{-1} ( E_{0} )$ starts in $U_{-}$ and ends in $U_{+}$. But, for $n$ large enough, $\rho_{n} ( t_{n}^{\pm} ) \in U_{\pm}$ and $\rho_{n} ( t )\in W$ for all $t \in [ t_{n}^{-} , t_{n}^{+} ]$. This gives a contradiction.
\end{proof}

\Subsection{Proof of the assertions appearing before Theorem \ref{a2} and Theorem \ref{i55}} \label{a71}

First, we obtain some properties on the Hamiltonian flow on $\Lambda_{+}^{0}$. Let
\begin{equation*}
H_{p}^{+} = 2 \nabla \varphi_{+} (x) \cdot \partial_{x} ,
\end{equation*}
be the Hamiltonian vector field restricted to $\Lambda_{+}^{0}$. This means that $( x (t) , \xi (t) )$ is a Hamiltonian curve in $\Lambda_{+}^{0}$ if and only if $x (t)$ is an integral curve of $H_{p}^{+}$ and $\xi (t) = \nabla \varphi_{+} ( x (t) )$. Note that \eqref{a57} implies that $2 \nabla \varphi_{+} (x) = \lambda x + \CO ( x^{2} )$. Then, from \cite[Remark 3.10]{HeSj85_01}, we get that
\begin{equation} \label{b25}
x^{+} ( t , x ) : = \exp ( t H_{p}^{+} ) (x) ,
\end{equation}
is expandible (as $t \to - \infty$). We write
\begin{equation} \label{a58}
x^{+} ( t , x ) = e^{\lambda t} x_{1}^{+} ( t , x ) + \widetilde{\CO} ( e^{2 \lambda t} ) ,
\end{equation}
where $x_{1}^{+}$ is polynomial in $t$ and $C^{\infty}$ in $x$. Moreover, from \cite[Proposition 6.11]{ALBoRa08_01}, we know that $x_{1}^{+}$ does not depend on $t$. In the previous equation, $f ( t , x ) = \widetilde{\CO} ( m ( t ) )$ means that, for all $( \alpha , \beta ) \in \N^{1 + n}$ and $\varepsilon > 0$, we have
\begin{equation*}
\big\vert \partial_{x}^{\alpha} \partial_{x}^{\beta} f ( t , x ) \big\vert \lesssim e^{\varepsilon t } m ( t ) ,
\end{equation*}
uniformly with respect to $t , x$. On the other hand, since all the $\lambda_{j}$'s are equal, there exists $\Phi$, a $C^{\infty}$ diffeomorphism between two neighborhoods of $0 \in \R^{n}$, such that
\begin{equation} \label{a62}
\Phi \big( \exp ( t H_{p}^{+} ) (x) \big) = e^{\lambda t} \Phi (x) ,
\end{equation}
$\Phi (0) = 0$ and $d \Phi (0) = Id_{\R^{n}}$ (see e.g. \cite[Theorem IX.12.1]{Ha82_01}). Then, \eqref{a58} gives
\begin{equation*}
e^{\lambda t} \Phi (x) = \Phi \big( e^{\lambda t} x_{1}^{+} ( x ) + \widetilde{\CO} ( e^{2 \lambda t} ) \big) = e^{\lambda t} x_{1}^{+} ( x ) + \widetilde{\CO} ( e^{2 \lambda t} ) .
\end{equation*}
This implies that $\Phi = x_{1}^{+}$. Moreover, from the definition of $g_{+}$, we have
\begin{equation} \label{a60}
\Phi (x) = g_{+} ( x , \nabla \varphi_{+} (x) ) .
\end{equation}

We will now use $x^{+}$ to prove the assertions appearing before Theorem \ref{a2}. For $\varepsilon > 0$ small enough, all the curves of $\CH$ pass through $\{ \vert x \vert = \varepsilon \} \cap \Lambda_{+}^{0}$. Then, using \eqref{a60},
\begin{equation*}
\CH_{\rm tang}^{\pm \infty} = \Big\{ \frac{\Phi (x)}{\vert \Phi (x) \vert} ; \ x \in K \Big\}  \quad \text{with} \quad K = \big\{ x \in \R^{n} ; \ \vert x \vert = \varepsilon \text{ and } ( x, \nabla \varphi_{+} (x) ) \in \CH_{\rm tang} \big\} .
\end{equation*}
Now, since $\Phi$ is a local diffeomorphism near $0$, $\Phi$ does not vanish on $K$ which is a compact set. This implies the following result.

\begin{proposition}\sl \label{a61}
The set $\CH_{\rm tang}^{\pm \infty}$ is a compact subset of $\S^{n-1}$.
\end{proposition}

Let $\varepsilon > 0$ be small enough. For $t \in \R$ and $\alpha \in \S^{n-1}$, we set
\begin{equation} \label{a73}
x_{+} ( t , \alpha ) = x^{+} \Big( t - \frac{\ln \varepsilon}{\lambda} , \Phi^{-1} ( \varepsilon \alpha ) \Big) ,
\end{equation}
and $\rho_{+} ( t , \alpha ) = ( x_{+} ( t , \alpha ) , \xi_{+} ( t , \alpha ) ) \in T^{*} \R^{n}$ the corresponding Hamiltonian curve. In fact, $x_{+} ( t , \alpha )$ does not depend on $\varepsilon$.

\begin{proposition}\sl \label{a59}
The function $x_{+} ( t , \alpha ) : \R \times \S^{n-1} \to \R^{n}$ is expandible (as $t \to - \infty$). Moreover, we have
\begin{equation} \label{a75}
x_{+} ( t , \alpha ) = e^{\lambda t} \alpha + \widetilde{\CO} ( e^{2 \lambda t} ) .
\end{equation}
Eventually, the application
\begin{equation*}
\begin{aligned}
\Psi : \\
{}^{}
\end{aligned}
\left\{ \begin{aligned}
&\R \times \S^{n - 1} &\longrightarrow &&&\Lambda_{+}^{0} \setminus \{ 0 \}  \\
& ( t , \alpha ) &\longmapsto &&&\big( x_{+} ( t , \alpha ) , \nabla \varphi_{+} ( x _{+} ( t , \alpha ) ) \big)
\end{aligned} \right.
\end{equation*}
is a bijection from a neighborhood of $\{ - \infty \} \times \S^{n-1}$ to a pointed neighborhood of $0$.
\end{proposition}

\begin{proof}
The fact that $x_{+}$ is expandible follows directly from the corresponding property of $x^{+}$. Moreover, from $\Phi = x_{1}^{+}$ and \eqref{a58}, we have
\begin{align*}
x_{+} ( t , \alpha ) &= x^{+} \Big( t - \frac{\ln \varepsilon}{\lambda} , \Phi^{-1} ( \varepsilon \alpha ) \Big)   \\
&= e^{\lambda ( t - \frac{\ln \varepsilon}{\lambda} )} x_{1}^{+} \big( \Phi^{-1} ( \varepsilon \alpha ) \big) + \widetilde{\CO} ( e^{2 \lambda t} )  \\
&= e^{\lambda t} \alpha + \widetilde{\CO} ( e^{2 \lambda t} ) .
\end{align*}
The fact that $\Psi$ is injective follows from the previous equation. On the other hand, for all Hamiltonian curve $( x (t) , \xi (t) )$ in $\Lambda^{0}_{+}$, \eqref{a62} yields
\begin{equation} \label{a63}
x (t) = x^{+} ( t , x (0) ) = x_{+} \Big( t + \frac{\ln \vert \Phi ( x (0) ) \vert}{\lambda} , \frac{\Phi ( x (0) )}{\vert \Phi ( x (0) ) \vert} \Big) ,
\end{equation}
and the surjectivity follows.
\end{proof}

Let $\beta \in \S^{n - 1}$ and $0 < \varepsilon \leq \varepsilon_{0}$, with $\varepsilon_{0} > 0$ small enough. As in \eqref{a63}, we have
\begin{equation} \label{a66}
x^{+} ( t , \varepsilon \beta ) = x_{+} \Big( t + \frac{\ln \vert \Phi ( \varepsilon \beta ) \vert}{\lambda} , \frac{\Phi ( \varepsilon \beta )}{\vert \Phi ( \varepsilon \beta ) \vert} \Big) .
\end{equation}
Since $\Phi$ is smooth near $0$, we can write
\begin{equation*}
\Phi ( \varepsilon \beta ) = \varepsilon \beta + \varepsilon^{2} f ( \varepsilon , \beta ) ,
\end{equation*}
with $f \in C^{\infty} ( [ - \varepsilon_{0} , \varepsilon_{0} ] \times \S^{n - 1} )$. In particular, for $( \varepsilon , \beta ) \in ] 0 , \varepsilon_{0} ] \times \S^{n - 1}$,
\begin{align}
\vert \Phi ( \varepsilon \beta ) \vert &= \vert \varepsilon \beta + \varepsilon^{2} f ( \varepsilon , \beta ) \vert  \nonumber \\
&= \varepsilon \sqrt{\big( \beta + \varepsilon f ( \varepsilon , \beta ) \big) \cdot \big( \beta + \varepsilon f ( \varepsilon , \beta ) \big)}   \nonumber \\
&= \varepsilon \big( 1 + \varepsilon g ( \varepsilon , \beta ) \big) ,  \label{a67}
\end{align}
with $g \in C^{\infty} ( [ - \varepsilon_{0} , \varepsilon_{0} ] \times \S^{n - 1} )$. This implies that
\begin{equation} \label{a64}
F ( \varepsilon , \beta ) : = \frac{\Phi ( \varepsilon \beta )}{\vert \Phi ( \varepsilon \beta ) \vert} = \beta + \varepsilon k ( \varepsilon , \beta ) ,
\end{equation}
with $k \in C^{\infty} ( [ - \varepsilon_{0} , \varepsilon_{0} ] \times \S^{n - 1} )$. Note that $F ( \varepsilon , \beta )$ is the normalized asymptotic direction of the Hamiltonian curve of $\Lambda_{+}^{0}$ passing through $\varepsilon \beta$. From \eqref{a64}, $\beta \longmapsto F ( \varepsilon , \beta )$ is a local diffeomorphism and is injective on $\S^{n - 1}$. Then, its image is closed (since it is the image of a compact set by a continuous function) and open (since the application is a local diffeomorphism) of the connected set $\S^{n - 1}$ if $n \geq 2$. Therefore, for all $\varepsilon \in [ - \varepsilon_{0} , \varepsilon_{0} ]$,
\begin{equation*}
\beta \longmapsto F ( \varepsilon , \beta ) ,
\end{equation*}
is a diffeomorphism from $\S^{n - 1}$ onto itself whose inverse $F^{-1} ( \varepsilon , \beta )$ satisfies, from \eqref{a64},
\begin{equation} \label{a68}
F^{-1} ( \varepsilon , \beta ) = \beta + \varepsilon \ell ( \varepsilon , \beta ) .
\end{equation}
with $\ell \in C^{\infty} ( [ - \varepsilon_{0} , \varepsilon_{0} ] \times \S^{n - 1} )$.

\begin{proposition}\sl \label{a65}
For all $\alpha \in \S^{n - 1}$ and $0 < \varepsilon \leq \varepsilon_{0}$, the characteristic curve $x_{+} ( t , \alpha )$ intersects $\{ \vert x \vert = \varepsilon \}$ at a unique negative time $t_{+}^{\varepsilon} ( \alpha )$ which verifies
\begin{equation*}
t_{+}^{\varepsilon} ( \alpha ) = \frac{\ln \varepsilon}{\lambda} + \varepsilon m ( \varepsilon , \alpha ) ,
\end{equation*}
with $m \in C^{\infty} ( [ - \varepsilon_{0} , \varepsilon_{0} ] \times \S^{n - 1} )$.
\end{proposition}

\begin{proof}
From \eqref{a66}, $x_{+} ( t , \alpha )$ belongs to $\{ \vert x \vert = \varepsilon \}$ with $t < 0$ if and only if there exists $\beta \in \S^{n - 1}$ such that
\begin{equation*}
\alpha = \frac{\Phi ( \varepsilon \beta )}{\vert \Phi ( \varepsilon \beta ) \vert}  \qquad \text{and} \qquad t = \frac{\ln \vert \Phi ( \varepsilon \beta ) \vert}{\lambda} ,
\end{equation*}
which is equivalent to
\begin{equation*}
\beta = F^{-1} ( \varepsilon , \alpha ) \qquad \text{and} \qquad t = t_{+}^{\varepsilon} ( \alpha ) : = \frac{\ln \vert \Phi ( \varepsilon F^{-1} ( \varepsilon , \alpha )) \vert}{\lambda} .
\end{equation*}
So, we have proved that the curve $x_{+} ( t , \alpha )$ meets the sphere of radius $\varepsilon$ at a unique negative time $t_{+}^{\varepsilon} ( \alpha )$ given by the previous formula. In particular, $( \varepsilon , \alpha ) \longmapsto t_{+}^{\varepsilon} ( \alpha ) \in C^{\infty} ( ] 0 , \varepsilon_{0} ] \times \S^{n - 1} )$. Moreover, \eqref{a67} and \eqref{a68} imply
\begin{equation} \label{a70}
t_{+}^{\varepsilon} ( \alpha ) = \frac{\ln \varepsilon}{\lambda} + \varepsilon m ( \varepsilon , \alpha ) ,
\end{equation}
with $m \in C^{\infty} ( [ - \varepsilon_{0} , \varepsilon_{0} ] \times \S^{n - 1} )$.
\end{proof}

In particular, the previous proposition directly implies that $( \varepsilon , \alpha ) \longmapsto t_{+}^{\varepsilon} ( \alpha )$ is continuous on $] 0 , \varepsilon_{0} ] \times \CH_{\rm tang}^{+ \infty}$. We will now prove that the same thing holds for $t_{-}^{\varepsilon} ( \alpha )$. As before, we can parametrize $\Lambda_{-}^{0}$ by $x_{-} ( t , \omega )$, for $t \in \R$ and $\omega \in \S^{n - 1}$, which satisfies mutatis mutandis \eqref{a73} and the properties stated in Proposition \ref{a59}. Moreover, since we consider a Schr\"{o}dinger operator, we have $x_{-} ( t , \omega ) = x_{+} ( - t , \omega )$. Then, Proposition \ref{a65} implies that, for all $\omega \in \S^{n - 1}$ and $0 < \varepsilon \leq \varepsilon_{0}$, the characteristic curve $x_{-} ( t , \omega )$ intersects $\{ \vert x \vert = \varepsilon \}$ at a unique positive time, namely at $t = - t_{+}^{\varepsilon} ( \omega )$. In particular, for all $\alpha \in \CH_{\rm tang}^{+ \infty}$ and $0 < \varepsilon < \varepsilon_{0}$, the characteristic curve $x_{+} ( t , \alpha )$ intersects $\{ \vert x \vert = \varepsilon \}$ at a unique positive time $t_{-}^{\varepsilon} ( \alpha )$.

Now, consider $\alpha_{0} \in \CH_{\rm tang}^{+ \infty}$. Since $x_{+}( t , \alpha_{0} )$ goes to $0$ as $t$ tends to $+ \infty$, there exists a positive time $T_{\alpha_{0}}$ such that $x_{+} ( T_{\alpha_{0}} , \alpha_{0} ) \in \{ \vert x \vert \leq \varepsilon_{0} / 4 \}$. Then, by continuity of the Hamiltonian flow, for all $\alpha$ in a neighborhood in $\S^{n - 1}$ of $\alpha_{0}$, we have  $x_{+} ( T_{\alpha_{0}} , \alpha ) \in \{ \vert x \vert \leq \varepsilon_{0} / 2 \}$. By Proposition \ref{a44} and the compactness of $\CH_{\rm tang}^{+ \infty}$, it follows that there exists $T > 0$ such that
\begin{equation*}
x_{+} ( T , \alpha ) \in \{ \vert x \vert \leq \varepsilon_{0} \} ,
\end{equation*}
for all $\alpha \in \CH_{\rm tang}^{+ \infty}$. From \eqref{a63}, we then have
\begin{align}
x_{+} ( t , \alpha ) &= x^{+} ( - t + T ,x_{+} ( T , \alpha ) )    \nonumber \\
&= x_{+} \Big( - t + T + \frac{\ln \vert \Phi ( x_{+} ( T , \alpha ) ) \vert}{\lambda} , \frac{\Phi ( x_{+} ( T , \alpha ) )}{\vert \Phi ( x_{+} ( T , \alpha ) ) \vert} \Big) .  \label{a80}
\end{align}
Hence, Proposition \ref{a65} yields
\begin{equation*}
t_{-}^{\varepsilon} ( \alpha ) = T + \frac{\ln \vert \Phi ( x_{+} ( T , \alpha ) ) \vert}{\lambda} - t_{+}^{\varepsilon} \Big( \frac{\Phi ( x_{+} ( T , \alpha ) )}{\vert \Phi ( x_{+} ( T , \alpha ) ) \vert} \Big) .
\end{equation*}
Using the continuity of $\alpha \longmapsto x_{+} ( T , \alpha )$ and Proposition \ref{a65}, we obtain the following proposition.

\begin{proposition}\sl \label{a74}
For all $\alpha \in \CH_{\rm tang}^{+ \infty}$ and $0 < \varepsilon \leq \varepsilon_{0}$, the characteristic curve $x_{+} ( t , \alpha )$ intersects $\{ \vert x \vert = \varepsilon \}$ at a unique positive time $t_{-}^{\varepsilon} ( \alpha )$ which verifies
\begin{equation*}
t_{-}^{\varepsilon} ( \alpha ) = - \frac{\ln \varepsilon}{\lambda} + q ( \alpha ) + \varepsilon r ( \varepsilon , \alpha ) ,
\end{equation*}
with $q \in C^{0} ( \CH_{\rm tang}^{+ \infty} )$ and $r , t_{-} \in C^{0} ( [ - \varepsilon_{0} , \varepsilon_{0} ] \times \CH_{\rm tang}^{+ \infty} )$.
\end{proposition}

From Proposition \ref{a65} and Proposition \ref{a74}, we deduce

\begin{corollary}\sl
For all $\alpha \in \CH_{\rm tang}^{+ \infty}$, the limit
\begin{equation*}
T ( \alpha ) = \lim_{\varepsilon \to 0} \big( t_{-}^{\varepsilon} ( \alpha ) - t_{+}^{\varepsilon} ( \alpha ) - 2 \vert \ln \varepsilon \vert / \lambda \big) ,
\end{equation*}
exists. Moreover, $T \in C^{0} ( \CH_{\rm tang}^{+ \infty} )$.
\end{corollary}

We now consider the Maslov's determinant appearing in \eqref{a81}. From \eqref{a75}, $\alpha \longmapsto x_{+} ( s , \alpha )$ is an immersion for $s < 0$ sufficiently large. Moreover, the Hamiltonian vector field $H_{p}$ is transverse to the manifold $\{ ( x_{+} ( s , \alpha ) , \nabla \varphi_{+} ( x_{+} ( s , \alpha ) ) ) ; \ \alpha \in \S^{n - 1} \}$. Then, from the standard theory (see e.g. the book of Maslov and Fedoriuk \cite{MaFe81_01}), the Jacobian
\begin{equation*}
J ( t , \alpha ) = \det \frac{\partial x_{+} ( t , \alpha )}{\partial ( t , \alpha )} ,
\end{equation*}
is a $C^{\infty}$ function which does not vanish near $( s , \alpha )$, $\alpha \in \S^{n - 1}$. Moreover, assume that a part of the evolution of $\Lambda_{+}$ by the Hamiltonian flow projects nicely on the $x$-space near $x_{+} ( t , \alpha )$,
\begin{equation*}
\Lambda_{+} = \{ ( x , \nabla \psi ( x ) ) ; \ x \text{ near } x_{+} ( t , \alpha ) \} ,
\end{equation*}
for some generating function $\psi \in C^{\infty}( \R^{n} )$ defined in a neighborhood of $x_{+} ( t , \alpha )$. Then, in this case, $J ( t , \alpha )$ does not vanish and
\begin{equation} \label{a76}
\partial_{t} \ln \vert J ( t , \alpha ) \vert = 2 \Delta \psi ( x_{+} ( t , \alpha ) ) .
\end{equation}

Let now $\alpha \in \CH_{\rm tang}^{+ \infty}$. Near $\rho_{+} ( t , \alpha )$, with $t \in [ t_{+}^{\varepsilon} ( \alpha ) , t_{+}^{\varepsilon_{0}} ( \alpha ) ]$, $\varphi_{+}$ is a generating function of $\Lambda_{+} = \Lambda_{+}^{0}$. Then, \eqref{a76} yields
\begin{equation} \label{a77}
J ( t_{+}^{\varepsilon} ( \alpha ) , \alpha ) = J ( t_{+}^{\varepsilon_{0}} ( \alpha ) , \alpha ) e^{- 2 \int_{t_{+}^{\varepsilon} ( \alpha )}^{t_{+}^{\varepsilon_{0}} ( \alpha )} \Delta \varphi_{+} ( x_{+} ( s , \alpha ) ) \, d s} ,
\end{equation}
and $J ( t_{+}^{\varepsilon} ( \alpha ) , \alpha )$ does not vanish. On the other hand, since $T_{\rho} \Lambda_{-} = T_{\rho} \Lambda_{+}$ on each point $\rho \in \CH_{\rm tang}$ and since $\Lambda_{-} = \Lambda_{-}^{0}$ projects nicely on the $x$-space near $0$, the evolution of $\Lambda_{+}$ projects nicely on the $x$-space near each point $\rho_{+} ( t , \alpha )$ for $t \in [ t_{-}^{\varepsilon_{0}} ( \alpha ) , t_{-}^{\varepsilon} ( \alpha ) ]$. Moreover, if $\psi$ is a generating function of the evolution of $\Lambda_{+}$, we have
\begin{equation*}
\Delta \psi ( x_{+} ( t , \alpha ) ) = \Delta \varphi_{-} ( x_{+} ( t , \alpha ) ) .
\end{equation*}
Therefore,  \eqref{a76} implies
\begin{equation} \label{a78}
J ( t_{-}^{\varepsilon} ( \alpha ) , \alpha ) = J ( t_{-}^{\varepsilon_{0}} ( \alpha ) , \alpha ) e^{2 \int_{t_{-}^{\varepsilon_{0}} ( \alpha )}^{t_{-}^{\varepsilon} ( \alpha )} \Delta \varphi_{-} ( x_{+} ( s , \alpha ) ) \, d s} ,
\end{equation}
and $J ( t_{-}^{\varepsilon} ( \alpha ) , \alpha )$ does not vanish. Then, using the notation of \eqref{a81} and combining \eqref{a77} and \eqref{a78}, we get
\begin{align}
\CM_{\varepsilon} ( \alpha ) &= \frac{\sqrt{J} ( t_{+}^{\varepsilon} ( \alpha ) , \alpha )}{\sqrt{J} ( t_{-}^{\varepsilon} ( \alpha ) , \alpha )}    \nonumber \\
&= \frac{\sqrt{J} ( t_{+}^{\varepsilon_{0}} ( \alpha ) , \alpha )}{\sqrt{J} ( t_{-}^{\varepsilon_{0}} ( \alpha ) , \alpha )} \frac{\sqrt{J} ( t_{+}^{\varepsilon} ( \alpha ) , \alpha )}{\sqrt{J} ( t_{+}^{\varepsilon_{0}} ( \alpha ) , \alpha )} \frac{\sqrt{J} ( t_{-}^{\varepsilon_{0}} ( \alpha ) , \alpha )}{\sqrt{J} ( t_{-}^{\varepsilon} ( \alpha ) , \alpha )}    \nonumber \\
&= \CM_{\varepsilon_{0}} ( \alpha ) e^{- \int_{t_{+}^{\varepsilon} ( \alpha )}^{t_{+}^{\varepsilon_{0}} ( \alpha )} \Delta \varphi_{+} ( x_{+} ( s , \alpha ) ) \, d s} e^{- \int_{t_{-}^{\varepsilon_{0}} ( \alpha )}^{t_{-}^{\varepsilon} ( \alpha )} \Delta \varphi_{-} ( x_{+} ( s , \alpha ) ) \, d s}   \nonumber  \\
&= \CM_{\varepsilon_{0}} ( \alpha ) e^{- G ( \varepsilon , \alpha )} ,  \label{a79}
\end{align}
where
\begin{equation*}
G ( \varepsilon , \alpha ) = \int_{t_{+}^{\varepsilon} ( \alpha )}^{t_{+}^{\varepsilon_{0}} ( \alpha )} \Delta \varphi_{+} ( x_{+} ( s , \alpha ) ) \, d s + \int_{t_{-}^{\varepsilon_{0}} ( \alpha )}^{t_{-}^{\varepsilon} ( \alpha )} \Delta \varphi_{-} ( x_{+} ( s , \alpha ) ) \, d s ,
\end{equation*}
which is continuous on $] 0 , \varepsilon_{0} ] \times \CH_{\rm tang}^{+ \infty}$ thanks to Proposition \ref{a65} and Proposition \ref{a74}.

From \eqref{a75} and \eqref{a80}, there exists $C > 0$ such that
\begin{equation*}
\vert x_{+} ( t , \alpha ) \vert \leq C e^{- \lambda \vert t \vert} ,
\end{equation*}
for all $t \in \R$ and $\alpha \in \CH_{\rm tang}^{+ \infty}$. Then, the asymptotic expansions of $t_{\pm}^{\bullet}$ given in Proposition \ref{a65} and Proposition \ref{a74} and that of $\varphi_{\pm}$ given in \eqref{a57} imply
\begin{align}
G ( \varepsilon , \alpha ) - G ( \widetilde{\varepsilon} , \alpha ) ={}& \int_{t_{+}^{\varepsilon} ( \alpha )}^{t_{+}^{\widetilde{\varepsilon}} ( \alpha )} \Delta \varphi_{+} ( x_{+} ( s , \alpha ) ) \, d s + \int_{t_{-}^{\widetilde{\varepsilon}} ( \alpha )}^{t_{-}^{\varepsilon} ( \alpha )} \Delta \varphi_{-} ( x_{+} ( s , \alpha ) ) \, d s  \nonumber \\
={}& \frac{n \lambda}{2} \big( t_{+}^{\widetilde{\varepsilon}} ( \alpha ) - t_{+}^{\varepsilon} ( \alpha ) \big) + \int_{t_{+}^{\varepsilon} ( \alpha )}^{t_{+}^{\widetilde{\varepsilon}} ( \alpha )} \CO ( \vert x_{+} ( s , \alpha ) \vert ) \, d s  \nonumber   \\
&- \frac{n \lambda}{2} \big( t_{-}^{\varepsilon} ( \alpha ) - t_{-}^{\widetilde{\varepsilon}} ( \alpha ) \big) + \int_{t_{-}^{\widetilde{\varepsilon}} ( \alpha )}^{t_{-}^{\varepsilon} ( \alpha )} \CO ( \vert x_{+} ( s , \alpha ) \vert ) \, d s   \nonumber  \\
={}& \CO ( \varepsilon ) + \CO ( \widetilde{\varepsilon} ) + \int_{t_{+}^{\varepsilon} ( \alpha )}^{t_{+}^{\widetilde{\varepsilon}} ( \alpha )} \CO ( e^{- \lambda \vert t \vert} ) \, d s + \int_{t_{-}^{\widetilde{\varepsilon}} ( \alpha )}^{t_{-}^{\varepsilon} ( \alpha )} \CO ( \vert e^{- \lambda \vert t \vert} ) \, d s \nonumber  \\
={}& \CO ( \varepsilon ) + \CO ( \widetilde{\varepsilon} ) ,
\end{align}
uniformly in $0 < \varepsilon , \widetilde{\varepsilon} \leq \varepsilon_{0}$ and $\alpha \in \CH_{\rm tang}^{+ \infty}$. By the Cauchy criterion, $G ( \varepsilon , \alpha )$ converges as $\varepsilon \to 0$, uniformly in $\alpha \in \CH_{\rm tang}^{+ \infty}$, to a function denoted by $G ( 0 , \alpha )$. Moreover, since $G ( \varepsilon , \alpha )$ is continuous for $\varepsilon \neq 0$, the application $G ( \varepsilon , \alpha )$ is continuous, and then bounded, on the compact set $[ 0 , \varepsilon_{0} ] \times \CH_{\rm tang}^{+ \infty}$. Combining with \eqref{a79}, we eventually obtain

\begin{proposition}\sl \label{b53}
The map $\CM_{\varepsilon} ( \alpha )$ is well-defined on $] 0 , \varepsilon_{0} ] \times \CH_{\rm tang}^{+ \infty}$ and converges as $\varepsilon \to 0$, uniformly in $\alpha \in \CH_{\rm tang}^{+ \infty}$, to a function $\CM_{0} ( \alpha )$. Moreover, the map
\begin{equation*}
( \varepsilon , \alpha ) \longmapsto \CM_{\varepsilon} ( \alpha ) ,
\end{equation*}
is continuous and does not vanish on $[ 0 , \varepsilon_{0} ] \times \CH_{\rm tang}^{+ \infty}$.
\end{proposition}

Recall that, for $\omega \in \CH_{\rm tang}^{- \infty}$, $\alpha ( \omega ) \in \CH_{\rm tang}^{+ \infty}$ has been defined in Section \ref{s2} as the normalized asymptotic direction of $x_{-}  ( t , \omega )$ as $t \to - \infty$ (see Figure \ref{f9}). From \eqref{a64} and the definition of $\CH_{\rm tang}^{\pm \infty}$, this application is well-defined. Furthermore, it satisfies

\begin{lemma}\sl \label{b67}
The map $\alpha : \CH_{\rm tang}^{- \infty} \longrightarrow \CH_{\rm tang}^{+ \infty}$ is a homeomorphism.
\end{lemma}

\begin{proof}
By definition, $\alpha$ is a bijection from $\CH_{\rm tang}^{- \infty}$ onto $\CH_{\rm tang}^{+ \infty}$. Moreover, using the discussion after \eqref{a70} and also \eqref{a80}, we can write
\begin{equation*}
\alpha ( \omega ) = \frac{\Phi ( x_{+} ( T , \omega ) )}{\vert \Phi ( x_{+} ( T , \omega ) ) \vert} .
\end{equation*}
Thus, $\alpha$ is continuous since $\Phi$ and $x_{+}$ are smooth. Eventually, since we consider a Schr\"{o}dinger operator, $\alpha$ is an involution and the lemma follows.
\end{proof}

It remains to explain Remark \ref{a82}. We consider the case where $\CH_{\rm tang}^{+ \infty}$ contains a neighborhood of $\alpha_{0}$. Then, from \eqref{a80},
\begin{equation*}
x_{+} ( t , \alpha ) = x_{+} \big( - t + T ( \alpha ) , \omega ( \alpha ) \big) ,
\end{equation*}
where
\begin{equation*}
T ( \alpha ) = T + \frac{\ln \vert \Phi ( x_{+} ( T , \alpha ) ) \vert}{\lambda} \qquad \text{and} \qquad \omega ( \alpha ) = \frac{\Phi ( x_{+} ( T , \alpha ) )}{\vert \Phi ( x_{+} ( T , \alpha ) ) \vert} ,
\end{equation*}
are $C^{\infty}$ near $\alpha_{0}$. Using now \eqref{a75} and the asymptotic of $t_{+}^{\varepsilon} ( \alpha )$ given in Proposition \ref{a65}, we obtain, in a neighborhood of $\alpha_{0}$,
\begin{equation*}
\left\{ \begin{aligned}
&\frac{\partial x_{+} ( t , \alpha )}{\partial t} \big\vert_{t = t_{+}^{\varepsilon} ( \alpha )} = \lambda \varepsilon \alpha + \CO ( \varepsilon^{3/2} ) ,   \\
&\frac{\partial x_{+} ( t , \alpha )}{\partial \alpha} \big\vert_{t = t_{+}^{\varepsilon} ( \alpha )} = \varepsilon Id_{T_{\alpha} \S^{n - 1}} + \CO ( \varepsilon^{3/2} ) ,
\end{aligned} \right.
\end{equation*}
and
\begin{equation*}
\left\{ \begin{aligned}
&\frac{\partial x_{+} ( t , \alpha )}{\partial t} \big\vert_{t = t_{-}^{\varepsilon} ( \alpha )} = - \lambda \varepsilon \omega ( \alpha ) + \CO ( \varepsilon^{3/2} ) ,   \\
&\frac{\partial x_{+} ( t , \alpha )}{\partial \alpha} \big\vert_{t = t_{-}^{\varepsilon} ( \alpha )} = \varepsilon \frac{\partial \omega ( \alpha )}{\partial \alpha} + \CO ( \varepsilon^{3/2} ) .
\end{aligned} \right.
\end{equation*}
So, taking the determinant and letting $\varepsilon$ go to $0$, we see that \eqref{a83} holds in that case.

\Subsection{Construction of potential bumps with large scattering angles} \label{s21}

In this part, we construct potential bumps with large scattering angle in dimension $n = 2$. They are used to build some examples with specific geometries. The idea is to truncate properly $E_{0} - \lambda^{2} x^{2} / 4$.

\begin{figure}
\begin{center}
\begin{picture}(0,0)%
\includegraphics{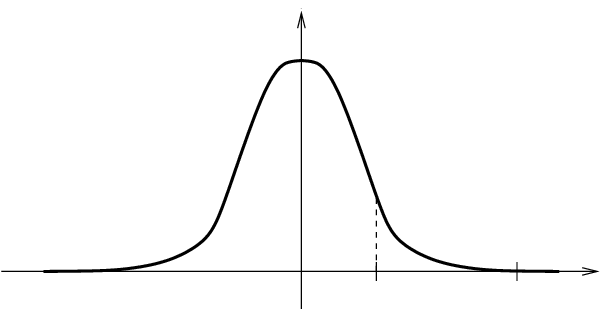}%
\end{picture}%
\setlength{\unitlength}{1184sp}%
\begingroup\makeatletter\ifx\SetFigFont\undefined%
\gdef\SetFigFont#1#2#3#4#5{%
  \reset@font\fontsize{#1}{#2pt}%
  \fontfamily{#3}\fontseries{#4}\fontshape{#5}%
  \selectfont}%
\fi\endgroup%
\begin{picture}(9644,5111)(1179,-6050)
\put(10501,-4936){\makebox(0,0)[lb]{\smash{{\SetFigFont{9}{10.8}{\rmdefault}{\mddefault}{\updefault}$x$}}}}
\put(6151,-1636){\makebox(0,0)[lb]{\smash{{\SetFigFont{9}{10.8}{\rmdefault}{\mddefault}{\updefault}$E_{0}$}}}}
\put(6901,-5986){\makebox(0,0)[b]{\smash{{\SetFigFont{9}{10.8}{\rmdefault}{\mddefault}{\updefault}$\frac{2}{\lambda} \sqrt{E_{0} - \frac{1}{N}}$}}}}
\put(9526,-5986){\makebox(0,0)[b]{\smash{{\SetFigFont{9}{10.8}{\rmdefault}{\mddefault}{\updefault}$\frac{2}{\lambda} \sqrt{E_{0}}$}}}}
\put(7201,-3286){\makebox(0,0)[lb]{\smash{{\SetFigFont{9}{10.8}{\rmdefault}{\mddefault}{\updefault}$V_{N} (x)$}}}}
\end{picture} $\qquad \qquad$ \begin{picture}(0,0)%
\includegraphics{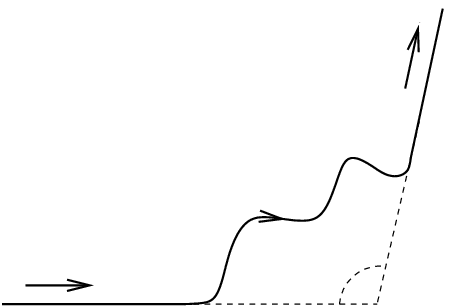}%
\end{picture}%
\setlength{\unitlength}{1184sp}%
\begingroup\makeatletter\ifx\SetFigFont\undefined%
\gdef\SetFigFont#1#2#3#4#5{%
  \reset@font\fontsize{#1}{#2pt}%
  \fontfamily{#3}\fontseries{#4}\fontshape{#5}%
  \selectfont}%
\fi\endgroup%
\begin{picture}(7116,4791)(3868,-3919)
\put(8626,-3211){\makebox(0,0)[lb]{\smash{{\SetFigFont{9}{10.8}{\rmdefault}{\mddefault}{\updefault}$\theta ( \rho )$}}}}
\put(8701, 14){\makebox(0,0)[lb]{\smash{{\SetFigFont{9}{10.8}{\rmdefault}{\mddefault}{\updefault}$\xi ( + \infty )$}}}}
\put(4051,-3286){\makebox(0,0)[lb]{\smash{{\SetFigFont{9}{10.8}{\rmdefault}{\mddefault}{\updefault}$\xi ( - \infty )$}}}}
\put(7126,-2086){\makebox(0,0)[lb]{\smash{{\SetFigFont{9}{10.8}{\rmdefault}{\mddefault}{\updefault}$\pi_{x} ( \rho (t) )$}}}}
\end{picture}%
\end{center}
\caption{The potential $V_{N}$ and the scattering angle $\theta ( \rho )$.} \label{f18}
\end{figure}

Let $\chi \in C^{\infty} ( \R ; [ 0 , 1 ] )$ be a non-increasing function such that $\chi = 1$ near $] - \infty , - 1 ]$ and $\chi = 0$ near $[ 0 , + \infty [$. For $N \geq 1$, the potential $V_{N} (x)$ is defined by
\begin{equation} \label{b89}
V_{N} (x) = \bigg( E_{0} - \frac{\lambda^{2}}{4} x^{2} \bigg) \chi \bigg( N \bigg( \frac{\lambda^{2}}{4} x^{2} - E_{0} \bigg) \bigg) .
\end{equation}
For $N \gg 1$, $V_{N} \in C^{\infty}_{0} ( \R^{2} ; [ 0 , E_{0} ] )$ is a radial function satisfying $x \cdot \nabla V_{N} (x) < 0$ for $x$ in the interior of $\supp V_{N} \setminus \{ 0 \}$ and 
\begin{equation*}
V_{N} (x) = E_{0} - \frac{\lambda^{2}}{4} x^{2} ,
\end{equation*}
near $0$ (see Figure \ref{f18}). We also define $p_{N} = \xi^{2} + V_{N} (x)$, the symbol of the Schr\"{o}dinger operator with potential $V_{N}$. From \eqref{d27}, any Hamiltonian trajectory $\rho (t) = ( x (t) , \xi (t) ) \subset p_{N}^{- 1} ( E_{0} )$ satisfies
\begin{equation} \label{b93}
\partial_{t} x^{2} (t) = 4 x (t) \cdot \xi (t) ,
\end{equation}
and
\begin{equation} \label{b92}
\partial_{t}^{2} x^{2} (t) = 8 \xi^{2} (t) - 4 x \cdot \nabla V_{N} (x) \geq 8 \xi^{2} (t) ,
\end{equation}
since $x \cdot \nabla V_{N} (x) \leq 0$. In particular, $x^{2} (t)$ is strictly convex when $\rho \neq ( 0 , 0 )$. Thus, any bicharacteristic curve which is bounded for positive or negative times goes to $( 0 , 0 )$. So,
\begin{equation*}
K ( E_{0} ) = \{ ( 0 , 0 ) \} ,
\end{equation*}
and any Hamiltonian trajectory in $p^{-1}_{N} ( E_{0} ) \setminus ( \Lambda_{-} \cup \Lambda_{+} )$ diverges as $t \to \pm \infty$. Moreover, $V_{N}$ satisfies the following estimate.

\begin{lemma}\sl \label{b90}
There exists $C > 0$ such that
\begin{equation*}
\vert \nabla V_{N} (x) \vert \leq C ,
\end{equation*}
for all $x \in \R^{2}$ and $N \geq 1$.
\end{lemma}

\begin{proof}
Indeed, a direct computation gives
\begin{align*}
\vert \nabla V_{N} (x) \vert &= \bigg\vert - \frac{\lambda^{2}}{2} x \chi \bigg( N \bigg( \frac{\lambda^{2}}{4} x^{2} - E_{0} \bigg) \bigg) + \frac{\lambda^{2}}{2} x N \bigg( E_{0} - \frac{\lambda^{2}}{4} x^{2} \bigg) \chi^{\prime} \bigg( N \bigg( \frac{\lambda^{2}}{4} x^{2} - E_{0} \bigg) \bigg) \bigg\vert  \\
&\leq \lambda \sqrt{E_{0}} \big( 1 + \Vert \chi^{\prime} \Vert_{L^{\infty}} \big) ,
\end{align*}
since $\vert x \vert \leq 2 \sqrt{E_{0}} / \lambda$ on the support of $V_{N}$.
\end{proof}

Let $\rho (t) = ( x (t) , \xi (t) )$ be a Hamiltonian trajectory in $p^{- 1}_{N} ( E_{0} ) \setminus ( \Lambda_{-} \cup \Lambda_{+} )$. Since $V_{N}$ is compactly supported and $x (t)$ goes to $\infty$ as $t \to \pm \infty$, $\xi (t)$ is constant equal to $\xi ( \pm \infty ) \in \sqrt{E_{0}} \S^{1}$ for $\pm t$ large enough. The scattering angle of the curve $\rho$ is then defined as
\begin{equation} \label{n14}
\theta ( \rho ) = \ang ( \xi ( - \infty ) , \xi ( + \infty ) ) \in [ 0 , \pi ] ,
\end{equation}
as illustrated in Figure \ref{f18}. The following result states that this scattering angle is essentially always bigger than $\pi/2$ for $N$ large enough.

\begin{proposition}\sl \label{a14}
For all $\varepsilon > 0$, there exists $N_{\varepsilon} \geq 1$ such that, for all $N \geq N_{\varepsilon}$, we have
\begin{equation*}
\theta ( \rho ) \geq \frac{\pi}{2} - \varepsilon ,
\end{equation*}
for all Hamiltonian curve $\rho \subset p^{-1}_{N} ( E_{0} ) \setminus ( \Lambda_{-} \cup \Lambda_{+} )$.
\end{proposition}

\begin{proof}
We now decompose $\R^{2}$ in three regions:
\begin{align*}
\CD_{1} &= \bigg\{ \vert x \vert \leq \frac{2}{\lambda} \sqrt{E_{0} - \frac{1}{N}} \bigg\} ,   \\
\CD_{2} &= \bigg\{ \frac{2}{\lambda} \sqrt{E_{0} - \frac{1}{N}} \leq \vert x \vert \leq \frac{2}{\lambda} \sqrt{E_{0}} \bigg\} ,   \\
\CD_{3} &= \Big\{ \frac{2}{\lambda} \sqrt{E_{0}} \leq \vert x \vert \Big\} ,
\end{align*}
as illustrated in Figure \ref{f19}. Of course, from \eqref{b89}, we have
\begin{equation*}
V_{N} ( x ) = \left\{ \begin{aligned}
&E_{0} - \frac{\lambda^{2}}{4} x^{2} &&\text{for } x \in \CD_{1} ,  \\
&0 &&\text{for } x \in \CD_{3} .
\end{aligned} \right.
\end{equation*}

\begin{figure}
\begin{center}
\begin{picture}(0,0)%
\includegraphics{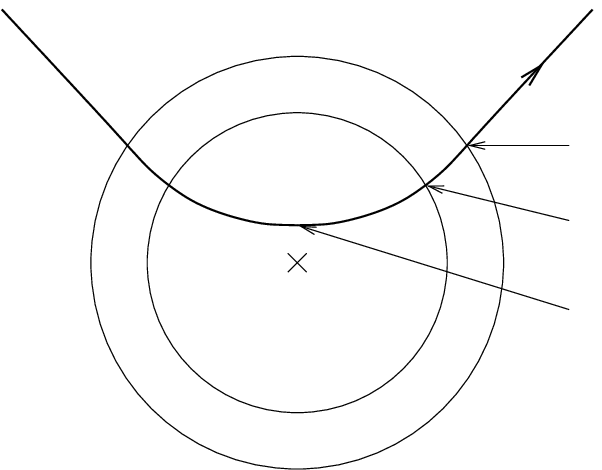}%
\end{picture}%
\setlength{\unitlength}{1184sp}%
\begingroup\makeatletter\ifx\SetFigFont\undefined%
\gdef\SetFigFont#1#2#3#4#5{%
  \reset@font\fontsize{#1}{#2pt}%
  \fontfamily{#3}\fontseries{#4}\fontshape{#5}%
  \selectfont}%
\fi\endgroup%
\begin{picture}(9516,7399)(1243,-8477)
\put(8851,-1861){\makebox(0,0)[lb]{\smash{{\SetFigFont{9}{10.8}{\rmdefault}{\mddefault}{\updefault}$x ( t )$}}}}
\put(7501,-7486){\makebox(0,0)[lb]{\smash{{\SetFigFont{9}{10.8}{\rmdefault}{\mddefault}{\updefault}$\CD_{2}$}}}}
\put(6451,-6961){\makebox(0,0)[lb]{\smash{{\SetFigFont{9}{10.8}{\rmdefault}{\mddefault}{\updefault}$\CD_{1}$}}}}
\put(8626,-7936){\makebox(0,0)[lb]{\smash{{\SetFigFont{9}{10.8}{\rmdefault}{\mddefault}{\updefault}$\CD_{3}$}}}}
\put(6001,-5686){\makebox(0,0)[b]{\smash{{\SetFigFont{9}{10.8}{\rmdefault}{\mddefault}{\updefault}$0$}}}}
\put(10501,-3361){\makebox(0,0)[lb]{\smash{{\SetFigFont{9}{10.8}{\rmdefault}{\mddefault}{\updefault}$x ( t_{2} )$}}}}
\put(10501,-4636){\makebox(0,0)[lb]{\smash{{\SetFigFont{9}{10.8}{\rmdefault}{\mddefault}{\updefault}$x ( t_{1} )$}}}}
\put(10501,-6061){\makebox(0,0)[lb]{\smash{{\SetFigFont{9}{10.8}{\rmdefault}{\mddefault}{\updefault}$x ( t_{0} )$}}}}
\end{picture}%
\end{center}
\caption{The geometric setting of the proof of Proposition \ref{a14}.} \label{f19}
\end{figure}

Since $x^{2} (t)$ is strictly convex (see \eqref{b92} and below), one can show that there exists a unique time, denoted by $t_{0}$ in the following, such that $\partial_{t} x^{2} ( t_{0} ) = 0$. Then, after a symplectic rotation, we can assume that
\begin{equation*}
x ( t_{0} ) = \left( \begin{gathered}
0 \\
\alpha
\end{gathered} \right)
\qquad \text{and} \qquad
\xi ( t_{0} ) = \left( \begin{gathered}
\beta \\
0
\end{gathered} \right) ,
\end{equation*}
for some $\alpha > 0$ and $\beta \in \R$. By symmetry, it is enough to consider the situation $\beta \geq 0$ (see Figure \ref{f19}).

We first assume that $x ( t_{0} ) \in \CD_{1}$. Then, $\beta = \alpha \lambda / 2$. As long as $x (t)$ stays in $\CD_{1}$, we have
\begin{equation*}
x (t) = \left( \begin{gathered}
\alpha \sinh \big( \lambda ( t - t_{0} ) \big)  \\
\alpha \cosh \big( \lambda ( t - t_{0} ) \big)
\end{gathered} \right)
\qquad \text{and} \qquad
\xi (t) = \left( \begin{gathered}
\frac{\lambda}{2} \alpha \cosh \big( \lambda ( t - t_{0} ) \big)  \\
\frac{\lambda}{2} \alpha \sinh \big( \lambda ( t - t_{0} ) \big)
\end{gathered} \right) .
\end{equation*}
In particular, there exists a first time $t_{1} \geq t_{0}$ such that $x ( t_{1} ) \in \partial \CD_{1}$. At this time, the previous equations give $x_{1} ( t_{1} ) , x_{2} ( t_{1} ) , \xi_{1} ( t_{1} ) , \xi_{2} ( t_{1} ) \geq 0$ and
\begin{equation} \label{b91}
\frac{\xi_{2} ( t_{1} )}{\xi_{1} ( t_{1} )} = \tanh \big( \lambda ( t_{1} - t_{0} ) \big) \leq 1 .
\end{equation}
Moreover, since $x ( t_{1} ) \in \partial \CD_{1}$ and $p_{N} ( \rho ( t_{1} ) ) = E_{0}$, we have
\begin{equation} \label{b95}
x^{2} ( t_{1} ) = \frac{4}{\lambda^{2}} \Big( E_{0} - \frac{1}{N} \Big) \qquad \text{and} \qquad \xi^{2} ( t_{1} ) = E_{0} - \frac{1}{N} .
\end{equation}
Combining with \eqref{b91}, it implies
\begin{equation} \label{b97}
\xi_{1}^{2} ( t_{1} ) \geq \frac{\xi^{2} (t)}{2} = \frac{E_{0}}{2} - \frac{1}{2 N} \geq \frac{E_{0}}{4} ,
\end{equation}
for $N$ large enough.

We will now show that there exists $t_{2} \geq t_{1}$ close to $t_{1}$ such that $x ( t ) \in \CD_{3}$ for all $t \geq t_{2}$. From Lemma \ref{b90}, we have
\begin{equation*}
\xi ( t_{1} + \delta ) = \xi ( t_{1} ) + \CO ( \delta ) .
\end{equation*}
Here and in the following, the $\CO$'s are uniform with respect to $N \geq 1$ and to the curve $\rho$. Then, \eqref{b92} yields
\begin{equation*}
\partial_{t} x^{2} ( t_{1} + \delta ) \geq \partial_{t} x^{2} ( t_{1} ) + 8 \xi^{2} ( t_{1} ) \delta + \CO ( \delta^{2} ) .
\end{equation*}
Using \eqref{b93} and the sign of the coordinates of $\rho ( t_{1} )$ given above \eqref{b91}, we get $\partial_{t} x^{2} ( t_{1} ) \geq 0$ and then
\begin{equation} \label{b94}
\partial_{t} x^{2} ( t_{1} + \delta ) \geq 8 \xi^{2} ( t_{1} ) \delta + \CO ( \delta^{2} ) .
\end{equation}
Integrating this estimate with respect to $\delta$ implies
\begin{equation*}
x^{2} ( t_{1} + \delta ) \geq x^{2} ( t_{1} ) + 4 \xi^{2} ( t_{1} ) \delta^{2} + \CO ( \delta^{3} ) .
\end{equation*}
Using \eqref{b95}, we obtain
\begin{equation*}
x^{2} ( t_{1} + \delta ) \geq \frac{4}{\lambda^{2}} E_{0} + 4 E_{0} \delta^{2} + \CO \big( N^{-1} + \delta^{2} N^{-1} + \delta^{3} \big) .
\end{equation*}
If we take $\delta = N^{- 1 / 3}$, we have proved
\begin{equation*}
x^{2} ( t_{1} + N^{- \frac{1}{3}} ) \geq \frac{4}{\lambda^{2}} E_{0} + 4 E_{0} N^{- \frac{2}{3}} + \CO ( N^{-1} ) > \frac{4}{\lambda^{2}} E_{0} ,
\end{equation*}
for $N$ large enough. Thus, for $N$ large enough, the theorem of the intermediate values provides a (first) time $t_{2} \in ] t_{1} , t_{1} + N^{- 1 / 3} [$ such that $x ( t_{2} ) \in \partial \CD_{3}$. Moreover, from \eqref{b94} and \eqref{b95}, we have
\begin{equation*}
\partial_{t} x^{2} ( t_{2}) \geq 0 ,
\end{equation*}
for $N$ large enough. Since $x^{2} (t)$ is convex from \eqref{b92}, it implies $x^{2} (t) \geq x^{2} ( t_{2} ) = 4 E_{0} / \lambda^{2}$ for all $t \geq t_{2}$. In other words, $x (t) \in \CD_{3}$ for all $t \geq t_{2}$. Since the Hamiltonian trajectories are uniform rectilinear in $\CD_{3}$, we get $\xi ( + \infty ) = \xi ( t_{2} )$. Thus, combining with Lemma \ref{b90}, it follows
\begin{equation} \label{b96}
\xi ( + \infty ) = \xi ( t_{1} ) + \CO ( \vert t_{2} - t_{1} \vert ) = \xi ( t_{1} ) + \CO \big( N^{- \frac{1}{3}} \big) .
\end{equation}

Summing up, we have
\begin{align*}
\theta ( \rho ) &= 2 \ang \bigg( \left( \begin{array}{c}
0 \\
1 
\end{array} \right) , \xi ( + \infty ) \bigg)  \\
&= \pi - 2 \arctan \Big( \frac{\xi_{2} ( + \infty )}{\xi_{1} ( + \infty )} \Big)   \\
&= \pi - 2 \arctan \Big( \frac{\xi_{2} ( t_{1} )}{\xi_{1} ( t_{1} )} \Big) + \CO \big( N^{- \frac{1}{3}} \big) ,
\end{align*}
from \eqref{b97} and \eqref{b96}. Using now \eqref{b91}, we eventually obtain
\begin{equation}
\theta ( \rho ) \geq \frac{\pi}{2} + \CO \big( N^{- \frac{1}{3}} \big) .
\end{equation}
which finishes the proof if $\rho ( t_{0} ) \in \CD_{1}$.

If $\rho ( t_{0} ) \in \CD_{2}$, we can proceed as before. More precisely, taking $t_{1} = t_{0}$, the proof is the same except that \eqref{b95} is replaced by
\begin{equation*}
x^{2} ( t_{1} ) = \frac{4}{\lambda^{2}} E_{0} + \CO \big( N^{- 1} \big) \qquad \text{and} \qquad \xi^{2} ( t_{1} ) = E_{0} + \CO \big( N^{- 1} \big) .
\end{equation*}
Finally, if $\rho ( t_{0} ) \in \CD_{3}$, we directly have $\theta ( \rho ) = \pi$.
\end{proof}

We end this part by constructing a potential with a small scattering angle. It is used to create a counterexample in Example \ref{k74}.

\begin{example}\rm \label{k76}
Let $F \in C^{\infty}_{0} ( \R^{2} ; [ 0 , + \infty [ )$ be a radial function such that $F ( 0 ) > E_{0} > 0$ and $x \cdot \nabla F ( x ) < 0$ for $x$ in the interior of $\supp F \setminus \{ 0 \}$. Such potential is illustrated in Figure \ref{f58}. We define
\begin{equation*}
q ( x , \xi ) = \xi^{2} + F ( x ) ,
\end{equation*}
the symbol of the Schr\"{o}dinger operator with symbol $F ( x )$. The previous hypotheses imply that $q$ is non-trapping at energy $E_{0}$ ($x \cdot \xi$ is an escape function).

\begin{figure}
\begin{center}
\begin{picture}(0,0)%
\includegraphics{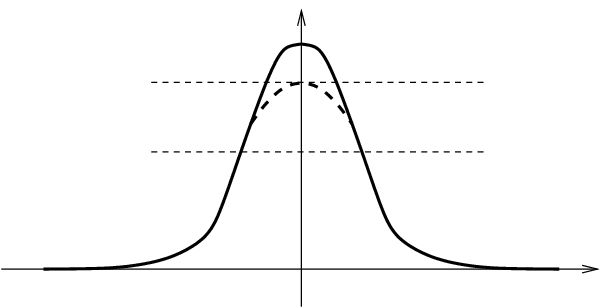}%
\end{picture}%
\setlength{\unitlength}{1184sp}%
\begingroup\makeatletter\ifx\SetFigFont\undefined%
\gdef\SetFigFont#1#2#3#4#5{%
  \reset@font\fontsize{#1}{#2pt}%
  \fontfamily{#3}\fontseries{#4}\fontshape{#5}%
  \selectfont}%
\fi\endgroup%
\begin{picture}(9644,4844)(1179,-5783)
\put(5626,-2761){\makebox(0,0)[lb]{\smash{{\SetFigFont{9}{10.8}{\rmdefault}{\mddefault}{\updefault}$W$}}}}
\put(10501,-4936){\makebox(0,0)[lb]{\smash{{\SetFigFont{9}{10.8}{\rmdefault}{\mddefault}{\updefault}$x$}}}}
\put(3451,-2311){\makebox(0,0)[rb]{\smash{{\SetFigFont{9}{10.8}{\rmdefault}{\mddefault}{\updefault}$E_{0}$}}}}
\put(3451,-3436){\makebox(0,0)[rb]{\smash{{\SetFigFont{9}{10.8}{\rmdefault}{\mddefault}{\updefault}$E_{0} - \varepsilon$}}}}
\put(6376,-1636){\makebox(0,0)[lb]{\smash{{\SetFigFont{9}{10.8}{\rmdefault}{\mddefault}{\updefault}$F$}}}}
\end{picture} $\qquad \qquad \qquad$ \begin{picture}(0,0)%
\includegraphics{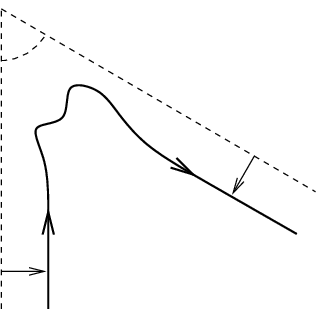}%
\end{picture}%
\setlength{\unitlength}{1184sp}%
\begingroup\makeatletter\ifx\SetFigFont\undefined%
\gdef\SetFigFont#1#2#3#4#5{%
  \reset@font\fontsize{#1}{#2pt}%
  \fontfamily{#3}\fontseries{#4}\fontshape{#5}%
  \selectfont}%
\fi\endgroup%
\begin{picture}(5069,4855)(2379,-3994)
\put(3751,-1861){\makebox(0,0)[b]{\smash{{\SetFigFont{9}{10.8}{\rmdefault}{\mddefault}{\updefault}$x ( t )$}}}}
\put(3001,-361){\makebox(0,0)[b]{\smash{{\SetFigFont{9}{10.8}{\rmdefault}{\mddefault}{\updefault}$\pi / 3$}}}}
\put(6676,-2086){\makebox(0,0)[b]{\smash{{\SetFigFont{9}{10.8}{\rmdefault}{\mddefault}{\updefault}$z_{+}$}}}}
\put(2776,-3811){\makebox(0,0)[b]{\smash{{\SetFigFont{9}{10.8}{\rmdefault}{\mddefault}{\updefault}$z_{-}$}}}}
\end{picture}%
\end{center}
\caption{The potentials $F , W$ and the impact parameters $z_{\pm}$.} \label{f58}
\end{figure}

Fix $R > 0$ large enough. For $z \in \R$, let $\theta ( z )$ be the scattering angle (see \eqref{n14}) of the Hamiltonian trajectory
\begin{equation*}
t \longmapsto \exp ( t H_{p} ) \big( z , - R , 0 , \sqrt{E_{0}} \big) .
\end{equation*}
The previous curve is the unique Hamiltonian trajectory of energy $E_{0}$ with incoming direction $( 0 , 1 )$ and impact parameter $z \in \R$. For $z > 0$ large enough, we have $\theta ( z ) = \pi$ since $F$ is compactly supported. On the other hand, $\theta ( 0 ) = 0$. Since $\theta$ is smooth and $\theta ( [ 0 , + \infty [ ) \subset [ 0 , \pi ]$, the mean value theorem implies that there exists $z_{-} \in ] 0 , + \infty [$ such that $\theta ( z_{-} ) = \pi / 3$. In the sequel, we note
\begin{equation*}
\rho ( t ) = ( x ( t ) , \xi ( t ) ) : = \exp ( t H_{p} ) \big( z_{-} , - R , 0 , \sqrt{E_{0}} \big) .
\end{equation*}

Assume that $\xi ( t_{0} ) = 0$ for some $t_{0} \in \R$. Then, the characteristic trajectory $s \mapsto x ( t_{0} + s )$ belongs to the half line $x ( t_{0} ) [ 0 , + \infty [$ with $ x ( t_{0} ) \neq 0$ since $F$ is radial. This is impossible because the scattering angle of $\rho$ is $\pi / 3$ and not $0$. Thus, we have $\xi ( t ) \neq 0$ for all $t \in \R$. Using $\vert \xi ( t ) \vert = \sqrt{E_{0}}$ for $\vert t \vert$ large enough, we deduce that
\begin{equation*}
\varepsilon : = \inf_{t \in \R} \vert \xi ( t ) \vert > 0 .
\end{equation*}
From $\xi^{2} ( t ) + F ( x ( t ) ) = E_{0}$, it yields
\begin{equation} \label{k78}
F ( x ( t ) ) \leq E_{0} - \varepsilon ,
\end{equation}
for all $t \in \R$.

We now modify the function $F$ outside $F^{- 1} ( [ 0 , E_{0} - \varepsilon ])$. Let $W \in C^{\infty}_{0} ( \R^{2} , [ 0 , E_{0} ] )$ be a radial function such that $x \cdot \nabla W ( x ) < 0$ for $x$ in the interior of $\supp W \setminus \{ 0 \}$ and
\begin{equation} \label{k77}
W ( x ) = \left\{ \begin{aligned}
&F ( x ) &&\text{ near } F^{- 1} ( [ 0 , E_{0} - \varepsilon ]) ,  \\
&E_{0} - \frac{\lambda^{2}}{2} x^{2} &&\text{ near } 0 ,
\end{aligned} \right.
\end{equation}
(see Figure \ref{f58}). Here, $\lambda$ is any given positive constant. By construction and \eqref{k78}, we have $W = F$ near the characteristic curve $\{ x ( t ) ; \ t \in \R \}$. This shows that $\rho ( t )$ is also a Hamiltonian trajectory of $p ( x , \xi ) = \xi^{2} + W ( x )$. In other words, we have found a reasonable bump $W$ such that
\begin{equation} \label{k80}
p \text{ has a trajectory of scattering angle } \pi / 3 \text{ in the energy surface } E_{0} .
\end{equation}

Finally, we compute the outgoing impact parameter of $\rho ( t )$ using a classical trick. We set $f ( t ) = x ( t ) \cdot J \xi ( t )$ with the matrix
\begin{equation*}
J = \left( \begin{array}{cc}
0 & 1 \\
- 1 & 0
\end{array} \right) .
\end{equation*}
For $\pm t > 0$ large enough, $f ( t ) = \sqrt{E_{0}} z_{\pm}$ where $z_{+}$ (resp. $z_{-}$) is the outgoing (resp. incoming) impact parameter (see Figure \ref{f58}). The Hamiltonian equations give
\begin{equation*}
\partial_{t} f ( t ) = 2 \xi ( t ) \cdot J \xi ( t ) - x ( t ) \cdot J \nabla W ( x ( t ) ) = 0 ,
\end{equation*}
since $W ( x )$ is radial. The previous properties imply
\begin{equation} \label{k79}
z_{-} = z_{+} .
\end{equation}
\end{example}

\Subsection{Proof of Proposition \ref{e24}} \label{s86}

In this section, we prove Proposition \ref{e24}. By assumption, there exists $\varepsilon > 0$ such that $p ( x , \xi ) = E_{0} + \xi^{2} - \sum \lambda_{j}^{2} x_{j}^{2} / 4$ in $B : = B ( 0 , \varepsilon ) \subset T^{*} \R^{n}$. Let $( x ( t ) , \xi ( t ) )$ be a Hamiltonian trajectory in $B$. A classical computation gives
\begin{equation*}
\partial_{t}^{2} x^{2} = 4 \partial_{t} ( x \cdot \xi ) = 8 \xi^{2} + 2 L^{2} x \cdot x ,
\end{equation*}
where $L = \diag ( \lambda_{1} , \ldots , \lambda_{n} )$. Thus, the function $t \mapsto x^{2} ( t )$ is strictly convex in $B \cap p^{-1} ( E )$ for all $E \neq E_{0}$. In particular, every Hamiltonian trajectory in $p^{- 1} ( E )$ with $E \neq E_{0}$ which enters in $B$ must leave $B$.

We prove the proposition using a contradiction argument. So, assume that all the energies just above $E_{0}$ are not non-trapped. Then, there exist a decreasing sequence of $E$ converging to $E_{0}$ and $\rho_{E} \in K ( E )$. From the previous discussion, we can always assume that $\rho_{E} \notin B$. On the other hand, $K ( [ E_{0} , E_{0} + 1] )$ is bounded from \ref{h1}. Thus, up to the extraction of a subsequence, $\rho_{E}$ tends to $\rho_{E_{0}}$ as $E$ goes to $E_{0}$. In particular, $p ( \rho_{E_{0}} ) = E_{0}$ and $\rho_{E_{0}} \neq 0$ since $\rho_{E_{0}} \notin B$. Lastly, the continuity of the Hamiltonian flow implies that $\rho_{E_{0}} \in K ( E_{0} )$.

\begin{figure}
\begin{center}
\begin{picture}(0,0)%
\includegraphics{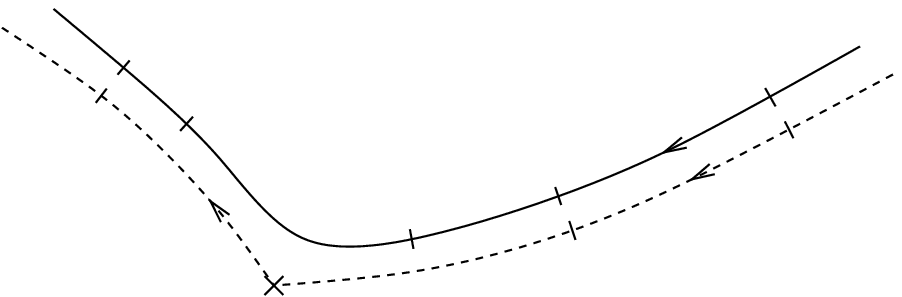}%
\end{picture}%
\setlength{\unitlength}{1184sp}%
\begingroup\makeatletter\ifx\SetFigFont\undefined%
\gdef\SetFigFont#1#2#3#4#5{%
  \reset@font\fontsize{#1}{#2pt}%
  \fontfamily{#3}\fontseries{#4}\fontshape{#5}%
  \selectfont}%
\fi\endgroup%
\begin{picture}(14316,4641)(1618,-6544)
\put(2326,-3961){\makebox(0,0)[lb]{\smash{{\SetFigFont{9}{10.8}{\rmdefault}{\mddefault}{\updefault}$\rho_{+}^{E_{0}}$}}}}
\put(4501,-6436){\makebox(0,0)[lb]{\smash{{\SetFigFont{9}{10.8}{\rmdefault}{\mddefault}{\updefault}$( 0 , 0 )$}}}}
\put(10201,-4636){\makebox(0,0)[lb]{\smash{{\SetFigFont{9}{10.8}{\rmdefault}{\mddefault}{\updefault}$\rho_{-}^{E}$}}}}
\put(13351,-2986){\makebox(0,0)[lb]{\smash{{\SetFigFont{9}{10.8}{\rmdefault}{\mddefault}{\updefault}$\rho_{E}$}}}}
\put(14326,-4411){\makebox(0,0)[lb]{\smash{{\SetFigFont{9}{10.8}{\rmdefault}{\mddefault}{\updefault}$\rho_{E_{0}}$}}}}
\put(10801,-6061){\makebox(0,0)[lb]{\smash{{\SetFigFont{9}{10.8}{\rmdefault}{\mddefault}{\updefault}$\rho_{-}^{E_{0}}$}}}}
\put(6751,-5011){\makebox(0,0)[lb]{\smash{{\SetFigFont{9}{10.8}{\rmdefault}{\mddefault}{\updefault}$\rho^{t , E} = \rho_{-}^{E} ( t )$}}}}
\put(4876,-3661){\makebox(0,0)[lb]{\smash{{\SetFigFont{9}{10.8}{\rmdefault}{\mddefault}{\updefault}$\rho_{+}^{E} ( - t )$}}}}
\put(3751,-2536){\makebox(0,0)[lb]{\smash{{\SetFigFont{9}{10.8}{\rmdefault}{\mddefault}{\updefault}$\rho_{+}^{E} = \rho_{-}^{E} ( T_{E} )$}}}}
\end{picture}%
\end{center}
\caption{The geometric setting of Section \ref{s86}.} \label{f46}
\end{figure}

For $\rho \in T^{*} \R^{n}$ and $t \in \R$, we note $\rho ( t ) = \exp ( t H_{p} ) ( \rho )$. Since $K ( E_{0} ) \setminus \{ ( 0 , 0 ) \} = \CH$ and since every homoclinic trajectory goes to $( 0 , 0 )$ as $t \to + \infty$, there exists $t_{-} \geq 0$ such that $\rho_{E_{0}} ( t_{-} ) \in B \cap \Lambda_{-}^{0}$. We define $\rho_{-}^{E} = \rho_{E} ( t_{-} )$ which lies in $B$ for $E$ close enough to $E_{0}$. Moreover, from the first paragraph of this proof, the curve $\rho_{-}^{E} ( t )$ leaves $B$ for $E \neq E_{0}$. Let $T_{E} > 0$ be the first time such that $\rho_{-}^{E} ( T_{E} ) \in \partial B$. Since $\rho_{-}^{E} = \rho_{-}^{E_{0}} + o_{E \to E_{0}} ( 1 )$ and $\rho_{-}^{E_{0}} \in \Lambda_{-}^{0}$, the continuity of the Hamiltonian flow implies that
\begin{equation} \label{j51}
T_{E} \longrightarrow + \infty ,
\end{equation}
as $E$ goes to $E_{0}$. We define $\rho_{+}^{E} = \rho_{-}^{E} ( T_{E} )$. Since $\partial B$ is a compact set, we can always assume up to the extraction of a subsequence that $\rho_{+}^{E}$ tends to $\rho_{+}^{E_{0}}$ as $E$ goes to $E_{0}$. As before, $\rho_{+}^{E_{0}} \in \CH \cap \partial B$. Furthermore, using one more time the continuity of the Hamiltonian flow, we deduce that $\rho_{+}^{E_{0}} ( t ) \in \overline{B}$ for all $t \leq 0$. Then, Proposition \ref{a44} gives $\rho_{+}^{E_{0}} \in \Lambda_{+}^{0}$. The situation is illustrated in Figure 
\ref{f46}.

Let $( x_{\bullet}^{\star} , \xi_{\bullet}^{\star} )$ denotes the coordinates of $\rho_{\bullet}^{\star}$. From \eqref{d2}, we have
\begin{align*}
x_{-}^{E_{0}} ( t ) &= g_{-} ( \rho_{-}^{E_{0}} ) e^{- \lambda_{1} t } + o_{t \to + \infty} ( e^{- \lambda_{1} t } ) ,   \\
x_{+}^{E_{0}} ( - t ) &= g_{+} ( \rho_{+}^{E_{0}} ) e^{- \lambda_{1} t } + o_{t \to + \infty} ( e^{- \lambda_{1} t } ) ,
\end{align*}
for $t > 0$. Then, the continuity of the Hamiltonian flow gives
\begin{align}
x_{-}^{E} ( t ) &= g_{-} ( \rho_{-}^{E_{0}} ) e^{- \lambda_{1} t } + o_{t \to + \infty} ( e^{- \lambda_{1} t } ) + o_{E \to E_{0}}^{t} ( 1 ) ,  \label{j54} \\
x_{+}^{E} ( - t ) &= g_{+} ( \rho_{+}^{E_{0}} ) e^{- \lambda_{1} t } + o_{t \to + \infty} ( e^{- \lambda_{1} t } ) + o_{E \to E_{0}}^{t} ( 1 ) .  \label{j53}
\end{align}
The notation $o_{E \to E_{0}}^{t} ( 1 )$ designs a function which tends to $0$ as $E$ goes to $E_{0}$ for $t \geq 0$ fixed. In particular, 
\begin{equation} \label{j58}
x_{-}^{E} ( t ) \cdot x_{+}^{E} ( - t ) = g_{-} ( \rho_{-}^{E_{0}} ) \cdot g_{+} ( \rho_{+}^{E_{0}} ) e^{- 2 \lambda_{1} t } + o_{t \to + \infty} ( e^{- 2 \lambda_{1} t } ) + o_{E \to E_{0}}^{t} ( 1 ) .
\end{equation}

We now define
\begin{equation*}
\rho^{t , E} = ( x^{t , E } , \xi^{t , E} ) : = \rho_{-}^{E} ( t ) .
\end{equation*}
Since $\Lambda_{-}^{0}$ is given by $\{ ( x , \xi ) ; \ \xi = - L x / 2 \}$ and $\rho^{t , E} = \rho_{-}^{E_{0}} ( t ) + o_{E \to E_{0}}^{t} ( 1 )$, we get
\begin{equation} \label{j55}
x_{-}^{E} ( t ) = L^{- 1} \Big( \frac{L}{2} x^{t , E} - \xi^{t , E} \Big) + o_{E \to E_{0}}^{t} ( 1 ) .
\end{equation}
On the other hand, since $V ( x ) = E_{0} - \sum \lambda_{j}^{2} x_{j}^{2} / 4$ in $B$, the Hamiltonian flow is explicit and we can write
\begin{align}
x_{+}^{E} ( - t ) &= x_{-}^{E} ( T_{E} - t ) =  x^{t , E} ( T_{E} - 2 t ) \nonumber \\
&= L^{- 1} e^{L ( T_{E} - 2 t )} \Big( \frac{L}{2} x^{t , E} + \xi^{t , E} \Big) + L^{- 1} e^{- L ( T_{E} - 2 t )} \Big( \frac{L}{2} x^{t , E} - \xi^{t , E} \Big) . \label{j52}
\end{align}
From \eqref{j51}, we deduce
\begin{equation*}
L^{- 1} e^{- L ( T_{E} - 2 t )} \Big( \frac{L}{2} x^{t , E} - \xi^{t , E} \Big) = o_{E \to E_{0}}^{t} ( 1 ) .
\end{equation*}
On the other hand, since the $j$-th coordinate of $g_{+} ( \rho_{+}^{E_{0}} )$ vanishes if $\lambda_{j} > \lambda_{1}$, \eqref{j53} and \eqref{j52} imply
\begin{equation*}
\lambda_{j}^{- 1} e^{\lambda_{j} ( T_{E} - 2 t )} \Big( \frac{\lambda_{j}}{2} x^{t , E} + \xi^{t , E} \Big)_{j} = o_{t \to + \infty} ( e^{- \lambda_{1} t } ) + o_{E \to E_{0}}^{t} ( 1 ) ,
\end{equation*}
if $\lambda_{j} > \lambda_{1}$. Since $T_{E} - 2 t \geq 0$ for $E$ close to $E_{0}$ (depending on $t$), the previous equation yields
\begin{equation*}
\lambda_{1}^{- 2} \lambda_{j} e^{\lambda_{1} ( T_{E} - 2 t )} \Big( \frac{\lambda_{j}}{2} x^{t , E} + \xi^{t , E} \Big)_{j} = o_{t \to + \infty} ( e^{- \lambda_{1} t } ) + o_{E \to E_{0}}^{t} ( 1 ) ,
\end{equation*}
if $\lambda_{j} > \lambda_{1}$. Combining the previous formulas, \eqref{j52} becomes
\begin{equation} \label{j56}
x_{+}^{E} ( - t ) = \lambda_{1}^{- 2} e^{\lambda_{1} ( T_{E} - 2 t )} L \Big( \frac{L}{2} x^{t , E} + \xi^{t , E} \Big) + o_{t \to + \infty} ( e^{- \lambda_{1} t } ) + o_{E \to E_{0}}^{t} ( 1 ) .
\end{equation}
Making the scalar product between \eqref{j55} and \eqref{j56}, we deduce
\begin{align}
x_{-}^{E} ( t ) \cdot x_{+}^{E} ( - t ) &= \lambda_{1}^{- 2} e^{\lambda_{1} ( T_{E} - 2 t )} \Big( \frac{L}{2} x^{t , E} - \xi^{t , E} \Big) \cdot \Big( \frac{L}{2} x^{t , E} + \xi^{t , E} \Big) + o_{t \to + \infty} ( e^{- 2 \lambda_{1} t } ) + o_{E \to E_{0}}^{t} ( 1 )  \nonumber \\
&= \lambda_{1}^{- 2} e^{\lambda_{1} ( T_{E} - 2 t )} \big( E_{0} - p ( \rho^{t , E} ) \big) + o_{t \to + \infty} ( e^{- 2 \lambda_{1} t } ) + o_{E \to E_{0}}^{t} ( 1 )  \nonumber \\
&= \lambda_{1}^{- 2} e^{\lambda_{1} ( T_{E} - 2 t )} ( E_{0} - E ) + o_{t \to + \infty} ( e^{- 2 \lambda_{1} t } ) + o_{E \to E_{0}}^{t} ( 1 )  \nonumber \\
&\leq o_{t \to + \infty} ( e^{- 2 \lambda_{1} t } ) + o_{E \to E_{0}}^{t} ( 1 ) , \label{j57}
\end{align}
since $E > E_{0}$. We have used here the explicit form of $p$ in $B$.

Comparing \eqref{j58} and \eqref{j57}, it comes
\begin{equation*}
g_{-} ( \rho_{-}^{E_{0}} ) \cdot g_{+} ( \rho_{+}^{E_{0}} ) \leq o_{t \to + \infty} ( 1 ) + o_{E \to E_{0}}^{t} ( 1 ) .
\end{equation*}
Taking $t$ large enough and then $E$ close enough to $E_{0}$, it gives
\begin{equation*}
g_{-} ( \rho_{-}^{E_{0}} ) \cdot g_{+} ( \rho_{+}^{E_{0}} ) \leq 0 ,
\end{equation*}
since the left hand side of the above equation is independent of $t , E$. Eventually, this is in contradiction with $g_{+}^{k} \cdot g_{-}^{\ell} > 0$ for all $k , \ell \in \{ 1 ,\ldots , K \}$, and the proposition follows.

\Subsection{Proof of Lemma \ref{b17}} \label{s20}

This proof relies on \eqref{b19} and we use the definitions of Section \ref{s52}. To keep the notations of Theorem \ref{a32}, we assume in the sequel that $y_{j} = ( \varepsilon , 0 , \ldots , 0 )$ and $H_{j} = \{ y_{1} = \varepsilon \}$. It can always be done by a linear change of coordinates. First, we estimate the behavior of the asymptotic directions. Let $x \in \varepsilon \S^{n-1} \cap \supp a_{+}$. Then, from \eqref{a60}, $g_{+} ( \rho^{+}_{x} ) = \Phi (x)$ where $\Phi \in C^{\infty}$, $\Phi (0) = 0$ and $d \Phi (0) = I d_{\R^{n}}$. In particular,
\begin{equation} \label{b21}
g_{+} ( \rho^{+}_{x} ) = x + \CO ( x^{2} ) = x + \CO ( \varepsilon^{2} ) .
\end{equation}
Now, let $y \in H_{j} \cap \supp a_{- , 0}^{j}$. From \eqref{b14}--\eqref{b15}, we deduce
\begin{equation} \label{n12}
y = Y ( y ) + \CO_{\varepsilon} ( \nu ) .
\end{equation}
Using the symmetry of the Hamiltonian curves $( x ( t ) , \xi ( t ) ) \longmapsto ( x ( - t ) , - \xi ( - t ) )$, we also have $g_{-} ( \rho^{-}_{y} ) = \Phi ( y )$. Therefore, we obtain as in \eqref{b21}
\begin{equation} \label{b22}
g_{-} ( \rho^{-}_{y} ) = y + \CO ( y^{2} ) = Y (y) + \CO ( \varepsilon^{2} ) + \CO_{\varepsilon} ( \nu ) ,
\end{equation}
since $\vert Y ( y ) \vert = \varepsilon$. In particular, \eqref{b21} and \eqref{b22} imply
\begin{equation} \label{b38}
\vert g_{-} ( \rho^{-}_{y} ) \vert = \varepsilon + \CO ( \varepsilon^{2} ) + \CO_{\varepsilon} ( \nu ) \qquad \text{and} \qquad \vert g_{+} ( \rho^{+}_{x} ) \vert = \varepsilon + \CO ( \varepsilon^{2} ) ,
\end{equation}
whereas, thanks to \eqref{n13},
\begin{align}
\big( i \lambda ( g_{+} ( \rho^{+}_{x} ) \cdot g_{-} ( \rho^{-}_{y} ) ) \big)^{- S ( z , h ) / \lambda} ={}& \big( i \lambda ( g_{+} ( \rho^{+}_{x} ) \cdot g_{-} ( \rho^{-}_{y} ) ) \big)^{- \frac{n}{2} + i \frac{z - E_{0}}{\lambda h}}   \nonumber \\
={}& \big( i \lambda x \cdot Y ( y ) \big)^{- \frac{n}{2} + i \frac{z - E_{0}}{\lambda h}} + \CO ( \varepsilon^{1 - n} ) + \CO_{\varepsilon} ( \nu ) . \label{j32}
\end{align}
Note that, in the previous equation, the quantities $( g_{+} ( \rho^{+}_{x} ) \cdot g_{-} ( \rho^{-}_{y} ) )$ and $( x \cdot Y (y) )$ do not vanish from \ref{h4} for $V_{\pm}^{\varepsilon}$ sufficiently close to $\pi_{x} ( \CH_{{\rm tang} , \pm}^{\varepsilon} )$. Taking the modulus and using $\im z = \CO ( h \vert \ln h \vert^{-1} )$, we get
\begin{align}
\Big\vert \big( i \lambda ( g_{+} ( \rho^{+}_{x} ) \cdot g_{-} ( \rho^{-}_{y} ) ) \big)^{- S ( z , h ) / \lambda} \Big\vert ={}& e^{- \frac{\pi}{2} \frac{\re z - E_{0}}{\lambda h} \sgn ( x \cdot Y (y) )} \big\vert ( \lambda x \cdot Y (y) ) \big\vert^{- \frac{n}{2} - \frac{\im z}{\lambda h}}   \nonumber \\
&\qquad \qquad \quad + \CO ( \varepsilon^{1 - n} ) + \CO_{\varepsilon} ( \nu ) \nonumber \\
={}& \lambda^{- \frac{n}{2}} e^{- \frac{\pi}{2} \frac{\re z - E_{0}}{\lambda h} \sgn ( x \cdot Y (y) )} \big\vert ( x \cdot Y (y) ) \big\vert^{- \frac{n}{2}}  \nonumber \\
&\qquad \qquad \quad + \CO ( \varepsilon^{1 - n} ) + \CO_{\varepsilon} ( \nu ) + \CO_{\varepsilon , \nu} \big( \vert \ln h \vert^{-1} \big) . \label{b39}
\end{align}
Again, using that $\im z = \CO ( h \vert \ln h \vert^{-1} )$, we immediately obtain
\begin{equation} \label{j33}
\big\vert \Gamma ( S ( z , h ) / \lambda ) \big\vert = \Big\vert \Gamma \Big( \frac{n}{2} - i \frac{\re z - E_{0}}{\lambda h} \Big) \Big\vert + \CO \big( \vert \ln h \vert^{-1} \big) .
\end{equation}

Otherwise, \eqref{a57} implies
\begin{equation} \label{b23}
\big\vert \det \nabla^{2}_{y^{\prime} , y^{\prime} } \varphi_{-} (y) \big\vert^{\frac{1}{2}} = \Big\vert \det \Big( - \frac{\lambda}{2} I d_{\R^{n-1}} + \CO (y) \Big) \Big\vert^{\frac{1}{2}} = \Big( \frac{\lambda}{2} \Big)^{\frac{n-1}{2}} + \CO ( \varepsilon ) + \CO ( \nu ) ,
\end{equation}
since $\vert y \vert = \varepsilon + \CO ( \nu )$. The same way, we obtain
\begin{equation} \label{b24}
\big\vert \partial_{\xi_{1}} p ( \rho_{y}^{-} ) \big\vert^{\frac{1}{2}} = \sqrt{2 \big\vert \partial_{x_{1}} \varphi_{-} (y) \big\vert} = \sqrt{\lambda \vert y_{1} \vert + \CO ( y^{2} )} = \sqrt{\lambda \varepsilon} + \CO ( \varepsilon^{\frac{3}{2}} ) + \CO_{\varepsilon} ( \nu ) .
\end{equation}
since $y = y_{j} + \CO ( \nu )$.

We now consider $x (t) = \pi_{x} ( \exp ( t H_{p} ) ( \rho_{x}^{+} ) )$. Then, there exists $C > 0$ such that $\vert x (t) \vert \leq C e^{\lambda t} \vert x \vert$ for $t \leq 0$ (this follows from \eqref{a62} for instance). Using \eqref{a57}, it yields
\begin{equation*}
\Delta \varphi_{+} ( x ( s ) ) - n \lambda / 2 = \CO ( e^{\lambda s} \varepsilon ) .
\end{equation*}
Then, we get
\begin{equation} \label{b20}
e^{\int_{0}^{- \infty} ( \Delta \varphi_{+} ( x ( s ) ) - n \lambda /2 ) \, d s} = 1 + \CO ( \varepsilon ) .
\end{equation}

It remains to study the last term in \eqref{b19}. For that, we will follow the approach of \cite[Lemma 2.1]{HeSj85_01} and \cite[Proposition C.1]{ALBoRa08_01}. First, using the notation $( y (t) , \eta (t) ) = \exp ( t H_{p} ) ( \rho_{y}^{-} )$, we have
\begin{equation*}
\frac{\partial y ( t , y^{\prime} , \eta^{\prime} )}{\partial t} \vert_{\eta^{\prime} = \partial_{y^{\prime}} \varphi_{-} (y)} = 2 \eta (t) = 2 \nabla \varphi_{-} ( y (t) ) .
\end{equation*}
Combining with \eqref{a57} and \eqref{d2}, this implies
\begin{equation*}
\frac{\partial y ( t , y^{\prime} , \eta^{\prime} )}{\partial t} \vert_{\eta^{\prime} = \partial_{y^{\prime}} \varphi_{-} (y)} = - \lambda y (t) + \CO ( y (t)^{2} ) = - \lambda g_{-} ( \rho_{y}^{-} ) e^{- \lambda t} + o_{t \to + \infty} ( e^{- \lambda t } ) .
\end{equation*}
Eventually, $g_{-} ( \rho^{-}_{y} ) = \Phi ( y )$ gives $g_{-} ( \rho^{-}_{y} ) = y + \CO ( y^{2} ) = y_{j} + \CO ( \varepsilon^{2} + \nu )$ and then
\begin{equation} \label{b26}
\frac{\partial y ( t , y^{\prime} , \eta^{\prime} )}{\partial t} \vert_{\eta^{\prime} = \partial_{y^{\prime}} \varphi_{-} (y)} = - \lambda
\left( \begin{array}{c}
\varepsilon  \\
0 \\
\vdots \\
0
\end{array} \right)
e^{- \lambda t} + \CO ( \varepsilon^{2} + \nu ) e^{- \lambda t} + o_{t \to + \infty} ( e^{- \lambda t } ) .
\end{equation}
We now compute the derivative of $y ( t , y^{\prime} , \eta^{\prime} )$ with respect to $y^{\prime}$ at $\eta^{\prime} = \partial_{y^{\prime}} \varphi_{-} (y)$. For that, it is equivalent to compute the evolution of a tangent vector $( \delta_{y} , \delta_{\eta} )$ of
\begin{equation*}
\Lambda_{\eta^{\prime}} = \big\{ ( \varepsilon , y^{\prime} , f_{-} ( \varepsilon , y^{\prime} , \eta^{\prime} ) , \eta^{\prime} ) ; \ y^{\prime} \in \R^{n-1} \big\} ,
\end{equation*}
along $\exp ( t H_{p} ) ( \rho_{y}^{-} )$. Indeed, $\Lambda_{\eta^{\prime}}$ is the space of the initial condition of $y ( t , y^{\prime} , \eta^{\prime} )$ given in Theorem \ref{a32}. From the definition of $\Lambda_{\eta^{\prime}}$, we have
\begin{equation} \label{b29}
\delta_{y} = ( 0 , \delta_{y^{\prime}} ) \qquad \text{and} \qquad \delta_{\eta} = \big( \partial_{y^{\prime}} f_{-} ( \varepsilon , y^{\prime} , \eta^{\prime} ) \cdot \delta_{y^{\prime}} , 0 \big) .
\end{equation}
We then follow the proof of \cite[Lemma 2.1]{HeSj85_01}. There exist symplectic local coordinates $( k , \kappa ) \in T^{*} \R^{n}$ centered at $( 0 , 0 )$ such that $\Lambda_{-}^{0}$ (resp. $\Lambda_{+}^{0}$) is given by $k = 0$ (resp. $\kappa =0$) and
\begin{align}
k =& \frac{1}{\sqrt{\lambda}} \Big( \eta + \frac{\lambda}{2} y \Big) + \CO \big( ( y , \eta )^{2} \big) ,   \label{b33} \\
\kappa =& \frac{1}{\sqrt{\lambda}} \Big( \eta - \frac{\lambda}{2} y \Big) + \CO \big( ( y , \eta )^{2} \big) . \label{b34}
\end{align}
Then, $p ( y , \eta ) = E_{0} + A ( k , \kappa ) k \cdot \kappa$ with $A (0,0) = \lambda I d_{\R^{n}}$. In these coordinates, \eqref{b29} becomes
\begin{align}
\delta_{k} &= \frac{1}{2 \sqrt{\lambda}} \big( 2 \partial_{y^{\prime}} f_{-} ( \varepsilon , y^{\prime} , \eta^{\prime} ) \cdot \delta_{y^{\prime}} , \lambda \delta_{y^{\prime}} \big) + \CO ( \varepsilon + \nu ) \vert \delta_{y^{\prime}} \vert ,    \label{b30} \\
\delta_{\kappa} &= \frac{1}{2 \sqrt{\lambda}} \big( 2 \partial_{y^{\prime}} f_{-} ( \varepsilon , y^{\prime} , \eta^{\prime} ) \cdot \delta_{y^{\prime}} , - \lambda \delta_{y^{\prime}} \big) + \CO ( \varepsilon + \nu ) \vert \delta_{y^{\prime}} \vert .   \label{b31}
\end{align}
Moreover, the evolution of a tangent vector $( \delta_{k} , \delta_{\kappa} )$ at $\exp ( t H_{p} ) ( \rho_{y}^{-} )$ satisfies the following equation
\begin{equation}  \label{b27}
\frac{d}{d t}
\left( \begin{array}{c}
\delta_{k} \\
\delta_{\kappa}
\end{array} \right)
=
\left( \begin{array}{cc}
\lambda + \CO \big( ( \varepsilon + \nu ) e^{-\lambda t} \big) & 0 \\
\CO \big( ( \varepsilon + \nu ) e^{-\lambda t} \big) & - \lambda + \CO \big( ( \varepsilon + \nu ) e^{-\lambda t} \big)
\end{array} \right) \left( \begin{array}{c}
\delta_{k} \\
\delta_{\kappa}
\end{array} \right) .
\end{equation}
The first equation
\begin{equation*}
\partial_{t} \delta_{k} = \big( \lambda + \CO \big( ( \varepsilon + \nu ) e^{-\lambda t} \big) \big) \delta_{k} ,
\end{equation*}
implies
\begin{equation} \label{b28}
\delta_{k} (t) = e^{\lambda t} \delta_{k} (0) + \CO ( \varepsilon + \nu ) e^{\lambda t} \vert \delta_{k} (0) \vert ,
\end{equation}
by the usual theory of the ordinary differential equations. On the other hand, $\Lambda_{-}$ and $\Lambda_{\eta^{\prime}}$ intersect transversally along a Hamiltonian curve from \cite[Lemma 5.6]{BoFuRaZe07_01}. Then, \cite[Lemma 2.1]{HeSj85_01} implies that, for $( \delta_{k} (0) , \delta_{\kappa} (0) ) \in T \Lambda_{\eta^{\prime}}$, we have
\begin{equation*}
\delta_{\kappa} (t) = B (t) \delta_{k} (t) ,
\end{equation*}
where $B (t) = \CO_{\varepsilon , \nu} ( e^{- \lambda t} )$. Then, \eqref{b28} gives
\begin{equation} \label{b32}
\delta_{\kappa} (t) = \CO_{\varepsilon , \nu} (1) \vert \delta_{k} (0) \vert .
\end{equation}
Inverting the change of variables \eqref{b33}--\eqref{b34}, the expansions \eqref{b28} and \eqref{b32} yield
\begin{equation} \label{g81}
\delta_{y} (t) = \frac{1}{\sqrt{\lambda}} e^{\lambda t} \delta_{k} (0) + \CO ( \varepsilon + \nu ) e^{\lambda t} \vert \delta_{k} (0) \vert + \CO_{\varepsilon , \nu} (1) \vert \delta_{k} (0) \vert .
\end{equation}
We will now estimate $\partial_{y^{\prime}} f_{-}$. By definition, we have
\begin{equation*}
f_{-} ( \varepsilon , y^{\prime} , \partial_{y^{\prime}} \varphi_{-} (y) ) = \partial_{y_{1}} \varphi_{-} (y) = - \frac{\lambda}{2} y_{1} + \CO ( y^{2} ) = - \frac{\lambda}{2} \varepsilon + \CO ( \varepsilon^{2} + \nu ) .
\end{equation*}
Furthermore, differentiating $p ( \varepsilon , y^{\prime} , f_{-} , \eta^{\prime} ) = E_{0}$, we get
\begin{equation*}
2 f_{-} ( \varepsilon , y^{\prime} , \partial_{y^{\prime}} \varphi_{-} (y) ) \partial_{y^{\prime}} f_{-} ( \varepsilon , y^{\prime} , \partial_{y^{\prime}} \varphi_{-} (y) ) = - \partial_{y^{\prime}} V ( \varepsilon , y^{\prime} ) = \frac{\lambda^{2}}{2} y^{\prime} + \CO ( ( \varepsilon , y^{\prime} )^{2} ) = \CO ( \varepsilon + \nu ) .
\end{equation*}
The two last estimates imply
\begin{equation} \label{b37}
\partial_{y^{\prime}} f_{-} ( \varepsilon , y^{\prime} , \partial_{y^{\prime}} \varphi_{-} (y) ) = \CO (1) + \CO_{\varepsilon} ( \nu ) .
\end{equation}
Then, \eqref{b30}, \eqref{g81} and \eqref{b37} give
\begin{equation*}
\delta_{y} (t) = \frac{1}{2 \lambda} e^{\lambda t} \big( 2 \partial_{y^{\prime}} f_{-} ( \varepsilon , y^{\prime} , \eta^{\prime} ) \cdot \delta_{y^{\prime}} , \lambda \delta_{y^{\prime}} \big) + \CO ( \varepsilon ) e^{\lambda t} \vert \delta_{y^{\prime}} \vert + \CO_{\varepsilon} ( \nu ) e^{\lambda t} \vert \delta_{y^{\prime}} \vert + \CO_{\varepsilon , \nu} (1) \vert \delta_{y^{\prime}} \vert .
\end{equation*}
In other words, we have proved
\begin{align}
\frac{\partial y ( t , y^{\prime} , \eta^{\prime} )}{\partial y^{\prime}} \vert_{\eta^{\prime} = \partial_{y^{\prime}} \varphi_{-} (y)} = \frac{1}{2 \lambda} e^{\lambda t} 
\left( \begin{array}{c}
2 \partial_{y^{\prime}} f_{-} ( \varepsilon , y^{\prime} , \partial_{y^{\prime}} \varphi_{-} (y) ) \\
\lambda I d_{\R^{n-1}}
\end{array} \right)
+ {} & \CO ( \varepsilon ) e^{\lambda t} \nonumber \\
&+ \CO_{\varepsilon} ( \nu ) e^{\lambda t} + \CO_{\varepsilon , \nu} (1) . \label{b35}
\end{align}
Combining \eqref{b26}, \eqref{b37} and \eqref{b35}, we get
\begin{equation*}
\Big\vert \det \frac{\partial y ( t , y^{\prime} , \eta^{\prime} )}{\partial ( t , y^{\prime} )}\vert_{\eta^{\prime} = \partial_{y^{\prime}} \varphi_{-} ( y )} \Big\vert = \frac{\lambda \varepsilon}{2^{n-1}} e^{( n - 2 ) \lambda t} + \big( \CO ( \varepsilon^{2} ) + \CO_{\varepsilon} ( \nu ) + o_{t \to + \infty}^{\varepsilon , \nu} ( 1 ) \big) e^{( n - 2 ) \lambda t} ,
\end{equation*}
which finally implies
\begin{equation} \label{b36}
\lim_{t\to + \infty} \frac{e^{\frac{n - 2}{2} \lambda t}}{{\sqrt{\Big\vert \det \frac{\partial y( t , y^{\prime} , \eta^{\prime} )}{\partial ( t , y^{\prime} )}\vert_{\eta^{\prime} = \partial_{y^{\prime}} \varphi_{-} ( y )} \Big\vert}}} = \frac{2^{\frac{n - 1}{2}}}{\sqrt{\lambda \varepsilon}} + \CO ( \sqrt{\varepsilon} ) + \CO_{\varepsilon} ( \nu ) .
\end{equation}

Eventually, \eqref{j2} follows from \eqref{b19} together with \eqref{b38}, \eqref{j32}, \eqref{b23}, \eqref{b24}, \eqref{b20} and \eqref{b36}. The same way, \eqref{j3} is a consequence of the estimates \eqref{b38}, \eqref{b39}, \eqref{j33}, \eqref{b23}, \eqref{b24}, \eqref{b20} and \eqref{b36}.

\Subsection{A self-adjoint operator for Example \ref{k87}} \label{s90}

\begin{figure}
\begin{center}
\begin{picture}(0,0)%
\includegraphics{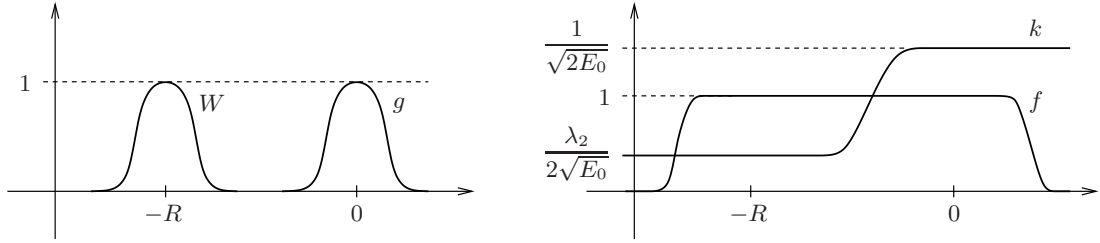}%
\end{picture}%
\setlength{\unitlength}{987sp}%
\begingroup\makeatletter\ifx\SetFigFont\undefined%
\gdef\SetFigFont#1#2#3#4#5{%
  \reset@font\fontsize{#1}{#2pt}%
  \fontfamily{#3}\fontseries{#4}\fontshape{#5}%
  \selectfont}%
\fi\endgroup%
\begin{picture}(27516,6066)(-27182,-5194)
\put(-23249,-4711){\makebox(0,0)[b]{\smash{{\SetFigFont{9}{10.8}{\rmdefault}{\mddefault}{\updefault}$- R$}}}}
\put(-17474,-1861){\makebox(0,0)[lb]{\smash{{\SetFigFont{9}{10.8}{\rmdefault}{\mddefault}{\updefault}$g$}}}}
\put(-22349,-1936){\makebox(0,0)[lb]{\smash{{\SetFigFont{9}{10.8}{\rmdefault}{\mddefault}{\updefault}$W$}}}}
\put(-26549,-1432){\makebox(0,0)[rb]{\smash{{\SetFigFont{9}{10.8}{\rmdefault}{\mddefault}{\updefault}$1$}}}}
\put(-1499,-61){\makebox(0,0)[lb]{\smash{{\SetFigFont{9}{10.8}{\rmdefault}{\mddefault}{\updefault}$k$}}}}
\put(-11999,-436){\makebox(0,0)[rb]{\smash{{\SetFigFont{9}{10.8}{\rmdefault}{\mddefault}{\updefault}$\displaystyle \frac{1}{\sqrt{2 E_{0}}}$}}}}
\put(-11999,-3136){\makebox(0,0)[rb]{\smash{{\SetFigFont{9}{10.8}{\rmdefault}{\mddefault}{\updefault}$\displaystyle \frac{\lambda_{2}}{2 \sqrt{E_{0}}}$}}}}
\put(-8549,-4711){\makebox(0,0)[b]{\smash{{\SetFigFont{9}{10.8}{\rmdefault}{\mddefault}{\updefault}$- R$}}}}
\put(-3374,-4711){\makebox(0,0)[b]{\smash{{\SetFigFont{9}{10.8}{\rmdefault}{\mddefault}{\updefault}$0$}}}}
\put(-1499,-1861){\makebox(0,0)[lb]{\smash{{\SetFigFont{9}{10.8}{\rmdefault}{\mddefault}{\updefault}$f$}}}}
\put(-11999,-1786){\makebox(0,0)[rb]{\smash{{\SetFigFont{9}{10.8}{\rmdefault}{\mddefault}{\updefault}$1$}}}}
\put(-18374,-4711){\makebox(0,0)[b]{\smash{{\SetFigFont{9}{10.8}{\rmdefault}{\mddefault}{\updefault}$0$}}}}
\end{picture}%
\end{center}
\caption{The functions $f , g , k , W$.} \label{f61}
\end{figure}

We construct here a self-adjoint pseudodifferential operator satisfying the assumptions of Example \ref{k87}. For that, we use the symbol build in Example \ref{b81}. Thus, we set
\begin{equation*}
q_{1} ( x_{1} , \xi_{1} ) = \psi ( x_{1} ) e^{- x_{1}^{2}}  \big( 1 - \xi_{1}^{2} ( \xi_{1} + 1 / 2 ) \big) e^{- \xi_{1}^{4} / 5} \in S ( 1 ) ,
\end{equation*}
where $\psi \in C^{\infty}_{0} ( \R ; [ 0 , 1 ] )$ is equal to $1$ on a sufficiently large neighborhood of $0$. In particular, there exists $C_{1} > 0$ such that
\begin{equation*}
q_{1} ( x_{1} , \xi_{1} ) + C_{1} > 0,
\end{equation*}
for all $( x_{1} , \xi_{1} ) \in T^{*} \R$. Next, let
\begin{equation*}
q_{2} ( x_{1} , \xi_{1} ) = W ( x_{1} ) + q_{1} ( x_{1} , \xi_{1} ) + C_{1} \in S ( 1 ) ,
\end{equation*}
where the potential $W \in C^{\infty}_{0} ( \R ; [ 0 , 1 ] )$ is a non-degenerate bump illustrated in Figure \ref{f61} and satisfying
\begin{equation*}
W ( x_{1} ) = 1 - \frac{\lambda_{2}^{2}}{4} ( x_{1} + R )^{2} + \CO \big( ( x_{1} + R )^{3} \big) ,
\end{equation*}
near $- R$ for some $R > 0$ large enough. Here, $\lambda_{2} > 0$ can be chosen arbitrarily. Note that $q_{2} ( x_{1} , \xi_{1} ) > 0$ for all $( x_{1} , \xi_{1} ) \in T^{*} \R$. We then define the energy $E_{0} = 1 + C_{1} > 0$ and the symbol
\begin{equation} \label{k96}
p_{0} ( x , \xi ) = \xi^{2} + q_{2} ( x_{1} , \xi_{1} ) f ( x_{1} ) g ( x_{2} k ( x_{1} ) ) ,
\end{equation}
for $( x, \xi ) \in T^{*} \R^{2}$. The graph of the smooth functions $f , g , k : \R \longmapsto \R$ is drawn in Figure \ref{f61}. Moreover, we assume that $g ( y ) = 1 - y^{2} + \CO ( y^{3} )$ near $0$. The symbol $p_{0}$ is real valued, belongs to $S ( \< \xi \>^{2} )$ and is classically elliptic in this class (this means $p_{0} ( x , \xi ) \geq \xi^{2} - C$). Moreover, $p_{0} ( x , \xi ) = \xi^{2}$ for $x$ outside a compact set.

\begin{lemma}\sl \label{k97}
For all $( x , \xi ) \in K ( E_{0} )$, we have $x_{2} = \xi_{2} = 0$.
\end{lemma}

\begin{proof}
Let $\rho ( t ) = ( x ( t ) , \xi ( t ) )$ be a trapped Hamiltonian trajectory in $p_{0}^{- 1} ( E_{0} )$. From \eqref{k96} and the definition of the Hamiltonian vector field \eqref{d27}, we get
\begin{equation} \label{k98}
\left\{ \begin{aligned}
&\partial_{t} x_{2} = \partial_{\xi_{2}} p_{0} ( x , \xi ) = 2 \xi_{2} \text{ has the sign of } \xi_{2} ,   \\
&\partial_{t} \xi_{2} = - \partial_{x_{2}} p_{0} ( x , \xi ) = - q_{2} ( x_{1} , \xi_{1} ) f ( x_{1} ) k ( x_{1} ) g^{\prime} ( x_{2} k ( x_{1} ) ) \text{ has the sign of } x_{2} .
\end{aligned} \right.
\end{equation}
We first assume that $x_{2} ( 0 ) > 0$.

In addition, we suppose that $\xi_{2} ( 0 ) > 0$. Then, \eqref{k98} implies that $x_{2} ( t ) > 0$ and $\xi_{2} ( t ) > 0$ for all $t \geq 0$. Another application of \eqref{k98} gives $\xi_{2} ( t ) \geq \xi_{2} ( 0 )$ and eventually $x_{2} ( t ) \geq x_{2} ( 0 ) + t \xi_{2} ( 0 )$ for all $t \geq 0$. It implies that this trajectory is not trapped and gives a contradiction. The same way, one can exclude the assertion $\xi_{2} ( 0 ) < 0$ considering the limit $t \to - \infty$. Thus, we have $\xi_{2} ( 0 ) = 0$. Assume that there exists $s > 0$ such that $\xi_{2} ( s ) \neq 0$. Let
\begin{equation*}
t = \inf \{ s > 0 ; \ \xi_{2} ( s ) \neq 0 \} .
\end{equation*}
From \eqref{k98}, we have $x_{2} ( t ) = x_{2} ( 0 ) > 0$. By continuity, there exists $s > t$ (close to $t$) such that $x_{2} ( s ) \geq x_{2} ( 0 ) / 2 > 0$ and $\xi_{2} ( s ) \neq 0$. By the previous paragraph, this is impossible. This shows that $\xi_{2} ( t ) = 0$ for all $t \geq 0$. Working symmetrically for $t < 0$, we obtain
\begin{equation} \label{k99}
\xi_{2} ( t ) = 0 \qquad \text{and} \qquad x_{2} ( t ) = x_{2} ( 0 ) ,
\end{equation}
for all $t \in \R$.

Combining \eqref{k98} and \eqref{k99}, we get
\begin{equation*}
0 = \partial_{t} \xi_{2} ( t ) = - q_{2} ( x_{1} , \xi_{1} ) f ( x_{1} ) k ( x_{1} ) g^{\prime} ( x_{2} k ( x_{1} ) ) ,
\end{equation*}
for all $t \in \R$. From the form of the functions $f , g , k , q_{2}$, it yields
\begin{equation} \label{m72}
f ( x_{1} ) = 0 \qquad \text{or} \qquad g ( x_{2} k ( x_{1} ) ) = 0 ,
\end{equation}
for all $t \in \R$. On the other hand, the Hamiltonian equations \eqref{d27} write
\begin{equation*}
\partial_{t} x_{1} = \partial_{\xi_{1}} p_{0} ( x , \xi ) = 2 \xi_{1} + ( \partial_{\xi_{1}} q_{2} ) ( x_{1} , \xi_{1} ) f ( x_{1} ) g ( x_{2} k ( x_{1} ) ) ,
\end{equation*}
and
\begin{align*}
\partial_{t} \xi_{1} = - \partial_{x_{1}} p_{0} ( x , \xi ) ={}& - ( \partial_{x_{1}} q_{2} ) ( x_{1} , \xi_{1} ) f ( x_{1} ) g ( x_{2} k ( x_{1} ) ) - q_{2} ( x_{1} , \xi_{1} ) f^{\prime} ( x_{1} ) g ( x_{2} k ( x_{1} ) ) \\
&- q_{2} ( x_{1} , \xi_{1} ) f ( x_{1} ) k^{\prime} ( x_{1} ) g^{\prime} ( x_{2} k ( x_{1} ) ) .
\end{align*}
Since $f ( y ) = 0$ (resp. $g ( y ) = 0$) implies $f^{\prime} ( y ) = 0$ (resp. $g^{\prime} ( y ) = 0$), \eqref{m72} gives $\partial_{t} x_{1} ( t ) = 2 \xi_{1} ( t )$ and $\partial_{t} \xi_{1} ( t ) = 0$ for all $t \in \R$. Thus,
\begin{equation} \label{m73}
\xi_{1} ( t ) = \xi_{1} ( 0 ) \qquad \text{and} \qquad x_{1} ( t ) = x_{1} ( 0 ) + t \xi_{1} ( 0 ) ,
\end{equation}
for all $t \in \R$. Lastly, \eqref{k96}, \eqref{k99}, \eqref{m72} and \eqref{m73} imply
\begin{equation*}
0 < E_{0} = p_{0} ( x , \xi ) = \xi^{2} = \xi_{1}^{2} ( 0 ) .
\end{equation*}
and then $\xi_{1} ( 0 ) \neq 0$. As consequence, \eqref{m73} shows that the trajectory $t \longmapsto \rho ( t )$ is non-trapped which is impossible.

Starting from the assumption $x_{2} ( 0 ) > 0$, we have just obtained a contradiction. Working the same way if $x_{2} ( 0 ) < 0$, we deduce $x_{2} ( 0 ) = 0$ and then $x_{2} ( t ) = 0$ for all $t \in \R$. Eventually, \eqref{d27} yields
\begin{equation*}
\xi_{2} ( t ) = \partial_{t} x_{2} ( t ) / 2 = 0 ,
\end{equation*}
and the lemma follows.
\end{proof}

\begin{figure}
\begin{center}
\begin{picture}(0,0)%
\includegraphics{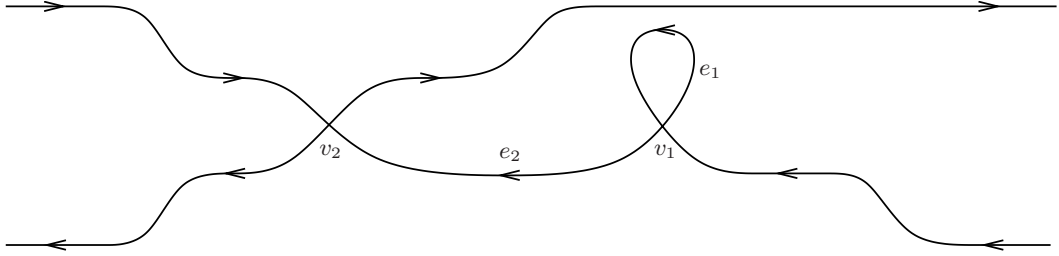}%
\end{picture}%
\setlength{\unitlength}{987sp}%
\begingroup\makeatletter\ifx\SetFigFont\undefined%
\gdef\SetFigFont#1#2#3#4#5{%
  \reset@font\fontsize{#1}{#2pt}%
  \fontfamily{#3}\fontseries{#4}\fontshape{#5}%
  \selectfont}%
\fi\endgroup%
\begin{picture}(26488,6264)(-34543,-7093)
\put(-16799,-2686){\makebox(0,0)[b]{\smash{{\SetFigFont{9}{10.8}{\rmdefault}{\mddefault}{\updefault}$e_{1}$}}}}
\put(-26324,-4711){\makebox(0,0)[b]{\smash{{\SetFigFont{9}{10.8}{\rmdefault}{\mddefault}{\updefault}$v_{2}$}}}}
\put(-17924,-4711){\makebox(0,0)[b]{\smash{{\SetFigFont{9}{10.8}{\rmdefault}{\mddefault}{\updefault}$v_{1}$}}}}
\put(-21824,-4786){\makebox(0,0)[b]{\smash{{\SetFigFont{9}{10.8}{\rmdefault}{\mddefault}{\updefault}$e_{2}$}}}}
\end{picture}%
\end{center}
\caption{The energy surface $p_{0}^{- 1} ( E_{0} )$ restricted to $x_{2} = \xi_{2} = 0$.} \label{f62}
\end{figure}

The symbol $p_{0}$ restricted to $x_{2} = \xi_{2} = 0$ takes the form
\begin{equation*}
p_{0} {}_{\vert_{x_{2} = \xi_{2} = 0}} ( x , \xi ) = \xi_{1}^{2} + q_{2} ( x_{1} , \xi_{1} ) f ( x_{1} ) .
\end{equation*}
The energy surface $p_{0}^{- 1} ( E_{0} )$ restricted to $x_{2} = \xi_{2} = 0$ is illustrated in Figure \ref{f62}. In particular, the trapped set $K ( E_{0} )$ consists of two fixed points $( v_{1} , 0 ) = ( 0 , 0 )$ and $( v_{2} , 0 ) = ( - R , 0 )$, a homoclinic trajectory $e_{1}$ around $v_{1}$ (which is that of Example \ref{b81} drawn in Figure \ref{f14}) and a heteroclinic trajectory $e_{2}$ from $v_{1}$ to $v_{2}$.

It remains to study the behavior of $p_{0}$ in a vicinity of the critical points. Near $( v_{1} , 0 ) = ( 0 , 0 )$, \eqref{k96} gives
\begin{equation*}
p_{0} ( x , \xi ) = E_{0} + \xi_{1}^{2} / 2 + \xi_{2}^{2} - x_{1}^{2} - x_{2}^{2} / 2 + \CO \big( ( x , \xi )^{3} \big) .
\end{equation*}
Thus, applying the symplectic homothety $( x , \xi ) \longmapsto ( y , \eta ) = ( \sqrt{2} x_{1} , x_{2} , \xi_{1} / \sqrt{2} , \xi_{2} )$, we get
\begin{equation*}
p_{0} ( x , \xi ) = E_{0} + \eta^{2} - y^{2} / 2 + \CO \big( ( y , \eta )^{3} \big) .
\end{equation*}
In other words, $( v_{1} , 0 )$ is an isotropic fixed point with coefficient $\lambda_{1} = \sqrt{2}$. On the other hand, we have near $( v_{2} , 0 ) = ( - R , 0 )$
\begin{equation*}
p_{0} ( x , \xi ) = E_{0} + \eta^{2} - \frac{\lambda_{2}^{2}}{4} y^{2} + \CO \big( ( y , \eta )^{3} \big) .
\end{equation*}
with $( y , \eta ) = ( x_{1} + R , x_{2} , \xi_{1} , \xi_{2} )$. Thus, $( v_{2} , 0 )$ is an isotropic fixed point with coefficient $\lambda_{2}$ which can be chosen arbitrarily small.

Summing up, we define
\begin{equation}
P = - ( 1 - \varphi ) ( x ) h^{2} \Delta ( 1 - \varphi ) ( x ) + \chi ( x ) \Op ( p_{0} ) \chi ( x ) ,
\end{equation}
where $\chi , \varphi \in C^{\infty}_{0} ( \R^{2} ; [ 0 , 1 ] )$ are such that $\chi = 1$ on a sufficiently large neighborhood of $0$ and $\chi^{2} ( x ) + ( 1 - \varphi )^{2} ( x ) = 1$ for all $x \in \R^{2}$. From \eqref{k96}, this operator is a pseudodifferential operator whose symbol $p \in S ( \< \xi \>^{2} )$ is real, classically elliptic and verifies
\begin{equation*}
p = p_{0} + S ( h^{2} \< x \>^{- \infty} ) .
\end{equation*}
Thus, we have constructed a self-adjoint operator on $H^{2} ( \R^{2} )$ satisfying the assumptions of Example \ref{k87} and such that $\lambda_{2}$ can be taken arbitrarily small. It is not a Schr\"{o}dinger operator, but it is within the scope of Remark \ref{c13}. In particular, the resonances of $P$ are well-defined using the theory of Sj\"{o}strand and Zworski \cite{SjZw91_01} since $P$ coincides with $- h^{2} \Delta$ outside a compact set. Note also that the transversality hypothesis \ref{h19} follows from the repulsive character of $p_{0}$ in the variable $x_{2}$ (see the proof of Lemma \ref{k97}).

\section{Spectral radius of $\CT_{0}$ and $\CT$} \label{s70}

In this part, we show that $\tau \longmapsto \CA_{0} ( \tau )$ is continuous on $\R$ and that $\Vert \CT_{0} ( \tau )^{k} \Vert^{\frac{1}{k}}$ converges uniformly to $\CA_{0} ( \tau )$ on the compact set in $\tau$. In the sequel, we will assume that $\mes_{\S^{n - 1}} ( \CH_{\rm tang}^{- \infty} ) > 0$ since all the considered operators vanish in the other case (see Remark \ref{b69} $i)$).

We will first work on the vicinity of a point $\tau_{0} \in \R$. From \ref{h4}, \eqref{b70} and Proposition \ref{b53}, there exists $C_{1} > 0$ such that
\begin{equation} \label{c2}
1 / C_{1} \leq \CT_{0} ( \tau ) ( \omega , \widetilde{\omega} ) \leq C_{1} ,
\end{equation}
for all $\vert \tau - \tau_{0} \vert \leq 1$ and $\omega , \widetilde{\omega} \in \CH_{\rm tang}^{- \infty}$. Moreover, from \eqref{b70}, there exists $C_{2} > 0$ such that
\begin{equation} \label{b99}
\CT_{0} ( \tau ) ( \omega , \widetilde{\omega} ) = \CT_{0} ( \sigma ) ( \omega , \widetilde{\omega} ) + \CR ( \tau , \sigma ) ( \omega , \widetilde{\omega} ) ,
\end{equation}
with
\begin{equation} \label{c1}
\big\vert \CR ( \tau , \sigma ) ( \omega , \widetilde{\omega} ) \big\vert \leq C_{2} \vert \tau - \sigma \vert ,
\end{equation}
for all $\tau , \sigma \in \tau_{0} + [ - 1 , 1 ]$ and $\omega , \widetilde{\omega} \in \CH_{\rm tang}^{- \infty}$.

For an operator $T$ with continuous kernel $0 \leq T ( \omega , \widetilde{\omega} ) \in L^{\infty}$, we have
\begin{equation*}
\vert T (f) ( \omega ) \vert \leq T ( \vert f \vert ) ( \omega ) ,
\end{equation*}
for all $f \in L^{\infty}$ and $\omega \in \CH_{\rm tang}^{- \infty}$. This implies
\begin{equation} \label{c6}
\Vert T \Vert_{L^{\infty} \to L^{\infty}} = \sup_{\Vert f \Vert_{L^{\infty}} = 1 , \ f \geq 0} \Vert T (f) \Vert_{L^{\infty}} .
\end{equation}
Thus, let $f \in L^{\infty} ( \CH_{\rm tang}^{- \infty} )$ be such that $\Vert f \Vert_{L^{\infty}} = 1$ and $f \geq 0$.

\begin{lemma}\sl  \label{b98}
We have
\begin{equation*}
\CT_{0}^{k} ( \tau ) (f) ( \omega ) \leq \big( 1 + C_{1} C_{2} \vert \tau - \sigma \vert \big)^{k} \CT_{0}^{k} ( \sigma ) (f) ( \omega ) ,
\end{equation*}
for all $k \in \N$, $\tau , \sigma \in \tau_{0} + [ - 1 , 1 ]$ and $\omega \in \CH_{\rm tang}^{- \infty}$.
\end{lemma}

\begin{proof}
We will prove this estimate by induction over $k$. It is clearly satisfied for $k = 0$. Assume now that it holds true for some $k \in \N$. From \eqref{b99}, we can write
\begin{equation} \label{c3}
\CT_{0}^{k + 1} ( \tau ) (f) ( \omega ) = \CT_{0} ( \sigma ) \CT_{0}^{k} ( \tau ) (f) ( \omega ) + \CR ( \tau , \sigma ) \CT_{0}^{k} ( \tau ) (f) ( \omega ) .
\end{equation}
Since $\CT_{0} ( \sigma ) ( \omega , \widetilde{\omega} ) \geq 0$, the induction assumption implies
\begin{equation} \label{c4}
\CT_{0} ( \sigma ) \CT_{0}^{k} ( \tau ) (f) ( \omega ) \leq \big( 1 + C_{1} C_{2}  \vert \tau - \sigma \vert \big)^{k} \CT_{0}^{k + 1} ( \sigma ) (f) ( \omega ) .
\end{equation}
On the other hand, using \eqref{c1} and $\CT_{0}^{k} ( \tau ) (f) ( \omega ) \geq 0$, we obtain
\begin{equation*}
\CR ( \tau , \sigma ) \CT_{0}^{k} ( \tau ) (f) ( \omega ) \leq C_{2} \vert \tau - \sigma \vert \int \CT_{0}^{k} ( \tau ) (f) ( \widetilde{\omega} ) \, d \widetilde{\omega} .
\end{equation*}
Combining with \eqref{c2} and \eqref{c4}, it yields
\begin{align}
\CR ( \tau , \sigma ) \CT_{0}^{k} ( \tau ) (f) ( \omega ) &\leq C_{1} C_{2} \vert \tau - \sigma \vert \int \CT_{0} ( \sigma ) ( \omega , \widetilde{\omega} ) \CT_{0}^{k} ( \tau ) (f) ( \widetilde{\omega} ) \, d \widetilde{\omega}    \nonumber \\
&= C_{1} C_{2} \vert \tau - \sigma \vert \CT_{0} ( \sigma ) \CT_{0}^{k} ( \tau ) (f) ( \omega )  \nonumber \\
&\leq C_{1} C_{2} \vert \tau - \sigma \vert \big( 1 + C_{1} C_{2}  \vert \tau - \sigma \vert \big)^{k} \CT_{0}^{k + 1} ( \sigma ) (f) ( \omega ) .    \label{c5}
\end{align}
Thus, \eqref{c3}, \eqref{c4} and \eqref{c5} give the required estimate for $k + 1$ and then for all $k \in \N$.
\end{proof}

From \eqref{c6} and Lemma \ref{b98}, we obtain
\begin{equation*}
\big\Vert \CT_{0}^{k} ( \tau ) \big\Vert \leq \big( 1 + C_{1} C_{2} \vert \tau - \sigma \vert \big)^{k} \big\Vert \CT_{0}^{k} ( \sigma ) \big\Vert ,
\end{equation*}
for all $k \in \N$, $\tau , \sigma \in \tau_{0} + [ - 1 , 1 ]$. Thus, we deduce
\begin{equation} \label{c10}
\big( 1 + C_{1} C_{2} \vert \tau - \tau_{0} \vert \big)^{- 1} \big\Vert \CT_{0}^{k} ( \tau_{0} ) \big\Vert^{\frac{1}{k}} \leq \big\Vert \CT_{0}^{k} ( \tau ) \big\Vert^{\frac{1}{k}} \leq \big( 1 + C_{1} C_{2} \vert \tau - \tau_{0} \vert \big) \big\Vert \CT_{0}^{k} ( \tau_{0} ) \big\Vert^{\frac{1}{k}} ,
\end{equation}
for all $k \in \N$ and $\vert \tau - \tau_{0} \vert \leq 1$. Taking the limit $k \to + \infty$, we get
\begin{equation*}
\big( 1 + C_{1} C_{2} \vert \tau - \tau_{0} \vert \big)^{- 1} \CA_{0} ( \tau_{0} ) \leq \CA_{0} ( \tau ) \leq \big( 1 + C_{1} C_{2} \vert \tau - \tau_{0} \vert \big) \CA_{0} ( \tau_{0} ) ,
\end{equation*}
and then

\begin{proposition}\sl \label{c7}
The function $\tau \longmapsto \CA_{0} ( \tau )$ is continuous on $\R$.
\end{proposition}

Moreover, we have the following result used in Section \ref{s52}.

\begin{proposition}\sl \label{c8}
Let $\rho > 1$ and $K \subset \R$ be a compact set. Then, there exists an integer $k_{0} \geq 1$ such that
\begin{equation*}
\big\Vert \CT_{0}^{k} ( \tau ) \big\Vert^{\frac{1}{k}} \leq \rho \CA_{0} ( \tau ) ,
\end{equation*}
for all $k \geq k_{0}$ and $\tau \in K$.
\end{proposition}

\begin{proof}
By compactness, it is enough to obtain this property near any point $\tau_{0} \in \R$. Let $\delta > 0$. Then, there exists an integer $k_{0} \geq 1$ such that
\begin{equation*}
\Big\vert \big\Vert \CT_{0}^{k} ( \tau_{0} ) \big\Vert^{\frac{1}{k}} - \CA_{0} ( \tau_{0} ) \Big\vert \leq \delta ,
\end{equation*}
for all $k \geq k_{0}$. On the other hand, from \eqref{c10}, there exists $c_{1} > 0$ such that
\begin{equation*}
\Big\vert \big\Vert \CT_{0}^{k} ( \tau ) \big\Vert^{\frac{1}{k}} - \big\Vert \CT_{0}^{k} ( \tau_{0} ) \big\Vert^{\frac{1}{k}} \Big\vert \leq \delta ,
\end{equation*}
for all $\vert \tau - \tau_{0} \vert \leq c_{1}$ and $k \geq 1$. Moreover, Proposition \ref{c7} provides $c_{2} > 0$ such that
\begin{equation} \label{c11}
\big\vert \CA_{0} ( \tau ) - \CA_{0} ( \tau_{0} ) \big\vert \leq \delta ,
\end{equation}
for all $\vert \tau - \tau_{0} \vert \leq c_{2}$. The three last estimates give
\begin{equation} \label{c9}
\big\Vert \CT_{0}^{k} ( \tau ) \big\Vert^{\frac{1}{k}} \leq \CA_{0} ( \tau ) + 3 \delta ,
\end{equation}
for all $\vert \tau - \tau_{0} \vert \leq \min ( c_{1} , c_{2} )$ and $k \geq k_{0}$.

From \eqref{c11} and $\CA_{0} ( \tau_{0} ) > 0$ (see Remark \ref{b69} $i)$ and the beginning of this section), we also have
\begin{equation*}
\CA_{0} ( \tau ) \geq \CA_{0} ( \tau_{0} ) - \delta \geq \CA_{0} ( \tau_{0} ) / 2 > 0,
\end{equation*}
for all $\delta$ small enough and $\vert \tau - \tau_{0} \vert \leq \min ( c_{1} , c_{2} )$. Thus, taking $\delta$ sufficiently small, \eqref{c9} yields
\begin{equation*}
\big\Vert \CT_{0}^{k} ( \tau ) \big\Vert^{\frac{1}{k}} \leq \rho \CA_{0} ( \tau ) ,
\end{equation*}
for all $\vert \tau - \tau_{0} \vert \leq \min ( c_{1} , c_{2} )$ and $k \geq k_{0}$.
\end{proof}

We finish this part by proving \eqref{i57}, that is

\begin{lemma}\sl \label{i59}
Under the assumptions of Section \ref{s26}, we have
\begin{equation*}
\CA ( \tau , h ) \leq \CA_{0} ( \tau ) ,
\end{equation*}
for all $\tau , h$.
\end{lemma}

\begin{proof}
For $k \in \N \setminus \{ 0 \}$, the kernel of $\CT_{\bullet}^{k}$ is given by
\begin{equation*}
\CT_{\bullet}^{k} ( \omega , \widetilde{\omega} ) = \int_{( \CH^{- \infty}_{\rm tang} )^{k - 1}} \CT_{\bullet} ( \omega , \omega_{1} ) \cdots \CT_{\bullet} ( \omega_{k - 1} , \widetilde{\omega} ) \, d \omega_{1} \cdots d \omega_{k - 1} .
\end{equation*}
Using \eqref{i56}, it yields
\begin{equation*}
\big\vert \CT^{k} ( \omega , \widetilde{\omega} ) \big\vert \leq \CT^{k}_{0} ( \omega , \widetilde{\omega} ) ,
\end{equation*}
for all $k \in \N \setminus \{ 0 \}$ and $\omega , \widetilde{\omega} \in \CH^{- \infty}_{\rm tang}$. It implies
\begin{equation}
\big\vert \CT^{k} u ( \omega ) \big\vert \leq \int \big\vert \CT^{k} ( \omega , \widetilde{\omega} ) u ( \widetilde{\omega} ) \big\vert \, d \widetilde{\omega} \leq \CT^{k}_{0} \vert u \vert ( \omega ) ,
\end{equation}
for all $u \in L^{2} ( \CH^{- \infty}_{\rm tang} )$ and $\omega \in \CH^{- \infty}_{\rm tang}$. In particular,
\begin{equation*}
\big\Vert \CT^{k} u \big\Vert \leq \big\Vert \CT^{k}_{0} \vert u \vert \big\Vert \leq \big\Vert \CT^{k}_{0} \big\Vert \Vert \vert u \vert \Vert = \big\Vert \CT^{k}_{0} \big\Vert \Vert u \Vert .
\end{equation*}
This gives $\Vert \CT^{k} \Vert \leq \Vert \CT^{k}_{0} \Vert$ and the lemma follows from the definition of the spectral radius.
\end{proof}

\section{Distorted and truncated estimates} \label{s41}

In this part, we compare distorted and truncated quantities. The first result, which is essentially in the folklore of the theory of resonances, shows the equivalence of the norms of the resolvent.

\begin{proposition}\sl \label{j66}
Assume \ref{h1} and let $E_{0} , R > 0$, $\theta = h \vert \ln h \vert$ and $\chi \in C^{\infty}_{0} ( \R^{n} ; [ 0 , 1 ])$ be such that $\chi = 1$ on a sufficiently large neighborhood of $0$. We also suppose that the distortion \eqref{j70} occurs outside the support of $\chi$. Then, there exists $C > 0$ such that
\begin{equation*}
\big\Vert \chi ( P - z )^{-1} \chi \big\Vert \leq \big\Vert ( P_{\theta} - z )^{-1} \big\Vert \leq h^{- C} \big\Vert \chi ( P - z )^{-1} \chi \big\Vert ,
\end{equation*}
for $h$ small enough and $z \in B ( E_{0} , R h )$ that is not a resonance.
\end{proposition}

In fact, if the distortion function $F$ is well-chosen (that is $H_{p} ( F ( x ) \cdot \xi ) \geq 0$ on the entire energy surface $p^{- 1} ( E_{0} )$) or if the distortion angle satisfies $\theta = \CO ( h )$, it could be possible to replace $h^{- C}$ by $C$ in the previous inequality.

\begin{proof}
Using \eqref{c12} and $\vert \chi \vert \leq 1$, we immediately obtain
\begin{equation} \label{j67}
\big\Vert \chi ( P - z )^{-1} \chi \big\Vert \leq \big\Vert ( P_{\theta} - z )^{-1} \big\Vert .
\end{equation}
For the second inequality, we construct an operator $Q$ which coincides with $P$ outside a compact set and which is non-trapping at energy $E_{0}$. This idea comes back to Robert and Tamura \cite{RoTa87_01}. Let $\varphi \in C^{\infty}_{0} ( \R^{n} ; [ 0 , 1 ] )$ be such that $\varphi ( x ) = 1$ if $\vert x \vert \leq 1$. For $r > 0$, we set
\begin{equation}
Q = - h^{2} \Delta + V_{r} ( x ) \qquad \text{with} \qquad V_{r} ( x ) = V ( x ) \Big( 1 - \varphi \Big( \frac{x}{r} \Big) \Big) .
\end{equation}
Its symbol is $q ( x , \xi ) = \xi^{2} + V_{r} ( x )$. From \ref{h1}, we have $\xi^{2} > E_{0} / 2$ on the energy surface $q^{- 1} ( E_{0} )$ for $r$ large enough. Moreover, a direct computation gives
\begin{equation*}
\{ q , x \cdot \xi \} = 2 \xi^{2} - x \cdot \nabla V ( x ) \Big( 1 - \varphi \Big( \frac{x}{r} \Big) \Big) + \frac{x}{r} \cdot \nabla \varphi \Big( \frac{x}{r} \Big) V ( x ) .
\end{equation*}
The Cauchy formula and \ref{h1} imply $x \cdot \nabla V ( x ) = o_{x \to + \infty} ( 1 )$. Thus, the previous equation becomes
\begin{equation} \label{l1}
\{ q , x \cdot \xi \} \geq E_{0} + o_{x \to + \infty} ( 1 ) \one_{\vert x \vert \geq r} \geq E_{0} / 2 ,
\end{equation}
for $r$ large enough and $( x , \xi ) \in q^{- 1} ( E_{0} )$. In particular, $x \cdot \xi$ is an escape function for $q$ at energy $E_{0}$ and then $Q$ is non-trapping at energy $E_{0}$. In the sequel, we fix $r$ large enough such that the previous property holds true.

We now assume that $\chi = 1$ near the support of $\varphi ( \cdot / r )$. In particular, $W = W \chi$ where $W = V - V_{r} \in C^{\infty}_{0} ( \R^{n} )$. Let $Q_{\theta}$ be the operator $Q$ distorted as in \eqref{a6}. From the resolvent estimate in the non-trapping case proved by Martinez \cite{Ma02_01}, there exists $M > 0$ such that
\begin{equation} \label{l3}
\big\Vert ( Q_{\theta} - z )^{-1} \big\Vert \lesssim h^{- M} ,
\end{equation}
uniformly for $h$ small enough and $z \in B ( E_{0} , R h )$. Using the resolvent identity and \eqref{c12}, we can write
\begin{align}
( P_{\theta} - z )^{-1} ={}& ( Q_{\theta} - z )^{-1} - ( P_{\theta} - z )^{-1} W ( Q_{\theta} - z )^{-1}   \nonumber \\
={}& ( Q_{\theta} - z )^{-1} - ( Q_{\theta} - z )^{-1} W ( Q_{\theta} - z )^{-1}    \nonumber \\
&\qquad \qquad \qquad + ( Q_{\theta} - z )^{-1} W ( P_{\theta} - z )^{-1} W ( Q_{\theta} - z )^{-1}   \nonumber \\
={}& ( Q_{\theta} - z )^{-1} - ( Q_{\theta} - z )^{-1} W ( Q_{\theta} - z )^{-1}    \nonumber \\
&\qquad \qquad \qquad + ( Q_{\theta} - z )^{-1} W \chi ( P - z )^{-1} \chi W ( Q_{\theta} - z )^{-1} ,  \label{j71}
\end{align}
for $z$ not a resonance. Combining with \eqref{l3}, this yields
\begin{equation} \label{j68}
\big\Vert ( P_{\theta} - z )^{-1} \big\Vert \lesssim h^{- 2 M} + h^{- 2 M} \big\Vert \chi ( P - z )^{-1} \chi \big\Vert ,
\end{equation}
uniformly for $z \in B ( E_{0} , R h )$ not a resonance. On the other hand, from \cite[Proposition 1.5]{BoPe13_01} of Petkov and the first author, we have
\begin{equation*}
\big\Vert \chi ( P - z )^{-1} \chi \big\Vert \gtrsim h^{- 1} ,
\end{equation*}
uniformly for $z \in B ( E_{0} , R h )$ not a resonance. Then, \eqref{j68} gives
\begin{equation} \label{j69}
\big\Vert ( P_{\theta} - z )^{-1} \big\Vert \leq h^{- C} \big\Vert \chi ( P - z )^{-1} \chi \big\Vert ,
\end{equation}
with $C = 2 M + 1$ for $h$ small enough and $z \in B ( E_{0} , R h )$ not a resonance. Summing up, the proposition follows from \eqref{j67} and \eqref{j69}.

If in addition the distortion function $F$ is well-chosen or if $\theta = \CO ( h )$ (see below the proposition), one could adapt the proof of \cite[Theorem 2.1]{Ma02_01} in order to obtain $M =1$ in \eqref{l3}. Moreover, using the commutator argument of \cite[Pages 437-438]{RoTa87_01} instead of the resolvent identity \eqref{j71}, one may gain a factor $h^{2}$ in \eqref{j68}. Thus, $h^{- C}$ could be replaced by $C$.
\end{proof}

We now prove the equivalence of the norms of the resonant states. Since Remark \ref{j48} is analogous to Proposition \ref{j66}, we use similar arguments.

\begin{proof}[Proof of Remark \ref{j48}]
Let $Q$ be the non-trapping operator build in the proof of Proposition \ref{j66}. We consider $\psi \in C^{\infty}_{0} ( \R^{n} )$ such that $\varphi ( \cdot / r ) \prec \psi \prec \chi$. By construction,
\begin{equation*}
( Q_{\theta} - z ) ( 1 - \psi ) = ( P_{\theta} - z ) ( 1 - \psi ) = ( 1 - \psi ) ( P_{\theta} - z ) - [ P , \psi ] .
\end{equation*}
In particular, if $u$ is a resonant state associated to $z$, this gives
\begin{equation} \label{l2}
( Q_{\theta} - z ) ( 1 - \psi ) u = - [ P , \psi ] u = - [ P , \psi ] \chi u .
\end{equation}
Using \eqref{l3} and the ellipticity of $Q_{\theta}$, we deduce
\begin{equation*}
\big\Vert ( Q_{\theta} - z )^{-1} \big\Vert_{H^{-2}_{h} \to L^{2}} \lesssim h^{- M} ,
\end{equation*}
uniformly for $h$ small enough and $z \in B ( E_{0} , R h )$. Then, \eqref{l2} gives
\begin{align}
\Vert ( 1 - \chi ) u \Vert &\leq \Vert ( 1 - \psi ) u \Vert   \nonumber \\
&= \big\Vert ( Q_{\theta} - z )^{-1} [ P , \psi ] \chi u \big\Vert  \nonumber \\
&\leq \big\Vert ( Q_{\theta} - z )^{-1} \big\Vert_{H^{- 2}_{h} \to L^{2}} \big\Vert [ P , \psi ] \big\Vert_{L^{2} \to H^{- 1}_{h}} \Vert \chi u \Vert  \nonumber \\
&\lesssim h^{1 - M} \Vert \chi u \Vert , \label{l4}
\end{align}
which implies the remark.

As in the end of the proof of Proposition \ref{j66}, one could obtain $M = 1$ in \eqref{l2} and \eqref{l4} under the assumptions below Remark \ref{j48} ($F$ well-chosen or $\theta = \CO ( h )$). Thus, $h^{- C}$ may be replaced by $C$ in the statement of the result.
\end{proof}

\section{Semiclassical principle} \label{b73}

This section is devoted to the proof of Corollary \ref{a7} and Corollary \ref{a8}. We follow Tang and Zworski \cite[Lemma 2]{TaZw98_01} in the general case, Nakamura, Stefanov and Zworski \cite[Proposition 3.1]{NaStZw03_01} in the non-trapping case and Burq \cite[Lemma 4.7]{Bu04_01} for hyperbolic trapping. The next lemma is a simple adaptation of these references.

\begin{lemma}\sl \label{b74}
Let $f ( z , \varepsilon )$ be holomorphic in $[ - 2 \varepsilon , 2 \varepsilon ] + i [ - \varepsilon , \varepsilon ]$ such that
\begin{equation*}
\vert f ( z , \varepsilon ) \vert \leq
\left\{ \begin{aligned}
&\varepsilon^{- C}  &&\text{for } z \in [ - 2 \varepsilon , 2 \varepsilon ] + i [ - \varepsilon , \varepsilon ] ,   \\
&\frac{1}{\vert \im z \vert} &&\text{for } \im z > 0,
\end{aligned} \right.
\end{equation*}
for some $C > 0$ and $\varepsilon$ small enough. Then, there exists $N > 0$, depending only on $C$, such that
\begin{equation*}
\forall z \in i [ - \varepsilon , \varepsilon ] , \qquad \vert f ( z , \varepsilon ) \vert \leq N \frac{\vert \ln \varepsilon \vert}{\varepsilon} e^{N \vert \im z \vert \vert \ln \varepsilon \vert / \varepsilon} ,
\end{equation*}
for $\varepsilon$ small enough.
\end{lemma}

\begin{proof}
We consider
\begin{equation*}
g ( z , \varepsilon ) = f ( z , \varepsilon ) e^{- C z^{2} \vert \ln \varepsilon \vert / \varepsilon^{2}} e^{- 2 C i z \vert \ln \varepsilon \vert /\varepsilon} .
\end{equation*}
For $z \in \pm 2 \varepsilon + i [ - \varepsilon , \varepsilon ]$, we have
\begin{align*}
\vert g ( z , \varepsilon ) \vert &\leq \vert f ( z , \varepsilon ) \vert e^{- C ( 4 \varepsilon^{2} - ( \im z )^{2} ) \vert \ln \varepsilon \vert / \varepsilon^{2}} e^{2 C \im z \vert \ln \varepsilon \vert /\varepsilon}   \\
&\leq \varepsilon^{- C} e^{- 3 C \vert \ln \varepsilon \vert} e^{2 C \vert \ln \varepsilon \vert} = 1 .
\end{align*}
The same way, for $z \in [ - 2 \varepsilon , 2 \varepsilon ] - i \varepsilon$, we get
\begin{align*}
\vert g ( z , \varepsilon ) \vert &\leq \vert f ( z , \varepsilon ) \vert e^{- C ( ( \re z )^{2} - \varepsilon^{2} ) \vert \ln \varepsilon \vert / \varepsilon^{2}} e^{- 2 C \vert \ln \varepsilon \vert}   \\
&\leq \varepsilon^{- C} e^{C \vert \ln \varepsilon \vert} e^{- 2 C \vert \ln \varepsilon \vert} = 1 .
\end{align*}
Eventually, for $z \in [ - 2 \varepsilon , 2 \varepsilon ] + i \varepsilon / \vert \ln \varepsilon \vert$, we can estimate
\begin{align*}
\vert g ( z , \varepsilon ) \vert &\leq \vert f ( z , \varepsilon ) \vert e^{- C ( ( \re z )^{2} - \varepsilon^{2} / \vert \ln \varepsilon \vert^{2} ) \vert \ln \varepsilon \vert / \varepsilon^{2}} e^{2 C}   \\
&\leq \frac{\vert \ln \varepsilon \vert}{\varepsilon} e^{C / \vert \ln \varepsilon \vert} e^{2 C} \leq N \frac{\vert \ln \varepsilon \vert}{\varepsilon} .
\end{align*}
From the maximum principle for holomorphic functions, we deduce that $\vert g ( z , \varepsilon ) \vert \leq N \vert \ln \varepsilon \vert / \varepsilon$ for all $z \in [ - 2 \varepsilon , 2 \varepsilon ] + i [ - \varepsilon , \varepsilon / \vert \ln \varepsilon \vert ]$. In particular,
\begin{equation*}
\vert f ( z , \varepsilon ) \vert \leq N \frac{\vert \ln \varepsilon \vert}{\varepsilon} e^{C ( ( \re z )^{2} - ( \im z )^{2} ) \vert \ln \varepsilon \vert / \varepsilon^{2}} e^{- 2 C \im z \vert \ln \varepsilon \vert /\varepsilon} \leq N \frac{\vert \ln \varepsilon \vert}{\varepsilon} e^{- 2 C \im z \vert \ln \varepsilon \vert /\varepsilon} ,
\end{equation*}
for all $z \in i [- \varepsilon , \varepsilon / \vert \ln \varepsilon \vert ]$. On the other hand, since the estimate is clear for $z \in i [ \varepsilon / \vert \ln \varepsilon \vert , \varepsilon ]$, the lemma follows.
\end{proof}

To prove Corollary \ref{a7} from Theorem \ref{a1}, it is enough to take $\varepsilon = \nu h$ with $\nu > 0$ small enough in Lemma \ref{b74}. To prove Corollary \ref{a8} from Theorem \ref{a2}, we choose $\varepsilon = \nu h \vert \ln h \vert^{- 1}$ with $\nu > 0$ small enough and use
\begin{equation*}
\frac{\vert \ln \varepsilon \vert}{\varepsilon} = \frac{\vert \ln  h \vert^{2}}{\nu h} ( 1 + o_{h \to 0} ( 1) ) .
\end{equation*}
Eventually, remark that Lemma \ref{b74} also gives the usual estimate
\begin{equation*}
\big\Vert \chi ( P -z )^{-1} \chi \big\Vert \lesssim\frac{1}{h} e^{N \vert \im z \vert / h} ,
\end{equation*}
in the non-trapping case (such an estimate follows from Proposition 3.1 of \cite{NaStZw03_01}). Indeed, to obtain this bound, it is enough to apply the result of Martinez \cite{Ma02_01} which provides a resonance free region of size $h \vert \ln h \vert$ with a polynomial estimate of the resolvent, to take $\varepsilon = \nu h \vert \ln h \vert$ with $\nu > 0$ small enough in Lemma \ref{b74} and to use
\begin{equation*}
\frac{\vert \ln \varepsilon \vert}{\varepsilon} = \frac{1}{\nu h} ( 1 + o_{h \to 0} ( 1) ) .
\end{equation*}
In other words, all these estimates (non-trapping case, hyperbolic trapping, Corollary \ref{a8}) follow from the same idea, provided one has a polynomial estimate of the resolvent at hand.

\newpage

\section*{Notations} \label{s34}
\begin{tabbing}
\sl{Notation} \hskip 30pt \= \sl{Description/definition} \hskip 180pt \= \sl{Reference} \\
\\
$\prec$ \> $f \prec g$ means that $g = 1$ near the support of $f$ \> \\
$\sim$ \> equivalent \> \\
$\simeq$  \> asymptotic expansion \> \\
$\approx$ \> close, has no precise meaning \> \\
$\oplus$ \> decomposition of a path \> Section \ref{s50} \\
$p_{1} \cup p_{2}$ \> concatenation, the path $p_{1}$ followed by $p_{2}$ \> Section \ref{s47} \\
$\CA_{\bullet}$ \> spectral radius of $\CT_{\bullet}$ \> \eqref{b71}, \eqref{j50} \\
$B ( c , r )$ \> open ball of center $c$ and radius $r$, $\{ x ; \ \vert x - c \vert < r \}$ \>  \\
$D_{0}$ \> minimal damping \> \eqref{m1}  \\
$\varphi_{\pm}$ \> generating function of $\Lambda_{\pm}^{0}$ \>  \eqref{a57}  \\
$\varphi_{+}^{1}$ \> generating function of $\Lambda_{+}^{1}$ \> Section \ref{s52} \\
$g_{\pm} ( \rho )$ \> asymptotic direction of the characteristic \> \eqref{d2} \\
\> trajectory starting from $\rho \in \Lambda_{\pm}$ \> \\
$\gamma_{k}$ \> homoclinic trajectory, $\gamma_{k} ( t ) = ( x_{k} ( t ) , \xi_{k} ( t ) )$  \> Section \ref{s61} \\
$\SG = ( \SV , \SE )$ \> graph whose $\SV$ is the set of fixed points and $\SE$ is \> Section \ref{s45} \\
\> the set of homoclinic/heteroclinic trajectories \> \\
$\SG^{*} = ( \SE , \SC )$ \> adjoint graph of $\SG$ \> Section \ref{s50} \\
$\Gamma_{0} ( h )$ \> leading term in the asymptotic of \> \eqref{m80} \\
\> the resonances generated by a barrier-top \> \\
$\Gamma ( h )$ \> exceptional set where the microlocal Cauchy \> Section \ref{m81} \\
\> problem at the barrier-top is not well-posed \> \\
$H_{p}$ \> Hamiltonian vector field \> \eqref{d27} \\
$\CH$ \> homoclinic set \> \ref{h3} \\
$\CH_{\rm tang}$ \> set of homoclinic trajectories along \> Section \ref{s77} \\
\> which $\Lambda_{\pm}$ are tangent \> \\
$\CH_{\rm trans}$ \> $\CH \setminus \CH_{\rm tang}$ \> Section \ref{s53} \\
$\CH_{\rm tang}^{\pm \infty}$ \> normalized asymptotic directions of $\CH_{\rm tang}$ \> Section \ref{s77} \\
\> as $t \to \mp \infty$ \> \\
$\CH_{\bullet , \pm}^{\varepsilon}$ \> $\CH_{\bullet} \cap \Lambda_{\pm}^{0} \cap \{ \vert x \vert = \varepsilon \}$ \> Sections \ref{s4} and \ref{s51} \\
$\CI ( \Lambda , m )$ \> space of semiclassical Lagrangian distributions \> Definition \ref{a30} \\
\> with manifold $\Lambda$ and symbol of order $m$ \> \\
$K$ \> number of transversal homoclinic trajectories \> Section \ref{s61} \\
$K ( E )$ \> trapped set at energy $E$ \> Section \ref{s2} \\
$\lambda_{j}$ \> $- \lambda_{j}^{2} / 2 $ are the eigenvalues of $\Hess V ( 0 )$ \> \ref{h2} \\
$\Lambda_{\pm}$ \> global incoming/outgoing Lagrangian manifold \> \eqref{j49} \\
$\Lambda_{\pm}^{0}$ \> local incoming/outgoing Lagrangian manifold near $0$ \> Section \ref{s2} \\
$\Lambda_{+}^{1}$ \> manifold $\Lambda_{+}^{0}$ after a turn along $\CH$ \> \eqref{a38} \\
$n$ \> dimension of the space \> Section \ref{s2} \\
$\N$ \> set of natural numbers, $\{ 0 , 1 , 2 , \ldots \}$ \> \\
$o_{b \to c}^{a} ( 1 )$ \> function which goes to $0$ as $b \to c$ for $a$ fixed \> \\
$\CO_{a} ( 1 )$ \> bounded function for $a$ fixed \> \\
$\pi_{x}$ \> base space projection, $\pi_{x} ( x , \xi ) = x$ \> \\
$P_{\theta}$ \> distorted operator of angle $\theta$ \> \eqref{a6} \\
$\CQ ( z , h )$\> $K \times K$ quantization matrix in Section \ref{s6} \> \eqref{d4} \\
$\widetilde{\CQ} , \widehat{\CQ}$ \> rewritings of $\CQ$ \> \eqref{d10}, \eqref{e18} \\
$\CQ^{2} ( z , h ) $\> quantization matrix for the second set of resonances \> \eqref{g91} \\
$\SQ ( z , h )$ \> global quantization matrix for multiple barriers \> \eqref{m2} \\
$\res ( P )$ \> set of the resonances of $P$ \> Section \ref{s2} \\
$\res_{0} ( P )$ \> set of the pseudo-resonances of $P$ which \> Definitions \ref{d1},  \\
\> depends on the nature of the trapped set \> \ref{g26}, \ref{i51} or \ref{j82} \\
$\res_{0}^{2} ( P )$ \> set of the pseudo-resonances of the second kind of $P$ \> Definition \ref{g85} \\
$\res_{\infty} ( P )$ \> set of the pseudo-resonances at infinite order of $P$ \> Definition \ref{e97} \\
$\sigma$\> rescaled spectral parameter, $\frac{z - E_{0}}{h}$ \> \eqref{d92} \\
$\spe ( \cdot )$ \> spectrum \> \\
$\spr ( \cdot )$ \> spectral radius \> \eqref{b71} \\
$S_{\pm}^{\varepsilon}$ \> $\Lambda_{\pm}^{0} \cap \{ \vert x \vert = \varepsilon \}$ \> Section \ref{s4} \\
$S ( m )$ \> space of symbols of order $m$ \> Section \ref{s35} \\
$S ( z , h )$ \> $\sum \frac{\lambda_{j}}{2}- i \frac{z - E_{0}}{h}$ \> \eqref{d5} or \eqref{n13} \\
$\tau$ \> $\re z = E_{0} + \tau h + o ( h )$ \> Section \ref{s61} \\
$T^{*} \R^{n}$\> cotangent bundle, $T^{*} \R^{n} = \R^{2 n}$ \> \\
$\CT$ \> quantization operator in Section \ref{s26} \> \eqref{c19} \\
$\widetilde{\CT} , \widehat{\CT} $\> rewritings of $\CT$ \> Section \ref{s53} \\
$\CT_{0}$ \> operator whose kernel is the modulus \> \eqref{b70} \\
\> of that of $\CT$ \> \\
$z_{q , k} ( \tau )$ \> asymptotic of the pseudo-resonances  \> \eqref{d95}
\end{tabbing}

\newpage

\setcounter{tocdepth}{2}
\tableofcontents

\newpage


\bibliographystyle{amsplain}
\providecommand{\bysame}{\leavevmode\hbox to3em{\hrulefill}\thinspace}
\providecommand{\MR}{\relax\ifhmode\unskip\space\fi MR }
\providecommand{\MRhref}[2]{%
  \href{http://www.ams.org/mathscinet-getitem?mr=#1}{#2}
}
\providecommand{\href}[2]{#2}

\end{document}